\documentclass[reqno]{amsart}
\usepackage[a4paper,total={170mm,257mm},left=20mm,top=20mm]{geometry}
\usepackage{amsmath,amssymb,amsthm,mathrsfs,booktabs,microtype,mathtools,array}
\usepackage[T1]{fontenc}
\usepackage{lmodern}
\usepackage{needspace}
\usepackage{tikz}
\usetikzlibrary{decorations.pathreplacing}
\usepackage[hypertexnames=false]{hyperref}
\hypersetup{hidelinks,bookmarksdepth=2,
 pdftitle={Onsager's Conjecture for Ideal MHD},
 pdfauthor={Matteo Giardi and L\'aszl\'o Sz\'ekelyhidi Jr.}}
\newtheorem{theorem}{Theorem}[section]
\newtheorem{proposition}[theorem]{Proposition}
\newtheorem{lemma}[theorem]{Lemma}
\newtheorem{corollary}[theorem]{Corollary}
\theoremstyle{definition}
\newtheorem{assumption}[theorem]{Assumption}
\newtheorem{definition}[theorem]{Definition}
\theoremstyle{remark}
\newtheorem{remark}[theorem]{Remark}
\newcommand{\sym}{\operatorname{sym}}
\newcommand{\Cov}{\operatorname{Cov}}
\newcommand{\R}{\mathbb R}
\newcommand{\rgood}{r_{\mathrm{good}}}
\newcommand{\dstar}{d_{\ast}}
\newcommand{\rprep}{r_{\mathrm{prep}}}
\newcommand{\rcut}{r_{\mathrm{cut}}}
\newcommand{\Nc}{N_{\mathrm c}}
\newcommand{\rfg}{r_{\mathrm{fg}}}
\newcommand{\gres}{\gamma_{\mathrm{res}}}
\newcommand{\dd}{\mathrm d}
\newcommand{\DD}{\mathrm D}
\newcommand{\ddiv}{\operatorname{div}}
\newcommand{\curl}{\operatorname{curl}}
\newcommand{\IId}{\mathrm{Id}}
\DeclareMathOperator{\supp}{supp}
\newcommand{\cn}{\cdot\nabla}
\numberwithin{equation}{section}
\makeatletter
\renewcommand\subsection{\@startsection{subsection}{2}%
 \z@{-\baselineskip}{-0.3em}%
 {\normalfont\bfseries}}
\renewcommand\subsubsection{\@startsection{subsubsection}{3}%
 \z@{-\baselineskip}{-0.3em}%
 {\normalfont\bfseries}}
\makeatother

\title[Onsager's conjecture for ideal MHD]{Onsager's Conjecture for Ideal Magnetohydrodynamics}
\author{Matteo Giardi}
\thanks{MG was supported by the International Max Planck Research School
Mathematics in the Sciences (IMPRS MiS)}
\author{L\'aszl\'o Sz\'ekelyhidi Jr.}
\thanks{LSz gratefully acknowledges the support of the Deutsche
Forschungsgemeinschaft (DFG, German Research Foundation) through
GZ SZ 325/2-1.}
\date{}

\begin{document}
\raggedbottom

\begin{abstract}
For any $0\le\gamma<1/3$ we construct weak solutions $(v,B,p)$ of the
ideal MHD equations with $v,B\in C^\gamma(\mathbb T^3\times\mathbb R)$,
which conserve neither the total energy nor the cross-helicity and have
nontrivial magnetic helicity. We also establish anisotropic H\"older bounds with distinct velocity
and magnetic exponents and stronger regularity along every magnetic
field line. The admissible exponents connect
the Goldreich--Sridhar \cite{Goldreich1995} spatial and parallel pair $1/3,1/2$
with the $1/4$ spatial scaling of the
Iroshnikov--Kraichnan weak-turbulence picture
\cite{Iroshnikov1963,Kraichnan1965}. In the spirit of
Arnold's formulation of ideal hydrodynamics, a solution is
regarded as a path of volume-preserving diffeomorphisms; the proof is
then based on the interplay between classical convex integration
techniques and geometric constructions at the level of the Lie algebra
of this Lie group. Our work substantially extends the recent result of Enciso,
Pe\~nafiel-Tom\'as and Peralta-Salas \cite{EnPePe} and can be used to reprove \cite{GR} of Giri and Radu without a Newton iteration.
\end{abstract}
\maketitle
{
\makeatletter
\providecommand\@dotsep{4.5}
\renewcommand{\l@section}{\@dottedtocline{1}{0em}{2.3em}}
\renewcommand{\l@subsection}{\@dottedtocline{2}{1.5em}{3em}}
\makeatother
\small\tableofcontents
}
\section{Introduction}

We study the ideal incompressible magnetohydrodynamic equations (MHD)
in three dimensions. This is a classical macroscopic model describing
the evolution of electrically conducting fluids such as plasmas and
liquid metals. The MHD system describes the simultaneous evolution of a
velocity field $v$ and a magnetic field $B$, both divergence free. The
evolution of $v$ is described by the Cauchy momentum equation with the
external force given by the Lorentz force induced by $B$;
incompressibility is ensured by a pressure gradient $\nabla p$. In
turn, the evolution of $B$ is described by the induction equation,
which couples the Maxwell--Faraday law with Ohm's law for the electric
field $E$ of a perfect conductor:
\begin{equation}\label{eq:faraday-ohm}
 \partial_tB+\curl E=0,\qquad \ddiv B=0,\qquad E+v\times B=0.
\end{equation}
The resulting system of equations is
\begin{equation}\tag{MHD}\label{eq:mhd-ideal}
 \begin{cases}
 \partial_tv+v\cn v+\nabla p-B\cn B=0,\\
 \partial_tB+\curl(B\times v)=0,\\
 \ddiv v=\ddiv B=0.
 \end{cases}
\end{equation}
Here $p$ denotes the total pressure, the sum of the fluid pressure and
the magnetic pressure $|B|^2/2$. We consider this system in the
periodic setting, in other words on the spatial domain
$\mathbb T^3=(\R/2\pi\mathbb Z)^3$, so that the unknowns are the vector
fields $v,B:\mathbb T^3\times\R\to\R^3$ and the scalar field
$p:\mathbb T^3\times\R\to\R$. In this paper we are concerned with weak
solutions of the MHD system, that is, triples $(v,B,p)$ satisfying
\eqref{eq:mhd-ideal} in the sense of distributions; the test functions
are smooth, periodic in space and compactly supported in time, and
pressures are normalized to have zero spatial mean.

For divergence-free fields one has
\[
 \curl(B\times v)=v\cn B-B\cn v=[v,B],
\]
the Lie bracket of the two vector fields. The induction equation states that $B$ is
transported as a vector field by the flow of $v$: the magnetic field
lines are frozen into the fluid. We use this form of the induction
equation throughout, together with the conventions
\[
 \begin{gathered}
 (\DD v)^{ij}=\partial_jv^i,\qquad (a\otimes b)^{ij}=a^ib^j,
 \qquad (\ddiv R)^i=\partial_jR^{ij},\\
 a\odot b=(a\otimes b+b\otimes a)/2.
 \end{gathered}
\]
Thus $\ddiv(a\otimes b)=b\cn a+(\ddiv b)a$, and the momentum equation reads
\[
 \partial_tv+\ddiv(v\otimes v-B\otimes B)+\nabla p=0.
\]

Our motivation comes from magnetohydrodynamic turbulence and the
search for an analogue of the famous K41 theory describing
hydrodynamic turbulence. In the hydrodynamic case a theory of weak
solutions of the incompressible Euler equations, describing ideal
hydrodynamics, arose from the work of Onsager in 1949 \cite{Onsager}.
This work has received considerable attention over the past decades
under the name of Onsager's conjecture, after its initial revival in
the early 1990s \cite{Eyink1994,CET}. The conjecture was eventually
resolved in the work of P.~Isett \cite{Isett2018}, see also
\cite{BuckmasterDeLellisSzekelyhidiVicol2019}, using the technique of
convex integration. For a comprehensive review of this subject we refer
to \cite{EyinkReview}. Subsequently, the technique has been adapted to
various other fluid-dynamic models, including the ideal MHD setting in
\cite{FLSz2021,BBV,FLSz2024,MiaoNieYe,EnPePe}. For a more in-depth
review of these works, see Subsection~\ref{ssec:literature} below. 

Our main result is the following.

\begin{theorem}[Main Theorem]\label{full:main}
For any $0\le\gamma<1/3$ there exist weak solutions $(v,B,p)$ of
\eqref{eq:mhd-ideal} with $v,B\in C^\gamma(\mathbb T^3\times\R)$ which
conserve neither the total energy nor the cross-helicity. The magnetic
field has zero spatial mean at every time, is nowhere zero, and has
nonzero constant magnetic helicity.
\end{theorem}
In fact, the construction gives more information on the separate regularity of
both fields as well as regularity along the magnetic field lines. Theorem \ref{full:main} is thus a special case of the following more general statement:

\begin{theorem}[Anisotropic regularity]\label{full:fieldline}
For any pair of exponents 
\begin{equation}\label{full:fieldline-exponents}
 0<\gamma_v<\tfrac13,\qquad
 0<\gamma_B<\tfrac12,\qquad 2\gamma_v+\gamma_B<1,
\end{equation}
there exists a weak solution $(v,B,p)$ of
\eqref{eq:mhd-ideal} with 
\[
 v\in C^{\gamma_v}(\mathbb T^3\times\R),\qquad
 B\in C^{\gamma_B}(\mathbb T^3\times\R)
\]
which conserves neither the total energy nor the cross-helicity but has nonzero constant magnetic helicity. The magnetic field has zero spatial mean at every time and is nowhere zero.

Moreover, along every magnetic field line, the solution satisfies
\begin{equation}\label{full:fieldline-holder}
 \begin{aligned}
 |v(t,\Gamma(s))-v(t,\Gamma(r))|
     &\le C|s-r|^{\gamma_v/(1-\gamma_B)},\\
 |B(t,\Gamma(s))-B(t,\Gamma(r))|
     &\le C|s-r|^{\gamma_B/(1-\gamma_B)}
 \end{aligned}
 \qquad (|s-r|\le1),
\end{equation}
with $C$ independent of $t$ and $\Gamma$. 
\end{theorem}

We recall that a magnetic field line at fixed time $t$ is an integral curve of the magnetic field $B(t,\cdot)$; that is, a curve satisfying
$\Gamma'(s)=B(t,\Gamma(s))$. In our setting $B$ is merely H\"older continuous, so magnetic field lines need not be unique, but existence is guaranteed by Peano's theorem. Further, since $B$ is non-vanishing, any field line $\Gamma$ is a regular $C^1$ curve. 

The connection with MHD turbulence theory, including the GS95
exponents and the Iroshnikov--Kraichnan weak-turbulence exponent
$1/4$, is discussed in Subsection~\ref{ssec:phenomenology}. For the moment observe that our conditions approach (but do not include) the endpoint case $\gamma_v\to\tfrac14$, $\gamma_B\to\tfrac12$, where the H\"older exponent of the magnetic field $B$ along field lines 
$\gamma_B/(1-\gamma_B)\to 1$.

\subsection{Context}
In order to explain the physical context and motivation of our work,
let us first recall the closely related hydrodynamic case and weak
solutions of the incompressible Euler equations. In a broad-brush
picture, ignoring corrections due to intermittency, there are two
classical questions defining the contours of the theory. The first is
about the scale-invariant power law characterizing the cascade of
kinetic energy; in the hydrodynamic context this is the famous
$k^{-5/3}$ law of Kolmogorov, the K41 theory \cite{Kolmogorov1941a}. The second is about the
critical regularity required for the validity of the ideal
conservation of energy; in the hydrodynamic context this is precisely
the question raised by L.~Onsager, with the by now verified $C^{1/3}$
threshold on the H\"older scale; let us call this the O49 theory
\cite{Onsager}. Sufficiency of H\"older regularity above $1/3$ was proved
in \cite{CET}, whereas sharpness below this threshold was established in
\cite{Isett2018}. It
appears to be a stroke of coincidence that the K41 exponent $-5/3$ and
the O49 exponent $1/3$ are both uniquely determined and agree in terms
of scaling. In particular, the K41 theory does not require any
dynamical insight into what is going on in the inertial range, but the
need for such an insight would arise if dimensional analysis failed to
give a unique answer, which is what happens for MHD turbulence, as we
now discuss. We also remark that the O49 exponent conveniently agrees with
the natural bound reachable by convex integration as introduced by
Nash \cite{Nash1954}; see \cite{DLSz2016} for a discussion of this
fact.

Now let us turn to MHD turbulence. Since any solution of the
incompressible Euler equations is trivially a solution of the ideal
MHD system with zero magnetic field, a fortiori both the $k^{-5/3}$
spectrum and the $C^{1/3}$ threshold appear as candidates. However, the
presence of a large-scale nonzero magnetic field $B_0$ substantially
changes the situation. Firstly, it introduces an inherent anisotropy,
in stark contrast with one of the basic assumptions underlying the K41
theory, leading to a distinction between $k_\parallel$ and $k_\perp$,
the components of the wave vector $k$ of the fluctuations of the
turbulent fields parallel and perpendicular to $B_0$. Secondly, it
leads to the presence of Alfv\'en waves, whose wave packets propagate
along $B_0$ with frequency $\omega=\pm B_0\cdot k$; see
Subsection~\ref{ssec:alfven} below. This introduces a second characteristic
time scale, the Alfv\'en time $\tau_A\sim(v_A|k_\parallel|)^{-1}$,
where $v_A=|B_0|$ in the normalization of \eqref{eq:mhd-ideal}. Both the anisotropy and the
additional characteristic scale mean that additional hypotheses are
needed to deduce scaling laws for the energy spectrum. Kraichnan
\cite{Kraichnan1965} argued that in homogeneous isotropic MHD
turbulence one should have equipartition of energy in the inertial
range between the kinetic and the magnetic energy spectra, and that the
scaling law should depend on $v_A$, leading to a $k^{-3/2}$ spectrum
(known as the Iroshnikov--Kraichnan spectrum). In contrast, Goldreich
and Sridhar \cite{Goldreich1995} argued that the strong anisotropy
should be reflected by different scaling laws for the $k_\parallel$
and $k_\perp$ spectra and, assuming that the Alfv\'en time equals the
eddy turnover time (the critical balance hypothesis), obtained
$E(k_\perp)\sim\epsilon^{2/3}k_\perp^{-5/3}$ for the perpendicular
energy spectrum, $\epsilon$ being the energy flux through the inertial
range. The nature of the inertial range in plasma turbulence is a
subject of intensive debate; we refer to the excellent survey
\cite{Schekochihin}. Our results provide mathematical evidence for the
admissibility of several proposed scalings at the level of anisotropic
regularity; see Subsection~\ref{ssec:phenomenology}.

In terms of analytic properties, it is known that a nonvanishing
large-scale magnetic field and the resulting hyperbolic structure can
have a regularizing effect on solutions; we refer to
\cite{BSS1988,ZHuali2026} and the work of the first author
\cite{Giardi2026Local}. Further, ideal MHD possesses two additional
ideal conserved quantities besides the total energy: the cross-helicity
and the magnetic helicity. These helicities generate the regular
Casimirs of the MHD coadjoint action \cite{KhesinPeraltaYang}. We recall that the total energy, the
cross-helicity and the magnetic helicity are given by
\begin{equation}\label{eq:invariants}
 \mathcal E(t)=\int_{\mathbb T^3}\bigl(|v|^2+|B|^2\bigr)\,\dd x,
 \qquad
 \mathcal H^\times(t)=\int_{\mathbb T^3}v\cdot B\,\dd x,
 \qquad
 \mathcal H^{\rm mag}(t)=\int_{\mathbb T^3}A\cdot B\,\dd x,
\end{equation}
with the fields evaluated at time $t$, where $A$ is a vector potential
of the magnetic field, $\curl A=B$.
For a divergence-free field $B$ with zero spatial mean, which is the
case for the solutions of Theorem~\ref{full:main}, the value of
$\mathcal H^{\rm mag}$ does not depend on the choice of the potential,
and we fix $A=\curl^{-1}B=-\Delta^{-1}\curl B$, the unique
divergence-free potential with zero mean. These functionals are well
defined for the continuous fields considered in this paper. Magnetic
helicity was first studied by Woltjer \cite{Woltjer} and interpreted
topologically as the linkage and twist of magnetic field lines in the
highly influential work of Moffatt \cite{Moffatt}; see also
\cite{ArnoldKhesinBook}. In analogy with the K41 hypothesis of
anomalous dissipation in the infinite Reynolds number limit, in MHD
turbulence the rate of total energy dissipation in viscous, resistive
MHD seems not to tend to zero when the Reynolds number and the
magnetic Reynolds number tend to infinity; see
\cite{DallasAlexakis,Linkmann2015nonuniversality,Mininni2009finite}.
On the other hand, magnetic helicity is a rather robust conserved
quantity even in turbulent regimes, and J.~B.~Taylor conjectured that
magnetic helicity is approximately conserved for small resistivities
\cite{Taylor1974}; this conjecture was recently verified in
\cite{FaracoLindberg-Taylor}.

The analogue of Onsager's conjecture, that is, the critical threshold
for the conservation of the above quantities, was first addressed by
Caflisch, Klapper and Steele \cite{CKS}, building upon the work of
Constantin, E and Titi \cite{CET}, and endpoint cases were subsequently
studied by Kang and Lee \cite{KangLee2007}. It turns out that the total
energy and the cross-helicity are conserved provided
$v,B\in C([0,T];C^\gamma(\mathbb T^3))$ with $\gamma>1/3$. More
precisely, for $v\in C([0,T];C^{\gamma_1}(\mathbb T^3))$ and
$B\in C([0,T];C^{\gamma_2}(\mathbb T^3))$ the total energy is conserved
if $\gamma_1>1/3$ and $\gamma_1+2\gamma_2>1$, and the cross-helicity is
conserved if $\gamma_2>1/3$ and $2\gamma_1+\gamma_2>1$, each condition
corresponding to the cubic terms in the respective balance. Conservation
of magnetic helicity, on the other hand, merely requires
$(v,B)\in L^3(\mathbb T^3\times(0,T))$
\cite{KangLee2007,FaracoLindberg-Taylor}. Thus a natural conjecture,
put forward for instance in \cite{rev2}, is the following direct
analogue of Onsager's conjecture: weak solutions
$v,B\in C([0,T];C^\gamma(\mathbb T^3))$ conserve the total energy and
the cross-helicity if $\gamma>1/3$, but may fail to do so if
$\gamma<1/3$. Theorem~\ref{full:main} settles the flexible half of this
conjecture. Note that, by the results just quoted, every continuous weak
solution of \eqref{eq:mhd-ideal} conserves magnetic helicity, so that
the constancy of the magnetic helicity asserted in
Theorem~\ref{full:main} is forced once the magnetic field is
continuous; the nontrivial statement is that it is nonzero.

\subsection{Consequences for the phenomenological theories}
\label{ssec:phenomenology}

Theorem~\ref{full:main} and Theorem~\ref{full:fieldline} capture
two features specific to MHD turbulence: transverse roughness and
greater regularity along the local magnetic field. The magnetic field
is everywhere nonzero and has nontrivial helicity, and the longitudinal
estimates follow its actual field lines. Thus the result addresses
the magnetic geometry of the fluctuations as well as the familiar
Onsager exponent of the Euler setting.

In the Goldreich--Sridhar theory \cite{Goldreich1995}, critical balance
gives, up to outer-scale factors,
\[
 k_\parallel\sim k_\perp^{2/3},\qquad
 E(k_\perp)\sim k_\perp^{-5/3},\qquad
 E(k_\parallel)\sim k_\parallel^{-2}.
\]
The parallel spectrum is measured relative to the local magnetic
field; see also \cite{NazarenkoSchekochihin2011}. Under the usual
spectral interpretation of second-order increments, these laws
correspond to exponents $1/3$ across the magnetic field and $1/2$ along it. According to Theorem~\ref{full:fieldline}, with $\gamma_v,\gamma_B \to 1/3$, the constructed solutions approach this pair
for both velocity and magnetic field. At fixed step $q$ of the convex integration scheme, this can be read as the critical balance hypothesis of GS95, namely
\[
    \chi\sim\frac{\text{Alfv\'en time} }{\text{eddy turnover time}}
    \sim\frac{\lambda_q\delta_q^{1/2}}{\lambda_q\delta_{B,q}^{1/2}}
    =\lambda_q^{\gamma_B-\gamma_v}\to1
\]
In fact, by choosing the parameters appropriately, we obtain the full ranges described in Theorem~\ref{full:fieldline}. As $\gamma_v$ increases from $1/4$ to $1/3$,
the velocity parallel exponent along magnetic field lines remains
$1/2$, while the limiting magnetic exponent decreases continuously
from $1$ to $1/2$. 

The value $\gamma_v=1/4$ is the exponent
associated with the $k^{-3/2}$ spectrum of
Iroshnikov--Kraichnan (IK) phenomenology \cite{Kraichnan1965},
based on weak Alfv\'enic interactions and an isotropy assumption. Theorem~\ref{full:fieldline} seems to suggest an anisotropic reinterpretation of this endpoint value. We finally remark that the IK picture should be distinguished from balanced
anisotropic weak-wave turbulence, whose perpendicular spectrum
is $k_\perp^{-2}$ \cite{GaltierEtAl2000}.

This continuous family provides mathematical evidence that the GS95
and $1/4$ velocity scalings, together with intermediate regimes, are all
possible at the level of the underlying PDE in terms of anisotropic regularity.

\subsection{Overview of the mathematical literature}\label{ssec:literature}
We start with a simple observation: any weak solution of the
three-dimensional incompressible Euler equations is trivially a weak
solution of \eqref{eq:mhd-ideal} with vanishing magnetic field
$B\equiv0$, so that Onsager's conjecture with exponent $1/3$ is valid
in this special case. More generally, weak $2\frac12$-dimensional
solutions have been constructed by Bronzi, Lopes Filho and Nussenzveig
Lopes \cite{Bronzi2015} in the class $L^\infty$ and, recently, by
Miao, Nie and Ye \cite{MiaoNieYe} in the class
$C^\gamma(\mathbb T^3\times\R)$ for every $\gamma<1/3$. These are
solutions of the form
\[
 \begin{aligned}
 v(x_1,x_2,x_3,t)
   &=(u^1(x_1,x_2,t),u^2(x_1,x_2,t),\theta(x_1,x_2,t)),\\
 B(x,t)&=(0,0,\theta(x_1,x_2,t)),
 \end{aligned}
\]
which in particular depend only on two coordinates orthogonal to the
magnetic field. Under this structural restriction, the
three-dimensional MHD system \eqref{eq:mhd-ideal} reduces to the
two-dimensional Euler equations for $u=(u^1,u^2)$ coupled with a
passive scalar equation for $\theta$:
\[
 \begin{cases}
 \partial_tu+u\cn u+\nabla p=0,\\
 \partial_t\theta+u\cn\theta=0,\\
 \ddiv u=0.
 \end{cases}
\]
For this special class of solutions, the authors of
\cite{MiaoNieYe} are able to reach the optimal $1/3$ exponent. However, in this setting the magnetic helicity is zero,
and its density vanishes in the following gauge. Under the assumption
$\int_{\mathbb T^3}\theta=0$, solve
\[
 \Delta\psi=\theta, \qquad \int_{\mathbb T^3}\psi=0.
\]
The solution $\psi$ is independent of $x_3$, since $\theta$ is.
Consequently,
\[
 A=(-\partial_2\psi,\partial_1\psi,0),\qquad
 \curl A=B,\qquad A\cdot B=0.
\]
This reduction excludes a nontrivial inverse cascade of magnetic
helicity \cite{Frisch1975,Muller2012}.

The first results for weak solutions which are not of this reduced
form appeared in the work of Beekie, Buckmaster and Vicol \cite{BBV}
and in the work of Faraco, Lindberg and the second author
\cite{FLSz2021}. In \cite{BBV} the authors constructed unbounded weak
solutions in the energy class $v,B\in L^\infty_tL^2_x$ which preserve
neither the magnetic helicity nor the energy or the cross-helicity.
Compared with the thresholds for the conservation of energy and
magnetic helicity obtained in \cite{CKS,KangLee2007} mentioned above,
we see that the space $L^\infty_tL^2_x$ is supercritical with respect
to magnetic helicity conservation. In view of Taylor's conjecture,
proved in \cite{FaracoLindberg-Taylor}, such solutions cannot arise as
weak ideal limits of Leray--Hopf solutions of viscous, resistive MHD. In \cite{FLSz2021} the
authors constructed bounded (but discontinuous) weak solutions which
do not preserve energy and cross-helicity; of course, a fortiori,
magnetic helicity is conserved, and indeed vanishes identically. In
the subsequent work \cite{FLSz2024} the authors extended their
construction to bounded weak solutions with nonvanishing (constant)
magnetic helicity. In both works \cite{FLSz2021,FLSz2024} a key point
was to introduce a relaxation of the MHD system \eqref{eq:mhd-ideal}
which decouples the effects of hydrodynamic turbulence (in the form of
the appearance of a Reynolds stress term in the momentum equation)
from possible small-scale effects in the Faraday--Ohm system
\eqref{eq:faraday-ohm}, in a manner consistent with the robustness and
conservation of magnetic helicity. In a nutshell, the relaxation
involves replacing Ohm's law (the third equation in
\eqref{eq:faraday-ohm}) by the nonlinear constraint $E\cdot B=0$.

Very recently, Enciso, Pe\~nafiel-Tom\'as and Peralta-Salas
\cite{EnPePe} succeeded in constructing H\"older continuous weak
solutions that do not conserve energy and cross-helicity and have
nontrivial magnetic helicity. The weak solutions obtained are in the
class $v,B\in C^\gamma(\mathbb T^3\times[0,T])$ with $\gamma\le1/200$.
Their construction is based on convex integration applied to a more
restrictive, partial relaxation of \eqref{eq:mhd-ideal}, in which the
Faraday--Ohm system \eqref{eq:faraday-ohm} is solved exactly, leaving
the appearance of a Reynolds stress term in the momentum equation as
the only small-scale effect; see \eqref{eq:mhd} below. In this way the
magnetic field $B$ is slaved to the velocity via Lie transport.
The present paper reaches every exponent below the Onsager threshold
$1/3$ within this partial relaxation. In the spirit of Arnold's
formulation of ideal hydrodynamics \cite{ArnoldKhesinBook}, we regard
a solution as a path of volume-preserving diffeomorphisms and construct
its perturbations in the Lie algebra of divergence-free vector fields.
Consequently, the induction equation is
satisfied exactly at every stage of the iteration, no symmetry
reduction is used, and the magnetic field is of order one, nowhere
zero and has nonzero magnetic helicity. The mechanism which makes
this possible, namely the separate estimate of material and magnetic
derivatives in charts adapted to the magnetic field, is explained in
Section~\ref{sec:proof-overview}. 

Finally, we point out an interesting connection to the local well-posedness theory in $H^s$ of ideal MHD. The first author's recent work \cite{Giardi2026Local}, extending the remarkable work of Huali Zhang \cite{ZHuali2026}, proves local well-posedness for the case of a nonzero constant initial magnetic field on $\R^n$ in the space $H^s$ with
\[
s> \frac{n+1}{2}.
\]
The constant magnetic background is subtracted in the Sobolev norm.
These solutions conserve renormalized energy and cross-helicity.
Since $s-n/2>1/2$, Sobolev embedding gives spatial H\"older exponents
strictly above $1/2$, while our longitudinal velocity exponents
approach $1/2$ from below - see Theorem \ref{full:fieldline}. At the level of directional H\"older
scaling, our flexibility result thus approaches this rigidity result
from below in the magnetic direction. 

\subsection{AI declaration} The authors made use of Anthropic's Claude Fable 5.1 and ChatGPT Astra in preparing this manuscript. In the following we specify their role in this work. 

First of all, this paper should be viewed as a technical extension of our previous work \cite{gmzl2026c15convexintegrationsolutions}, which itself was prepared entirely without the use of AI and arose from the first author's PhD thesis, ``On Onsager's Conjecture for Ideal MHD'' \cite{Giardi2026Thesis}. The main result in \cite{gmzl2026c15convexintegrationsolutions} is analogous to Theorem \ref{full:main} but with the restriction $\gamma<1/5$. In this paper we improve this to the optimal exponent $\gamma<1/3$. However, we wish to emphasize that this improvement is essentially a technical rather than a conceptual advancement: the main steps in the actual construction of the weak solutions remain the same. Indeed, as already pointed out in \cite{gmzl2026c15convexintegrationsolutions} (see its Section 1.3 on the Regularity Cap), the $1/5$ restriction arises because we did not perform the full iteration as has been implemented in \cite{GR,GiriKwonNovack2026}, as doing this would have required significant additional work, in particular on:
\begin{enumerate}
    \item Higher-order transport derivatives with corresponding antiderivative gain mechanisms entering the coefficient construction and the Galbrun solution. See \S\ref{ssec:principalperturbationoverview} and \S\ref{ssec:grintro} below for the role of these in the construction.
    \item Tracking the additional smallness of the parallel gradients. See \S\ref{prep:old-mixed}.
\end{enumerate}
In addition to implementing these using the AI assistance following precisely targeted prompts, the current manuscript also includes some substantial simplifications in the bookkeeping and organization of how to obtain the estimates, including for instance the formalization of the class of building blocks into derivative classes - these give an efficient repackaging of the sections "Lie-Taylor Series Expansion" and "Inductive Lemma" of \cite{gmzl2026c15convexintegrationsolutions}. The basic idea of this simplification was first developed by the authors, and pressured by the current times we then instructed AI to implement it in a rigorous proof. The same efficient reorganization was achieved with AI-assistance in transferring the first-order transport estimates in e.g. \cite[(1.16)]{gmzl2026c15convexintegrationsolutions} to their higher-order counterparts required here; see \eqref{prep:old-mixed}-\eqref{prep:old-stress-time} below. Finally, tables of notation and pictures are also AI-generated, e.g. Figure \ref{fig:window}.

We provide no Lean verification as we believe none is needed. We are fully aware of the content of the manuscript, and the mistakes are our own. Any comment or correction is welcome and the reader is invited to contact us directly.

\subsection{Organization of the paper}
Section~\ref{sec:proof-overview} explains the mechanism of the
construction and the parameter constraints in an informal way, starting
with the constant-background Alfv\'en waves and their exact local
realization through the magnetic charts.
Section~\ref{sec:inductive-scheme} introduces the scales, states the
inductive hypotheses and formulates the iterative step,
Proposition~\ref{full:step}.
Sections~\ref{sec:preparation}--\ref{sec:adapted-perturbation} construct
the local backgrounds, the charts and amplitudes, and the complete
deformation, ending with the exact endpoint stress identity of
Theorem~\ref{mom:complete-ledger}.
Sections~\ref{sec:generator-field-estimates}--\ref{sec:endpoint-stress}
estimate the Lagrangian displacement fields, velocity and magnetic fields, and stress, retaining the scale
inequality required by each error. Section~\ref{sec:parameter-choice}
then solves these inequalities and fixes the parameters, and
Section~\ref{sec:iteration-closure} closes the induction;
Section~\ref{sec:initialization-and-limit} supplies the initial approximate solution,
passes to the limit and proves Theorem~\ref{full:main} and
Theorem~\ref{full:fieldline}. The appendices contain the path
calculus and the properties of its derivative classes
(Appendix~\ref{sec:path}), the MHD Frobenius theorem
(Appendix~\ref{sec:mhd-frobenius}), and the higher-order estimates,
including the symmetric and oscillatory inverse divergence
(Appendix~\ref{sec:analytic-estimates},
Subsection~\ref{sec:inverse-divergence}). The forced Galbrun equation
is treated in Appendix~\ref{sec:galbrun}.

\section{Overview of the Proof}\label{sec:proof-overview}

The proof uses a convex integration scheme: at each step, we add perturbations satisfying a multiscale ansatz and use their averaged quadratic interactions to correct the stress from the previous iterate. The constant-background Alfv\'en calculation in \S\ref{ssec:alfven}
identifies the oscillations needed to balance this stress. We then explain how to obtain these by means of
divergence-free Lagrangian displacement fields (LDFs), deforming the flow map while keeping induction exact. After giving a flow chart of the proof in \S\ref{ssec:flowchart}, we discuss each step and in particular explain how the MHD Frobenius charts realise the same constant transport geometry
locally (\S\ref{ssec:frameintro}), how a calculus adapted to these deformations yields the required estimates (\S\ref{ssec:pathcalculusintro}), and how fast temporal oscillations enter the construction (\S\ref{ssec:principalperturbationoverview}). The constraints on the exponents are explained separately (\S\ref{ssec:overview-finite-derivatives}).

\subsection{Alfv\'en waves and partial relaxation}\label{ssec:alfven}
To motivate the construction, consider small fluctuations
$( v+w, B+b)$ around a constant state $( v, B)$ and
linearize \eqref{eq:mhd-ideal}:
\[
 \partial_tw+ v\cn w- B\cn b+\nabla p=0,\qquad
 \partial_tb+ v\cn b- B\cn w=0.
\]
Equivalently, in the Els\"asser variables $z^\pm=w\pm b$,
\[
 \partial_tz^\pm+( v\mp B)\cn z^\pm+\nabla p=0,
\]
which admit the plane-wave solutions
\begin{equation}\label{eq:ansatz0}
    z^\pm=\delta^{1/2}e^{i\lambda(x-t v\pm t B)\cdot k}\zeta^\pm
\end{equation}
with amplitude $\delta^{1/2}$, wave vector $k$ and state vectors
$\zeta^\pm$ with $\zeta^\pm\cdot k=0$. Whilst one can assume $ v=0$
by Galilean invariance, the differing wave frequencies $\mp\lambda B\cdot k$
reflect fundamental properties of Alfv\'en wave packets, most
pronounced in the case $k\parallel B$. On the other hand, the
nonlinear term in Els\"asser variables takes the form
$\ddiv(z^+\otimes z^-)$, and a key aspect of convex integration
constructions is that spatial averages of the nonlinear term have to
balance the Reynolds stress: we would like to have
$\langle z^+\otimes z^-\rangle=R$, uniformly in time; see \eqref{eq:relaxationintro}. For these single-phase waves, a nonzero time-independent spatial
average requires $k\cdot B=0$. This leads to the ansatz
\begin{equation}\label{eq:ansatz}
 w(x,t)=\delta^{1/2}\varphi\bigl(\lambda(x-t v)\cdot k\bigr)\zeta,
 \qquad k\cdot B=0,\quad k\cdot\zeta=0.
\end{equation}
The first step in any convex integration construction is to identify a
suitable relaxed system of equations. This is usually achieved by
local averaging or filtering, and is also a key step in identifying
the threshold regularity for energy conservation
\cite{CET,Eyink1994,EyinkReview,AluieEyink}. Such a filtering process
detects additional defect terms, arising from the effect of
small-scale fluctuations on the large-scale dynamics via nonlinear
interactions. In the case of \eqref{eq:mhd-ideal} the relaxation takes
the form
\begin{equation}\label{eq:relaxationintro}
 \begin{cases}
     \partial_t v+\ddiv( v\otimes v- B\otimes B)
   +\nabla p=\ddiv R,\\
 \partial_t B+\curl( B\times v)=\curl M,\\
 \ddiv v=\ddiv B=0,
 \end{cases}
\end{equation}
where $R$ is a symmetric two-tensor (usually called the Reynolds
stress) and $M$ is a vector field (called the subscale electromotive
force in \cite{AluieEyink}). An additional feature, not seen by pure
filtering but resulting from the differential (div-curl) structure,
first observed by L.~Tartar \cite{Tartar}, is that, at least in the
presence of a large scale separation, $M$ necessarily satisfies a
geometric constraint of the type $M\cdot B=0$. This constraint is
a local expression of magnetic helicity conservation \cite{FLSz2024}.
Following \cite{EnPePe}, we impose $M\equiv0$. At step $q$,
the relaxed fields $(v_q,B_q,p_q,R_q)$ satisfy
\begin{equation}\label{eq:mhd}
 \begin{cases}
 \partial_t v_q+\ddiv(v_q\otimes v_q-B_q\otimes B_q)
        +\nabla p_q=\ddiv R_q,\\
 \partial_t B_q+[v_q,B_q]=0,\\
 \ddiv v_q=\ddiv B_q=0.
 \end{cases}
\end{equation}
Only the momentum equation is relaxed. The induction equation says
that the magnetic field is the Lie transport of its initial value;
in the preceding notation,
$ B(\cdot,t)=( X_t)_* B|_{t=0}$, where $ X_t$ is the
Lagrangian flow of $ v$. To keep this constraint exact, as in \cite{EnPePe} we are therefore forced to work directly at the level of the flow map, transforming the
velocity and magnetic field together. In contrast with \cite{EnPePe}, however, we perturb the particle
trajectories by volume-preserving continuous deformations driven by an element of the Lie algebra, that is, a divergence-free vector field.

\textbf{Field deformations.} More precisely, for a volume-preserving map $X_s(t,\cdot)$ with
$X_0=\IId$, set
\begin{equation}\label{eq:sfieldsintro}
    v_s=\partial_tX_s\circ X_s^{-1}+(X_s)_*v_q,
 \qquad B_s=(X_s)_*B_q.
\end{equation}
Here $s\in[0,1]$ is the deformation parameter, while physical time
remains $t$. The velocity formula is the transformation law of the
spacetime vector $\partial_t+v_q$, so equivariance of the Lie bracket under pushforward gives
\[
 [\partial_t+v_s,B_s]=(X_s)_*[\partial_t+v_q,B_q]=0.
\]
Both fields remain
divergence free, and the perturbation creates no induction defect.
Thus the same deformation preserves precisely the transport constraint
needed by the construction.

We specify the deformation through its divergence-free generator
$\xi_s=\partial_sX_s\circ X_s^{-1}$. This is a generalised version of the Lagrangian displacement field (LDF)
of Newcomb \cite{Newcomb} and Vladimirov, Moffatt and Ilin \cite{VMI}:
it moves the construction from the group of volume-preserving
maps to its Lie algebra. The first variations are
\begin{equation}\label{eq:firstvariationoverview}
    \partial_sv_s=\mathcal D_{t,s}\xi_s,
 \qquad \partial_sB_s=\mathcal L_{B_s}\xi_s,
\end{equation}
where $\mathcal D_{t,s}=\partial_t+\mathcal L_{v_s}$ is the Lie material
derivative. At $s=0$, a displacement $\xi$ therefore has velocity
and magnetic first variations
$\mathcal D_{t,q}\xi=(\partial_t+\mathcal L_{v_q})\xi$ and
$\mathcal L_{B_q}\xi$.
These are the first-order approximations of our velocity and magnetic increments, as a Lie-Taylor expansion at $s=0$ shows; see \eqref{eq:taylor}.

\textbf{Momentum deformation.} The corresponding momentum change has a linear and a quadratic part.
Set $w=v_1-v_q$ and $b=B_1-B_q$. For the momentum expression
\[
    \mathcal M(v,B)=\partial_tv+\ddiv(v\otimes v-B\otimes B),
\]
the exact $s=1$ endpoint identity at step $q$ is
\begin{equation}\label{eq:momentumvariationintro}
    \mathcal M(v_q+w,B_q+b)-\mathcal M(v_q,B_q)
 =\mathcal F_{v_q,B_q}(w,b)+\ddiv(w\otimes w-b\otimes b),
\end{equation}
where $\mathcal F_{v_q,B_q}$ is the linearized force. Substituting the two
first variations into this force gives a second-order equation for the
displacement, the MHD version of the Galbrun equation of acoustics, sometimes called the Frieman-Rotenberg equation
(introduced in \cite{frieman1960hydromagnetic} by Frieman and Rotenberg in their study of stability of MHD equilibria; see also Lindblad \cite{Lindblad2005FreeBoundary} for its relevance in the free boundary
problem for the Euler equation). This will play an important role and is addressed in \S\ref{ssec:grintro} below. 

\subsection{Basic heuristics on the exponents}\label{ssec:overview-heuristics}

Let us explain the heuristics behind the basic scheme, our inductive assumptions in Section \ref{ssec:fixed-scales} below, and how they lead to our claimed $C^{1/3-}$ regularity. Recall that our construction proceeds by building a sequence of approximations $(v_q,B_q,p_q,R_q)$ which solve the partially relaxed system \eqref{eq:mhd}, where the defect $R_q$ is a symmetric tensor. The starting point, as in numerous works on convex integration in the context of Onsager's conjecture \cite{DS1,DLSz2016,Isett2018,BuckmasterDeLellisSzekelyhidiVicol2019,GR,burczak_szekelyhidi_wu_2023} is a super-exponentially growing sequence of frequencies
$\lambda_q\sim a^{b^q}$ for some large $a\gg1 $ and $b>1$ with $b\sim 1$, and a corresponding sequence of amplitudes $\delta_q^{1/2}=\lambda_q^{-\beta}$ for some $\beta>0$ - see Section \ref{ssec:fixed-scales} for the precise definitions. The sequence of velocity and magnetic increments satisfies
\begin{align*}
\|v_{q+1}-v_q\|_{0}&\lesssim \delta_q^{1/2}\,,\\
\|\nabla v_q\|_{0}&\lesssim \delta_q^{1/2}\lambda_q\,,\\
\|B_{q+1}-B_q\|_{0}&\lesssim \delta_{B,q}^{1/2}\,,\\
\|\nabla B_q\|_{0}&\lesssim \delta_{B,q}^{1/2}\lambda_q\,,	
\end{align*}
with $\delta_{B,q}$ denoting the amplitudes of magnetic field increments. As is easy to see from these assumptions by interpolation, $\beta>0$ denotes the borderline H\"older exponent reachable in the limit of such an approximating sequence. The size of the defect $R_q$ is measured as 
\begin{align*}
\|R_q\|_0\lesssim \delta_{q+1},	
\end{align*}
the choice being motivated by the fact that $R_q$ formally arises from high-high to low interaction terms from small-scale velocity (Reynolds-) and magnetic (Lorentz-)stress, which with some abuse of notation we may denote as
\begin{align*}
R_q\sim \langle (v_{q+1}-v_q)\otimes (v_{q+1}-v_q)-(B_{q+1}-B_q)\otimes (B_{q+1}-B_q)\rangle+{\rm l.o.t.}  	
\end{align*}
In fact it turns out that magnetic stress contributions are negligible, $\delta_{B,q}\ll \delta_q$, cf. Section \ref{ssec:fixed-scales}. Thus, in a sense we treat the MHD system \eqref{eq:mhd-ideal} and its partial relaxation \eqref{eq:mhd} as a small, nonlinear perturbation of the Euler equations and their relaxation, usually called the Euler-Reynolds system \cite{DS1}. Thus, in a sense this generalizes the case of $2\frac12$ dimensional flow constructions in \cite{Bronzi2015,MiaoNieYe}. In this connection we recall that the $1/3$ borderline in \cite{BuckmasterDeLellisSzekelyhidiVicol2019,MiaoNieYe} arises from the linear interaction terms, namely from the interaction of large-scale velocity gradients with small-scale fluctuations (called ``Nash-term'') and from the Lagrangian deformation of small-scale fluctuations by the large-scale flow (called ``transport term''). In particular the latter requires control almost up to the natural (``eddy turnover'') time-scale $\tau_c\sim (\delta_q^{1/2}\lambda_q)^{-1}$ (in fact there is some room, see Section \ref{ssec:fixed-scales}). In the present paper the corresponding two terms are denoted by $R^{\mathrm p}_{\rm na}$ and $R^{\mathrm p}_{\rm tr}$, referred to as principal Nash stress and principal transport stress - see \eqref{mom:principal-interaction-stress}-\eqref{mom:principal-transport-stress}.

If we assume for the moment that the worst residual error arises from a temporal cutoff on time-scale $\tau_c$, we are led to the estimate 
\begin{equation*}
\|R_{q+1}\|_{0}\lesssim \frac{\delta_{q+1}^{1/2}}{\lambda_{q+1}\tau_c},	
\end{equation*}
which, in combination with the inductive assumption above, leads to the required exponent $1/3$ - see for instance Chapter 7 in \cite{DLSz2016} for a quick heuristic explanation. 
Although the original resolution of Onsager's conjecture in \cite{Isett2018} (see also \cite{BuckmasterDeLellisSzekelyhidiVicol2019}) used the size of the residual linear interaction terms to determine the optimal exponent $\beta$, in later works (see for instance \cite{GR,GiriKwonNovack2026,burczak_szekelyhidi_wu_2023}) it was realized that it suffices to make sure these terms are qualitatively small, because then to reach quantitative smallness one can use higher order expansions and/or an inner iteration. In our work this small gain is related to the ratio $\tau_a/\tau_c$, where $\tau_a$ is the time-scale of temporal oscillations - see Section \ref{ssec:fixed-scales}, in analogy to the work of Giri-Radu \cite{GR}. In our previous work \cite{gmzl2026c15convexintegrationsolutions} this smallness was the main restriction, leading to the exponent $1/5$. To remove this technical restriction requires higher order expansions and much more careful bookkeeping, although the conceptual framework of the construction remains the same. In particular, our actual inductive assumptions in \eqref{prep:old-ordinary}-\eqref{full:stress-support} include separate bounds on pure spatial, transport and magnetic derivatives, each to a certain high order ($j_1$ for transport and magnetic, $r_{\rm good}$ for pure spatial), chosen sufficiently large to be able to account for various derivative losses. 


\subsection{Ansatz and flowchart of the proof}\label{ssec:flowchart} 
The Lagrangian displacement field (LDF) used in the construction has the form 
\[
    \xi_s=\xi^{\mathrm p}_s+\xi^{\mathrm c}.
\]
The principal part $\xi^{\mathrm p}_s$ contains the fast space and fast time oscillations. Its
quadratic interactions correct $R_q$ after space--time averaging;
see Subsection~\ref{ssec:principalperturbationoverview}. The corrector
$\xi^{\mathrm c}$ varies slowly in space and cancels the remaining
spatially averaged forcing, whose fast temporal profile has zero mean. In fact, it rewrites the stress in advance,
as in the construction of Giri--Radu \cite{GR}; see
Subsection~\ref{ssec:grintro}. 

In the remaining part of \S\ref{sec:proof-overview} we explain its construction, which consists of the following main steps:
\begin{enumerate}
    \item \textit{Smoothing and background preparation:} evolve the mollified data with the smoothed stress on time intervals of length comparable to $\tau_c$, recovering the induction equation and hence the Alfv\'en commutation structure. Estimate the differences from the original fields. Output: $v_{\ell,n},B_{\ell,n}$.
    \item \textit{Stress decomposition:} construct the Lagrangian charts $\Psi_I$ adapted to the magnetic field, decompose the smoothed stress into rank-one tensors with positive coefficients $a_I$, and form the localised approximate antiderivatives $\mathfrak a_I$ using the fast time profiles $\alpha_I$. Output: $a_I, \alpha_I,\mathfrak a_I,\Psi_I$.
    \item \textit{Construction of the corrector:} solve the linearised equation in the Lie algebra with forcing built from the zero-mean fast temporal profile $1-\alpha_I^2$, then localise the solution using source-free time cut off regions. This prepares the stress for the principal perturbation, à la Giri--Radu \cite{GR}. Output: $\xi_n^{\mathrm c}=\curl \Theta_n^{\mathrm c}$.
    \item \textit{Construction of the principal part:} form the oscillatory potentials from the charts and coefficients, take their curls, and push forward the complete potentials by the corrector flow. Their quadratic interactions balance the remaining stress while the induction equation remains exact. Output: $\xi_{I,s}^{\mathrm p}=\curl \Theta_{I,s}^{\mathrm p}$.
    \item \textit{Space-time averaging:} cancel the old stress with the averaged quadratic terms. Estimate the oscillatory terms, cutoff errors and background differences to obtain the new stress and close the field estimates.
    \item \textit{Initialisation of the scheme and conclusion:} choose the initial fields with different energy and cross-helicity at two times outside the perturbation supports. Iterate the construction and pass to the limit using the summability of the increments.
\end{enumerate}
We now explain these steps.

\subsection{Preparation: non-linear smoothing and local backgrounds}\label{ssec:deeppreparation} Mollification is a standard ingredient in convex integration schemes. In our setting, it does not preserve the induction equation because it does not commute with the quadratic Lie bracket. To circumvent this issue we partition time at the scale $\tau_c$, obtained by slightly shortening the eddy turnover time~$(\lambda_q\delta_q^{1/2})^{-1}$, and then evolve the mollified data by the relaxed system with
the smoothed stress. This step schematically reads
\[
    (v_q,B_q)\leadsto (v_{\ell}, B_\ell)\leadsto (v_{\ell,n},B_{\ell,n}),
\]
where the spatial mollification length $\ell$ satisfies
\[
    \lambda_q< \ell^{-1}\ll\lambda_{q+1},
\]
and $n$ indexes the intervals in the time partition. The resulting smooth local background on each interval satisfies
the induction equation exactly, and approximation estimates
compare it with the original fields; these differences, the
background gaps, 
\begin{equation}\label{eq:gapsintro}
    (v_q, B_q)-(v_{\ell,n},B_{\ell,n})
\end{equation}
remain explicit inputs to the final stress estimate. Choosing the mollifier with sufficiently many vanishing moments makes these differences sufficiently small even when the exponent $\gamma_\ell$ in $\ell\lambda_q=\varepsilon_{q+1}^{\gamma_\ell}<1$ is arbitrarily small; see Proposition~\ref{moll:spatial-basic} and \cite{gromov1986partial}.

\subsection{Preparation: local charts and frames} \label{ssec:frameintro}
We can now return to \eqref{eq:ansatz} on a nonconstant background. Suppressing indices $n$ and $\ell$, the differential-geometric structure that makes this passage possible, at
least locally in space-time, is the commutation relation
\begin{equation}\label{eq:commutation}
 [\partial_t+ v\cn, B\cn]
 =(\partial_t B+ v\cn B- B\cn v)\cn=0,
\end{equation}
a consequence of the induction equation, which the local background construction
satisfies exactly. Equivalently, the two Alfv\'en
transport operators 
\[
\mathcal{A}^\pm=\partial_t+( v\pm B)\cn
\]
commute. This lets us straighten the magnetic field locally in space on
one time slice and
then transport the coordinates by the velocity flow, preserving the
straightening identity throughout the slow-time window. This is the MHD Frobenius theorem of
Appendix~\ref{sec:mhd-frobenius}. We obtain charts $\Psi_I$ satisfying
\[
 D_t\Psi_I=0,\qquad D_B\Psi_I\parallel e_3,
 \qquad \det\DD\Psi_I=1.
\]
The index $I$ records the space--time
localization and a positively oriented orthonormal triplet
$(k,\nu,\zeta)$, with $k\times\nu=\zeta$. The chart itself does not depend on this triplet. For scalars written in the coordinates $y=\Psi_I(t,x)$, these identities
give the exact correspondence
\[
 D_t=\partial_t\big|_y,\qquad
 D_B=c_I^{-1}\partial_{y^3}.
\]
Thus the two transport directions are those of a constant magnetic
background in Lagrangian coordinates. For an orthonormal triplet as above, in the directional coordinates
\[
    y_I^\eta=\Psi_I\cdot\eta \qquad \eta\in\{k_I,\nu_I,\zeta_I\}
\]
choose the oscillation direction $k_I$ with $k_I\cdot e_3=0$. For any profile $\varphi:\R\to \R$ we then have
\[
    D_t[\varphi(\lambda y_I^{k_I})]=D_B[\varphi(\lambda y_I^{k_I})]=0.
\]
Volume preservation finally gives divergence-free vector fields
\[
 \frac{\partial}{\partial y_I^{\zeta_I}}=dy_I^{k_I}\times dy_I^{\nu_I},
\]
and similar transport identities hold when we replace the transport derivatives with Lie derivatives for this frame and coframe. With this at hand, we can rewrite \eqref{eq:ansatz} as follows:
$(x-t v)\cdot k$ becomes $y_I^{k_I}$ and the constant state vector
$\zeta$ becomes $\frac{\partial}{\partial y_I^{\zeta_I}}$, namely
\[
    w=\delta^{1/2}\varphi\bigl(\lambda(x-t v)\cdot k\bigr)\zeta \qquad \leadsto \qquad w=\delta^{1/2}\varphi\bigl(\lambda y_I^{k_I}\bigr)\frac{\partial}{\partial y_I^{\zeta_I}}
\]
This realises the
constant-background geometry exactly in each chart of the local background. The invariant properties and geometric identities that these objects enjoy are crucial in the estimates below.

The full magnetic gradient of the iterate is smaller than the velocity
gradient, since $[B_q]_1\lesssim\lambda_q\delta_{B,q}^{1/2}$ while
$[v_q]_1\lesssim\lambda_q\delta_q^{1/2}$ and $\delta_{B,q}\ll \delta_q $.
This permits localization on spatial scale $\lambda_\parallel^{-1}$ and time scale $\tau_c$ while satisfying the deformation bound
\[
 \lambda_\parallel^{-1}[B]_1+\tau_c[v]_1
 =o(1),
 \qquad \tau_c\ll\lambda_\parallel^{-1},
\]
and preserving the anisotropic scaling required by the construction. This is the job of the cutoffs $\eta_I$ in time and $\chi_I$ in space introduced below. The full gradient bound on $B_q$ is essential for this construction.

\subsection{From the LDF to the fields: the role of fast time}\label{ssec:principalperturbationoverview} 
The tensor $\delta_{q+1}\IId-R_q$ is close to a positive multiple of
the identity. After smoothing in space and along the commuting
velocity and magnetic flows in spacetime, a fixed geometric decomposition expresses
it as a positive sum of squares of coordinate vectors, namely
\begin{equation}\label{eq:decompositionintro}
    \varrho^2\sum_n\eta_n^2\mathcal J_n(\delta_{q+1}\IId-R_q)=\sum_{I}a_I^2
   \partial_{y_I^{\zeta_I}}\otimes\partial_{y_I^{\zeta_I}}.
\end{equation}
The coefficients
of this decomposition determine the amplitudes directly:
\[
 a_I=\varrho\eta_I\chi_I
 \left[q_{\zeta_I}\!\left((\DD\Psi_I)
   \mathcal J_n(\delta_{q+1}\IId-R_q)(\DD\Psi_I)^T\right)\right]^{1/2}.
\]
Here $q_\zeta$ is a fixed linear functional on symmetric tensors, positive near the
identity, $\mathcal J_n$ is the tensor smoothing, $\varrho$ the outer
time cutoff, $\eta_I$ the slow time partition and $\chi_I$ the spatial cutoff. The difference between the smoothed and the
original target is part of the new stress. The fast time mentioned in \S\ref{ssec:overview-heuristics} enters the construction in the following two mechanisms, which can be traced back to \cite{GR} and \cite{EnPePe}. 

\textit{1) Disjointness of the local contributions.} The chart construction above reduces the oscillations to plane waves
transverse to the magnetic direction. Waves associated with different directions
or charts may have overlapping supports. We group the indices into finitely
many types, each recording the time parity, spatial color and direction,
as in Lemma~\ref{amplitude:profiles}. Local contributions of the same
type have disjoint supports; different types are separated by periodic
time profiles at scale $\tau_a$. Fix smooth unit-periodic profiles $\alpha_I$ with
$\int_0^1\alpha_I^2\,\dd\tau=1$. With $\varphi(u)=\sqrt2\sin u$, a natural ansatz for the velocity perturbation $w$ of the form \eqref{eq:ansatz} is then
\[
 w\sim a_I\alpha_I(t/\tau_a)\varphi'(\lambda_{q+1}y_I^{k_I})
       \frac{\partial}{\partial y_I^{\zeta_I}}.
\]
The resulting local contributions have disjoint space-time supports: distinct types occupy separated time subintervals,
while local contributions of the same type already have disjoint supports. 

\textit{2) The fast-time gain in the LDF.}  By \eqref{eq:firstvariationoverview}, to realize a velocity variation
$w=v_1-v$ of size $a_I$ through an LDF, we need to solve approximately
\begin{equation}\label{eq:firststateintro}
    w=(\partial_t+\mathcal{L}_{v})\xi_0+\dots
\end{equation}
At the scalar level this requires an approximate antiderivative of
$\alpha_I(t/\tau_a)a_I$ along the velocity flow. A normalized $(j_0+1)$-st derivative of a
bump function has $j_0$ successive primitives supported in the same time subinterval
(Lemma~\ref{amplitude:profiles}). We choose $\alpha_I$ in this way and set
\[
 \mathfrak a_I=
 \sum_{j=0}^{j_0-1}(-1)^j\tau_a^{j+1}
       \alpha_I^{[j+1]}(t/\tau_a)D_t^ja_I.
\]
The sum telescopes under one material derivative,
\begin{equation}\label{eq:telescopeintro}
    D_t\mathfrak a_I=\alpha_I(t/\tau_a)a_I+a_I^{\rm c},
 \qquad
 a_I^{\rm c}=(-1)^{j_0-1}\tau_a^{j_0}
       \alpha_I^{[j_0]}(t/\tau_a)D_t^{j_0}a_I,
\end{equation}
so the primitive has size $\tau_a\delta_{q+1}^{1/2}$, is supported
within the support of $a_I$, and has an error of size
$\varepsilon_\tau^{j_0}\delta_{q+1}^{1/2}$ at low derivative orders, where $\varepsilon_\tau=\tau_a/\tau_c$ is the ratio of the fast and slow time scales. Lemma~\ref{principal:material-primitive} gives the construction and estimates.

The principal LDF is then defined through its one-form potential 
\begin{equation}\label{eq:principalpreparation}
    \Theta_I^{\rm p}=\lambda_{q+1}^{-1}\mathfrak a_I
     \varphi(\lambda_{q+1}y_I^{k_I})dy_I^{\nu_I},
 \qquad \xi_I^{\rm p}=\curl\Theta_I^{\rm p}, \qquad \xi^{\rm p}=\sum_I\xi_I^{\rm p}.
\end{equation}
This definition preserves the disjoint supports of the local contributions and makes $\xi^{\rm p}$ divergence free,
with the coordinate directions oriented so that the leading curl points
along $\frac{\partial}{\partial y_I^{\zeta_I}}$. We first estimate this LDF and
then explain why the exact field variations produced by its flow, namely \eqref{eq:sfieldsintro}, have the same size bounds as their first-order Taylor expansion at $s=0$; see \eqref{eq:firststateintro}. 

\subsection{Estimates on the perturbation: path calculus.}\label{ssec:pathcalculusintro} Fix an index $I=(n,J,\zeta)$, with local background $(v_{\ell,I},B_{\ell,I})=(v_{\ell,n},B_{\ell,n})$, and let $X_s^{\rm p}$ be the flow of $\xi^{\rm p}$ at fixed physical time. Write
$\mathcal U_s^{\rm p}=(X_s^{\rm p})_*$ for its pushforward action and
\begin{equation}\label{eq:fieldsrepresentationintro0}
    w_{I,s}=(\partial_tX_s^{\rm p})\circ(X_s^{\rm p})^{-1}
                   +\mathcal U_s^{\rm p}v_{\ell,I}-v_{\ell,I},
 \qquad b_{I,s}=\mathcal U_s^{\rm p}B_{\ell,I}-B_{\ell,I}
\end{equation}
for the induced field variations on the support of the corresponding local contribution, before the corrector pushforward; recall \eqref{eq:sfieldsintro} above. We think of $s\mapsto X_s^{\rm p}$ as a path in the volume-preserving diffeomorphisms starting at the identity for each fixed physical time. The first-variation identities of Lemma~\ref{lem:path-transport}
give the representations
\begin{equation}\label{eq:fieldsrepresentationintro1}
    w_{I,s}=
      \int_0^s\mathcal U_r^{\rm p}\mathcal D_{t,I}\xi_I^{\rm p}\,\dd r, \qquad
 b_{I,s}=\int_0^s\mathcal U_r^{\rm p}\mathcal L_{B_{\ell,I}}\xi_I^{\rm p}\,\dd r.
\end{equation}
This can be seen as an average along the adjoint action, but we will not pursue this point of view here.  To estimate these integrals, we first bound the LDF. We then introduce norms that track the spatial, material and magnetic derivatives of fields expressed in the local frames of \S\ref{ssec:frameintro}, and prove that these norms are controlled under the principal pushforward. This applies in particular to $\mathcal D_{t,I}\xi_I^{\rm p}$ and $\mathcal L_{B_{\ell,I}}\xi_I^{\rm p}$ in \eqref{eq:fieldsrepresentationintro1}, and gives a path calculus adapted to the anisotropy of the construction.

\textbf{Estimates on the LDF.} Write $\lambda=\lambda_{q+1}$ and
\[
 f_I=\lambda^{-1}\mathfrak a_I\varphi(\lambda y_I^{k_I}),
 \qquad \Theta_I^{\rm p}=f_I\,dy_I^{\nu_I}.
\]
The curl product rule and volume preservation give
\[
 \begin{aligned}
 \xi_I^{\rm p}&=\nabla f_I\times dy_I^{\nu_I}\\
 &=\mathfrak a_I\varphi'(\lambda y_I^{k_I})\frac{\partial}{\partial y_I^{\zeta_I}}
   +\lambda^{-1}\varphi(\lambda y_I^{k_I})\,
       \nabla\mathfrak a_I\times dy_I^{\nu_I}.
 \end{aligned}
\]
The first term is orthogonal to the phase gradient. The second is the
correction required by localization; its derivative falls on the slow
amplitude, so it is smaller by $(\ell\lambda)^{-1}$.

Since the local background fields
are divergence-free, their Lie derivatives commute with curl,
viewed as exterior differentiation from one-forms to vector fields.
From the transport properties highlighted earlier, we then deduce  
\[
 \mathcal D_{t,I}^j\mathcal L_{B_{\ell,I}}^m\xi_I^{\rm p}
 =\curl\!\left[\lambda^{-1}(D_{t,I}^jD_{B,I}^m\mathfrak a_I)
                   \varphi(\lambda y_I^{k_I})dy_I^{\nu_I}\right].
\]
Thus neither transport operator differentiates the fast spatial oscillation. One curl cancels the prefactor
$\lambda^{-1}$, and each further spatial derivative costs at most
$\lambda$, since $\ell^{-1}\le\lambda$. The primitive bounds give
\[
 \|\mathcal D_{t,I}^j\mathcal L_{B_{\ell,I}}^m\xi_I^{\rm p}\|_r
 \lesssim \tau_a\delta_{q+1}^{1/2}
                  \lambda^r\tau_a^{-j}\lambda_\parallel^m,
\]
for $j+m\le j_1$ within the supplied derivative range; a spatial
H\"older exponent $\alpha$ adds $\lambda^\alpha$. This proves the
vector estimate directly from the scalar amplitude estimate. Proposition~\ref{principal:generator-estimates} and
Lemma~\ref{high:component-lie-conversion} give the derivative bounds
and their conversion to transport derivatives $D_{t,I},D_{B,I}$.

\textbf{Path calculus.} We first reduce the problem of estimating objects of the form $\mathcal U_r^{\rm p}F$ for a vector field $F=\curl(g\,dy_I^{\nu_I})$, with $g$ vanishing outside the support of the local contribution, to a scalar case and then deal with that through a class norm and a Gr\"onwall argument.  

\textit{Reduction to the scalar case.} The Euclidean gradient of $\xi_I^{\rm p}$ can be as large as
$\tau_a\lambda\delta_{q+1}^{1/2}$, so an estimate exponential in this
quantity would be useless. However, $\xi_I^{\rm p}=\nabla f_I\times dy_I^{\nu_I}$ gives
\[
 \xi_I^{\rm p}\cn f_I=\xi_I^{\rm p}\cn y_I^{\nu_I}=0,\qquad
 \mathcal L_{\xi_I^{\rm p}}(g\,dy_I^{\nu_I})=(\xi_I^{\rm p}\cn g)\,dy_I^{\nu_I}.
\]
The
flow $X_s^{\rm p}$ preserves $f_I$, $dy_I^{\nu_I}$ and $df_I$, and volume preservation gives
\begin{equation}\label{eq:generatorcommuteintro}
    \mathcal U_s^{\rm p}\curl(g\,dy_I^{\nu_I})
       =\curl\bigl((\mathcal U_s^{\rm p}g)\,dy_I^{\nu_I}\bigr).
\end{equation}
For vector fields of this form, it therefore suffices to estimate a
transported scalar and then take one curl.

\textit{Class norm.} The oscillation coordinate $y_I^{k_I}$ itself is not exactly invariant under this flow:
only the curl corrector term in $\xi_I^{\rm p}$ can move it, and the curl formula
above gives
\begin{equation}\label{eq:phaseerrorprelimintro}
    \xi_I^{\rm p}\cn(\lambda y_I^{k_I})
       =-\varphi(\lambda y_I^{k_I})\partial_{y_I^{\zeta_I}}\mathfrak a_I.
\end{equation}
The derivative on the right is slow. To measure this, distinguish the
oscillatory coordinate $y_I^{k_I}$ from $y_I^\perp=(y_I^{\nu_I},y_I^{\zeta_I})$ and use
\[
 \lambda^{-1}\partial_{y_I^{k_I}},\qquad
 \ell\partial_{y_I^{\nu_I}},\qquad \ell\partial_{y_I^{\zeta_I}},\qquad
 \tau_aD_{t,I},\qquad \lambda_\parallel^{-1}D_{B,I}.
\]
These derivatives commute in the chart of the local background. For derivative order $N$
and mixed transport derivative order $J$, the corresponding norm is
\[
 \|u\|_{N,J}
 =\max_{\substack{p+|\sigma|+h+m\le N\\h+m\le J}}
   \lambda^{-p}\ell^{|\sigma|}\tau_a^h\lambda_\parallel^{-m}
   \|\partial_{y_I^{k_I}}^p\partial_{y_I^\perp}^\sigma D_{t,I}^hD_{B,I}^mu\|_\infty.
\]
The supremum in $\|\cdot\|_\infty$ is taken over time and space on the support of the local contribution indexed by $I$. The weights
remove the expected cost of each derivative: the two nonoscillating
coordinate derivatives cost only $\ell^{-1}$. This is important because the objects we estimate then become bounded independently of the iteration index $q$, so the exponential factor in the Gr\"onwall estimate is uniformly bounded and can be absorbed into the implicit constant. These are norms of scalar coefficients; the physical
vector fields are recovered using their chart frames and curl.

\textit{Gr\"onwall estimates in the class norm.} Motivated by \eqref{eq:generatorcommuteintro} and \eqref{eq:phaseerrorprelimintro} we want to estimate in the norm $\|\cdot\|_{N,J}$ solutions $u_s$ of the transport equation on the support of each local contribution:
\[
    (\partial_s+\xi_I^{\rm p}\cn)u_s=g_s.
\]
The coefficients of $\xi_I^{\rm p}\cn$ in the three normalised spatial
directions are
\[
 -\varphi\partial_{y_I^{\zeta_I}}\mathfrak a_I,\qquad 0,\qquad
 \ell^{-1}\mathfrak a_I\varphi'
       +(\ell\lambda)^{-1}\varphi\partial_{y_I^{k_I}}\mathfrak a_I.
\]
Their $\|\cdot\|_{N,J}$ norms are bounded by
$C\tau_a\ell^{-1}\delta_{q+1}^{1/2}\ll1$, using one additional
spatial derivative of $\mathfrak a_I$. This makes precise why the large
spatial gradient is harmless here: differentiation with respect to $y_I^{k_I}$ is
normalised by $\lambda^{-1}$, while motion across surfaces of constant phase
already has that small factor in the operator $\xi_I^{\rm p}\cn$. 

Commuting the normalised derivatives with
$(\partial_s+\xi_I^{\rm p}\cn)u_s=g_s$ differentiates one of these small
coefficients and leaves at most $N$ derivatives on $u_s$. Integration
along characteristics and Gr\"onwall therefore give
\[
 \sup_{0\le s\le1}\|u_s\|_{N,J}
 \lesssim\|u_0\|_{N,J}+\int_0^1\|g_s\|_{N,J}\,\dd s,
\]
with no dependence on the large spatial gradient. On scalars, differentiation along $\xi_I^{\rm p}$ satisfies
\[
 \|\xi_I^{\rm p}\cn u\|_{N,J}
 \lesssim\tau_a\ell^{-1}\delta_{q+1}^{1/2}\|u\|_{N+1,J}.
\]
The additional derivative is spatial; the mixed transport derivative order stays $J$.
This is the small factor later used in the variations of the quadratic
tensor where an additional expansion is needed. Lemma~\ref{principal:finite-transport} proves these estimates.

\textit{Example: the phase error under $\mathcal U_s^{\rm p}$.} If $Y_{I,s}$ is the principal flow in chart coordinates, the phase
displacement
\[
 \psi_{I,s}=\lambda(y_I^{k_I}\circ Y_{I,s}^{-1}-y_I^{k_I})
\]
satisfies, by \eqref{eq:phaseerrorprelimintro},
\[
 (\partial_s+\xi_I^{\rm p}\cn)\psi_{I,s}
   =\varphi(\lambda y_I^{k_I})\partial_{y_I^{\zeta_I}}\mathfrak a_I,
 \qquad \psi_{I,0}=0.
\]
The scalar estimate bounds $\|\psi_{I,s}\|_{N,J}$ by
$C\tau_a\ell^{-1}\delta_{q+1}^{1/2}$. The chain rule then controls
$\varphi(\lambda y_I^{k_I}+\psi_{I,s})$ in the same norm, so the transported
profile retains the spatial oscillation scale. The transverse displacement, divided by $\ell$, satisfies the same bound.
Corollary~\ref{principal:phase-flow}
proves these scalar, phase and frame estimates and their path-integral
versions.

\textbf{Field estimates.} Rewriting \eqref{eq:fieldsrepresentationintro1} by means of \eqref{eq:generatorcommuteintro} we obtain
\[
 \begin{aligned}
 w_{I,s}&=\curl\!\left[
      \left(\int_0^s\mathcal U_r^{\rm p}(D_{t,I}f_I)\,\dd r\right)dy_I^{\nu_I}\right], \qquad b_{I,s}=\curl\!\left[
      \left(\int_0^s\mathcal U_r^{\rm p}(D_{B,I}f_I)\,\dd r\right)dy_I^{\nu_I}\right],\\
 D_{t,I}f_I&=\lambda^{-1}(\alpha_I(t/\tau_a)a_I+a_I^{\rm c})\varphi(\lambda y_I^{k_I}),\qquad
 D_{B,I}f_I=\lambda^{-1}(D_{B,I}\mathfrak a_I)\varphi(\lambda y_I^{k_I}).
 \end{aligned}
\]
The scalar path estimate preserves the source norms, with amplitudes
$\lambda^{-1}\delta_{q+1}^{1/2}$ and
$\lambda^{-1}\tau_a\lambda_\parallel\delta_{q+1}^{1/2}$. The scalar estimates and the invariant chart coframe therefore give
\[
 \sup_{0\le s\le1}\|w_{I,s}\|_r\lesssim
                  \lambda^r\delta_{q+1}^{1/2},\qquad
 \sup_{0\le s\le1}\|b_{I,s}\|_r\lesssim
                  \lambda^r\tau_a\lambda_\parallel\delta_{q+1}^{1/2}.
\]
Estimates for material and magnetic derivatives follow from the same norm with $J>0$.
The extra magnetic derivative in $D_{B,I}\mathfrak a_I$ has a loss at the
highest transport derivative order. The passage to derivative bounds for the next iterate and higher
spatial orders is detailed in Section~\ref{ssec:derivative-closure}; its
role in reaching the exponent $1/3$ is explained in
Subsection~\ref{ssec:overview-finite-derivatives}.
Proposition~\ref{principal:nonprincipal-field-smallness}
states the full bounds, including the corrector pushforward. 

\textbf{Background gaps.} For a general vector $F$ on which we have no geometric control, and in particular for the background differences in \eqref{eq:gapsintro}, the principal
pushforward in chart coordinates is
\[
 (Y_{I,s})_*F=(D_{y_I}Y_{I,s}\circ Y_{I,s}^{-1})(F\circ Y_{I,s}^{-1}),\qquad
 \|D_{y_I}Y_{I,s}-\IId\|_0\lesssim
       \tau_a\lambda_{q+1}\delta_{q+1}^{1/2}.
\]
The product rule introduces at most the factor
$1+\tau_a\lambda_{q+1}\delta_{q+1}^{1/2}$, with no additional
derivative of $F$. The smallness gained in the mollification error absorbs this large factor
for the background gaps, including the higher derivative orders covered by the hypotheses;
see Corollary~\ref{principal:vector-transport}.

\textbf{Quadratic interaction estimates: the antiderivative and magnetic errors.} The leading velocity and magnetic coefficients in \eqref{eq:fieldsrepresentationintro1} are
$D_{t,I}\mathfrak a_I$ and $D_{B,I}\mathfrak a_I$. The corresponding coefficient in the quadratic momentum tensor
\[
    w_{I,1}\otimes w_{I,1}-b_{I,1}\otimes b_{I,1}
\]
of \eqref{eq:momentumvariationintro} is therefore
\begin{equation}\label{eq:leadingquadraticintro}
    (D_{t,I}\mathfrak a_I)^2-(D_{B,I}\mathfrak a_I)^2
 =\alpha_I^2a_I^2+2\alpha_Ia_Ia_I^{\rm c}+(a_I^{\rm c})^2-(D_{B,I}\mathfrak a_I)^2.
\end{equation}
Since $\varphi(u)=\sqrt2\sin u$, the spatial average of $(\varphi')^2$
is one. Thus the first term supplies the tensor in
\eqref{eq:decompositionintro} after averaging the fast space and time
profiles. The middle terms
form the primitive error $R^{\rm p}_{\rm time}$ and the negative
magnetic square forms $R^{\rm p}_{\rm mag}$. Their smallness determines the primitive expansion order $j_0$ and the admissible
magnetic exponent $\gamma_\parallel$; see
Subsection~\ref{ssec:overview-finite-derivatives}.

For the principal path integrals, \eqref{principal:exact-square}
separates the initial quadratic tensor, its first Lie variation, and
second-order terms. The contributions estimated below enter the principal
quadratic stress $R^{\rm p}_{\rm quad}$, defined in
\eqref{mom:N-remainder}. The leading oscillatory tensors have vector factors orthogonal to the
phase gradient, since $\frac{\partial}{\partial y_I^{\zeta_I}}\cdot\nabla y_I^{k_I}=0$.
Their divergence therefore differentiates only slow coefficients and
frames; periodic antiderivatives in Lemma~\ref{osc:parametrix} then
gain $(\ell\lambda_{q+1})^{-1}$. This inversion is also needed for
the tangential part of the first Lie variation, whose Fourier modes in the fast phase
are $\pm1,\pm3$ by \eqref{principal:first-lie-mean}. Products containing
a localisation gradient corrector term already contain the factor $\lambda_{q+1}^{-1}$ and are estimated
directly. Thus the complete first variation, including these products,
has the required gain. The second variation and covariance contain
two applications of the principal LDF to scalars and are bounded directly by
$\tau_a^2\ell^{-2}\delta_{q+1}^2$; see
\eqref{principal:quadratic-second-bound}. Applying oscillatory inversion
only to the initial tensor would leave the first variation too large
near $1/3$ (Remark~\ref{iter:first-order-warning}).

\subsection{Time-averaging and the Galbrun equation}\label{ssec:grintro}
The space-time average of the leading quadratic term cancels $R_q$
through \eqref{eq:decompositionintro}, up to mollification errors. Spatial
oscillations are treated by oscillatory inverse divergence. The remaining tensor
source is
\[
 \mathsf F_n
 =\sum_{I=(n,J,\zeta)}[1-\alpha_I(t/\tau_a)^2]a_I^2
    \partial_{y_I^{\zeta_I}}\otimes\partial_{y_I^{\zeta_I}}.
\]
It is slow in space, while each temporal profile $1-\alpha_I^2$
has zero mean. We use these temporal cancellations to construct a
small displacement that compensates for $\ddiv\mathsf F_n$. As in \cite{GR}, we construct this displacement before the principal perturbation.

On a fixed time interval, write $(v,B)=(v_{\ell,n},B_{\ell,n})$
and suppress the index. The Lie-algebra linearization of
momentum is
\[
 \mathcal F_{v,B}(\mathcal D_t\xi,\mathcal L_B\xi)
   =(D_t^2-D_B^2)\xi-\DD(D_tv-D_BB)\xi.
\]
Writing $\xi=\curl\Theta$, the Galbrun operator on one-form
potentials is
\[
 \begin{aligned}
 \mathscr G\Theta
 &=(D_t^2-D_B^2)\Theta+\mathscr H_tD_t\Theta
       +\mathscr H_BD_B\Theta+\mathscr H_2\Theta,\\
 \curl\mathscr G\Theta
 &=\mathbb P\mathcal F_{v,B}(\mathcal D_t\xi,\mathcal L_B\xi).
 \end{aligned}
\]
Here $\mathbb P$ is the Leray projection, and
$\mathscr H_t,\mathscr H_B,\mathscr H_2$ are operators of spatial
order zero formed from $\nabla v$, $\nabla B$ and
$\nabla(\ddiv R-\nabla p)$;
see \eqref{galbrun:separate-lower-operators}--\eqref{mom:H2}.
Thus we solve $\mathscr G\Theta=\mathscr T\mathsf F$ with vanishing
initial data, where
$\mathscr T=\curl^{-1}\mathbb P\ddiv$ has order zero.

For a mean-zero periodic profile, the superscript $[j]$ here denotes its
$j$th mean-zero periodic primitive. The first explicit potential is
\[
 \Theta_{\mathrm{exp},1}
 =\tau_a^2\sum_I(1-\alpha_I^2)^{[2]}(t/\tau_a)
   \mathscr T\bigl(a_I^2
    \partial_{y_I^{\zeta_I}}\otimes
    \partial_{y_I^{\zeta_I}}\bigr),
\]
and direct differentiation gives
\[
 \begin{aligned}
 \mathscr G\Theta_{\mathrm{exp},1}
 &=\mathscr T\mathsf F+E_1,\\
 E_1
 &=\sum_I\Bigl[
   \tau_a(1-\alpha_I^2)^{[1]}(t/\tau_a)
       (2D_t+\mathscr H_t)\\
 &\hspace{3.2em}
   +\tau_a^2(1-\alpha_I^2)^{[2]}(t/\tau_a)\mathscr G
   \Bigr]\mathscr T\bigl(a_I^2
    \partial_{y_I^{\zeta_I}}\otimes
    \partial_{y_I^{\zeta_I}}\bigr).
 \end{aligned}
\]
Cancelling successive residual terms by further temporal primitives
gives the iterative splitting
\[
 \begin{aligned}
 \Theta&=\Theta_{\mathrm{exp},N}+Z_N,\\
 \mathscr G\Theta_{\mathrm{exp},N}&=\mathscr T\mathsf F+E_N,
 \qquad \mathscr GZ_N=-E_N,\\
 Z_N(t_-)&=D_tZ_N(t_-)=0.
 \end{aligned}
\]
Every term of $\Theta_{\mathrm{exp},N}$ is supported in time where
the slow coefficients are supported, and each cancellation gains
$\varepsilon_\tau=\tau_a/\tau_c$. The exact solution $\Theta$ does
not depend on $N$. The characteristic equations for
$D_tZ_N\pm D_BZ_N$, followed by one integration for $Z_N$, give
over a time interval of length $O(\tau_c)$
\[
 \|E_N\|_\alpha
 +\tau_c^{-2}\|Z_N\|_\alpha
 +\tau_c^{-1}\|D_tZ_N\|_\alpha
 \lesssim_N\ell^{-\alpha}\varepsilon_\tau^N\delta_{q+1},
\]
on the derivative ranges in Appendix~\ref{sec:galbrun}.
In the source region, $N=2$ already gives the potential size
$\ell^{-\alpha}\tau_a^2\delta_{q+1}$. Magnetic derivatives act only
on slow coefficients and cost only $\lambda_\parallel$.

\begin{remark}
The expansion $\Theta_{\mathrm{exp},N}$ uses the same antiderivative
mechanism as the coefficient construction in \eqref{eq:telescopeintro}:
successive errors are cancelled by further temporal primitives, and
each cancellation gains $\varepsilon_\tau$. One could therefore use
this expansion to construct an approximate Galbrun solution, retaining
$\mathcal R\curl E_N$ as an additional stress error. We instead solve
the Galbrun equation exactly and use the expansion to estimate its
solution and the error introduced by the time cutoff.
In the two-dimensional Euler setting, these estimates can also be
used to reprove the result of Giri--Radu~\cite{GR} without a Newton
iteration.
\end{remark}

The exact solution can persist after the source vanishes. To keep the
perturbation compact in time, restore the window index and set
\[
 \begin{gathered}
 \Theta_n^{\rm c}=\tilde\eta_n\Theta_n,\qquad
 \xi_n^{\rm c}=\curl\Theta_n^{\rm c},\\
 \Theta^{\rm c}=\sum_n\Theta_n^{\rm c},\qquad
 \xi^{\rm c}=\curl\Theta^{\rm c},
 \end{gathered}
\]
where $\tilde\eta_n=1$ near the source support and changes only in
source-free cut off regions, with
$|\tilde\eta_n^{(h)}|\lesssim_h\tau_c^{-h}$. In either cut off region the
explicit part vanishes. Thus, suppressing $n$ again, the cutoff error
depends only on the remainder:
\[
 \begin{aligned}
 [\mathscr G,\tilde\eta]\Theta
 &=[\mathscr G,\tilde\eta]Z_N\\
 &=\tilde\eta''Z_N+2\tilde\eta'D_tZ_N
                  +\tilde\eta'\mathscr H_tZ_N.
 \end{aligned}
\]
The factors $\tau_c^{-2}$ and $\tau_c^{-1}$ cancel the factors obtained
from the two integrations over an interval of length $O(\tau_c)$ above. Since $\mathcal R\curl$ has order zero,
\begin{equation}\label{eq:Rcutintro}
    \|R^{\mathrm{cut}}\|_0
 \lesssim_N\ell^{-\alpha}\varepsilon_\tau^N\delta_{q+1},
 \qquad 1\le N\le j_1-2,
\end{equation}
with the ranges of material and magnetic derivatives in
\eqref{galbrun:genericoutgoingtailA}. The incoming contribution is
exactly zero by uniqueness; the outgoing contribution gains
$\varepsilon_\tau^N$ for every admissible $N$. These slow cut off regions differ from the
gaps between fast principal profiles: $1-\alpha_I^2=1$ where
$\alpha_I=0$, so the corrector need not vanish there.
Figure~\ref{fig:window} shows the forcing and outgoing cutoff-error supports.

Taking a curl costs the slow spatial frequency $\lambda_q$.
Consequently, for each local background, with its material and magnetic Lie derivatives,
\[
 \begin{aligned}
 \|\xi_n^{\rm c}\|_0
   &\lesssim\ell^{-\alpha}\tau_a^2\lambda_q\delta_{q+1},\\
 \|\mathcal D_{t,n}\xi_n^{\rm c}\|_0
   &\lesssim\ell^{-\alpha}\tau_a\lambda_q\delta_{q+1},\\
 \|\mathcal L_{B_{\ell,n}}\xi_n^{\rm c}\|_0
   &\lesssim\ell^{-\alpha}\tau_a^2\lambda_\parallel\lambda_q\delta_{q+1}.
 \end{aligned}
\]
Each corrector first variation is smaller than its principal counterpart
by $\ell^{-\alpha}\tau_a\lambda_q\delta_{q+1}^{1/2}\ll1$.
Moreover,
\[
 \|\DD\xi^{\rm c}\|_0\lesssim
 \ell^{-\alpha}\tau_a^2\lambda_q^2\delta_{q+1}\ll1,
\]
so the corrector flow is close to the identity in $C^1$.
Proposition~\ref{gg:physical} and
\eqref{prep:corrector-final-comparison} give the precise field
estimates and their comparison with the actual background.

\begin{figure}[htbp]
\centering
\begin{tikzpicture}[x=1.6cm,y=0.8cm,font=\small]
 \fill[gray!8] (-2,0) rectangle (-1,1.48);
 \fill[gray!8] (1,0) rectangle (2,1.48);
 \fill[gray!18] (-1,0) rectangle (1,1.48);
 \fill[gray!40] (1.12,0) rectangle (1.85,1.48);
 \draw[gray!50,densely dotted] (-2,0) -- (-2,1.48);
 \draw[gray!50,densely dotted] (2,0) -- (2,1.48);
 \draw[->] (-4.18,0) -- (4.28,0);
 \foreach \x in {-4,-3,-2,-1,0,1,2,3,4}
   \draw (\x,0.06) -- (\x,-0.06) node[below=2pt] {$\x$};
 \node[anchor=east] at (4.2,-0.8) {$(t-t_n)/\tau_c$};
 \draw[thick]
   (-2,0) -- (-1.85,0)
   .. controls (-1.60,0) and (-1.48,1.2) .. (-1.12,1.2)
   -- (1.12,1.2)
   .. controls (1.48,1.2) and (1.60,0) .. (1.85,0) -- (2,0);
 \node[anchor=west] at (1.62,1.2) {$\tilde\eta_n$};
 \node[align=left,anchor=west] at (2.15,0.66)
   {outgoing\\cutoff error};
 \draw[->] (2.10,0.42) -- (1.62,0.42);
 \draw[thick,dashed]
   (-0.73,0) .. controls (-0.53,0) and (-0.46,1.2) .. (-0.25,1.2)
   -- (0.25,1.2)
   .. controls (0.46,1.2) and (0.53,0) .. (0.73,0);
 \node at (0,0.52) {$\eta_n$};
 \draw[thick] (-3,0) -- (-3,0.63);
 \node[align=center,anchor=south] at (-3,0.72)
   {zero initial data\\$t=t_n-3\tau_c$};
 \draw[densely dotted,thick] (-4,1.85)
   .. controls (-1.5,1.85) and (1.5,1.79) .. (4,1.77);
 \node[anchor=west,fill=white,inner sep=2pt] at (-3.85,1.94)
   {$\varrho$\quad(scale $\delta_q^{1/2}$)};
 \draw[thick] (-3,2.48) -- (3,2.48);
 \fill (-3,2.48) circle (1.3pt);
 \fill (3,2.48) circle (1.3pt);
 \node[fill=white,inner sep=3pt] at (0,2.48)
   {$\mathscr G_n\Theta_n=\mathscr T\mathsf F_n$};
 \draw[decorate,decoration={brace,amplitude=4pt}]
   (-4,3.03) -- (4,3.03)
   node[midway,above=5pt] {$I_n^{\rm bg}$};
 \draw[line width=2pt] (-0.73,-0.65) -- (0.73,-0.65);
 \node[below=4pt] at (0,-0.65)
   {forcing: $[1-\alpha_I(t/\tau_a)^2]a_I^2$};
 \draw[decorate,decoration={brace,mirror,amplitude=4pt}]
   (-1,-1.42) -- (1,-1.42)
   node[midway,below=5pt] {$I_n^{\rm src}$};
 \draw[decorate,decoration={brace,mirror,amplitude=4pt}]
   (-2,-2.10) -- (2,-2.10)
   node[midway,below=5pt] {$I_n^{\rm col}$};
\end{tikzpicture}
\caption{Source and cutoff-error supports on a Galbrun window.
The bar marks the slow support of the forcing, contained in
$\supp(\varrho\eta_n)\Subset I_n^{\rm src}$; its fast factors
$1-\alpha_I(t/\tau_a)^2$ need not vanish between principal profiles.
The darker outgoing transition contains the support of
$[\mathscr G_n,\tilde\eta_n]\Theta_n$, where only the remainder $Z_N$
survives. The incoming cutoff error is zero. The Cauchy problem with zero initial data starts
at $t_n-3\tau_c$; the windows and outer cutoff are specified in
\eqref{prep:slow-windows} and \eqref{prep:outer-cutoff}.}
\label{fig:window}
\end{figure}

\subsection{Final ansatz}
The potential in \eqref{eq:principalpreparation} in fact defines only the initial
principal LDF $\xi_0^{\rm p}$. Let $X_s^{\rm c}$ and
$X_s^{\rm p}$ be the flows of $\xi^{\rm c}$ and $\xi_0^{\rm p}$.
We compose them in the order
\[
 X_s=X_s^{\rm c}\circ X_s^{\rm p},\qquad
 \xi_s=\xi^{\rm c}+\mathcal U_s^{\rm c}\xi_0^{\rm p},\qquad
 \mathcal U_s^{\rm c}=(X_s^{\rm c})_*.
\]
Thus the corrector flow transports the entire principal potential:
\[
 \Theta_{I,s}^{\rm p}
 =\lambda_{q+1}^{-1}(\mathcal U_s^{\rm c}\mathfrak a_I)
   \varphi\bigl(\lambda_{q+1}\mathcal U_s^{\rm c}y_I^{k_I}\bigr)
   d(\mathcal U_s^{\rm c}y_I^{\nu_I}),\qquad
 \xi_s^{\rm p}=\sum_I\curl\Theta_{I,s}^{\rm p}.
\]

In a chart of the local background, suppressing the local background indices, denote the $X_s^{\rm c}$ deformation of $(v,B)$ by
$(v_s^{\rm c},B_s^{\rm c})$. Intertwining of Lie derivatives with the corrector pushforward gives the exact identities
\[
 \begin{aligned}
 D_{t,s}^{\rm c}(\mathcal U_s^{\rm c}y_I^{k_I})
   &=\mathcal U_s^{\rm c}D_ty_I^{k_I}=0,\\
 D_{B_s^{\rm c}}(\mathcal U_s^{\rm c}y_I^{k_I})
   &=\mathcal U_s^{\rm c}D_By_I^{k_I}=0.
 \end{aligned}
\]
Keeping the original phase would instead produce
\[
 \begin{aligned}
 D_{t,s}^{\rm c}\varphi(\lambda_{q+1}y_I^{k_I})
 &=\lambda_{q+1}\varphi'(\lambda_{q+1}y_I^{k_I})
            (v_s^{\rm c}-v)\cn y_I^{k_I},\\
 D_{B_s^{\rm c}}\varphi(\lambda_{q+1}y_I^{k_I})
 &=\lambda_{q+1}\varphi'(\lambda_{q+1}y_I^{k_I})
            (B_s^{\rm c}-B)\cn y_I^{k_I}.
 \end{aligned}
\]
The pushforward removes these terms with their large factor
$\lambda_{q+1}$. The same intertwining applies to the coframe and
all Lie derivatives of the principal LDF; see
\eqref{direct:differentiated-transport}. By
\eqref{path:factorized-field-increments}, the principal increments
relative to $(v_s^{\rm c},B_s^{\rm c})$ are precisely
$\mathcal U_s^{\rm c}$ applied to the path integrals
\eqref{eq:fieldsrepresentationintro1}. Hence the preceding scalar
estimates apply in the chart transported by the corrector flow. This avoids the additional chart construction used after the corrector perturbation in \cite{gmzl2026c15convexintegrationsolutions}. The small principal phase drift
is measured from $\lambda_{q+1}\mathcal U_s^{\rm c}y_I^{k_I}$,
as in \eqref{principal:phase-corrector-composition}.

The corrector flow preserves the disjointness of the principal LDF. In corrector coordinates, $\xi_0^{\rm p}$ vanishes in the original
cut-off regions, so its trajectories cannot cross between the supports of distinct local contributions
(Lemma~\ref{mom:separation}). This eliminates quadratic interactions between principal perturbations with distinct indices $I$.

\subsection{Parameter constraints}\label{ssec:overview-finite-derivatives}
To approach $1/3$, the fast-time exponent $\gamma_a$ and the
spatial mollification exponent $\gamma_\ell$ must be arbitrarily small, with
$b>1$ close to one. The constructions above compensate for this through the approximation gains from vanishing moments in the mollifier (\S\ref{ssec:deeppreparation}) and repeated gains of $\tau_a/\tau_c$ (\S\ref{ssec:principalperturbationoverview}, \S\ref{ssec:grintro}). This leaves the following main errors.

\textit{Temporal errors.}
The principal transport stress
$R^{\mathrm p}_{\rm tr}$ imposes the
upper bound on $\gamma_a$. In the linearised force
\eqref{eq:momentumvariationintro}, differentiating the fast profile
costs $\tau_a^{-1}$, while the potential
\eqref{eq:principalpreparation} supplies $\lambda_{q+1}^{-1}$.
Its contribution relative to each local background therefore satisfies the bound
\eqref{principal:natural-stresses}, of size
\[
 \frac{\delta_{q+1}^{1/2}}{\lambda_{q+1}\tau_a}
 =\varepsilon_{q+1}^{1-\beta-\gamma_a-\gamma_\ell}\delta_{q+1}.
\]
The competing lower bound comes from the antiderivative stress
$R^{\mathrm p}_{\rm time}$ and the cutoff stress $R^{\mathrm{cut}}$ in the Galbrun construction. Both constructions gain powers of the ratio of the fast and slow time scales, so the two stresses satisfy the same bound
$\varepsilon_{q+1}^{j_0\gamma_a}\delta_{q+1}$; see
\eqref{aniso:stress-material-error} in \S\ref{ssec:pathcalculusintro} and \eqref{eq:Rcutintro} in \S\ref{ssec:grintro}.

The flow-smoothing part of $R^{\rm moll}$ has this size by the
order-$j_0$ approximation estimate \eqref{stress:mollification-zero}.
Thus we require
\[
 j_0\gamma_a>2b\beta,
 \qquad \gamma_a+\gamma_\ell<1-(2b+1)\beta.
\]
A large $j_0$ makes the temporal remainders sufficiently small without
forcing $\gamma_a$ away from zero.

\textit{Quadratic path errors.}
The path calculus of Subsection~\ref{ssec:pathcalculusintro} separates
two kinds of contributions to $R^{\mathrm p}_{\rm quad}$ relative to the local backgrounds. Spatial oscillation
and the first Lie variation are controlled by oscillatory inverse
divergence, with size $(\ell\lambda_{q+1})^{-1}\delta_{q+1}$.
The corresponding estimates follow from
\eqref{osc:principal-classes} and \eqref{principal:first-lie-mean}.
Differentiating the slow coefficients along the principal LDF
contributes a factor $\tau_a\ell^{-1}\delta_{q+1}^{1/2}$.
The second variation and covariance each contain two such factors,
giving size
\[
 (\tau_a\ell^{-1}\delta_{q+1}^{1/2})^2\delta_{q+1}=\varepsilon_{q+1}^{2(\gamma_a+\beta)}\delta_{q+1}.
\]
See \eqref{principal:quadratic-second-bound}. This also dominates the
corrector square $R^{\mathrm c}_{\rm quad}$ and the chart-transfer
stress $R^{\mathrm p}_{\rm chart}$, whose size is
$\tau_a^2\lambda_q^2\delta_{q+1}^2$, since $\lambda_q\le\ell^{-1}$.
The resulting restriction is
\begin{equation}\label{eq:quadraticpathintro}
    \gamma_a>(b-1)\beta
\end{equation}
matching exactly the bound for \cite[$R_{q+1}^{\rm Newton}$]{GR} in our notation. The inverse divergence treatment of the first Lie variation is essential: a direct size estimate would impose
$\gamma_a+\beta>2b\beta$ and reduce the limiting range to $1/4$;
see Remark~\ref{iter:first-order-warning}.

\textit{Magnetic square and charts.}
The magnetic stress $R^{\mathrm p}_{\rm mag}$ is the negative square
retained in \eqref{eq:leadingquadraticintro}. Phase and coframe
invariance leave the magnetic derivative acting on the slow primitive
of size $\tau_a\delta_{q+1}^{1/2}$, with differentiation cost $\lambda_\parallel$.
As shown in \eqref{aniso:stress-magnetic-error}, its size is
\begin{equation}\label{eq:magneticlowerintro}
    (\tau_a\lambda_\parallel)^2\delta_{q+1}
 =\varepsilon_{q+1}^{2(\gamma_a+\gamma_\parallel)}\delta_{q+1},
 \qquad \gamma_a+\gamma_\parallel>b\beta.
\end{equation}
On the other hand, the charts of Subsection~\ref{ssec:frameintro}
have radius $\lambda_\parallel^{-1}$. Their small-distortion condition
\eqref{part:flow-smallness} uses
\[
 \lambda_\parallel^{-1}[B_q]_1
 \lesssim\varepsilon_{q+1}^{\gamma_\ell+
       [\gamma_a-(b-1)\gamma_\parallel]/b},
\]
see \eqref{setup:magnetic-rate}. We impose $0<(b-1)\gamma_\parallel<\gamma_a$. Together with the
magnetic lower bound \eqref{eq:magneticlowerintro}, this requires precisely
$\gamma_a>(b-1)\beta$, matching the quadratic path error \eqref{eq:quadraticpathintro}.

\textit{Reaching $1/3$.}
Ignoring $\gamma_\ell$, the decisive interval is therefore
\[
 (b-1)\beta<\gamma_a<1-(2b+1)\beta,
\]
which is nonempty precisely when
\[
 1-(2b+1)\beta-(b-1)\beta=1-3b\beta>0,
\]
that is, $\beta<1/(3b)$. For any $\beta<1/3$, we
first choose $b>1$ sufficiently close to one, then admissible
$\gamma_a,\gamma_\parallel$, and finally $j_0$ large enough that
$j_0\gamma_a>2b\beta$. Increasing $j_1$ supplies the remaining
transport derivatives. At the highest transport derivative orders, the ratios of the
local derivative costs to those of the next iterate absorb the coarse
Galbrun and primitive losses, preserving the required smallness gains;
see Lemma~\ref{aniso:parameter-budgets}, in particular
\eqref{aniso:corrector-final-comparison} and
\eqref{aniso:collar-depth-comparison}. We then take $\gamma_\ell>0$
small enough to preserve the strict margins and finally fix the Schauder loss $\gamma_S$ sufficiently
small with $0<3\gamma_S<\gamma_\ell$. The mollification and background comparison estimates control the
spatial mollification error and the background-gap contributions,
including the large factor arising from the derivative matrix of the principal
flow; see \eqref{iter:accuracy-requirement} and
\eqref{iter:principal-gap-budgets}. The remaining comparisons of derivative costs
and derivative reserves are verified in Section~\ref{sec:parameter-choice}.

\subsection{Regularity along magnetic field lines}\label{ssec:parallelregintro}
The gain in Theorem~\ref{full:fieldline} comes from slower variation
along the magnetic field. In the charts of
Subsection~\ref{ssec:frameintro}, the fast spatial phase is constant
in the local magnetic direction. A magnetic derivative therefore
acts on the slow coefficients, without differentiating the fast
oscillation. The path calculus of
Subsection~\ref{ssec:pathcalculusintro} preserves this structure
under the deformation. These directional estimates give information
beyond the joint H\"older regularity of $v$ and $B$.

To explain the exponents, we write the estimates schematically using
$\gamma_v,\gamma_B$ from the corollary. At spatial frequency
$\lambda_q$, the smooth increments have amplitudes
\[
 \|v_q-v_{q-1}\|_0\lesssim\lambda_q^{-\gamma_v},
 \qquad
 \|B_q-B_{q-1}\|_0\lesssim\lambda_q^{-\gamma_B}.
\]
The inductive magnetic derivative bounds
\eqref{prep:old-mixed} assign each magnetic derivative the cost
$\lambda_q^{1-\gamma_B}$, corresponding to the size of the
magnetic gradient, rather than the full spatial differentiation cost $\lambda_q$.
In particular, they give schematically
\[
 \|B_q\cn v_q\|_0\lesssim\lambda_q^{1-\gamma_B-\gamma_v},
 \qquad
 \|B_q\cn B_q\|_0\lesssim\lambda_q^{1-2\gamma_B}.
\]

The curves in the corollary follow the limiting field $B$.
To pass to this direction, write
\[
 B\cn v_q=B_q\cn v_q+(B-B_q)\cn v_q.
\]
The magnetic tail has size $\|B-B_q\|_0\lesssim\lambda_q^{-\gamma_B}$,
while $[v_q]_1\lesssim\lambda_q^{1-\gamma_v}$. Thus the second
term has the same bound $\lambda_q^{1-\gamma_B-\gamma_v}$ as
the first. The same argument for $B_q$ gives
$\|B\cn B_q\|_0\lesssim\lambda_q^{1-2\gamma_B}$.
The improved directional bounds therefore hold along the actual
field lines and also apply to differences of successive approximations;
see \eqref{iter:limiting-magnetic-derivative}.

Fix $t$ and any $C^1$ curve $\Gamma'=B(t,\Gamma)$, and let
$h=|s-r|\le1$. For each smooth increment, there are two ways to
bound its change along the curve: use its amplitude, or integrate
its derivative along $B$. For the velocity and magnetic increments,
respectively, this gives
\[
 \min\{\lambda_q^{-\gamma_v},
             h\lambda_q^{1-\gamma_B-\gamma_v}\},
 \qquad
 \min\{\lambda_q^{-\gamma_B},
             h\lambda_q^{1-2\gamma_B}\}.
\]
Thus an oscillation at spatial scale $\lambda_q^{-1}$ changes
appreciably only over the longer field-line scale
$\lambda_q^{-(1-\gamma_B)}$. For a given $h$, the derivative
bound controls frequencies below the transition
$h\lambda_q^{1-\gamma_B}\sim1$, and the amplitude bound controls
those above it. Both sums are dominated by their terms nearest
this transition. Evaluating the amplitudes at the transition gives
\[
 \begin{aligned}
 |v(t,\Gamma(s))-v(t,\Gamma(r))|
     &\lesssim h^{\gamma_v/(1-\gamma_B)},\\
 |B(t,\Gamma(s))-B(t,\Gamma(r))|
     &\lesssim h^{\gamma_B/(1-\gamma_B)}.
 \end{aligned}
\]
The proof in Subsection~\ref{ssec:fieldlinesproof} chooses slightly
stronger construction exponents to obtain precisely the exponents
stated in the corollary. The chain rule is used only for the smooth
increments, so the argument applies to every such curve without
requiring uniqueness. Since $B$ is uniformly bounded and bounded
away from zero, arclength gives the same exponents.

As $(\gamma_v,\gamma_B)$ approaches $(1/4,1/2)$ within the
corollary's range, the field-line exponents approach $1/2$ for
$v$ and $1$ for $B$. In this way the magnetic field can be almost
Lipschitz along its own field lines while the stated joint magnetic
exponent approaches $1/2$.

\section{Inductive Scheme}\label{sec:inductive-scheme}

To prove Theorem~\ref{full:main}, we construct smooth solutions
$(v_q,B_q,p_q,R_q)$, $q\ge0$, of the relaxed system
\begin{equation}\label{iter:relaxed-system}
 \begin{cases}
 \partial_t v_q+\ddiv(v_q\otimes v_q-B_q\otimes B_q)
        +\nabla p_q=\ddiv R_q,\\
 \partial_t B_q+[v_q,B_q]=0,\\
 \ddiv v_q=\ddiv B_q=0.
 \end{cases}
\end{equation}
on $\mathbb T^3\times\R$, with $R_q$ symmetric. We normalize $p_q$ to
have zero spatial mean, so that it is uniquely determined by the
Poisson equation associated with \eqref{iter:relaxed-system}. Thus the
induction equation and the divergence constraints hold exactly at each
stage. The bounds below and the increment estimates of
Proposition~\ref{full:step} imply convergence of $v_q,B_q$ in the
required H\"older spaces and convergence of $R_q$ to zero.

\subsection{Inductive assumptions and iterative proposition}

For $q\ge0$, we define the frequencies, their ratios and the squared
increment amplitudes by
\[
 \begin{gathered}
 \lambda_q=a^{b^q},\qquad
 \varepsilon_0=a^{1/b-1},\qquad
 \varepsilon_{q+1}=\frac{\lambda_q}{\lambda_{q+1}},\\
 \delta_q=\lambda_q^{-2\beta},\qquad
 \delta_{B,q}=\delta_q\varepsilon_q^{2(\gamma_a+\gamma_\parallel)}.
 \end{gathered}
\]
Here $a>1$ will be taken sufficiently large, $b>1$ is close to one,
and $0<\beta<1/3$ determines the velocity regularity. The positive
exponents $\gamma_a,\gamma_\parallel$ determine the additional
magnetic regularity through this choice of amplitude. For the pair
at stage $q$ write
\[
 D_{t,q}=\partial_t+v_q\cn,\qquad D_{B_q}=B_q\cn,\qquad
 \mathcal D_{t,q}=\partial_t+\mathcal L_{v_q}.
\]
The Lie derivative conventions are recalled in
Section~\ref{ssec:path-conventions}. All norms are spatial norms,
uniform in physical time.

Fix a small Schauder-loss parameter $\gamma_S>0$ and set
$\alpha=(b-1)\gamma_S/b$. Choose a constant $C_0>1$, a bound $j_1$
on the sum of the material and magnetic derivative orders,
and a maximum derivative order $\rgood$, as specified in
Proposition~\ref{full:step}. We impose the following
bounds for nonnegative integers $r,h,k,m$ in the indicated ranges.

The ordinary field derivatives satisfy, for $h\le j_1+1$ and
$1\le r+h\le\rgood+4$,
\begin{equation}\label{prep:old-ordinary}
 \begin{aligned}
 \|\partial_t^h v_q\|_r&\le C_0^{h+1}\lambda_q^{r+h}\delta_q^{1/2},\\
 \|\partial_t^h B_q\|_r&\le C_0^{h+1}\lambda_q^{r+h}\delta_{B,q}^{1/2}.
 \end{aligned}
\end{equation}
The material and magnetic derivatives satisfy, for $k+m\le j_1$ and
$1\le r+k+m\le\rgood$,
\begin{equation}\label{prep:old-mixed}
 \begin{aligned}
 \|D_{t,q}^kD_{B_q}^m v_q\|_r
 &\le C_0^{k+m+1}\lambda_q^r\delta_q^{1/2}
       (\lambda_q\delta_q^{1/2})^k(\lambda_q\delta_{B,q}^{1/2})^m,\\
 \|D_{t,q}^kD_{B_q}^m B_q\|_r
 &\le C_0^{k+m+1}\lambda_q^r\delta_{B,q}^{1/2}
       (\lambda_q\delta_q^{1/2})^k(\lambda_q\delta_{B,q}^{1/2})^m.
 \end{aligned}
\end{equation}
For the pressure, with $r\ge1$, $k+m\le j_1$ and
$r+k+m\le\rgood$, we assume
\begin{equation}\label{prep:old-pressure}
 \|D_{t,q}^kD_{B_q}^m p_q\|_r
 \le C_0^{k+m+2}\lambda_q^r\delta_q
       (\lambda_q\delta_q^{1/2})^k(\lambda_q\delta_{B,q}^{1/2})^m.
\end{equation}
The stress satisfies, for $k+m\le j_1$ and
$r+k+m\le\rgood-\mathbf1_{\{k+m>0\}}$,
\begin{equation}\label{prep:old-stress}
 \|D_{t,q}^kD_{B_q}^mR_q\|_r+\|\mathcal D_{t,q}^k\mathcal L_{B_q}^mR_q\|_r
 \le\lambda_q^{r-\alpha}\delta_{q+1}
       (\lambda_q\delta_q^{1/2})^k(\lambda_q\delta_{B,q}^{1/2})^m,
\end{equation}
and, for $h\le j_1$ and $r+h\le\rgood$,
\begin{equation}\label{prep:old-stress-time}
 \|\partial_t^hR_q\|_r
 \le C_0^h\lambda_q^{r+h-\alpha}\delta_{q+1}.
\end{equation}
The pointwise bounds and time support are
\begin{equation}\label{full:pointwise-induction}
 \|v_q\|_0,\|B_q\|_0\le 2(1-\delta_q^{1/2}),\qquad
 |B_q|\ge\tfrac14(1+\delta_q^{1/2}),
\end{equation}
\begin{equation}\label{full:stress-support}
 \supp_tR_q\subset
 \bigl((1+\delta_q^{1/2})/2,\,5(1-\delta_q^{1/2})/2\bigr).
\end{equation}

The iterative step is given by the following proposition. We choose its
parameters in Section~\ref{sec:parameter-choice}, in the order
$\beta,b$, then $\gamma_a,\gamma_\parallel$, then the transport orders,
$\gamma_\ell$, $\gamma_S$, and finally the approximation orders and
maximum derivative orders. We then fix $C_0$ and take $a$ sufficiently
large. Each choice depends only on the preceding parameters, the fixed
profiles and the geometric data. In particular, the choices are uniform
in $q$ and in the solution of \eqref{iter:relaxed-system}. We use
throughout the construction the scale inequalities and derivative
bounds established in that section.

\begin{proposition}[The iterative step]\label{full:step}
Fix $0<\beta<1/3$ and $1<b<9/8$ with $3b\beta<1$, and let
$\gamma_a,\gamma_\parallel$ satisfy \eqref{iter:exponent-window}.
There exist integers $j_0,j_1$ such that the following holds for every
sufficiently small $\gamma_\ell\in(0,1)$ and, depending on
$\gamma_\ell$, every sufficiently small
$\gamma_S\in(0,\gamma_\ell/3)$. There are approximation orders
$m_0,\rfg$, a maximum derivative order $\rgood$, and constants
$C_0>1$, $a_0>1$ such that, for any fixed $a\ge a_0$ and every $q\ge0$,
if $(v_q,B_q,p_q,R_q)$ is a smooth solution of
\eqref{iter:relaxed-system} satisfying
\eqref{prep:old-ordinary}--\eqref{full:stress-support} with
$\alpha=(b-1)\gamma_S/b$, then there exists a smooth solution
$(v_{q+1},B_{q+1},p_{q+1},R_{q+1})$ of the same system satisfying
these assumptions with $q+1$ in place of $q$.
Moreover,
\begin{equation}\label{full:increment-bound}
 \begin{aligned}
 \|\partial_t^h(v_{q+1}-v_q)\|_r
 &\le C_0^{h+1}\lambda_{q+1}^{r+h}\delta_{q+1}^{1/2},\\
 \|\partial_t^h(B_{q+1}-B_q)\|_r
 &\le C_0^{h+1}\lambda_{q+1}^{r+h}\delta_{B,q+1}^{1/2},
 \end{aligned}
\end{equation}
for $h\le j_1+1$ and $r+h\le\rgood+4$, and both field increments
are supported in $(1/2,5/2)$.
\end{proposition}
The derivative and pointwise bounds for the fields and their increments
are verified in
Subsection~\ref{ssec:fields-pressure-support}, using
Proposition~\ref{full:endpoint-derivative-closure}.
Proposition~\ref{stress:complete-bounds} gives the stress bounds.
The pressure bounds are verified in Subsection~\ref{ssec:iteration-pressure},
and the support conditions and exact equations in
Subsection~\ref{ssec:iteration-support}.

\subsection{Parameters and notation}\label{ssec:fixed-scales}

\subsubsection{Scales}
At stage $q$, the construction uses a mollification length $\ell$,
a preparation time $\tau_c$, a fast time $\tau_a$ and a spatial chart
radius $\lambda_\parallel^{-1}$. With the exponents chosen in
Proposition~\ref{full:step}, we set
\begin{equation}\label{setup:scales}
 \begin{gathered}
 \ell^{-1}=\lambda_q^{1+(b-1)\gamma_\ell},\qquad
 \alpha=\frac{b-1}{b}\gamma_S,\\
 \tau_c=\varepsilon_{q+1}^{\gamma_\ell}\frac{1}{\lambda_q\delta_q^{1/2}},\qquad
 \tau_a=\varepsilon_{q+1}^{\gamma_a}\tau_c,\\
 \varepsilon_\tau=\frac{\tau_a}{\tau_c},\qquad
 \lambda_\parallel=\varepsilon_{q+1}^{\gamma_\parallel}\tau_c^{-1}.
 \end{gathered}
\end{equation}
The exponent $\gamma_\ell$ measures both the mollification gain and the
shortening of the eddy turnover time, while $\gamma_S$ measures the
Schauder loss. We choose $3\gamma_S<\gamma_\ell$.
The stage index is suppressed in $\ell,\tau_c,\tau_a,\lambda_\parallel$.
Slow and fast material derivatives have costs $\tau_c^{-1}$ and
$\tau_a^{-1}$, respectively, and magnetic derivatives have cost
$\lambda_\parallel$. The charts have spatial radius
$\lambda_\parallel^{-1}$ in every direction and time intervals of
length $O(\tau_c)$; the transverse oscillation frequency is
$\lambda_{q+1}$.

The parameter choices ensure
\[
 \lambda_q<\ell^{-1}<\lambda_{q+1},\qquad
 \tau_c\lambda_q\delta_q^{1/2}\ll1,\qquad
 \lambda_\parallel^{-1}\lambda_q\delta_{B,q}^{1/2}\ll1,
 \qquad \tau_a<\tau_c.
\]
These inequalities separate the spatial frequencies, control the
chart distortions under the velocity and magnetic flows, and provide the gain in
the supported approximate antiderivative. In particular, the magnetic differentiation cost
of the iterate is related to the corresponding cost in the charts by
\begin{equation}\label{setup:magnetic-rate}
 \lambda_q\delta_{B,q}^{1/2}
 =\lambda_\parallel
   \varepsilon_{q+1}^{\gamma_\ell+(\gamma_a-(b-1)\gamma_\parallel)/b}.
\end{equation}
The frequency ratios satisfy
\[
 \varepsilon_{q+1}=\varepsilon_q^b
 =\lambda_q^{1-b}\qquad(q\ge0),
\]
so the exponents of $\lambda_q$ in $\delta_q^{1/2}$ and
$\delta_{B,q}^{1/2}$ differ by
$(b-1)(\gamma_a+\gamma_\parallel)/b$. The stress amplitude at stage
$q$ is $\lambda_q^{-\alpha}\delta_{q+1}$.
The definitions also give
\begin{equation}\label{setup:scale-identities}
 \begin{gathered}
 \ell\lambda_q=\varepsilon_{q+1}^{\gamma_\ell},\qquad
 (\ell\lambda_{q+1})^{-1}=\varepsilon_{q+1}^{1-\gamma_\ell},\\
 \lambda_{q+1}^{\alpha}=\varepsilon_{q+1}^{-\gamma_S},\qquad
 \varepsilon_\tau=\varepsilon_{q+1}^{\gamma_a},\\
 \tau_c\lambda_\parallel=\varepsilon_{q+1}^{\gamma_\parallel},\qquad
 \Bigl(\frac{\delta_{B,q}}{\delta_q}\Bigr)^{1/2}
       =\varepsilon_q^{\gamma_a+\gamma_\parallel},\\
 \Bigl(\frac{\delta_{B,q+1}}{\delta_{q+1}}\Bigr)^{1/2}
       =\varepsilon_{q+1}^{\gamma_a+\gamma_\parallel}
       =\tau_a\lambda_\parallel.
 \end{gathered}
\end{equation}
Thus each fixed Schauder factor is bounded by
$\varepsilon_{q+1}^{-\gamma_S}$, while the spatial separation is
measured by $\gamma_\ell$. The comparisons of derivative costs used throughout the
construction are
\begin{equation}\label{setup:rate-comparisons}
 \begin{gathered}
 \tau_a^{-1}\le\lambda_q,\qquad
 \frac{\tau_a^{-1}}{\lambda_{q+1}\delta_{q+1}^{1/2}}
 =\frac{\lambda_\parallel}{\lambda_{q+1}\delta_{B,q+1}^{1/2}}
 =\varepsilon_{q+1}^{1-\beta-\gamma_a-\gamma_\ell}
 \le\varepsilon_{q+1}^{\gamma_a},\\
 \varepsilon_{q+1}\le\varepsilon_{q+1}^{\gamma_a},\qquad
 \ell^{-\alpha}\le\lambda_{q+1}^{\alpha},\qquad
 \tau_a\ell^{-1}\delta_{q+1}^{1/2}=o(1),\qquad
 \tau_a^2\lambda_q^2\ell^{-\alpha}\delta_{q+1}=o(1).
 \end{gathered}
\end{equation}
Lemma~\ref{iter:parameter-lemma} proves these comparisons. The last
two quantities are the normalized coefficient size of the principal
LDF acting on scalars and the gradient size of the corrector,
respectively.

\subsubsection{Derivative ranges}
We will use the sharp, lossy and ordinary classes defined below.
In the inductive hypotheses, the material and magnetic derivative
orders satisfy $k+m\le j_1$ and $r+k+m\le\rgood$. The ordinary field
bounds hold for $h\le j_1+1$ and $1\le r+h\le\rgood+4$, whereas the
ordinary stress bounds hold for $h\le j_1$ and $r+h\le\rgood$. Interpolation
between consecutive spatial orders gives sharp bounds for the old
iterate through order $\rgood-1$, and lossy bounds through order
$\rgood+3$ for the fields and $\rgood-1$ for the stress. We do not
assume estimates beyond these ranges.

After the preparation in Section~\ref{sec:preparation}, the sharp
estimates hold through order
\begin{equation}\label{setup:prepared-ceiling}
 \rprep=\rgood-\dstar,\qquad \dstar=2m_0+8.
\end{equation}
The estimates in
Sections~\ref{sec:charts-and-amplitudes}--\ref{sec:endpoint-stress}
require further spatial derivatives; we specify the number at each
application. In particular, the sharp estimates for every stress term
hold through order $\rcut$, defined in
Section~\ref{ssec:derivative-reserve}. The constructed quantities also
satisfy lossy and ordinary bounds. The lossy bounds hold at every
spatial order, within the stated material and magnetic derivative
ranges, except for quantities containing a difference between the old
iterate and a local background.
Under \eqref{setup:reserve-inequality}, Lemma~\ref{setup:calculus}(vii)
gives estimates at the remaining spatial orders allowed by these lossy
bounds, with the derivative costs of the next iterate.
Section~\ref{ssec:derivative-reserve} specifies the additional
derivatives needed for this application. We first estimate each stress
with the local derivative costs, then apply this lemma to obtain the
inductive bounds. At the highest material and magnetic derivative
orders, Theorem~\ref{galbrun:linear-theorem}(f) estimates the Galbrun
solution with the sharp amplitudes by using the coarse estimates.

\paragraph{Field and pressure bounds.}
The inductive field bounds give $[B_q]_1\lesssim\lambda_q\delta_{B,q}^{1/2}$. We need
this bound on the full gradient: a bound on $D_{B_q}B_q$ alone would
control neither the magnetic gradient nor the distortion of the charts.
We assume four additional orders of spatial and ordinary time
derivatives for the following reason. The new stress contains the
differences $v_q-v_{\ell,n}$ and $B_q-B_{\ell,n}$, differentiated
at most twice beyond the order of the stress estimate. Interpolation
for the H\"older estimate requires one further spatial derivative.
Finally, changing to the new material and magnetic derivatives requires
stress estimates one order beyond the range of the inductive material
and magnetic bounds. We prove the increment estimates throughout the
resulting range of spatial and ordinary time derivatives.
The induction identity $D_{t,q}B_q=D_{B_q}v_q$ is consistent with
these scales, since
\[
 \delta_{B,q}^{1/2}(\lambda_q\delta_q^{1/2})
 =\delta_q^{1/2}(\lambda_q\delta_{B,q}^{1/2}).
\]
By Lemma~\ref{high:transport-gradient}, commuting a spatial gradient
with $D_{t,q}^kD_{B_q}^m$ introduces gradients with material and
magnetic orders summing to less than $k+m$. Induction on $k+m$
therefore gives the gradient estimates without increasing the required
derivative order.

We also propagate the material and magnetic derivative bounds for the
pressure. When estimating Lie derivatives of the stress, we include
the terms involving the background gradients. In particular, for the
contravariant identity tensor,
$\mathcal L_{B_q}\IId=-(\nabla B_q+\nabla B_q^{\mathsf T})$;
its magnetic Lie derivative need not vanish.

\subsubsection{Notation for differential operators}\label{ssec:path-conventions}
For a background $(v,B)$, the material and magnetic derivatives are
\begin{equation}
 D_t=\partial_t+v\cn,\qquad D_B=B\cn.
\end{equation}
They act componentwise on vectors and tensors.
For a vector $F$, a one-form $\Theta$ and a contravariant two-tensor $T$,
\begin{equation}
 \mathcal L_zF=z\cn F-(\DD z)F,\qquad
 \mathcal L_z\Theta=z\cn\Theta+(\DD z)^T\Theta,\qquad
 \mathcal L_zT=z\cn T-(\DD z)T-T(\DD z)^T,
\end{equation}
and $\mathcal L_z=z\cn$ on scalars. The material and magnetic Lie derivatives are
$\mathcal D_t=\partial_t+\mathcal L_v$ and $\mathcal L_B$. The induction
equation gives
\[
 [D_t,D_B]=0,\qquad [\mathcal D_t,\mathcal L_B]=0.
\]
On a path $(v_s,B_s)$, we write $D_{t,s}$ and $D_{B_s}$.
We use $k$ for the material derivative order, $m$ for the magnetic
derivative order, and $h$ for the ordinary time derivative order.
The Els\"asser combinations
$\mathcal A^\pm=D_t\pm D_B$ and $\mathscr A^\pm=\mathcal D_t\pm\mathcal L_B$
occur only in exact algebraic identities. The norm $\|f\|_r$ is the
spatial $C^r$ norm, uniform in physical time, and $\|f\|_{r+\alpha}$ the
corresponding H\"older norm; see Appendix~\ref{ssec:holder-calculus}.
We use italic $d$ for exterior differentiation and upright $\dd$ for
integration measures and ordinary derivatives. A one-form is identified
with its Euclidean coefficients when writing its curl, with the
convention
\[
 d\Theta=\iota_{\curl\Theta}\mathrm{vol},\qquad
 \mathrm{vol}=d x_1\wedge d x_2\wedge d x_3.
\]

\subsubsection{Derivative classes}
We group the derivative estimates by amplitude, spatial frequency, and
material and magnetic derivative costs. We include integer and
H\"older estimates separately, as required for the product, transport
and singular integral estimates below. The three pairs of costs are
\begin{equation}\label{setup:rate-pairs}
 \begin{aligned}
 \mathrm c&=(\tau_c^{-1},\lambda_\parallel),\\
 \mathrm a&=(\tau_a^{-1},\lambda_\parallel),\\
 \mathrm f&=(\lambda_{q+1}\delta_{q+1}^{1/2},
              \lambda_{q+1}\delta_{B,q+1}^{1/2}).
 \end{aligned}
\end{equation}
The first two pairs give the local derivative costs, with material
derivative costs $\tau_c^{-1}$ and $\tau_a^{-1}$, respectively. The last pair gives the derivative costs of the next
iterate, which dominate the local derivative costs by
\eqref{setup:rate-comparisons}.

\begin{definition}[Derivative classes]
Fix a background and its transports $D_t,D_B$, a spatial frequency
$\Lambda>0$, a pair of derivative costs $\rho=(\mathrm a_t,\mathrm a_B)$ with positive
components, a maximum derivative order
$N\in\{0,1,\ldots\}\cup\{\infty\}$, a bound $0\le J\le j_1+1$
on the sum of the material and magnetic derivative orders,
and amplitudes $E,E'>0$. The components $\mathrm a_t,\mathrm a_B$
always refer to the chosen pair: for $\rho=\mathrm a$, for example,
\[
 \mathrm a=(\mathrm a_t,\mathrm a_B),\qquad
 \mathrm a_t=\tau_a^{-1},\qquad \mathrm a_B=\lambda_\parallel.
\]
The transport derivatives act componentwise unless Lie derivatives
are displayed, and belong to the background named in the statement.
The norms include the supremum over the time interval and the chart
neighborhood in question. All classes are taken componentwise on tensors.

The \emph{sharp class} is defined by
\begin{equation}\label{setup:class}
 F\in\mathcal C_{N,J}(E;\Lambda,\rho)
 \quad\Longleftrightarrow\quad
 \begin{cases}
 \|D_t^kD_B^mF\|_r\le E\Lambda^r\mathrm a_t^k\mathrm a_B^m,\\
 \|D_t^kD_B^mF\|_{r+\alpha}
 \le\max\{\ell^{-1},\Lambda\}^{\alpha}
       E\Lambda^r\mathrm a_t^k\mathrm a_B^m,
 \end{cases}
\end{equation}
where both inequalities hold for all nonnegative integers $r,k,m$
with $k+m\le J$ and $r+k+m\le N$.
We write $\mathcal C_N=\mathcal C_{N,j_1}$.

The \emph{lossy class} assigns the cost $\ell^{-1}$ to every
derivative:
\begin{equation}
 F\in\mathcal K_{N,J}(E')
 \quad\Longleftrightarrow\quad
 \begin{cases}
 \|D_t^kD_B^mF\|_r\le E'\ell^{-(r+k+m)},\\
 \|D_t^kD_B^mF\|_{r+\alpha}
       \le E'\ell^{-(r+k+m)-\alpha},
 \end{cases}
\end{equation}
for the same range $k+m\le J$, $r+k+m\le N$.
We write $\mathcal K_J=\mathcal K_{\infty,J}$ and
$\mathcal K=\mathcal K_{j_1}$.

For a nonnegative integer $H$, the \emph{ordinary class} records
spatial derivatives and ordinary time derivatives through order $H$:
\begin{equation}\label{setup:ordinary-class}
 F\in\mathcal O_{N,H}(E;\Lambda)
 \quad\Longleftrightarrow\quad
 \begin{cases}
 \|\partial_t^hF\|_r\le E\Lambda^{r+h},\\
 \|\partial_t^hF\|_{r+\alpha}
       \le\max\{\ell^{-1},\Lambda\}^{\alpha}E\Lambda^{r+h},
 \end{cases}
\end{equation}
for all nonnegative integers $r,h$ with $h\le H$ and $r+h\le N$.
Thus the second index of $\mathcal O_{N,H}$ bounds the number of
ordinary time derivatives, whereas the second index of
$\mathcal C_{N,J}$ and $\mathcal K_{N,J}$ bounds the sum of the
material and magnetic derivative orders.

If a quantity belongs to both $\mathcal C_{N,J}(E;\Lambda,\rho)$ and
$\mathcal K_{N',J}(E')$, we call $E'/E$ the \emph{loss} of the pair
of estimates.

For $N,J\ge1$, a background is called
$(\Lambda,\rho)$-\emph{adapted} through $(N,J)$ if
\[
 \nabla v\in\mathcal C_{N-1,J-1}(C\mathrm a_t;\Lambda,\rho),\qquad
 \nabla B\in\mathcal C_{N-1,J-1}(C\mathrm a_B;\Lambda,\rho).
\]
\end{definition}

It follows directly from the definition that
\[
 \mathcal C_{N,J}(E;\Lambda,\rho)
 \subset\mathcal C_{N,J}(E;\Lambda',\rho')
\]
whenever $\Lambda\le\Lambda'$ and each component of $\rho$ is at
most the corresponding component of $\rho'$. To apply the calculus
for the lossy classes, we also require
\[
 \nabla v,\nabla B\in\mathcal K_{N-1,J-1}(C\ell^{-1}),
\]
which every background of the construction satisfies through the
maximum derivative order of its lossy bounds. The class calculus is stated and proved
in Lemma~\ref{setup:calculus} of
Section~\ref{ssec:class-properties}.

The H\"older factor in \eqref{setup:class} is $\ell^{-\alpha}$
for a slow quantity with $\Lambda\le\ell^{-1}$, and
\begin{equation}\label{stress:holder-unit}
 \lambda_{q+1}^{\alpha}=\varepsilon_{q+1}^{-\gamma_S}
\end{equation}
for a fast quantity with $\Lambda=\lambda_{q+1}$. The slow quantities
are smoothed at scale $\ell$ or solve transport equations with such
data. Integer estimates through a given derivative order yield the H\"older
estimates one order lower, by interpolation:
\[
 \|F\|_{r+\alpha}
 \le\|F\|_r^{1-\alpha}\|F\|_{r+1}^{\alpha},\qquad
 \Lambda^\alpha\le\max\{\ell^{-1},\Lambda\}^{\alpha}.
\]
The product estimate \eqref{setup:product} uses only one factor in
the H\"older norm, and hence introduces only one H\"older factor.
An integer bound deduced from a H\"older bound retains that factor
in its amplitude. This occurs in
Appendix~\ref{ssec:component-transport-estimates} and in the Galbrun
estimates. We include the factor in the sharp amplitude; we omit it
from a lossy amplitude only when a separate lossy estimate is available.

Products multiply the losses of their factors, up to fixed constants.
Estimating a transport derivative with the lossy factor $\ell^{-1}$ in
both the sharp and lossy estimates leaves their ratio unchanged.
The operations in Lemma~\ref{setup:calculus}(iii)--(vi) preserve
losses up to fixed constants and fixed powers of
$\lambda_{q+1}^{\alpha}$. In particular, commuting a transport
derivative through an order-zero singular integral introduces a
further H\"older factor; see Lemma~\ref{setup:calculus}(v) and
Appendix~\ref{ssec:holder-calculus}. Each such factor is at most
$\lambda_{q+1}^{\alpha}$, and
\begin{equation}\label{setup:holder-constant}
 c=64+14j_0+8j_1
\end{equation}
bounds their total number in every stress term, as recorded in
Section~\ref{ssec:derivative-closure}. The choice of parameters absorbs
these fixed losses.

\smallskip
\noindent\emph{Indices and profiles.}
The index $n$ specifies a slow time interval and $J$ a spatial chart;
$I=(n,J,\zeta)$ is the combined index for a local contribution. Temporal and spatial profiles are
written $\alpha_I(t/\tau_a)$ and $\varphi(\lambda_{q+1}y_I^{k_I})$;
the variable $\tau$ is used only as the argument of fixed temporal
profiles and their primitives, and $\alpha$ without an index is the
H\"older-loss exponent. The superscript $[j]$ denotes a normalized
one-variable primitive.

\section{Local Backgrounds}\label{sec:preparation}

We first mollify the old iterate in space and then solve
\eqref{prep:local-equation}, with force $\ddiv R_\ell$ and the
mollified fields as initial data. This gives smooth
local backgrounds satisfying the induction equation. The projected
transport estimates in Appendix~\ref{ssec:component-transport-estimates}
bound their differences from the mollified fields in spatial and
ordinary time derivatives. We then apply the comparison lemmas of
Appendix~\ref{ssec:finite-derivative-conversion} to estimate material,
magnetic and Lie derivatives of these differences. In the passage to
Lie derivatives, the background gradients are differentiated at most
$k+m-1$ times in the material and magnetic directions. Thus we can
estimate one more such derivative of the differences without assuming
it for the old iterate.

The approximation order $m_0$ satisfies \eqref{iter:ordinary-accuracy}, so
that
\begin{equation}\label{prep:accuracy-inequality}
 \varepsilon_\ell=(\ell\lambda_q)^{m_0}
 \le\ell^{2\alpha}\varepsilon_{q+1}^2
        \delta_{B,q}^{(j_1+3)/2}.
\end{equation}
The inequalities $\tau_c^{-1}\ge\lambda_q\delta_{B,q}^{1/2}$ and
$\lambda_\parallel\ge\lambda_q\delta_{B,q}^{1/2}$, together with
\eqref{prep:accuracy-inequality}, will allow us to replace the ordinary
derivative costs in the comparison estimates by the local material and
magnetic costs, including at the smaller magnetic amplitude. We use
only the derivatives provided by
\eqref{prep:old-ordinary}--\eqref{prep:old-stress-time}. The estimates
hold in the classes of Section~\ref{ssec:fixed-scales} through order
$\rprep=\rgood-\dstar$; the choice \eqref{setup:prepared-ceiling}
leaves the derivatives needed for the comparisons below. We also use
the scale inequalities \eqref{setup:rate-comparisons} throughout.

\subsection{Spatial smoothing}

Choose a real compactly supported smooth convolution kernel $\rho$
with integral one and vanishing moments through degree $2m_0+2$; the
construction in Appendix~\ref{ssec:spatial-mollification} gives such a
kernel, which need not be nonnegative. Set
\[
 \mathcal J_\ell f=\rho_\ell*f,\qquad
 (v_\ell,B_\ell,p_\ell,R_\ell)
   =\mathcal J_\ell(v_q,B_q,p_q,R_q).
\]
The convolution commutes with ordinary derivatives and with divergence.
Define
\begin{align*}
 R_\ell^c&=v_\ell\otimes v_\ell-B_\ell\otimes B_\ell
       -\mathcal J_\ell(v_q\otimes v_q-B_q\otimes B_q),\\
 M_\ell&=B_\ell\times v_\ell-\mathcal J_\ell(B_q\times v_q),
\end{align*}
so that
\begin{align}
 \partial_tv_\ell+\ddiv(v_\ell\otimes v_\ell-B_\ell\otimes B_\ell)
       +\nabla p_\ell&=\ddiv(R_\ell+R_\ell^c),\notag\\
 \partial_tB_\ell+\curl(B_\ell\times v_\ell)&=\curl M_\ell.
 \label{prep:mollified-equations}
\end{align}
Apply Proposition~\ref{moll:spatial-basic} with $d=2m_0$ to the
ordinary inductive bounds, for $k\le j_1+1$ and
$r+k+2m_0\le\rgood$. In each quadratic commutator term, both factors
are differentiated and therefore have amplitude at most
$\delta_q^{1/2}$. After distributing the time derivatives by Leibniz'
rule, each factor is differentiated at most $r+k+2m_0$ times.
Thus the approximation estimate and the estimates for $R_\ell^c$
and $M_\ell$ all retain the factor $\varepsilon_\ell^2$.

The convolution is bounded on the integer norms. Its H\"older estimate
in Proposition~\ref{moll:spatial-basic},
\[
 \|\mathcal J_\ell f\|_{r+\alpha}\le C\ell^{-\alpha}\|f\|_r,
\]
gives the factor $\ell^{-\alpha}$ whenever the integer norm on the
right-hand side is controlled, without requiring another derivative.
Consequently,
\eqref{prep:old-stress-time} and rule (vi) of
Lemma~\ref{setup:calculus} show that the mollified stress satisfies
\begin{equation}\label{prep:mollified-stress-classes}
 R_\ell\in\mathcal O_{\rgood,j_1}(C\lambda_q^{-\alpha}\delta_{q+1};\lambda_q)
 \cap\mathcal O_{\infty,j_1}(C\lambda_q^{-\alpha}\delta_{q+1};\ell^{-1}).
\end{equation}
The mollified fields satisfy
\[
 \begin{aligned}
 \nabla v_\ell&\in\mathcal O_{\rgood+2,j_1+1}(C\lambda_q\delta_q^{1/2};\lambda_q)
 \cap\mathcal O_{\infty,j_1+1}(C\lambda_q\delta_q^{1/2};\ell^{-1}),\\
 \nabla B_\ell&\in\mathcal O_{\rgood+2,j_1+1}(C\lambda_q\delta_{B,q}^{1/2};\lambda_q)
 \cap\mathcal O_{\infty,j_1+1}(C\lambda_q\delta_{B,q}^{1/2};\ell^{-1}).
 \end{aligned}
\]
At spatial orders beyond the sharp range, the derivatives fall on the
convolution kernel and give the stated lossy bounds.

\subsection{Local correction}

We cover time by intervals on the slow time scale. Set $t_n=n\tau_c$ and
choose a fixed nonnegative bump $\eta^0\in C_c^\infty((-3/4,3/4))$ that
is positive on $[-1/2,1/2]$. Its square-normalized translates are
\begin{equation}
 \eta_n(t)=\frac{\eta^0(t/\tau_c-n)}
 {\bigl(\sum_{j\in\mathbb Z}\eta^0(t/\tau_c-j)^2\bigr)^{1/2}},
 \qquad \sum_n\eta_n^2=1,\qquad
 |\partial_t^k\eta_n|\le C_k\tau_c^{-k}.
\end{equation}
At every time, the argument of at least one translate belongs to $[-1/2,1/2]$;
the denominator is therefore bounded below by a fixed positive constant.
At most two cutoffs are nonzero, and the supports of distinct cutoffs
whose indices have the same parity are disjoint. Choose the common outer cutoff
$\varrho\in C_c^\infty(\R)$, $0\le\varrho\le1$, with
\begin{equation}\label{prep:outer-cutoff}
 \begin{gathered}
 \varrho=1\quad\hbox{on }
 \left[\frac12+\frac{\delta_q^{1/2}}2,
             \frac52-\frac{5\delta_q^{1/2}}2\right],\qquad
 \supp\varrho\Subset
 \left(\frac12+\frac{\delta_q^{1/2}}4,
             \frac52-\frac{\delta_q^{1/2}}4\right),
 \qquad |\varrho^{(k)}|\le C_k\delta_q^{-k/2}.
 \end{gathered}
\end{equation}
A product of two fixed smooth transition functions, rescaled by
$\delta_q^{1/2}$, gives such a cutoff. In particular, $\varrho=1$ on the time
support of $R_q$ and of its spatial convolution. Only the
finitely many intervals with
$\operatorname{dist}(t_n,\supp\varrho)\le2\tau_c$ enter the
construction, and all sums over time intervals below run over these
indices. Their quadratic partition of unity equals one on the fixed neighborhood
\[
 \{t:\operatorname{dist}(t,\supp\varrho)<\tau_c\}.
\]
Indeed, every nonzero translate on this neighborhood has its center at
distance less than $7\tau_c/4$ from $\supp\varrho$ and hence belongs to
the selected family. The temporal smoothing operations have length
$C\tau_a$ and remain within this neighborhood for sufficiently large $a$.
The quadratic partition of unity and the stress decomposition for this finite family hold
there, and hold globally after multiplication by $\varrho^2$.
The cut off regions for the correction lie at distance at most $4\tau_c$ from
$\supp\varrho$. Since $\tau_c/\delta_q^{1/2}\to0$ and
$\delta_{q+1}^{1/2}/\delta_q^{1/2}\to0$, we may require
\[
 4\tau_c+\tfrac52\delta_{q+1}^{1/2}<\tfrac14\delta_q^{1/2},
\]
which together with \eqref{prep:outer-cutoff} places every outgoing
cut off region strictly inside the next stress-support interval
\eqref{full:stress-support} with $q+1$ in place of $q$.

For the source, the correction cutoff and the local background we
use the nested intervals
\begin{equation}\label{prep:slow-windows}
 I_n^{\rm src}=(t_n-\tau_c,t_n+\tau_c),\qquad
 I_n^{\rm col}=(t_n-2\tau_c,t_n+2\tau_c),\qquad
 I_n^{\rm bg}=(t_n-4\tau_c,t_n+4\tau_c).
\end{equation}
The support of $\eta_n$ lies at distance at least $\tau_c/4$ from the boundary of
$I_n^{\rm src}$; a smooth time cutoff equal to one on $I_n^{\rm src}$
and supported in $I_n^{\rm col}$ has derivatives bounded by
$C_k\tau_c^{-k}$; all local equations and charts are available on
$I_n^{\rm bg}$, which extends before and after $I_n^{\rm col}$.
Figure~\ref{fig:window} shows these intervals together with the
correction cutoff and the forcing and outgoing cutoff-error supports.

For each reference time $t_n$, consider the Cauchy problem
\begin{equation}\label{prep:local-equation}
 \begin{cases}
 \partial_t v_{\ell,n}+\ddiv(v_{\ell,n}\otimes v_{\ell,n}-B_{\ell,n}\otimes B_{\ell,n})+\nabla p_{\ell,n}
      =\ddiv R_\ell,\\
 \partial_tB_{\ell,n}+[v_{\ell,n},B_{\ell,n}]=0,\qquad \ddiv v_{\ell,n}=\ddiv B_{\ell,n}=0,\\
 (v_{\ell,n},B_{\ell,n})(t_n)=(v_\ell,B_\ell)(t_n),\qquad \displaystyle\int p_{\ell,n}=0.
 \end{cases}
\end{equation}
We require the solution on $I_n^{\rm bg}$; any fixed enlargement of
this interval affects only the constants. Write
$D_{t,n}=\partial_t+v_{\ell,n}\cn$, $D_{B,n}=B_{\ell,n}\cn$ and
$\mathcal D_{t,n}=\partial_t+\mathcal L_{v_{\ell,n}}$, and define the comparison fields
\[
 \Delta_n^v=v_\ell-v_{\ell,n},\qquad
 \Delta_n^B=B_\ell-B_{\ell,n},\qquad
 G_n^v=v_q-v_{\ell,n},\qquad
 G_n^B=B_q-B_{\ell,n}.
\]
With $z_{\ell,n}^\pm=v_{\ell,n}\pm B_{\ell,n}$, their Els\"asser
combinations are
\[
 \Delta_n^\pm=\Delta_n^v\pm\Delta_n^B=z_\ell^\pm-z_{\ell,n}^\pm,\qquad
 G_n^\pm=G_n^v\pm G_n^B=z_q^\pm-z_{\ell,n}^\pm.
\]

\begin{lemma}[Local evolution and comparison estimates]
\label{prep:ordinary-comparison}
The solution of \eqref{prep:local-equation} exists on $I_n^{\rm bg}$.
Its fields are uniformly bounded, and the magnitude of its magnetic
field is bounded below by a positive constant. For $k\le j_1+1$ and
$r+k\le\rgood-2m_0-3$,
\begin{align}
 \|\partial_t^k \Delta_n^\pm\|_{r+\alpha}
 &\le C\ell^{-\alpha}\varepsilon_\ell^2
       \lambda_q^{r+k}\varepsilon_{q+1}^{\gamma_\ell}\delta_q^{1/2},
 \label{prep:correction-comparison}\\
 \|\partial_t^kG_n^\pm\|_{r+\alpha}
 &\le C\ell^{-\alpha}\varepsilon_\ell^2
       \lambda_q^{r+k}\delta_q^{1/2}.
 \label{prep:complete-comparison}
\end{align}
For $k\le j_1+1$ and every $r+k\le\rgood+3$,
\[
 \|\partial_t^kG_n^\pm\|_{r+\alpha}\le C\ell^{-\alpha}\ell^{-(r+k)}\delta_q^{1/2}.
\]
The evolved fields satisfy
\begin{equation}\label{prep:whole-patch-fields}
 \|\partial_t^k z_{\ell,n}^\pm\|_{r+\alpha}
 \le C\ell^{-\alpha}\lambda_q^{r+k}\delta_q^{1/2},
 \quad k\le j_1+1,\quad 1\le r+k\le\rgood-3,
\end{equation}
and the lossy bound with $\lambda_q\ell^{-(r+k-1)}$ in place of
$\lambda_q^{r+k}$ for $k\le j_1+1$ and every $r+k\ge1$. In the classes of
Section~\ref{ssec:fixed-scales},
\begin{equation}\label{prep:local-field-classes}
 \nabla z_{\ell,n}^\pm\in
 \mathcal O_{\rgood-4,j_1+1}(C\ell^{-\alpha}\lambda_q\delta_q^{1/2};\lambda_q)
 \cap\mathcal O_{\infty,j_1+1}(C\ell^{-\alpha}\lambda_q\delta_q^{1/2};\ell^{-1}).
\end{equation}
\end{lemma}
The lossy comparison follows from the triangle inequality and
\eqref{prep:old-ordinary}, without the approximation factor
$\varepsilon_\ell^2$. For the evolved fields, the first derivative has
cost $\lambda_q$ and subsequent derivatives have cost $\ell^{-1}$.
We prove the lemma at the end of
Section~\ref{ssec:component-transport-estimates}, where the projected
transport estimates are available. Applied to
\eqref{prep:mollified-equations}, they control the local evolution error;
adding the spatial smoothing error gives the comparison with the old iterate.

\begin{corollary}[Magnetic comparison]
\label{prep:magnetic-ordinary-comparison}
For nonnegative integers $r,k$ with
$k\le j_1+1$ and $r+k\le\rgood-2m_0-3$,
\begin{equation}
 \|\partial_t^kG_n^B\|_{r+\alpha}
 \le C\ell^{-\alpha}\varepsilon_\ell
       \lambda_q^{r+k}\delta_{B,q}^{1/2},\qquad
 \|\partial_t^k\Delta_n^B\|_{r+\alpha}
 \le C\ell^{-\alpha}\varepsilon_\ell
       \lambda_q^{r+k}\varepsilon_{q+1}^{\gamma_\ell}
       \delta_{B,q}^{1/2}.
\end{equation}
For $k\le j_1+1$ and $r+k\le\rgood+3$, the lossy bound is
\[
 \|\partial_t^kG_n^B\|_{r+\alpha}\le C\ell^{-\alpha}\ell^{-(r+k)}\delta_q^{1/2},
\]
whose loss is $\varepsilon_\ell^{-1}(\delta_q/\delta_{B,q})^{1/2}$.
\end{corollary}
\begin{proof}
Apply Lemma~\ref{prep:ordinary-comparison} to
$G_n^B=(G_n^+-G_n^-)/2$ and
$\Delta_n^B=(\Delta_n^+-\Delta_n^-)/2$. The two sharp estimates
follow from \eqref{prep:accuracy-inequality}, since
\[
 \varepsilon_\ell^2\delta_q^{1/2}\le\varepsilon_\ell\delta_{B,q}^{1/2}.
\]
The factor $\varepsilon_{q+1}^{\gamma_\ell}$ in the estimate for
$\Delta_n^B$ is unchanged. The same half-difference and the lossy
estimate of that lemma give the stated lossy bound.

The comparison used here is for the coupled velocity and magnetic
equations. Indeed,
\[
 D_{t,n}\Delta_n^B-D_{B,n}\Delta_n^v
 =\curl M_\ell+\Delta_n^B\cn v_\ell-\Delta_n^v\cn B_\ell,
 \qquad \Delta_n^B(t_n)=0.
\]
Thus the estimate for $\Delta_n^B$ uses the velocity comparison to
control $D_{B,n}\Delta_n^v$. Lemma~\ref{prep:ordinary-comparison}
provides the factor $\varepsilon_\ell^2$ using spatial and ordinary
time derivatives of the old fields through order $r+k+2m_0+2$; the
preceding inequality then gives the magnetic amplitude.
\end{proof}

\subsection{Separate transports}

The local induction equation gives $[D_{t,n},D_{B,n}]=0$, and the
comparison identities give
\begin{equation}\label{prep:frame-comparison}
 D_{t,n}=D_{t,q}-G_n^v\cn,\qquad
 D_{B,n}=D_{B_q}-G_n^B\cn.
\end{equation}
We use \eqref{prep:frame-comparison} to deduce estimates for
$D_{t,n}$ and $D_{B,n}$ from \eqref{prep:old-mixed}. In particular,
the estimates for derivatives of $B_{\ell,n}$ retain the amplitude
$\delta_{B,q}^{1/2}$.

\begin{lemma}[Local material and magnetic derivative bounds]
\label{prep:finite-conversion}
With respect to the local transports $D_{t,n},D_{B,n}$,
\begin{equation}\label{prep:local-gradient}
 \nabla v_{\ell,n}\in\mathcal C_{\rprep}(C\lambda_q\delta_q^{1/2};\lambda_q,\mathrm c),
 \qquad
 \nabla B_{\ell,n}\in\mathcal C_{\rprep}(C\lambda_q\delta_{B,q}^{1/2};\lambda_q,\mathrm c),
\end{equation}
so that the local background is $(\lambda_q,\mathrm c)$-adapted through
$(\rprep+1,j_1+1)$ with the ordinary bounds
\eqref{prep:local-field-classes}. We also have
\begin{equation}\label{prep:local-acceleration}
 \nabla(\ddiv R_\ell-\nabla p_{\ell,n})
 \in\mathcal C_{\rprep}(C\lambda_q^2\delta_q;\lambda_q,\mathrm c),
\end{equation}
and the stress obeys
\begin{equation}\label{prep:stress-good}
 R_\ell\in\mathcal C_{\rprep}(C\lambda_q^{-\alpha}\delta_{q+1};\lambda_q,\mathrm c)
 \cap\mathcal K(C\lambda_q^{-\alpha}\delta_{q+1}),
\end{equation}
in the classes defined by transport derivatives and in those defined by Lie derivatives.
\end{lemma}
\begin{proof}
We apply the comparison lemmas of
Appendix~\ref{ssec:finite-derivative-conversion} to
\eqref{prep:frame-comparison}. We use both their integer and H\"older
estimates. In particular, a comparison estimated by
\eqref{prep:complete-comparison} retains the factor $\ell^{-\alpha}$
even in the integer norm.

\emph{1. Gradients and approximation estimates.}
By Lemma~\ref{high:transport-gradient} and \eqref{prep:old-mixed},
the old gradients satisfy the required estimates for
$r+k+m\le\rgood-1$. In the induction on $k+m$, the commutator terms
contain gradients with fewer than $k+m$ material and magnetic
derivatives in sum. The estimate for $\nabla B_q$ therefore retains
its magnetic amplitude. Interpolating consecutive spatial estimates
gives the H\"older bound with factor
$\lambda_q^\alpha\le\ell^{-\alpha}$.

We next apply Lemma~\ref{high:ordinary-transport-comparison} with
$\Lambda=\lambda_q$. Its coefficient hypotheses follow from
\eqref{prep:old-ordinary}, \eqref{prep:whole-patch-fields} and
$\ell^{-\alpha}\delta_q^{1/2}\le1$, which bounds the amplitudes of
the differentiated coefficients by one. For an expression containing
$s$ spatial or ordinary time derivatives of a field,
Lemma~\ref{prep:ordinary-comparison} gives the comparison error
\begin{equation}\label{prep:raw-conversion-error}
 C\ell^{-\alpha}\varepsilon_\ell^2
              \lambda_q^{r+k+m+s}\delta_q^{1/2};
\end{equation}
in both norms. Here the factor $\ell^{-\alpha}$ comes from
\eqref{prep:complete-comparison}. We use $s=1$ for a field gradient;
for a gap, comparison with zero gives $s=0$. The lemma requires at
most $k$ additional time derivatives of this expression. Including its
$s$ previous derivatives, the required ordinary derivatives of the
field have order at most $r+k+m+s$. For $k+m\le j_1+1$, we have
\begin{equation}\label{prep:conversion-margin}
 \frac{\lambda_q^{k+m}}{\tau_c^{-k}\lambda_\parallel^m}
       \le\delta_{B,q}^{-(j_1+1)/2},\qquad
 \ell^{-\alpha}\varepsilon_\ell^2\delta_{B,q}^{-(j_1+2)/2}
 \le\ell^\alpha\varepsilon_\ell\varepsilon_{q+1}^2\delta_{B,q}^{1/2}
 \le\varepsilon_{q+1}^2,
\end{equation}
where the second inequality follows from
\eqref{prep:accuracy-inequality}.
Since $0<\delta_{B,q}^{1/2}\le\delta_q^{1/2}\le1$, the additional factor
$\delta_{B,q}^{-1/2}$ bounds both the ratio
$(\delta_q/\delta_{B,q})^{1/2}$ needed for the magnetic amplitude and
the factor $\delta_q^{-1/2}$ needed for
\eqref{prep:local-acceleration}. Thus the error
\eqref{prep:raw-conversion-error} is bounded at both field amplitudes,
which proves \eqref{prep:local-gradient}. The largest sum of spatial and ordinary time derivative orders
used is $\rprep+2=\rgood-2m_0-6$, within the range of
Lemma~\ref{prep:ordinary-comparison}. The same argument, with the sum of
the material and magnetic derivative orders at most $j_1+1$ and the
sum including spatial derivatives at most $\rprep+1$,
proves adaptedness through $(\rprep+1,j_1+1)$.

\emph{2. Stress.}
Lemma~\ref{high:component-lie-conversion}, with the old gradient
bounds, first expresses \eqref{prep:old-stress} in transport derivatives.
We then apply Proposition~\ref{moll:differentiated-spatial} with costs
$\tau_c^{-1}$ and $\lambda_\parallel$ to obtain
\eqref{prep:stress-good} for the old transports. Its coefficient
hypotheses are the old gradient bounds proved above. The lossy assertion
follows from \eqref{prep:mollified-stress-classes} by expanding
$D_{t,n}=\partial_t+v_{\ell,n}\cn$ and using the lossy field bounds.

To replace the old transports by the local ones on $R_\ell$, we use
Lemma~\ref{high:ordinary-transport-comparison}. The ordinary stress
amplitude is $\lambda_q^{-\alpha}\delta_{q+1}$, so the comparison error
is bounded by
\[
 C\ell^{-\alpha}\varepsilon_\ell^2\delta_q^{1/2}
 \lambda_q^{r+k+m-\alpha}\delta_{q+1}.
\]
The inequality \eqref{prep:conversion-margin} absorbs this error, proving
\eqref{prep:stress-good} for the local transports. Finally,
Lemma~\ref{high:component-lie-conversion} applies to the local gradients
with costs $\tau_c^{-1},\lambda_\parallel$ and gives the Lie derivative
bounds. Both conversions require derivatives of the gradients for which
the sum of the material and magnetic derivative orders is at most $k+m-1$.

\emph{3. The momentum equation and derivative orders.}
Compare the first-order expressions
\[
 \ddiv R_\ell-\nabla p_{\ell,n}=D_{t,n}v_{\ell,n}-D_{B,n}B_{\ell,n},\qquad
 \ddiv R_q-\nabla p_q=D_{t,q}v_q-D_{B_q}B_q.
\]
The difference of these expressions contains a gap and one ordinary
derivative. After taking a gradient,
Lemma~\ref{high:ordinary-transport-comparison} therefore gives
\eqref{prep:raw-conversion-error} with $s=2$. The gradient
$\nabla(\ddiv R_q-\nabla p_q)$ is controlled by the bounds for
material and magnetic derivatives of the stress and pressure, with spatial
derivatives commuted by Lemma~\ref{high:transport-gradient}. Including
the derivative used for interpolation, this requires stress derivatives
through spatial order $r+k+m+3$ and compositions of material and magnetic
derivatives of the pressure with $k+m\le j_1$.
The inequality \eqref{prep:conversion-margin} absorbs the division by
$\delta_q^{1/2}$ and proves \eqref{prep:local-acceleration}.
The terms $D_{t,n}v_{\ell,n}$ and $D_{t,q}v_q$ require ordinary
time derivatives of the fields through $k+1\le j_1+1$; the stress is
differentiated only through $k\le j_1$.

The Taylor estimates for the difference between the two expressions
above require spatial and
ordinary time derivatives of the old fields through order $\rprep+2m_0+4$.
This is the largest such order used in the proof, and is at most
$\rgood$ by \eqref{setup:prepared-ceiling}. The estimates for material and magnetic derivatives
require at most $\rprep+3\le\rgood$ derivatives, of which at most
$j_1$ are material or magnetic, so all the preceding estimates
lie within the inductive ranges.
\end{proof}

\begin{corollary}[Comparison in material and magnetic derivatives]
\label{prep:mixed-comparison}
With respect to the local transports, each
$F\in\{G_n^v,G_n^\pm,\Delta_n^v,\Delta_n^\pm\}$ satisfies
\begin{equation}\label{prep:local-gap-good}
 F\in\mathcal C_{\rprep+1,j_1+1}(C\varepsilon_{q+1}^2\delta_q^{1/2};\lambda_q,\mathrm c)
 \cap\mathcal K_{\rgood+3,j_1+1}(C\delta_q^{1/2}),
\end{equation}
and the magnetic differences satisfy
\[
 G_n^B,\,\Delta_n^B\in
 \mathcal C_{\rprep+1,j_1+1}
 (C\varepsilon_{q+1}^2\delta_{B,q}^{1/2};\lambda_q,\mathrm c)
 \cap\mathcal K_{\rgood+3,j_1+1}(C\delta_q^{1/2}).
\]
Both assertions hold in the
vector classes defined by transport derivatives and in those defined by Lie derivatives. One material or magnetic Lie derivative of $F$ lies in
\[
 \mathcal C_{\rprep}(C\varepsilon_{q+1}^2\delta_q^{1/2}\tau_c^{-1};\lambda_q,\mathrm c)
 \quad\hbox{or}\quad
 \mathcal C_{\rprep}(C\varepsilon_{q+1}^2\delta_q^{1/2}\lambda_\parallel;\lambda_q,\mathrm c),
\]
according to its type, with the magnetic amplitude for
$F=G_n^B,\Delta_n^B$.
The same assertions hold for differences of local backgrounds on an
overlap, with the corresponding magnetic improvement, and differences
involving only local evolutions retain the additional factor
$\varepsilon_{q+1}^{\gamma_\ell}$.
\end{corollary}
In the sequel we use the factor $\varepsilon_{q+1}^2$ in place of
$\varepsilon_\ell^2$ from the ordinary comparison, or
$\varepsilon_\ell$ from
Corollary~\ref{prep:magnetic-ordinary-comparison}. Both replacements
follow from \eqref{prep:accuracy-inequality} and remove the dependence
on $m_0$ from the later estimates. The respective losses are
$\varepsilon_{q+1}^{-2}$ and $\varepsilon_{q+1}^{-2}(\delta_q/\delta_{B,q})^{1/2}$.
Their maximum derivative order $\rgood+3$ follows from the old ordinary
bounds \eqref{prep:old-ordinary}, allowing one additional spatial
derivative for the H\"older estimate.

\begin{proof}
Each gap is a half-sum, half-difference or difference of the quantities
in Lemma~\ref{prep:ordinary-comparison}, and hence satisfies its
ordinary estimate with amplitude
$C\varepsilon_\ell^2\delta_q^{1/2}$. Apply
Lemma~\ref{high:ordinary-transport-comparison} to this gap and zero.
By \eqref{prep:raw-conversion-error} with $s=0$ and
\eqref{prep:conversion-margin}, we obtain the sharp class in
\eqref{prep:local-gap-good}, including the smaller amplitude for the
magnetic gaps. This application uses ordinary derivatives only through
order $\rprep+1$.

To obtain the lossy class, we expand the local transports in ordinary
derivatives and use Lemma~\ref{prep:ordinary-comparison} and
Corollary~\ref{prep:magnetic-ordinary-comparison}. A derivative with
$r+k+m\le\rgood+3$ uses spatial and ordinary time derivatives of the gap through order $r+k+m$, with time
order at most $k$. The coefficients in this expansion are derivatives
of the local fields, whose lossy bounds hold at every spatial order.

For Lie derivatives, apply
Lemma~\ref{high:component-lie-conversion} with $J=j_1+1$ and
$N=\rprep+1$. The required gradient bounds are
\eqref{prep:local-gradient}, with $k+m\le j_1$ and
$r+k+m\le\rprep$. Every additional term contains a gap and thus has
the same amplitude and approximation factor. The claims for one further
material or magnetic Lie derivative follow from
Lemma~\ref{setup:calculus}(ii)--(iii). On an overlap, we may use
either local background and repeat the argument. When both backgrounds
are local evolutions, their difference is a difference of
$\Delta_n^\pm$, so the estimates retain
$\varepsilon_{q+1}^{\gamma_\ell}$.
\end{proof}

\subsection{Background fields and pressure}

\begin{proposition}[Local background bounds]
\label{prep:background-verification}
The local background satisfies
Assumption~\ref{galbrun:linear-assumptions}(i) and (ii) through
derivative order $\rprep$. Its pressure satisfies
\begin{equation}\label{prep:local-pressure-good}
 \begin{aligned}
 \nabla p_{\ell,n}&\in\mathcal C_{\rprep-1}(C\lambda_q\delta_q;\lambda_q,\mathrm c)
 \cap\mathcal K(C\lambda_q\ell^{-\alpha}\delta_q)\\
 &\quad\cap\mathcal O_{\infty,j_1}(C\lambda_q\ell^{-\alpha}\delta_q;\ell^{-1}),
 \end{aligned}
\end{equation}
where the $\mathcal C$ and $\mathcal K$ classes use the local transports.
The ordinary bounds for the field gradients are given in
\eqref{prep:local-field-classes}.
\end{proposition}
\begin{proof}
All transports in this proof belong to the local background
$(v_{\ell,n},B_{\ell,n})$.

\emph{1. Local coefficients.}
The local terms in \eqref{galbrun:separate-lower-operators} and
\eqref{mom:H2} contain $\nabla v_{\ell,n}$, $\nabla B_{\ell,n}$,
$\nabla(\ddiv R_\ell-\nabla p_{\ell,n})$, the stress, and products
of two field gradients. Apply
Lemma~\ref{prep:finite-conversion} to these coefficients and use the
product rule. The comparison in its proof, with $s=0$, also bounds
derivatives of the fields of positive order with amplitudes
$\delta_q^{1/2}$ for
$v_{\ell,n}$ and $\delta_{B,q}^{1/2}$ for $B_{\ell,n}$. The
undifferentiated magnetic field is only bounded by a constant.
Lemma~\ref{high:component-lie-conversion} gives the Lie derivative
estimates using coefficient derivatives with material and magnetic
orders summing to at most $k+m-1$. In particular,
\begin{align*}
 \|D_{t,n}^kD_{B,n}^m\nabla v_{\ell,n}\|_{r+\alpha}
 &\le C\ell^{-\alpha}\varepsilon_{q+1}^{\gamma_\ell}
       \lambda_q^r\tau_c^{-k-1}\lambda_\parallel^m,\\
 \|D_{t,n}^kD_{B,n}^m\nabla B_{\ell,n}\|_{r+\alpha}
 &\le C\ell^{-\alpha}\varepsilon_{q+1}^{\gamma_\ell}
       \varepsilon_{q+1}^{(\gamma_a-(b-1)\gamma_\parallel)/b}
       \lambda_q^r\tau_c^{-k}\lambda_\parallel^{m+1},\\
 \|D_{t,n}^kD_{B,n}^m\nabla(\ddiv R_\ell-\nabla p_{\ell,n})\|_{r+\alpha}
 &\le C\ell^{-\alpha}\varepsilon_{q+1}^{2\gamma_\ell}
       \lambda_q^r\tau_c^{-k-2}\lambda_\parallel^m,
\end{align*}
Since $\gamma_\ell>3\gamma_S$, the factors
$\varepsilon_{q+1}^{\gamma_\ell}$ and
$\varepsilon_{q+1}^{2\gamma_\ell}$ are bounded by
$\lambda_{q+1}^{-3\alpha}$ and $\lambda_{q+1}^{-6\alpha}$,
respectively. These give the small factors required in
\eqref{galbrun:background-low} and in the nonlocal terms below.

\emph{2. Bounds for the nonlocal operators.}
Apply Proposition~\ref{high:cz-iterated} to $(1,v_{\ell,n})$ and
$(0,B_{\ell,n})$. By Lemma~\ref{principal:finite-multiplier}, the
$k$ material and $m$ magnetic commutators have their respective costs;
in particular, a magnetic derivative of a coefficient already
differentiated in the material direction still costs
$\lambda_\parallel$. For a term containing $p$ gradients,
\eqref{high:cz-iterated-estimate} and Step 1 give the factor
$(\ell^{-\alpha}\lambda_{q+1}^{-3\alpha})^p$. This absorbs the
additional H\"older factors in \eqref{high:cz-scaled-output}, whose
product is at most $\lambda_{q+1}^{p\alpha}$, because
\[
 (\ell^{-\alpha}\lambda_{q+1}^{-3\alpha})^p\lambda_{q+1}^{p\alpha}
 \le\lambda_{q+1}^{-p\alpha}.
\]
The term involving $\nabla(\ddiv R_\ell-\nabla p_{\ell,n})$ has the cost of two material
derivatives and the small factor $\ell^{-\alpha}\lambda_{q+1}^{-6\alpha}$.
This factor is no larger than the square of the gradient factor, so
the preceding argument also controls $\mathscr H_2$. Consequently the
source keeps its prescribed H\"older weight without a further loss from
the coefficients.

For ordinary time derivatives, differentiate the Cartesian coefficient
formulas at fixed $x$. The spatial multipliers commute with
$\partial_t$. In terms that also contain material or magnetic
derivatives, use $\partial_t=D_{t,n}-v_{\ell,n}\cn$ and the product
rule. At time order $h$, this introduces at most $h$ additional material
and spatial derivatives in sum. The coefficient estimates therefore
apply when
\[
 h+k+m\le j_1,\qquad r+h+k+m+3\le\rprep.
\]
These are the derivative restrictions for the nonlocal estimates.

\emph{3. Material and magnetic derivatives of the pressure.}
Subtracting the two momentum equations and
using \eqref{prep:frame-comparison} gives the exact identity
\begin{equation}
 \nabla(p_{\ell,n}-p_q)
 =\ddiv(R_\ell-R_q)+D_{t,n}G_n^v+G_n^v\cn v_q
    -D_{B,n}G_n^B-G_n^B\cn B_q.
\end{equation}
For $r\ge1$, Lemma~\ref{high:ordinary-transport-comparison} estimates
the right-hand side using the ordinary gap bounds and the moment
estimate for $R_\ell-R_q$. Every field term contains a gap, giving
\[
 \|D_{t,n}^kD_{B,n}^m\nabla(p_{\ell,n}-p_q)\|_{r-1+\alpha}
 \le C\ell^{-\alpha}\varepsilon_\ell^2
                 \lambda_q^{r+k+m}\delta_q^{1/2},
\]
and the same bound holds in the integer norm.
The material derivative in the pressure identity uses ordinary time
derivatives through $k+1\le j_1+1$ for the fields and $k\le j_1$ for
the stress. The Taylor estimate uses spatial and ordinary time derivatives of order at most
$r+k+m+2m_0+1$, which is also available. To compare the local and old
transports of $\nabla p_q$, we use the momentum equation in the form
\[
 \nabla p_q=\ddiv R_q-\partial_tv_q-v_q\cn v_q+B_q\cn B_q.
\]
At time order $h\le j_1$, this formula uses the fields through order
$h+1$ and the stress through order $h$. Thus
Lemma~\ref{high:ordinary-transport-comparison} gives the same error for
the change of transports. With amplitude $\delta_q$, the inequality
\eqref{prep:conversion-margin} absorbs both errors into
$C\lambda_q^r\delta_q\tau_c^{-k}\lambda_\parallel^m$. The remaining
term is controlled at this scale by the old bounds for material and magnetic derivatives of the pressure,
Lemma~\ref{high:transport-gradient}, and interpolation between
consecutive spatial orders. Applying the reverse direction of
Lemma~\ref{high:transport-gradient} with the local gradient bounds
\eqref{prep:local-gradient} yields
\[
 \|D_{t,n}^kD_{B,n}^m\nabla p_{\ell,n}\|_{r-1}
 \le C\lambda_q^r\delta_q\tau_c^{-k}\lambda_\parallel^m,
 \qquad r+k+m\le\rprep.
\]
The H\"older bound has the factor $\ell^{-\alpha}$, which gives the
sharp class in \eqref{prep:local-pressure-good}. For every periodic
scalar $f$, divergence freedom implies
\[
 \int D_{t,n}f\,\dd x=\partial_t\int f\,\dd x,\qquad
 \int D_{B,n}f\,\dd x=0.
\]
Hence every composition of material and magnetic derivatives of the normalized pressure has zero
mean. Poincar\'e's inequality then bounds the pressure itself by its
differential.

\emph{4. Ordinary pressure and lossy bounds.}
Divergence of the local momentum equation gives
\begin{equation}
 \Delta p_{\ell,n}=\partial_i\partial_j(R_\ell)_{ij}
 -\partial_i(v_{\ell,n})_j\partial_j(v_{\ell,n})_i
 +\partial_i(B_{\ell,n})_j\partial_j(B_{\ell,n})_i,
\end{equation}
in which both factors in each quadratic term are gradients.
The field estimate below follows from the lossy version of
\eqref{prep:whole-patch-fields}. For the pressure and
$\nabla(\ddiv R_\ell-\nabla p_{\ell,n})$, differentiate the Poisson equation $k\le j_1$ times in time and use
the field estimates, the smoothing bounds for $R_\ell$, and the
order-zero estimate for $\nabla^2\Delta^{-1}$. For every $r\ge0$,
we obtain
\[
 \begin{aligned}
 \|\partial_t^k\nabla(v_{\ell,n},B_{\ell,n})\|_{r+\alpha}
 &\le C\lambda_q\ell^{-r-k-\alpha}\delta_q^{1/2},
 &&k\le j_1+1,\\
 \|\partial_t^k\nabla(D_{t,n}v_{\ell,n}-D_{B,n}B_{\ell,n})\|_{r+\alpha}
 &\le C\lambda_q^2\ell^{-r-k-2\alpha}\delta_q,
 &&k\le j_1.
 \end{aligned}
\]
The mean-zero normalization of the pressure fixes its additive
constant. These estimates prove the ordinary coefficient bounds in (ii);
no additional time derivative of the stress is needed. Expanding the
local transports in ordinary derivatives gives
\[
 \begin{aligned}
 \nabla v_{\ell,n},\,\nabla B_{\ell,n}
 &\in\mathcal K_{j_1+1}(C\lambda_q\ell^{-\alpha}\delta_q^{1/2}),\\
 \nabla(\ddiv R_\ell-\nabla p_{\ell,n})
 &\in\mathcal K(C\lambda_q^2\ell^{-2\alpha}\delta_q),\\
 \nabla p_{\ell,n}
 &\in\mathcal K(C\lambda_q\ell^{-\alpha}\delta_q),
 \end{aligned}
\]
which also proves the $\mathcal K$ inclusion in
\eqref{prep:local-pressure-good}.
Conversely, expanding
$\partial_t=D_{t,n}-v_{\ell,n}\cn$ in the pressure estimate gives
its $\mathcal O_{\infty,j_1}$ bound. At time order $h$, this uses
at most $h$ material derivatives of the pressure and ordinary
derivatives of the velocity through order $h-1$. The velocity
bounds and $\ell^{1-\alpha}\lambda_q\delta_q^{1/2}\le1$ give the
cost $\ell^{-1}$ for each derivative, with the same amplitude
$C\lambda_q\ell^{-\alpha}\delta_q$.
\end{proof}

\begin{remark}[Inductive pressure bounds]
The proposition gives the local pressure bounds from the inductive
pressure bounds. We prove the corresponding bounds for the next iterate
in Section~\ref{sec:iteration-closure}, using the pressure differential
and the comparison of the old and new material and magnetic derivatives.
\end{remark}
\section{Charts and Amplitudes}
\label{sec:charts-and-amplitudes}

We use the commuting local velocity and magnetic flows to construct
coordinates on spatial patches of size $\lambda_\parallel^{-1}$ and
time intervals of size $\tau_c$. A quadratic partition of unity and a
decomposition in six fixed directions then express a smoothing of
$\delta_{q+1}\IId-R_q$ as a sum of positive rank-one tensors. We
define the amplitudes from their coefficients and choose temporal
profiles for which the Galbrun forcing has zero mean in the fast
variable.

\subsection{Charts and partition}
\label{ssec:material-partition}

For each slow time interval indexed by $n$, the local fields
$(v_{\ell,n},B_{\ell,n})$ satisfy the induction equation and hence
$[D_{t,n},D_{B,n}]=0$. Moreover, \eqref{prep:local-gradient} gives
\begin{equation}\label{part:flow-smallness}
 \begin{aligned}
 \tau_c[v_{\ell,n}]_1
 &\le C\ell^{-\alpha}\varepsilon_{q+1}^{\gamma_\ell}=o(1),\\
 \lambda_\parallel^{-1}[B_{\ell,n}]_1
 &\le C\ell^{-\alpha}\varepsilon_{q+1}^{\gamma_\ell+
              (\gamma_a-(b-1)\gamma_\parallel)/b}=o(1),\\
 \tau_c&\le\tfrac14\lambda_\parallel^{-1}.
 \end{aligned}
\end{equation}
Proposition~\ref{frobenius:partition}, on the time interval $[t_n-4\tau_c,t_n+4\tau_c]$,
therefore provides volume-preserving charts on spatial patches of
radius comparable to $\lambda_\parallel^{-1}$, together with a
subordinate quadratic partition of unity, satisfying
\begin{equation}
 D_{t,n}\Psi_{n,J}=0,\qquad D_{B,n}\Psi_{n,J}=c_{n,J}^{-1}e_3,\qquad
 \det\DD\Psi_{n,J}=1,\qquad \tfrac13\le c_{n,J}\le8,
\end{equation}
with $\DD\Psi_{n,J}$ uniformly close to a rotation. Here the fixed bounds
on $c_{n,J}$ follow from $1/8\le|B_{\ell,n}|\le3$, obtained from
\eqref{full:pointwise-induction} and the preparation error. Also,
\begin{equation}
 D_{t,n}\chi_{n,J}=0,\qquad \sum_J\chi_{n,J}^2=1.
\end{equation}
The supports have bounded overlap and eight spatial colors, and fixed
cut off regions of patches of the same color are disjoint.

To obtain the chart and partition bounds, we apply
Proposition~\ref{frobenius:estimates} and the corresponding partition
estimates with
\[
 (\Lambda,N)=(\lambda_q,\rprep+3)\quad\hbox{and}\quad(\ell^{-1},N)
 \hbox{ for every }N,\qquad k_{\max}=j_1+2,
\]
For the material and magnetic derivatives, the parameters are
\[
 \mathrm a_t=\tau_c^{-1},\qquad \mathrm a_B=\lambda_\parallel,\qquad
 j_{\max}=j_1+1,\qquad N_{\rm char}=\rprep+1.
\]
The ordinary derivative hypotheses follow from
\eqref{prep:whole-patch-fields}, since
$\rgood-3\ge\rprep+3$ and the lossy bounds hold at every spatial
order. For the background gradients, the hypotheses are
$k+m\le j_{\max}-1=j_1$ and $r+k+m\le\rprep$, so
\eqref{prep:local-gradient} applies. Consequently,
$\DD\Psi_{n,J}$ and its inverse belong to
\[
 \mathcal O_{\rprep+2,j_1+1}(C;\lambda_q)\cap\mathcal O_{\infty,j_1+1}(C;\ell^{-1})
 \quad\hbox{and}\quad
 \mathcal C_{\rprep+1,j_1+1}(C;\lambda_q,\mathrm c).
\]
The maps satisfy the corresponding bounds with one further spatial
derivative. We chose $N$ one order
above the required ordinary derivative range so that interpolation
gives the H\"older estimates, with factor $\ell^{-\alpha}$.
The transport identities for the cutoffs give
\begin{equation}
 D_{t,n}^kD_{B,n}^m(\eta_n\chi_{n,J})=\eta_n^{(k)}D_{B,n}^m\chi_{n,J}
\end{equation}
and, on the fixed neighborhood of $\supp\varrho$,
\begin{equation}
 \sum_{n,J}\eta_n^2\chi_{n,J}^2=1 .
\end{equation}
The cut off regions are wide enough for every smoothing trajectory.
Smoothing along the velocity flow changes time by
$O(\tau_a)=o(\tau_c)$ and fixes $\Psi_{n,J}$. Magnetic smoothing changes
only its third coordinate, by
$O(\varepsilon_\tau\lambda_\parallel^{-1})=o(\lambda_\parallel^{-1})$,
while spatial convolution displaces points by
$O(\ell)=o(\lambda_\parallel^{-1})$. The bounded chart differential
therefore ensures that the zero extensions used in these operations
are well defined along all the trajectories.

\subsection{Directions and profiles}

\begin{lemma}[Positive rank-one stress decomposition]
\label{prep:geometric-decomposition}
There are six unit directions $\mathscr Z$ and linear forms
$q_\zeta$ on $\sym_3$ such that
\begin{equation}\label{prep:geometric-data}
 q_\zeta(S)=\tfrac54\zeta\cdot S\zeta-\tfrac14\operatorname{tr}S,
 \qquad S=\sum_{\zeta\in\mathscr Z}q_\zeta(S)\zeta\otimes\zeta.
\end{equation}
They satisfy $q_\zeta(\IId)=1/2$ and $q_\zeta(S)\ge1/4$ in a fixed
neighborhood of $\IId$. Each $\zeta$ has an oriented orthonormal
completion $(k_\zeta,\nu_\zeta,\zeta)$ with $k_\zeta\cdot e_3=0$.
\end{lemma}
\begin{proof}
Put $g=(1+\sqrt5)/2$ and normalize the six vectors
\[
 (0,1,g),\ (0,-1,g),\ (1,g,0),\ (-1,g,0),\ (g,0,1),\ (g,0,-1)
\]
by $(1+g^2)^{-1/2}$. Since $g^2=g+1$, their tensors
$\Pi_\zeta=\zeta\otimes\zeta$ satisfy $\Pi_\zeta:\Pi_{\zeta'}=1$ for
$\zeta=\zeta'$ and $=1/5$ otherwise, and $\sum_\zeta\Pi_\zeta=2\IId$.
Their Gram matrix is
\[
 (4/5)\IId_6+(1/5)\mathbf1\mathbf1^T.
\]
Its eigenvalues are $4/5$ and $2$, so it is invertible and the six
tensors form a basis of $\sym_3$. The displayed linear forms satisfy
\[
 q_\zeta(\Pi_{\zeta'})=\delta_{\zeta,\zeta'},\qquad
 q_\zeta(\IId)=1/2.
\]
The first identity proves the decomposition. The second gives positivity
in a fixed neighborhood of the identity, by continuity of the six
linear forms.
Choose a unit $k_\zeta\in e_3^\perp\cap\zeta^\perp$ and put
$\nu_\zeta=\zeta\times k_\zeta$.
\end{proof}

For $I=(n,J,\zeta)$, use the corresponding time interval, chart and direction:
\[
 \begin{gathered}
 \Psi_I=\Psi_{n,J},\qquad \chi_I=\chi_{n,J},\qquad \eta_I=\eta_n,\qquad c_I=c_{n,J},
 \qquad
 (k_I,\nu_I,\zeta_I)=(k_\zeta,\nu_\zeta,\zeta),\\
 y_I^\eta=\Psi_I\cdot\eta,\qquad
 \frac{\partial}{\partial y_I^\eta}=(\DD\Psi_I)^{-1}\eta,
 \qquad \eta\in\{k_I,\nu_I,\zeta_I\}.
 \end{gathered}
\]
The same index on a local background denotes its time interval $n$;
in particular,
\[
 \begin{gathered}
 v_{\ell,I}=v_{\ell,n},\qquad B_{\ell,I}=B_{\ell,n},\\
 D_{t,I}=D_{t,n},\qquad D_{B,I}=D_{B,n},\qquad
 \mathcal D_{t,I}=\mathcal D_{t,n}.
 \end{gathered}
\]
For a fixed chart $(n,J)$, write $M_{n,J}=\DD\Psi_{n,J}$.
Applying \eqref{prep:geometric-data} to $M_{n,J}SM_{n,J}^T$
and transforming back to physical coordinates gives the decomposition
\[
 S=\sum_{\zeta\in\mathscr Z}q_\zeta(M_{n,J}SM_{n,J}^T)
       (M_{n,J}^{-1}\zeta)\otimes(M_{n,J}^{-1}\zeta).
\]
The two identities in \eqref{frobenius:matrix-mixeds} give
\[
 \begin{aligned}
 D_{t,n}M_{n,J}&=-M_{n,J}\DD v_{\ell,n},\\
 D_{B,n}M_{n,J}&=-M_{n,J}\DD B_{\ell,n},\\
 D_{t,n}(M_{n,J}SM_{n,J}^T)
 &=M_{n,J}\bigl(D_{t,n}S-(\DD v_{\ell,n})S
             -S(\DD v_{\ell,n})^T\bigr)M_{n,J}^T,\\
 D_{B,n}(M_{n,J}SM_{n,J}^T)
 &=M_{n,J}\bigl(D_{B,n}S-(\DD B_{\ell,n})S
             -S(\DD B_{\ell,n})^T\bigr)M_{n,J}^T.
 \end{aligned}
\]
The expressions in parentheses are the material and magnetic Lie
derivatives of the contravariant tensor $S$. Linearity of $q_\zeta$
therefore gives, with $I=(n,J,\zeta)$,
\begin{equation}\label{prep:dual-naturality}
 \begin{aligned}
 D_{t,n}q_\zeta(M_{n,J}SM_{n,J}^T)
 &=q_\zeta(M_{n,J}\mathcal D_{t,n}S M_{n,J}^T),\\
 D_{B,n}q_\zeta(M_{n,J}SM_{n,J}^T)
 &=q_\zeta(M_{n,J}\mathcal L_{B_{\ell,n}}S M_{n,J}^T),\\
 \mathcal D_{t,I}\bigl(\partial_{y_I^{\zeta_I}}\otimes\partial_{y_I^{\zeta_I}}\bigr)
 &=\mathcal L_{B_{\ell,I}}\bigl(\partial_{y_I^{\zeta_I}}\otimes\partial_{y_I^{\zeta_I}}\bigr)=0.
 \end{aligned}
\end{equation}
These identities iterate in the two commuting transports.

We assign a temporal profile to each combination of time parity,
spatial color and direction, giving $2\cdot8\cdot6=96$ types. Contributions
of the same type have disjoint complete supports. For distinct types,
we place the profiles in disjoint time subintervals of one period, as
in the following lemma.

\begin{lemma}[Separated profiles]\label{amplitude:profiles}
For every fixed integer $j_0\ge1$, these types admit profiles
$\alpha_I$, constant within each type, whose normalized primitives
through order $j_0$ are compactly supported in the assigned time subinterval and obey
\begin{equation}\label{amplitude:profile-identities}
 \begin{gathered}
 \partial_\tau\alpha_I^{[j+1]}=\alpha_I^{[j]}\quad(0\le j<j_0),\qquad
 \alpha_I^{[0]}=\alpha_I,\qquad
 \int_0^1\alpha_I^2\,\dd\tau=1,\\
 \int_0^1\alpha_I^{[j]}\,\dd\tau=0
 \quad(0\le j\le j_0).
 \end{gathered}
\end{equation}
For each prescribed derivative order and number of further normalized
periodic primitives, the resulting profiles are uniformly bounded
over the types.
\end{lemma}
\begin{proof}
Choose a nonzero smooth bump $\psi_I$ in its assigned time subinterval and set
\[
 \alpha_I^{[j]}=\psi_I^{(j_0+1-j)}/\|\psi_I^{(j_0+1)}\|_{L^2},
 \qquad 0\le j\le j_0.
\]
The recurrence follows by differentiation, and the denominator gives
the required square normalization. Since every numerator is a
derivative of a compactly supported bump, all the prescribed primitives
have zero mean and remain supported in the assigned subinterval.
For any mean-zero periodic function $f$,
\[
 \partial_\tau^{-1}f(\tau)
 =\int_0^\tau f(s)\,\dd s-\int_0^1\int_0^u f(s)\,\dd s\,\dd u,
 \qquad \|\partial_\tau^{-1}f\|_\infty\le2\|f\|_\infty,
\]
and iteration gives the remaining finite bounds.
\end{proof}

\subsection{Coefficients and forcing}

We smooth $\delta_{q+1}\IId-R_q$ by the pullback averages of
Proposition~\ref{moll:mixed-calculus}, choosing kernels with vanishing
moments through order $j_0$ and setting
\begin{equation}
 \mathcal J_n=\mathcal J_{s_t,s_B}^{(n)}\mathcal J_\ell,
 \qquad s_t=\tau_a,\qquad s_B=\varepsilon_\tau\lambda_\parallel^{-1}.
\end{equation}
The spatial convolution is applied once, before the smoothing by tensor pullbacks along the flows:
\[
 \mathcal J_n(\delta_{q+1}\IId-R_q)
 =\mathcal J_{s_t,s_B}^{(n)}(\delta_{q+1}\IId-R_\ell).
\]
Define the complete amplitudes by
\begin{equation}\label{amplitude:coefficient-formula}
 a_I=\varrho\eta_I\chi_I
 \left[q_{\zeta_I}\!\left((\DD\Psi_I)
    \mathcal J_n(\delta_{q+1}\IId-R_q)(\DD\Psi_I)^T\right)\right]^{1/2}.
\end{equation}
By the kernel bounds and the estimates for tensor pullbacks,
\[
 \|\mathcal J_n\IId-\IId\|_0
 \le C(s_t[v_{\ell,n}]_1+s_B[B_{\ell,n}]_1)=o(1),\qquad
 \|\mathcal J_nR_q\|_0\le C\|R_q\|_0=o(\delta_{q+1}),
\]
where the second estimate uses the inductive stress bound.
Since $\DD\Psi_I$ is close to a rotation, the tensor in
\eqref{amplitude:coefficient-formula}, divided by $\delta_{q+1}$,
belongs to the neighborhood of $\IId$ in
Lemma~\ref{prep:geometric-decomposition}. By linearity of
$q_{\zeta_I}$, the expression under the square root is therefore
between $c\delta_{q+1}$ and $C\delta_{q+1}$. Hence $a_I$ is smooth.
The argument requires no positivity of the smoothing kernels.

For a fixed $n$, the notation $\sum_{I=(n,J,\zeta)}$ denotes summation
over $J$ and $\zeta$, whereas $\sum_I$ includes all indices. The
geometric decomposition and quadratic partition of unity imply
\begin{equation}
 \sum_{I=(n,J,\zeta)}a_I^2
        \partial_{y_I^\zeta}\otimes\partial_{y_I^\zeta}
 =\varrho^2\eta_n^2\mathcal J_n(\delta_{q+1}\IId-R_q),
\end{equation}
with the cutoffs outside the smoothing operator. The prescribed
Galbrun forcing is
\begin{equation}\label{amplitude:centered-source}
 \mathsf F_n(t,x)=\sum_{I=(n,J,\zeta)}
       [1-\alpha_I(t/\tau_a)^2]a_I(t,x)^2
              \partial_{y_I^\zeta}\otimes\partial_{y_I^\zeta},
\end{equation}
whose temporal factors have zero mean by
\eqref{amplitude:profile-identities}.

\subsection{Amplitude bounds}

Smoothing along the flows gives estimates for material and magnetic
derivatives of the tensor, and hence of the amplitude, when
$p+m>j_1$. By Lemma~\ref{aniso:amplitude-overflow}, each derivative
beyond this range costs an additional factor
$\varepsilon_\tau^{-1}$. We include these factors in the estimates
below. They require no material or magnetic derivatives of the old
stress beyond the inductive range.

\begin{corollary}[Amplitude bounds]
For each index $I$, with transports $D_{t,I},D_{B,I}$, and for
$h\le j_1+1$, $p+m\le j_0+j_1+2$ and $r+h+\min\{p+m,j_1\}\le\rprep$,
\begin{equation}\label{amplitude:amplitude-bounds}
 \|\partial_t^hD_{t,I}^pD_{B,I}^ma_I\|_{r}
 \le C\delta_{q+1}^{1/2}\lambda_q^{r+h}\tau_c^{-p}\lambda_\parallel^m
       \varepsilon_\tau^{-[p+m-j_1]^+},
\end{equation}
and for $h\le j_1+1$, $p+m\le j_0+j_1+2$ and every $r$,
\begin{equation}\label{amplitude:crude-amplitude-bounds}
 \|\partial_t^hD_{t,I}^pD_{B,I}^ma_I\|_{r}
 \le C\delta_{q+1}^{1/2}\ell^{-(r+h)}\tau_a^{-p}s_B^{-m}
 \le C\delta_{q+1}^{1/2}\ell^{-(r+h+p+m)},
\end{equation}
and both hold in the H\"older norm $\|\cdot\|_{r+\alpha}$ with the factor
$\ell^{-\alpha}$. In particular
\begin{equation}\label{amplitude:classes}
 a_I\in\mathcal C_{\rprep}(C\delta_{q+1}^{1/2};\lambda_q,\mathrm c)
 \cap\mathcal K(C\delta_{q+1}^{1/2})
 \cap\mathcal O_{\rprep,j_1+1}(C\delta_{q+1}^{1/2};\lambda_q)
 \cap\mathcal O_{\infty,j_1+1}(C\delta_{q+1}^{1/2};\ell^{-1}),
\end{equation}
with no loss.
\end{corollary}

\begin{proof}
The estimates follow from the flow smoothing in
\eqref{amplitude:coefficient-formula} and the chain rule. In this
argument the ordinary time derivative order of the old stress never
exceeds $j_1$, and neither does the sum of its material and magnetic
derivative orders.

\emph{1. Material and magnetic derivatives.}
The scalar under the square root in \eqref{amplitude:coefficient-formula}
is comparable to $\delta_{q+1}$. On this interval,
\[
 |\dd^hs^{1/2}/\dd s^h|\le C_h\delta_{q+1}^{1/2-h}.
\]
We apply Lemma~\ref{aniso:amplitude-overflow} to
$\delta_{q+1}\IId-R_\ell$. The required Lie derivative estimates
follow from \eqref{prep:stress-good} and the identities
\[
 \begin{aligned}
 \mathcal D_{t,I}(\delta_{q+1}\IId)&=-\delta_{q+1}(\DD v_{\ell,I}+(\DD v_{\ell,I})^T),\\
 \mathcal L_{B_{\ell,I}}(\delta_{q+1}\IId)&=-\delta_{q+1}(\DD B_{\ell,I}+(\DD B_{\ell,I})^T).
 \end{aligned}
\]
Take $e=\varepsilon_\tau$ in the lemma, since
$s_t=\varepsilon_\tau\tau_c$ and
$s_B=\varepsilon_\tau\lambda_\parallel^{-1}$. It gives the material
and magnetic costs $\tau_c^{-p}\lambda_\parallel^m$ and the loss
$\varepsilon_\tau^{-[p+m-j_1]^+}$. The derivatives of the chart
components satisfy the same estimates by
\eqref{prep:dual-naturality} and
Proposition~\ref{frobenius:estimates}.
In the chain rule, the material orders on the factors sum to $p$ and
the magnetic orders sum to $m$. The amplitude is consequently
$\delta_{q+1}^{1/2}$, with the same derivative costs and at most the
stated loss. For the cutoffs, we use
\[
 D_{B,I}\varrho=D_{B,I}\eta_I=0,\qquad D_{t,I}\chi_I=0,\qquad
 |\varrho^{(h)}|\lesssim\tau_c^{-h}.
\]
The identities in Section~\ref{ssec:material-partition} give a factor
$\lambda_\parallel$ for each magnetic derivative of $\chi_I$.
Leibniz' rule and the tame product estimate,
with one factor in the highest spatial norm, now give
\eqref{amplitude:amplitude-bounds} for $h=0$. The smoothing lemma
uses at most $\min\{p+m,j_1\}$ material and magnetic derivatives in
sum of $\delta_{q+1}\IId-R_\ell$, which gives the stated spatial
range.

\emph{2. Ordinary derivatives.}
To estimate ordinary time derivatives, use
$\partial_t=D_{t,I}-v_{\ell,I}\cn$ and
$[D_{t,I},\partial_i]=-(\partial_iv_{\ell,I})\cn$.
At order $h$, expanding $\partial_t^hD_{t,I}^pD_{B,I}^ma_I$ gives products
of velocity derivatives and terms
$\DD^sD_{t,I}^{p+u}D_{B,I}^ma_I$ with $s+u\le h$.
The local background bounds apply to the velocity factors: their
derivative order is at most $r+h-1$, with time order at most
$h-1\le j_1$. Each of the $s$ additional spatial derivatives costs
$\lambda_q$ in the sharp estimate. For the $u$ additional material
derivatives, Step 1 and $\tau_a^{-1}\le\lambda_q$ give
\[
 \tau_c^{-u}\varepsilon_\tau^{-[p+u+m-j_1]^+}
 \le\tau_c^{-u}\varepsilon_\tau^{-u}\varepsilon_\tau^{-[p+m-j_1]^+}
 =\tau_a^{-u}\varepsilon_\tau^{-[p+m-j_1]^+}
 \le\lambda_q^u\varepsilon_\tau^{-[p+m-j_1]^+},
\]
Thus each term is bounded as in
\eqref{amplitude:amplitude-bounds}, since $s+u\le h$. The same
expansion applies to the lossy estimates, with $\ell^{-1}$ in place
of $\lambda_q$ for both additional spatial and material derivatives.

\emph{3. Lossy bounds.}
Apply \eqref{moll:mixed-transfer} to transfer all material and
magnetic derivatives of the smoothed tensor, including its identity
term, to the two kernels. The remaining chart components contain only
spatial derivatives of $\delta_{q+1}\IId-R_\ell$ and of the chart
maps. Their lossy estimates hold at every spatial order. The kernel
costs satisfy
\[
 s_B^{-1}=\tau_a^{-1}\varepsilon_{q+1}^{\gamma_\parallel}
 \le\tau_a^{-1}.
\]
The chain rule for the square root and the cutoff identities therefore give
\[
 \|D_{t,I}^uD_{B,I}^wa_I\|_s\le C\delta_{q+1}^{1/2}\ell^{-s}\tau_a^{-u}s_B^{-w}.
\]
Here we also used $\tau_c^{-1}\le\tau_a^{-1}$ and
$\lambda_\parallel\le s_B^{-1}$ for the cutoff derivatives.
Interpolation with the next spatial order gives the H\"older
estimate. Apply the expansion in Step 2 with the lossy velocity bounds
to obtain \eqref{amplitude:crude-amplitude-bounds}; its second
inequality follows from
$s_B^{-1}\le\tau_a^{-1}\le\lambda_q\le\ell^{-1}$.
This argument requires spatial derivatives of the chart matrices
through order $r+h+1$, and of the maps through order $r+h+2$.
These estimates hold at every spatial order. It does not require an
ordinary time derivative of order $j_1+1$ of the old stress.
Finally, taking $p+m\le j_1$ in the two bounds gives
\eqref{amplitude:classes}, since the additional loss is then absent.
\end{proof}

The amplitude bounds apply to the approximate antiderivative constructed in
Lemma~\ref{principal:material-primitive}. A material derivative of
$D_{t,I}^ja_I$ in that primitive is estimated by
\eqref{amplitude:amplitude-bounds} with $(p,m)=(j+k,m)$, retaining
the stated additional derivative loss.
\section{The Complete Deformation}\label{sec:adapted-perturbation}

Let $(v_q,B_q,p_q,R_q)$ be the relaxed solution at stage $q$.
We construct the next pair of fields by composing the corrector flow
with a principal flow defined in the local background charts. The
corrector solves the equation with the centered source, while the
principal perturbation produces the prescribed quadratic tensor.
To compute the new stress, we also determine how the composed flow
acts on the differences between the actual and local background fields.

\subsection{Construction of the Corrector}

Fix $n$. We apply Theorem~\ref{galbrun:linear-theorem} to the source
\eqref{amplitude:centered-source}, grouping the local contributions
according to their $96$ fixed profile types. Each profile
$1-\alpha_I^2$ has zero mean. Its derivative and normalized periodic
primitive bounds follow from Lemma~\ref{amplitude:profiles}.
The corresponding slow tensor is the sum of
$a_I^2\partial_{y_I^{\zeta_I}}\otimes\partial_{y_I^{\zeta_I}}$ over
indices of that type. Since the coordinate frame has zero material
and magnetic Lie derivatives,
\[
 \begin{aligned}
 &\mathcal D_{t,n}^k\mathcal L_{B_{\ell,n}}^m
   (a_I^2\partial_{y_I^{\zeta_I}}\otimes\partial_{y_I^{\zeta_I}})\\
 &\qquad=D_{t,n}^kD_{B,n}^m(a_I^2)
       \partial_{y_I^{\zeta_I}}\otimes\partial_{y_I^{\zeta_I}}.
 \end{aligned}
\]
Thus \eqref{amplitude:classes}, the frame estimates and
Lemma~\ref{setup:calculus}(i),(iii) give, for every summand,
\[
 a_I^2\partial_{y_I^{\zeta_I}}\otimes\partial_{y_I^{\zeta_I}}
 \in\mathcal C_{\rprep}(C\delta_{q+1};\lambda_q,\mathrm c).
\]
The supports have bounded overlap, also in the surrounding cut off
regions. For a H\"older difference, only summands supported at one of
the two points contribute. Thus each grouped tensor $\mathsf S_b$
satisfies \eqref{galbrun:source-bounds} with
$M_{\rm src}=C\delta_{q+1}$ through derivative order $\rprep$, uniformly
in the number of patches. We group the tensors before applying any
spatial Fourier multiplier. The resulting source is supported in
$\supp(\varrho\eta_n)\Subset I_n^{\rm src}$; its remaining bounds
\eqref{galbrun:crude-source-bounds} follow from the lossy and ordinary
estimates in \eqref{amplitude:classes}. Finally,
Proposition~\ref{prep:background-verification} verifies the background
hypotheses on the full interval, so the theorem applies.

The theorem therefore determines a smooth potential on
$[t_n-3\tau_c,t_n+3\tau_c]\Subset I_n^{\rm bg}$ by
\begin{equation}\label{amplitude:forced-equation}
 \mathscr G_n\Theta_n=\mathscr T\mathsf F_n,
 \qquad \Theta_n(t_n-3\tau_c)=D_{t,n}\Theta_n(t_n-3\tau_c)=0.
\end{equation}
Choose $\tilde\eta_n\in C_c^\infty(I_n^{\rm col})$, equal to one
near $\overline{I_n^{\rm src}}$, with
$|\tilde\eta_n^{(k)}|\le C_k\tau_c^{-k}$. Set
\[
 \begin{gathered}
 \Theta_n^{\rm c}=\tilde\eta_n\Theta_n,\qquad
 \xi_n^{\rm c}=\curl\Theta_n^{\rm c},\\
 \Theta^{\rm c}=\sum_n\Theta_n^{\rm c},\qquad
 \xi^{\rm c}=\sum_n\xi_n^{\rm c}=\curl\Theta^{\rm c}.
 \end{gathered}
\]
Each summand extends smoothly by zero. The localized potential and its
first material and magnetic derivatives satisfy the sharp, coarse and
lossy estimates of Theorem~\ref{galbrun:linear-theorem}(a)--(d),
including in the cut off region after the source support.

\paragraph{\textbf{Corrector flow.}}

The corrector flow and its tensor pushforward are
\[
 \partial_sX_s^{\rm c}=\xi^{\rm c}\circ X_s^{\rm c},\qquad
 X_0^{\rm c}=\IId,\qquad \mathcal U_s^{\rm c}=(X_s^{\rm c})_*.
\]
The LDF is independent of $s$ and may depend on physical time.
The corrector acts on the actual fields by
\begin{equation}\label{direct:reference-fields}
 v_s^{\rm c}=(\partial_tX_s^{\rm c})\circ(X_s^{\rm c})^{-1}
                         +\mathcal U_s^{\rm c} v_q,
 \qquad B_s^{\rm c}=\mathcal U_s^{\rm c} B_q,
\end{equation}
and on each local background by
\[
 v_{\ell,n,s}^{\rm c}
 = (\partial_tX_s^{\rm c})\circ(X_s^{\rm c})^{-1}
      +\mathcal U_s^{\rm c}v_{\ell,n},\qquad
 B_{\ell,n,s}^{\rm c}=\mathcal U_s^{\rm c}B_{\ell,n}.
\]
The Lie derivatives for the corrected local background are
$\mathcal D_{t,n,s}^{\rm c}=\partial_t+\mathcal L_{v_{\ell,n,s}^{\rm c}}$
and $\mathcal L_{B_{\ell,n,s}^{\rm c}}$.
Set
\[
 \mathcal P^{\rm c}f=\int_0^1\mathcal U_s^{\rm c}f\,\dd s.
\]

\paragraph{\textbf{Corrector fields.}}

Since $\tilde\eta_n=1$ on its source support,
\[
 \mathscr G_n\Theta_n^{\rm c}
 =\mathscr T\mathsf F_n
        +[\mathscr G_n,\tilde\eta_n]\Theta_n.
\]
The commutator is supported where the cutoff varies, by
\eqref{galbrun:cutoff-identity}, and vanishes in the cut off
region before the source support because the potential does. To separate
the linear solution
from the flow correction and the terms containing the initial gaps
$G_n^\pm$, define
\[
 w_n^{\mathrm c,\mathrm{main}}\pm b_n^{\mathrm c,\mathrm{main}}
       =\curl\mathscr A_{\ell,n}^\pm\Theta_n^{\rm c},
\]
and their one-form remainders by
\begin{equation}\label{mom:remainder-potentials}
 \begin{aligned}
 \Theta_n^{\mathrm c,\mathrm{rem},\pm}
   &=(\mathcal P^{\rm c}-\IId)\mathscr A_{\ell,n}^\pm\Theta_n^{\rm c}
       +\mathcal P^{\rm c}\mathcal L_{G_n^\pm}\Theta_n^{\rm c},\\
 \Theta_{w,n}^{\mathrm c,\mathrm{rem}}
   &=\tfrac12(\Theta_n^{\mathrm c,\mathrm{rem},+}+\Theta_n^{\mathrm c,\mathrm{rem},-}),\\
 \Theta_{b,n}^{\mathrm c,\mathrm{rem}}
   &=\tfrac12(\Theta_n^{\mathrm c,\mathrm{rem},+}-\Theta_n^{\mathrm c,\mathrm{rem},-}).
 \end{aligned}
\end{equation}
Their curls are denoted by
\[
 w_n^{\mathrm c,\mathrm{rem}}=\curl\Theta_{w,n}^{\mathrm c,\mathrm{rem}},\qquad
 b_n^{\mathrm c,\mathrm{rem}}=\curl\Theta_{b,n}^{\mathrm c,\mathrm{rem}}.
\]

The complete corrector increments are
\[
 \begin{aligned}
 w_n&=w_n^{\mathrm c,\mathrm{main}}+w_n^{\mathrm c,\mathrm{rem}},
 &b_n&=b_n^{\mathrm c,\mathrm{main}}+b_n^{\mathrm c,\mathrm{rem}},\\
 w^{\rm c}&=\sum_nw_n=v_1^{\rm c}-v_q,
 &b^{\rm c}&=\sum_nb_n=B_1^{\rm c}-B_q.
 \end{aligned}
\]
Indeed, \eqref{eq:signed-increments} gives
\[
 \begin{aligned}
 (v_1^{\rm c}-v_q)\pm(B_1^{\rm c}-B_q)
 &=\curl\sum_n\mathcal P^{\rm c}
       \bigl(\mathscr A_{\ell,n}^\pm
                +\mathcal L_{G_n^\pm}\bigr)\Theta_n^{\rm c}\\
 &=\curl\sum_n\Bigl[\mathscr A_{\ell,n}^\pm\Theta_n^{\rm c}
       +(\mathcal P^{\rm c}-\IId)\mathscr A_{\ell,n}^\pm\Theta_n^{\rm c}
       +\mathcal P^{\rm c}\mathcal L_{G_n^\pm}\Theta_n^{\rm c}\Bigr]\\
 &=\sum_n\bigl(w_n^{\mathrm c,\mathrm{main}}
                     \pm b_n^{\mathrm c,\mathrm{main}}\bigr)
      +\curl\sum_n\Theta_n^{\mathrm c,\mathrm{rem},\pm}.
 \end{aligned}
\]
The last line is \eqref{mom:remainder-potentials}, which proves the
decomposition. For the term containing the background difference,
Cartan's formula gives
\[
 \mathcal L_G\Theta=\iota_Gd\Theta+d(\iota_G\Theta),
\]
The exact form has zero curl. We may therefore estimate the contraction
of $G$ with $d\Theta$, which does not differentiate the initial gap.

\paragraph{\textbf{Pressure.}}

The mean-zero pressure of the linear corrector is
\begin{equation}\label{gg:pressure}
 \begin{aligned}
 \pi_n=\Delta^{-1}\ddiv\ddiv\bigl(&\mathsf F_n
       -2(\mathcal D_{t,n}\Theta_n^{\rm c})\times\nabla v_{\ell,n}\\
       &+2(\mathcal L_{B_{\ell,n}}\Theta_n^{\rm c})\times\nabla B_{\ell,n}\bigr),
 \end{aligned}
\end{equation}
with the matrices $Y\times\nabla z$ of \eqref{mom:Kdefinition}. This is
the potential formula in Corollary~\ref{galbrun:split-pressure} with
source $\mathsf F_n$. The cutoff force is a curl and
contributes no pressure.

\subsection{Background gaps}

Recall the background differences defined in the preparation:
\[
 G_n^v=v_q-v_{\ell,n},\qquad
 G_n^B=B_q-B_{\ell,n},\qquad
 G_n^\pm=G_n^v\pm G_n^B.
\]

Subtracting the affine actions on the actual fields and the local backgrounds in
\eqref{direct:reference-fields}, the common time derivative cancels, giving
\begin{equation}\label{mom:transported-gaps}
 v_s^{\rm c}-v_{\ell,n,s}^{\rm c}=\mathcal U_s^{\rm c}G_n^v,
 \qquad
 B_s^{\rm c}-B_{\ell,n,s}^{\rm c}=\mathcal U_s^{\rm c}G_n^B.
\end{equation}
Volume preservation and the path identities imply that both pairs
satisfy incompressibility and induction. Their Lie derivatives
intertwine with $\mathcal U_s^{\rm c}$ by Lemma~\ref{direct:flatness}.
Below we compute the additional change in the gaps under the principal flow.

\subsection{Construction of the Principal Perturbation}
\label{ssec:exact-phase-geometry}
\paragraph{\textbf{Physical chart convention.}}

The construction uses local coordinates, but all spatial estimates
differentiate Cartesian components. In particular, the coordinate
vectors $\frac{\partial}{\partial y_I^\eta}$ and dual spatial one-forms
$d y_I^\eta$ are variable fields in these estimates.

Fix an index $I$. In the corresponding local background chart,
the two physical transport operators are constant-coefficient
coordinate derivatives:
\begin{equation}\label{part:physical-chart-derivatives}
 D_{t,I}=\partial_t\big|_{y_I},\qquad
 D_{B,I}=c_I^{-1}e_3\cdot\nabla_{y_I},\qquad
 D_{t,I}y_I^{k_I}=D_{B,I}y_I^{k_I}=0.
\end{equation}
Here ordinary $\partial_t$ in Cartesian coordinates fixes $x$;
$\partial_t|_{y_I}$ fixes the chart coordinate. The differentials
$d y_I^\eta$ in the potentials are spatial one-forms.

Define the spatial profile by
\begin{equation}
 \varphi(u)=\sqrt2\sin u,\qquad
 \varphi'(u)^2=1+\cos(2u),\qquad
 \varphi(u)^2=1-\cos(2u),
\end{equation}
so that $\varphi'^2$ and $\varphi^2$ have mean one over a period and
$\varphi\varphi'$ has mean zero. Its Fourier modes in the fast phase variable will be used
in the quadratic stress.

All fields below are functions of physical $(t,x)$, equivalently
of $(t,y_I)$. Their profile factors are $\alpha_I(t/\tau_a)$ and
$\varphi(\lambda_{q+1}y_I^{k_I})$; we sometimes suppress these
arguments in local formulas. Primes always differentiate the fixed
one-variable profile. Thus, for a physical scalar $A$, differentiating
the temporal factor gives
\begin{equation}
 \begin{aligned}
 D_{t,I}\bigl(\alpha_I(t/\tau_a)A\bigr)
 &=\tau_a^{-1}\alpha_I'(t/\tau_a)A+\alpha_I(t/\tau_a)D_{t,I}A,\\
 D_{B,I}\bigl(\alpha_I(t/\tau_a)A\bigr)&=\alpha_I(t/\tau_a)D_{B,I}A.
 \end{aligned}
\end{equation}
Both transports annihilate the spatial profile factor. The chart
coordinate derivatives commute with each other and with these two
transports. Thus for any physical scalar $f(t,y_I)$,
\begin{align}
 d f&=(\partial_{y_I^{k_I}}f)d y_I^{k_I}
      +(\partial_{y_I^{\nu_I}}f)d y_I^{\nu_I}
      +(\partial_{y_I^{\zeta_I}}f)d y_I^{\zeta_I},
 \\
 \curl(f d y_I^{\nu_I})
 &=-\partial_{y_I^{\zeta_I}}f\,\frac{\partial}{\partial y_I^{k_I}}
     +\partial_{y_I^{k_I}}f\,\frac{\partial}{\partial y_I^{\zeta_I}}.
 \label{part:physical-curl}
\end{align}
The last identity uses the volume-preserving, oriented chart. In
particular, differentiating $A\varphi(\lambda_{q+1}y_I^{k_I})$ in the
$y_I^{k_I}$ direction gives
\[
 (\partial_{y_I^{k_I}}A)\varphi(\lambda_{q+1}y_I^{k_I})
 +\lambda_{q+1}A\varphi'(\lambda_{q+1}y_I^{k_I}).
\]
For tensors, we apply these identities to the chart components and use
the vanishing Lie derivatives of the coordinate frame.

\paragraph{\textbf{Principal potential and composition.}}

The principal potential uses a localized approximate antiderivative of
$\alpha_Ia_I$ along the velocity flow.
The mean-zero periodic primitives of Lemma~\ref{amplitude:profiles} satisfy
\begin{equation}
 \frac{\dd}{\dd t}\alpha_I^{[j]}(t/\tau_a)
   =\tau_a^{-1}\alpha_I^{[j-1]}(t/\tau_a),
 \qquad 1\le j\le j_0.
\end{equation}
Define
\begin{equation}
 \begin{aligned}
 \mathfrak a_I&=\sum_{j=0}^{j_0-1}(-1)^j\tau_a^{j+1}
                  \alpha_I^{[j+1]}(t/\tau_a)D_{t,I}^ja_I,\\
 a_I^{\rm c}&=(-1)^{j_0-1}\tau_a^{j_0}
                  \alpha_I^{[j_0]}(t/\tau_a)D_{t,I}^{j_0}a_I.
 \end{aligned}
\end{equation}
Lemma~\ref{principal:material-primitive} gives
\begin{equation}\label{principal:primitive-identity}
 D_{t,I}\mathfrak a_I=\alpha_I(t/\tau_a)a_I+a_I^{\rm c}.
\end{equation}
Thus $a_I^{\rm c}$ is the error in the antiderivative identity.
Differentiation does not enlarge the amplitude support, so the sum
remains supported in the same local contribution. The estimates for
$\mathfrak a_I$ and $a_I^{\rm c}$ are proved in
\eqref{principal:primitive-sharp}. The smaller gain at the highest
magnetic derivative orders is compensated by
\eqref{aniso:finite-depth-budget} when using the derivative scales
of the next iterate.

Define the principal potential before the corrector pushforward and its flow by
\begin{equation}
 \begin{gathered}
 \Theta_{I,0}^{\rm p}
 =\frac{1}{\lambda_{q+1}}\mathfrak a_I
       \varphi(\lambda_{q+1}y_I^{k_I}) d y_I^{\nu_I},
 \qquad \xi_0^{\rm p}=\sum_I\curl\Theta_{I,0}^{\rm p},\\
 \partial_sX_s^{\rm p}=\xi_0^{\rm p}\circ X_s^{\rm p},\qquad
 X_0^{\rm p}=\IId.
 \end{gathered}
\end{equation}
The complete principal potentials extend smoothly by zero across
the cut off regions of their charts. The principal LDF is independent of $s$
and may depend on physical time. Write
$\mathcal U_s^{\rm p}=(X_s^{\rm p})_*$ for its tensor pushforward.

We now transport each complete potential and tensor by
$\mathcal U_s^{\rm c}$. Composing the principal flow with the corrector
flow gives the full deformation:
\begin{equation}\label{path:full-generator}
 \begin{gathered}
 \Theta_{I,s}^{\rm p}=\mathcal U_s^{\rm c}\Theta_{I,0}^{\rm p}
 =\lambda_{q+1}^{-1}\mathcal U_s^{\rm c}
  \bigl(\mathfrak a_I\varphi(\lambda_{q+1}y_I^{k_I})d y_I^{\nu_I}\bigr),
 \qquad \xi_{I,s}^{\rm p}=\curl\Theta_{I,s}^{\rm p},\\
 X_s=X_s^{\rm c}\circ X_s^{\rm p},\qquad
 \xi_s=\xi^{\rm c}+\xi_s^{\rm p},\qquad
 \xi_s^{\rm p}=\sum_I\xi_{I,s}^{\rm p}=\mathcal U_s^{\rm c}\xi_0^{\rm p},
 \qquad \partial_sX_s=\xi_s\circ X_s.
 \end{gathered}
\end{equation}

The LDF and propagator of the composed flow are given by
Lemma~\ref{direct:flatness}. By \eqref{direct:differentiated-transport},
the material and magnetic Lie derivatives intertwine with
$\mathcal U_s^{\rm c}$.

For any tensor $T$ independent of $s$, these identities read
\begin{equation}\label{direct:stationary-data}
 \begin{aligned}
 (\partial_s+\mathcal L_{\xi^{\rm c}})\mathcal U_s^{\rm c}T&=0,\\
 \mathcal D_{t,I,s}^{\rm c}\mathcal U_s^{\rm c}T
   &=\mathcal U_s^{\rm c}\mathcal D_{t,I}T,\\
 \mathcal L_{B_{\ell,I,s}^{\rm c}}\mathcal U_s^{\rm c}T
   &=\mathcal U_s^{\rm c}\mathcal L_{B_{\ell,I}}T.
 \end{aligned}
\end{equation}
It follows from \eqref{direct:stationary-data} that the transported
phase, coframe and dual vector have zero material and magnetic Lie
derivatives with respect to the corrected background. The same
identity applies to $a_I$, whose material derivative need not vanish.

Only in endpoint formulas do we use the coordinate shorthand
\[
 y_{I,1}^{\eta}=\mathcal U_1^{\rm c}y_I^{\eta},\qquad
 d y_{I,1}^{\eta}=\mathcal U_1^{\rm c}d y_I^{\eta},\qquad
 \frac{\partial}{\partial y_{I,1}^{\eta}}=\mathcal U_1^{\rm c}\frac{\partial}{\partial y_I^{\eta}}.
\]
These are the coordinates and dual frames of the volume-preserving map
$\Psi_I\circ(X_1^{\rm c})^{-1}$. In
particular,
\[
 \frac{\partial}{\partial y_{I,1}^{\zeta_I}}
 =\nabla y_{I,1}^{k_I}\times d y_{I,1}^{\nu_I}.
\]
The coordinates $y_I^\eta$ are defined before the corrector pushforward; the endpoint coordinates
$y_{I,1}^\eta$ include the common corrector pushforward.

Since every LDF in the construction is a curl, the full path preserves
volume. Consequently, the actual fields
\begin{equation}\label{path:full-fields}
 v_s=(\partial_tX_s)\circ X_s^{-1}+(X_s)_*v_q,
 \qquad B_s=(X_s)_*B_q
\end{equation}
satisfy incompressibility and induction. To determine the change in
the momentum equation, we express the principal increments relative
to the local backgrounds and the corrected fields and apply the path
identities.

The time intervals \eqref{prep:slow-windows} and the profile assignment
preceding Lemma~\ref{amplitude:profiles} give the following support
property.

\begin{lemma}[Separation and bounded overlap]\label{mom:separation}
The corrector time intervals have uniformly bounded overlap. The
complete principal supports are disjoint and remain so after the
common pushforward $\mathcal U_s^{\rm c}$, whose LDF is the sum of the
corrector LDFs. Corrector time intervals may intersect a principal
support. The full flow does not transport a principal increment
across a transported cut off region.
\end{lemma}
\begin{proof}
For every $t$,
\[
 \#\{n:|t-t_n|\le2\tau_c\}\le5,
\]
since $t_{n+1}-t_n=\tau_c$. To prove disjointness of the principal
supports, consider first distinct profile types: their periodic
supports have disjoint neighborhoods. Contributions of the same type
have separated complete supports by construction. In both cases,
the approximate antiderivative preserves disjointness, since
differentiation does not enlarge the amplitude support and the
periodic primitives vanish in the same cut off region as the profiles.
Thus distinct local contributions satisfy
\[
 \supp\Theta_{I,0}^{\rm p}\cap\supp\Theta_{J,0}^{\rm p}=\varnothing,
 \qquad
 \supp\Theta_{I,s}^{\rm p}
       =X_s^{\rm c}(\supp\Theta_{I,0}^{\rm p}).
\]
By the second identity, all supports are mapped by the same diffeomorphism,
so their images remain disjoint. Smooth extension by zero also gives
\[
 \supp\partial^\beta\Theta_{I,0}^{\rm p}
 \subseteq\supp\Theta_{I,0}^{\rm p}
\]
for every derivative. In the corrector coordinates, the full flow is
$X_s^{\rm p}$. Its LDF vanishes in the cut off regions, so points
there are fixed by this flow. Uniqueness for the characteristic
equation therefore prevents a trajectory from crossing one of these regions.
\end{proof}

Fix $I$. Until the endpoint formulas below, we work on its local
support in the corresponding background chart; material and magnetic
derivatives before the corrector pushforward refer to this background.
The operators $\mathcal U_s^{\rm p}$ and $\mathcal U_s^{\rm c}$ still
denote the common pushforwards.

\paragraph{\textbf{Invariant potential.}}

On this support, volume preservation gives
$\mathrm{vol}=d y_I^{k_I}\wedge d y_I^{\nu_I}\wedge d y_I^{\zeta_I}$.
Separation implies that all other principal LDFs vanish on a
neighborhood of this support. The common principal flow therefore agrees there with the flow of
the local LDF. At fixed physical time, the latter is
\[
 f_I(t,y_I)=\lambda_{q+1}^{-1}\mathfrak a_I(t,y_I)
                       \varphi(\lambda_{q+1}y_I^{k_I}),\qquad
 \xi_{I,0}^{\rm p}=\curl(f_I\,dy_I^{\nu_I}).
\]
In local formulas we abbreviate $\varphi^{(j)}(\lambda_{q+1}y_I^{k_I})$
by $\varphi^{(j)}$ and $\alpha_I^{[j]}(t/\tau_a)$ by
$\alpha_I^{[j]}$, with the same convention for $\alpha_I$.
The scalar $\mathfrak a_I$ includes the complete localized amplitude.

Lemma~\ref{phase:prepared-identities} proves that the principal flow
preserves $f_I$, $dy_I^{\nu_I}$ and $df_I$, and that its tensor pushforward on
$g\,dy_I^{\nu_I}$ reduces to scalar transport of $g$. In particular,
\[
 \mathcal U_r^{\rm p}\curl(g\,dy_I^{\nu_I})
   =\curl((\mathcal U_r^{\rm p}g)\,dy_I^{\nu_I}).
\]
Although the potential is invariant, the oscillatory coordinate varies
according to \eqref{path:intrinsic-phase-drift}. We estimate this
variation, the frame and the path integrals in
Subsection~\ref{ssec:path-estimates}, Corollary~\ref{principal:phase-flow}.

By \eqref{path:factorized-transport}, the corresponding identities
on the full path are the corrector pushforwards of these scalar transport identities.

\paragraph{\textbf{Increments relative to the local background.}}

The first Lie derivatives of the potential on this support, before the
corrector pushforward, are
\begin{equation}
 \begin{aligned}
 \mathcal D_{t,I}\Theta_{I,0}^{\rm p}
   &=\lambda_{q+1}^{-1}(\alpha_I a_I+a_I^{\rm c})\varphi\,dy_I^{\nu_I},\\
 \mathcal L_{B_{\ell,I}}\Theta_{I,0}^{\rm p}
   &=\lambda_{q+1}^{-1}(D_{B,I}\mathfrak a_I)\varphi\,dy_I^{\nu_I},\\
 F_{I,0}^\pm&=(\mathcal D_{t,I}\pm\mathcal L_{B_{\ell,I}})\Theta_{I,0}^{\rm p}.
 \end{aligned}
\end{equation}
These formulas follow from \eqref{principal:primitive-identity},
since the local transports annihilate the phase and coframe. The
material Lie derivative incorporates both affine velocity terms in a
single potential, as made explicit in
\eqref{principal:velocity-duhamel}. Differentiating once more and using
commutation of the local transports gives
\begin{equation}
 \begin{aligned}
 (\mathcal D_{t,I}^2-\mathcal L_{B_{\ell,I}}^2)\Theta_{I,0}^{\rm p}
 &=\lambda_{q+1}^{-1}(D_{t,I}^2\mathfrak a_I-D_{B,I}^2\mathfrak a_I)
          \varphi(\lambda_{q+1}y_I^{k_I})d y_I^{\nu_I}\\
 &=\lambda_{q+1}^{-1}
  \bigl[\tau_a^{-1}\alpha_I' a_I+\alpha_I D_{t,I}a_I+D_{t,I}a_I^{\rm c}-D_{B,I}^2\mathfrak a_I\bigr]
          \varphi(\lambda_{q+1}y_I^{k_I})d y_I^{\nu_I}.
 \end{aligned}
\end{equation}
The second equality follows by differentiating the same antiderivative identity.
Every curl of these periodic potentials is divergence free and has
zero spatial mean. Define
\begin{equation}\label{path:prepared-increments}
 \begin{aligned}
 \Theta_{I,s}^{{\rm p},\pm}
 &=\mathcal U_s^{\rm c}\int_0^s\mathcal U_r^{\rm p}F_{I,0}^\pm\,\dd r,\\
 \Theta_{w,I,s}^{\rm p}&=\tfrac12(\Theta_{I,s}^{{\rm p},+}+\Theta_{I,s}^{{\rm p},-}),\qquad
 \Theta_{b,I,s}^{\rm p}=\tfrac12(\Theta_{I,s}^{{\rm p},+}-\Theta_{I,s}^{{\rm p},-}).
 \end{aligned}
\end{equation}
Before the corrector pushforward, this integral solves the transport
equation with autonomous LDF, source $F_{I,0}^\pm$ and zero initial value.

\begin{lemma}[Principal field increments]
Let the local background fields transported by the full map be
\begin{equation}
 v_{\ell,I,s}=(\partial_tX_s)\circ X_s^{-1}+(X_s)_*v_{\ell,I},
 \qquad B_{\ell,I,s}=(X_s)_*B_{\ell,I}.
\end{equation}
Set $w_{I,s}=v_{\ell,I,s}-v_{\ell,I,s}^{\rm c}$ and
$b_{I,s}=B_{\ell,I,s}-B_{\ell,I,s}^{\rm c}$. On the image of the full
local support under $X_s^{\rm c}$, we have
\begin{equation}\label{principal:exact-local-increments}
 w_{I,s}\pm b_{I,s}=\curl\Theta_{I,s}^{{\rm p},\pm},\qquad
 w_{I,s}=\curl\Theta_{w,I,s}^{\rm p},\qquad b_{I,s}=\curl\Theta_{b,I,s}^{\rm p}.
\end{equation}
At the endpoint,
\begin{equation}\label{path:direct-path-average}
 \Theta_{I,1}^{{\rm p},\pm}=\mathcal U_1^{\rm c}\mathcal P^{\rm p}F_{I,0}^\pm,
 \qquad \mathcal P^{\rm p}=\int_0^1\mathcal U_r^{\rm p}\,\dd r,
 \qquad \mathcal U_r^{\rm p}=e^{-r\mathcal L_{\xi_0^{\rm p}}}.
\end{equation}
\end{lemma}
\begin{proof}
Apply \eqref{path:factorized-field-increments} to the local background
pair. On this support the principal sources are
$\curl\mathcal D_{t,I}\Theta_{I,0}^{\rm p}$ and
$\curl\mathcal L_{B_{\ell,I}}\Theta_{I,0}^{\rm p}$.
The pushforwards commute with curl because the flows preserve volume.
Hence
\[
 w_{I,s}\pm b_{I,s}
 =\mathcal U_s^{\rm c}\curl\int_0^s
             \mathcal U_r^{\rm p}F_{I,0}^\pm\,\dd r
 =\curl\Theta_{I,s}^{{\rm p},\pm}.
\]
Taking the half-sum and half-difference proves
\eqref{principal:exact-local-increments}; setting $s=1$ gives
\eqref{path:direct-path-average}. This use of the path lemma allows
the full LDF to depend on $s$. Only the principal LDF before the
corrector pushforward is independent of $s$.
\end{proof}

The exact scalar representations are now
\begin{equation}
 \begin{aligned}
 \Theta_{w,I,s}^{\rm p}&=\lambda_{q+1}^{-1}\mathcal U_s^{\rm c}
       \int_0^s\mathcal U_r^{\rm p}
           ((\alpha_I a_I+a_I^{\rm c})\varphi\,dy_I^{\nu_I})\,\dd r,\\
 \Theta_{b,I,s}^{\rm p}&=\lambda_{q+1}^{-1}\mathcal U_s^{\rm c}
       \int_0^s\mathcal U_r^{\rm p}
           ((D_{B,I}\mathfrak a_I)\varphi\,dy_I^{\nu_I})\,\dd r.
 \end{aligned}
\end{equation}
The sources include all derivatives of the outer cutoff in the
approximate antiderivative. Although $\alpha_I a_I\varphi$ need not
be invariant, these formulas allow us to estimate its variation at
zeros of $\mathfrak a_I$ and the temporal profiles, since no division
by these factors is involved.

Formula \eqref{principal:scalar-action} for the directional derivative
$\xi_{I,0}^{\rm p}\cn$ and its normalized coefficient bounds are proved in
Subsection~\ref{ssec:path-estimates}.
Writing
\[
 (\mathcal U_s^{\rm c})^{-1}\Theta_{w,I,s}^{\rm p}=f_{w,I,s}dy_I^{\nu_I},
\]
and similarly for the magnetic potential, gives the physical equations
\begin{equation}\label{principal:profile-equation}
 \begin{aligned}
 (\partial_s+\xi_{I,0}^{\rm p}\cn)f_{w,I,s}
    &=\lambda_{q+1}^{-1}(\alpha_I a_I+a_I^{\rm c})\varphi,
 &f_{w,I,0}&=0,\\
 (\partial_s+\xi_{I,0}^{\rm p}\cn)f_{b,I,s}
    &=\lambda_{q+1}^{-1}(D_{B,I}\mathfrak a_I)\varphi,
 &f_{b,I,0}&=0.
 \end{aligned}
\end{equation}
Indeed, for an LDF independent of $s$,
\[
 (\partial_s+\mathcal L_{\xi_0^{\rm p}})
 \int_0^s\mathcal U_r^{\rm p}F\,\dd r=F.
\]
The normalized coefficient bounds of
Subsection~\ref{ssec:path-estimates} show that the flow differential
is close to the identity in anisotropic coordinates. We estimate the
potentials in these coordinates. In Euclidean coordinates, the flow
differential may contain the large factor
$\tau_a\lambda_{q+1}\delta_{q+1}^{1/2}$.

\paragraph{\textbf{Principal gap increments.}}

Fix an index $I=(n,J,\zeta)$ and denote the corresponding initial
background differences by
\[
 G_I^v=v_q-v_{\ell,I},\qquad G_I^B=B_q-B_{\ell,I},\qquad G_I^\pm=G_I^v\pm G_I^B.
\]
On the support of this local contribution, the common pushforward
$\mathcal U_s^{\rm p}$ agrees with pushforward by the flow of the local LDF
$\xi_{I,0}^{\rm p}$ before the corrector pushforward. The factorization \eqref{path:full-generator}
and the affine velocity action give the three comparisons
\begin{align}
 G_{I,s}^{\rm c,\pm}
  &:=v_s^{\rm c}\pm B_s^{\rm c}
       -(v_{\ell,I,s}^{\rm c}\pm B_{\ell,I,s}^{\rm c})
    =\mathcal U_s^{\rm c}G_I^\pm,
                                      \\
 G_{I,s}^\pm&=\mathcal U_s^{\rm c}\mathcal U_s^{\rm p}G_I^\pm,
                                      \label{adapt:full-block-gap}\\
 G_{I,s}^\pm-G_{I,s}^{\rm c,\pm}
   &=\mathcal U_s^{\rm c}(\mathcal U_s^{\rm p}-\IId)G_I^\pm
     =\curl\Theta_{I,s}^{\mathrm p,\mathrm{gap},\pm}.
                                      \label{adapt:gap-curl}
\end{align}
Here $G_{I,s}^\pm$ is the signed difference between the actual fields
and the local fields transported by the full path. Subtracting the
corrector pushforward of the initial gap gives the additional
principal increment in the last line. The identities hold on the
full transported local support, where every other principal LDF
vanishes, and refer to background $I$ even when the slow time
intervals overlap.

To prove the potential identity, work in the chart of the local background indexed by $I$.
With the scalar $f_I$ above, Lemma~\ref{phase:prepared-identities} gives
\[
 \Theta_{I,0}^{\rm p}=f_I d y_I^{\nu_I}.
\]
The invariants in \eqref{phase:prepared-invariance} reduce the gap
increment to scalar equations as well. The relevant contraction is
\[
 \iota_{G_I^\pm} d\Theta_{I,0}^{\rm p}
 =(G_I^\pm\cn f_I)\, d y_I^{\nu_I}-(G_I^\pm\cn y_I^{\nu_I})\, d f_I.
\]
Consequently the exact gap potential is
\begin{equation}\label{adapt:prepared-gap-potential}
 \begin{split}
 \Theta_{I,s}^{\mathrm p,\mathrm{gap},\pm}
 &=\mathcal U_s^{\rm c}\int_0^s\mathcal U_r^{\rm p}
                 (\iota_{G_I^\pm} d\Theta_{I,0}^{\rm p})\,\dd r\\
 &=\mathcal U_s^{\rm c}\left[
      \left(\int_0^s\mathcal U_r^{\rm p}(G_I^\pm\cn f_I)\,\dd r\right) d y_I^{\nu_I}
     -\left(\int_0^s\mathcal U_r^{\rm p}(G_I^\pm\cn y_I^{\nu_I})\,\dd r\right) d f_I
                              \right].
 \end{split}
\end{equation}
The factors $dy_I^{\nu_I}$ and $df_I$ are invariant under the principal
pushforward, which justifies taking them outside the integrals. Both
remaining integrals involve only scalar transport. Cartan's identity and
volume preservation give
\[
 \begin{aligned}
 \curl(\iota_{G_I^\pm} d\Theta_{I,0}^{\rm p})
   &=\mathcal L_{G_I^\pm}\xi_{I,0}^{\rm p}
     =-\mathcal L_{\xi_{I,0}^{\rm p}}G_I^\pm,\\
 \curl\Theta_{I,s}^{\mathrm p,\mathrm{gap},\pm}
   &=-\mathcal U_s^{\rm c}\int_0^s
                 \mathcal U_r^{\rm p}\mathcal L_{\xi_{I,0}^{\rm p}}G_I^\pm\,\dd r
     =\mathcal U_s^{\rm c}(\mathcal U_s^{\rm p}-\IId)G_I^\pm.
 \end{aligned}
\]
This proves \eqref{adapt:gap-curl}.

For the global definition, solve the scalar equation with source
$G_I^\pm\cn y_I^{\nu_I}$ on a neighborhood of the local support in
the chart, and extend its product with $d f_I$ by zero. The product
vanishes in the cut off region, although the scalar factor need not
be compactly supported. The other source, $G_I^\pm\cn f_I$, is
already supported in $\supp d f_I$. Since $f_I$ and $y_I^{\nu_I}$
are invariant, the principal flow preserves the supports of both
complete products. Their images under the common corrector flow
therefore satisfy Lemma~\ref{mom:separation}. No division by $a_I$
or $\mathfrak a_I$ is used, so the formulas also hold at their zeros.

\paragraph{\textbf{Endpoint fields.}}

For the potentials relative to the local backgrounds in \eqref{path:prepared-increments}, the
additional principal increments are
\[
 w_{I,s}^{\mathrm p,\mathrm{gap}}\pm b_{I,s}^{\mathrm p,\mathrm{gap}}
 =G_{I,s}^\pm-G_{I,s}^{\rm c,\pm}
 =\curl\Theta_{I,s}^{\mathrm p,\mathrm{gap},\pm}.
\]
The complete principal increments for each local contribution, extended by zero outside its
transported support, are
\[
 w_I^{\rm p}=\curl\Theta_{w,I,1}^{\rm p}+w_{I,1}^{\mathrm p,\mathrm{gap}},
 \qquad
 b_I^{\rm p}=\curl\Theta_{b,I,1}^{\rm p}+b_{I,1}^{\mathrm p,\mathrm{gap}}.
\]
Together with the corrector increments, their sums give the actual endpoint:
\begin{equation}\label{direct:actual-endpoint}
 \begin{aligned}
 w^{\rm p}&=\sum_Iw_I^{\rm p},
 &b^{\rm p}&=\sum_Ib_I^{\rm p},\\
 v_1&=v_q+w^{\rm c}+w^{\rm p},
 &B_1&=B_q+b^{\rm c}+b^{\rm p}.
 \end{aligned}
\end{equation}
Indeed, on the full support of each local principal contribution, the affine factorization gives
\[
 \begin{aligned}
 (v_1-v_1^{\rm c})\pm(B_1-B_1^{\rm c})
 &=(v_{\ell,I,1}-v_{\ell,I,1}^{\rm c})
       \pm(B_{\ell,I,1}-B_{\ell,I,1}^{\rm c})
       +G_{I,1}^\pm-G_{I,1}^{\rm c,\pm}\\
 &=\curl\Theta_{I,1}^{\mathrm p,\pm}
       +\curl\Theta_{I,1}^{\mathrm p,\mathrm{gap},\pm}
   =w_I^{\rm p}\pm b_I^{\rm p}.
 \end{aligned}
\]
The local field identity \eqref{principal:exact-local-increments} and
the gap identity \eqref{adapt:gap-curl} give the last line. Extending
the complete curl increments by zero and adding the corrector now
gives the global decomposition at the level of potentials:
\begin{equation}\label{direct:full-field-split}
 \begin{aligned}
 (v_1-v_q)\pm(B_1-B_q)
 &=w^{\rm c}\pm b^{\rm c}+\sum_I(w_I^{\rm p}\pm b_I^{\rm p})\\
 &=\sum_n\bigl(w_n^{\mathrm c,\mathrm{main}}
                  \pm b_n^{\mathrm c,\mathrm{main}}\bigr)\\
 &\quad+\curl\Bigl[\sum_n\Theta_n^{\mathrm c,\mathrm{rem},\pm}
       +\sum_I\bigl(\Theta_{I,1}^{\mathrm p,\pm}
                    +\Theta_{I,1}^{\mathrm p,\mathrm{gap},\pm}\bigr)\Bigr].
 \end{aligned}
\end{equation}
The half-sum and half-difference give the velocity and magnetic
increments. In the quadratic expansion, the principal sum has disjoint
supports, whereas the corrector terms may overlap in time. The initial
gaps are not mollified. In Section~\ref{sec:generator-field-estimates},
Proposition~\ref{gap:field-bounds} estimates them using
\eqref{prep:local-gap-good} and the ordinary bounds
\eqref{prep:old-ordinary}.

\subsection{Quadratic Terms}

\paragraph{\textbf{Quadratic decomposition.}}

The field decomposition \eqref{direct:full-field-split} gives
\begin{equation}\label{direct:full-quadratic-split}
 \begin{aligned}
 &(v_1-v_q)\otimes(v_1-v_q)-(B_1-B_q)\otimes(B_1-B_q)\\
 &\quad=w^{\rm c}\otimes w^{\rm c}-b^{\rm c}\otimes b^{\rm c}
       +2\bigl(w^{\rm c}\odot w^{\rm p}-b^{\rm c}\odot b^{\rm p}\bigr)\\
 &\qquad\quad+\sum_I
       \bigl(w_I^{\rm p}\otimes w_I^{\rm p}-b_I^{\rm p}\otimes b_I^{\rm p}\bigr).
 \end{aligned}
\end{equation}
Disjointness eliminates products of distinct principal contributions.
The corrector quadratic term still contains products from overlapping
time intervals, and the cross term is included in the force linearized
about $(v_1^{\rm c},B_1^{\rm c})$ in \eqref{direct:cross-grouping}.
For each principal contribution, expanding the gap terms gives
\begin{equation}\label{principal:gap-square}
 \begin{aligned}
 &w_I^{\rm p}\otimes w_I^{\rm p}-b_I^{\rm p}\otimes b_I^{\rm p}\\
 &\quad=\curl\bigl(\Theta_{I,1}^{\mathrm p,+}
                     +\Theta_{I,1}^{\mathrm p,\mathrm{gap},+}\bigr)
         \odot\curl\bigl(\Theta_{I,1}^{\mathrm p,-}
                     +\Theta_{I,1}^{\mathrm p,\mathrm{gap},-}\bigr)\\
 &\quad=\curl\Theta_{w,I,1}^{\rm p}\otimes\curl\Theta_{w,I,1}^{\rm p}
       -\curl\Theta_{b,I,1}^{\rm p}\otimes\curl\Theta_{b,I,1}^{\rm p}\\
 &\qquad+2\bigl(\curl\Theta_{w,I,1}^{\rm p}\odot w_{I,1}^{\mathrm p,\mathrm{gap}}
             -\curl\Theta_{b,I,1}^{\rm p}\odot b_{I,1}^{\mathrm p,\mathrm{gap}}\bigr)\\
 &\qquad+w_{I,1}^{\mathrm p,\mathrm{gap}}\otimes w_{I,1}^{\mathrm p,\mathrm{gap}}
          -b_{I,1}^{\mathrm p,\mathrm{gap}}\otimes b_{I,1}^{\mathrm p,\mathrm{gap}}.
 \end{aligned}
\end{equation}
The last two lines are the quadratic background corrections. We next
express the preceding quadratic term using the path average.

\paragraph{\textbf{Difference of squares.}}

Fix an index $I$ and use the corresponding endpoint increments relative to the local background
$w_{I,1},b_{I,1}$ from \eqref{principal:exact-local-increments}. Put
\[
 \begin{gathered}
 \mathbb Q_{I,0}=\curl F_{I,0}^+\odot\curl F_{I,0}^-,\qquad
 \mathbb Q_{I,1}=\mathcal U_1^{\rm c}\mathbb Q_{I,0},\\
 \mathbb Q_{I,2}=\mathcal U_1^{\rm c}\Bigl(\frac12\int_0^1(1-r)^2
       \mathcal U_r^{\rm p}\mathcal L_{\xi_{I,0}^{\rm p}}^2\mathbb Q_{I,0}\,\dd r
       -\Cov_{X^{\rm p}}(\curl F_{I,0}^+,\curl F_{I,0}^-)\Bigr),
 \end{gathered}
\]
where the last tensor collects the second variation and covariance
from \eqref{eq:quadratic-remainder}, followed by the corrector
pushforward. The map $\mathcal U_1^{\rm c}$ preserves symmetric
products and satisfies
$\mathcal U_1^{\rm c}\mathcal L_{\xi_{I,0}^{\rm p}}
=\mathcal L_{\xi_{I,1}^{\rm p}}\mathcal U_1^{\rm c}$.
Thus the path-product identity and the finite expansion
\eqref{eq:quadratic-remainder}, applied to
\eqref{path:direct-path-average}, express the endpoint square as
\begin{equation}\label{principal:exact-square}
 \begin{aligned}
 w_{I,1}\otimes w_{I,1}-b_{I,1}\otimes b_{I,1}
 &=(w_{I,1}+b_{I,1})\odot(w_{I,1}-b_{I,1})\\
 &=\mathcal U_1^{\rm c}\bigl[
       (\mathcal P^{\rm p}\curl F_{I,0}^+)\odot
       (\mathcal P^{\rm p}\curl F_{I,0}^-)\bigr]\\
 &=\mathcal U_1^{\rm c}\bigl[\mathcal P^{\rm p}\mathbb Q_{I,0}
                -\Cov_{X^{\rm p}}(\curl F_{I,0}^+,\curl F_{I,0}^-)\bigr]\\
 &=\mathbb Q_{I,1}-\frac12\mathcal L_{\xi_{I,1}^{\rm p}}\mathbb Q_{I,1}+\mathbb Q_{I,2}.
 \end{aligned}
\end{equation}

The first Lie variation of the quadratic tensor has the exact finite Fourier expansion in the fast phase variable
\begin{equation}\label{principal:first-lie-mean}
 \mathcal L_{\xi_{I,1}^{\rm p}}\mathbb Q_{I,1}
 =\sum_{h\in\{-3,-1,1,3\}} C_{I,h}(t,x)
          e^{ih\lambda_{q+1}y_{I,1}^{k_I}}.
\end{equation}
The coefficients $C_{I,h}$ are sums of products of transported
amplitudes and chart factors. To check the Fourier indices, note that $\varphi$ and its
derivatives have only the modes $\pm1$, as do $\curl F_{I,0}^\pm$
and $\xi_{I,0}^{\rm p}$. Thus $\mathbb Q_{I,0}$ has modes $0,\pm2$.
Spatial differentiation, whether of a coefficient or an exponential,
preserves these indices. Each term in
$\mathcal L_{\xi_{I,0}^{\rm p}}\mathbb Q_{I,0}$ consequently has
index $\pm1$ or $\pm3$. Composing the coefficients and phase with the
corrector map proves \eqref{principal:first-lie-mean}. This eliminates
the constant coefficient in the fast Fourier expansion, but does not
imply that the spatial integral of the tensor vanishes.

\paragraph{\textbf{The principal quadratic tensor.}}

For the same index $I$ at $s=1$, write its transported coefficients as
\[
 a_{I,1}=\mathcal U_1^{\rm c}a_I,\qquad
 a_{I,1}^{\rm c}=\mathcal U_1^{\rm c}a_I^{\rm c},\qquad
 \mathfrak a_{I,1}=\mathcal U_1^{\rm c}\mathfrak a_I.
\]
Use the corrected coordinates $y_{I,1}^\eta$ and the magnetic derivative
$D_{B_{\ell,I,1}^{\rm c}}$. Intertwining gives
\[
 D_{B_{\ell,I,1}^{\rm c}}\mathfrak a_{I,1}=\mathcal U_1^{\rm c}(D_{B,I}\mathfrak a_I).
\]
The profile arguments remain
$\alpha_I=\alpha_I(t/\tau_a)$ and
$\varphi^{(j)}=\varphi^{(j)}(\lambda_{q+1}y_{I,1}^{k_I})$.
Apply the curl product rule to the material and magnetic Lie
derivatives of the principal potential. The derivatives
$\mathcal D_{t,I}$ and $\mathcal L_{B_{\ell,I}}$ inside
$\mathcal U_1^{\rm c}$ are those of the local background, so
\begin{align}
 \curl(\mathcal U_1^{\rm c}\mathcal D_{t,I}\Theta_{I,0}^{\rm p})
 &= (\alpha_I a_{I,1}+a_{I,1}^{\rm c})\varphi'\frac{\partial}{\partial y_{I,1}^{\zeta_I}}
     +\lambda_{q+1}^{-1}\varphi\,
           \nabla(\alpha_I a_{I,1}+a_{I,1}^{\rm c})\times dy_{I,1}^{\nu_I},
                       \label{osc:velocity-split}\\
 \curl(\mathcal U_1^{\rm c}\mathcal L_{B_{\ell,I}}\Theta_{I,0}^{\rm p})
 &= (D_{B_{\ell,I,1}^{\rm c}}\mathfrak a_{I,1})\varphi'\frac{\partial}{\partial y_{I,1}^{\zeta_I}}
     +\lambda_{q+1}^{-1}\varphi\,
           \nabla(D_{B_{\ell,I,1}^{\rm c}}\mathfrak a_{I,1})\times dy_{I,1}^{\nu_I}.
                       \label{osc:magnetic-split}
\end{align}
Before averaging along the principal path, the difference of these tensor
squares is $\mathbb Q_{I,1}$. Expanding it gives the prescribed leading
term, the antiderivative correction, the magnetic quadratic term and the
terms in which curl differentiates a coefficient:
\begin{equation}\label{osc:complete-square}
 \begin{aligned}
 \mathbb Q_{I,1}
 ={}&\alpha_I^2a_{I,1}^2\bigl(1+\cos(2\lambda_{q+1}y_{I,1}^{k_I})\bigr)
                       \partial_{y_{I,1}^{\zeta_I}}\otimes\partial_{y_{I,1}^{\zeta_I}}\\
 &+\bigl[2\alpha_I a_{I,1}a_{I,1}^{\rm c}+(a_{I,1}^{\rm c})^2-(D_{B_{\ell,I,1}^{\rm c}}\mathfrak a_{I,1})^2\bigr]
                  (\varphi')^2\partial_{y_{I,1}^{\zeta_I}}\otimes\partial_{y_{I,1}^{\zeta_I}}\\
 &+2\lambda_{q+1}^{-1}\varphi\varphi'\partial_{y_{I,1}^{\zeta_I}}\odot
    \Bigl[(\alpha_I a_{I,1}+a_{I,1}^{\rm c})
            \bigl(\nabla(\alpha_I a_{I,1}+a_{I,1}^{\rm c})\times dy_{I,1}^{\nu_I}\bigr)\\
 &\hspace{49mm}-(D_{B_{\ell,I,1}^{\rm c}}\mathfrak a_{I,1})
                    \bigl(\nabla D_{B_{\ell,I,1}^{\rm c}}\mathfrak a_{I,1}\times dy_{I,1}^{\nu_I}\bigr)\Bigr]\\
 &+\lambda_{q+1}^{-2}\varphi^2\Bigl[
       \bigl(\nabla(\alpha_I a_{I,1}+a_{I,1}^{\rm c})\times dy_{I,1}^{\nu_I}\bigr)
         \otimes\bigl(\nabla(\alpha_I a_{I,1}+a_{I,1}^{\rm c})\times dy_{I,1}^{\nu_I}\bigr)\\
 &\hspace{29mm}-\bigl(\nabla D_{B_{\ell,I,1}^{\rm c}}\mathfrak a_{I,1}\times dy_{I,1}^{\nu_I}\bigr)
                    \otimes\bigl(\nabla D_{B_{\ell,I,1}^{\rm c}}\mathfrak a_{I,1}\times dy_{I,1}^{\nu_I}\bigr)\Bigr].
 \end{aligned}
\end{equation}
Here $\nabla$ differentiates the displayed coefficient.
The products $\varphi\varphi'$, $(\varphi')^2$ and $\varphi^2$
have constant Fourier coefficients $0$, $1$ and $1$ in the fast phase variable,
respectively. Thus the constant coefficient in that fast Fourier expansion of $\mathbb Q_{I,1}$ is
\begin{equation}\label{principal:complete-principal-average}
 \begin{split}
 &\big[(\alpha_I a_{I,1}+a_{I,1}^{\rm c})^2-(D_{B_{\ell,I,1}^{\rm c}}\mathfrak a_{I,1})^2\big]
                     (\partial_{y_{I,1}^{\zeta_I}}\otimes\partial_{y_{I,1}^{\zeta_I}})\\
 &+\lambda_{q+1}^{-2}\big[
    (\nabla(\alpha_I a_{I,1}+a_{I,1}^{\rm c})\times dy_{I,1}^{\nu_I})
                \otimes(\nabla(\alpha_I a_{I,1}+a_{I,1}^{\rm c})\times dy_{I,1}^{\nu_I})\\
 &\hspace{42mm}-(\nabla D_{B_{\ell,I,1}^{\rm c}}\mathfrak a_{I,1}\times dy_{I,1}^{\nu_I})
                \otimes(\nabla D_{B_{\ell,I,1}^{\rm c}}\mathfrak a_{I,1}\times dy_{I,1}^{\nu_I})\big].
 \end{split}
\end{equation}
The prescribed coefficient is $\alpha_I^2a_{I,1}^2$, and the corresponding
rank-one quadratic tensor in the direction $\frac{\partial}{\partial y_{I,1}^{\zeta_I}}$ is
\begin{equation}\label{osc:principal-square}
 \alpha_I^2a_{I,1}^2(1+\cos(2\lambda_{q+1}y_{I,1}^{k_I}))
                 (\partial_{y_{I,1}^{\zeta_I}}\otimes\partial_{y_{I,1}^{\zeta_I}}).
\end{equation}
It is convenient to define the error tensors using coefficients in
the original charts and then apply the corrector pushforward. With
$\varphi'=\varphi'(\lambda_{q+1}y_I^{k_I})$ in the $I$th summand,
the antiderivative and magnetic errors, including their oscillatory parts,
are
\begin{equation}\label{osc:primitive-magnetic-tensors}
 \begin{aligned}
 R^{\mathrm p}_{\rm time}
 &=\sum_I\mathcal U_1^{\rm c}
   \left[(2\alpha_Ia_Ia_I^{\rm c}+(a_I^{\rm c})^2)
       (\varphi')^2(\partial_{y_I^{\zeta_I}}\otimes\partial_{y_I^{\zeta_I}})\right],\\
 R^{\mathrm p}_{\rm mag}
 &=-\sum_I\mathcal U_1^{\rm c}
   \left[(D_{B,I}\mathfrak a_I)^2
       (\varphi')^2(\partial_{y_I^{\zeta_I}}\otimes\partial_{y_I^{\zeta_I}})\right].
 \end{aligned}
\end{equation}
The coefficients in the brackets are physical functions, and
$D_{t,I}^{j_0}a_I$ and $\mathfrak a_I$ include the derivatives of
$\varrho$. The two rank-one tensors therefore include their outer cutoff
contributions, while the curl products remain in
\eqref{principal:complete-principal-average}.

\subsection{Definition of the Stresses}

We start by defining the cutoff and mollification errors as 
\begin{equation}\label{mom:collar-stress}
 R^{\mathrm{cut}}=\sum_n\mathcal R\curl
                              [\mathscr G_n,\tilde\eta_n]\Theta_n
\end{equation}
and
\begin{equation}\label{mom:mollification-stress}
 \begin{aligned}
 R_n^{\rm moll}
   &=\varrho^2\eta_n^2(\mathcal J_n-\IId)
                         (\delta_{q+1}\IId-R_q),\\
 R^{\rm moll}&=\sum_nR_n^{\rm moll}.
 \end{aligned}
\end{equation}
The tensor $R^{\rm moll}$ accounts for mollification of
$\delta_{q+1}\IId-R_q$. The errors from regularizing the velocity and
magnetic fields are included in the terms defined below.

For the perturbations $(\sum_nw_n,\sum_nb_n)$, the linearized force is
\begin{align}
 \mathcal F_{v_q,B_q}\Bigl(\sum_nw_n,\sum_nb_n\Bigr)
 ={}&\sum_n\bigl((\partial_t+\mathcal L_{v_{\ell,n}})w_n
                                      -\mathcal L_{B_{\ell,n}}b_n\bigr)
       \notag\\
 &+2\sum_n(w_n\cn v_{\ell,n}-b_n\cn B_{\ell,n})
       \notag\\
 &+2\ddiv\sum_n\bigl(((v_q-v_{\ell,n})\odot w_n)
                         -((B_q-B_{\ell,n})\odot b_n)\bigr).
 \label{mom:linear-transfer}
\end{align}
Accordingly, we define the linear error terms, with the subscripts $\mathrm{tr}$, $\mathrm{na}$ and $\mathrm{moll}$ referring to the transport, Nash and mollification errors, by
\begin{align}
 R^{\mathrm c}_{\rm tr}&=\sum_n\mathcal R\bigl((\partial_t+\mathcal L_{v_{\ell,n}})w_n^{\mathrm c,\mathrm{rem}}
                                      -\mathcal L_{B_{\ell,n}}b_n^{\mathrm c,\mathrm{rem}}\bigr),
                                      \label{mom:corrector-transport-stress}\\
 R^{\mathrm c}_{\mathrm{na}}&=2\sum_n\mathcal R(w_n^{\mathrm c,\mathrm{rem}}\cn v_{\ell,n}
                                      -b_n^{\mathrm c,\mathrm{rem}}\cn B_{\ell,n}),
                                      \label{mom:corrector-interaction-stress}\\
 R^{\mathrm c}_{\rm moll}&=2\sum_n\bigl(((v_q-v_{\ell,n})\odot w_n)
                              -((B_q-B_{\ell,n})\odot b_n)\bigr),
                                      \label{mom:G-gap}\\
 R^{\mathrm p}_{\rm tr}&=\sum_I\mathcal R\bigl(
       (\partial_t+\mathcal L_{v_{\ell,I,1}^{\rm c}})w_I^{\rm p}
                   -\mathcal L_{B_{\ell,I,1}^{\rm c}}b_I^{\rm p}\bigr),
                                      \label{mom:principal-transport-stress}\\
 R^{\mathrm p}_{\mathrm{na}}&=2\sum_I\mathcal R\bigl(w_I^{\rm p}\cn v_{\ell,I,1}^{\rm c}
                         -b_I^{\rm p}\cn B_{\ell,I,1}^{\rm c}\bigr),
                                      \label{mom:principal-interaction-stress}\\
 R^{\mathrm p}_{\rm moll}&=2\sum_I\bigl(
       ((v_1^{\rm c}-v_{\ell,I,1}^{\rm c})\odot w_I^{\rm p})
       -((B_1^{\rm c}-B_{\ell,I,1}^{\rm c})\odot b_I^{\rm p})\bigr).
                                      \label{mom:principal-gap-stress}
\end{align}

Finally, the quadratic error term is formed from the full corrector increments:
\begin{equation}\label{mom:corrector-quadratic-stress}
 R^{\mathrm c}_{\rm quad}=w^{\rm c}\otimes w^{\rm c}-b^{\rm c}\otimes b^{\rm c}.
\end{equation}
Every field to which $\mathcal R$ is applied above has zero spatial mean:
it is a curl, a time derivative of a curl, a bracket of divergence-free
fields, or a divergence $r\cn z=\ddiv(z\otimes r)$.

For the leading principal square, sum over all selected indices:
\begin{equation}
 \mathsf T=\sum_I\alpha_I(t/\tau_a)^2a_I^2
       \partial_{y_I^{\zeta_I}}\otimes\partial_{y_I^{\zeta_I}},
\end{equation}
so that 
\begin{equation}\label{mom:target-balance}
 \sum_n\mathsf F_n+\mathsf T
 =\varrho^2(\delta_{q+1}\IId-R_q)+R^{\rm moll}.
\end{equation}
Since pushforward preserves tensor products, the leading quadratic tensor
at the endpoint is $\mathcal U_1^{\rm c}\mathsf T$, with the outer cutoff
included in the coefficients. We therefore define the chart stress by
\begin{equation}\label{mom:chart-stress}
 R^{\mathrm p}_{\rm chart}=(\mathcal U_1^{\rm c}-\IId)\mathsf T.
\end{equation}
This is an equality of tensors. Its estimate uses the Lie transport
equation, without commuting Euclidean divergence with
$\mathcal U_1^{\rm c}$.

We define $R^{\mathrm p}_{\rm quad}$ from the expansion
\eqref{principal:gap-square} on each full local support, using
\eqref{principal:exact-square} for the terms relative to the local
background and \eqref{osc:complete-square} for $\mathbb Q_{I,1}$.
For the cosine term in \eqref{osc:principal-square}, use
Lemma~\ref{osc:parametrix}. Apply the same lemma to each nonzero fast
Fourier mode of the component in the tangent direction of
$-\tfrac12\mathcal L_{\xi_{I,1}^{\rm p}}\mathbb Q_{I,1}$; its constant
mode vanishes by \eqref{principal:first-lie-mean}.

The remaining summands are the transverse components, the curl
products, $\mathbb Q_{I,2}$ and the last two lines of
\eqref{principal:gap-square}. The tensors in
\eqref{osc:primitive-magnetic-tensors} are excluded. Summing these
symmetric tensors over the local contributions defines the quadratic
stress. Wherever the construction requires a spatial inverse
divergence, we first sum the supported forces and then apply the
operator; the resulting tensor need not have the same support.

Each tensor constructed in this way has the same divergence as the
term it replaces. Thus
\begin{equation}\label{mom:N-remainder}
 \ddiv R^{\mathrm p}_{\rm quad}
 =\ddiv\bigl(w^{\rm p}\otimes w^{\rm p}
               -b^{\rm p}\otimes b^{\rm p}-\mathcal U_1^{\rm c}\mathsf T
               -R^{\mathrm p}_{\rm time}-R^{\mathrm p}_{\rm mag}\bigr).
\end{equation}
\begin{theorem}[Endpoint stress]\label{mom:complete-ledger}
The fields $(v_{q+1},B_{q+1})=(v_1,B_1)$ from
\eqref{direct:actual-endpoint}, with pressure
$p_{q+1}=p_q+\sum_n\pi_n$, satisfy the relaxed MHD equations with
symmetric stress
\begin{equation}\label{mom:full-stress}
 \begin{split}
 R_{q+1}={}&R^{\mathrm{cut}}+R^{\rm moll}
       +R^{\mathrm c}_{\rm tr}+R^{\mathrm c}_{\mathrm{na}}
       +R^{\mathrm c}_{\rm moll}+R^{\mathrm c}_{\rm quad}\\
       &+R^{\mathrm p}_{\rm tr}+R^{\mathrm p}_{\mathrm{na}}
       +R^{\mathrm p}_{\rm moll}+R^{\mathrm p}_{\rm chart}
       +R^{\mathrm p}_{\rm quad}+R^{\mathrm p}_{\rm time}+R^{\mathrm p}_{\rm mag}.
 \end{split}
\end{equation}
Every term is compactly supported in time within the selected time intervals.
\end{theorem}
\begin{proof}
The Galbrun potential identity applied to \eqref{amplitude:forced-equation}
gives the local linear equation for
$(w_n^{\mathrm c,\mathrm{main}},b_n^{\mathrm c,\mathrm{main}},\pi_n)$.
The cutoff terms satisfy
\[
 \sum_n\curl[\mathscr G_n,\tilde\eta_n]\Theta_n
       =\ddiv R^{\mathrm{cut}}.
\]
Using \eqref{mom:linear-transfer} and the definitions of the corrector
remainders therefore gives
\[
 \mathcal F_{v_q,B_q}(w^{\rm c},b^{\rm c})+\nabla\sum_n\pi_n
 =\ddiv\bigl(\sum_n\mathsf F_n+R^{\mathrm{cut}}
       +R^{\mathrm c}_{\rm tr}+R^{\mathrm c}_{\mathrm{na}}
       +R^{\mathrm c}_{\rm moll}\bigr).
\]
Applying \eqref{eq:endpoint-polynomial} to the corrector increment about
the old pair, and then to the principal increment about the corrected
pair, gives
\begin{equation}\label{direct:cross-grouping}
 \begin{split}
 \mathcal M(v_1,B_1)-\mathcal M(v_q,B_q)
 ={}&\mathcal F_{v_q,B_q}(w^{\rm c},b^{\rm c})
          +\ddiv R^{\mathrm c}_{\rm quad}\\
 &+\mathcal F_{v_1^{\rm c},B_1^{\rm c}}(w^{\rm p},b^{\rm p})
       +\ddiv(w^{\rm p}\otimes w^{\rm p}-b^{\rm p}\otimes b^{\rm p}).
 \end{split}
\end{equation}
In particular, the cross products occur exactly once, since
\[
 \begin{aligned}
 \mathcal F_{v_1^{\rm c},B_1^{\rm c}}(w^{\rm p},b^{\rm p})
       -\mathcal F_{v_q,B_q}(w^{\rm p},b^{\rm p})
   &=2\ddiv(w^{\rm c}\odot w^{\rm p}-b^{\rm c}\odot b^{\rm p}),\\
 \mathcal F_{v_1^{\rm c},B_1^{\rm c}}(w^{\rm p},b^{\rm p})
   &=\ddiv(R^{\mathrm p}_{\rm tr}+R^{\mathrm p}_{\mathrm{na}}
                                     +R^{\mathrm p}_{\rm moll}).
 \end{aligned}
\]
For the remaining tensors, \eqref{mom:N-remainder} and
\eqref{mom:chart-stress} give
\[
 \ddiv\bigl(w^{\rm p}\otimes w^{\rm p}-b^{\rm p}\otimes b^{\rm p}
                      -\mathsf T\bigr)
   =\ddiv(R^{\mathrm p}_{\rm quad}+R^{\mathrm p}_{\rm chart}
          +R^{\mathrm p}_{\rm time}+R^{\mathrm p}_{\rm mag}).
\]
Moreover, \eqref{mom:target-balance} gives
\[
 \begin{aligned}
 R_q+\sum_n\mathsf F_n+\mathsf T-R^{\rm moll}
   &=(1-\varrho^2)R_q+\varrho^2\delta_{q+1}\IId,\\
 \ddiv\!\left(R_q+\sum_n\mathsf F_n+\mathsf T\right)
   &=\ddiv R^{\rm moll}.
 \end{aligned}
\]
The second equality follows from
$(1-\varrho^2)R_q=0$, which holds by \eqref{prep:outer-cutoff}, and
from the fact that the remaining isotropic tensor is spatially constant. Substituting these identities into
\eqref{direct:cross-grouping} and using the old momentum equation
proves \eqref{mom:full-stress} with $p_{q+1}=p_q+\sum_n\pi_n$.
Every force to which inverse divergence is applied has zero mean and
has not been projected, so no additional pressure term appears. Volume preservation and
\eqref{path:full-fields} give incompressibility and induction.
Finally, all operators defining the stress are spatial or act in $s$
at fixed physical time; they preserve the selected time supports.
\end{proof}

The next two sections estimate the LDFs, fields and this stress.
\section{Estimates: LDFs and Fields}\label{sec:generator-field-estimates}

In this section we estimate the corrector, the principal perturbation
and the background gaps, and deduce the field bounds at the next
stage. For the gaps, the estimates must include the possibly large
Euclidean differential of the principal flow. We use $k$ and $m$ for
the material and magnetic derivative orders. When passing to the new
background, we convert Lie derivatives into derivatives of Cartesian
components. The required parameter inequalities are collected in
Section~\ref{sec:parameter-choice}.

The classes are those of Section~\ref{ssec:fixed-scales}, with the
reference fields specified in each statement. We first prove the
estimates at the local derivative scales. They remain valid with the
same amplitudes at the scales of the next iterate: for the Galbrun
equation this is Theorem~\ref{galbrun:linear-theorem}(f), for $k$
material and $m$ magnetic derivatives with $k+m\le j_1$; for the
other estimates it follows from \eqref{setup:rate-comparisons}.
Throughout this section and the next, set
\begin{equation}\label{principal:ceilings}
 \Nc=\rprep-j_0-12,
\end{equation}
where $\rprep$ is the order provided by Section~\ref{sec:preparation}.
The sharp estimates below hold through order $\Nc$. Assuming the
inequality of Section~\ref{ssec:derivative-reserve},
Lemma~\ref{setup:calculus}(vii) extends them to every spatial order
for which the corresponding lossy estimate is available, without
changing the restrictions on material and magnetic derivatives.
This gives all spatial orders for the smooth constructed quantities.
For terms containing a background gap, the derivative orders are those
of Proposition~\ref{gap:field-bounds} and \eqref{gg:crude}.

\subsection{Corrector bounds}

\begin{proposition}[Corrector LDF]\label{gg:physical}
For each $n$, take Lie derivatives with respect to the local
background on the corresponding time interval. Then
\begin{equation}\label{gg:common-generator-bounds}
 \begin{aligned}
 \Theta_n^{\rm c}&\in
 \mathcal C_{\rprep-10,j_1-5}(C\ell^{-\alpha}\tau_a^2\delta_{q+1};\lambda_q,\mathrm a)
 \cap\mathcal C_{\rprep-10}(C\ell^{-\alpha}\tau_a^2\delta_{q+1};\lambda_q,\mathrm f),\\
 \xi_n^{\rm c}&\in
 \mathcal C_{\rprep-11,j_1-5}(C\ell^{-\alpha}\tau_a^2\lambda_q\delta_{q+1};\lambda_q,\mathrm a)
 \cap\mathcal C_{\rprep-11}(C\ell^{-\alpha}\tau_a^2\lambda_q\delta_{q+1};\lambda_q,\mathrm f),\\
 \mathcal D_{t,n}\xi_n^{\rm c}&\in
 \mathcal C_{\rprep-11,j_1-5}(C\ell^{-\alpha}\tau_a\lambda_q\delta_{q+1};\lambda_q,\mathrm a)
 \cap\mathcal C_{\rprep-11}(C\ell^{-\alpha}\tau_a\lambda_q\delta_{q+1};\lambda_q,\mathrm f),\\
 \mathcal L_{B_{\ell,n}}\xi_n^{\rm c}&\in
 \mathcal C_{\rprep-11,j_1-5}(C\ell^{-\alpha}\tau_a^2\lambda_\parallel\lambda_q\delta_{q+1};\lambda_q,\mathrm a)
 \cap\mathcal C_{\rprep-11}(C\ell^{-\alpha}\tau_a^2\lambda_\parallel\lambda_q\delta_{q+1};\lambda_q,\mathrm f),
 \end{aligned}
\end{equation}
and their lossy and ordinary classes are
\[
 \begin{aligned}
 \Theta_n^{\rm c}&\in\mathcal K(C\tau_c^2\delta_{q+1})
 \cap\mathcal O_{\infty,j_1+2}(C\tau_c^2\delta_{q+1};\ell^{-1}),\\
 \xi_n^{\rm c}&\in\mathcal K(C\tau_c^2\ell^{-1}\delta_{q+1})
 \cap\mathcal O_{\infty,j_1+1}(C\tau_c^2\ell^{-1}\delta_{q+1};\ell^{-1}),\\
 \mathcal D_{t,n}\xi_n^{\rm c}&\in\mathcal K(C\tau_c\ell^{-1}\delta_{q+1})
 \cap\mathcal O_{\infty,j_1+1}(C\tau_c\ell^{-1}\delta_{q+1};\ell^{-1}),\\
 \mathcal L_{B_{\ell,n}}\xi_n^{\rm c}&\in\mathcal K(C\tau_c\ell^{-1}\delta_{q+1})
 \cap\mathcal O_{\infty,j_1+1}(C\tau_c\ell^{-1}\delta_{q+1};\ell^{-1}).
 \end{aligned}
\]
In particular, for every integer $0\le r\le\rprep-11$,
\begin{equation}\label{gg:physical-scales}
 \|\xi^{\rm c}\|_{r+\alpha}
 \le C\tau_a^2\lambda_q^{r+1}\ell^{-\alpha}\delta_{q+1},\qquad
 \|\DD\xi^{\rm c}\|_\alpha
 \le C\tau_a^2\lambda_q^2\ell^{-\alpha}\delta_{q+1}.
\end{equation}
The difference of the second material and magnetic Lie derivatives satisfies
\begin{equation}\label{gg:high-mixed-state}
 (\mathcal D_{t,n}^2-\mathcal L_{B_{\ell,n}}^2)\Theta_n^{\rm c}\in
 \mathcal C_{\rprep-10,j_1-5}(C\lambda_{q+1}^{(j_1+1)\alpha}\delta_{q+1};\lambda_q,\mathrm a)
 \cap\mathcal C_{\rprep-10}(C\lambda_{q+1}^{(j_1+1)\alpha}\delta_{q+1};\lambda_q,\mathrm f)
 \cap\mathcal K(C\lambda_{q+1}^{(j_1+1)\alpha}\delta_{q+1}),
\end{equation}
and the same estimates, including the lossy and ordinary bounds,
hold after summing over $n$, relative to any fixed active background.
\end{proposition}
\begin{proof}
Section~\ref{sec:adapted-perturbation} verifies the hypotheses of
Theorem~\ref{galbrun:linear-theorem} with $M_{\rm src}=C\delta_{q+1}$
and slow source tensors bounded through derivative order $\rprep$.
Parts (a), (c) and (f) give the estimates for $\Theta_n^{\rm c}$ and
its first material and magnetic derivatives through order $\rprep-10$.
The lossy estimates follow by expanding the transport operators in
ordinary derivatives and using the lossy bounds for the fields.
Since the background fields preserve volume, curl commutes with their
Lie derivatives, and hence
\[
 \mathcal D_{t,n}^k\mathcal L_{B_{\ell,n}}^m\xi_n^{\rm c}
 =\curl\mathcal D_{t,n}^k\mathcal L_{B_{\ell,n}}^m\Theta_n^{\rm c}.
\]
Lemma~\ref{setup:calculus}(ii) now gives the asserted estimates for
$\xi_n^{\rm c}$ and its first variations, with one additional spatial
derivative of the potential. For the spatial H\"older norm, use
Proposition~\ref{galbrun:fast-interior} and the proof of
Theorem~\ref{galbrun:linear-theorem}. Since the cutoff depends only
on time, they give
\[
 \|\Theta_n^{\rm c}\|_{r+\alpha}
 \le C\tau_a^2\lambda_q^r\ell^{-\alpha}\delta_{q+1},
 \qquad 0\le r\le\rprep-10.
\]
Taking curl and using this estimate at order $r+1$ proves the first
inequality in \eqref{gg:physical-scales} for $0\le r\le\rprep-11$.
The case $r=1$ gives the second inequality, also recorded in
\eqref{galbrun:low-deformation}. The sum over $n$ satisfies the same
bounds because the time supports have bounded overlap.

For \eqref{gg:high-mixed-state}, write the local equation in terms of Lie derivatives as
\[
 (\mathcal D_{t,n}^2-\mathcal L_{B_{\ell,n}}^2)\Theta_n^{\rm c}
 =\mathscr T\mathsf F_n+\mathscr C^{\rm cut}[\tilde\eta_n,\Theta_n]
 -2\mathscr T\bigl((\mathcal D_{t,n}\Theta_n^{\rm c})\times\nabla v_{\ell,n}
       -(\mathcal L_{B_{\ell,n}}\Theta_n^{\rm c})\times\nabla B_{\ell,n}\bigr),
\]
as given by Proposition~\ref{mom:lower-structure}.
By \eqref{amplitude:classes}, the source satisfies
\[
 \mathsf F_n\in
 \mathcal C_{\rprep}(C\delta_{q+1};\lambda_q,\mathrm a)
 \cap\mathcal K(C\delta_{q+1}).
\]
A material derivative of a fast profile costs $\tau_a^{-1}$, so
Lemma~\ref{setup:calculus}(v) gives the required bounds for
$\mathscr T\mathsf F_n$ in \eqref{gg:high-mixed-state}. The cutoff
commutator is controlled by Theorem~\ref{galbrun:linear-theorem}(d),
which places it in
\[
 \mathcal C_{\rprep-2}(C\ell^{-\alpha}\delta_{q+1};\lambda_q,\mathrm a)
 \cap\mathcal K(C\delta_{q+1})
\]
for $k+m\le j_1$. These bounds imply the classes in
\eqref{gg:high-mixed-state} because the local derivative costs are at
most those of the next iterate and
$\ell^{-\alpha}\le\lambda_{q+1}^{\alpha}$.

For the remaining terms, apply Lemma~\ref{setup:calculus}(i) to the
background gradients in \eqref{prep:local-gradient} and the first
derivatives of the potential. The two products inside $\mathscr T$
have amplitudes
\[
 \begin{aligned}
 C\ell^{-\alpha}\lambda_q\delta_q^{1/2}\tau_a\delta_{q+1}
 &=C\ell^{-\alpha}\varepsilon_\tau\varepsilon_{q+1}^{\gamma_\ell}\delta_{q+1},\\
 C\ell^{-\alpha}\lambda_q\delta_{B,q}^{1/2}\tau_a^2\lambda_\parallel\delta_{q+1}
 &\le C\ell^{-\alpha}(\tau_a\lambda_\parallel)^2\varepsilon_{q+1}^{\gamma_\ell}\delta_{q+1}.
 \end{aligned}
\]
Both amplitudes are bounded by $C\delta_{q+1}$; applying rule (v)
therefore yields the asserted classes. The same products belong to
\[
 \mathcal K(C\ell^{-\alpha}\varepsilon_{q+1}^{\gamma_\ell}\delta_{q+1})
 \subset\mathcal K(C\delta_{q+1}),
\]
since $C\ell^{-\alpha}\varepsilon_{q+1}^{\gamma_\ell}\delta_{q+1}
\le C\delta_{q+1}$, using $\ell^{-\alpha}\le\lambda_{q+1}^{\alpha}$.
Using the sharp bounds from part (a) of the theorem in these products
gives the local costs $\mathrm a$ for compositions of at most $j_1-5$
material and magnetic derivatives.

For distinct local backgrounds, we have
\begin{equation}\label{gg:adjacent}
 \begin{aligned}
 \mathcal D_{t,n}\xi_j^{\rm c}
   &=\mathcal D_{t,j}\xi_j^{\rm c}
                +\mathcal L_{v_{\ell,n}-v_{\ell,j}}\xi_j^{\rm c},\\
 \mathcal L_{B_{\ell,n}}\xi_j^{\rm c}
   &=\mathcal L_{B_{\ell,j}}\xi_j^{\rm c}
                +\mathcal L_{B_{\ell,n}-B_{\ell,j}}\xi_j^{\rm c}.
 \end{aligned}
\end{equation}
The additional Lie derivatives in these identities are bounded if the
background differences, normalized by the local derivative costs,
satisfy
\[
 \tau_a\lambda_q\|v_{\ell,n}-v_{\ell,j}\|_0\le C,
 \qquad
 \frac{\lambda_q\|B_{\ell,n}-B_{\ell,j}\|_0}{\lambda_\parallel}\le C,
\]
and the corresponding differentiated estimates, including one
additional spatial derivative for tensor Lie derivatives. These
bounds follow from Corollary~\ref{prep:mixed-comparison}.
Lemma~\ref{setup:calculus}(iv), with $\Lambda=\lambda_q$, therefore
compares the two backgrounds through order $\rprep$ on their
overlapping time intervals, for either pair of derivative costs.
The constants are uniform, so bounded temporal overlap gives the
claimed estimates for the sum relative to any fixed active background.
\end{proof}

\paragraph{\textbf{Corrector pushforward.}}

We apply Lemma~\ref{gg:transport-classes} with $\xi=\xi^{\rm c}$ and
$(v,B)=(v_{\ell,n},B_{\ell,n})$. The Lie derivatives satisfy
\begin{equation}\label{gg:intertwining}
 \mathcal D_{t,n,s}^{\rm c}\mathcal U_s^{\rm c}=\mathcal U_s^{\rm c}\mathcal D_{t,n},\qquad
 \mathcal L_{B_{\ell,n,s}^{\rm c}}\mathcal U_s^{\rm c}=\mathcal U_s^{\rm c}\mathcal L_{B_{\ell,n}}.
\end{equation}
At $\Lambda=\lambda_q$, the spatial hypotheses follow from
\eqref{gg:physical-scales} through $N'+1\le\rprep-11$, with size
$C\tau_a^2\lambda_q^2\ell^{-\alpha}\delta_{q+1}$. At the higher
frequency $\Lambda=\lambda_{q+1}$, the lossy bounds control the higher
spatial orders. Indeed, for $r\ge2$, they give
\[
 \begin{aligned}
 \|\xi^{\rm c}\|_{r+\alpha}
 &\le C\tau_c^2\ell^{-1-r-\alpha}\delta_{q+1}\\
 &=C\tau_c^2\ell^{-2}\delta_{q+1}\,\ell^{-\alpha}
   (\ell\lambda_{q+1})^{-(r-1)}\lambda_{q+1}^{r-1}
 \le C\lambda_{q+1}^{r-1}.
 \end{aligned}
\]
Here $\tau_c^2\ell^{-2}\delta_{q+1}=\varepsilon_{q+1}^{2\beta}\le1$ and
$\ell^{-\alpha}\le\lambda_{q+1}^{\alpha}\le\ell\lambda_{q+1}$;
the $r=1$ bound is \eqref{gg:physical-scales}.
Thus pushforward by the corrector preserves the following classes, with
Lie derivatives taken along the corrected background:
\begin{equation}\label{gg:corrector-class-ranges}
 \begin{gathered}
 \mathcal C_{N',J}(E;\Lambda,\rho),\qquad
 \rho\in\{\mathrm a,\mathrm f\},\quad J\le j_1,\\
 \begin{cases}
 N'\le\rprep-12,&\Lambda=\lambda_q,\\
 N'\le\infty,&\Lambda=\lambda_{q+1},
 \end{cases}
 \qquad
 \mathcal K_{N'',J}(E'),\quad N''\le\infty,
 \end{gathered}
\end{equation}
up to fixed constants. Applying the ordinary derivative bounds for
the LDF gives the same conclusion for
$\mathcal O_{N,H}(E_{\rm o};\Lambda_{\rm o})$, where
$\Lambda_{\rm o}\in\{\ell^{-1},\lambda_{q+1}\}$,
$N\le\infty$ and $0\le H\le j_1+1$, independently of $N',N''$.
The spatial estimates involve derivatives of the LDF; the corrector
displacement need not be small at the oscillation scale.

For estimates relative to a fixed background, the lemma also requires
bounds for the first variations, with the material and magnetic Lie
derivative orders summing to at most $J-1$ and the spatial, material
and magnetic orders summing to at most $N'+1$.
These follow from \eqref{gg:common-generator-bounds} and
\eqref{gg:adjacent} for $J\le j_1$ when $\rho=\mathrm f$, and for
$J\le j_1-4$ when $\rho=\mathrm a$, in both cases with
$N'\le\rprep-12$. Using the identities
\begin{equation}\label{prep:corrector-final-comparison}
 \begin{aligned}
 c_{v,n,s}=v_{\ell,n,s}^{\rm c}-v_{\ell,n}
   &=\int_0^s\mathcal U_u^{\rm c}\mathcal D_{t,n}\xi^{\rm c}\,\dd u,\\
 c_{B,n,s}=B_{\ell,n,s}^{\rm c}-B_{\ell,n}
   &=\int_0^s\mathcal U_u^{\rm c}\mathcal L_{B_{\ell,n}}\xi^{\rm c}\,\dd u.
 \end{aligned}
\end{equation}
we therefore obtain, relative to the fixed local background,
\begin{equation}\label{gg:corrector-comparison-classes}
 \begin{aligned}
 c_{v,n,s}&\in\mathcal C_{\rprep-12}
 (C\ell^{-\alpha}\tau_a\lambda_q\delta_{q+1};\lambda_q,\mathrm f),\\
 c_{B,n,s}&\in\mathcal C_{\rprep-12}
 (C\ell^{-\alpha}\tau_a^2\lambda_\parallel\lambda_q\delta_{q+1};\lambda_q,\mathrm f),\\
 \nabla c_{v,n,s}&\in\mathcal C_{\rprep-13}
 (C\ell^{-\alpha}\tau_a\lambda_q^2\delta_{q+1};\lambda_q,\mathrm f),\\
 \nabla c_{B,n,s}&\in\mathcal C_{\rprep-13}
 (C\ell^{-\alpha}\tau_a^2\lambda_\parallel\lambda_q^2\delta_{q+1};\lambda_q,\mathrm f).
 \end{aligned}
\end{equation}
The same bounds hold at the fast local derivative costs for
compositions of at most $j_1-5$ material and magnetic derivatives.
The lossy classes are
\[
 \begin{aligned}
 c_{v,n,s},c_{B,n,s}&\in\mathcal K(C\tau_c\ell^{-1}\delta_{q+1}),\\
 \nabla c_{v,n,s},\nabla c_{B,n,s}&\in\mathcal K(C\tau_c\ell^{-2}\delta_{q+1}).
 \end{aligned}
\]
To pass to transport derivatives, we use the adaptedness of the local background.

After normalization by the derivative costs of the next iterate, the
largest factors in the change of transport occur at
$\Lambda=\lambda_{q+1}$ and satisfy
\begin{equation}\label{gg:corrector-rate-comparisons}
 \begin{gathered}
 \ell^{-\alpha}\tau_a\lambda_q\delta_{q+1}^{1/2}
 \le\varepsilon_{q+1}^{\beta+\gamma_a+\gamma_\ell-\gamma_S},\\
 \ell^{-\alpha}(\tau_a\lambda_\parallel)
       (\tau_a\lambda_q\delta_{q+1}^{1/2})
       \Bigl(\frac{\delta_{q+1}}{\delta_{B,q+1}}\Bigr)^{1/2}
 \le\varepsilon_{q+1}^{\beta+\gamma_a+\gamma_\ell-\gamma_S}.
 \end{gathered}
\end{equation}
Both exponents are positive. At $\Lambda=\lambda_q$ each quotient
has the additional factor $\varepsilon_{q+1}$. At the fast local
derivative costs and $\Lambda=\lambda_q$, both factors are at most
$C\tau_a^2\lambda_q^2\ell^{-\alpha}\delta_{q+1}$; at the lossy
costs they are at most $C\tau_c\ell^{-1}\delta_{q+1}\le C$.
The gradient amplitudes in \eqref{gg:corrector-comparison-classes}
are at most $C\tau_a^{-1}$ and $C\lambda_\parallel$, respectively.
Adding the original gradients from \eqref{prep:local-gradient} and
using the local-background derivative expansion in
Lemma~\ref{aniso:old-frame-operator-change}, we obtain
\begin{equation}\label{gg:corrected-adaptedness}
 (v_{\ell,n,s}^{\rm c},B_{\ell,n,s}^{\rm c})\text{ is }
 \begin{cases}
 (\lambda_q,\mathrm f)\text{-adapted through }(\rprep-12,j_1),\\
 (\lambda_q,\mathrm a)\text{-adapted through }(\rprep-12,j_1-4).
 \end{cases}
\end{equation}
The range follows from the definition of adaptedness: at order
$(N,J)$ it uses gradient bounds only through $(N-1,J-1)$.
Lemma~\ref{setup:calculus}(iv) then compares derivatives relative to
the local and corrected backgrounds in both directions. At the costs
of the next iterate, it applies through order $\rprep-12$ at either
spatial frequency. At frequency $\lambda_q$ and the fast local costs,
it applies to at most $j_1-4$ material and magnetic derivatives in
all. The $\mathcal K_{N'',J}(E')$ estimates in
\eqref{gg:corrector-class-ranges} are likewise equivalent, with fixed
constants.

For $F$ in the sharp class at frequency $\Lambda=\lambda_q$ or in
the lossy class of \eqref{gg:corrector-class-ranges},
Lemma~\ref{gg:transport-classes} gives, with either background fixed,
\begin{equation}\label{gg:corrector-defects}
 \begin{aligned}
 (\mathcal U_s^{\rm c}-\IId)F&\in\mathcal C_{N'-1,J}
 (C\tau_a^2\lambda_q^2\ell^{-\alpha}\delta_{q+1}E;\lambda_q,\rho),\\
 (\mathcal U_s^{\rm c}-\IId)F&\in\mathcal K_{N''-1,J}
 (C\tau_c^2\ell^{-2}\delta_{q+1}E'),
 \end{aligned}
\end{equation}
where $N'\le\rprep-12$, with $J\le j_1$ for $\rho=\mathrm f$
and $J\le j_1-4$ for $\rho=\mathrm a$. At the lossy scale, the normalized field increments are bounded by
the same factor, since $\tau_c\ell^{-1}\ge1$. Integrating with the
background fixed gives the same estimates for $\mathcal P^{\rm c}$ and
its difference from the identity.

Finally, composition of the chart with the inverse corrector flow and tensor pushforward of its dual frames give
\begin{equation}\label{gg:corrected-chart-bounds}
 \begin{gathered}
 y_{I,s}^\eta=y_I^\eta\circ(X_s^{\rm c})^{-1},\qquad
 dy_{I,s}^\eta=\mathcal U_s^{\rm c}dy_I^\eta,\qquad
 \frac{\partial}{\partial y_{I,s}^\eta}=\mathcal U_s^{\rm c}\frac{\partial}{\partial y_I^\eta}.
 \end{gathered}
\end{equation}
Apply the tensor pushforward estimate to the frames and the ordinary
flow estimate to the coordinates. This gives the chart bounds of
Section~\ref{ssec:material-partition} through order $\rprep-13$ for
the corrected charts. Their material and magnetic derivatives are
given by \eqref{direct:stationary-data}. The transported profile
formula at the end of Section~\ref{ssec:transported-charts} also
applies, giving the estimates for higher derivatives.
\subsection{The approximate antiderivative and scalar transport}
\label{ssec:principal-spatial-continuation}

\emph{The approximate antiderivative.}
We apply
Lemma~\ref{principal:material-primitive} with
\[
 N=j_0,\qquad J=j_1,\qquad F=\delta_{q+1}^{1/2},\qquad
 \Lambda=\lambda_q,\qquad \mathrm a_t=\tau_c^{-1},\qquad \mathrm a_B=\lambda_\parallel,
\]
and \eqref{amplitude:amplitude-bounds}. For a fixed index $I$,
applying $D_{t,I}^kD_{B,I}^m$ to the approximate antiderivative
requires the spatial norm of order $r$ of
$D_{t,I}^{k+j}D_{B,I}^ma_I$, with $j\le j_0$.
These amplitude estimates are available when
\[
 r+\min\{k+m+j_0,j_1\}\le \rprep,
\]
in particular when $r+k+m\le \rprep-j_0$. Accounting also for the derivative loss in the amplitude estimate, the
lemma gives the following bounds at the fast local derivative costs:
\begin{equation}\label{principal:primitive-sharp}
 \begin{aligned}
 \|D_{t,I}^kD_{B,I}^m\mathfrak a_I\|_{r}
 &\le C\delta_{q+1}^{1/2}\lambda_q^r\tau_a^{1-k}\lambda_\parallel^m
        \varepsilon_\tau^{-[m-j_1]^+},\\
 \|D_{t,I}^kD_{B,I}^ma_I^{\rm c}\|_{r}
 &\le C\delta_{q+1}^{1/2}\lambda_q^r\tau_a^{-k}\lambda_\parallel^m
        \varepsilon_\tau^{\min\{j_0,j_1-m\}},\\
 \|D_{t,I}^kD_{B,I}^m(\mathfrak a_I-\tau_a\alpha_I^{[1]}a_I)\|_{r}
 &\le C\delta_{q+1}^{1/2}\lambda_q^r\tau_a^{1-k}\lambda_\parallel^m
        \varepsilon_\tau^{\min\{1,j_1-m\}},
 \end{aligned}
\end{equation}
for
\[
 k+m\le j_1+2,\qquad
 r+\min\{k+m+j_0,j_1\}\le \rprep.
\]
The corresponding H\"older bounds include the factor $\ell^{-\alpha}$.
For the lossy bounds, valid for $r\ge0$ and $k+m\le j_1+2$,
replace $\lambda_q^r$ by
$\ell^{-r}$ and $\lambda_\parallel^m$ by
$s_B^{-m}=\varepsilon_\tau^{-m}\lambda_\parallel^m$, omitting the gain
in the second and third lines. The ordinary derivatives
follow from \eqref{amplitude:crude-amplitude-bounds}: if $b$ of the
$h$ additional ordinary time derivatives fall on the amplitude in the $j$th
summand of $\mathfrak a_I$,
\[
 \tau_a^{j+1}\tau_a^{-(h-b)}\|\partial_t^bD_{t,I}^ja_I\|_{r}
 \le C\tau_a\delta_{q+1}^{1/2}\tau_a^{-(h-b)}\ell^{-r-b}
 \le C\tau_a\delta_{q+1}^{1/2}\ell^{-r-h}.
\]
The factor $\tau_a^j$ thus compensates for differentiating the
amplitude in the $j$th summand, also at $h=j_1+1$.
Combining these estimates with \eqref{aniso:finite-depth-budget} gives
\begin{equation}\label{principal:primitive-classes}
 \begin{aligned}
 \mathfrak a_I&\in\mathcal C_{\rprep-j_0}(C\tau_a\delta_{q+1}^{1/2};\lambda_q,\mathrm a)
 \cap\mathcal K(C\tau_a\delta_{q+1}^{1/2})\cap\mathcal O_{\infty,j_1+1}(C\tau_a\delta_{q+1}^{1/2};\ell^{-1}),\\
 a_I^{\rm c}&\in\mathcal C_{\rprep-j_0}(C\varepsilon_\tau^{j_0}\delta_{q+1}^{1/2};\lambda_q,\mathrm f)
 \cap\mathcal C_{\rprep-j_0,j_1-j_0}(C\varepsilon_\tau^{j_0}\delta_{q+1}^{1/2};\lambda_q,\mathrm a)
 \cap\mathcal K(C\delta_{q+1}^{1/2}),\\
 \mathfrak a_I-\tau_a\alpha_I^{[1]}a_I&\in
 \mathcal C_{\rprep-j_0}(C\varepsilon_\tau\tau_a\delta_{q+1}^{1/2};\lambda_q,\mathrm f)
 \cap\mathcal K(C\tau_a\delta_{q+1}^{1/2}),\\
 D_{B,I}\mathfrak a_I&\in\mathcal C_{\rprep-j_0-1}(C\tau_a\lambda_\parallel\delta_{q+1}^{1/2};\lambda_q,\mathrm f)
 \cap\mathcal C_{\rprep-j_0-1,j_1-1}(C\tau_a\lambda_\parallel\delta_{q+1}^{1/2};\lambda_q,\mathrm a)
 \cap\mathcal K(C\varepsilon_\tau^{-1}\tau_a\lambda_\parallel\delta_{q+1}^{1/2}),
 \end{aligned}
\end{equation}
with the ordinary classes $\mathcal O_{\infty,j_1+1}(\,\cdot\,;\ell^{-1})$ of
the same amplitudes as the lossy ones. The ratios of the lossy amplitudes to the sharp amplitudes are
$\varepsilon_\tau^{-j_0}$ for $a_I^{\rm c}$, $\varepsilon_\tau^{-1}$ for the
difference and for $D_{B,I}\mathfrak a_I$, and one for $\mathfrak a_I$.
Indeed, for $d\in\{1,j_0\}$ and $m\le j_1$,
\[
 \varepsilon_\tau^{\min\{d,j_1-m\}}
 \Bigl(\frac{\tau_a^{-1}}{\lambda_{q+1}\delta_{q+1}^{1/2}}\Bigr)^k
 \Bigl(\frac{\lambda_\parallel}{\lambda_{q+1}\delta_{B,q+1}^{1/2}}\Bigr)^m
 \le\varepsilon_\tau^d ;
\]
for $d=j_0$ this follows from \eqref{aniso:finite-depth-budget} with
$s=0$, since both ratios of derivative costs are
$\varepsilon_{q+1}^{1-\beta-\gamma_a-\gamma_\ell}\le\varepsilon_\tau$, and
for $d=1$ it follows from $\min\{1,j_1-m\}+m\ge1$ and
\eqref{setup:rate-comparisons}. The first magnetic derivative incurs the loss
$\varepsilon_\tau^{-1}$ only at $m=j_1$. At that order,
\[
 \varepsilon_\tau^{-1}
 (\lambda_\parallel/(\lambda_{q+1}\delta_{B,q+1}^{1/2}))^{j_1}
 \le\varepsilon_\tau^{j_1-1}\le1.
\]

\emph{The transported approximate antiderivative and scalar transport.}
The principal LDF and its corrector pushforwards are
\[
 \xi_0^{\rm p}=\sum_I\curl\bigl(\lambda_{q+1}^{-1}\mathfrak a_I
       \varphi(\lambda_{q+1}y_I^{k_I})\,dy_I^{\nu_I}\bigr),
 \qquad \xi_s^{\rm p}=\mathcal U_s^{\rm c}\xi_0^{\rm p}.
\]
The coefficients in these finite profile sums consist of frames and
first spatial derivatives of $\mathfrak a_I$, so
Lemma~\ref{setup:calculus}(vii) applies. In the original chart on
support $I$, formula \eqref{principal:scalar-action} and
\eqref{principal:primitive-classes} bound the normalized coefficients
of scalar transport by $C\tau_a\ell^{-1}\delta_{q+1}^{1/2}$.
The coefficient estimates require one additional spatial derivative
of $\mathfrak a_I$ at the same material and magnetic derivative orders. In particular, the transport
estimates of Subsection~\ref{ssec:path-estimates} through order $N$
use at most $N+2$ spatial derivatives of $\mathfrak a_I$.
The choice $\Nc=\rprep-j_0-12$ therefore suffices for these estimates
and the corrector pushforward \eqref{gg:corrector-class-ranges}.
For higher spatial orders we use
Corollary~\ref{principal:rescaled-profile-continuation}, which requires
\[
 \ell^{-1}\le\Lambda_1\le\lambda_{q+1},\qquad
 \tau_a\delta_{q+1}^{1/2}\Lambda_1\le1.
\]
By \eqref{setup:rate-comparisons},
$\tau_a\lambda_{q+1}\delta_{q+1}^{1/2}\ge1$ and
$\tau_a\ell^{-1}\delta_{q+1}^{1/2}=o(1)$. We therefore take the largest admissible transverse frequency,
which gives the strongest gain for the slow coefficients:
\begin{equation}\label{principal:continuation-data}
 \Lambda_1=(\tau_a\delta_{q+1}^{1/2})^{-1},\qquad
 \vartheta=(\ell\Lambda_1)^{-1}=\tau_a\ell^{-1}\delta_{q+1}^{1/2}.
\end{equation}
With this choice, \eqref{iter:general-final-reserve} compensates for
the losses listed below and preserves the sharp amplitudes. Both
potentials have zero initial data. In the last row, the source is
that of $f_{w,I,s}-s\lambda_{q+1}^{-1}\alpha_I a_I\varphi$; its
equation is derived in Proposition~\ref{principal:nonprincipal-field-smallness}.
\begingroup
\renewcommand{\arraystretch}{1.15}
\[
\begin{array}{c|c|c}
 \text{coefficient or source}&\text{amplitude }E&\text{loss }E'/E\\ \hline
 \xi_{I,0}^{\rm p}\cn\text{ in the normalized basis}&\tau_a\ell^{-1}\delta_{q+1}^{1/2}&1\\
 \lambda_{q+1}^{-1}(\alpha_I a_I+a_I^{\rm c})\varphi
       &\lambda_{q+1}^{-1}\delta_{q+1}^{1/2}&1\\
 \lambda_{q+1}^{-1}(D_{B,I}\mathfrak a_I)\varphi
       &\lambda_{q+1}^{-1}\tau_a\lambda_\parallel\delta_{q+1}^{1/2}&\varepsilon_\tau^{-1}\\
 \lambda_{q+1}^{-1}(a_I^{\rm c}\varphi-s\,\xi_{I,0}^{\rm p}\cn(\alpha_I a_I\varphi))
       &\lambda_{q+1}^{-1}\delta_{q+1}^{1/2}
           (\varepsilon_\tau^{j_0}+\tau_a\ell^{-1}\delta_{q+1}^{1/2}\varepsilon_\tau)
       &\varepsilon_\tau^{-j_0}
\end{array}
\]
\endgroup
In each row, the slow coefficients, with the displayed spatial profiles
factored out, belong to $\mathcal K(CE')$, where $E'$ is determined
by the last two columns. The comparison above gives the sharp amplitudes
also at the derivative costs of the next iterate, with the losses from
\eqref{principal:primitive-classes}. The inverse characteristics and
the profiles composed with them satisfy the corresponding estimates
by Corollary~\ref{principal:phase-flow}.

\subsection{Principal LDF}

\begin{proposition}[Principal LDF]\label{principal:generator-estimates}
Fix $I$ and take derivatives with respect to the corrected background
$(v_{\ell,I,s}^{\rm c},B_{\ell,I,s}^{\rm c})$. On the image of the
local support under $X_s^{\rm c}$, we have
\begin{equation}
 \xi_{I,s}^{\rm p}\in
 \mathcal C_{\Nc}(C\tau_a\delta_{q+1}^{1/2};\lambda_{q+1},\mathrm a)
 \cap\mathcal C_{\Nc}(C\tau_a\delta_{q+1}^{1/2};\lambda_{q+1},\mathrm f)
 \cap\mathcal O_{\infty,j_1+1}(C\tau_a\delta_{q+1}^{1/2};\lambda_{q+1}),
\end{equation}
uniformly in $0\le s\le1$. The slow coefficients in the finite profile
expansion belong to $\mathcal K(C\tau_a\delta_{q+1}^{1/2})$.
Thus their sharp and lossy amplitudes agree up to fixed constants,
in the sense of Lemma~\ref{setup:calculus}(vii).
The same estimates hold for the sum
of the local contributions, whose supports are disjoint.
\end{proposition}
\begin{proof}
The phase and $dy_I^{\nu_I}$ have zero material and magnetic Lie
derivatives, so
\[
 (\mathcal D_{t,I}^k\mathcal L_{B_{\ell,I}}^m)\Theta_{I,0}^{\rm p}
 =\lambda_{q+1}^{-1}(D_{t,I}^kD_{B,I}^m\mathfrak a_I)
                   \varphi(\lambda_{q+1}y_I^{k_I})dy_I^{\nu_I} .
\]
Use \eqref{principal:primitive-classes} and take curl. Since
$\ell^{-1}\le\lambda_{q+1}$, curl and each further spatial derivative
cost at most $\lambda_{q+1}$; the H\"older factor
$\lambda_{q+1}^{\alpha}$ also bounds $\ell^{-\alpha}$.
This proves the estimate at $s=0$.
Lemma~\ref{gg:transport-classes} gives the estimate for every $s$.
For ordinary derivatives, apply the product and chain rules, using
the lossy bounds for the chart and corrector map. Finally, all
derivatives are supported in the same disjoint local supports, so
summing changes only a constant depending on the derivative order.
\end{proof}

\subsection{Principal fields relative to the local background}

\begin{proposition}[Principal fields relative to the local background]
\label{principal:nonprincipal-field-smallness}
The exact potentials and fields in
\eqref{principal:exact-local-increments} satisfy, with respect to the
corrected transports $\mathcal D_{t,I,s}^{\rm c}$ and
$\mathcal L_{B_{\ell,I,s}^{\rm c}}$, uniformly for $0\le s\le1$,
\begin{equation}\label{principal:potentials-bound}
 \begin{aligned}
 \Theta_{w,I,s}^{\rm p}&\in\mathcal C_{\Nc}(C\lambda_{q+1}^{-1}\delta_{q+1}^{1/2};\lambda_{q+1},\mathrm f),\\
 \Theta_{b,I,s}^{\rm p}&\in\mathcal C_{\Nc}(C\lambda_{q+1}^{-1}\tau_a\lambda_\parallel\delta_{q+1}^{1/2};\lambda_{q+1},\mathrm f).
 \end{aligned}
\end{equation}
The corresponding fields satisfy
\begin{equation}\label{principal:fields-bound}
 w_{I,s}\in\mathcal C_{\Nc}(C\delta_{q+1}^{1/2};\lambda_{q+1},\mathrm f),
 \qquad
 b_{I,s}\in\mathcal C_{\Nc}(C\tau_a\lambda_\parallel\delta_{q+1}^{1/2};\lambda_{q+1},\mathrm f),
\end{equation}
and the same estimates hold at the local derivative costs $\mathrm a$
for $\Theta_{w,I,s}^{\rm p},w_{I,s}$ after at most $j_1$ material
and magnetic Lie derivatives in all, and for
$\Theta_{b,I,s}^{\rm p},b_{I,s}$ after at most $j_1-1$ such
derivatives. Set
\[
 w_I^{\rm lead}=\lambda_{q+1}^{-1}\curl\mathcal U_1^{\rm c}
  \bigl(\alpha_I(t/\tau_a)a_I
         \varphi(\lambda_{q+1}y_I^{k_I})dy_I^{\nu_I}\bigr).
\]
Then
\begin{equation}\label{principal:full-magnetic-natural}
 w_{I,1}-w_I^{\rm lead}\in\mathcal C_{\Nc}\Bigl(C\bigl(\varepsilon_\tau^{j_0}
       +\tau_a\ell^{-1}\delta_{q+1}^{1/2}\varepsilon_\tau\bigr)
 \delta_{q+1}^{1/2};\lambda_{q+1},\mathrm f\Bigr),
\end{equation}
and the ordinary classes are
\begin{equation}\label{principal:full-magnetic-time}
 \begin{aligned}
 w_{I,1}-w_I^{\rm lead}&\in\mathcal O_{\infty,j_1+1}\Bigl(C(\varepsilon_\tau^{j_0}+\tau_a\ell^{-1}\delta_{q+1}^{1/2}\varepsilon_\tau)
            \delta_{q+1}^{1/2};\lambda_{q+1}\Bigr),\\
 b_{I,s}&\in\mathcal O_{\infty,j_1+1}(C\tau_a\lambda_\parallel\delta_{q+1}^{1/2};\lambda_{q+1}).
 \end{aligned}
\end{equation}
For $w_{I,s}$, $b_{I,s}$ and $w_{I,1}-w_I^{\rm lead}$, the slow
coefficients in the source formulas, before composition with the
principal flow, have losses $1$, $\varepsilon_\tau^{-1}$ and
$\varepsilon_\tau^{-j_0}$, respectively.
\end{proposition}
\begin{proof}
Work first in the original chart and apply
Lemma~\ref{principal:finite-transport} to the two scalar equations
\eqref{principal:profile-equation}. Its coefficient bound follows
from \eqref{principal:scalar-action} and
\eqref{principal:primitive-classes}: the normalized coefficients of
$\xi_{I,0}^{\rm p}\cn$ have size
$C\tau_a\ell^{-1}\delta_{q+1}^{1/2}$. The source estimates are given
in Section~\ref{ssec:principal-spatial-continuation}, with one extra
magnetic derivative of $\mathfrak a_I$ for the magnetic source.
Thus the lemma applies through order $\Nc+1$ at both the original
chart scales and the derivative costs of the next iterate.

By \eqref{principal:phase-path-average}, composing the sources with
the characteristics preserves their amplitudes and derivative
weights. Converting chart derivatives to physical derivatives gives
\eqref{principal:potentials-bound}, in fact through order $\Nc+1$.
Taking curl gives \eqref{principal:fields-bound} through order $\Nc$.
Finally apply the corrector pushforward. It intertwines the material
and magnetic Lie derivatives, transforms the phase by
\eqref{principal:phase-corrector-composition}, and preserves the
estimates by Lemma~\ref{gg:transport-classes}.

To estimate the difference from the leading velocity, we use the
following identity, valid for any slow scalar $g$:
\begin{equation}\label{principal:principal-constant-bound}
 \xi_{I,0}^{\rm p}\cn(g\varphi)
 =\varphi\varphi'(\mathfrak a_I\partial_{y_I^{\zeta_I}}g-g\partial_{y_I^{\zeta_I}}\mathfrak a_I)
 +\lambda_{q+1}^{-1}\varphi^2
       (\partial_{y_I^{k_I}}\mathfrak a_I\,\partial_{y_I^{\zeta_I}}g
              -\partial_{y_I^{\zeta_I}}\mathfrak a_I\,\partial_{y_I^{k_I}}g).
\end{equation}
For $g=\alpha_I a_I$, substitute
\[
 \mathfrak a_I=\tau_a\alpha_I^{[1]}a_I
                  +(\mathfrak a_I-\tau_a\alpha_I^{[1]}a_I).
\]
In each parenthesis, the contribution of the first summand cancels.
The remaining terms contain the difference estimated in
\eqref{principal:primitive-classes}, so
$\xi_{I,0}^{\rm p}\cn(\alpha_I a_I\varphi)$ has amplitude
$C\tau_a\ell^{-1}\varepsilon_\tau\delta_{q+1}$ at the costs of the
next iterate. This cancellation controls the change of the leading
velocity along the principal flow. Indeed, subtracting
$s\lambda_{q+1}^{-1}\alpha_I a_I\varphi$ from the first equation in
\eqref{principal:profile-equation} leaves the source
\[
 \lambda_{q+1}^{-1}a_I^{\rm c}\varphi
       -s\lambda_{q+1}^{-1}\xi_{I,0}^{\rm p}\cn(\alpha_I a_I\varphi).
\]
The amplitude of this source is the bound in
\eqref{principal:full-magnetic-natural} divided by $\lambda_{q+1}$.
The scalar transport estimate followed by curl proves the claim.
The cancellation is polynomial, so it is also valid at zeros of
$a_I$ or the profiles. Differentiating the same products in ordinary
coordinates and using the ordinary estimates for the approximate
antiderivative proves \eqref{principal:full-magnetic-time}, with the
losses listed in Section~\ref{ssec:principal-spatial-continuation}.
\end{proof}

The preceding calculation also estimates the Lie derivatives needed
for the quadratic tensor. In the chart norms at the original scales,
through derivative order $\Nc$, we have
\begin{equation}\label{principal:scalar-action-bounds}
 \begin{aligned}
 \|\mathcal L_{\xi_{I,0}^{\rm p}}^j\mathcal D_{t,I}\Theta_{I,0}^{\rm p}\|_{N,J}
 &\le C\lambda_{q+1}^{-1}\delta_{q+1}^{1/2}
       \tau_a^j\ell^{-j}\delta_{q+1}^{j/2}
       \bigl(\varepsilon_\tau^{\min\{1,j_1-J\}}+\varepsilon_\tau^{\min\{j_0,j_1-J\}}\bigr),
       \qquad j=1,2,\\
 \|\mathcal L_{\xi_{I,0}^{\rm p}}^j\mathcal L_{B_{\ell,I}}\Theta_{I,0}^{\rm p}\|_{N,J}
 &\le C\lambda_{q+1}^{-1}\delta_{q+1}^{1/2}
       \tau_a^j\ell^{-j}\delta_{q+1}^{j/2}\tau_a\lambda_\parallel
       \varepsilon_\tau^{-[J+1-j_1]^+},\qquad j=0,1,2,
 \end{aligned}
\end{equation}
where $\|\cdot\|_{N,J}$ is the chart norm
\eqref{principal:profile-norm} of the scalar coefficient of
$dy_I^{\nu_I}$, with at most $J\le j_1$ material and magnetic
derivatives in all. For $N\le\Nc$, the coefficient derivatives in
\eqref{principal:scalar-action} are bounded by
\eqref{principal:primitive-sharp}, since $\Nc=\rprep-j_0-12$.
In the first line, the gain follows
from the estimate for the approximate antiderivative minus its leading
term and from the estimate for $a_I^{\rm c}$, each differentiated
$J$ times in the magnetic direction. It disappears when $J=j_1$.
Accordingly, Section~\ref{sec:endpoint-stress} uses the estimate
without this gain. In the second line, the loss comes from
$D_{B,I}\mathfrak a_I$ at magnetic order $J+1$. At $J=j_1$ the
prefactor absorbs it, because
$\tau_a\lambda_\parallel\varepsilon_\tau^{-1}=\tau_c\lambda_\parallel\le1$.
For $j=1$, we use \eqref{principal:principal-constant-bound} for the
material Lie derivative of the principal potential and the source bound
for its magnetic Lie derivative.
Applying \eqref{principal:scalar-action} once more gives
$j=2$ and requires one additional coefficient derivative. In both
applications, we first estimate the Lie derivative in the basis at
the original chart scales
\[
 \lambda_{q+1}^{-1}\partial_{y_I^{k_I}},\qquad
 \ell\partial_{y_I^{\nu_I}},\qquad \ell\partial_{y_I^{\zeta_I}},
\]
whose normalized coefficients have size
$C\tau_a\ell^{-1}\delta_{q+1}^{1/2}$. We then estimate the resulting expressions at the rescaled spatial
frequencies, with the same factor
$\tau_a^j\ell^{-j}\delta_{q+1}^{j/2}$ in the amplitude. Since curl
commutes with these Lie derivatives, the same argument controls the
field variations. Section~\ref{sec:endpoint-stress} uses these bounds
for the second variation and covariance, and the nonzero fast Fourier
modes in \eqref{principal:first-lie-mean} for the first variation.

\subsection{Background gaps}\label{ssec:background-gap-fields}

For the initial gaps $G_I^v=v_q-v_{\ell,I}$ and
$G_I^B=B_q-B_{\ell,I}$, take the half-sum and half-difference in
\eqref{mom:transported-gaps} and \eqref{adapt:full-block-gap}.
Together with the principal increments relative to the local
background, this gives, on the full local support,
\begin{equation}
 \begin{aligned}
 v_s-v_s^{\rm c}
 &=\mathcal U_s^{\rm c}\left[
     \curl(f_{w,I,s}\,dy_I^{\nu_I})+(\mathcal U_s^{\rm p}-\IId)G_I^v\right],\\
 B_s-B_s^{\rm c}
 &=\mathcal U_s^{\rm c}\left[
     \curl(f_{b,I,s}\,dy_I^{\nu_I})+(\mathcal U_s^{\rm p}-\IId)G_I^B\right],
 \end{aligned}
\end{equation}
where $f_{w,I,s},f_{b,I,s}$ are the scalar potentials from
\eqref{principal:profile-equation}, before the corrector pushforward.
It remains to estimate the two gap terms. We retain the smaller
magnetic amplitude and include the factor from the large spatial
differential of the principal flow.

For $G=G_I^v$ or $G=G_I^B$, write $G_s,G_s^{\rm c}$ and
$\Theta_s^{\mathrm p,\mathrm{gap}}$ for the corresponding half-sums
or half-differences of the signed comparisons and potentials above.

\begin{proposition}[Local gap bounds]\label{gap:field-bounds}
Fix an index $I$ and let $G$ denote either $G_I^v$ or $G_I^B$, satisfying
\[
 G\in\mathcal C_{\rprep,j_1}(A_G;\ell^{-1},\mathrm a)
 \cap\mathcal K_{\rgood+3,j_1}(A_G')
\]
in chart components with respect to the local background transports, with amplitudes
\[
 \begin{gathered}
 A_{G_I^v}=C\varepsilon_{q+1}^2\delta_q^{1/2},\qquad
 A_{G_I^B}=C\varepsilon_{q+1}^2\delta_{B,q}^{1/2},\\
 A_{G_I^v}'=A_{G_I^B}'=C\delta_q^{1/2}.
 \end{gathered}
\]
Assume also the following bounds for ordinary derivatives of $G$
in integer norms:
\[
 \|\partial_t^hG\|_r\le C A_G'\ell^{-r-h},
 \qquad h\le j_1+1,\quad r+h\le\rgood+4,
\]
in Cartesian and chart components.
Then, with respect to the corrected transports and uniformly in $s$,
\begin{align}
 \Theta_s^{\mathrm p,\mathrm{gap}}
 &\in\mathcal C_{\Nc}(C\tau_a\delta_{q+1}^{1/2}A_G;\lambda_{q+1},\mathrm f),
                     \label{adapt:gap-natural}\\
 G_s-G_s^{\rm c}
 &\in\mathcal C_{\Nc}(C\tau_a\lambda_{q+1}\delta_{q+1}^{1/2}A_G;\lambda_{q+1},\mathrm f),
                     \label{adapt:gap-field-natural}\\
 G_s&\in\mathcal C_{\Nc}(C(1+\tau_a\lambda_{q+1}\delta_{q+1}^{1/2})A_G;\lambda_{q+1},\mathrm f),
                     \label{adapt:full-gap-natural}
\end{align}
The same estimates hold at the local derivative costs for
compositions of at most $j_1$ material and magnetic derivatives,
and in the ordinary classes
$\mathcal O_{\rgood+2,j_1+1}(\,\cdot\,;\lambda_{q+1})$ with the same
amplitudes. The coefficient losses are $A_G'/A_G$, with lossy estimates
through order $\rgood+3$. In addition, the following estimates in integer norms hold for
$h\le j_1+1$ and $r+h\le\rgood+4$:
\[
 \begin{aligned}
 \|\partial_t^h(G_s-G_s^{\rm c})\|_r
 &\le C\tau_a\lambda_{q+1}\delta_{q+1}^{1/2}A_G
                                      \lambda_{q+1}^{r+h},\\
 \|\partial_t^hG_s\|_r
 &\le C(1+\tau_a\lambda_{q+1}\delta_{q+1}^{1/2})A_G
                                      \lambda_{q+1}^{r+h}.
 \end{aligned}
\]
For the magnetic field $A_G=A_{G_I^B}$ throughout.
\end{proposition}
\begin{proof}
By Corollary~\ref{prep:mixed-comparison}, the gaps have the stated
amplitudes at frequency $\lambda_q$ and costs $\mathrm c$, through
order $\rprep+1$ and with at most $j_1+1$ material and magnetic
derivatives in all. Since $\lambda_q\le\ell^{-1}$ and
$\tau_c^{-1}\le\tau_a^{-1}$, these imply the
$\mathcal C_{\rprep,j_1}$ and $\mathcal K_{\rgood+3,j_1}$ hypotheses
in the statement. The estimates hold in both Cartesian and chart
components: by Section~\ref{ssec:material-partition}, multiplication
by the chart matrices preserves the sharp bounds through order
$\rprep+1$ and the $\mathcal K_{\rgood+3,j_1}$ bound.

In \eqref{adapt:prepared-gap-potential}, the factors $dy_I^{\nu_I}$
and $df_I$ are invariant and the scalar sources satisfy
\[
 |G\cn f_I|\le C\tau_a\delta_{q+1}^{1/2}|G|,\qquad
 |G\cn y_I^{\nu_I}|\,|df_I|\le C\tau_a\delta_{q+1}^{1/2}|G|.
\]
Apply
\eqref{principal:structured-duhamel}--\eqref{principal:structured-coefficient-bound}
with $E_1=dy_I^{\nu_I}$, $E_2=df_I$ and scalar coefficients
$G\cn f_I$, $-G\cn y_I^{\nu_I}$. The product rule and the bounds
for derivatives of $df_I$ give \eqref{adapt:gap-natural}. In the
second term, only the product with $df_I$ is extended by zero;
the scalar coefficient need not be compactly supported.

For the fields, apply Corollary~\ref{principal:vector-transport} to
the chart components of $G$. To estimate $\|\cdot\|_{N,J}$, its vector
pushforward formula uses derivatives of $G$ only for $k+m\le J$ and
$r+k+m\le N$. The extra derivative of the flow differential falls on
the smooth
approximate antiderivative. By \eqref{principal:scalar-action}, an
estimate of order $N$ uses at most $N+2$ spatial derivatives of
$\mathfrak a_I$, with the same number of material and magnetic
derivatives. If $Y_s$ denotes the principal flow in the coordinates
$y_I$, the antiderivative estimates give
\begin{equation}
 \sup_s\|D_{y_I}Y_s-\IId\|_{N,J}
       \le C\tau_a\lambda_{q+1}\delta_{q+1}^{1/2}
\end{equation}
by \eqref{principal:phase-jacobian-bound}, also in the chart norms at the rescaled spatial frequencies
since $(\vartheta\ell)^{-1}\le\lambda_{q+1}$.
Thus the two terms in \eqref{adapt:gap-stretching-composition}
have amplitudes $C\tau_a\lambda_{q+1}\delta_{q+1}^{1/2}A_G$ and
$CA_G$, with the stated derivative weights. Their sum proves
\eqref{adapt:gap-field-natural} and \eqref{adapt:full-gap-natural},
since $\tau_a\lambda_{q+1}\delta_{q+1}^{1/2}\ge1$ by
\eqref{setup:rate-comparisons}. Apply
\eqref{gg:corrector-class-ranges} and the chart estimates following
\eqref{gg:corrected-chart-bounds} to return to Cartesian components
after the corrector pushforward.

We now prove the integer norm estimates. Fix
$1\le r+h\le\rgood+4$ and $h\le j_1+1$.
Subtracting the local fields from the old fields, and applying
\eqref{prep:old-ordinary} and \eqref{prep:whole-patch-fields}, gives
\[
 \|\partial_t^hG\|_r
 \le C\delta_q^{1/2}
   \bigl(\lambda_q^{r+h}
      +\ell^{-\alpha}\lambda_q\ell^{-(r+h-1)}\bigr)
 \le C\delta_q^{1/2}\ell^{-r-h}.
\]
The last inequality uses $\lambda_q\le\ell^{-1}$ and
$\ell^{-\alpha}\ell\lambda_q\le
\varepsilon_{q+1}^{\gamma_\ell-\gamma_S}\le1$, since
$\gamma_S<\gamma_\ell$ by \eqref{iter:small-gammaS}.
At order zero, use \eqref{prep:local-gap-good}. The chart matrices and
their inverses satisfy lossy estimates for all spatial derivatives
and ordinary time derivatives through order $j_1+1$. The product
and chain rules therefore give the same amplitude in chart components,
for both ordinary physical time and $D_{t,I}=\partial_t|_{y_I}$.
Indeed, a derivative of order $r+h$ in Lagrangian coordinates uses
derivatives of $G$ whose spatial and time orders sum to at most $r+h$,
with at most $h$ ordinary time derivatives. All other factors are
derivatives of the smooth chart or background velocity. This proves
the gap estimates through order $\rgood+4$ in these coordinates.

For the sharp estimates, apply \eqref{prep:local-gap-good} with
$m=0$ and use the coordinate frame pushed forward by the velocity
flow. This gives amplitude $A_G$ for either gap through derivative
order $\rprep$, with time order at most $j_1+1$. For higher orders,
apply Lemma~\ref{app:finite-spatial-reserve} in the rescaled integer
norms before transporting by the principal flow. When $r+h>\rprep$,
the ratio of the lossy bound to the required bound is at most
\[
 \frac{A_G'}{A_G}(\ell\Lambda_1)^{-(r+h)}
 \le\frac{A_G'}{A_G}(\ell\Lambda_1)^{-(\rprep+1)}\le C.
\]
Indeed, with the transverse scales in \eqref{principal:continuation-data},
we have $(\ell\Lambda_1)^{-1}=\vartheta\le1$, and the two ratios of
lossy to sharp amplitudes are
\[
 \frac{A_{G_I^v}'}{A_{G_I^v}}=C\varepsilon_{q+1}^{-2},\qquad
 \frac{A_{G_I^B}'}{A_{G_I^B}}
 =C\varepsilon_{q+1}^{-2}\Bigl(\frac{\delta_q}{\delta_{B,q}}\Bigr)^{1/2}
 =C\varepsilon_{q+1}^{-2-(\gamma_a+\gamma_\parallel)/b},
\]
each bounded by $C\varepsilon_{q+1}^{-\gres}$ according to
\eqref{stress:loss-maximum}. By \eqref{stress:ceiling-minimum},
$\rprep+1\ge\rcut-j_1-1$, so the preceding quotient is at most
$C\varepsilon_{q+1}^{-\gres}\vartheta^{\rcut-j_1-1}$ and is bounded
by \eqref{iter:general-final-reserve}. A longitudinal or material
derivative gives the same conclusion: its normalization uses
$\lambda_{q+1}\ge\Lambda_1$, while its lossy cost is $\ell^{-1}$,
and hence it contributes at most $(\ell\lambda_{q+1})^{-1}\le\vartheta$.
Thus every distribution of the $r+h$ derivatives has the required
bound. Together with the direct sharp estimate through order
$\rprep$, this proves the rescaled estimate through $r+h=\rgood+4$.

Apply Corollary~\ref{principal:vector-transport} in the rescaled
Lagrangian norms. Composition requires no further derivatives of $G$;
the flow differential contributes
$1+\tau_a\lambda_{q+1}\delta_{q+1}^{1/2}$. Next expand
$\partial_t=D_{t,I}-v_{\ell,I}\cn$ to recover ordinary time
derivatives. This does not increase the number of derivatives of the
gap and uses time derivatives of the local velocity only through
order $h-1$. For the corrector pushforward, we need one further
spatial derivative of its flow and time derivatives through order $h$.

For completeness, at $h\le j_1+1$ and $r+h\le\rgood+4$ the
required gap derivatives have spatial and ordinary time orders
summing to at most $r+h$, with time order at most $h$; these are
provided by \eqref{prep:old-ordinary}. For the local velocity, the
sum is again at most $r+h$ and the time order is at most
$h-1\le j_1$. The chart matrices are needed through spatial order
$r+h+1$ and time order $h$. Their bounds, and those for the smooth
local velocity, follow from Sections~\ref{sec:preparation}
and~\ref{ssec:material-partition} at all these orders.

It remains to check the derivatives of the LDFs. With
$(N,J)=(r+h,0)$, Corollary~\ref{principal:vector-transport} uses
the normalized principal coefficients through $(r+h+1,0)$.
By their definition, this requires $\mathfrak a_I$ through spatial
order $r+h+2$ and ordinary time order $h$. These bounds follow
from \eqref{principal:primitive-classes} and the preceding extension
to higher derivatives, with the same sharp and lossy amplitudes,
up to fixed constants, after at most $j_1+1$ additional ordinary time
derivatives. Proposition~\ref{gg:physical} supplies the corrector
LDF bounds through spatial order $r+h+1$ and time order $h\le j_1+1$.
We have therefore proved the integer norm estimates through order
$\rgood+4$. Interpolating with one additional spatial derivative
gives the stated ordinary H\"older estimates. The only loss in
either argument is the loss of the initial gap $G$.

\end{proof}

For either $G=G_n^v$ or $G=G_n^B$, the corrector has the exact
endpoint formula
\begin{equation}\label{gg:unmollified-endpoint}
 \curl\int_0^1\mathcal U_s^{\rm c}\iota_Gd\Theta^{\rm c}\,\dd s
       =(\mathcal U_1^{\rm c}-\IId)G,
\end{equation}
since Cartan's identity identifies the curl of the integrand with
$-\mathcal U_s^{\rm c}\mathcal L_{\xi^{\rm c}}G$, whose integral is
the endpoint difference. Estimating the two terms on the right-hand side
separately uses the ordinary gap bounds of
Proposition~\ref{gap:field-bounds}, with $h\le j_1+1$ and
$r+h\le\rgood+4$.

\subsection{A common background}

Spatial Fourier multipliers do not preserve local supports. We
therefore fix one active local background on each time neighborhood
and estimate all supported terms relative to that background before
summing and applying the multiplier. For overlapping time intervals $n,j$ at the same
deformation parameter $s$, the leading comparison factors at
frequency $\lambda_{q+1}$ are
\begin{equation}\label{stress:common-window}
 \tau_a\lambda_{q+1}\|v_{\ell,n,s}^{\rm c}-v_{\ell,j,s}^{\rm c}\|_0,
 \qquad
 \frac{\lambda_{q+1}\|B_{\ell,n,s}^{\rm c}-B_{\ell,j,s}^{\rm c}\|_0}{\lambda_\parallel}.
\end{equation}
The first factor is bounded by the velocity gap estimate. For the
second, we use the improved magnetic estimate below. For tensor Lie
derivatives, we also need one spatial derivative of each difference.
In taking these derivatives, the reference background is fixed
throughout the time neighborhood.

\begin{lemma}[Magnetic gap comparison]\label{aniso:magnetic-gap-closure}
Under the improved magnetic comparison of
Corollary~\ref{prep:magnetic-ordinary-comparison}, the second factor in
\eqref{stress:common-window} is bounded, and the complete principal
magnetic gap is smaller than the bound required for the next magnetic
increment in the classes of Proposition~\ref{gap:field-bounds}.
\end{lemma}
\begin{proof}
By \eqref{prep:local-gap-good}, the magnetic gap has amplitude
$A_{G_I^B}=C\varepsilon_{q+1}^2\delta_{B,q}^{1/2}$. The corrector
pushforward preserves this bound and its material and magnetic Lie
derivative estimates. Consequently,
\[
 \frac{\lambda_{q+1}A_{G_I^B}}{\lambda_\parallel}
 =C\varepsilon_{q+1}^{\gamma_\ell}
   \varepsilon_{q+1}^{1+(\gamma_a-(b-1)\gamma_\parallel)/b}=o(1),
\]
and the same quotient bounds the gradient difference after one
additional spatial derivative. For the complete principal gap,
\eqref{adapt:full-gap-natural} gives
\[
 \frac{(1+\tau_a\lambda_{q+1}\delta_{q+1}^{1/2})A_{G_I^B}}
       {\delta_{B,q+1}^{1/2}}
 \le C\left[
       \varepsilon_{q+1}^{2-\beta-\frac{b-1}{b}(\gamma_a+\gamma_\parallel)}
       +\varepsilon_{q+1}^{\gamma_\ell}
          \varepsilon_{q+1}^{1+\gamma_a-\frac{b-1}{b}(\gamma_a+\gamma_\parallel)}\right]
       =o(1),
\]
since both exponents are positive. The same gain holds after material
and magnetic differentiation by \eqref{prep:local-gap-good} and
Lemma~\ref{setup:calculus}(ii), and the local derivative costs are
bounded by those of the next iterate by
\eqref{setup:rate-comparisons}. For ordinary derivatives, apply
\eqref{adapt:gap-stretching-composition} and
Proposition~\ref{gap:field-bounds} to the principal gap, and
\eqref{gg:unmollified-endpoint} to the corrector gap. At spatial order
$r$ and ordinary time order $h$, these formulas use derivatives of
$B_q$ only with $h\le j_1+1$ and $r+h\le\rgood+4$, as provided by
\eqref{prep:old-ordinary}.
\end{proof}

\subsection{Endpoint derivatives}

The new transports are
\[
 D_{t,q+1}=D_{t,*}+u\cn,\qquad D_{B,q+1}=D_{B,*}+b\cn,
\]
where $(v_*,B_*)$ is the selected corrected local reference and
$u=v_{q+1}-v_*$, $b=B_{q+1}-B_*$. These identities apply to
Cartesian components; for tensors we use the corresponding Lie
derivatives. We first estimate the increments relative to the
reference fields and then apply these identities.

\begin{proposition}[Endpoint derivatives]
\label{full:endpoint-derivative-closure}
Assume \eqref{iter:general-final-reserve}. The
constructed fields satisfy the split field bounds at stage $q+1$. In
particular, with respect to the transports $D_{t,*},D_{B,*}$,
\begin{equation}\label{full:endpoint-increment-bounds}
 \begin{gathered}
 u\in\mathcal C_{\rgood}(C\delta_{q+1}^{1/2};\lambda_{q+1},\mathrm f)
 \cap\mathcal O_{\rgood+2,j_1+1}(C\delta_{q+1}^{1/2};\lambda_{q+1}),\\
 b\in\mathcal C_{\rgood}(C\delta_{B,q+1}^{1/2};\lambda_{q+1},\mathrm f)
 \cap\mathcal O_{\rgood+2,j_1+1}(C\delta_{B,q+1}^{1/2};\lambda_{q+1}),
 \end{gathered}
\end{equation}
and the $\mathcal C_{\rgood}$ estimates hold also with respect to
$D_{t,q+1},D_{B,q+1}$, up to fixed multiplicative constants, for
$k+m\le j_1$ and $r+k+m\le\rgood$. For sufficiently large $a$, the
constant $C$ in
\eqref{full:endpoint-increment-bounds} depends only on the fixed
orders, profiles and geometric data, independently of the inductive
constant $C_0$. The selected reference is
$(\lambda_{q+1},\mathrm f)$-adapted through $(\infty,j_1)$:
\begin{equation}\label{full:endpoint-gradient-bounds}
 \begin{aligned}
 \nabla v_*&\in\mathcal C_{\infty}(C\lambda_{q+1}\delta_{q+1}^{1/2};\lambda_{q+1},\mathrm f),\\
 \nabla B_*&\in\mathcal C_{\infty}(C\lambda_{q+1}\delta_{B,q+1}^{1/2};\lambda_{q+1},\mathrm f),
 \end{aligned}
\end{equation}
and the new background $(v_{q+1},B_{q+1})$ is adapted through
$(\rgood,j_1)$. Consequently, by Lemma~\ref{setup:calculus}(iii),(iv),
the bounds defining a Lie or transport derivative class
$\mathcal C_{N',J}(E;\lambda_{q+1},\mathrm f)$ with $N'\le\rgood$ and
$J\le j_1$ are equivalent up to fixed multiplicative constants with
respect to the corrected reference, the selected reference and the
actual endpoint transports.
\end{proposition}
\begin{proof}
We use the field estimates proved above. Steps~1--3 establish the
bounds for the increments and gradients; Step~4 proves that $C$
is independent of $C_0$ and gives the inductive constants
$C_0^{k+m+1}$.

\emph{1. The magnetic size.}
The principal magnetic increment relative to the local background has size
$C\tau_a\lambda_\parallel\delta_{q+1}^{1/2}$, and its quotient by the
required bound for the next magnetic increment is
\begin{equation}\label{full:field-reserve-quotient}
 \tau_a\lambda_\parallel\Bigl(\frac{\delta_{q+1}}{\delta_{B,q+1}}\Bigr)^{1/2}
 =1.
\end{equation}
Integrating \eqref{gg:common-generator-bounds} along the corrector
path gives the same bound for the corrector magnetic increment,
with the additional factor
$\tau_a\lambda_q\ell^{-\alpha}\delta_{q+1}^{1/2}$. The full magnetic
gap is smaller by Lemma~\ref{aniso:magnetic-gap-closure}.
Together these prove the bound for the magnetic increment and its
first spatial derivative.

\emph{2. Material and magnetic derivatives of increments relative to the local reference.}
Relative to the corrected reference, $u,b$ consist of the principal
fields \eqref{principal:fields-bound}, the principal gap fields
\eqref{adapt:gap-field-natural}, and the transported initial gaps
$\mathcal U_1^{\rm c}G_n^v$ and $\mathcal U_1^{\rm c}G_n^B$; the
corrector is already included in the reference. Apply
Propositions~\ref{gg:physical},
\ref{principal:nonprincipal-field-smallness} and
\ref{gap:field-bounds}, Lemma~\ref{gg:transport-classes}, and Step~1.
At the derivative costs of the next iterate, each velocity term has
amplitude $C\delta_{q+1}^{1/2}$ and each magnetic term has amplitude
$C\delta_{B,q+1}^{1/2}$ through at least order $\Nc$. For the
antiderivative error, this follows from the factor
$\varepsilon_\tau^{\min\{j_0,j_1-m\}}$ for $k+m\le j_1$.

Each term is a finite sum of slow coefficients multiplied by profiles
of the transported coordinate. Apply Lemma~\ref{setup:calculus}(vii),
using the losses in Section~\ref{ssec:derivative-closure}, with
$\Lambda_1=\lambda_{q+1}$ or, for transport by the principal flow,
with \eqref{principal:continuation-data}. In either case its
hypothesis is \eqref{iter:general-final-reserve}. The sharp estimates
therefore hold up to the spatial orders allowed by the lossy
estimates, with one additional derivative for summing in the charts.
For smooth coefficients this includes every spatial order; for the
gaps it includes $r+k+m\le\rgood+2$ by \eqref{prep:local-gap-good}.
In both cases, $k+m\le j_1$. This proves
\eqref{full:endpoint-increment-bounds} for Lie
derivatives relative to the corrected reference.

To pass to transport derivatives, we need the reference gradient
bounds. By \eqref{gg:corrected-adaptedness}, the reference is
$(\lambda_q,\mathrm f)$-adapted through $(\rprep-12,j_1)$, hence
also $(\lambda_{q+1},\mathrm f)$-adapted there. For higher spatial
orders, apply Lemma~\ref{setup:calculus}(vii) with
$\Lambda_1=\lambda_{q+1}$ to the bounds
\eqref{prep:local-gradient} and \eqref{prep:local-field-classes}:
\[
 \begin{aligned}
 \nabla v_{\ell,n}&\in
 \mathcal C_{\rprep}(C\lambda_q\delta_q^{1/2};\lambda_q,\mathrm c)
 \cap\mathcal K(C\ell^{-\alpha}\lambda_q\delta_q^{1/2}),\\
 \nabla B_{\ell,n}&\in
 \mathcal C_{\rprep}(C\lambda_q\delta_{B,q}^{1/2};\lambda_q,\mathrm c)
 \cap\mathcal K(C\ell^{-\alpha}\lambda_q\delta_q^{1/2}),
 \end{aligned}
\]
with respective losses $\ell^{-\alpha}$ and
$\ell^{-\alpha}(\delta_q/\delta_{B,q})^{1/2}$. Apply the same lemma to
$\nabla c_{v,n,1},\nabla c_{B,n,1}$, whose classes are
\eqref{gg:corrector-comparison-classes}, with respective losses
$\varepsilon_\tau^{-1}(\ell\lambda_q)^{-2}$ and
$\varepsilon_\tau^{-2}\varepsilon_{q+1}^{-\gamma_\parallel}
 (\ell\lambda_q)^{-2}$.
Each loss is at most $\varepsilon_{q+1}^{-\gres}$, and all four sharp
estimates hold through order $\rcut$. The assumed inequality therefore
gives
\[
 \begin{aligned}
 \nabla v_{\ell,n}&\in\mathcal C_\infty(C\lambda_q\delta_q^{1/2};\lambda_{q+1},\mathrm f),\\
 \nabla B_{\ell,n}&\in\mathcal C_\infty(C\lambda_q\delta_{B,q}^{1/2};\lambda_{q+1},\mathrm f),\\
 \nabla c_{v,n,1}&\in\mathcal C_\infty(C\ell^{-\alpha}\tau_a\lambda_q^2\delta_{q+1};\lambda_{q+1},\mathrm f),\\
 \nabla c_{B,n,1}&\in\mathcal C_\infty(C\ell^{-\alpha}\tau_a^2\lambda_\parallel\lambda_q^2\delta_{q+1};\lambda_{q+1},\mathrm f).
 \end{aligned}
\]
The comparisons of derivative costs \eqref{gg:corrector-rate-comparisons} now give
\eqref{full:endpoint-gradient-bounds}. Lemma~\ref{setup:calculus}(iii)
converts the increment classes defined by Lie derivatives into classes defined by transport derivatives.

\emph{3. Changing transports.}
Apply Lemma~\ref{aniso:old-frame-operator-change}, with
$\Lambda=\lambda_{q+1}$, to the Cartesian components of the
increments, using \eqref{full:endpoint-increment-bounds} and
\eqref{full:endpoint-gradient-bounds}. The required material and
magnetic orders sum to at most $j_1-1$, and the spatial, material
and magnetic orders sum to at most $\rgood-1$.
The lemma expresses derivatives in $D_{t,q+1},D_{B,q+1}$ in terms
of derivatives relative to the reference fields. This expression is
linear in the differentiated field, so the magnetic amplitude is
preserved.

Apply the same expansion to $\nabla v_*+\nabla u$ and
$\nabla B_*+\nabla b$, with spatial, material and magnetic orders
summing to at most $\rgood-1$ and material and magnetic orders
summing to at most $j_1-1$. The estimates relative to the reference
follow from its gradient bounds and one spatial derivative of the
increments. Thus the new gradients have amplitudes
$C\lambda_{q+1}\delta_{q+1}^{1/2}$ and
$C\lambda_{q+1}\delta_{B,q+1}^{1/2}$, respectively, proving
adaptedness through $(\rgood,j_1)$. Lemma~\ref{setup:calculus}(iii)
now gives the Lie derivative estimates, and part (iv) compares them
in both directions with the reference estimates.

\emph{4. The complete fields.}
The undifferentiated sizes of $v_*$ and $B_*$ need not be small, so we apply
the preceding expansion to first derivatives. The exact identities are
\begin{equation}\label{full:first-mixed-expansion}
 \begin{aligned}
 D_{t,q+1}v_{q+1}&=D_{t,*}v_*+u\cn v_*+D_{t,*}u+u\cn u,\\
 D_{B,q+1}v_{q+1}&=D_{B,*}v_*+b\cn v_*+D_{B,*}u+b\cn u,\\
 D_{t,q+1}B_{q+1}&=D_{t,*}B_*+u\cn B_*+D_{t,*}b+u\cn b,\\
 D_{B,q+1}B_{q+1}&=D_{B,*}B_*+b\cn B_*+D_{B,*}b+b\cn b.
 \end{aligned}
\end{equation}
The background fields on the right-hand sides are all differentiated,
so their terms can be estimated by the inductive derivative bounds.
The magnetic amplitude is retained: $u\cn B_*$ uses the bound for
$DB_*$, and $b\cn v_*$ uses the amplitude of $b$.
Since the induction equation is preserved,
$D_{t,q+1}B_{q+1}=D_{B,q+1}v_{q+1}$, and both sides require the
bound $\lambda_{q+1}(\delta_{q+1}\delta_{B,q+1})^{1/2}$.
For the four old-field contributions, the ratios to their respective
new inductive bounds are
\[
 \varepsilon_{q+1}^{1-2\beta},\qquad
 \varepsilon_{q+1}^{1-2\beta-\frac{b-1}{b}(\gamma_a+\gamma_\parallel)},\qquad
 \varepsilon_{q+1}^{1-2\beta-\frac{b-1}{b}(\gamma_a+\gamma_\parallel)},\qquad
 \varepsilon_{q+1}^{1-2\beta-2\frac{b-1}{b}(\gamma_a+\gamma_\parallel)},
\]
and all tend to zero because
$\beta+\frac{b-1}{b}(\gamma_a+\gamma_\parallel)<1/2$ by
\eqref{iter:exponent-window}. This proves the magnetic derivative
bound at spatial order zero, with the amplitude for the differentiated
magnetic field. For a product of an increment with a background
gradient, the corresponding ratio is bounded by
$\lambda_q\delta_q^{1/2}/(\lambda_{q+1}\delta_{q+1}^{1/2})$ or
$\lambda_q\delta_{B,q}^{1/2}/(\lambda_{q+1}\delta_{B,q+1}^{1/2})$,
multiplied by the bounded increment quotient from Step~1. Each
further material or magnetic derivative relative to the reference
fields contributes its old cost divided by its new cost.
Corrector terms satisfy the same estimates, with the additional
small factor in Step~1.

For the leading velocity interaction, the polarization is tangent
to the surfaces of constant phase. Since the local supports are
disjoint, we obtain
\[
 \|\sum_I w_I^{\rm lead}\cn w_I^{\rm lead}\|_r
 \le C\ell^{-1}\delta_{q+1}\lambda_{q+1}^r.
\]
Indeed, the tangent vector differentiates only the slow amplitude
or frame. The curl correction has size
$(\ell\lambda_{q+1})^{-1}\delta_{q+1}^{1/2}$, which compensates
for a derivative of frequency $\lambda_{q+1}$ on the other factor.
Every other product contains the small local velocity remainder,
a magnetic increment, or a gap. Using their previously proved
amplitudes and the derivative costs of the next iterate gives
\begin{equation}
 \|D_{B,q+1}u\|_0
 \le C\lambda_{q+1}(\delta_{q+1}\delta_{B,q+1})^{1/2},\qquad
 \|D_{t,q+1}u\|_0\le C\lambda_{q+1}\delta_{q+1}.
\end{equation}
Differentiate the four identities repeatedly and apply
Lemma~\ref{setup:calculus}(i),(ii). This proves the field estimates
for $1\le k+m\le j_1$ and $r+k+m\le\rgood$, retaining the smaller
amplitude in products containing $b$.
At the highest magnetic derivative orders of the antiderivative
remainder, the gain from the ratio of old to new magnetic derivative
costs compensates for the loss of local smallness, by
\eqref{aniso:finite-depth-budget}.

We now verify the dependence on $C_0$ required for
\eqref{prep:old-mixed} at stage $q+1$, with constants $C_0^{k+m+1}$.
The derivative orders, the profiles in Lemma~\ref{amplitude:profiles}
and the geometry in Lemma~\ref{prep:geometric-decomposition} are
fixed. The chart and frame bounds of
Section~\ref{ssec:material-partition} are also fixed; in particular,
the chart differential differs from a rotation by $o(1)$.
Constants depending on $C_0$ can therefore enter only through the
old inductive estimates. We check that these terms contain a
positive power of $\varepsilon_{q+1}$, so increasing $a$ gives
\eqref{full:endpoint-increment-bounds} with constants independent
of $C_0$.
For $\sum_I\curl\Theta_{w,I,1}^{\rm p}$ without derivatives, use
the fixed profile and frame bounds and $a_I^2\le C\delta_{q+1}$
from Lemma~\ref{prep:geometric-decomposition}.
Corollary~\ref{principal:phase-flow} bounds the principal scalar
flow, and \eqref{gg:physical-scales} bounds the corrector
differential, with normalized coefficients $o(1)$.
Thus the principal velocity has the required amplitude with a fixed constant.

For $D_{B,I}\mathfrak a_I$, consider the leading term
$\tau_a\alpha_I^{[1]}D_{B,I}a_I$. If the derivative falls on the
spatial cutoff, use the fixed cutoff profiles, the bounded chart
differential, and $\|B_{\ell,I}\|_0\le2+o(1)$ from
\eqref{full:pointwise-induction}; the resulting normalized constant
is independent of $C_0$. If it falls on the stress or chart
matrices, apply the local background estimates. The chart bounds
or the ratios of old to local derivative costs then give a positive
power in every term whose constant depends on $C_0$. The other
terms in the approximate antiderivative have an additional
$\varepsilon_\tau$. The principal magnetic amplitude is therefore
$C\delta_{B,q+1}^{1/2}$ with fixed $C$.

By Steps~1 and~2, the corrector fields, gap fields, transported
mollification gaps and $\sum_I(w_{I,1}-w_I^{\rm lead})$ have
amplitudes $o(1)\delta_{q+1}^{1/2}$ or $o(1)\delta_{B,q+1}^{1/2}$.
Here $o(1)$ is a positive power of $\varepsilon_{q+1}$ times a
constant that may depend on $C_0$.

After differentiation, a spatial oscillatory profile contributes
$\lambda_{q+1}$ with a fixed constant. A derivative of a slow
coefficient, frame, temporal profile or local background gradient,
normalized by the costs of the next iterate, contributes at most
$(\ell\lambda_{q+1})^{-1}$, $\varepsilon_\tau$ or
$\lambda_q\delta_q^{1/2}/(\lambda_{q+1}\delta_{q+1}^{1/2})$.
These are positive powers of $\varepsilon_{q+1}$, with constants
that may depend on $C_0$. Taking $a$ large absorbs all such terms
and proves the asserted independence of $C_0$ in
\eqref{full:endpoint-increment-bounds}.

Finally, expand $D_{t,q+1}^kD_{B_{q+1}}^mv_{q+1}$ and
$D_{t,q+1}^kD_{B_{q+1}}^mB_{q+1}$ using
\eqref{full:first-mixed-expansion} and Leibniz' rule. There are
finitely many terms, with combinatorial constants depending only on
$k+m\le j_1$. A term containing no reference field or its gradient
is a product of at most $k+m+1$ factors from $u$, $b$ and their
spatial, material and magnetic derivatives relative to the reference.
By Lemma~\ref{setup:calculus}(i),(ii), its ratio to the required
inductive bound is at most a product of $k+m+1$ fixed amplitude
constants, hence a fixed constant to the power $k+m+1$.

Every other term contains a spatial, material or magnetic derivative
of $v_*$ or $B_*$ relative to the reference fields. Dividing by the
required bound gives one of the positive powers of
$\varepsilon_{q+1}$ above. Further material or magnetic derivatives
improve the ratio by their old costs divided by their new costs;
the remaining constant may depend on $C_0$. Choose $C_0$ larger
than the fixed constants, and then choose $a$ so large that the sum
of these terms is a fixed small fraction of the inductive bound.
This proves \eqref{prep:old-mixed} at stage $q+1$ with constants
$C_0^{k+m+1}$, in the order prescribed in
Section~\ref{ssec:fields-pressure-support}. The tangency estimate
in fact makes the leading velocity interaction $o(1)$ relative to
its inductive bound, although the fixed bound already suffices for
$C_0^{k+m+1}$.

\emph{5. Ordinary derivative bounds.}
First consider the smooth constructed increments. Their membership
in $\mathcal O_{\infty,j_1+1}(\,\cdot\,;\lambda_{q+1})$, with the amplitudes
of Step~2, follows from Propositions~\ref{gg:physical},
\ref{principal:generator-estimates} and
\ref{principal:nonprincipal-field-smallness}, the ordinary estimates in
\eqref{principal:primitive-classes}, and Lemma~\ref{gg:transport-classes}.
In differentiating their source formulas, the Galbrun potential is
needed through time order $j_1+2$, and its first derivatives through
$j_1+1$; these estimates have already been proved.

For terms containing a principal gap, apply
Proposition~\ref{gap:field-bounds} in integer norms with
$h\le j_1+1$ and $r+h\le\rgood+4$. The vector pushforward after
rescaling uses gap derivatives whose spatial and ordinary time orders
sum to at most $r+h$, with time order at most $h$. It retains both
the magnetic amplitude and the single factor from the large spatial
differential. Hence Step~1 and \eqref{iter:principal-gap-budgets}
apply also when $r+h=\rgood+3$ or $r+h=\rgood+4$.

For the corrector gaps, fix one active local reference on a time
neighborhood. The global corrector flow gives exactly
\[
 \begin{aligned}
 v_1^{\rm c}-v_q
 &=\int_0^1\mathcal U_s^{\rm c}
                 \mathcal D_{t,n}\xi^{\rm c}\,\dd s
          +(\mathcal U_1^{\rm c}-\IId)G_n^v,\\
 B_1^{\rm c}-B_q
 &=\int_0^1\mathcal U_s^{\rm c}
                 \mathcal L_{B_{\ell,n}}\xi^{\rm c}\,\dd s
          +(\mathcal U_1^{\rm c}-\IId)G_n^B.
 \end{aligned}
\]
Estimate the source integrals by \eqref{gg:adjacent}; they contain
only smooth local background fields and the complete corrector LDF.
For the two endpoint differences, apply the integer norm argument
of Proposition~\ref{gap:field-bounds} at the oscillation scale.
The corrector differential is bounded, and its pushforward preserves
ordinary norms by Lemma~\ref{gg:transport-classes}. Bounding the
two endpoint values separately uses only $r+h$ derivatives of
$G_n^v,G_n^B$ and gives amplitudes $A_{G^v},A_{G^B}$, both smaller
than the required increment amplitudes. The extra spatial derivative
falls on the smooth corrector map, and the time order remains
$h\le j_1+1$. The selected reference is fixed on the time neighborhood
throughout this differentiation.

Summing the local increments, with their bounded overlap, now gives
\[
 \begin{aligned}
 \|\partial_t^h(v_{q+1}-v_q)\|_r
 &\le C_h\lambda_{q+1}^{r+h}\delta_{q+1}^{1/2},\\
 \|\partial_t^h(B_{q+1}-B_q)\|_r
 &\le C_h\lambda_{q+1}^{r+h}\delta_{B,q+1}^{1/2},
 \end{aligned}
\]
for $h\le j_1+1$ and $r+h\le\rgood+4$. When $r+h=\rgood+3$ or
$r+h=\rgood+4$, use the preceding integer norm estimates; all
derivatives of the
unmollified fields are provided by \eqref{prep:old-ordinary}.
When $r+h\ge1$, the ratios of the old ordinary bounds to the new
inductive bounds are $\varepsilon_{q+1}^{r+h-\beta}$ and
$\varepsilon_{q+1}^{r+h-\beta-\frac{b-1}{b}(\gamma_a+\gamma_\parallel)}$.
Both tend to zero, so adding the old fields proves the ordinary
inductive estimates through $\rgood+4$, with a fixed fraction of
the inductive constants to spare for sufficiently large $a$.

To check the constants, expand $\partial_t=D_{t,*}-v_*\cn$.
The undifferentiated reference velocity satisfies
$\|v_*\|_0\le\|v_q\|_0+o(1)\le2+o(1)$ by
\eqref{full:pointwise-induction}, independently of $C_0$.
Derivatives of the reference velocity and every other term with a
$C_0$-dependent constant have a positive power of
$\varepsilon_{q+1}$, as in Step~4. The derivative orders satisfy
$h\le j_1+1$ and $r+h\le\rgood+4$. Choose $C_0$ above the
corresponding fixed constants and then increase $a$ to obtain
\eqref{full:increment-bound} with constants $C_0^{h+1}$.
\end{proof}
\section{Estimates: Stresses}\label{sec:endpoint-stress}

In this section we estimate the stress defined in \eqref{mom:full-stress}.
We first treat the linear terms, then the quadratic terms, including
the antiderivative error and the magnetic square. We prove the
estimates for spatial, material, and magnetic derivatives separately
where their costs differ. We use $\Nc=\rprep-j_0-12$ from
\eqref{principal:ceilings}. By \eqref{stress:holder-unit},
$\lambda_{q+1}^{\alpha}$ bounds the H\"older factor
$\ell^{-\alpha}$; Lemma~\ref{setup:calculus}(v) introduces only fixed
powers of $\lambda_{q+1}^{\alpha}$. The parameter inequalities used
below are proved in Section~\ref{sec:parameter-choice}.

\subsection{The complete bound}

\begin{proposition}[Stress bounds]\label{stress:complete-bounds}
Choose the parameters as in Section~\ref{sec:parameter-choice}.
Then, for all sufficiently large $a$ and uniformly in $q$, the
stress defined by \eqref{mom:full-stress} satisfies
\begin{gather}
 \|R_{q+1}\|_r
 \le\lambda_{q+1}^{r-\alpha}\delta_{q+2},
       \qquad 0\le r\le\rgood,\label{full:stress-spatial}\\
 \|\partial_t^hR_{q+1}\|_r
 \le C_0^h\lambda_{q+1}^{r+h-\alpha}\delta_{q+2},
       \qquad h\le j_1,\quad r+h\le\rgood,\\
 \|D_{t,q+1}^kD_{B_{q+1}}^mR_{q+1}\|_r
 +\|\mathcal D_{t,q+1}^k\mathcal L_{B_{q+1}}^mR_{q+1}\|_r
 \le \lambda_{q+1}^{r-\alpha}
  (\lambda_{q+1}\delta_{q+1}^{1/2})^k
  (\lambda_{q+1}\delta_{B,q+1}^{1/2})^m\delta_{q+2},
       \label{full:stress-transport}
\end{gather}
where $1\le k+m\le j_1$ and $r+k+m\le\rgood-1$ in the last line.
\end{proposition}

We prove the proposition in
Sections~\ref{ssec:linear-terms}--\ref{ssec:derivative-closure}.
Recall the decomposition
\begin{equation}\label{stress:linear-quadratic-split}
 \begin{aligned}
 R^{\rm lin}={}&R^{\mathrm{cut}}+R^{\rm moll}
   +R^{\mathrm c}_{\rm tr}+R^{\mathrm c}_{\mathrm{na}}+R^{\mathrm c}_{\rm moll}
   +R^{\mathrm p}_{\rm tr}+R^{\mathrm p}_{\mathrm{na}}+R^{\mathrm p}_{\rm moll},\\
 R^{\rm quad}={}&R^{\mathrm c}_{\rm quad}+R^{\mathrm p}_{\rm chart}
   +R^{\mathrm p}_{\rm time}+R^{\mathrm p}_{\rm mag}
   +R^{\mathrm p}_{\rm quad},\qquad R_{q+1}=R^{\rm lin}+R^{\rm quad}.
 \end{aligned}
\end{equation}
The principal linear force is computed relative to the corrector
endpoint. It already contains the products between the corrector and
principal increments, so these products are absent from $R^{\rm quad}$.

We use the sharp and lossy estimates of
Section~\ref{sec:generator-field-estimates}. The \emph{size} of a
tensor will mean its sharp amplitude at the derivative costs of the
next iterate; its \emph{loss} will mean the ratio of its lossy and
sharp amplitudes. We also specify the derivative orders for which
the sharp estimate holds. Lemma~\ref{setup:calculus}(i),(ii) estimates
each product through the derivative range common to its factors.
The local derivative costs determine its amplitude, and the
comparison with the costs of the next iterate gives the estimates
used below.

For the oscillatory terms, we apply
Lemma~\ref{setup:calculus}(vii) to the slow coefficients of the
finite expansion in profiles of the transported chart coordinate.
We then apply Lemma~\ref{setup:calculus}(v) to the Fourier
multiplier. Its hypotheses hold because the corrected background is
$(\lambda_{q+1},\mathrm f)$-adapted through $(\infty,j_1)$ by
Proposition~\ref{full:endpoint-derivative-closure}. The spatial
operators and the corrector pushforward preserve these estimates
at the orders stated in Lemmas~\ref{setup:calculus}(v)
and~\ref{gg:transport-classes}. In
Section~\ref{ssec:derivative-closure} we compare these derivatives
with those of the endpoint background and complete the proof.
The following table lists the sizes to be estimated. The strict
parameter inequalities absorb their fixed powers of
$\lambda_{q+1}^{\alpha}$.
\begin{center}
\small
\setlength{\tabcolsep}{4pt}
\renewcommand{\arraystretch}{1.18}
\begin{tabular}{@{}>{\raggedright\arraybackslash}p{0.40\textwidth}>{\raggedright\arraybackslash}p{0.56\textwidth}@{}}
\toprule
Stress (definition) & Contribution and size estimate \\
\midrule
$R^{\mathrm{cut}},R^{\rm moll}$
 \eqref{mom:collar-stress}, \eqref{mom:mollification-stress}
 & Cutoff error and flow-smoothing part of $R^{\rm moll}$: $\varepsilon_\tau^{j_0}\delta_{q+1}$ \\
$R^{\rm moll}$ \eqref{mom:mollification-stress}
 & Spatial part: $\varepsilon_{q+1}^3\delta_{q+1}$ \\
$R^{\mathrm c}_{\rm tr},R^{\mathrm c}_{\mathrm{na}},R^{\mathrm c}_{\rm moll}$
 \eqref{mom:corrector-transport-stress}, \eqref{mom:corrector-interaction-stress}, \eqref{mom:G-gap}
 & $[\tau_a^2\lambda_q^2\delta_{q+1}
    +\varepsilon_\tau\varepsilon_{q+1}^{2+\gamma_\ell}]
       \delta_{q+1}$ \\
$R^{\mathrm p}_{\rm tr},R^{\mathrm p}_{\mathrm{na}}$
 \eqref{mom:principal-transport-stress}, \eqref{mom:principal-interaction-stress}
 & Fields relative to the local background: $\delta_{q+1}^{1/2}/(\lambda_{q+1}\tau_a)$ \\
$R^{\mathrm c}_{\rm quad},R^{\mathrm p}_{\rm chart}$
 \eqref{mom:corrector-quadratic-stress}, \eqref{mom:chart-stress}
 & $\tau_a^2\lambda_q^2\delta_{q+1}^2$ \\
$R^{\mathrm p}_{\rm time}$ \eqref{osc:primitive-magnetic-tensors}
 & $\varepsilon_\tau^{j_0}\delta_{q+1}$ \\
$R^{\mathrm p}_{\rm mag}$ \eqref{osc:primitive-magnetic-tensors}
 & $\delta_{q+1}(\tau_a\lambda_\parallel)^2$ \\
$R^{\mathrm p}_{\rm quad}$ \eqref{mom:N-remainder}
 & Oscillations and first Lie variation of the quadratic tensor (for the latter, see
 Proposition~\ref{stress:first-lie-variation}):
 $(\ell\lambda_{q+1})^{-1}\delta_{q+1}$ \\
$R^{\mathrm p}_{\rm quad}$ \eqref{mom:N-remainder}
 & Second variation and covariance: $\tau_a^2\ell^{-2}\delta_{q+1}^2$ \\
$R^{\mathrm p}_{\rm tr},R^{\mathrm p}_{\mathrm{na}},R^{\mathrm p}_{\rm moll},R^{\mathrm p}_{\rm quad}$
 \eqref{mom:principal-transport-stress}, \eqref{mom:principal-interaction-stress}, \eqref{mom:principal-gap-stress}, \eqref{mom:N-remainder}; \eqref{principal:gap-square}
 & Background differences:
 \newline $\varepsilon_{q+1}^2(\delta_q\delta_{q+1})^{1/2}
    (1+\tau_a\lambda_{q+1}\delta_{q+1}^{1/2})$ \\
\bottomrule
\end{tabular}
\end{center}
The complete amplitude contains $\varrho$. Derivatives of this
cutoff therefore contribute to $R^{\mathrm p}_{\rm time}$
and to the slow coefficients in the table. For $a_I^{\rm c}$, the
Galbrun cutoff error, and
$\mathfrak a_I-\tau_a\alpha_I^{[1]}a_I$, we use
\eqref{principal:primitive-classes} and
Theorem~\ref{galbrun:linear-theorem}(f) to retain the cancellation
factors after differentiation at the costs of the next iterate.
For the unmollified terms $R_q-R_\ell$ and $G_n$, we use the
derivative bounds in \eqref{prep:old-stress-time} and
\eqref{prep:local-gap-good}, respectively. Accordingly,
Section~\ref{ssec:derivative-closure} treats $R^{\rm moll}$
separately and uses the lossy estimates for the other terms,
including the additional spatial derivative required for the
H\"older bounds. This gives the stress estimates through order
$\rgood$.

\Needspace{6\baselineskip}
\subsection{Linear terms}\label{ssec:linear-terms}

\paragraph{Cutoff error.}
Applying Theorem~\ref{galbrun:linear-theorem}(d),(f) to the
localized Galbrun solution gives
\begin{equation}
 R^{\mathrm{cut}}\in
 \mathcal C_{\rprep-j_0-7}(C\ell^{-\alpha}\varepsilon_\tau^{j_0}\delta_{q+1};\lambda_q,\mathrm f)
 \cap\mathcal K(C\delta_{q+1})\cap\mathcal O_{\infty,j_1+1}(C\delta_{q+1};\ell^{-1}),
\end{equation}
with respect to the local background transports; the loss is $\varepsilon_\tau^{-j_0}$.

\paragraph{Smoothing.}
The definition of $\mathcal J_n$ gives the exact decomposition
\begin{equation}
 R_n^{\rm moll}=\varrho^2\eta_n^2
 \left[(R_q-R_\ell)
  +(\mathcal J_{s_t,s_B}^{(n)}-\IId)
                         (\delta_{q+1}\IId-R_\ell)\right].
\end{equation}
For the spatial part, the vanishing moments of the kernel and
\eqref{moll:spatial-accuracy} with $d=m_0$ give
\begin{equation}
 \|R_q-\mathcal J_\ell R_q\|_0
 \le C\ell^{m_0}\|R_q\|_{m_0}
 \le C\varepsilon_\ell\lambda_q^{-\alpha}\delta_{q+1}
 \le C\varepsilon_{q+1}^{3}\delta_{q+1},
\end{equation}
using $m_0\le\rgood$ spatial derivatives of $R_q$ and no time
derivatives. For the flow part,
\eqref{moll:mixed-accuracy} with $d=j_0$ applies to
$\delta_{q+1}\IId-R_\ell$: the required Lie derivatives, through
$j_0\le j_1$ transport derivatives, follow from
\eqref{prep:stress-good} and the symmetric-gradient identities
\[
 \mathcal D_{t,n}\IId=-2\sym\DD v_{\ell,n},\qquad
 \mathcal L_{B_{\ell,n}}\IId=-2\sym\DD B_{\ell,n}.
\]
For each term with $j_0$ material and magnetic derivatives in all,
the scale factors satisfy
\[
 s_t^is_B^j\tau_c^{-i}\lambda_\parallel^j=\varepsilon_\tau^{j_0},
 \qquad i+j=j_0.
\]
Thus the flow part has size $C\varepsilon_\tau^{j_0}\delta_{q+1}$.
Summing over the time intervals and using their bounded overlap gives
\begin{equation}\label{stress:mollification-zero}
 \|R^{\rm moll}\|_0\le C(\varepsilon_{q+1}^3+\varepsilon_\tau^{j_0})\delta_{q+1}.
\end{equation}
We will estimate the positive-order derivatives of $R^{\rm moll}$
in Section~\ref{ssec:derivative-closure}. In addition to the preceding
order-zero estimate, we need
\begin{equation}\label{stress:mollification-classes}
 R_n^{\rm moll}\in\mathcal C_{\rprep}(C\lambda_q^{-\alpha}\delta_{q+1};\lambda_q,\mathrm c)
 \cap\mathcal K_{\rgood-1,j_1}(C\lambda_q^{-\alpha}\delta_{q+1}),
\end{equation}
To prove the sharp estimate, apply \eqref{prep:stress-good} and
\eqref{moll:mixed-transfer} to the smoothed terms, and use
\eqref{prep:frame-comparison} and
Corollary~\ref{prep:mixed-comparison} to express derivatives of the
old stress relative to the local background. Differentiating a time
cutoff contributes $\tau_c^{-1}$. Flow smoothing commutes with the
local transports, and the commutators with $\mathcal J_\ell$
satisfy the same amplitude bounds.

For the lossy estimate, expand the local material and magnetic
derivatives into ordinary derivatives. Applied to $R_q$, a
composition whose spatial, material, and magnetic derivative orders
sum to $r+k+m$ involves ordinary spatial and time derivatives
whose orders sum to at most $r+k+m$, with at most $k\le j_1$
time derivatives. The remaining factors are derivatives of the local
fields and satisfy the lossy estimates. We then use
\eqref{prep:old-stress-time}, with one additional spatial derivative
to obtain the H\"older bound. Thus $r+k+m\le\rgood-1$, which
proves the second estimate.

\paragraph{Potential identities.}
Let $\Theta_w,\Theta_b$ be a local pair of remainder or increment
potentials. For the background $(v,B)$ associated with these potentials, the transport and Nash stresses are
\begin{align}
 R_{\rm tr}
 &=\mathcal R\curl(\mathcal D_t\Theta_w-\mathcal L_B\Theta_b),
                        \label{principal:local-transport-stress}\\
 R_{\mathrm{na}}
 &=2\mathcal R\ddiv(\Theta_w\times\nabla v-\Theta_b\times\nabla B).
                        \label{principal:local-interaction-stress}
\end{align}
Commutation of curl with the corresponding Lie derivatives gives the
first identity, and \eqref{mom:Kdiv} gives the second. In both cases,
the fixed spatial operator acts on the sum of the forces over the
torus, as prescribed by \eqref{stress:common-window}.

\paragraph{Corrector source.}
For a fixed local background $n$, put
\[
 \begin{gathered}
 F_{w,n}=\mathcal D_{t,n}\Theta_n^{\rm c},\qquad
 F_{b,n}=\mathcal L_{B_{\ell,n}}\Theta_n^{\rm c},\\
 G_n^v=v_q-v_{\ell,n},\qquad G_n^B=B_q-B_{\ell,n}.
 \end{gathered}
\]
Then \eqref{mom:remainder-potentials} gives
\[
 \begin{aligned}
 \Theta_{w,n}^{\mathrm c,\mathrm{rem}}&=(\mathcal P^{\rm c}-\IId)F_{w,n}
                         +\mathcal P^{\rm c}\mathcal L_{G_n^v}\Theta_n^{\rm c},\\
 \Theta_{b,n}^{\mathrm c,\mathrm{rem}}&=(\mathcal P^{\rm c}-\IId)F_{b,n}
                         +\mathcal P^{\rm c}\mathcal L_{G_n^B}\Theta_n^{\rm c},
 \end{aligned}
\]
and the exact source decomposition is
\begin{equation}\label{gg:five-terms}
 \begin{aligned}
 &\mathcal D_{t,n}\Theta_{w,n}^{\mathrm c,\mathrm{rem}}
          -\mathcal L_{B_{\ell,n}}\Theta_{b,n}^{\mathrm c,\mathrm{rem}}\\
 &\quad=(\mathcal P^{\rm c}-\IId)
                 (\mathcal D_{t,n}^2-\mathcal L_{B_{\ell,n}}^2)\Theta_n^{\rm c}\\
 &\qquad+[\mathcal D_{t,n},\mathcal P^{\rm c}]F_{w,n}
                       -[\mathcal L_{B_{\ell,n}},\mathcal P^{\rm c}]F_{b,n}\\
 &\qquad+\mathcal P^{\rm c}(\mathcal L_{G_n^v}F_{w,n}-\mathcal L_{G_n^B}F_{b,n})\\
 &\qquad+\mathcal P^{\rm c}(\mathcal L_{\mathcal D_{t,n}G_n^v}\Theta_n^{\rm c}
                     -\mathcal L_{\mathcal L_{B_{\ell,n}}G_n^B}\Theta_n^{\rm c})\\
 &\qquad+[\mathcal D_{t,n},\mathcal P^{\rm c}]\mathcal L_{G_n^v}\Theta_n^{\rm c}
       -[\mathcal L_{B_{\ell,n}},\mathcal P^{\rm c}]\mathcal L_{G_n^B}\Theta_n^{\rm c}.
 \end{aligned}
\end{equation}
Here the derivatives refer to the local background:
\[
 \mathcal D_{t,n}=\partial_t+\mathcal L_{v_{\ell,n}},\qquad
 \mathcal L_{B_{\ell,n}}G_n^B=[B_{\ell,n},G_n^B].
\]
In particular,
\[
 \mathcal D_{t,n}G_n^v=\partial_tG_n^v+[v_{\ell,n},G_n^v].
\]
The commutators with the path average are obtained from those with the
pushforward. For $\mathscr D=\mathcal D_{t,n}$ or
$\mathcal L_{B_{\ell,n}}$, Duhamel's formula gives
\begin{equation}
 [\mathscr D,\mathcal U_s^{\rm c}]
 =-\int_0^s\mathcal U_{s-u}^{\rm c}
          \mathcal L_{\mathscr D\xi^{\rm c}}\mathcal U_u^{\rm c}\,\dd u,
\end{equation}
Here $\xi^{\rm c}$, $\mathcal U_s^{\rm c}$, and
$\mathcal P^{\rm c}$ are formed from the sum of the correctors on
all time intervals. To estimate $\mathscr D\xi^{\rm c}$, commute
curl with the first Lie derivative of
$\Theta^{\rm c}=\sum_j\Theta_j^{\rm c}$. Proposition~\ref{gg:physical}
bounds this derivative of the potential. Applying
Lemma~\ref{gg:transport-classes} to the two pushforwards in Duhamel's
formula then bounds the commutator by a product with one spatial
derivative on one of its factors.

We now estimate \eqref{gg:five-terms}. We use
Corollary~\ref{prep:mixed-comparison} for the gaps and
Lemma~\ref{gg:transport-classes} for their pushforwards. For the Lie
derivative of a one-form, the required product estimate is
\[
 \|\mathcal L_XY\|_\alpha
 \le C\bigl(\|X\|_\alpha\|Y\|_{1+\alpha}
             +\|X\|_{1+\alpha}\|Y\|_\alpha\bigr).
\]
For the order-zero estimate, we use the following bounds, with $r=0,1,2$:
\[
 \begin{gathered}
 \|F_{w,n}\|_{r+\alpha}\le C\tau_a\lambda_q^r\ell^{-\alpha}\delta_{q+1},\qquad
 \|F_{b,n}\|_{r+\alpha}\le C\tau_a^2\lambda_q^r\ell^{-\alpha}\lambda_\parallel\delta_{q+1},\\
 \|\mathcal D_{t,n}\xi^{\rm c}\|_{r+\alpha}\le C\tau_a\lambda_q^{r+1}\ell^{-\alpha}\delta_{q+1},\qquad
 \|\mathcal L_{B_{\ell,n}}\xi^{\rm c}\|_{r+\alpha}\le C\tau_a^2\lambda_q^{r+1}\ell^{-\alpha}\lambda_\parallel\delta_{q+1},
 \end{gathered}
\]
Here \eqref{gg:adjacent} allows us to estimate the summed corrector
relative to a common background. Together with
\eqref{gg:high-mixed-state}, these bounds give
\[
 \|(\mathcal P^{\rm c}-\IId)(\mathcal D_{t,n}^2-\mathcal L_{B_{\ell,n}}^2)\Theta_n^{\rm c}\|_\alpha
 +\|[\mathcal D_{t,n},\mathcal P^{\rm c}]F_{w,n}\|_\alpha
 +\|[\mathcal L_{B_{\ell,n}},\mathcal P^{\rm c}]F_{b,n}\|_\alpha
 \le C\lambda_{q+1}^{(j_1+3)\alpha}\tau_a^2\lambda_q^2\delta_{q+1}^2.
\]
For the first term, apply \eqref{gg:high-mixed-state}, which
contributes $\lambda_{q+1}^{(j_1+1)\alpha}$, and estimate the
difference between the path average and the identity, which
contributes $\ell^{-\alpha}$. The H\"older estimate introduces
one further factor $\ell^{-\alpha}\le\lambda_{q+1}^{\alpha}$.
For the commutators, use the products of first Lie derivatives
above; the magnetic commutator contains the additional factor
$(\tau_a\lambda_\parallel)^2\le1$.

It remains to estimate the three terms containing a background gap.
Applying \eqref{prep:local-gap-good}, also to the transport
derivative of the gap in the second of these terms, gives
\[
 \begin{aligned}
 \|\mathcal P^{\rm c}(\mathcal L_{G_n^v}F_{w,n}-\mathcal L_{G_n^B}F_{b,n})\|_\alpha
 &\le C\tau_a\lambda_q\ell^{-2\alpha}\varepsilon_{q+1}^2\delta_q^{1/2}\delta_{q+1},\\
 \|\mathcal P^{\rm c}(\mathcal L_{\mathcal D_{t,n}G_n^v}\Theta_n^{\rm c}
                    -\mathcal L_{\mathcal L_{B_{\ell,n}}G_n^B}\Theta_n^{\rm c})\|_\alpha
 &\le C\varepsilon_\tau\tau_a\lambda_q\ell^{-2\alpha}\varepsilon_{q+1}^2\delta_q^{1/2}\delta_{q+1},\\
 \|[\mathcal D_{t,n},\mathcal P^{\rm c}]\mathcal L_{G_n^v}\Theta_n^{\rm c}
        -[\mathcal L_{B_{\ell,n}},\mathcal P^{\rm c}]\mathcal L_{G_n^B}\Theta_n^{\rm c}\|_\alpha
 &\le C\tau_a^3\lambda_q^3\ell^{-3\alpha}\varepsilon_{q+1}^2\delta_q^{1/2}\delta_{q+1}^2,
 \end{aligned}
\]
where the second line uses $\tau_a\tau_c^{-1}=\varepsilon_\tau$ and
$\tau_a\lambda_\parallel=\varepsilon_\tau\varepsilon_{q+1}^{\gamma_\parallel}\le\varepsilon_\tau$,
and the third has two corrector factors contributing
$\tau_a^2\lambda_q^2\ell^{-\alpha}\delta_{q+1}$. We sum the local
forces over $n$ before applying the fixed spatial operators. The
operator $\mathcal R\curl$ has order zero and contributes the factor
$\lambda_{q+1}^{(j_1+1)\alpha}$ by
Lemma~\ref{setup:calculus}(v); multiplying the remainder potential by
a background gradient gives the Nash term, and \eqref{mom:G-gap}
gives the comparison term, which contains no nonlocal Fourier multiplier. We obtain
\begin{equation}\label{stress:corrector-linear}
 \begin{aligned}
 \|R^{\mathrm c}_{\rm tr}\|_\alpha
 &\le C\lambda_{q+1}^{(2j_1+4)\alpha}\delta_{q+1}
       (\tau_a^2\lambda_q^2\delta_{q+1}
            +\tau_a\lambda_q\varepsilon_{q+1}^2\delta_q^{1/2}),\\
 \|R^{\mathrm c}_{\mathrm{na}}\|_\alpha
 &\le C\varepsilon_\tau\lambda_{q+1}^{(2j_1+4)\alpha}\delta_{q+1}
       (\tau_a^2\lambda_q^2\delta_{q+1}
            +\tau_a\lambda_q\varepsilon_{q+1}^2\delta_q^{1/2}),\\
 \|R^{\mathrm c}_{\rm moll}\|_\alpha
 &\le C\lambda_{q+1}^{\alpha}\tau_a\lambda_q
                      \varepsilon_{q+1}^2\delta_q^{1/2}\delta_{q+1}.
 \end{aligned}
\end{equation}
Since $\tau_a\lambda_q\delta_q^{1/2}=\varepsilon_\tau\varepsilon_{q+1}^{\gamma_\ell}$,
their sum obeys
\begin{equation}\label{gg:actual-stress}
 \|R^{\mathrm c}_{\rm tr}\|_\alpha+
 \|R^{\mathrm c}_{\mathrm{na}}\|_\alpha+
 \|R^{\mathrm c}_{\rm moll}\|_\alpha
 \le C\lambda_{q+1}^{(2j_1+4)\alpha}
 [\tau_a^2\lambda_q^2\delta_{q+1}
  +\varepsilon_\tau\varepsilon_{q+1}^{2+\gamma_\ell}]
                                                   \delta_{q+1}.
\end{equation}
Apply spatial, material, and magnetic derivatives to the same
identities and use the product rule. This gives
\eqref{stress:corrector-linear}, with the derivative costs of the
next iterate, when $k+m\le j_1$ and the spatial, material and
magnetic orders sum to at most $\rprep-14$. Indeed,
\eqref{gg:high-mixed-state} bounds
$(\mathcal D_{t,n}^2-\mathcal L_{B_{\ell,n}}^2)\Theta_n^{\rm c}$
through order $\rprep-10$. A Lie derivative along an LDF, a curl,
or a commutator with the path average uses one additional spatial
derivative.
The gap estimates hold through order $\rprep+1$, whereas the
corrector pushforward requires derivatives of its argument only
through order $\rprep-12$, and one more for
$\mathcal P^{\rm c}-\IId$. Thus all the differentiated terms
satisfy the asserted bounds. We prove the lossy estimates in
Section~\ref{ssec:derivative-closure}.

\paragraph{Principal transport.}
In the chart corresponding to $I$, before the corrector pushforward,
we apply the local material derivative $D_{t,I}$ to the velocity
equation in \eqref{principal:profile-equation} and subtract the local
magnetic derivative $D_{B,I}$ of the magnetic equation. This gives
\begin{equation}\label{principal:exact-transport-potential}
 \begin{aligned}
 (\partial_s+\xi_{I,0}^{\rm p}\cn)(D_{t,I}f_{w,I,s}-D_{B,I}f_{b,I,s})
 ={}&\lambda_{q+1}^{-1}(D_{t,I}^2-D_{B,I}^2)\mathfrak a_I\,
                           \varphi(\lambda_{q+1}y_I^{k_I})\\
 &-[D_{t,I},\xi_{I,0}^{\rm p}\cn]f_{w,I,s}
       +[D_{B,I},\xi_{I,0}^{\rm p}\cn]f_{b,I,s},
 \end{aligned}
\end{equation}
with initial value $D_{t,I}f_{w,I,0}-D_{B,I}f_{b,I,0}=0$.
Here $\xi_{I,0}^{\rm p}\cn$ acts on physical scalars. The local background
derivatives commute with the chart directions and annihilate $y_I^{k_I}$,
so differentiating the scalar equations introduces exactly the two
displayed commutators. The source contains the complete approximate
antiderivative,
\[
 D_{t,I}^2\mathfrak a_I=\tau_a^{-1}\alpha_I'(t/\tau_a)a_I
                 +\alpha_I(t/\tau_a)D_{t,I}a_I+D_{t,I}a_I^{\rm c},
 \qquad D_{B,I}^2\mathfrak a_I=D_{B,I}(D_{B,I}\mathfrak a_I),
\]
and the bounds for the approximate antiderivative imply
\[
 \begin{aligned}
 \|(D_{t,I}^2-D_{B,I}^2)\mathfrak a_I\|_0
 &\le C\tau_a^{-1}\delta_{q+1}^{1/2}
       [1+\varepsilon_\tau+\varepsilon_\tau^{j_0}
          +(\tau_a\lambda_\parallel)^2]\\
 &\le C\tau_a^{-1}\delta_{q+1}^{1/2}.
 \end{aligned}
\]
We next estimate compositions of up to $j_1$ material and
magnetic derivatives of the source. Apply
\eqref{principal:primitive-sharp} with the material and magnetic
orders summing to at most $j_1+2$, and use the costs of the next
iterate. For the term with two material derivatives, this gives
\[
 D_{t,I}^2\mathfrak a_I\in
 \mathcal C_{\rprep-j_0-2}
 (C\tau_a^{-1}\delta_{q+1}^{1/2};\lambda_q,\mathrm f)
\]
without additional derivative loss. For the magnetic term, the corresponding estimate is
\[
 D_{B,I}^2\mathfrak a_I\in
 \mathcal C_{\rprep-j_0-2}
 (C\tau_a\lambda_\parallel^2\delta_{q+1}^{1/2};\lambda_q,\mathrm f).
\]
For the two highest orders of magnetic differentiation, the loss and
the ratio of the magnetic derivative costs satisfy
\[
 \varepsilon_\tau^{-[m+2-j_1]^+}\le\varepsilon_\tau^{-2},\qquad
 \left(\frac{\lambda_\parallel}
 {\lambda_{q+1}\delta_{B,q+1}^{1/2}}\right)^m\le\varepsilon_\tau^m.
\]
For $m\le j_1-2$, there is no loss. For $m=j_1-1$ and
$m=j_1$, the losses are $\varepsilon_\tau^{-1}$ and
$\varepsilon_\tau^{-2}$, respectively. Since $j_1\ge2$, the ratio of
the derivative costs absorbs these factors in both cases, as in
\eqref{principal:primitive-classes}. The magnetic amplitude is no larger than
the amplitude of the material-derivative term, since
\[
 \tau_a^2\lambda_\parallel^2
 \le(\tau_c\lambda_\parallel)^2=\varepsilon_{q+1}^{2\gamma_\parallel}\le1.
\]
To estimate the two commutators in
\eqref{principal:exact-transport-potential}, use the chart norms
\eqref{principal:profile-norm}. The coefficients of $[D_{t,I},\xi_{I,0}^{\rm p}\cn]$ have size $C\ell^{-1}\delta_{q+1}^{1/2}$
and $f_{w,I,s}$ has size $C\lambda_{q+1}^{-1}\delta_{q+1}^{1/2}$, while
$[D_{B,I},\xi_{I,0}^{\rm p}\cn]$ has size $C\tau_a\ell^{-1}\delta_{q+1}^{1/2}\lambda_\parallel$
and $f_{b,I,s}$ has size $C\lambda_{q+1}^{-1}\tau_a\lambda_\parallel\delta_{q+1}^{1/2}$,
each of the latter already containing one magnetic derivative.
Applying the product rule for $m\le j_1$ further magnetic derivatives
introduces an additional loss of at most $\varepsilon_\tau^{-2}$.
Thus the ratio
of their product to the source amplitude is bounded by
\[
 C\tau_a\ell^{-1}\delta_{q+1}^{1/2}
   (\tau_a\lambda_\parallel)^2\varepsilon_\tau^{-2}
 \le C\tau_a\ell^{-1}\delta_{q+1}^{1/2},
\]
because
\[
 \tau_a\lambda_\parallel\varepsilon_\tau^{-1}
 =\tau_c\lambda_\parallel\le1.
\]
The spatial coefficient estimates in
\eqref{principal:scalar-action} and
\eqref{principal:primitive-classes} apply to the principal LDF in
these commutators. No further magnetic derivative is introduced.
We can therefore apply
\eqref{principal:structured-coefficient-bound} in the invariant
coframe $dy_I^{\nu_I}$, with the source estimates just proved.
This bounds
$(D_{t,I}f_{w,I,1}-D_{B,I}f_{b,I,1})dy_I^{\nu_I}$
relative to the local background.

The principal transport potential is its corrector pushforward.
Indeed, \eqref{principal:profile-equation} identifies
$\Theta_{w,I,1}^{\rm p}=\mathcal U_1^{\rm c}(f_{w,I,1}dy_I^{\nu_I})$ and
$\Theta_{b,I,1}^{\rm p}=\mathcal U_1^{\rm c}(f_{b,I,1}dy_I^{\nu_I})$, so
\eqref{direct:differentiated-transport} gives the exact identity
\[
 \begin{aligned}
 &(\partial_t+\mathcal L_{v_{\ell,I,1}^{\rm c}})\Theta_{w,I,1}^{\rm p}
      -\mathcal L_{B_{\ell,I,1}^{\rm c}}\Theta_{b,I,1}^{\rm p}\\
 &\qquad=\mathcal U_1^{\rm c}
              \bigl[(D_{t,I}f_{w,I,1}-D_{B,I}f_{b,I,1})dy_I^{\nu_I}\bigr].
 \end{aligned}
\]
The operators on the left are those of the corrected background, and
$D_{t,I},D_{B,I}$ on the right are those of the local background.
This identity incorporates the change induced by the corrector
pushforward. Lemma~\ref{gg:transport-classes} now yields
\begin{equation}\label{principal:transport-source-bound}
 \begin{aligned}
 \mathcal D_{t,I,1}^{\rm c}\Theta_{w,I,1}^{\rm p}-\mathcal L_{B_{\ell,I,1}^{\rm c}}\Theta_{b,I,1}^{\rm p}
 &\in\mathcal C_{\Nc}(C\lambda_{q+1}^{-1}\tau_a^{-1}\delta_{q+1}^{1/2};\lambda_{q+1},\mathrm f)\\
 &\quad\cap\mathcal O_{\infty,j_1+1}(C\lambda_{q+1}^{-1}\tau_a^{-1}\delta_{q+1}^{1/2};\lambda_{q+1}),
 \end{aligned}
\end{equation}
for compositions of at most $j_1$ material and magnetic derivatives
with respect to the corrected background. Before composition with the
principal flow, the slow coefficients in the source representation
belong to $\mathcal K(C\lambda_{q+1}^{-1}\ell^{-1}\delta_{q+1}^{1/2})$.
Indeed, the source
is $\lambda_{q+1}^{-1}$ times one material derivative of
$D_{t,I}\mathfrak a_I=\alpha_I a_I+a_I^{\rm c}\in\mathcal K(C\delta_{q+1}^{1/2})$
or one magnetic derivative of $D_{B,I}\mathfrak a_I$. The slow
coefficients used to estimate the commutator terms, before composition with the principal flow, belong
to the same $\mathcal K$ class by
\eqref{principal:primitive-classes}. The loss is therefore
$\tau_a\ell^{-1}$.

The identities \eqref{principal:local-transport-stress} and
\eqref{principal:local-interaction-stress} turn these potential estimates
into stress estimates. We sum the local forces relative to a common
corrected background before applying the spatial operators. For the
Nash stress, the force is the divergence of the product of a principal
potential and a background gradient. Lemma~\ref{setup:calculus}(v)
therefore gives
\begin{equation}\label{principal:natural-stresses}
 R^{\mathrm p}_{\rm tr}-R_{\rm tr}^{\rm gap},\quad R^{\mathrm p}_{\mathrm{na}}-R_{\mathrm{na}}^{\rm gap}\in
 \mathcal C_{\Nc}\Bigl(C\lambda_{q+1}^{(j_1+2)\alpha}
              \frac{\delta_{q+1}^{1/2}}{\lambda_{q+1}\tau_a};\lambda_{q+1},\mathrm f\Bigr)
 \cap\mathcal O_{\infty,j_1+1}\Bigl(C\lambda_{q+1}^{\alpha}\frac{\delta_{q+1}^{1/2}}{\lambda_{q+1}\tau_a};\lambda_{q+1}\Bigr),
\end{equation}
where the subtracted tensors are the gap contributions estimated
below, and the estimates are relative to a common corrected background.
For the Nash stress, the uncorrected background gradient has amplitude
$C\tau_c^{-1}\le C\tau_a^{-1}$ by \eqref{prep:local-gradient}, and the
corrector increment of the gradient has amplitude
$C\ell^{-\alpha}\tau_a\lambda_q^2\delta_{q+1}\le C\tau_a^{-1}$ by
\eqref{gg:corrector-comparison-classes}, so each gradient contribution has amplitude at most
$C\tau_a^{-1}$. Its product with the principal potential therefore
has the same size as the transport stress. For the lossy estimates,
\[
 \begin{aligned}
 \nabla v_{\ell,I},\nabla B_{\ell,I}
 &\in\mathcal K(C\lambda_q\ell^{-\alpha}\delta_q^{1/2}),\\
 \nabla c_{v,I,1},\nabla c_{B,I,1}
 &\in\mathcal K(C\tau_c\ell^{-2}\delta_{q+1}).
 \end{aligned}
\]
Their amplitudes satisfy
\[
 C\lambda_q\ell^{-\alpha}\delta_q^{1/2}\le C\tau_a^{-1},\qquad
 C\tau_c\ell^{-2}\delta_{q+1}\le C\tau_a^{-1}.
\]
The second comparison follows from
\[
 \tau_a\tau_c\ell^{-2}\delta_{q+1}
 =\varepsilon_\tau\delta_{q+1}/\delta_q
 =\varepsilon_{q+1}^{\gamma_a+2\beta}\le1.
\]
Thus, before composition with the principal flow, the slow coefficients
in the representation of the Nash product belong to
$\mathcal K(C\lambda_{q+1}^{-1}\ell^{-1}\delta_{q+1}^{1/2})$,
with the same loss $\tau_a\ell^{-1}$ as the transport force.

\paragraph{Principal background differences.}
We apply the same two operators to the gap potentials, again summing
the local forces first. The extra material derivative removes the
factor $\tau_a$ from their amplitude; a magnetic derivative retains its
smaller cost. Consequently Proposition~\ref{gap:field-bounds} yields
\begin{equation}\label{stress:principal-linear-gap}
 R_{\rm tr}^{\rm gap},\,R_{\mathrm{na}}^{\rm gap}\in
 \mathcal C_{\Nc}\bigl(C\lambda_{q+1}^{(j_1+2)\alpha}\delta_{q+1}^{1/2}\varepsilon_{q+1}^2\delta_q^{1/2};\lambda_{q+1},\mathrm f\bigr).
\end{equation}
The comparison tensor has the corrected old gap as its first factor
and the complete principal increment as its second, hence
\begin{equation}\label{stress:principal-comparison}
 R^{\mathrm p}_{\rm moll}\in
 \mathcal C_{\Nc}\bigl(C\lambda_{q+1}^{(j_1+2)\alpha}\varepsilon_{q+1}^2\delta_q^{1/2}
  \delta_{q+1}^{1/2}
       (1+\tau_a\lambda_{q+1}\varepsilon_{q+1}^2\delta_q^{1/2});\lambda_{q+1},\mathrm f\bigr).
\end{equation}
When differentiating these products, the derivative of the gap
can have one higher material or magnetic order than the stress
being estimated. We use \eqref{prep:local-gap-good} with the
material and magnetic orders summing to at most $j_1+1$.
Apply Lemma~\ref{setup:calculus}(vii) to $G$ with at most $j_1+1$
material and magnetic derivatives in all. The differentiated gap
is then bounded with at most $j_1$ further such derivatives,
as required in the product. A material derivative now contributes
\[
 \lambda_{q+1}\delta_{q+1}^{1/2}
 =\tau_a^{-1}(\tau_a\lambda_{q+1}\delta_{q+1}^{1/2}).
\]
This gives the additional factor in
\eqref{adapt:all-gap-stresses} below. The differentiated gap has the
same loss as $G$. Multiplying the loss for the gap,
$\varepsilon_{q+1}^{-2}(\delta_q/\delta_{B,q})^{1/2}$, by the loss
for the principal increment, at most $\varepsilon_\tau^{-1}$, gives
the asserted coefficient loss.

\subsection{Quadratic terms}

\paragraph{Corrector square and pushforward of the leading quadratic tensor.}
The quadratic corrector stress involves the global increments, so the
corrector contributions must first be summed. On a fixed active time
interval $n$, these increments satisfy
\[
 \begin{aligned}
 w^{\rm c}&=v_{\ell,n,1}^{\rm c}-v_{\ell,n}
                  +(\mathcal U_1^{\rm c}-\IId)G_n^v,\\
 b^{\rm c}&=B_{\ell,n,1}^{\rm c}-B_{\ell,n}
                  +(\mathcal U_1^{\rm c}-\IId)G_n^B.
 \end{aligned}
\]
The increments of the local background are estimated by
\eqref{gg:corrector-comparison-classes}. For the gap terms, use
\[
 (\mathcal U_1^{\rm c}-\IId)G_n^v
 =-\int_0^1\mathcal U_s^{\rm c}\mathcal L_{\xi^{\rm c}}G_n^v\,\dd s,
\]
and the analogous identity for $G_n^B$. By Proposition~\ref{gg:physical}
and Corollary~\ref{prep:mixed-comparison}, the sources have sharp
bounds through derivative order $\rprep-12$ and lossy bounds through
$\rgood+2$, for compositions of at most $j_1$ material and magnetic
derivatives.
Estimating the Lie derivative requires each factor through one
additional spatial derivative. Lemma~\ref{gg:transport-classes}
preserves the resulting bounds with respect to the fixed local background.
The ratios of the gap contributions to the amplitudes in the sharp
and lossy estimates for the increments are bounded, respectively, by
\[
 \begin{aligned}
 \tau_a\lambda_q\varepsilon_{q+1}^2\delta_q^{1/2}
 &=\varepsilon_\tau\varepsilon_{q+1}^{2+\gamma_\ell}\le1,\\
 \tau_c\ell^{-1}\delta_q^{1/2}
 &=1.
 \end{aligned}
\]
Consequently,
\[
 w^{\rm c},b^{\rm c}\in
 \mathcal C_{\rprep-12}
 (C\ell^{-\alpha}\tau_a\lambda_q\delta_{q+1};\lambda_q,\mathrm f)
 \cap\mathcal K_{\rgood+2,j_1}(C\tau_c\ell^{-1}\delta_{q+1}).
\]
The product estimate therefore gives
\begin{equation}\label{stress:corrector-square}
 R^{\mathrm c}_{\rm quad}\in
 \mathcal C_{\rprep-12}(C\ell^{-2\alpha}\tau_a^2\lambda_q^2\delta_{q+1}^2;\lambda_q,\mathrm f)
 \cap\mathcal K_{\rgood+2,j_1}(C\tau_c^2\ell^{-2}\delta_{q+1}^2),
\end{equation}
with loss $\varepsilon_\tau^{-2}(\ell\lambda_q)^{-2}$. For the corrector pushforward of the leading quadratic tensor, use the identity
\[
 (\mathcal U_1^{\rm c}-\IId)\mathsf T
 =-\int_0^1\mathcal U_s^{\rm c}\mathcal L_{\xi^{\rm c}}\mathsf T\,\dd s.
\]
The tensor has sharp amplitude $\delta_{q+1}$, up to a fixed
constant, and $\mathsf T\in\mathcal K(C\delta_{q+1})$.
The integral identity and Lemma~\ref{gg:transport-classes} consequently give
\begin{equation}
 R^{\mathrm p}_{\rm chart}\in
 \mathcal C_{\Nc}(C\ell^{-\alpha}\tau_a^2\lambda_q^2\delta_{q+1}^2;\lambda_{q+1},\mathrm f),
\end{equation}
and
\[
 R^{\mathrm p}_{\rm chart}\in
 \mathcal K(C\ell^{-2\alpha}\tau_c^2\ell^{-2}\delta_{q+1}^2),
\]
with coefficient loss $\varepsilon_\tau^{-2}(\ell\lambda_q)^{-2}\ell^{-\alpha}$.
In the product rule, derivatives of the tensor and of the LDF are
estimated at their respective costs.

\paragraph{Antiderivative error.}
The identity
\[
 D_{t,I}\mathfrak a_I=\alpha_I(t/\tau_a)a_I+a_I^{\rm c}
\]
shows that subtracting $\alpha_I(t/\tau_a)^2a_I^2$ from the
coefficient of the rank-one quadratic tensor leaves
\[
 2\alpha_I(t/\tau_a)a_Ia_I^{\rm c}+(a_I^{\rm c})^2.
\]
The corrected frame and the oscillatory profile are bounded.
Since the complete principal supports are disjoint,
\eqref{principal:primitive-classes} gives
\begin{equation}\label{aniso:stress-material-error}
 R^{\mathrm p}_{\rm time}\in
 \mathcal C_{\Nc}(C\varepsilon_\tau^{j_0}\delta_{q+1};\lambda_{q+1},\mathrm f),
\end{equation}
Its slow coefficients belong to $\mathcal K(C\delta_{q+1})$,
with loss $\varepsilon_\tau^{-j_0}$. This includes all derivatives of
$\varrho$, which occur inside $D_{t,I}^{j_0}a_I$.

\paragraph{Magnetic square.}
Commutation $[D_{t,I},D_{B,I}]=0$ on each local background gives
\[
 D_{B,I}\mathfrak a_I=\sum_{j=0}^{j_0-1}(-1)^j\tau_a^{j+1}
 \alpha_I^{[j+1]}(t/\tau_a)D_{t,I}^jD_{B,I}a_I,
 \qquad
 \|D_{B,I}\mathfrak a_I\|_0\le C\tau_a\lambda_\parallel\delta_{q+1}^{1/2}
 \sum_{j<j_0}\varepsilon_\tau^j\le C\tau_a\lambda_\parallel\delta_{q+1}^{1/2},
\]
by the amplitude bounds and the finite geometric sum. Therefore, by
\eqref{principal:primitive-classes},
\begin{equation}\label{aniso:stress-magnetic-error}
 R^{\mathrm p}_{\rm mag}\in
 \mathcal C_{\Nc}(C\delta_{q+1}(\tau_a\lambda_\parallel)^2;\lambda_{q+1},\mathrm f),
\end{equation}
Its slow coefficients belong to
$\mathcal K(C\varepsilon_\tau^{-2}\delta_{q+1}(\tau_a\lambda_\parallel)^2)$,
with loss $\varepsilon_\tau^{-2}$. We retain this tensor,
including the constant Fourier mode in the fast phase variable.
Its estimate follows directly from the class of
$D_{B,I}\mathfrak a_I$ at the costs of the next iterate, including
magnetic order $j_1$. The parameter inequalities make both this bound and
\eqref{aniso:stress-material-error} smaller than the required bound
for the next stress.

\paragraph{Oscillatory square.}
The remaining rank-one quadratic tensor, tangent to the surfaces of
constant phase, is the corrector pushforward of the explicit tensor
\[
 \alpha_I(t/\tau_a)^2a_I^2\bigl[\varphi'(\lambda_{q+1}y_I^{k_I})^2-1\bigr]
 \partial_{y_I^{\zeta_I}}\otimes\partial_{y_I^{\zeta_I}}.
\]
The fixed trigonometric polynomial $\varphi'^2-1$ has no constant
Fourier mode in the fast phase variable. In the divergence of the
displayed tensor, only the slow coefficients and frame are
differentiated; applying the finite oscillatory inverse divergence gives
the size estimate
\begin{equation}
 C(\ell\lambda_{q+1})^{-1}\delta_{q+1}.
\end{equation}
For every other curl product, a spatial derivative falls on a slow
coefficient or frame and leaves the factor $\lambda_{q+1}^{-1}$.
The product estimate therefore gives the same bound directly.
We keep the constant Fourier modes of all transverse products in
the stress.

For the nonzero modes, apply Lemma~\ref{osc:parametrix} at the
derivative costs of the next iterate. We verify its coefficient
hypotheses using the estimates for $a_I^2$ and the corrected
frames through order $\Nc$. Its background hypotheses follow
from \eqref{gg:corrected-adaptedness}: the corrected background
is $(\lambda_q,\mathrm f)$-adapted through
$(\rprep-12,j_1)$, and the slow spatial scale satisfies
$\ell^{-1}\ge\lambda_q$. These gradient estimates hold on the
whole torus. Only the phase, coefficients, and frames are restricted
to a chart neighborhood of the support. The lemma uses
$\rfg+3$ additional spatial derivatives, and one more gives the
H\"older estimate. For ordinary time derivatives,
\eqref{amplitude:crude-amplitude-bounds} and the corrector pushforward
estimates following \eqref{gg:corrected-chart-bounds} give the
coefficient, frame and coframe bounds through time order $j_1+1$.
The corrected reference fields satisfy the corresponding component
bounds through time order $j_1$. Thus the ordinary derivative
estimate in Lemma~\ref{osc:parametrix} applies through $j_1+1$,
with the same spatial derivative restrictions. The oscillatory-square
contribution to $R^{\mathrm p}_{\rm quad}$ therefore belongs to
\begin{equation}\label{osc:principal-classes}
 \mathcal C_{\Nc-\rfg-4}(C\lambda_{q+1}^{\alpha}(\ell\lambda_{q+1})^{-1}\delta_{q+1};\lambda_{q+1},\mathrm f)
 \cap\mathcal O_{\infty,j_1+1}(C\lambda_{q+1}^{\alpha}(\ell\lambda_{q+1})^{-1}\delta_{q+1};\lambda_{q+1}),
\end{equation}
Each slow coefficient with sharp amplitude $E$ belongs to
$\mathcal K(CE)$, since $a_I^2\in\mathcal K(C\delta_{q+1})$
and the frames belong to $\mathcal K(C)$.
To construct the parametrix, perform the finite expansion in
each invariant chart; the corrected local transports preserve its
phase and coframe. Express the resulting remainders relative to
one fixed global background, sum them according to
\eqref{stress:common-window}, and then apply the global inverse
divergence. The commutator estimate in Lemma~\ref{osc:parametrix}
gives the integer derivative bounds at the costs of the next
iterate. It requires no phase invariance in the common background.

For higher spatial derivatives, first apply
Lemma~\ref{setup:calculus}(vii) to the coefficients. Then use
\eqref{full:endpoint-gradient-bounds} at scale $\lambda_{q+1}$
in the inverse-divergence estimate for the remainder.

\paragraph{First Lie variation.}
\begin{proposition}[Stress from the first Lie variation]\label{stress:first-lie-variation}
Assume the coefficient and frame estimates of
Section~\ref{sec:generator-field-estimates}. Then the contribution
to $R^{\mathrm p}_{\rm quad}$ of
$-\tfrac12\mathcal L_{\xi_{I,1}^{\rm p}}\mathbb Q_{I,1}$ belongs
to the following class at the derivative costs of the next iterate:
\[
 \mathcal C_{\Nc-\rfg-4}
 \bigl(C\lambda_{q+1}^{\alpha}(\ell\lambda_{q+1})^{-1}
       \delta_{q+1};\lambda_{q+1},\mathrm f\bigr)
\]
It also belongs to the ordinary class in
\eqref{osc:principal-classes}. Each slow coefficient in the original
chart, after factoring out the fast profiles, belongs to
$\mathcal K(CE)$ when its sharp amplitude is $E$. More precisely, its
order-zero norm is bounded by
\[
 C\lambda_{q+1}^{\alpha}\tau_a\ell^{-2}
       \lambda_{q+1}^{-1}\delta_{q+1}^{3/2}.
\]
\end{proposition}
\begin{proof}
\emph{1. The leading tensor and the curl terms.}
We apply the oscillatory inverse divergence to the leading rank-one
term to gain $\lambda_{q+1}^{-1}$. The terms in which curl
differentiates a slow coefficient or frame already contain this
factor and can be retained directly.

By \eqref{principal:first-lie-mean}, the only Fourier modes of the
first Lie variation are $\pm1,\pm3$ in the fast phase variable.
Use \eqref{osc:velocity-split}--\eqref{osc:magnetic-split} to
separate each field into its leading term, tangent to the surfaces
of constant phase, and the terms in which curl differentiates a
slow coefficient or frame. For
$F_{I,0}^\pm=f_I^\pm dy_I^{\nu_I}$ and
$W_I^\pm=\curl F_{I,0}^\pm$, commute curl with the Lie derivative
to obtain
\[
 \mathcal L_{\xi_{I,0}^{\rm p}}W_I^\pm
 =\curl((\xi_{I,0}^{\rm p}\cn f_I^\pm)\,dy_I^{\nu_I}).
\]
The estimate for differentiation of scalars along
$\xi_{I,0}^{\rm p}$ bounds the leading term of each differentiated
field by $C\tau_a\ell^{-1}\delta_{q+1}$ and the curl terms
differentiating a slow coefficient or frame by
$C\tau_a\ell^{-2}\lambda_{q+1}^{-1}\delta_{q+1}$.
Multiplying by the undifferentiated field bounds the coefficients of
the leading rank-one tensor by
$C\tau_a\ell^{-1}\delta_{q+1}^{3/2}$.
Every product containing a curl term differentiating a slow coefficient
or frame already has the smaller size
$C\tau_a\ell^{-2}\lambda_{q+1}^{-1}\delta_{q+1}^{3/2}$.
Commutation of curl with the corrector pushforward gives the same
decomposition after the corrector pushforward.

For a slow coefficient $g$, the chart identities give
\[
 \partial_{y_{I,1}^{\zeta_I}}[g\varphi(\lambda_{q+1}y_{I,1}^{k_I})]
 =(\partial_{y_{I,1}^{\zeta_I}}g)\varphi(\lambda_{q+1}y_{I,1}^{k_I}),\qquad
 \ddiv\frac{\partial}{\partial y_{I,1}^{\zeta_I}}=0.
\]
Apply these identities to the divergence of the leading tensor.
The derivative falls only on a slow coefficient or on the frame;
the latter gives
$(\frac{\partial}{\partial y_{I,1}^{\zeta_I}}\cn)\frac{\partial}{\partial y_{I,1}^{\zeta_I}}$
from \eqref{osc:tangent-divergence}. Each such derivative contributes
$\ell^{-1}$. Apply the finite parametrix separately to the four
nonzero Fourier modes to gain $\lambda_{q+1}^{-1}$. The products
containing a curl term satisfy the resulting bound directly and
remain in the stress.

The coefficient estimates for differentiation along the principal
LDF hold through order $\Nc-2$. All these calculations use the
complete $\mathfrak a_I$, without division by an amplitude or a
profile. In particular, after factoring out the fast profiles from
$\xi_{I,0}^{\rm p}\cn$ in the original basis, each slow coefficient
with sharp amplitude $E$ belongs to $\mathcal K(CE)$. The product
estimates retain this property.

\smallskip
\noindent\emph{2. Spatial, material, and magnetic derivatives.}
To prove the estimates in \eqref{osc:principal-classes}, fix
derivative orders satisfying
\[
 k+m\le j_1,\qquad r+k+m\le\Nc-\rfg-4,
\]
Apply $(D_{t,I,1}^{\rm c})^kD_{B_{\ell,I,1}^{\rm c}}^m$ to the
scalar coefficients and the corresponding Lie derivatives to the
corrected frames. The phase is annihilated and the coframe is
preserved by the corrected transports. Consequently, the product
rule and \eqref{principal:primitive-classes} give
\[
 (\lambda_{q+1}\delta_{q+1}^{1/2})^k
 (\lambda_{q+1}\delta_{B,q+1}^{1/2})^m;
\]
the estimates for the approximate antiderivative retain these factors
also at the highest magnetic orders. The derivative along the
principal LDF uses one more coefficient derivative. Apply the
parametrix to each differentiated coefficient, using the
$\rfg+4$ additional spatial derivatives in
\eqref{osc:principal-classes}, one of which is needed for the
H\"older estimate. This proves the asserted bound for material
and magnetic derivatives. Applying ordinary derivatives instead,
and using the corresponding estimates for the coefficients and
frames, proves the ordinary derivative bound. These ordinary
coefficient bounds hold through time order $j_1+1$, by
\eqref{amplitude:crude-amplitude-bounds}, so the same ordinary
derivative estimate in Lemma~\ref{osc:parametrix} applies.
Lemma~\ref{gg:transport-classes} then gives both estimates after
the corrector pushforward.

\smallskip
\noindent\emph{3. The amplitude bound.}
The order-zero estimate from the first step satisfies
\[
 \begin{aligned}
 C\lambda_{q+1}^{\alpha}\tau_a\ell^{-2}
       \lambda_{q+1}^{-1}\delta_{q+1}^{3/2}
 &=C\lambda_{q+1}^{\alpha}
   (\tau_a\ell^{-1}\delta_{q+1}^{1/2})
   (\ell\lambda_{q+1})^{-1}\delta_{q+1}\\
 &\le C\lambda_{q+1}^{\alpha}
   (\ell\lambda_{q+1})^{-1}\delta_{q+1},
 \end{aligned}
\]
where the last inequality uses only
$\tau_a\ell^{-1}\delta_{q+1}^{1/2}\le1$.
\end{proof}

This estimate removes the first-variation restriction on the
H\"older exponent. The stress comparisons and the choice of temporal
parameters in \eqref{iter:exponent-window} are established in
Section~\ref{sec:parameter-choice}. At $\gamma_\ell=\gamma_S=0$, a direct bound
of the first Lie variation of the quadratic tensor by $C\tau_a\ell^{-1}\delta_{q+1}^{3/2}$
would instead require $\gamma_a+\beta>2b\beta$. Combining this with
the principal linear-stress condition
$\gamma_a<1-(2b+1)\beta$ gives $\beta<1/(4b)$. The absence of the constant Fourier mode in the fast phase variable
and the divergence identity for the leading rank-one tensor remove
this additional restriction.

\paragraph{Second variation and covariance.}
For $W_I^\pm=\curl F_{I,0}^\pm$, the estimates for differentiation along the principal LDF on scalars give
\[
 \begin{aligned}
 \|W_I^\pm\|_0&\le C\delta_{q+1}^{1/2},\\
 \|\mathcal L_{\xi_{I,0}^{\rm p}}W_I^\pm\|_0
   &\le C\tau_a\ell^{-1}\delta_{q+1},\\
 \|\mathcal L_{\xi_{I,0}^{\rm p}}^2W_I^\pm\|_0
   &\le C\tau_a^2\ell^{-2}\delta_{q+1}^{3/2}.
 \end{aligned}
\]
The Leibniz rule for the second Lie derivative produces two terms with
both derivatives on one factor and twice the product of the first Lie
derivatives. Each has size $C\tau_a^2\ell^{-2}\delta_{q+1}^2$, so
\[
 \|\mathcal L_{\xi_{I,0}^{\rm p}}^2(W_I^+\odot W_I^-)\|_0
 \le C\tau_a^2\ell^{-2}\delta_{q+1}^2.
\]
To estimate the covariance, use the exact difference formula
\[
 \mathcal U_r^{\rm p}W_I^\pm-\mathcal U_u^{\rm p}W_I^\pm
 =-\int_u^r\mathcal U_s^{\rm p}\mathcal L_{\xi_{I,0}^{\rm p}}W_I^\pm\,\dd s.
\]
To bound the integrand, write it as the curl of
$\mathcal U_s^{\rm p}\mathcal L_{\xi_{I,0}^{\rm p}}F_{I,0}^\pm$.
Its scalar coefficient solves the transport equation in the invariant
coframe. Apply \eqref{principal:structured-coefficient-bound},
taking one additional spatial derivative for curl. The resulting
bound is $C\tau_a\ell^{-1}\delta_{q+1}$, uniformly in $s$.
This argument estimates the scalar coefficient directly and does
not require a uniform bound on the differential of the principal
flow. Integration gives
$C|r-u|\tau_a\ell^{-1}\delta_{q+1}$ for each difference. Their
product, and hence the covariance, is bounded by
$C\tau_a^2\ell^{-2}\delta_{q+1}^2$.

For the second-variation integrand, apply the same argument to
$(\mathcal U_s^{\rm p}\mathcal L_{\xi_{I,0}^{\rm p}}^2W_I^+)\odot
 (\mathcal U_s^{\rm p}W_I^-)$, to the term with the signs exchanged,
and to twice the product of the transported first Lie derivatives.
Each term satisfies the same bound. To include derivatives, recall
that the one-forms $\mathcal L_{\xi_{I,0}^{\rm p}}^jF_{I,0}^\pm$,
$j=0,1,2$, are finite sums of profiles multiplied by slow
coefficients. By \eqref{principal:scalar-action-bounds}, their
amplitudes contain $\tau_a^j\ell^{-j}\delta_{q+1}^{j/2}$. Before
composition with the principal flow, each of these slow coefficients
with sharp amplitude $E$ belongs to $\mathcal K(CE)$. Apply
Corollary~\ref{principal:rescaled-profile-continuation}, followed by
\eqref{principal:structured-coefficient-bound}, with
\eqref{principal:continuation-data}. The product rule gives
\begin{equation}\label{principal:quadratic-second-bound}
 \mathbb Q_{I,2}\in
 \mathcal C_{\Nc-2}(C\tau_a^2\ell^{-2}\delta_{q+1}^2;\lambda_{q+1},\mathrm f)
 \cap\mathcal O_{\infty,j_1+1}(C\tau_a^2\ell^{-2}\delta_{q+1}^2;\lambda_{q+1}).
\end{equation}
Here the first magnetic derivative of the potential and each of its Lie derivatives along the principal LDF
are estimated before taking curls and products. At the local derivative costs,
the factor associated with magnetic differentiation at the endpoint satisfies
\[
 \tau_a\lambda_\parallel\varepsilon_\tau^{-1}=\tau_c\lambda_\parallel
 =\varepsilon_{q+1}^{\gamma_\parallel}\le1.
\]
Thus the magnetic terms have the common amplitude in
\eqref{principal:quadratic-second-bound}. Since the tensor is retained
directly, this coefficient estimate completes its bound.

\paragraph{Quadratic background differences.}
The gap terms in \eqref{principal:gap-square} are the products of
$\curl\Theta_{w,I,1}^{\rm p}$ and $\curl\Theta_{b,I,1}^{\rm p}$ with
$w_{I,1}^{\mathrm p,\mathrm{gap}}$ and $b_{I,1}^{\mathrm p,\mathrm{gap}}$,
and the difference of the tensor squares of the latter fields. By
Proposition~\ref{gap:field-bounds}, their norm is at most
\[
 C\tau_a\lambda_{q+1}\varepsilon_{q+1}^2\delta_q^{1/2}\delta_{q+1}
 (1+\tau_a\lambda_{q+1}\varepsilon_{q+1}^2\delta_q^{1/2}).
\]
Combined with
\eqref{stress:principal-linear-gap} and
\eqref{stress:principal-comparison}, and using
$\tau_a\lambda_{q+1}\varepsilon_{q+1}^2\delta_q^{1/2}=o(1)$,
all principal background-difference stresses lie in
\begin{equation}\label{adapt:all-gap-stresses}
 \mathcal C_{\Nc}\bigl(C\lambda_{q+1}^{(j_1+2)\alpha}
       (1+\tau_a\lambda_{q+1}\delta_{q+1}^{1/2})\varepsilon_{q+1}^2\delta_q^{1/2}
                                             \delta_{q+1}^{1/2};\lambda_{q+1},\mathrm f\bigr),
\end{equation}
Before composition with the principal flow, each slow coefficient
in these gap representations, with sharp amplitude $E$, belongs to
\[
 \mathcal K_{\rgood+1,j_1}
 \bigl(CE\varepsilon_{q+1}^{-2}(\delta_q/\delta_{B,q})^{1/2}
               \varepsilon_\tau^{-1}\bigr).
\]
Thus the coefficient loss is at most
$\varepsilon_{q+1}^{-2}(\delta_q/\delta_{B,q})^{1/2}\varepsilon_\tau^{-1}$.
We have estimated the amplitude of every term in
\eqref{stress:linear-quadratic-split}. It remains to prove these
estimates through spatial order $\rgood$, and for material and
magnetic derivatives as in \eqref{full:stress-transport}, and to
compare the amplitudes with the inductive bound.

\subsection{Estimates for higher derivatives}\label{ssec:derivative-closure}

\paragraph{Comparison of transport operators.}
We compare the transport operators in the estimates above with
those of the endpoint. The terms $R^{\mathrm{cut}}$, $R^{\rm moll}$
and the corrector stresses have been estimated relative to the local
background; the other terms have been estimated relative to the
corrected background.

We apply Lemma~\ref{setup:calculus}(iv). By the gradient estimates
in Proposition~\ref{full:endpoint-derivative-closure}, the local and
corrected backgrounds are $(\lambda_{q+1},\mathrm f)$-adapted
through $(\infty,j_1)$. The selected background satisfies the
same bounds, and the endpoint background satisfies them through
$(\rgood,j_1)$. The differences between these backgrounds consist
of the corrector integrals $c_{v,n,1},c_{B,n,1}$ from
\eqref{prep:corrector-final-comparison} and the remaining
increments. The former belong to
$\mathcal C_\infty(\,\cdot\,;\lambda_{q+1},\mathrm f)$, with
amplitudes at most $\delta_{q+1}^{1/2}$ and
$\delta_{B,q+1}^{1/2}$, by the estimates of that proposition. For the
remaining increments $u,b$, use
\eqref{full:endpoint-increment-bounds}.

All hypotheses of the lemma thus hold with
$\Lambda=\lambda_{q+1}$. It follows that the
$\mathcal C_{N',J}(E;\lambda_{q+1},\mathrm f)$ estimates are
equivalent, up to fixed constants, in these four backgrounds whenever
$N'\le\rgood$ and $J\le j_1$. The argument applies to transport
and Lie derivatives of the same fixed tensor. The material and
magnetic orders on a coefficient sum to at most $j_1-1$. We may
therefore prove Proposition~\ref{stress:complete-bounds} in any of these
backgrounds.

\paragraph{Antiderivative and magnetic errors.}
For the antiderivative error $R^{\mathrm p}_{\rm time}$ and the magnetic
square $R^{\mathrm p}_{\rm mag}$, the full gains in
\eqref{aniso:stress-material-error} and
\eqref{aniso:stress-magnetic-error} persist at the costs of the next
iterate. To see this, use Corollary~\ref{aniso:endpoint-residual-gain}
with $\Lambda=\lambda_{q+1}$,
$\mathrm a_t=\lambda_{q+1}\delta_{q+1}^{1/2}$,
$\mathrm a_B=\lambda_{q+1}\delta_{B,q+1}^{1/2}$,
$\mathrm a_{t,0}=\tau_a^{-1}$, $\mathrm a_{B,0}=\lambda_\parallel$ and
$e=\varepsilon_\tau$. The stress terms under consideration contain
at most two magnetic derivatives of the approximate antiderivative
before further differentiation. We apply
\eqref{principal:primitive-overflow} to these terms. The comparisons
\eqref{setup:rate-comparisons} supply the required derivative costs, and
$j_1\ge j_0+2$ by \eqref{iter:depth-anisotropy}. The analogous gain
for the cutoff error follows from Theorem~\ref{galbrun:linear-theorem}(f).

\paragraph{Smoothing error.}
We estimate the derivatives of the mollification error using
\eqref{stress:mollification-classes}. Fix $(r,k,m)$ such that
\[
 r+k+m\ge1,\qquad k+m\le j_1,\qquad r+k+m\le\rgood-1.
\]
Expand
$\mathcal D_{t,q+1}^k\mathcal L_{B_{q+1}}^mR_n^{\rm moll}$
using Lemma~\ref{setup:calculus}(iv), then apply the spatial
derivatives and the product rule. First consider the terms in which
the stress $R_n^{\rm moll}$ is undifferentiated. These are products of background
gradients of degree $(k,m)$ and $R_n^{\rm moll}$. The gradient
bounds give
\[
 C\lambda_{q+1}^r(\lambda_{q+1}\delta_{q+1}^{1/2})^k
 (\lambda_{q+1}\delta_{B,q+1}^{1/2})^m\|R^{\rm moll}\|_0.
\]
The factor involving the undifferentiated stress is small by
\eqref{stress:mollification-zero}.

In each remaining term, the spatial, material and magnetic derivatives
of $R_n^{\rm moll}$ have orders $(r',k',m')$, respectively, where
\[
 r'+k'+m'\ge1,\qquad k'\le k,\qquad m'\le m,\qquad
 r'+k'+m'\le r+k+m.
\]
\emph{Case 1: $r'+k'+m'\le\rprep$.}
Apply the sharp estimate in \eqref{stress:mollification-classes}
and compare it with
\[
 \lambda_{q+1}^{r'-\alpha}
 (\lambda_{q+1}\delta_{q+1}^{1/2})^{k'}
 (\lambda_{q+1}\delta_{B,q+1}^{1/2})^{m'}\delta_{q+2}.
\]
By \eqref{setup:rate-comparisons}, the ratio of the sharp bound to
the required bound is at most
\[
 \varepsilon_{q+1}^{r'+k'(1-\beta-\gamma_\ell)+m'(1-\beta-\gamma_a-\gamma_\ell)-2b\beta-\gamma_S}
 \le\varepsilon_{q+1}^{1-\beta-\gamma_a-\gamma_\ell-\gamma_S-2b\beta}=o(1)
\]
by the last inequality in \eqref{iter:auxiliary-margins}. The
remaining gradient factors have exactly the required derivative costs.
\smallskip
\noindent\emph{Case 2: $r'+k'+m'>\rprep$.}
Apply instead the lossy estimate in the same display. After
division by the required bound, we obtain
\[
 \varepsilon_{q+1}^{-2b\beta-2\gamma_S}(\ell\lambda_{q+1})^{-r'}
 \max\Bigl\{1,\frac{\ell^{-1}}{\lambda_{q+1}\delta_{B,q+1}^{1/2}}\Bigr\}^{j_1}
 \le\varepsilon_{q+1}^{(1-\gamma_\ell)(\rcut-j_1+1)-2b\beta-2\gamma_S
   -j_1[b\beta/(b-1)+\gamma_a+\gamma_\parallel-(1-\gamma_\ell)]^+}=o(1).
\]
Indeed, the restriction on the number of material and magnetic
derivatives gives
\[
 r'\ge\rprep-j_1+1\ge\rcut-j_1+1,
\]
and the positivity required here follows from
\eqref{iter:general-final-reserve} and $\gres\ge2b\beta+2\gamma_S$.
To compare the exponents, note that the cited inequality uses
$\rcut-j_1-1$ in place of $\rcut-j_1+1$, subtracts $\gres$ in place
of $2b\beta+2\gamma_S$, and allows $j_1+1$ in place of $j_1$
lossy transport derivatives. Its positive coefficient
$\gamma_a+\beta$ is also smaller than $1-\gamma_\ell$.
For ordinary time derivatives, \eqref{prep:old-stress-time} gives the
same argument with cost at most $\ell^{-1}\le\lambda_{q+1}$ per
derivative. Thus, for $r+h\ge1$, the ratio satisfies
\[
 \varepsilon_{q+1}^{(1-\gamma_\ell)(r+h)-2b\beta-2\gamma_S}=o(1).
\]
The two cases give the required bound for all derivatives of
$R^{\rm moll}$ of positive order. The approximation gain is used only at
order zero. Since this tensor involves no nonlocal Fourier multiplier,
the ordinary time derivatives of $R_q$ have order at most $j_1$
and its spatial and ordinary time orders sum to at most $\rgood$,
as in \eqref{prep:old-stress-time}.

\paragraph{Lossy classes.}
For the corrector stresses, we apply the product estimate in
Lemma~\ref{setup:calculus}(i) to the decomposition
\eqref{gg:five-terms}. The potentials, fields, and gaps satisfy
\begin{equation}
 \begin{gathered}
 \Theta_n^{\rm c}\in\mathcal K(C\tau_c^2\delta_{q+1}),\qquad
 F_{w,n},F_{b,n}\in\mathcal K(C\tau_c\delta_{q+1}),\\
 \mathcal D_{t,n}\xi^{\rm c},\mathcal L_{B_{\ell,n}}\xi^{\rm c}\in\mathcal K(C\tau_c\ell^{-1}\delta_{q+1}),\\
 (\mathcal D_{t,n}^2-\mathcal L_{B_{\ell,n}}^2)\Theta_n^{\rm c}\in\mathcal K(C\lambda_{q+1}^{(j_1+1)\alpha}\delta_{q+1}),\\
 G_n^v,G_n^B\in\mathcal K_{\rgood+3,j_1+1}(C\delta_q^{1/2}),\qquad
 \nabla v_{\ell,n},\nabla B_{\ell,n}\in\mathcal K_{j_1}(C\lambda_q\ell^{-\alpha}\delta_q^{1/2}),
 \end{gathered}
\end{equation}
These follow from Proposition~\ref{gg:physical},
Corollary~\ref{prep:mixed-comparison} and
\eqref{prep:local-field-classes}. We also use
\eqref{gg:corrector-class-ranges} and \eqref{gg:corrector-defects}:
the operators $\mathcal U_s^{\rm c}$ and $\mathcal P^{\rm c}$ are
bounded by a fixed constant in the lossy classes, whereas
$\mathcal P^{\rm c}-\IId$ multiplies the amplitude by
$\tau_c^2\ell^{-2}\delta_{q+1}=\varepsilon_{q+1}^{2\beta}\le1$
and requires one additional derivative. The five products in
\eqref{gg:five-terms}, in the order displayed there, belong to
\[
 \begin{gathered}
 \mathcal K(C\lambda_{q+1}^{(j_1+1)\alpha}\tau_c^2\ell^{-2}\delta_{q+1}^2),\qquad
 \mathcal K(C\tau_c^2\ell^{-2}\delta_{q+1}^2),\\
 \mathcal K_{\rgood+2,j_1}(C\tau_c\ell^{-1}\delta_q^{1/2}\delta_{q+1}),\\
 \mathcal K_{\rgood+1,j_1}(C\tau_c^2\ell^{-2}\delta_q^{1/2}\delta_{q+1}),\qquad
 \mathcal K_{\rgood+1,j_1}(C\tau_c^3\ell^{-3}\delta_q^{1/2}\delta_{q+1}^2),
 \end{gathered}
\]
respectively. The operator $\mathcal R\curl$ introduces the factor
$\lambda_{q+1}^{(j_1+1)\alpha}$. The same argument gives the
$\mathcal K_{\rgood+1,j_1}$ bound in \eqref{gg:crude} for the Nash
and comparison tensors: the background gradient contributes at most
$C\lambda_q\ell^{-\alpha}\delta_q^{1/2}\le C\ell^{-1}$ for one
of the derivatives. Since $\tau_c\ell^{-1}\ge1$ and
$\tau_c\ell^{-1}\delta_{q+1}\le1$, the fourth amplitude dominates the
third and the fifth, and
\begin{equation}\label{gg:crude}
 R^{\mathrm c}_{\rm tr},\,R^{\mathrm c}_{\mathrm{na}},\,R^{\mathrm c}_{\rm moll}
 \in\mathcal K_{\rgood+1,j_1}\bigl(C\lambda_{q+1}^{(2j_1+2)\alpha}
 [\tau_c^2\ell^{-2}\delta_{q+1}+\tau_c^2\ell^{-2}\delta_q^{1/2}]\delta_{q+1}\bigr).
\end{equation}
The bound on the sum of the spatial, material and magnetic orders
is two less than for the gap: one derivative forms
$\mathcal D_{t,n}G_n^v$ and a second enters the Lie derivative.
For the quadratic corrector stress, \eqref{stress:corrector-square}
and the same estimates give
$R^{\mathrm c}_{\rm quad}\in\mathcal K_{\rgood+2,j_1}(C\tau_c^2\ell^{-2}\delta_{q+1}^2)$,
which extends beyond the required order $\rgood+1$.

For fast tensors, we use the coefficient estimates listed below.
After factoring out the fast profiles, and before any composition
with the principal flow, each slow coefficient with lossy amplitude
$E'$ belongs to
$\mathcal K_{\rgood+1,j_1}(CE')$ when a gap occurs, and to
$\mathcal K(CE')$ otherwise.
Lemma~\ref{setup:calculus}(vii) gives these derivative ranges. For a tensor transported by the
principal flow, we first apply
Corollary~\ref{principal:rescaled-profile-continuation} to the
coefficients, source, and initial datum of the transport equation.
Its hypothesis follows from \eqref{iter:general-final-reserve}.

\paragraph{Summary of the estimates.}
The estimates for the tensors in \eqref{stress:linear-quadratic-split}
are collected below. For each size $E$, the table gives the order $N'$
through which the sharp estimate holds at the costs of the next iterate,
and an exponent $g$ with $E'/E\le\varepsilon_{q+1}^{-g}$ for the
$\mathcal K$ amplitude. For a slow tensor the $\mathcal K$ bound
applies to the tensor itself; for a fast tensor it applies to the
original slow coefficients, after factoring out the fast profiles
and before any composition with the principal flow. These
$\mathcal K$ bounds refer to the local background, or its corrector
pushforward, used in the corresponding estimate above. The bound is in
$\mathcal K(CE')$ except for terms containing gaps, where we use
$\mathcal K_{\rgood+1,j_1}(CE')$, and the smoothing error, where
\eqref{stress:mollification-classes} gives
$\mathcal K_{\rgood-1,j_1}(C\lambda_q^{-\alpha}\delta_{q+1})$.
These exponents follow from \eqref{setup:scale-identities}
and the identities
\[
 \begin{gathered}
 \tau_a\lambda_\parallel\le\varepsilon_\tau,\qquad
 \tau_c=\varepsilon_\tau^{-1}\tau_a,\qquad
 \ell\lambda_q=\varepsilon_{q+1}^{\gamma_\ell},\\
 (\delta_q/\delta_{B,q})^{1/2}
   =\varepsilon_{q+1}^{-(\gamma_a+\gamma_\parallel)/b},\\
 \tau_a\ell^{-1}
   =\varepsilon_{q+1}^{\gamma_a-\beta/(b-1)}
   \le\varepsilon_{q+1}^{-\beta/(b-1)}.
 \end{gathered}
\]
For a sum, we use the largest termwise quotient as a loss bound;
for a product, we use the product of the factor quotients.
\begin{center}
\footnotesize
\setlength{\tabcolsep}{3pt}
\renewcommand{\arraystretch}{1.25}
\begin{tabular}{@{}>{\raggedright\arraybackslash}p{30mm}>{\raggedright\arraybackslash}p{76mm}>{\centering\arraybackslash}p{20mm}>{\raggedright\arraybackslash}p{37mm}@{}}
\toprule
Stress (definition) & Contribution and size $E$ & Maximum derivative order & Loss exponent $g$ in the $\mathcal K$ bound\\
\midrule
$R^{\mathrm{cut}}$ \eqref{mom:collar-stress}
 & $\ell^{-\alpha}\varepsilon_\tau^{j_0}\delta_{q+1}$
 & $\rprep-j_0-7$ & $j_0\gamma_a$\\
$R^{\rm moll}$ \eqref{mom:mollification-stress}
 & Positive orders: $\lambda_{q+1}^{-\alpha}\delta_{q+2}$
 & $\rprep$ & $2b\beta+\gamma_S$\\
$R^{\mathrm c}_{\rm tr},R^{\mathrm c}_{\rm moll}$\newline
 \eqref{mom:corrector-transport-stress}, \eqref{mom:G-gap}
 & $\lambda_{q+1}^{(2j_1+4)\alpha}[\tau_a^2\lambda_q^2\delta_{q+1}+\tau_a\lambda_q\varepsilon_{q+1}^2\delta_q^{1/2}]\delta_{q+1}$
 & $\rprep-14$
 & $2+2\gamma_a+\beta/(b-1)+\gamma_\ell$\\
$R^{\mathrm c}_{\mathrm{na}}$ \eqref{mom:corrector-interaction-stress}
 & $\varepsilon_\tau\times$ the preceding
 & $\rprep-14$
 & $2+3\gamma_a+\beta/(b-1)+\gamma_\ell$\\
$R^{\mathrm p}_{\rm tr},R^{\mathrm p}_{\mathrm{na}}$\newline
 \eqref{mom:principal-transport-stress}, \eqref{mom:principal-interaction-stress}
 & Fields relative to the local background: $\lambda_{q+1}^{(j_1+2)\alpha}\delta_{q+1}^{1/2}/(\lambda_{q+1}\tau_a)$
 & $\Nc$ & $\beta/(b-1)+2\gamma_S$\\
$R^{\mathrm p}_{\rm tr},R^{\mathrm p}_{\mathrm{na}}$\newline
$R^{\mathrm p}_{\rm moll},R^{\mathrm p}_{\rm quad}$\newline
 \eqref{mom:principal-gap-stress}, \eqref{principal:gap-square},
 \eqref{mom:principal-transport-stress}, \eqref{mom:principal-interaction-stress}, \eqref{mom:N-remainder}
 & Comparison and gap terms: $\lambda_{q+1}^{(j_1+2)\alpha}(1+\tau_a\lambda_{q+1}\delta_{q+1}^{1/2})\varepsilon_{q+1}^2(\delta_q\delta_{q+1})^{1/2}$
 & $\Nc$
 & $2+\gamma_a+\frac{\gamma_a+\gamma_\parallel}{b}+2\gamma_S$\\
$R^{\mathrm c}_{\rm quad}$ \eqref{mom:corrector-quadratic-stress}
 & $\ell^{-2\alpha}\tau_a^2\lambda_q^2\delta_{q+1}^2$
 & $\rprep-12$ & $2\gamma_a+2\gamma_\ell$\\
$R^{\mathrm p}_{\rm chart}$ \eqref{mom:chart-stress}
 & $\ell^{-\alpha}\tau_a^2\lambda_q^2\delta_{q+1}^2$
 & $\Nc$ & $2\gamma_a+2\gamma_\ell+\gamma_S$\\
$R^{\mathrm p}_{\rm time}$ \eqref{osc:primitive-magnetic-tensors}
 & $\varepsilon_\tau^{j_0}\delta_{q+1}$
 & $\Nc$ & $j_0\gamma_a$\\
$R^{\mathrm p}_{\rm mag}$ \eqref{osc:primitive-magnetic-tensors}
 & $(\tau_a\lambda_\parallel)^2\delta_{q+1}$
 & $\Nc$ & $2\gamma_a$\\
$R^{\mathrm p}_{\rm quad}$ \eqref{mom:N-remainder}
 & Oscillations and first Lie variation: $\lambda_{q+1}^{\alpha}(\ell\lambda_{q+1})^{-1}\delta_{q+1}$
 & $\Nc-\rfg-4$ & $\gamma_S$\\
$R^{\mathrm p}_{\rm quad}$ \eqref{mom:N-remainder}
 & Second variation, covariance: $\tau_a^2\ell^{-2}\delta_{q+1}^2$
 & $\Nc-2$ & $0$\\
\bottomrule
\end{tabular}
\end{center}
The losses in the third and fourth rows follow by dividing the
lossy amplitudes in \eqref{gg:crude} by the corresponding sharp
amplitudes in \eqref{stress:corrector-linear}. The two bracketed terms
give, respectively,
\[
 \begin{aligned}
 \frac{\tau_c^2\ell^{-2}}{\tau_a^2\lambda_q^2}
   &=\varepsilon_\tau^{-2}(\ell\lambda_q)^{-2},\\
 \frac{\tau_c^2\ell^{-2}}{\tau_a\lambda_q\varepsilon_{q+1}^2}
   &=\varepsilon_\tau^{-2}(\tau_a\ell^{-1})
       (\ell\lambda_q)^{-1}\varepsilon_{q+1}^{-2}.
 \end{aligned}
\]
The lossy H\"older factor $\lambda_{q+1}^{(2j_1+2)\alpha}$ is smaller
than the sharp one and therefore does not increase either quotient.
The entries of the last column are therefore bounded by
\begin{equation}\label{stress:loss-maximum}
 \begin{aligned}
 \gres:={}&2+3\gamma_a+\frac{\beta}{b-1}
       +\frac{\gamma_a+\gamma_\parallel}{b}+j_0\gamma_a+2b\beta\\
       &+4\gamma_\ell+(c+3)\gamma_S,
 \end{aligned}
\end{equation}
where the term $c\gamma_S$ includes the H\"older factors suppressed in
the lossy amplitudes; as shown below, there are at most $c$ such factors.
Thus this exponent bounds the loss in every row. The smallest sharp
derivative range in the table ends at
\begin{equation}\label{stress:ceiling-minimum}
 \rcut:=\Nc-\rfg-4=\rprep-j_0-\rfg-16,
\end{equation}
Thus, after the changes of transport, every tensor has a sharp
class through at least $\rcut$. Apart from $R^{\rm moll}$,
the slow tensors and the original slow coefficients of fast tensors
with sharp amplitude $E$, before composition with the principal
flow, belong to
$\mathcal K_{\rgood+1,j_1}(CE\varepsilon_{q+1}^{-\gres})$
relative to the backgrounds used in their estimates above.
The derivatives of $R^{\rm moll}$ were estimated separately. We use these bounds in
Section~\ref{ssec:derivative-reserve} to extend the sharp estimates.

\paragraph{H\"older factors.}
In the classes of Section~\ref{ssec:fixed-scales} the H\"older bound
of a product contains the larger of the H\"older factors of its
factors, whereas the fixed factors $\ell^{-\alpha}$ and $\lambda_{q+1}^{\alpha}$
already contained in integer amplitudes multiply, by rule (i). For the corrector stresses, the proof of
\eqref{gg:high-mixed-state} absorbs the H\"older factors in the
lower-order Galbrun products before applying $\mathscr T$; the resulting
class amplitude contains only $\lambda_{q+1}^{(j_1+1)\alpha}$.
The difference between the path average and the identity contributes
at most one further factor $\ell^{-\alpha}$, so the forces in
\eqref{gg:five-terms} have amplitudes containing at most $j_1+2$
such factors.
The final $\mathcal R\curl$ contributes at most $j_1+1$ by
Lemma~\ref{setup:calculus}(v), and the defining H\"older bound adds
one common factor. Thus the exponent $2j_1+4$ in
\eqref{stress:corrector-linear} and the table also covers the
$C^\alpha$ norms. The commutator and gap contributions use products
of the first-derivative and gap bounds and satisfy the same count;
the Nash contribution retains the additional factor $\varepsilon_\tau$.
The lossy calculation \eqref{gg:crude} uses at most $2j_1+2$
factors. Products with gaps or local gradients elsewhere add at most
two, and the comparison of transport operators in
Lemma~\ref{setup:calculus}(iii),(iv) introduces no further factor.
The constant $c=64+14j_0+8j_1$ in
\eqref{setup:holder-constant} therefore covers all these losses.

\paragraph{Completion of the stress estimates.}
We first compare the amplitudes with the desired stress bound.
For every term in the table except $R^{\rm moll}$, the sharp
amplitude is smaller than
$\lambda_{q+1}^{-\alpha}\delta_{q+2}$ by a factor tending to zero.
Indeed, the parameter comparisons contain
$\lambda_{q+1}^{-c\alpha}$, which absorbs the H\"older factors
proved above. For the spatial part of the order-zero smoothing
error use \eqref{iter:mollifier-budget}; for its flow part use
\eqref{iter:scalar-budgets}. The cutoff error also satisfies
\eqref{iter:scalar-budgets}. The gap products
satisfy \eqref{iter:principal-gap-budgets}. The four bounds in
\eqref{iter:scalar-budgets} also control the corrector products,
principal transport, oscillatory inverse divergence, second
variation, and the antiderivative error
$\varepsilon_\tau^{j_0}\delta_{q+1}$. Finally,
\eqref{aniso:magnetic-stress-budget} controls the magnetic square.

We next prove the derivative estimates, beginning with
$R^{\mathrm p}_{\rm tr}-R_{\rm tr}^{\rm gap}$, the term that
imposes the principal linear-stress condition. Here
$R_{\rm tr}^{\rm gap}$ is the sum of the gap terms in
\eqref{stress:principal-linear-gap}. Recall from
\eqref{stress:common-window} that the spatial multiplier is
applied after summing the local terms. The local transport
potential on each complete support is
$\mathcal U_1^{\rm c}[(D_{t,I}f_{w,I,1}-D_{B,I}f_{b,I,1})dy_I^{\nu_I}]$,
where $D_{t,I}f_{w,I,1}-D_{B,I}f_{b,I,1}$ solves
\eqref{principal:exact-transport-potential}.
We apply Corollary~\ref{principal:rescaled-profile-continuation}
to this equation. The estimates following
\eqref{principal:exact-transport-potential} bound its source by
the sharp amplitude
$\lambda_{q+1}^{-1}\tau_a^{-1}\delta_{q+1}^{1/2}$ through order
$\rprep-j_0-2$, at the derivative costs of the next iterate.
These estimates apply separately to
$D_{t,I}^2\mathfrak a_I$, $D_{B,I}^2\mathfrak a_I$, and the two
commutators. The original slow coefficients in this source
representation, before composition with the principal flow, belong to
$\mathcal K(C\lambda_{q+1}^{-1}\ell^{-1}\delta_{q+1}^{1/2})$.
The ratio of this amplitude to the sharp amplitude is
$\tau_a\ell^{-1}\le\varepsilon_{q+1}^{-\gres}$. For the
slow coefficients of $\xi_{I,0}^{\rm p}\cn$, after factoring out
the fast profiles in the original chart, a sharp amplitude $E$
gives membership in $\mathcal K(CE)$.

Choose $\Lambda_1=(\tau_a\delta_{q+1}^{1/2})^{-1}$, so that
$(\ell\Lambda_1)^{-1}=\vartheta$ as in
\eqref{principal:continuation-data}. To verify
\eqref{setup:reserve-inequality}, we estimate
\[
 \begin{aligned}
 &\frac{E'}{E}(\ell\Lambda_1)^{-(\rprep-j_0-2-j_1)}
 \max\Bigl\{1,\frac{\ell^{-1}}{\lambda_{q+1}\delta_{q+1}^{1/2}},
             \frac{\ell^{-1}}{\lambda_{q+1}\delta_{B,q+1}^{1/2}}\Bigr\}^{j_1}\\
 &\qquad\le\varepsilon_{q+1}^{-\gres}\vartheta^{\rcut-j_1-1}
 \max\Bigl\{1,\frac{\ell^{-1}}{\lambda_{q+1}\delta_{q+1}^{1/2}},
             \frac{\ell^{-1}}{\lambda_{q+1}\delta_{B,q+1}^{1/2}}\Bigr\}^{j_1+1},
 \end{aligned}
\]
since $\rprep-j_0-2\ge\rcut$ and $\vartheta\le1$, and the right side
is $o(1)$ by \eqref{iter:transported-profile-comparison}.
Corollary~\ref{principal:rescaled-profile-continuation} therefore
gives $(D_{t,I}f_{w,I,1}-D_{B,I}f_{b,I,1})dy_I^{\nu_I}\in
\mathcal C_\infty(C\lambda_{q+1}^{-1}\tau_a^{-1}\delta_{q+1}^{1/2};\lambda_{q+1},\mathrm f)$
with respect to the local background transports. Thus
\eqref{principal:transport-source-bound} holds for all spatial
derivatives with the same amplitude.
The corrector pushforward preserves this class at every order with respect
to the corrected background, by the case $\Lambda=\lambda_{q+1}$ of
\eqref{gg:corrector-class-ranges}. The corrected background is
$(\lambda_{q+1},\mathrm f)$-adapted through $(\infty,j_1)$ by
\eqref{full:endpoint-gradient-bounds}, so Lemma~\ref{setup:calculus}(v)
places $R^{\mathrm p}_{\rm tr}-R_{\rm tr}^{\rm gap}$ in
$\mathcal C_\infty(C\lambda_{q+1}^{(j_1+2)\alpha}\lambda_{q+1}^{-1}\tau_a^{-1}\delta_{q+1}^{1/2};\lambda_{q+1},\mathrm f)$,
with the amplitude of \eqref{principal:natural-stresses} for all
spatial derivatives. The comparison of transport operators at the
start of this subsection gives the same class relative to the
endpoint background through order $\rgood$. Consequently, for $k+m\le j_1$ and
$r+k+m\le\rgood-1$,
\[
 \frac{\|\mathcal D_{t,q+1}^k\mathcal L_{B_{q+1}}^m
       (R^{\mathrm p}_{\rm tr}-R_{\rm tr}^{\rm gap})\|_r}
      {\lambda_{q+1}^{r-\alpha}(\lambda_{q+1}\delta_{q+1}^{1/2})^k
       (\lambda_{q+1}\delta_{B,q+1}^{1/2})^m\delta_{q+2}}
 \le C\lambda_{q+1}^{(j_1+3)\alpha}
      \frac{\tau_a^{-1}\lambda_{q+1}^{-1}\delta_{q+1}^{1/2}}{\delta_{q+2}}
 =C\varepsilon_{q+1}^{1-(2b+1)\beta-\gamma_a-\gamma_\ell-(j_1+3)\gamma_S},
\]
by \eqref{setup:scale-identities}, and the exponent is positive because
$1-\beta-\gamma_a-\gamma_\ell>2b\beta+c\gamma_S$ with $c\ge j_1+3$, by
Lemma~\ref{iter:parameter-lemma}. This is the second entry of
\eqref{iter:scalar-budgets} together with its H\"older factors.

For $R_{\rm tr}^{\rm gap}$, apply
\eqref{stress:principal-linear-gap} and
\eqref{adapt:all-gap-stresses}. The lossy estimates hold through
at least $\rgood+1$, so the same application of the coefficient
lemma and comparison of transport operators gives the sharp
estimate through $\rgood$. The amplitude satisfies
\eqref{iter:principal-gap-budgets}. Adding this estimate to the
one just proved gives the required bound for
$R^{\mathrm p}_{\rm tr}$. The gap estimates through $\rgood+1$
include every derivative used in this argument.

For the remaining slow tensors, apply
Lemma~\ref{setup:calculus}(vii) with
$\Lambda_1=\lambda_{q+1}$. For tensors transported by the principal
flow, use \eqref{principal:continuation-data}, as in the preceding
argument. The parameter inequality
\eqref{iter:general-final-reserve} applies with the loss from
\eqref{stress:loss-maximum} and the derivative order from
\eqref{stress:ceiling-minimum}. It gives the sharp estimates
through $\rgood$. The lossy coefficient estimates hold through
at least $\rgood+1$, including the additional spatial derivative
used to obtain the H\"older norm in physical coordinates from the
chart estimates. Together with the separate estimate for
the smoothing error, this proves the required derivative estimates
for every tensor.

The tensors retained directly in the stress, including
$R^{\mathrm p}_{\rm time}$ and $R^{\mathrm p}_{\rm mag}$, require
no Fourier multiplier. In the slow case, the choice
$\Lambda_1=\lambda_{q+1}$ in Lemma~\ref{setup:calculus}(vii)
satisfies a weaker hypothesis than the transported case, since
$(\ell\lambda_{q+1})^{-1}\le\vartheta$. In each case the ratio of
the sharp amplitude to
$\lambda_{q+1}^{-\alpha}\delta_{q+2}$ tends to zero by the
corresponding parameter inequality, including the H\"older
factors. The ratio of the lossy amplitude to the sharp amplitude
is at most $\varepsilon_{q+1}^{-\gres}$ by
\eqref{stress:loss-maximum}.

Ordinary derivatives follow from these classes through time order
$j_1$. Expanding $\partial_t=D_{t,q+1}-v_{q+1}\cn$ in a class with
$j_1$ transport derivatives costs at most $C\lambda_{q+1}$ per
ordinary derivative, by Proposition~\ref{full:endpoint-derivative-closure}.
The linear and quadratic sums consequently satisfy
\begin{equation}\label{stress:prepared-complete}
 \begin{aligned}
 \|\partial_t^hR^{\rm lin}\|_r+\|\partial_t^hR^{\rm quad}\|_r
 &=o(\lambda_{q+1}^{r+h-\alpha}\delta_{q+2}),\\
 \|\mathcal D_{t,*}^k\mathcal L_{B_*}^mR^{\rm lin}\|_r
 +\|\mathcal D_{t,*}^k\mathcal L_{B_*}^mR^{\rm quad}\|_r
 &=o\bigl(\lambda_{q+1}^{r-\alpha}
 (\lambda_{q+1}\delta_{q+1}^{1/2})^k
                       (\lambda_{q+1}\delta_{B,q+1}^{1/2})^m\delta_{q+2}\bigr),
 \end{aligned}
\end{equation}
The first line holds for
\[
 h\le j_1,\qquad r+h\le\rgood,
\]
and the second holds for
\[
 k+m\le j_1,\qquad r+k+m\le\rgood-1;
\]
the additional spatial derivative serves the H\"older estimate and the
comparison between Lie and transport derivatives. The comparison of transport operators proved at the start of the
subsection gives these bounds relative to the endpoint background. Since the $o(1)$ bounds are uniform
in $q$ as $a\to\infty$, increasing $a$ absorbs the finite sum and all
operator constants. In particular, the Lie derivative estimate has
constant one. This proves
\eqref{full:stress-spatial}--\eqref{full:stress-transport} and hence
Proposition~\ref{stress:complete-bounds}.\qed

\subsection{Pressure}

Corollary~\ref{galbrun:split-pressure} applies to the prescribed source
$\mathsf F_n$ in \eqref{amplitude:forced-equation}, whose amplitude is
$C\delta_{q+1}$ by the coefficient construction. It gives the local
pressure bounds for $k+m\le j_1$ and $r+k+m+4\le\rprep$,
including magnetic order $j_1$. The ordinary estimate
\eqref{galbrun:ordinary-pressure-bound} also gives
\begin{equation}\label{gg:pressure-bound}
 \pi_n\in\mathcal O_{\infty,j_1}(C\ell^{-\alpha}\delta_{q+1};\ell^{-1}),
\end{equation}
We express the transport derivatives of $dp_{q+1}$ relative to the endpoint
background in Section~\ref{sec:iteration-closure}.
\section{Choice of Parameters}\label{sec:parameter-choice}

We now choose the parameters left unspecified in the preceding
sections. They must make the stresses in
Section~\ref{sec:endpoint-stress} smaller than $\delta_{q+2}$ and
bound the magnetic increments by
$C\delta_{B,q+1}^{1/2}$, with all derivatives prescribed in
Section~\ref{ssec:fixed-scales}.

The order of the choices is as follows. First choose $b$,
$\gamma_a$, and $\gamma_\parallel$, then the integers $j_0,j_1$.
Next choose the spatial mollification exponent $\gamma_\ell$
and then the H\"older-loss exponent $\gamma_S$ sufficiently small.
Choose the approximation orders for the local construction and the
oscillatory inverse divergence, followed by the number of spatial
derivatives. Finally choose $a$ sufficiently large to absorb all
fixed constants. We verify each of these choices below.

\subsection{Scales and stress sizes}

The background distortions in \eqref{part:flow-smallness} are controlled by
\[
 \tau_c\lambda_q\delta_q^{1/2}=\varepsilon_{q+1}^{\gamma_\ell},\qquad
 \lambda_\parallel^{-1}\lambda_q\delta_{B,q}^{1/2}
 =\varepsilon_{q+1}^{\gamma_\ell+[\gamma_a-(b-1)\gamma_\parallel]/b}.
\]
The conditions $(b-1)\gamma_\parallel<\gamma_a$ and
$0<3\gamma_S<\gamma_\ell$ make these distortions small even after
H\"older losses: since $\ell^{-\alpha}\le\varepsilon_{q+1}^{-\gamma_S}$,
these quantities still tend to zero after multiplication by
$\ell^{-\alpha}$, as required by the chart and Galbrun estimates. Together with \eqref{setup:scales}, this
gives an admissible range of chart radii.

We first determine the conditions on
$\tau_a=\varepsilon_{q+1}^{\gamma_a}\tau_c$ imposed by the
antiderivative error and by a fast time derivative of the principal
perturbation. Using \eqref{setup:scale-identities}, write the
stress bounds as powers of $\varepsilon_{q+1}$. We first set
$\gamma_\ell=\gamma_S=0$, and then choose $\gamma_\ell>0$
and $\gamma_S>0$, in this order, sufficiently small to preserve
the strict inequalities.

The proof of Proposition~\ref{stress:complete-bounds} requires
the following four quantities to be smaller than
$\lambda_{q+1}^{-c\alpha}\delta_{q+2}$:
\begin{equation}\label{iter:scalar-budgets}
 \varepsilon_{q+1}^{j_0\gamma_a}\delta_{q+1},\qquad
 \frac{\tau_a^{-1}\delta_{q+1}^{1/2}}{\lambda_{q+1}},\qquad
 (\ell\lambda_{q+1})^{-1}\delta_{q+1},\qquad
 \ell^{-2\alpha}\tau_a^2\ell^{-2}\delta_{q+1}^2.
\end{equation}
These are, respectively, the cutoff, antiderivative and
flow-smoothing errors; the principal linear stress; spatial oscillation
and the first Lie variation (Proposition~\ref{stress:first-lie-variation});
and the second variation and covariance.
In the linear stress, the gain from inverse divergence compensates for
the fast time derivative. The last size is quadratic in
$\tau_a\ell^{-1}\delta_{q+1}^{1/2}$ and also covers the corrector
deformation $\ell^{-2\alpha}\tau_a^2\lambda_q^2\delta_{q+1}^2$ in
\eqref{gg:actual-stress}, since $\lambda_q\le\ell^{-1}$.
The factor $\lambda_{q+1}^{-c\alpha}$ absorbs the losses from Fourier
multipliers and from the comparison of transport and Lie derivatives
in \eqref{stress:prepared-complete}.

By \eqref{setup:scale-identities},
\[
 \frac{\lambda_{q+1}^{-c\alpha}\delta_{q+2}}{\delta_{q+1}}
 =\varepsilon_{q+1}^{2b\beta+c\gamma_S},\qquad
 \ell^{-\alpha}\le\varepsilon_{q+1}^{-\gamma_S},\qquad
 \tau_a\ell^{-1}\delta_{q+1}^{1/2}
 =\varepsilon_{q+1}^{\gamma_a+\beta}.
\]
After division by $\delta_{q+1}$, the four quantities in
\eqref{iter:scalar-budgets} are bounded by powers of
$\varepsilon_{q+1}$ with exponents
\begin{equation}\label{iter:exponent-list}
 j_0\gamma_a,\qquad 1-\beta-\gamma_a-\gamma_\ell,\qquad
 1-\gamma_\ell,\qquad 2\gamma_a+2\beta-2\gamma_S.
\end{equation}
At $\gamma_\ell=\gamma_S=0$, the first two errors require the
strict inequalities
\begin{equation}\label{iter:principal-margins}
 j_0\gamma_a-2b\beta>0,\qquad
 1-(2b+1)\beta-\gamma_a>0.
\end{equation}
Equivalently,
\[
 \frac{2b\beta}{j_0}<\gamma_a<1-(2b+1)\beta.
\]
This interval is nonempty precisely when
\[
 (2b+1+2b/j_0)\beta<1.
\] Its limiting range as $j_0\to\infty$ and
$b\downarrow1$ is $\beta<1/3$. The spatial-oscillation and quadratic
errors impose, respectively,
\[
 2b\beta<1,\qquad (b-1)\beta<\gamma_a.
\]
In particular, the upper bound on the fast-time exponent comes from
the linear stress, whereas the antiderivative error gives a lower
bound on its product with $j_0$.

The same parameters must separate the local derivative costs from
those of the next iterate. The relevant ratios are
\begin{align*}
 \frac{\tau_a^{-1}}{\lambda_q}
 &=\lambda_q^{-\beta+(b-1)(\gamma_a+\gamma_\ell)},\\
 \frac{\tau_a^{-1}}{\lambda_{q+1}\delta_{q+1}^{1/2}}
 &=\varepsilon_{q+1}^{1-\beta-\gamma_a-\gamma_\ell},\\
 \frac{\lambda_\parallel}{\lambda_{q+1}\delta_{B,q+1}^{1/2}}
 &=\varepsilon_{q+1}^{1-\beta-\gamma_a-\gamma_\ell}.
\end{align*}
For the first ratio to vanish and for each of the last two to
decay faster than $\varepsilon_{q+1}^{\gamma_a}$ before including
H\"older losses, the conditions at $\gamma_\ell=\gamma_S=0$ are
\[
 (b-1)\gamma_a<\beta,\qquad 2\gamma_a<1-\beta.
\]
The last two ratios of derivative costs coincide because
$\delta_{B,q+1}^{1/2}=\tau_a\lambda_\parallel\delta_{q+1}^{1/2}$. The remaining small-deformation
comparisons follow directly from
\begin{align*}
 \tau_a\ell^{-1}\ell^{-\alpha}\delta_q^{1/2}
 &\le\varepsilon_{q+1}^{\gamma_a-\gamma_S},\\
 \tau_a^2\lambda_q^2\ell^{-\alpha}\delta_{q+1}
 &=\ell^{-\alpha}\varepsilon_{q+1}^{2(\gamma_a+\gamma_\ell)}
       \frac{\delta_{q+1}}{\delta_q}.
\end{align*}

\begin{lemma}[Scale comparisons]\label{iter:parameter-lemma}
Let $0<\beta<1/3$, $b>1$, $0<(b-1)\gamma_\parallel<\gamma_a$ and
an integer $j_0\ge3$ satisfy \eqref{iter:principal-margins} and
\[
 2b\beta<1,\qquad (b-1)\beta<\gamma_a,\qquad
 (b-1)\gamma_a<\beta,\qquad 2\gamma_a<1-\beta.
\]
Fix an integer $j_1\ge0$ and set $c=64+14j_0+8j_1$.
For every sufficiently small $\gamma_\ell>0$, and then every
sufficiently small $\gamma_S>0$, the scales satisfy
$\lambda_q<\ell^{-1}<\lambda_{q+1}$, the derivative-cost comparisons
\eqref{setup:rate-comparisons}, and
\begin{equation}
 \tau_a\ell^{-1}\ell^{-\alpha}\delta_q^{1/2}=o(1).
\end{equation}
Each stress norm bound in \eqref{iter:scalar-budgets} is
$o(\lambda_{q+1}^{-c\alpha}\delta_{q+2})$. All limits are uniform
in $q$ as $a\to\infty$.
\end{lemma}
\begin{proof}
Set $\gamma_\ell=\gamma_S=0$. By the hypotheses, every exponent
in \eqref{iter:exponent-list} is strictly greater than $2b\beta$.
The first derivative-cost ratio is a negative power of
$\lambda_q$, and the exponents in the other two ratios are
strictly greater than $\gamma_a$. Choose $\gamma_\ell>0$
sufficiently small to preserve these inequalities, and then choose
$\gamma_S>0$ sufficiently small that the stress exponents remain
greater than $2b\beta+c\gamma_S$. This proves the stress
comparisons and the stated bounds for the three derivative-cost
ratios. The two deformation estimates above also tend to zero.

The remaining inequalities in the derivative-cost comparisons are
$\varepsilon_{q+1}\le\varepsilon_{q+1}^{\gamma_a}$ and
$\ell^{-\alpha}\le\lambda_{q+1}^\alpha$. They follow from
$\gamma_a<1$ and \eqref{setup:scale-identities}. Taking
$0<\gamma_\ell<1$ gives the strict ordering of the spatial scales.

For use in the local construction and the Galbrun estimates, we
make the same choices sufficiently small to satisfy also
\begin{equation}\label{iter:small-gammaS}
 \begin{gathered}
 1-\beta-\gamma_a-\gamma_\ell\ge\tfrac12(1-\beta-\gamma_a),\qquad
 \gamma_a>4\gamma_\ell,\\
 \gamma_\ell\le\frac1{2(j_1+2)},\qquad
 4\gamma_\ell+(c+10)\gamma_S<1,\\
 0<3\gamma_S<\gamma_\ell.
 \end{gathered}
\end{equation}
Indeed, $\gamma_\ell$ may first be chosen to satisfy the first
four inequalities at $\gamma_S=0$. The hypotheses on the fixed
parameters allow such a choice. A sufficiently small positive
$\gamma_S$ then gives all the stated bounds. In particular, $\gamma_S<\gamma_\ell$ absorbs the Schauder
factor in the mollification gain:
\[
 \ell^{1-\alpha}\lambda_q,
 \quad \tau_c\ell^{-\alpha}\lambda_q\delta_q^{1/2}
 \le\varepsilon_{q+1}^{\gamma_\ell-\gamma_S}=o(1).
\]
The bounds $\gamma_S<\gamma_\ell<\gamma_a$ and
$(b-1)\gamma_a<\beta$, together with
$\gamma_\ell+\gamma_S<1$, also give
\[
 \ell^{-\alpha}\delta_q^{1/2}
 \le\varepsilon_{q+1}^{\beta/(b-1)-\gamma_S}=o(1),\qquad
 \frac{\lambda_{q+1}^{\alpha}}{\ell\lambda_{q+1}}
 =\varepsilon_{q+1}^{1-\gamma_\ell-\gamma_S}=o(1).
\]
These inequalities give the field bounds at integer derivative orders,
the local flow estimates, and the estimates at the spatial frequency
of the next iterate.
Since $\varepsilon_{q+1}\le a^{-(b-1)}$, all the limits are uniform in $q$.
\end{proof}

\begin{remark}\label{iter:first-order-warning}
The bound for the second variation and covariance is essential. Replacing
$\tau_a^2\ell^{-2}\delta_{q+1}^2$ by $\tau_a\ell^{-1}\delta_{q+1}^{3/2}$
would require, at zero H\"older loss, $\gamma_a+\beta\ge2b\beta$;
together with the fast time comparison this gives $\beta<1/(4b)$,
whose limit as $b\downarrow1$ is $1/4$. A first Lie variation must
therefore be cancelled or estimated as an oscillatory divergence before
the scale comparisons are applied.
\end{remark}

\subsection{Conditions on the magnetic scale}

The magnetic increment has amplitude $\delta_{B,q+1}^{1/2}$, so its
square imposes an additional condition on the momentum stress. At
the highest orders of transport differentiation, the estimates for the
approximate antiderivative give less cancellation. The smaller local
derivative costs compensate for this loss when the estimates are
expressed at the costs of the next iterate. We choose the magnetic scale and the derivative orders to
satisfy both requirements.

\begin{lemma}[Anisotropic comparisons]\label{aniso:parameter-budgets}
Fix $0<\beta<1/3$ and $1<b<9/8$ with $3b\beta<1$, and let
$\gamma_a,\gamma_\parallel$ satisfy \eqref{iter:exponent-window}.
There exist integers $j_0\ge3$ and $j_1\ge j_0+2$ such that, with
$c=64+14j_0+8j_1$, for every sufficiently small
$\gamma_\ell\in(0,1)$ and then every sufficiently small
$\gamma_S\in(0,\gamma_\ell/3)$, the conclusions of
Lemma~\ref{iter:parameter-lemma} hold together with the following
comparisons:
\begin{align}
 \delta_{q+1}(\tau_a\lambda_\parallel)^2
 &=o(\lambda_{q+1}^{-c\alpha}\delta_{q+2}),
 \label{aniso:magnetic-stress-budget}\\
 \tau_a\lambda_\parallel\delta_{q+1}^{1/2}
 &=\delta_{B,q+1}^{1/2}.
 \label{aniso:magnetic-field-budget}
\end{align}
For every $0\le m\le j_1$ and $s\in\{0,1,2\}$,
\begin{equation}\label{aniso:finite-depth-budget}
 \min\{j_0,j_1-m-s\}\gamma_a+m\gamma_a
 \ge j_0\gamma_a.
\end{equation}
For $h=k+m\ge j_1-4$, we have
\begin{equation}\label{aniso:corrector-final-comparison}
 \varepsilon_{q+1}^{-2\gamma_a-\gamma_\parallel}
 \Bigl(\frac{\tau_a^{-1}}{\lambda_{q+1}\delta_{q+1}^{1/2}}\Bigr)^k
 \Bigl(\frac{\lambda_\parallel}{\lambda_{q+1}\delta_{B,q+1}^{1/2}}\Bigr)^m
 \le\varepsilon_{q+1}^{j_0\gamma_a+1}.
\end{equation}
For $0\le h=k+m\le j_1$ and
$N=\min\{j_0,\max\{0,j_1-h-2\}\}$,
\begin{equation}\label{aniso:collar-depth-comparison}
 \varepsilon_{q+1}^{N\gamma_a}
 \Bigl(\frac{\tau_a^{-1}}{\lambda_{q+1}\delta_{q+1}^{1/2}}\Bigr)^k
 \Bigl(\frac{\lambda_\parallel}{\lambda_{q+1}\delta_{B,q+1}^{1/2}}\Bigr)^m
 \le\varepsilon_{q+1}^{j_0\gamma_a}.
\end{equation}
Finally, the four auxiliary inequalities
\begin{equation}\label{iter:auxiliary-margins}
 \begin{gathered}
 2-\gamma_S>2b\beta+c\gamma_S,\qquad
 2-\beta-\gamma_S>2b\beta+c\gamma_S,\\
 1+\gamma_a+\gamma_\ell-\gamma_S>2b\beta+c\gamma_S,\\
 1-\beta-\gamma_a-\gamma_\ell-\gamma_S>2b\beta+c\gamma_S
 \end{gathered}
\end{equation}
hold.
All limits are uniform in $q$ as $a\to\infty$.
\end{lemma}
\begin{proof}
\emph{1. Conditions on the magnetic scale.}
Fix the exponents $\beta,b,\gamma_a,\gamma_\parallel$ in the
statement. The definitions give
\[
 \tau_a\lambda_\parallel=\varepsilon_{q+1}^{\gamma_a+\gamma_\parallel},
 \qquad
 \frac{\tau_a\lambda_\parallel\delta_{q+1}^{1/2}}
      {\delta_{B,q+1}^{1/2}}=1.
\]
Thus the magnetic square requires
$\gamma_a+\gamma_\parallel>b\beta$ at zero H\"older loss.
The second identity is the normalization of the magnetic increment;
its integer derivative estimates retain this amplitude, with fixed
constants. The H\"older estimates include the factor in
\eqref{setup:class}.

The first magnetic derivative of the old magnetic field has size
\[
 \lambda_q\delta_{B,q}
 =\lambda_q^{1-2\beta-2(b-1)(\gamma_a+\gamma_\parallel)/b}.
\]
For its old contribution to be smaller than the next inductive bound
in Step~4 of Proposition~\ref{full:endpoint-derivative-closure}, we
require
\[
 \beta+\frac{b-1}{b}(\gamma_a+\gamma_\parallel)<\frac12.
\]
The same condition makes each old contribution to the first material and magnetic
derivatives of the fields smaller than its new bound.

\smallskip
\noindent\emph{2. Choice of the number of transport derivatives.}
By \eqref{setup:rate-comparisons}, both ratios of local derivative
costs to those of the next iterate are
$\varepsilon_{q+1}^{1-\beta-\gamma_a-\gamma_\ell}$, bounded by
$\varepsilon_{q+1}^{\gamma_a}$. To prove the estimate for the
approximate antiderivative, choose
\begin{equation}\label{iter:depth-anisotropy}
 j_1\ge j_0+2.
\end{equation}
Indeed, if $j_1-m-s\ge j_0$, then
\eqref{aniso:finite-depth-budget} is immediate; otherwise its left-hand
side equals $(j_1-s)\gamma_a\ge j_0\gamma_a$. The exponent inequality
also covers $j_1-m-s<0$.

For compositions of $h\ge j_1-4$ material and magnetic derivatives,
the coarse Galbrun estimate loses a factor
\[
 \frac{\tau_c}{\tau_a^2\lambda_\parallel}
 =\varepsilon_{q+1}^{-2\gamma_a-\gamma_\parallel}
\]
compared with the sharp bound for the first magnetic derivative of
the Galbrun potential. Expressing each derivative at the costs of
the next iterate gives an exponent at least $(1-\beta-\gamma_a)/2$, by
\eqref{iter:small-gammaS}. We require
\begin{equation}\label{aniso:corrector-derivative-reserve}
 (j_1-4)\frac{1-\beta-\gamma_a}{2}
 >\frac{\beta}{b-1}+\frac{\gamma_a+\gamma_\parallel}{b}
       +10+(j_0+6)\gamma_a+2b\beta.
\end{equation}
Since $(\gamma_a+\gamma_\parallel)/b>\gamma_\parallel$ follows
from $(b-1)\gamma_\parallel<\gamma_a$, this implies
\eqref{aniso:corrector-final-comparison} and also covers the
magnetic gap loss in the derivative estimates. Increasing $j_1$
provides the transport derivatives required by
Theorem~\ref{galbrun:linear-theorem}(f).

\smallskip
\noindent\emph{3. The auxiliary inequalities.}
For spatial mollification and the products of gaps with increments,
including those containing the large spatial gradient of the
principal LDF, we must verify
\eqref{iter:mollifier-budget}--\eqref{iter:principal-gap-budgets}.
The differentiated mollification error in
Section~\ref{ssec:derivative-closure} also requires the last
auxiliary inequality. Set $\gamma_\ell=\gamma_S=0$.
The inequality $2b\beta<1$ and \eqref{iter:principal-margins}
give all four comparisons:
\[
 2>2b\beta,\qquad 2-\beta>2b\beta,\qquad
 1+\gamma_a>2b\beta,\qquad
 1-(2b+1)\beta-\gamma_a>0.
\]

\smallskip
\noindent\textbf{Choice of exponents.}
The fixed value of $b$ satisfies
\[
 1<b<\frac98,\qquad 3b\beta<1.
\]
The fixed exponents $\gamma_a,\gamma_\parallel$ satisfy
\begin{equation}\label{iter:exponent-window}
 \begin{gathered}
 (b-1)\beta<\gamma_a<\min\Bigl\{1-(2b+1)\beta,\ \frac{1-\beta}2\Bigr\},
 \qquad (b-1)\gamma_a<\beta,\\
 \gamma_a+\gamma_\parallel>b\beta,\qquad
 0<(b-1)\gamma_\parallel<\gamma_a,\\
 (b-1)(\gamma_a+\gamma_\parallel)<b(1/2-\beta).
 \end{gathered}
\end{equation}
We verify that such exponents exist. The interval for $\gamma_a$
is nonempty by $3b\beta<1$ and $b<2$. A choice sufficiently
close to its lower endpoint also satisfies
$(b-1)\gamma_a<\beta$. To choose $\gamma_\parallel$, observe
that
\[
 \max\Bigl\{(b-1)\beta,\frac{b-1}{b}\gamma_a\Bigr\}
 <\min\{\gamma_a,1/2-\beta\}.
\]
The first lower bound is smaller than both upper bounds by
$\gamma_a>(b-1)\beta$ and $b\beta<1/2$; the second is smaller than
$\gamma_a$, and is smaller than $1/2-\beta$ since $b<9/8$ and
$\gamma_a<(1-\beta)/2$. Choosing
$(b-1)(\gamma_a+\gamma_\parallel)/b$ between these bounds gives
precisely the magnetic conditions in \eqref{iter:exponent-window}.
Thus the same exponents satisfy the stress, chart and derivative-cost
conditions, including those for the estimates of material and magnetic derivatives of the fields.

Every positive magnetic regularity gain below
$\min\{1-3\beta,1/2-\beta\}$ is attainable. Choose $b-1$ sufficiently
small at the outset, then $\gamma_a$ above the prescribed gain and
below the upper bounds in \eqref{iter:exponent-window}. The two lower
bounds in the preceding display tend to zero as $b\downarrow1$,
so $\gamma_\parallel$ can be chosen with the prescribed value of
$(b-1)(\gamma_a+\gamma_\parallel)/b$.

We now make the remaining choices for the fixed exponents in the
statement. Choose an integer $j_0\ge3$ such that
$j_0\gamma_a>2b\beta$. Choose $j_1$ sufficiently large to satisfy
both \eqref{iter:depth-anisotropy} and
\eqref{aniso:corrector-derivative-reserve}. This is possible since
$1-\beta-\gamma_a>0$. Fix $c=64+14j_0+8j_1$. By
Lemma~\ref{iter:parameter-lemma}, we can next choose
$\gamma_\ell>0$, and then $\gamma_S>0$, sufficiently small
that \eqref{iter:small-gammaS},
\eqref{iter:auxiliary-margins}, and
\[
 2(\gamma_a+\gamma_\parallel)>2b\beta+c\gamma_S
\]
hold.
The stress and derivative cost inequalities are strict at
$\gamma_\ell=\gamma_S=0$. The preparation and Schauder inequalities
in \eqref{iter:small-gammaS} hold by taking $\gamma_S$ sufficiently
small after $\gamma_\ell$. They give
\eqref{aniso:magnetic-stress-budget}, while
\eqref{aniso:magnetic-field-budget} follows from the definition of
$\delta_{B,q+1}$ without a loss.

For \eqref{aniso:corrector-final-comparison}, subtract the coarse
Galbrun loss from \eqref{aniso:corrector-derivative-reserve}:
\[
 (j_1-4)\frac{1-\beta-\gamma_a}{2}
 -2\gamma_a-\gamma_\parallel
 >\frac{\beta}{b-1}+\frac{\gamma_a-(b-1)\gamma_\parallel}{b}
       +10+(j_0+4)\gamma_a+2b\beta
 >j_0\gamma_a+1.
\]
For the cutoff error, \eqref{iter:depth-anisotropy} gives $j_1\ge j_0+2$, so
\[
 N+h\ge\min\{j_0,j_1-2\}=j_0.
\]
Each ratio of local derivative costs to those of the next iterate is at most $\varepsilon_{q+1}^{\gamma_a}$. Thus the
left-hand side of \eqref{aniso:collar-depth-comparison} is bounded by
$\varepsilon_{q+1}^{(N+h)\gamma_a}$, whose exponent satisfies
\[
 (N+h)\gamma_a\ge j_0\gamma_a.
\]
This proves \eqref{aniso:collar-depth-comparison}. The choices of
$\gamma_\ell$ and $\gamma_S$ give all the asserted inequalities.
The approximation orders and the number of spatial derivatives
will be chosen below; they do not change these inequalities.
\end{proof}

\begin{remark}[Admissible magnetic regularity exponents]
The magnetic regularity gain
$(b-1)(\gamma_a+\gamma_\parallel)/b$ is subject to the conditions on charts,
material and magnetic derivatives of the fields, and stresses in \eqref{iter:exponent-window}.
A prescribed positive gain as $b\downarrow1$ requires
$\gamma_\parallel$ of order $(b-1)^{-1}$. This dependence is
included in \eqref{aniso:corrector-derivative-reserve} and
\eqref{stress:loss-maximum}, which determine the required number of
derivatives. The regularity
exponents in space and time and along magnetic field lines follow from these amplitude bounds and derivative costs in
Section~\ref{sec:initialization-and-limit}.
\end{remark}

\begin{remark}[Magnetic derivatives]
To apply these scalar comparisons, the differentiated estimates must
keep track of the number of magnetic derivatives on each remainder.
Lemma~\ref{aniso:old-frame-operator-change} separates the gradient
factors in the tensor Lie derivatives from derivatives of the scalar
coefficient.
\end{remark}

\subsection{Approximation error bounds}

We first choose the order of the spatial mollifier used in the
local construction. To estimate material and magnetic derivatives
with their separate slow costs and to retain the smaller magnetic
amplitude, we impose \eqref{prep:accuracy-inequality}:
\begin{equation}\label{iter:accuracy-requirement}
 \varepsilon_\ell:=(\ell\lambda_q)^{m_0}
 \le\ell^{2\alpha}\varepsilon_{q+1}^2
        \delta_{B,q}^{(j_1+3)/2}.
\end{equation}
To check the exponent on the right, divide each slow material or
magnetic derivative cost by $\lambda_q$. Each ratio is at least
$\delta_{B,q}^{1/2}$. In passing from ordinary derivative
estimates to local transport estimates,
Lemma~\ref{prep:finite-conversion} and
Corollary~\ref{prep:mixed-comparison} use at most $j_1+1$ such
ratios. One further half power bounds the magnetic amplitude ratio
or the normalization in \eqref{prep:local-acceleration}, and the last
half power gives the strict inequality in
\eqref{prep:conversion-margin}. Substitute
$\ell\lambda_q=\varepsilon_{q+1}^{\gamma_\ell}$ and
$\delta_{B,q}^{1/2}=\varepsilon_{q+1}^{\beta/(b-1)+(\gamma_a+\gamma_\parallel)/b}$
and choose
\begin{equation}\label{iter:ordinary-accuracy}
 m_0=1+\left\lceil
 \frac{3+(\beta/(b-1)+(\gamma_a+\gamma_\parallel)/b)(j_1+3)}{\gamma_\ell}\right\rceil .
\end{equation}
Indeed, this choice gives
\[
 \varepsilon_\ell=\varepsilon_{q+1}^{m_0\gamma_\ell}
 \le\varepsilon_{q+1}^3\delta_{B,q}^{(j_1+3)/2}.
\]
The remaining factor $\varepsilon_{q+1}$ is sufficient for the
factor $\ell^{2\alpha}$ in \eqref{iter:accuracy-requirement}, since
\[
 \ell^{2\alpha}
 =\varepsilon_{q+1}^{2\gamma_S[1+(b-1)\gamma_\ell]/b}
\]
and the exponent on the right is less than one by
\eqref{iter:small-gammaS}.
This proves the mollification condition with the same Taylor order
as the local correction. In
Lemma~\ref{prep:ordinary-comparison}, the Taylor estimate must
first give $\varepsilon_\ell^2$; we therefore choose the kernel
with vanishing moments through $2m_0+2$. One approximation factor
gives the improved magnetic comparison, and the second controls
the material and magnetic derivatives of the coefficients and gaps.

Take $\dstar=2m_0+8$ as in \eqref{setup:prepared-ceiling}.
This is sufficient for all additional ordinary derivatives used
in Section~\ref{sec:preparation}: the Taylor estimate for the
mollification errors in Step~3 of the proof of
Lemma~\ref{prep:finite-conversion} uses $2m_0+4$ additional
derivatives; the refined comparison of
Lemma~\ref{prep:ordinary-comparison}, used for material and
magnetic derivatives of the gaps, uses $2m_0+5$. The remaining
three derivatives give the H\"older estimates for the local
coefficients.

The mollification error enters the stress only through its order-zero
size, and the background gaps through their products with the
increments. The required scalar comparisons are
\begin{equation}\label{iter:mollifier-budget}
 \ell^{-\alpha}\varepsilon_{q+1}^2\delta_{q+1}
 =o(\lambda_{q+1}^{-c\alpha}\delta_{q+2})
\end{equation}
and, since Corollary~\ref{principal:vector-transport} bounds the
principal pushforward of an unmollified gap of size
$\ell^{-\alpha}\varepsilon_{q+1}^2\delta_q^{1/2}$ by the single
additional factor $\tau_a\lambda_{q+1}\delta_{q+1}^{1/2}$,
\begin{equation}\label{iter:principal-gap-budgets}
 \begin{aligned}
 (1+\tau_a\lambda_{q+1}\delta_{q+1}^{1/2})\delta_{q+1}^{1/2}
       \ell^{-\alpha}\varepsilon_{q+1}^2\delta_q^{1/2}
   &=o(\lambda_{q+1}^{-c\alpha}\delta_{q+2}),\\
 \tau_a\lambda_{q+1}\delta_{q+1}^{1/2}\ell^{-\alpha}
       \varepsilon_{q+1}^2\delta_q^{1/2}
   &=o(\delta_{q+1}^{1/2}).
 \end{aligned}
\end{equation}
After division by $\delta_{q+1}$, the exponents are respectively at
least $2-\gamma_S$, $2-\beta-\gamma_S$ (without the additional factor) and
$1+\gamma_a+\gamma_\ell-\gamma_S$ (with it), all exceeding $2b\beta+c\gamma_S$ by
\eqref{iter:auxiliary-margins}; the field comparison, divided by
$\delta_{q+1}^{1/2}$, has the last positive exponent. These are used in
\eqref{stress:prepared-complete} and
\eqref{full:endpoint-increment-bounds}.

Next choose the number of terms in the oscillatory parametrix.
Each integration of the phase profile in Lemma~\ref{osc:parametrix}
gives $(\ell\lambda_{q+1})^{-1}$. For the final remainder,
use the integrable kernel of the inverse divergence and the
commutator estimate in that lemma. This retains the material and
magnetic derivative costs, but does not give
$\lambda_{q+1}^{-1}$. Thus comparison with the desired stress
bound introduces $\lambda_{q+1}$, and we impose
\eqref{osc:reserve}:
\begin{equation}\label{iter:oscillatory-reserve}
 \lambda_{q+1}(\ell\lambda_{q+1})^{-\rfg}
 =o(1).
\end{equation}
By \eqref{setup:scale-identities}, the left-hand side is
$\lambda_q^{b-\rfg(b-1)(1-\gamma_\ell)}$. Thus choose the integer
$\rfg$ so that
\begin{equation}
 \rfg(b-1)(1-\gamma_\ell)>b,
\end{equation}
which makes \eqref{iter:oscillatory-reserve} tend to zero. The
parametrix requires $\rfg+3$ additional spatial derivatives. With $m_0$
and $\rfg$ fixed, we can now choose the maximum derivative orders.

\subsection{Choice of spatial derivative orders}\label{ssec:derivative-reserve}

We finish by choosing the number of spatial derivatives so that
Lemma~\ref{setup:calculus}(vii) applies to every stress term
requiring an extension of its sharp estimates. The lossy estimates
hold at all spatial orders, subject to the stated restrictions
on material and magnetic derivatives, except when a background
gap occurs. In that case they hold through at least
$\rgood+1$. The smoothing error needs no such argument, since
its derivatives were estimated separately in
Section~\ref{ssec:derivative-closure}.

To extend the sharp estimates through $\rgood$, it is enough
to verify \eqref{setup:reserve-inequality} for the smallest of the
derivative orders through which sharp estimates have been proved,
$\rcut$, and the largest loss $\varepsilon_{q+1}^{-\gres}$. By the table in
Section~\ref{ssec:derivative-closure}, these are
\begin{equation}\label{iter:common-switch}
 \begin{aligned}
 \rcut&=\rprep-j_0-\rfg-16,\\
 \gres&=2+3\gamma_a+\frac{\beta}{b-1}+\frac{\gamma_a+\gamma_\parallel}{b}
        +j_0\gamma_a+2b\beta+4\gamma_\ell+(c+3)\gamma_S.
 \end{aligned}
\end{equation}
These are \eqref{stress:ceiling-minimum} and \eqref{stress:loss-maximum}.
The derivative ranges follow by tracing the additional derivatives in
Lemma~\ref{setup:calculus} and the appendices. Preparation uses
$\dstar=2m_0+8$ derivatives of the old iterate beyond the asserted
range. The Galbrun estimates require ten further derivatives of the
coefficients and forcing, while curl and the corrector pushforward
each require one spatial derivative. Thus the corrector bounds hold
through $\rprep-12$.

The approximate antiderivative requires $j_0$ additional material
derivatives of its coefficient. The transport lemmas in charts in
Subsection~\ref{ssec:path-estimates} require two further spatial
derivatives. All principal quantities can therefore be estimated
through order
\[
 \Nc=\rprep-j_0-12\le\rprep-12.
\]
Finally, the oscillatory inverse divergence requires $\rfg+4$
additional spatial derivatives, including the one used to obtain the H\"older estimate. This gives
\[
 \rcut=\Nc-\rfg-4.
\]
We also check the loss by dividing the lossy estimate of each
factor by its sharp estimate in
Sections~\ref{sec:preparation}--\ref{sec:endpoint-stress}.
The antiderivative and magnetic cancellation estimates give
$\varepsilon_{q+1}^{-j_0\gamma_a}$ and
$\varepsilon_{q+1}^{-2\gamma_a}$; the coarse Galbrun estimate
gives $\varepsilon_{q+1}^{-2\gamma_a-\gamma_\parallel}$.
For a gap, the approximation factor and the magnetic amplitude
give $\varepsilon_{q+1}^{-2}(\delta_q/\delta_{B,q})^{1/2}$.
The ratios of spatial derivative costs contribute
$(\ell\lambda_q)^{-1}$ per derivative, and the ratio of lossy
and fast material derivative costs contributes $\tau_a\ell^{-1}$.
For the smoothing error, comparison with the required amplitude
gives $\varepsilon_{q+1}^{-2b\beta}$. Finally, each fixed
H\"older factor contributes $\varepsilon_{q+1}^{-\gamma_S}$.
There are at most $c$ such factors, giving $c\gamma_S$ in the
exponent. These are all the losses. Their exponents do not increase
with the spatial derivative order being estimated: the tame
product, chain-rule, and flow estimates put only one factor at the
highest order and estimate all other factors at fixed lower orders.

\smallskip
\noindent\emph{Absorbing the losses.}
We apply \eqref{setup:reserve-inequality} when the sum of the spatial,
material, and magnetic derivative orders is at most $\rcut$, the
sum of the material and magnetic derivative orders is at most
$J\le j_1+1$, and the loss is $\varepsilon_{q+1}^{-\gres}$. The factor
$(\ell\Lambda_1)^{-(\rcut-J)}$ must absorb both this loss and the
ratios of the lossy derivative costs to the material and magnetic
derivative costs of the next iterate. These ratios are
\[
 \frac{\ell^{-1}}{\lambda_{q+1}\delta_{q+1}^{1/2}}
 =\varepsilon_{q+1}^{1-\gamma_\ell-b\beta/(b-1)},\qquad
 \frac{\ell^{-1}}{\lambda_{q+1}\delta_{B,q+1}^{1/2}}
 =\varepsilon_{q+1}^{1-\gamma_\ell-b\beta/(b-1)-\gamma_a-\gamma_\parallel}.
\]
The magnetic factor is the larger one, so it controls every split
$k+m\le J\le j_1+1$. Two transverse scales occur. For expressions that
are not transported by the principal flow, $\Lambda_1=\lambda_{q+1}$
and $(\ell\lambda_{q+1})^{-1}=\varepsilon_{q+1}^{1-\gamma_\ell}$. For
expressions transported by the principal flow,
Corollary~\ref{principal:rescaled-profile-continuation} measures the
derivatives of the slow coefficients, source, and initial datum at
scale $\lambda_{q+1}^{-1}$ in the oscillation direction and at scale
$\vartheta\ell$ in the transverse directions. The transport bound requires
\[
 \frac{\tau_a\delta_{q+1}^{1/2}}{\vartheta\ell}\le1.
\]
Accordingly, \eqref{principal:continuation-data} takes the smallest
admissible $\vartheta$. By \eqref{setup:scale-identities},
\[
 \vartheta=\tau_a\ell^{-1}\delta_{q+1}^{1/2}
 =\varepsilon_{q+1}^{\gamma_a+\beta}
 \ge(\ell\lambda_{q+1})^{-1},
\]
where the last comparison follows from
$\tau_a\lambda_{q+1}\delta_{q+1}^{1/2}\ge1$ in
\eqref{setup:rate-comparisons}. A derivative falling on a slow
coefficient therefore gains $\vartheta$, and both cases are covered
by the single condition
\begin{equation}\label{iter:general-final-reserve}
 (\rcut-j_1-1)(\gamma_a+\beta)
 >\gres+
 (j_1+1)\left[\frac{b\beta}{b-1}+\gamma_a+\gamma_\parallel-(1-\gamma_\ell)\right]^+,
\end{equation}
This condition implies $\rcut\ge j_1+2$ and gives
\begin{equation}\label{iter:transported-profile-comparison}
 \varepsilon_{q+1}^{-\gres}
 \vartheta^{\rcut-j_1-1}
 \max\Bigl\{1,\frac{\ell^{-1}}{\lambda_{q+1}\delta_{q+1}^{1/2}},
 \frac{\ell^{-1}}{\lambda_{q+1}\delta_{B,q+1}^{1/2}}\Bigr\}^{j_1+1}
 =o(1),
\end{equation}
which is \eqref{setup:reserve-inequality} with
$\Lambda_1=(\ell\vartheta)^{-1}$ after increasing $a$; the case
$\Lambda_1=\lambda_{q+1}$ follows because $\vartheta\ge(\ell\lambda_{q+1})^{-1}$.

Since the coefficient of $\rcut$ in \eqref{iter:general-final-reserve}
is positive and its right-hand side is fixed, an admissible integer
exists. We take $\rcut$ to be the smallest such integer and set
\begin{equation}
 \rgood=\rcut+j_0+\rfg+16+\dstar,\qquad \dstar=2m_0+8,
\end{equation}
so that $\rprep=\rgood-\dstar$ and \eqref{iter:common-switch} hold. We
impose the ordinary field bounds through $\rgood+4$. No further real
parameter is needed. Every strict comparison retains
a positive power of $\varepsilon_{q+1}$. Since
$\varepsilon_{q+1}\le a^{-(b-1)}$, a sufficiently large $a$ absorbs
all fixed constants uniformly in $q$.
\section{Closing the Iteration}\label{sec:iteration-closure}

We now prove the induction step by combining the field and pressure
estimates with the support properties of the construction. The stress
bound is supplied by Proposition~\ref{stress:complete-bounds}.

\subsection{Velocity and magnetic field}
\label{ssec:fields-pressure-support}

We apply Proposition~\ref{full:endpoint-derivative-closure} to the ordinary derivatives
and to compositions of material and magnetic derivatives. For the magnetic
increment, the proposition separates the principal increment relative to the
local background, the corrector increment and the difference between the
backgrounds. Each has the required size by \eqref{full:field-reserve-quotient}
and Lemma~\ref{aniso:magnetic-gap-closure}.

For $r+h\ge1$, the old velocity and magnetic bounds are smaller than
the new inductive bounds by the respective factors
$\varepsilon_{q+1}^{r+h-\beta}$ and
$\varepsilon_{q+1}^{r+h-\beta-\frac{b-1}{b}(\gamma_a+\gamma_\parallel)}$.
This is the comparison in Step~5 of
Proposition~\ref{full:endpoint-derivative-closure}; the same proposition
provides the corresponding comparisons for compositions of material and magnetic derivatives.

We next choose $C_0$ and $a$. By Steps~4 and~5 of
Proposition~\ref{full:endpoint-derivative-closure}, the constants in the
principal estimates are independent of $C_0$. In particular, for large $a$,
the constants in \eqref{full:endpoint-increment-bounds} depend only on the
fixed derivative orders, profiles and geometric data. The product estimates
of Lemma~\ref{setup:calculus} use derivatives whose orders have sum at most
$\rgood+4$, with at most $j_1$ material and magnetic derivatives combined;
their constants and those from Leibniz' rule depend only on these fixed orders.
The reference fields have size $2+o(1)$, so the leading magnetic source
and the ordinary time derivatives also have constants independent of $C_0$.
All terms whose constants depend on $C_0$ contain a positive power of
$\varepsilon_{q+1}$. We may therefore first choose $C_0$ larger than a fixed
multiple of the constants in the principal estimates, and then choose $a$
so large that each $o(1)$ term is a fixed small fraction of the corresponding
inductive bound.

The triangle inequality now gives the field derivative bounds and
\eqref{full:increment-bound}. To obtain the pointwise bounds, increase
$a$ further so that
\[
 \begin{aligned}
 \|v_{q+1}-v_q\|_0,\ \|B_{q+1}-B_q\|_0
 &\le2(\delta_q^{1/2}-\delta_{q+1}^{1/2}),\\
 \|B_{q+1}-B_q\|_0
 &\le\tfrac14(\delta_q^{1/2}-\delta_{q+1}^{1/2}).
 \end{aligned}
\]
These inequalities follow uniformly in $q$ from the uniform increment
bounds and $\delta_{q+1}^{1/2}/\delta_q^{1/2}
=\varepsilon_{q+1}^\beta\to0$. Applying the triangle and reverse
triangle inequalities to \eqref{full:pointwise-induction} now gives
\[
 \|v_{q+1}\|_0,\ \|B_{q+1}\|_0
 \le2(1-\delta_{q+1}^{1/2}),\qquad
 |B_{q+1}|\ge\tfrac14(1+\delta_{q+1}^{1/2}).
\]

\subsection{Pressure}\label{ssec:iteration-pressure}

Write $p_{q+1}=p_q+\sum_n\pi_n$. We estimate its differential to avoid
requiring an additional derivative of the old pressure. Since exterior
differentiation commutes with Lie differentiation, at stage $q$ or $q+1$ we have
\[
 \mathcal D_t^k\mathcal L_B^m(dp)=d(D_t^kD_B^mp),\qquad
 \mathcal L_u(dp)=d(u\cn p),
\]
and each term introduced by changing the material or magnetic derivative
contains a spatial derivative of the pressure. Thus, to prove
\eqref{prep:old-pressure} at stage $q+1$, we estimate
$d(D_{t,q+1}^kD_{B_{q+1}}^mp_{q+1})$ in the norm of order $r-1$.
The spatial, material and magnetic derivative orders applied to
$dp_{q+1}$ have sum at most $\rgood-1$. We treat $dp_q$ and the
local pressure differentials separately;
for the latter, we will also use the bounds at all spatial derivative orders
and Lemma~\ref{setup:calculus}(vii).

By \eqref{prep:old-pressure}, $dp_q$ has amplitude
$C\lambda_q\delta_q$ and spatial frequency $\lambda_q$, where $C$ may
depend on $C_0$. This estimate holds when the spatial, material and magnetic
derivative orders have sum at most $\rgood-1$, and the material and magnetic
orders have sum at most $j_1$. Its material and magnetic derivative costs,
$\lambda_q\delta_q^{1/2}$ and $\lambda_q\delta_{B,q}^{1/2}$, are bounded
by those at the next iterate. Moreover, \eqref{prep:old-mixed} and
Lemma~\ref{high:transport-gradient} show that the old background is
$(\lambda_{q+1},\mathrm f)$-adapted through $(\rgood-1,j_1)$.
To compare its derivatives with those at the new iterate, we use
\[
 \begin{aligned}
 v_{q+1}-v_q&\in
 \mathcal C_{\rgood-2,j_1-1}(C\delta_{q+1}^{1/2};\lambda_{q+1},\mathrm f),\\
 B_{q+1}-B_q&\in
 \mathcal C_{\rgood-2,j_1-1}(C\delta_{B,q+1}^{1/2};\lambda_{q+1},\mathrm f)
 \end{aligned}
\]
with derivatives taken relative to the old background. Indeed, the
corresponding bounds relative to the new background follow by decomposing
the increments as in Proposition~\ref{full:endpoint-derivative-closure}.
The corrector integrals $c_{v,n,1},c_{B,n,1}$ and the gaps $G_n^v,G_n^B$
satisfy the same bounds by Lemma~\ref{setup:calculus}(vii) and
Section~\ref{ssec:derivative-closure}. Apply Lemma~\ref{setup:calculus}(iv)
with the new background fixed and increments $v_q-v_{q+1}$ and
$B_q-B_{q+1}$ to obtain the two displayed inclusions.

We now estimate $dp_q$ relative to the new background. First apply
Lemma~\ref{setup:calculus}(iii) to pass from the old Lie derivative
bounds to bounds for the componentwise material and magnetic
derivatives. Apply part (iv) with the old background fixed, then
part (iii) to obtain the new Lie derivative bounds. This gives
\[
 \begin{aligned}
 \|\mathcal D_{t,q+1}^k\mathcal L_{B_{q+1}}^m(dp_q)\|_{r-1}
 &\le C\lambda_q\delta_q\lambda_{q+1}^{r-1}
 (\lambda_{q+1}\delta_{q+1}^{1/2})^k
 (\lambda_{q+1}\delta_{B,q+1}^{1/2})^m,
 \end{aligned}
\]
for $r\ge1$, $k+m\le j_1$ and $r+k+m\le\rgood$. In this application,
the derivative orders on the increments have sum at most $\rgood-2$,
with at most $j_1-1$ material and magnetic derivatives combined. The
bounds for the old pressure are therefore used only at the orders
provided by the induction hypotheses.

For each local pressure, apply
Corollary~\ref{galbrun:split-pressure} at one higher spatial order
in \eqref{galbrun:split-pressure-bound}. Commuting exterior
differentiation with Lie derivatives gives, relative to the local
background on the time interval indexed by $n$,
\[
 d\pi_n\in\mathcal C_{\rprep-5}(C\ell^{-\alpha}\lambda_q\delta_{q+1};\lambda_q,\mathrm a)
 \cap\mathcal K(C\ell^{-1}\delta_{q+1}),
\]
where the second inclusion follows from \eqref{gg:pressure-bound}
at one higher spatial order. Expanding the material and magnetic
derivatives into ordinary derivatives uses at most $k$ time derivatives
and gives the relative loss
$(\ell\lambda_q)^{-1}\ell^{\alpha}\le\varepsilon_{q+1}^{-\gamma_\ell}$.
This is bounded by $\varepsilon_{q+1}^{-\gres}$, and $\rprep-5$ is at
least $\rcut$. Thus Lemma~\ref{setup:calculus}(vii) applies and gives
\[
 d\pi_n\in
 \mathcal C_\infty(C\ell^{-\alpha}\lambda_q\delta_{q+1};\lambda_{q+1},\mathrm f).
\]
By Section~\ref{ssec:derivative-closure}, the same estimates hold relative
to the new material and magnetic derivatives when the spatial, material
and magnetic orders sum to at most $\rgood-1$.
Both $dp_q$ and the local pressure differentials are smaller than the
required inductive bound, since
\[
 \frac{\lambda_q\delta_q}{\lambda_{q+1}\delta_{q+1}}
             =\varepsilon_{q+1}^{1-2\beta}=o(1),\qquad
 \frac{\ell^{-\alpha}\lambda_q\delta_{q+1}}{\lambda_{q+1}\delta_{q+1}}
             \le\varepsilon_{q+1}^{1-\gamma_S}=o(1).
\]
Here the H\"older factor is bounded by
$\ell^{-\alpha}\le\varepsilon_{q+1}^{-\gamma_S}$ and $\gamma_S<1$.
Since the time supports have bounded overlap, we obtain
\[
 \|d(D_{t,q+1}^kD_{B_{q+1}}^mp_{q+1})\|_{r-1}
 \le C_0^{k+m+2}\lambda_{q+1}^r\delta_{q+1}
 (\lambda_{q+1}\delta_{q+1}^{1/2})^k
                      (\lambda_{q+1}\delta_{B,q+1}^{1/2})^m,
\]
for $r\ge1$, $k+m\le j_1$ and $r+k+m\le\rgood$. This is
\eqref{prep:old-pressure} at stage $q+1$, with the pressure differential
estimated in the integer norm of order $r-1$.

Both vector fields are divergence free, so
\[
 \int D_{B_{q+1}}f\,\dd x=0,\qquad
 \int D_{t,q+1}f\,\dd x=\partial_t\int f\,\dd x.
\]
Since $p_q$ and each $\pi_n$ have zero mean, every
$D_{t,q+1}^kD_{B_{q+1}}^mp_{q+1}$ has zero mean. Poincar\'e's inequality gives
the required pressure norm from its differential.

\subsection{Support and exact equations}\label{ssec:iteration-support}

Every correction is supported within $4\tau_c$ of
$\supp\varrho$. The potentials vanish smoothly at the time endpoints,
and the elliptic operators act only in space. By
\eqref{prep:outer-cutoff}, the support of every increment and every new
stress term is therefore contained in
\[
 \left(\tfrac12+\tfrac14\delta_q^{1/2}-4\tau_c,
       \tfrac52-\tfrac14\delta_q^{1/2}+4\tau_c\right).
\]
The difference between the old stress and its spatial mollification is supported
in the old time interval, which is contained in this interval. The cutoff
construction in Section~\ref{sec:preparation} already imposes
\[
 4\tau_c+\tfrac52\delta_{q+1}^{1/2}
                                 <\tfrac14\delta_q^{1/2}.
\]
Indeed, $\tau_c/\delta_q^{1/2}=\ell/\delta_q
=\lambda_q^{-1+2\beta-(b-1)\gamma_\ell}\to0$ and
$\delta_{q+1}^{1/2}/\delta_q^{1/2}=\varepsilon_{q+1}^{\beta}\to0$.
The preceding inequality shows that the support interval above
is compactly contained in
\[
 \left(\tfrac12+\tfrac12\delta_{q+1}^{1/2},
       \tfrac52-\tfrac52\delta_{q+1}^{1/2}\right),
\]
which proves the next stress-support condition and the asserted
support of the increments in $(1/2,5/2)$.
The endpoint identity gives the relaxed momentum equation, while the
volume-preserving spacetime pushforward preserves incompressibility
and the induction equation. Together with the stress, field and
pressure estimates, this completes the proof of
Proposition~\ref{full:step}.

\section{Initialization and Passage to the Limit}\label{sec:initialization-and-limit}

The increment estimates give convergence to a weak solution and
determine its regularity. We also construct an initial relaxed solution
whose values outside the perturbation interval force the limiting total
energy and cross-helicity to vary.

\begin{proposition}[Passage from the iteration to a weak solution]
\label{iter:limit-proposition}
Let $(v_q,B_q,p_q,R_q)$ be smooth solutions of \eqref{eq:mhd} with
zero spatial mean pressure and
$v_0,B_0\in C_b^1(\mathbb T^3\times\R)$. Suppose a constant $C$,
independent of $q$ and $a$, satisfies
\begin{align}
 \|v_{q+1}-v_q\|_0
 &\le C\delta_{q+1}^{1/2},\qquad
 \|B_{q+1}-B_q\|_0\le C\delta_{B,q+1}^{1/2},
 \label{iter:increment-zero}\\
 \|\DD_{x,t}(v_{q+1}-v_q)\|_0
 &\le C\lambda_{q+1}\delta_{q+1}^{1/2},\qquad
 \|\DD_{x,t}(B_{q+1}-B_q)\|_0
 \le C\lambda_{q+1}\delta_{B,q+1}^{1/2},
 \label{iter:increment-one}\\
 \|R_q\|_0+\lambda_q^{-1}[R_q]_1
 &\le C\lambda_q^{-\alpha}\delta_{q+1}.
 \label{iter:stress-limit}
\end{align}
Assume also that every increment is supported in $1/2<t<5/2$.
Then $v_q,B_q$ converge in $C^\gamma(\mathbb T^3\times\R)$ for every
$0\le\gamma<\beta$ to a weak ideal MHD solution. The pressure
converges in $C_tC_x^\gamma$ on the same range. Moreover,
$B_q$ converges in $C^\gamma(\mathbb T^3\times\R)$ for every
$0\le\gamma<\beta+\frac{b-1}{b}(\gamma_a+\gamma_\parallel)$.
If the initial relaxed solution is the one constructed in
Lemma~\ref{iter:initializer-lemma}, then, for sufficiently large $a$,
the limiting magnetic field is nowhere zero and has zero spatial
mean and nonzero constant magnetic helicity. The limiting solution
fails to conserve both total energy and cross-helicity.
\end{proposition}

\subsection{An explicit initial iterate}

\begin{lemma}\label{iter:initializer-lemma}
For all sufficiently large $a$, the fields and stress in
\eqref{iter:initializer-fields}--\eqref{iter:initializer-stress}
satisfy every induction hypothesis at $q=0$.
\end{lemma}
\begin{proof}
Choose a fixed smooth function $\chi_{\rm seed}:\R\to[0,1]$ with
$\chi_{\rm seed}=1$ on $t\le2/3$ and $\chi_{\rm seed}=0$ on $t\ge1$.
We choose a magnetic field of fixed frequency and a velocity with two
components. The high frequency component makes the stress produced by the time cutoff
small after spatial antidifferentiation. The nonzero low frequency component,
parallel to the magnetic field, makes the cross-helicities at the two
endpoints different. More precisely, define
\begin{equation}\label{iter:initializer-fields}
 \begin{aligned}
 k_{\mathrm{in}}&=\lfloor\lambda_0\delta_0^{1/2}\rfloor,
 &c_{\mathrm{in}}&=\delta_1^2/k_{\mathrm{in}}^2,\\
 \bar B&=(\cos x_3,\sin x_3,0),
 &\bar v&=(\cos(k_{\mathrm{in}}x_3),\sin(k_{\mathrm{in}}x_3),0)
                       +c_{\mathrm{in}}\bar B,\\
 v_0(t,x)&=\chi_{\rm seed}(t)\bar v(x),
 &B_0(t,x)&=\bar B(x),\qquad p_0=0.
 \end{aligned}
\end{equation}
For large $a$, $k_{\mathrm{in}}\ge2$ and
$\frac12\lambda_0^{1-\beta}\le k_{\mathrm{in}}\le\lambda_0^{1-\beta}$.
Define the horizontal vector
\[
 A(x_3)=k_{\mathrm{in}}^{-1}
       (\sin(k_{\mathrm{in}}x_3),-\cos(k_{\mathrm{in}}x_3),0)
       +c_{\mathrm{in}}(\sin x_3,-\cos x_3,0).
\]
Thus $\partial_3A=\bar v$. Put
\begin{equation}\label{iter:initializer-stress}
 R_0=\chi_{\rm seed}'(A\otimes e_3+e_3\otimes A).
\end{equation}
Both vector fields in \eqref{iter:initializer-fields} are horizontal,
and their spatial dependence is only through $x_3$. Consequently
\[
 \ddiv v_0=\ddiv B_0=0,\qquad
 (v_0\cn)v_0=(v_0\cn)B_0=(B_0\cn)v_0=(B_0\cn)B_0=0,
\]
and
\[
 \ddiv R_0=\chi_{\rm seed}'\bar v=\partial_tv_0,\qquad
 \partial_tB_0=0.
\]
This proves the relaxed equations. The stress is symmetric and
trace-free, and its time support is contained in $[2/3,1]$.

\smallskip
\noindent\emph{Derivative bounds.}
On any component depending only on $(t,x_3)$ one has
$D_{t,0}=\partial_t$ and $D_{B_0}=0$. In particular,
\begin{equation}
 \begin{aligned}
 D_{t,0}^kD_{B_0}^m v_0
 &=\begin{cases}\chi_{\rm seed}^{(k)}\bar v,&m=0,\\0,&m\ge1,\end{cases}\\
 D_{t,0}^kD_{B_0}^m B_0&=0\quad(k+m\ge1),\\
 D_{t,0}^kD_{B_0}^m R_0
 &=\begin{cases}\partial_t^kR_0,&m=0,\\0,&m\ge1.\end{cases}
 \end{aligned}
\end{equation}
For the stress, the Lie derivatives also differentiate the tensor
factors. With primes on $A,\bar v,\bar B$ denoting $\partial_3$, the
first material and magnetic Lie derivatives are
\[
 \begin{aligned}
 \mathcal D_{t,0}R_0
 &=\chi_{\rm seed}''(A\otimes e_3+e_3\otimes A)
   -\chi_{\rm seed}\chi_{\rm seed}'
                      (\bar v'\otimes A+A\otimes\bar v'),\\
 \mathcal L_{B_0}R_0
 &=-\chi_{\rm seed}'(\bar B'\otimes A+A\otimes\bar B').
 \end{aligned}
\]
A tensor with two horizontal factors depending only on $(t,x_3)$ has
zero Lie derivative along either horizontal field, apart from the
ordinary time derivative in $\mathcal D_{t,0}$. Thus the horizontal
tensors produced above satisfy
\[
 \mathcal L_{B_0}^2R_0=0,\qquad
 \mathcal D_{t,0}^k\mathcal L_{B_0}R_0
 =-\chi_{\rm seed}^{(k+1)}(\bar B'\otimes A+A\otimes\bar B').
\]
More explicitly, put $S=A\otimes e_3+e_3\otimes A$ and
$H=\bar v'\otimes A+A\otimes\bar v'$. Then
$\mathcal D_{t,0}S=-\chi_{\rm seed}H$ and
$\mathcal D_{t,0}H=0$. For $k\ge1$,
\[
 \mathcal D_{t,0}^kR_0
 =\chi_{\rm seed}^{(k+1)}S
 -\sum_{j=0}^{k-1}\partial_t^{k-1-j}
      (\chi_{\rm seed}\chi_{\rm seed}^{(j+1)})H.
\]
This identity follows by induction. Together with
\eqref{iter:initializer-fields} and \eqref{iter:initializer-stress},
it gives the following bounds for all nonnegative integers $r,k$,
with constants depending on $r,k$ and the fixed cutoff $\chi_{\rm seed}$:
\begin{align*}
 \|\partial_t^kv_0\|_r&\le C_{r,k}k_{\mathrm{in}}^r,
 &\|\nabla B_0\|_r&\le C_r,\\
 \|\partial_t^kR_0\|_r&\le C_{r,k}(k_{\mathrm{in}}^{r-1}+c_{\mathrm{in}}),
 &\|\mathcal D_{t,0}^kR_0\|_r&\le C_{r,k}k_{\mathrm{in}}^r\quad(k\ge1),\\
 \|\mathcal D_{t,0}^k\mathcal L_{B_0}R_0\|_r
 &\le C_{r,k}k_{\mathrm{in}}^{r-1}.
\end{align*}
In these estimates, the product of the high frequency velocity gradient
with the low frequency part of $A$ is bounded using
$c_{\mathrm{in}}k_{\mathrm{in}}\le1$. This is why the coefficient of the
velocity component parallel to the magnetic field contains
$k_{\mathrm{in}}^{-2}$.

To verify the inductive stress bounds, it suffices to show that the
following three ratios tend to zero:
\[
 \begin{aligned}
 \frac{k_{\mathrm{in}}^{r-1}}
      {\lambda_0^{r-\alpha}\delta_1}
 &\lesssim\lambda_0^{-1+\beta+\alpha+2b\beta-\beta r},\\
 \frac{k_{\mathrm{in}}^r}
      {\lambda_0^{r-\alpha}\delta_1
                  (\lambda_0\delta_0^{1/2})^k}
 &\lesssim\lambda_0^{\alpha+2b\beta-\beta r-k(1-\beta)},
                         \qquad k\ge1,\\
 \frac{k_{\mathrm{in}}^{r-1}}
      {\lambda_0^{r+1-\alpha}\delta_{B,0}^{1/2}\delta_1}
 &\lesssim\lambda_0^{-2+2\beta+\frac{b-1}{b}(\gamma_a+\gamma_\parallel)
                                  +\alpha+2b\beta-\beta r}.
 \end{aligned}
\]
The first two tend to zero because
\[
 1-(2b+1)\beta-\alpha
 >\gamma_a-\alpha>0
\]
by \eqref{iter:principal-margins} and the chosen smallness of
$\gamma_S$. The third ratio also tends to zero: the inequality
$0<\frac{b-1}{b}(\gamma_a+\gamma_\parallel)<1-\beta$ makes its exponent smaller than that
of the first ratio.
Additional material derivatives only increase the denominators;
additional magnetic derivatives vanish. For the low frequency part of the stress, we have
\[
 \frac{c_{\mathrm{in}}}{\lambda_0^{-\alpha}\delta_1}
 =k_{\mathrm{in}}^{-2}\lambda_0^{\alpha}\delta_1=o(1).
\]
These comparisons prove \eqref{prep:old-stress} for both transport and
Lie derivatives, and \eqref{prep:old-stress-time} for ordinary time
derivatives, at $q=0$.

The spatial bounds for $B_0$ follow separately from
\[
 \|\nabla B_0\|_r\le C_r
       =o(\lambda_0^{r+1}\delta_{B,0}^{1/2}).
\]
All positive-order compositions of material and magnetic derivatives of $B_0$ vanish, as do its
Lie derivatives along the velocity flow and along $B_0$.
For $v_0$, divide the bounds for ordinary derivatives and for compositions of material and magnetic derivatives by their respective
inductive bounds. The ratios are bounded by
\[
 \begin{aligned}
 \frac{k_{\mathrm{in}}^r}{\lambda_0^{r+h}\delta_0^{1/2}}
 &\le\lambda_0^{\beta(1-r)-h},&&r+h\ge1,\\
 \frac{k_{\mathrm{in}}^r}
 {\lambda_0^r\delta_0^{1/2}(\lambda_0\delta_0^{1/2})^k}
 &\le\lambda_0^{\beta(1-r)-k(1-\beta)},&&r+k\ge1.
 \end{aligned}
\]
The first exponent is negative except at $(r,h)=(1,0)$, and the second
except at $(r,k)=(1,0)$. At these two orders, choose $C_0$ to absorb
the bounded ratios. With $C_0$ fixed, increasing $a$ gives the estimates
at every other derivative order in \eqref{prep:old-ordinary} and
\eqref{prep:old-mixed}. Since $p_0=0$, the pressure bounds hold as well.

For the pointwise estimates, $|B_0|=1$ and
$|v_0|\le1+c_{\mathrm{in}}$. Since $c_{\mathrm{in}}\to0$ and
$\delta_0^{1/2}\to0$ as $a\to\infty$, the two upper bounds are below
$2(1-\delta_0^{1/2})$ for large $a$, and the lower magnetic bound in
\eqref{full:pointwise-induction} holds as well.

Finally, for large $a$, the interval $[2/3,1]$ is compactly contained
in the stress-support interval \eqref{full:stress-support} at $q=0$.
The inclusion $\supp R_0\subset[2/3,1]$ gives the initial support
condition. The cutoff bounds for subsequent iterates
were established in Section~\ref{sec:preparation}.
\end{proof}

\subsection{Convergence and conserved quantities}

\begin{proof}[Proof of Proposition~\ref{iter:limit-proposition}]
Fix $0<\gamma<\beta$. Interpolation of
\eqref{iter:increment-zero} and \eqref{iter:increment-one} gives
\[
 \|v_{q+1}-v_q\|_{C^\gamma_{x,t}}
 +\|B_{q+1}-B_q\|_{C^\gamma_{x,t}}
 \le C_\gamma\lambda_{q+1}^{\gamma-\beta}.
\]
Since $\gamma<\beta$,
\[
 \sum_{q\ge0}\lambda_{q+1}^{\gamma-\beta}
   =\sum_{q\ge0}a^{-(\beta-\gamma)b^{q+1}}<\infty.
\]
Hence $(v_q,B_q)$ converges in $C^\gamma_{x,t}$. Interpolating
\eqref{iter:stress-limit}, we also obtain
\[
 \|R_q\|_{C_tC_x^\gamma}
 \le C_\gamma\lambda_q^{\gamma-\alpha-2b\beta}
 \longrightarrow0,
\]
since $\gamma<\beta<2b\beta+\alpha$.
Taking divergence of the momentum equation gives the normalized
pressure formula
\[
 p_q=\Delta^{-1}\ddiv\ddiv
       (R_q-v_q\otimes v_q+B_q\otimes B_q).
\]
The products on the right converge in $C_tC_x^\gamma$, since
$C^\gamma$ is an algebra. Boundedness of $\Delta^{-1}\ddiv\ddiv$ on
$C^\gamma$ for $\gamma>0$ therefore gives
\[
 p_q\longrightarrow
 p=\Delta^{-1}\ddiv\ddiv(-v\otimes v+B\otimes B)
 \quad\hbox{in }C_tC_x^\gamma.
\]
The case $\gamma=0$ follows from convergence at any positive
exponent below $\beta$.

Interpolating the magnetic bounds in \eqref{iter:increment-zero} and
\eqref{iter:increment-one} gives, for every
$0\le\gamma<\beta+\frac{b-1}{b}(\gamma_a+\gamma_\parallel)<1/2$,
\[
 \|B_{q+1}-B_q\|_{C^\gamma_{x,t}}
 \le C_\gamma\lambda_{q+1}^{\gamma}\delta_{B,q+1}^{1/2}
 =C_\gamma\lambda_{q+1}^{\gamma-\beta-\frac{b-1}{b}(\gamma_a+\gamma_\parallel)},
 \qquad
 \sum_{q\ge0}\lambda_{q+1}^{\gamma-\beta-\frac{b-1}{b}(\gamma_a+\gamma_\parallel)}<\infty.
\]
For $\gamma=0$, this is the uniform increment bound. Summing over $q$
shows that $B_q$ converges to $B$ in each of the stated space--time
H\"older norms.

In particular,
\[
 v_q\otimes v_q-B_q\otimes B_q-R_q\longrightarrow v\otimes v-B\otimes B,
 \qquad v_q\times B_q\longrightarrow v\times B
\]
uniformly. All terms in the relaxed equations therefore pass to the
limit against smooth compactly supported test functions, giving the
ideal MHD equations in distributions.

For the initial relaxed solution of Lemma~\ref{iter:initializer-lemma},
the magnetic field has unit length. Choose $a$ sufficiently large that
\[
 C\sum_{q\ge0}\delta_{B,q+1}^{1/2}<1/2.
\]
Then, uniformly in $q,t,x$,
\[
 |B_q|\ge |B_0|-\sum_{i<q}\|B_{i+1}-B_i\|_0>\tfrac12,
 \qquad |B|\ge\tfrac12.
\]

The perturbations vanish outside $[1/2,5/2]$, so the limiting fields
equal $(v_0,B_0)$ at $t=0$ and $t=3$. Evaluating the total energy and
cross-helicity by \eqref{eq:invariants} for the initial relaxed solution gives
\[
 \begin{aligned}
 \mathcal E_0(t)&=8\pi^3[1+(1+c_{\mathrm{in}}^2)\chi_{\rm seed}(t)^2],\\
 \mathcal H_0^\times(t)&=8\pi^3c_{\mathrm{in}}\chi_{\rm seed}(t).
 \end{aligned}
\]
Consequently
\[
 \mathcal E(0)-\mathcal E(3)=8\pi^3(1+c_{\mathrm{in}}^2),\qquad
 \mathcal H^\times(0)-\mathcal H^\times(3)
   =8\pi^3c_{\mathrm{in}}>0.
\]

Magnetic helicity is instead constant along the entire iteration.
First, the exact induction equation gives
\[
 \partial_t\int_{\mathbb T^3}B_q\,\dd x
   =\int_{\mathbb T^3}\curl(v_q\times B_q)\,\dd x=0.
\]
The field $B_q(0)=\bar B$ has zero mean, so every $B_q$ remains
mean zero. Its vector potential is
$A_q=\curl^{-1}B_q=-\Delta^{-1}\curl B_q$, as in
\eqref{eq:invariants}. On mean-zero divergence-free fields,
$\curl^2=-\Delta$ and hence $\curl A_q=B_q$. The operator
$\curl^{-1}$ is self-adjoint in $L^2$ by integration by parts; its
nonzero multipliers have size $|k|^{-1}\le1$, so it is also bounded on
$L^2$. Using $\partial_tB_q=\curl(v_q\times B_q)$ and integrating by
parts, we obtain
\[
 \frac{\dd}{\dd t}\int_{\mathbb T^3}A_q\cdot B_q\,\dd x
 =2\int_{\mathbb T^3}A_q\cdot\curl(v_q\times B_q)\,\dd x
 =2\int_{\mathbb T^3}B_q\cdot(v_q\times B_q)\,\dd x=0.
\]
At $t=0$, all magnetic iterates $B_q$ equal $\bar B$, and
$\curl\bar B=-\bar B$. Hence their common magnetic helicity is
\[
 \int_{\mathbb T^3}(\curl^{-1}\bar B)\cdot\bar B\,\dd x
 =-8\pi^3\ne0.
\]
To pass this invariant to the limit, use $B_q\to B$ in $L^2$,
uniformly in time, and the boundedness of $\curl^{-1}$ on $L^2$ in
this subspace:
\[
 \left|\int(\curl^{-1}B_q)\cdot B_q
              -\int(\curl^{-1}B)\cdot B\right|
 \le C\|B_q-B\|_2(\|B_q\|_2+\|B\|_2)\longrightarrow0.
\]
The limiting helicity is therefore the same nonzero constant at
every time, as claimed.
\end{proof}

\begin{proof}[Proof of Theorem~\ref{full:main}]
Given $0\le\gamma<1/3$, choose $\gamma<\beta<1/3$ and then the
parameters as in Section~\ref{sec:parameter-choice}.
Lemma~\ref{iter:initializer-lemma} supplies the initial iterate for
sufficiently large $a$. Proposition~\ref{full:step} then constructs the
sequence with the separate increment bounds
\eqref{full:increment-bound}, and Proposition~\ref{iter:limit-proposition}
gives the claimed $C^\gamma$ regularity and the remaining properties.
\end{proof}

\subsection{Anisotropic regularity}\label{ssec:fieldlinesproof}

\begin{proof}[Proof of Theorem~\ref{full:fieldline}]
Fix $\gamma_v,\gamma_B$ as in \eqref{full:fieldline-exponents}.
Choose $\gamma_v<\beta<1/3$ sufficiently close to $\gamma_v$, and then
choose the parameters so that
\[
 \max\{\beta,\gamma_B\}
 <\beta+\frac{b-1}{b}(\gamma_a+\gamma_\parallel)
 <\min\{\tfrac12,1-2\beta\}.
\]
These choices are possible by \eqref{iter:exponent-window} and its proof. They imply
\begin{equation}\label{iter:fieldline-choice}
 \begin{aligned}
 \frac{\beta}{1-\beta-\frac{b-1}{b}(\gamma_a+\gamma_\parallel)}
     &>\frac{\gamma_v}{1-\gamma_B},\\
 \frac{\beta+\frac{b-1}{b}(\gamma_a+\gamma_\parallel)}
      {1-\beta-\frac{b-1}{b}(\gamma_a+\gamma_\parallel)}
     &>\frac{\gamma_B}{1-\gamma_B}.
 \end{aligned}
\end{equation}
Choose the remaining parameters as in
Section~\ref{sec:parameter-choice}. The proof of
Theorem~\ref{full:main} with these parameters gives a weak solution
with all the stated properties. Proposition~\ref{iter:limit-proposition}
gives $v\in C^{\gamma_v}(\mathbb T^3\times\R)$ and
$B\in C^{\gamma_B}(\mathbb T^3\times\R)$.
It remains to estimate the fields along an integral curve of $B$.
We prove the stronger estimates with the exponents on the left of
\eqref{iter:fieldline-choice}. To treat both fields at once, set
\[
 \begin{aligned}
 f&=v,& f_q&=v_q,& \rho&=\beta,\\
 f&=B,& f_q&=B_q,& \rho&=\beta+\frac{b-1}{b}(\gamma_a+\gamma_\parallel).
 \end{aligned}
\] Since
$\beta+\frac{b-1}{b}(\gamma_a+\gamma_\parallel)<1/2$,
\[
 0<\rho\le\beta+\frac{b-1}{b}(\gamma_a+\gamma_\parallel)<\tfrac12,
 \qquad 1-\beta-\frac{b-1}{b}(\gamma_a+\gamma_\parallel)-\rho\ge1-2\beta-2\frac{b-1}{b}(\gamma_a+\gamma_\parallel)>0.
\]
All constants below are uniform in physical time.

The full increment bound~\eqref{full:increment-bound} gives
\[
 \|B-B_q\|_0\le C\sum_{j\ge q+1}\delta_{B,j}^{1/2}
       \le C\delta_{B,q+1}^{1/2}
       =C\lambda_{q+1}^{-\beta-\frac{b-1}{b}(\gamma_a+\gamma_\parallel)}.
\]
The sums are controlled by geometric series since
$\lambda_{j+1}/\lambda_j\ge a^{b-1}>1$. The ordinary derivative bounds in
\eqref{prep:old-ordinary} and the first magnetic derivative bound
in \eqref{prep:old-mixed} give the following estimate along the limiting
magnetic field:
\begin{equation}\label{iter:limiting-magnetic-derivative}
 \begin{split}
 \|B\cn f_q\|_0
 &\le\|B_q\cn f_q\|_0
                +C\|B-B_q\|_0[f_q]_1\\
 &\le C\lambda_q^{1-\beta-\frac{b-1}{b}(\gamma_a+\gamma_\parallel)-\rho}
       +C\lambda_{q+1}^{-\beta-\frac{b-1}{b}(\gamma_a+\gamma_\parallel)}
                   \lambda_q^{1-\rho}\\
 &\le C\lambda_q^{1-\beta-\frac{b-1}{b}(\gamma_a+\gamma_\parallel)-\rho}.
 \end{split}
\end{equation}
For $q\ge1$, set $d_q=f_q-f_{q-1}$. The increment bound and
\eqref{iter:limiting-magnetic-derivative}, applied at $q$ and
$q-1$, yield
\[
 \|d_q\|_0\le C\lambda_q^{-\rho},\qquad
 \|B\cn d_q\|_0
       \le C\lambda_q^{1-\beta-\frac{b-1}{b}(\gamma_a+\gamma_\parallel)-\rho}.
\]
Here the exponent in the second bound is positive, so its value at
$q-1$ is bounded by its value at $q$.

Fix $t$ and an integral curve $\Gamma'=B(t,\Gamma)$. The curve is
$C^1$ and each $d_q$ is smooth; thus the chain rule applies even though
$B$ need not be differentiable. For $h=|s-r|\in(0,1]$, the two bounds
above give
\[
 |d_q(t,\Gamma(s))-d_q(t,\Gamma(r))|
 \le C\min\{\lambda_q^{-\rho},
                   h\lambda_q^{1-\beta-\frac{b-1}{b}(\gamma_a+\gamma_\parallel)-\rho}\}.
\]
Let $Q$ be the largest integer $q\ge1$ for which
$\lambda_q^{1-\beta-\frac{b-1}{b}(\gamma_a+\gamma_\parallel)}\le h^{-1}$, and set $Q=0$ if there is
none. Since $f_0$ has a uniformly bounded spatial gradient and $B$
is bounded, its contribution is at most $Ch$. Summing the uniformly convergent
series $f=f_0+\sum_{q\ge1}d_q$ gives
\begin{equation}\label{iter:fieldline-summation}
 \begin{split}
 |f(t,\Gamma(s))-f(t,\Gamma(r))|
 &\le Ch+C h\sum_{1\le q\le Q}
           \lambda_q^{1-\beta-\frac{b-1}{b}(\gamma_a+\gamma_\parallel)-\rho}
          +C\sum_{q>Q}\lambda_q^{-\rho}\\
 &\le C h^{\rho/(1-\beta-\frac{b-1}{b}(\gamma_a+\gamma_\parallel))}.
 \end{split}
\end{equation}
Indeed, the first sum is controlled by its last term and
$\lambda_Q^{1-\beta-\frac{b-1}{b}(\gamma_a+\gamma_\parallel)}\le h^{-1}$; it is empty
when $Q=0$. The second is controlled by its first term and
$\lambda_{Q+1}^{1-\beta-\frac{b-1}{b}(\gamma_a+\gamma_\parallel)}>h^{-1}$.
The contribution of $f_0$ is also bounded by the displayed power since its
exponent lies in $(0,1)$.
This proves the estimates with the exponents in
\eqref{iter:fieldline-choice}, hence \eqref{full:fieldline-holder}.

By continuity of $B(t,\cdot)$, the existence theorem for ordinary
differential equations gives an integral curve through every point.
Boundedness of $B$ on the torus extends each curve to all parameter values.
No uniqueness is needed, since the argument above applies to every such
curve. Passing to the limit in \eqref{full:pointwise-induction} gives
$1/4\le|B|\le2$, and hence the change to arclength is uniformly bi-Lipschitz. For arclength
differences at most $1/4$, the estimate follows from
\eqref{iter:fieldline-summation} under this change of parameter. For
larger differences at most one, boundedness of $v$ and $B$ gives the
same assertion with a uniform constant.
\end{proof}

\appendix
\section{Lie Algebra Driven Perturbations}
\label{sec:path}

The induction equation identifies the magnetic field with the
pushforward of its initial value by the velocity flow. A
volume-preserving deformation of this flow therefore determines a
simultaneous perturbation of the velocity and magnetic field that
preserves induction exactly. This appendix develops the resulting
path identities and estimates. After deriving the identities for an
arbitrary Lagrangian displacement field, we compose the principal and
corrector flows and compare their action on different backgrounds.
The principal flow can have a large spatial differential even when
its phase displacement is small.
We use the phase displacement and invariant frame structure to estimate the field
increments, and then prove the properties of the derivative classes
needed for these estimates. We follow the tensor and derivative
conventions of Section~\ref{ssec:path-conventions}.

For a volume-preserving diffeomorphism $Y$, tensor pushforward
acts on a scalar $f$, vector $F$, one-form $\Theta$ and
contravariant two-tensor $T$ by
\[
 \begin{aligned}
 Y_*f&=f\circ Y^{-1},&
 Y_*F&=(\DD Y\,F)\circ Y^{-1},\\
 Y_*\Theta&=((\DD Y)^{-T}\Theta)\circ Y^{-1},&
 Y_*T&=(\DD Y\,T\,(\DD Y)^T)\circ Y^{-1}.
 \end{aligned}
\]
We write $Y^*=(Y^{-1})_*$ for the corresponding pullback. Since
$Y$ preserves $\mathrm{vol}$, pushforward of a two-form corresponds
to pushforward of a vector under the identification
$F\mapsto\iota_F\mathrm{vol}$.

Let $(v,B,p,R)$ be a smooth relaxed MHD solution of \eqref{eq:mhd},
and let $\Phi_{t\leftarrow t_0}$ be the volume-preserving flow of $v$.
The identity
\[
 \partial_t(\Phi_{t\leftarrow t_0}^*B(t))
 =\Phi_{t\leftarrow t_0}^*(\partial_tB+[v,B])
\]
shows that induction is equivalent to the pushforward representation
\[
 B(t)=(\Phi_{t\leftarrow t_0})_*B(t_0).
\]
Under scalar composition with this flow, $D_t$ becomes
$\partial_t|_y$ and $D_B$ becomes $B(t_0,y)\cn_y$; thus the two
operators $\mathcal A^\pm=D_t\pm D_B$ become
$\partial_t|_y\pm B(t_0,y)\cn_y$. Cartesian tensor components transform by this scalar composition,
whereas tensor pullback also includes the differential matrices above. The same pushforward identity conjugates
the magnetic flow:
\[
 e^{aB(t)}=\Phi_{t\leftarrow t_0}\circ e^{aB(t_0)}
                            \circ\Phi_{t_0\leftarrow t}.
\]
Thus integrals along the magnetic flow may be evaluated at the reference time
and then pushed forward by the velocity flow. For ordinary time derivatives,
we use $\partial_t=D_t-v\cn$.

We replace the flow map by
$X_s(t)\circ\Phi_{t\leftarrow t_0}$, where $X_0=\IId$ and every
$X_s(t)$ preserves volume. Differentiation in physical time and
transport of the magnetic field give
\begin{equation}\label{eq:path-fields}
 v_s=\partial_tX_s\circ X_s^{-1}+(X_s)_*v,\qquad
 B_s=(X_s)_*B.
\end{equation}
The term $\partial_tX_s\circ X_s^{-1}$ is the velocity contributed
by the time dependence of the deformation. The physical flow normalized
at $t_0$ is
\[
 \Phi^{(s)}_{t\leftarrow t_0}
 =X_s(t)\circ\Phi_{t\leftarrow t_0}\circ X_s(t_0)^{-1},
 \qquad
 B_s(t)=(\Phi^{(s)}_{t\leftarrow t_0})_*B_s(t_0).
\]
Thus induction holds along the whole deformation; if $X_s(t_0)=\IId$,
the reference magnetic data are preserved as well.

We construct the deformation from its
\emph{Lagrangian displacement field} (LDF),
\[
 \xi_s=\partial_sX_s\circ X_s^{-1},\qquad
 \partial_sX_s(t,x)=\xi_s(t,X_s(t,x)),\qquad X_0(t,x)=x.
\]
Here physical time $t$ is a parameter in the ordinary differential
equation, and $0\le s\le1$ parametrizes the deformation. We allow the LDF to depend on $s$
because the principal part will be pushed forward by the corrector flow.
For a smooth divergence-free $\xi_s$ on the torus, the
flow exists on the whole parameter interval, and
\[
 \partial_s\det\DD X_s
 =((\ddiv\xi_s)\circ X_s)\det\DD X_s=0.
\]
Conversely, this identity shows that the LDF of a smooth
volume-preserving path is divergence free. For $0\le r,s\le1$,
define the two-parameter propagator
\begin{equation}
 \mathcal U_{s\leftarrow r}:=(X_s\circ X_r^{-1})_*.
\end{equation}
The propagator acts by pushforward according to the tensor type of its argument.
We write $\mathcal D_{t,s}=\partial_t+\mathcal L_{v_s}$ and
$\mathcal L_{B_s}$ for the material and magnetic Lie derivatives on the path.

\subsection{Abstract Lagrangian displacement field}
\label{ssec:path-calculus}

We first derive the transport and variation identities for an
arbitrary $s$-dependent LDF $\xi_s$. The initial fields $v,B$ are fixed
in $s$. Only the finite Lie--Taylor formula at the end of the subsection
requires an autonomous LDF.

We shall call
\begin{equation}\label{eq:average}
 \mathcal P_X[F_\bullet]
       :=\int_0^1\mathcal U_{1\leftarrow r}F_r\,\dd r
\end{equation}
the \emph{path operator}; the dot denotes a family $r\mapsto F_r$.

\begin{lemma}[Path transport]
\label{lem:path-transport}
On each of these tensor types,
\begin{equation}\label{eq:generator}
 \begin{gathered}
 \mathcal U_{s\leftarrow r}\mathcal U_{r\leftarrow a}=\mathcal U_{s\leftarrow a},\qquad
 \mathcal U_{r\leftarrow r}=\IId,\qquad
 \mathcal U_{s\leftarrow r}^{-1}=\mathcal U_{r\leftarrow s},\\
 \partial_s\mathcal U_{s\leftarrow r}=-\mathcal L_{\xi_s}\mathcal U_{s\leftarrow r},
 \qquad
 \partial_r\mathcal U_{s\leftarrow r}=\mathcal U_{s\leftarrow r}\mathcal L_{\xi_r}.
 \end{gathered}
\end{equation}
Consequently, the solution of
$(\partial_s+\mathcal L_{\xi_s})F_s=G_s$ with value $F_r$ at $r$ is
\begin{equation}\label{eq:duhamel}
 F_s=\mathcal U_{s\leftarrow r}F_r+
            \int_r^s\mathcal U_{s\leftarrow a}G_a\,\dd a.
\end{equation}
Moreover, the propagator preserves tensor products and contractions
and commutes with exterior differentiation. In particular
\begin{equation}\label{eq:curl}
 \mathcal U_{s\leftarrow r}\curl\Theta
      =\curl(\mathcal U_{s\leftarrow r}\Theta).
\end{equation}

The propagator along the path intertwines the material and magnetic
Lie derivatives separately:
\begin{equation}
 \mathcal D_{t,s}\mathcal U_{s\leftarrow r}F
       =\mathcal U_{s\leftarrow r}\mathcal D_{t,r}F,
 \qquad
 \mathcal L_{B_s}\mathcal U_{s\leftarrow r}F
       =\mathcal U_{s\leftarrow r}\mathcal L_{B_r}F.
\end{equation}
Consequently, for all nonnegative integers $k,m$,
\begin{equation}\label{path:mixed-intertwining}
 \mathcal D_{t,s}^{k}\mathcal L_{B_s}^{m}\mathcal U_{s\leftarrow r}F
 =\mathcal U_{s\leftarrow r}\mathcal D_{t,r}^{k}\mathcal L_{B_r}^{m}F.
\end{equation}

Differentiation with respect to a fixed background gives a commutator
term. More precisely, let $\mathscr A=\partial_t+\mathcal L_z$, where
$z(t,x)$ is independent of the deformation parameter. For $r\le s$,
\begin{equation}\label{eq:transport}
 \mathscr A\mathcal U_{s\leftarrow r}F
 =\mathcal U_{s\leftarrow r}\mathscr AF
   -\int_r^s\mathcal U_{s\leftarrow a}
         \mathcal L_{\mathscr A\xi_a}\mathcal U_{a\leftarrow r}F\,\dd a.
\end{equation}
In particular,
\begin{equation}\label{eq:transport-average}
 \begin{aligned}
 \mathscr A\mathcal P_X[F_\bullet]
 ={}&\mathcal P_X[\mathscr AF_\bullet]\\
 &-\int_0^1\int_r^1\mathcal U_{1\leftarrow a}
       \mathcal L_{\mathscr A\xi_a}\mathcal U_{a\leftarrow r}F_r
                          \,\dd a\dd r.
 \end{aligned}
\end{equation}
The same identities hold with $\mathscr A=\mathcal L_z$.
\end{lemma}
\begin{proof}
For every smooth tensor family $F_s$, differentiation of its tensor
pullback gives
\[
 \partial_s(X_s^*F_s)=X_s^*(\partial_sF_s+\mathcal L_{\xi_s}F_s).
\]
For vectors and one-forms this follows by differentiating the factors
involving the differential matrix and its inverse or transpose in their
transformation laws; tensor products give the remaining cases. Thus
\[
 (\partial_s+\mathcal L_{\xi_s})\mathcal U_{s\leftarrow r}=0,
 \qquad
 \mathcal U_{s\leftarrow r}\mathcal U_{r\leftarrow a}=(X_s\circ X_a^{-1})_*
                              =\mathcal U_{s\leftarrow a}.
\]
The composition law gives the first line of \eqref{eq:generator}.
Differentiating it in $r$ gives
\[
 \partial_r\mathcal U_{s\leftarrow r}=\mathcal U_{s\leftarrow r}\mathcal L_{\xi_r}.
\]
For the forced equation, the pullback identity can be integrated in the
deformation parameter:
\[
 X_s^*F_s-X_r^*F_r=\int_r^sX_a^*G_a\,\dd a.
\]
Applying $(X_s)_*$ proves \eqref{eq:duhamel} and uniqueness.

Commutation of exterior differentiation with pushforward and preservation
of tensor contractions, together with $Y_*\mathrm{vol}=\mathrm{vol}$, give
\[
 d(Y_*\Theta)=Y_*d\Theta
 =Y_*(\iota_{\curl\Theta}\mathrm{vol})
 =\iota_{Y_*\curl\Theta}\mathrm{vol}.
\]
This is \eqref{eq:curl}. Differentiating in a volume-preserving flow
also gives $\curl\mathcal L_z=\mathcal L_z\curl$ when $\ddiv z=0$;
the fixed volume form likewise gives $[\partial_t,\curl]=0$.

To prove intertwining, let $Y=X_s\circ X_r^{-1}$ and
$u^Y=\partial_tY\circ Y^{-1}$. Then
\[
 \partial_t(Y_*F)=Y_*\partial_tF-\mathcal L_{u^Y}(Y_*F),
 \qquad v_s=u^Y+Y_*v_r,\qquad B_s=Y_*B_r.
\]
Substituting the velocity transformation law cancels the terms
containing $u^Y$. Equivariance of Lie differentiation under pushforward then gives
\[
 \begin{aligned}
 (\partial_t+\mathcal L_{v_s})(Y_*F)
 &=Y_*\partial_tF+\mathcal L_{Y_*v_r}(Y_*F)\\
 &=Y_*\mathcal D_{t,r}F,\\
 \mathcal L_{B_s}(Y_*F)&=Y_*\mathcal L_{B_r}F.
 \end{aligned}
\]
Iterating these two identities in the displayed order proves
\eqref{path:mixed-intertwining}.

For derivatives relative to a fixed background, the Jacobi identity gives
\[
 [\mathscr A,\mathcal L_{\xi_s}]
 =\mathcal L_{\partial_t\xi_s+[z,\xi_s]}
 =\mathcal L_{\mathscr A\xi_s}.
\]
Applying $\mathscr A$ to the homogeneous Lie transport equation for
$\mathcal U_{s\leftarrow r}F$ and using the independence of
$\mathscr A$ from $s$, we find that
$H_s=\mathscr A\mathcal U_{s\leftarrow r}F$ solves
\[
 (\partial_s+\mathcal L_{\xi_s})H_s
 =-[\mathscr A,\mathcal L_{\xi_s}]\mathcal U_{s\leftarrow r}F
 =-\mathcal L_{\mathscr A\xi_s}\mathcal U_{s\leftarrow r}F,
 \qquad H_r=\mathscr AF.
\]
Applying \eqref{eq:duhamel} gives \eqref{eq:transport}; integrating
that identity at $s=1$ in $r$ gives \eqref{eq:transport-average}.
The calculation with $\mathscr A=\mathcal L_z$ is identical after
omitting $\partial_t$.
\end{proof}

For $\mathscr A=\mathcal D_t$ or $\mathscr A=\mathcal L_B$,
formula~\eqref{eq:transport} expresses the commutator in terms of
$\mathcal D_t\xi_a$ or $\mathcal L_B\xi_a$, respectively. These terms
arise because the background is fixed in $s$, whereas
\eqref{path:mixed-intertwining} uses the fields on the path.
The two formulas can be iterated separately for material and magnetic
derivatives.

\paragraph{\textbf{Field increments.}}

We now express the field increments in terms of the LDF using
\eqref{eq:path-fields} and Lemma~\ref{lem:path-transport}.
First, writing $u_s^X=\partial_tX_s\circ X_s^{-1}$, volume preservation gives
\[
 \ddiv u_s^X=\partial_t\log\det\DD X_s\circ X_s^{-1}=0,
 \qquad
 \ddiv((X_s)_*F)=(\ddiv F)\circ X_s^{-1},
\]
so $v_s,B_s$ are divergence free. The spacetime map
$\widetilde X_s(t,x)=(t,X_s(t,x))$ satisfies
\[
 \begin{gathered}
 (\widetilde X_s)_*(\partial_t+v)=\partial_t+v_s,
 \qquad (\widetilde X_s)_*B=B_s,\\
 [\partial_t+v_s,B_s]
       =(\widetilde X_s)_*[\partial_t+v,B]=0,
 \end{gathered}
\]
which preserves induction along the path.
Commutation of differentiation in $s$ and $t$, together with \eqref{eq:generator}, implies
\[
 \partial_su_s^X=\partial_t\xi_s+[u_s^X,\xi_s],
 \qquad \partial_s((X_s)_*F)=[(X_s)_*F,\xi_s].
\]
Substitution in \eqref{eq:path-fields} gives the field variations
\begin{equation}\label{eq:path-linear-fields}
 \partial_sv_s=(\partial_t+\mathcal L_{v_s})\xi_s,
 \qquad \partial_sB_s=\mathcal L_{B_s}\xi_s.
\end{equation}
Put $w_s=v_s-v$ and $b_s=B_s-B$, and write $w=w_1$, $b=b_1$.
Subtracting the initial fields and moving the bracket terms to the left gives:
\[
 \begin{gathered}
 (\partial_s+\mathcal L_{\xi_s})w_s=\mathcal D_t\xi_s,
 \qquad
 (\partial_s+\mathcal L_{\xi_s})b_s=\mathcal L_B\xi_s,\\
 w_0=b_0=0.
 \end{gathered}
\]
The operators on the right are defined by the initial fields $v,B$, fixed in $s$.
Applying \eqref{eq:duhamel} gives
\begin{equation}\label{eq:path-velocity-duhamel}
 \begin{aligned}
 w_s&=\int_0^s\mathcal U_{s\leftarrow r}
             \bigl(\partial_t\xi_r+[v,\xi_r]\bigr)\,\dd r,\\
 b_s&=\int_0^s\mathcal U_{s\leftarrow r}[B,\xi_r]\,\dd r.
 \end{aligned}
\end{equation}
At the endpoint, these identities and their sum and difference read
\begin{equation}\label{eq:signed-increments}
 \begin{aligned}
 w&=\mathcal P_X[\mathcal D_t\xi_\bullet],&
 b&=\mathcal P_X[\mathcal L_B\xi_\bullet],\\
 w\pm b&=\mathcal P_X[(\mathcal D_t\pm\mathcal L_B)\xi_\bullet].
 \end{aligned}
\end{equation}
If the LDF has a one-form potential $\xi_s=\curl\Theta_s$,
curl commutes with both Lie derivatives and with the propagator by
\eqref{eq:curl}. Thus \eqref{eq:signed-increments} gives
\begin{equation}\label{eq:potential-increments}
 \begin{aligned}
 w&=\curl\mathcal P_X[\mathcal D_t\Theta_\bullet],&
 b&=\curl\mathcal P_X[\mathcal L_B\Theta_\bullet],\\
 w\pm b&=\curl\mathcal P_X[(\mathcal D_t\pm\mathcal L_B)\Theta_\bullet].
 \end{aligned}
\end{equation}
The representation~\eqref{eq:potential-increments} requires a global
potential for the chosen LDF; a general divergence-free vector field
on the torus need not admit one. If
$\xi_s(t_*,\cdot)=0$ for every $s$, then $X_s(t_*,\cdot)=\IId$
and the magnetic field at $t_*$ is unchanged. If the LDF
vanishes on a time interval, both increments vanish there.

\paragraph{\textbf{Momentum variation and force.}}

Write the momentum expression along the path without its pressure as
\[
 \begin{gathered}
 \mathcal M_s=\partial_tv_s+
               \ddiv(v_s\otimes v_s-B_s\otimes B_s),\\
 D_{t,s}=\partial_t+v_s\cn,\qquad D_{B_s}=B_s\cn.
 \end{gathered}
\]
For divergence-free increments $(w,b)$, the linear part of the
momentum change at $(v,B)$ is
\begin{equation}\label{eq:linear-force}
 \mathcal F_{v,B}(w,b)
 =D_tw-D_Bb+(\DD v)w-(\DD B)b.
\end{equation}
Differentiating $\mathcal M_s$ with respect to $s$ gives
\[
 \partial_s\mathcal M_s
 =\mathcal F_{v_s,B_s}(\partial_sv_s,\partial_sB_s).
\]
Expanding $\mathcal M_1-\mathcal M_0$ in the endpoint increments
$w,b$ from \eqref{eq:signed-increments} gives the exact identities
\begin{equation}\label{eq:endpoint-polynomial}
 \begin{aligned}
 \mathcal M_1-\mathcal M_0
   &=\mathcal F_{v,B}(w,b)
                  +\ddiv(w\otimes w-b\otimes b),\\
 w\otimes w-b\otimes b
   &=(w+b)\odot(w-b).
 \end{aligned}
\end{equation}
Writing this change as $\ddiv S-\nabla\pi$ gives the new stress
$R+S$ and pressure $p+\pi$.

To express the linear force through the displacement, take a
divergence-free vector field $\xi$ and its tangent fields
\begin{equation}\label{eq:linear-fields}
 w=(\partial_t+\mathcal L_v)\xi
      =D_t\xi-(\DD v)\xi,\qquad
 b=\mathcal L_B\xi=D_B\xi-(\DD B)\xi.
\end{equation}
The bracket preserves divergence-free vector fields, so $w,b$ are
divergence free. Induction gives
\[
 [\partial_t+\mathcal L_v,\mathcal L_B]=0,
 \qquad (\partial_t+\mathcal L_v)b=\mathcal L_Bw=-[w,B],
\]
which is the linearized induction equation.
Substituting \eqref{eq:linear-fields} in \eqref{eq:linear-force},
the terms containing $D_t\xi$ or $D_B\xi$ cancel, leaving
\begin{equation}\label{eq:displacement-force}
 \begin{aligned}
 \mathcal F_{v,B}(w,b)
 &=(D_t^2-D_B^2)\xi
   +[-D_t\DD v-(\DD v)^2+D_B\DD B+(\DD B)^2]\xi\\
 &=(D_t^2-D_B^2)\xi-\DD(D_tv-D_BB)\xi.
 \end{aligned}
\end{equation}
The second equality follows from the component product identities
\[
 D_t\DD v=\DD(D_tv)-(\DD v)^2,\qquad
 D_B\DD B=\DD(D_BB)-(\DD B)^2.
\]
The commuting scalar transports also give
$D_t^2-D_B^2=\mathcal A^-\mathcal A^+$.
By \eqref{eq:path-linear-fields}, the path variations are precisely
these tangent fields on the current background, so
\[
 \partial_s\mathcal M_s
 =(D_{t,s}^2-D_{B_s}^2)\xi_s-(\DD\mathcal M_s)\xi_s.
\]
Here $D_{t,s}v_s-D_{B_s}B_s=\mathcal M_s$. Integrating gives
\begin{equation}\label{eq:path-force-integral}
 \mathcal M_1-\mathcal M_0
 =\int_0^1\left[
     (D_{t,s}^2-D_{B_s}^2)\xi_s-(\DD\mathcal M_s)\xi_s
                 \right]\dd s.
\end{equation}

The quadratic term in \eqref{eq:endpoint-polynomial} is, by
\eqref{eq:signed-increments}, a product of path integrals. The next
identity separates this product into a transported tensor and a
covariance. Its difference formula will give a small factor from each
change along the path.

\begin{lemma}[Path products]
For smooth vector families $F_r,G_r$,
\begin{equation}\label{eq:product}
 (\mathcal P_X[F_\bullet])\odot
       (\mathcal P_X[G_\bullet])
 =\mathcal P_X[F_\bullet\odot G_\bullet]
                   -\Cov_X(F_\bullet,G_\bullet),
\end{equation}
where
\begin{equation}\label{eq:covariance}
 \begin{aligned}
 \Cov_X(F_\bullet,G_\bullet)
 :=\frac12\int_0^1\!\!\int_0^1
 &\bigl(\mathcal U_{1\leftarrow s}F_s-\mathcal U_{1\leftarrow r}F_r\bigr)\odot\\[-2pt]
 &\bigl(\mathcal U_{1\leftarrow s}G_s-\mathcal U_{1\leftarrow r}G_r\bigr)
                      \,\dd s\dd r.
 \end{aligned}
\end{equation}
The difference in either factor has the exact representation
\begin{equation}\label{eq:covariance-difference}
 \mathcal U_{1\leftarrow s}F_s-\mathcal U_{1\leftarrow r}F_r
 =\int_r^s\mathcal U_{1\leftarrow a}
        (\partial_aF_a+\mathcal L_{\xi_a}F_a)\,\dd a.
\end{equation}
\end{lemma}
\begin{proof}
Set $f_s=\mathcal U_{1\leftarrow s}F_s$ and $g_s=\mathcal U_{1\leftarrow s}G_s$.
Preservation of tensor products by pushforward gives
\[
 f_s\odot g_s=\mathcal U_{1\leftarrow s}(F_s\odot G_s).
\]
Expanding the double integral, interchanging $r$ and $s$, and applying
Fubini's theorem, we obtain
\begin{align*}
 \frac12\int_0^1\!\!\int_0^1
 (f_s-f_r)\odot(g_s-g_r)\,\dd s\dd r
 &=\int_0^1 f_s\odot g_s\,\dd s
   -\int_0^1\!\!\int_0^1f_s\odot g_r\,\dd s\dd r\\
 &=\mathcal P_X[F_\bullet\odot G_\bullet]
   -\mathcal P_X[F_\bullet]\odot\mathcal P_X[G_\bullet].
\end{align*}
This proves \eqref{eq:product}. Moreover, \eqref{eq:generator} gives
\[
 \partial_a(\mathcal U_{1\leftarrow a}F_a)
  =\mathcal U_{1\leftarrow a}(\partial_aF_a+\mathcal L_{\xi_a}F_a).
\]
Integrating from $r$ to $s$, with the oriented integral if $s<r$,
proves \eqref{eq:covariance-difference}.
\end{proof}

For a family depending on the deformation parameter,
\eqref{eq:covariance-difference} shows that the covariance estimate
must bound both $\partial_aF_a$ and $\mathcal L_{\xi_a}F_a$ with
the required small factor.

\paragraph{\textbf{Finite Lie--Taylor formula.}}

Suppose that $\xi_s=\xi$ is independent of the deformation
parameter, while still allowed to depend on physical time. At each
fixed $t$, $X_s$ is the flow of $\xi(t)$. For a tensor $F$ independent
of the deformation parameter, the substitution $s=1-r$ in
\eqref{eq:average} gives
\[
 \mathcal P_X[F]=\int_0^1\mathcal U_{s\leftarrow 0}F\,\dd s.
\]
Taylor's theorem in $s$ gives a finite expansion with an integral
remainder, without any assumption on convergence of the Lie series.

\begin{proposition}[Finite Lie--Taylor formula]
Let $\xi=\xi(t,x)$ be a smooth divergence-free vector field
independent of the deformation parameter, and let $X_s$ be its
flow as above. For every integer $K\ge0$ and every smooth tensor
$F$ independent of that parameter, the path operator
$\mathcal P_X$ satisfies
\begin{equation}\label{eq:taylor}
 \mathcal P_X[F]
 =\sum_{j=0}^{K}\frac{(-1)^j}{(j+1)!}\mathcal L_\xi^jF
 +\frac{(-1)^{K+1}}{(K+1)!}
    \int_0^1(1-s)^{K+1}\mathcal U_{s\leftarrow 0}
                              \mathcal L_\xi^{K+1}F\,\dd s.
\end{equation}
If $F$ and $G$ are smooth vector fields independent of the
deformation parameter and $T=F\odot G$, then
\begin{equation}\label{eq:quadratic-remainder}
 \begin{aligned}
 (\mathcal P_X[F])\odot(\mathcal P_X[G])
 ={}&T-\tfrac12\mathcal L_\xi T-\Cov_X(F,G)\\
 &+\tfrac12\int_0^1(1-s)^2\mathcal U_{s\leftarrow 0}
                                   \mathcal L_\xi^2T\,\dd s.
 \end{aligned}
\end{equation}
\end{proposition}
\begin{proof}
For an autonomous LDF, $(X_s)_*\xi=\xi$, so
\[
 \partial_s^j\mathcal U_{s\leftarrow 0}F
      =(-1)^j\mathcal U_{s\leftarrow 0}\mathcal L_\xi^jF.
\]
Taylor's formula with integral remainder gives
\[
 \mathcal U_{s\leftarrow 0}F
 =\sum_{j=0}^K\frac{(-s)^j}{j!}\mathcal L_\xi^jF
 +\frac{(-1)^{K+1}}{K!}\int_0^s(s-a)^K
                   \mathcal U_{a\leftarrow 0}\mathcal L_\xi^{K+1}F\,\dd a.
\]
Integrate in $s\in[0,1]$ and use
\[
 \int_0^1\frac{s^j}{j!}\,\dd s=\frac1{(j+1)!},
 \qquad \int_a^1\frac{(s-a)^K}{K!}\,\dd s
                      =\frac{(1-a)^{K+1}}{(K+1)!}.
\]
Interchanging the remainder integrals and using these identities proves
\eqref{eq:taylor}. For $T=F\odot G$, the expansion
with $K=1$ in the first term on the right of \eqref{eq:product}
gives \eqref{eq:quadratic-remainder}.
\end{proof}

\begin{remark}[Geometric motivation]
Formula~\eqref{eq:taylor} expresses the linear displacement and its
successive Lie corrections as terms of the same exact perturbation.
The quadratic identity~\eqref{eq:quadratic-remainder} separates the
first variation of $F\odot G$ from the quadratic Lie remainder and
the covariance. The invariant potential in
Subsection~\ref{ssec:exact-phase-geometry} will reduce the Lie
derivatives in these terms to scalar differentiation. After composition
with the corrector, the LDF depends on the deformation parameter, and
the exact identities~\eqref{eq:signed-increments} and
\eqref{eq:path-force-integral} continue to apply.
\end{remark}

\subsection{Composition of the principal and corrector flows and background comparison}
\label{ssec:path-composition}

The principal potential is oscillatory in the Lagrangian chart of the local
background. When we apply the corrector flow, we push forward both its
phase and its one-form factor. The principal LDF therefore depends on
$s$, but pulling it back by the corrector recovers the fixed field
$\xi_0^{\rm p}$. The following factorization allows us to estimate
its flow in the original chart.

\paragraph{\textbf{Composition.}}
Let $X_s^{\rm c},X_s^{\rm p}$ be the flows of smooth divergence-free
fields $\xi^{\rm c},\xi_0^{\rm p}$ independent of $s$. Write
\[
 \mathcal U_s^{\rm c}=(X_s^{\rm c})_*,\qquad
 \mathcal U_s^{\rm p}=(X_s^{\rm p})_*,
\]
and put
\[
 \nabla_s^{\rm c}=\partial_s+\mathcal L_{\xi^{\rm c}},\qquad
 \xi_s^{\rm p}=\mathcal U_s^{\rm c}\xi_0^{\rm p}.
\]

\begin{lemma}[Composition]\label{direct:flatness}
The path $X_s=X_s^{\rm c}\circ X_s^{\rm p}$ has LDF
$\xi_s=\xi^{\rm c}+\mathcal U_s^{\rm c}\xi_0^{\rm p}$ and propagator
\begin{equation}\label{path:factorized-transport}
 \mathcal U_{s\leftarrow r}
 =\mathcal U_s^{\rm c}\mathcal U_{s-r}^{\rm p}
                         (\mathcal U_r^{\rm c})^{-1},\qquad
 \mathcal U_{s\leftarrow 0}^{-1}\mathcal U_s^{\rm c}
                         =(\mathcal U_s^{\rm p})^{-1}.
\end{equation}
If $\xi_0^{\rm p}=\curl\Theta_0^{\rm p}$ and
$\Theta_s^{\rm p}=\mathcal U_s^{\rm c}\Theta_0^{\rm p}$, then
\begin{equation}\label{path:principal-transport}
 \nabla_s^{\rm c}\Theta_s^{\rm p}
 =\nabla_s^{\rm c}\xi_s^{\rm p}=0,\qquad
 [\nabla_s^{\rm c},\mathcal L_{\xi_s^{\rm p}}]=0.
\end{equation}
For any background satisfying induction, let $(v_s^{\rm c},B_s^{\rm c})$
be its corrector path fields and set
$\mathcal D_{t,s}^{\rm c}=\partial_t+\mathcal L_{v_s^{\rm c}}$. Then
\begin{equation}\label{direct:differentiated-transport}
 (\mathcal D_{t,s}^{\rm c})^k\mathcal L_{B_s^{\rm c}}^m
       \mathcal U_s^{\rm c}F
 =\mathcal U_s^{\rm c}\mathcal D_t^k\mathcal L_B^mF.
\end{equation}
The three operators $\nabla_s^{\rm c}$, $\mathcal D_{t,s}^{\rm c}$ and
$\mathcal L_{B_s^{\rm c}}$ commute pairwise.
If $(v_s,B_s)$ are the full path fields of the same initial pair, then
\begin{equation}\label{path:factorized-field-increments}
 \begin{aligned}
 v_s-v_s^{\rm c}
   &=\mathcal U_s^{\rm c}\int_0^s
             \mathcal U_r^{\rm p}\mathcal D_t\xi_0^{\rm p}\,\dd r,\\
 B_s-B_s^{\rm c}
   &=\mathcal U_s^{\rm c}\int_0^s
             \mathcal U_r^{\rm p}\mathcal L_B\xi_0^{\rm p}\,\dd r.
 \end{aligned}
\end{equation}
For a principal potential, curl may be moved outside either integral.
\end{lemma}
\begin{proof}
The chain rule and the group property of the principal flow give
\[
 \begin{aligned}
 \partial_sX_s
 &=\bigl(\xi^{\rm c}+\mathcal U_s^{\rm c}\xi_0^{\rm p}\bigr)\circ X_s,\\
 X_s\circ X_r^{-1}
 &=X_s^{\rm c}\circ X_{s-r}^{\rm p}\circ(X_r^{\rm c})^{-1}.
 \end{aligned}
\]
Taking the corresponding tensor pushforwards proves \eqref{path:factorized-transport}.
Equations~\eqref{eq:generator} and \eqref{eq:curl} give
\[
 \nabla_s^{\rm c}\Theta_s^{\rm p}
 =\nabla_s^{\rm c}\xi_s^{\rm p}=0.
\]
The Jacobi identity gives
\[
 [\nabla_s^{\rm c},\mathcal L_F]
 =\mathcal L_{\nabla_s^{\rm c}F},
\]
which proves the commutator identity in
\eqref{path:principal-transport}. The intertwining formula is
\eqref{path:mixed-intertwining} applied on the corrector path.
By the field variations~\eqref{eq:path-linear-fields},
$\nabla_s^{\rm c}$ commutes with the material and magnetic Lie
derivatives on the corrector path. The induction equation on that
path gives the remaining commutator.

For the field increments, subtract the corrector variations from
the full variations. Their differences solve
\[
 \begin{aligned}
 (\partial_s+\mathcal L_{\xi_s})(v_s-v_s^{\rm c})
   &=\mathcal D_{t,s}^{\rm c}\xi_s^{\rm p}
     =\mathcal U_s^{\rm c}\mathcal D_t\xi_0^{\rm p},\\
 (\partial_s+\mathcal L_{\xi_s})(B_s-B_s^{\rm c})
   &=\mathcal L_{B_s^{\rm c}}\xi_s^{\rm p}
     =\mathcal U_s^{\rm c}\mathcal L_B\xi_0^{\rm p},
 \end{aligned}
\]
with zero initial data. For the velocity, the bracket
$[v_s-v_s^{\rm c},\xi_s]$ has been moved to the left; the magnetic
equation follows in the same way. Thus, for either source $F$
relative to the initial background, Duhamel's formula gives
\[
 \begin{aligned}
 \int_0^s\mathcal U_{s\leftarrow r}\mathcal U_r^{\rm c}F\,\dd r
 &=\mathcal U_s^{\rm c}\int_0^s\mathcal U_{s-r}^{\rm p}F\,\dd r\\
 &=\mathcal U_s^{\rm c}\int_0^s\mathcal U_r^{\rm p}F\,\dd r.
 \end{aligned}
\]
Here the first equality follows from
\eqref{path:factorized-transport}, and the second is a change of
variable, since $F$ is independent of the integration parameter.
This proves \eqref{path:factorized-field-increments}. The potential
formula follows because curl commutes with the Lie derivatives
and with both pushforwards.
\end{proof}

\paragraph{\textbf{The chart after the corrector pushforward.}}
The factorization also describes the evolution of the Lagrangian
chart. Let $\Psi(t)$ be a volume-preserving Lagrangian chart before the corrector
acts, with coordinates $y^\eta$ as in
\eqref{part:physical-chart-derivatives}. The corrected chart is
$\Psi(t)\circ(X_s^{\rm c}(t))^{-1}$, with coordinates
$\mathcal U_s^{\rm c}y^\eta$. Pushforward commutes with exterior
differentiation and preserves scalar multiplication, so
\[
 d(\mathcal U_s^{\rm c}y^\eta)=\mathcal U_s^{\rm c}dy^\eta,
 \qquad
 \mathcal U_s^{\rm c}\bigl(f\,dy^\nu\bigr)
   =(\mathcal U_s^{\rm c}f)\,d(\mathcal U_s^{\rm c}y^\nu).
\]
Thus the transported coordinates, coframe and dual frame form a
volume-preserving chart. Formula~\eqref{direct:differentiated-transport}
identifies their material and magnetic Lie derivatives in the corrected
background with the pushforwards of those in the original background.

For the principal potential
$f\,dy^\nu=\lambda_{q+1}^{-1}\mathfrak a
\varphi(\lambda_{q+1}y^k)\,dy^\nu$, both its amplitude and its
oscillatory coordinate are transported by the corrector flow and expressed
in this chart. In these coordinates, the
principal flow is generated by $\xi_0^{\rm p}$.
The estimates below use invariance of $f$ and $dy^\nu$ under the
principal flow, and the small phase displacement, to reduce
tensor pushforward estimates to scalar estimates. The principal phase change is
measured from $\lambda_{q+1}\mathcal U_s^{\rm c}y^k$, the base phase
after the corrector pushforward, as in
\eqref{principal:phase-corrector-composition}.

\paragraph{\textbf{Background comparison.}}
The deformation constructed from the local background also acts on
the actual iterate. Since the same deformation acts on both, their
field gaps satisfy homogeneous Lie transport equations. These
equations give an exact comparison of the endpoint stresses, which
will be used in Section~\ref{ssec:background-gap-fields}.

Let $(v,B,p,R)$ and $(\bar v,\bar B,\bar p,\bar R)$ be smooth
relaxed MHD solutions on the same space--time domain. Apply
\eqref{eq:path-fields} to each with the same $X_s$, and use bars
for the fields arising from the second initial pair.
Their full gaps are
\[
 G_s^v=v_s-\bar v_s,\qquad G_s^B=B_s-\bar B_s.
\]
Subtracting the two field-variation formulas cancels the common
term containing the time derivative of the deformation. The gaps therefore satisfy
\[
 \begin{gathered}
 (\partial_s+\mathcal L_{\xi_s})G_s^{v}=0,\qquad
 (\partial_s+\mathcal L_{\xi_s})G_s^{B}=0,\\
 G_0^{v}=v-\bar v,\qquad G_0^{B}=B-\bar B.
 \end{gathered}
\]
We use the symmetric inverse divergence $\mathcal R$ of
Subsection~\ref{sec:inverse-divergence}.
\begin{lemma}[Exact endpoint transfer]
If the local endpoint has pressure $\bar p_1$ and stress
$\bar R_1$, the actual endpoint has pressure
$p_1=\bar p_1+p-\bar p$ and stress
\begin{equation}\label{eq:transfer-stress}
 \begin{aligned}
 R_1&=\bar R_1+R-\bar R
       +\mathcal R\partial_t(G_1^{v}-G_0^{v})+(T_1-T_0),\\
 T_s&=v_s\otimes v_s-B_s\otimes B_s
       -\bar v_s\otimes\bar v_s+\bar B_s\otimes\bar B_s.
 \end{aligned}
\end{equation}
\end{lemma}
\begin{proof}
The homogeneous gap equations imply
\[
 \frac{\dd}{\dd s}\int_{\mathbb T^3}G_s^v\,\dd x=0,
 \qquad
 \int_{\mathbb T^3}\partial_t(G_1^v-G_0^v)\,\dd x=0.
\]
Indeed, integration by parts on the torus shows that the bracket of
the divergence-free fields has zero integral. The spatial mean is
therefore independent of $s$, and differentiation of the endpoint
difference in physical time gives the second identity. Hence
$\mathcal R$ applies to the time derivative of the gap increment.
Subtracting the initial momentum difference from the endpoint
difference, we obtain
\begin{align*}
 &\mathcal M_1-\overline{\mathcal M}_1
       +\nabla(p-\bar p)\\
 &\qquad=\ddiv(R-\bar R)
       +\partial_t(G_1^v-G_0^v)+\ddiv(T_1-T_0)\\
 &\qquad=\ddiv\bigl(R-\bar R
       +\mathcal R\partial_t(G_1^v-G_0^v)+T_1-T_0\bigr).
\end{align*}
Adding the local endpoint equation proves \eqref{eq:transfer-stress}
with the stated pressure.
\end{proof}

\subsection{Perturbations in the adapted class and path estimates}
\label{ssec:path-estimates}

By \eqref{path:factorized-transport}, it suffices to estimate the principal
flow in the original chart and then apply the corrector pushforward.
We first prove a scalar transport estimate in normalized chart derivatives.
We then use the phase displacement and invariant frame to estimate the
velocity increment, the magnetic increment and general path integrals.
The factorization applies to the full $s$-dependent LDF.

\paragraph{\textbf{Phase flow and path integrals.}}\label{ssec:phase}

Fix an adapted Lagrangian chart as in
\eqref{part:physical-chart-derivatives}. In this chart,
$D_t=\partial_t|_y$ and $D_B$ has constant coefficients and no
$\partial_{y^k}$ component. We regard all functions below as functions
of physical time $t$ and Lagrangian coordinates $y$. Write
$\lambda=\lambda_{q+1}$ and $y^\perp=(y^\nu,y^\zeta)$. For a transverse
scale $0<\vartheta\le1$ and a pair $\rho$ of derivative costs $\mathrm a_t,\mathrm a_B$, the
commuting normalized derivatives are
\[
 \lambda^{-1}\partial_{y^k},\qquad
 \vartheta\ell\partial_{y^\nu},\qquad \vartheta\ell\partial_{y^\zeta},
 \qquad \mathrm a_t^{-1}D_t,\qquad \mathrm a_B^{-1}D_B,
\]
and for integers $N,J\ge0$ we define the anisotropic chart norm
\begin{equation}\label{principal:profile-norm}
 \|F\|_{N,J;\vartheta,\rho}=
 \max_{p+|\sigma|+h+m\le N,\ h+m\le J}
 \lambda^{-p}(\vartheta\ell)^{|\sigma|}\mathrm a_t^{-h}\mathrm a_B^{-m}
 \|\partial_{y^k}^p\partial_{y^\perp}^{\sigma}D_t^hD_B^mF\|_\infty .
\end{equation}
The supremum is taken over the given physical time interval and chart
neighborhood, and the norm of a vector of scalar coefficients is taken
componentwise. As in \eqref{setup:class}, each derivative is divided by
its corresponding frequency: $\lambda$ in the $y^k$ direction and
$(\vartheta\ell)^{-1}$ in the transverse directions.

We shall use three choices of these norms. The \emph{original} family
$\|F\|_{N,J}=\|F\|_{N,J;1,\mathrm a}$ uses transverse length $\ell$
and the fast material and magnetic derivative costs of the local background.
The \emph{rescaled} family $\|F\|_{N,J;\vartheta,\mathrm f}$ uses
$\vartheta<1$ and the derivative costs of the next iterate. For the
\emph{ordinary} family, we first work in Lagrangian coordinates and set
$m=0$ and $\mathrm a_t=\lambda$ in \eqref{principal:profile-norm}.
Thus the time derivative order $h$ is at most $J$, and
$p+|\sigma|+h\le N$. Differentiation in the oscillation direction
includes the derivatives of the factor $\varphi(\lambda y^k)$.

\begin{lemma}[Anisotropic chart transport]
\label{principal:finite-transport}
Suppose $0<\tau_a\ell^{-1}\delta_{q+1}^{1/2}\le1$ and $\lambda\ell\ge1$.
Fix nonnegative integers $N,J$ and the norm
\eqref{principal:profile-norm} with the scales specified above. Let
$V$ be a spatial vector field expressed in the normalized basis
\[
 \lambda^{-1}\frac{\partial}{\partial y^k},\qquad
 \vartheta\ell\frac{\partial}{\partial y^\nu},\qquad
 \vartheta\ell\frac{\partial}{\partial y^\zeta}.
\]
In the norm $\|\cdot\|_{N,J;\vartheta,\rho}$, assume that the
longitudinal coefficient is bounded by
$C\tau_a\ell^{-1}\delta_{q+1}^{1/2}$ and each transverse
coefficient by $C\tau_a\delta_{q+1}^{1/2}/(\vartheta\ell)\le C$.
Assume that the characteristics of
$V$ remain in the specified chart neighborhood for $0\le s\le1$. Consider
\[
 \partial_su_s+V\cn u_s+\mathsf Z_su_s=g_s,
 \qquad u_0=u,
\]
where $\mathsf Z_s=0$ for a scalar equation, or $\mathsf Z_s$ is a
matrix with $\sup_s\|\mathsf Z_s\|_{N,J;\vartheta,\rho}\le C$. Then
\begin{equation}\label{principal:scalar-homogeneous}
 \sup_{0\le s\le1}\|u_s\|_{N,J;\vartheta,\rho}
 \le C\left(\|u\|_{N,J;\vartheta,\rho}
                    +\int_0^1\|g_s\|_{N,J;\vartheta,\rho}\,\dd s\right).
\end{equation}
Suppose, in addition, that the norm belongs to the original family and
that all coefficients of $V$ and $\mathsf Z_s$ are bounded in that norm by
$C\tau_a\ell^{-1}\delta_{q+1}^{1/2}$. Then
\begin{equation}\label{principal:scalar-two-actions}
 \|(V\cn)^jF\|_{N,J}
          \le C\tau_a^j\ell^{-j}\delta_{q+1}^{j/2}\|F\|_{N+j,J},\qquad j=1,2,
\end{equation}
provided the coefficients of $V$ satisfy the same bounds through
$(N+j-1,J)$ and $F$ is bounded through $(N+j,J)$, for $j=1,2$.
If the solution is
homogeneous, then
\begin{equation}\label{principal:scalar-defect}
 \sup_{0\le s\le1}\|u_s-u\|_{N,J}
                    \le C\tau_a\ell^{-1}\delta_{q+1}^{1/2}\|u\|_{N+1,J}.
\end{equation}
In both estimates the additional derivatives are spatial, and the sum
of the material and magnetic derivative orders remains at most $J$.
The constants depend only on $N,J$ and the
coefficient bounds in the normalized basis; in particular, they do not
depend on the possibly large spatial gradient
$\tau_a\lambda\delta_{q+1}^{1/2}$.
\end{lemma}
\begin{proof}
Let $\partial_*$ be a composition of the normalized derivatives in
\eqref{principal:profile-norm}, with the orders prescribed there. Since
these derivatives commute, applying $\partial_*$ to the equation gives
\[
 (\partial_s+V\cn+\mathsf Z_s)\partial_*u_s
 =\partial_*g_s-[\partial_*,V\cn]u_s-[\partial_*,\mathsf Z_s]u_s .
\]
Each term in the first commutator has at least one derivative on a
coefficient of $V$. The remaining derivatives of $u_s$, including the
spatial derivative from $V\cn$, therefore have sum at most $N$.
Their material and magnetic orders still have sum at most $J$, because
the derivative from $V\cn$ is spatial. The commutator with $\mathsf Z_s$
contains no additional derivative of $u_s$. Thus Leibniz' rule and the
coefficient hypotheses bound both commutators by
$C\|u_s\|_{N,J;\vartheta,\rho}$. Integrating along the characteristics
of $V$ and taking the maximum over $\partial_*$ gives
\[
 \|u_s\|_{N,J;\vartheta,\rho}\le\|u\|_{N,J;\vartheta,\rho}
 +\int_0^s\|g_a\|_{N,J;\vartheta,\rho}\,\dd a
 +C\int_0^s\|u_a\|_{N,J;\vartheta,\rho}\,\dd a ,
\]
and Gronwall's inequality proves \eqref{principal:scalar-homogeneous}.
At the original scales, the small coefficient bounds give the
stronger estimate
\[
 \|V\cn F\|_{N,J}
 \le C\tau_a\ell^{-1}\delta_{q+1}^{1/2}\|F\|_{N+1,J}.
\]
Apply this estimate twice, using one additional spatial derivative of
the coefficients in the second application, to obtain
\eqref{principal:scalar-two-actions}. If the solution is homogeneous,
subtracting its initial datum shows that $u_s-u$ satisfies the same
transport equation with zero initial data and source
$-(V\cn+\mathsf Z_s)u$. The source is bounded by
$C\tau_a\ell^{-1}\delta_{q+1}^{1/2}\|u\|_{N+1,J}$.
Equation~\eqref{principal:scalar-homogeneous} then proves
\eqref{principal:scalar-defect} through $(N,J)$, using the initial
datum through $(N+1,J)$.
\end{proof}

The transverse coefficients need not be small in
$\|\cdot\|_{N,J;\vartheta,\mathrm f}$. To retain the small factors in
\eqref{principal:scalar-two-actions}, we first estimate the one or two
scalar transport derivatives in the original family, and then apply
Lemma~\ref{app:finite-spatial-reserve}. This gives the following rescaled
estimate, also when the coefficients contain temporal profiles.

\begin{corollary}[Rescaled estimates]
\label{principal:rescaled-profile-continuation}
Let
\[
 \ell^{-1}\le\Lambda_1\le\min\{\lambda,(\tau_a\delta_{q+1}^{1/2})^{-1}\},
 \qquad \vartheta=(\ell\Lambda_1)^{-1},
\]
so that
\[
 \vartheta\ge\max\{(\ell\lambda)^{-1},\tau_a\ell^{-1}\delta_{q+1}^{1/2}\}.
\]
Suppose that every coefficient of $V$ in the original
normalized basis, every entry of $\mathsf Z_s$, the initial datum and
the source can be written as finite sums $\sum_hb_h(t,y)H_h(\lambda y^k)$,
where the profiles are fixed bounded functions and the coefficients satisfy
\[
 b_h\in\mathcal C_{N,J}(E;\ell^{-1},\rho)\cap\mathcal K_{N',J}(E').
\]
Assume that $N,J,E,E'$ satisfy \eqref{app:reserve-inequality} for
this $\Lambda_1$. Here $E=\tau_a\ell^{-1}\delta_{q+1}^{1/2}$ for
the coefficients of $V$ and $\mathsf Z_s$, $E=E_0$ for the initial
datum, and $E=E(s)$ for the source. All lossy bounds are assumed
through the same order $N'\le\infty$. Then
\begin{equation}\label{principal:rescaled-profile-output}
 \sup_{0\le s\le1}\|u_s\|_{\bar N,J;\vartheta,\mathrm f}
 \le C\Bigl(E_0+\int_0^1E(s)\,\dd s\Bigr)
\end{equation}
for every $\bar N\le N'$. In physical coordinates the solution satisfies
\[
 u_s\in\mathcal C_{N'-1,J}(C(E_0+\int_0^1E);\lambda,\mathrm f)
\]
uniformly in $s$, where $N'-1=\infty$ when $N'=\infty$, provided
the chart and tensor frames satisfy the bounds in
Section~\ref{ssec:material-partition}. For initial data or sources
obtained from a previously estimated solution, it suffices to assume
the rescaled bounds directly. The same conclusions hold in the ordinary
family, with $\lambda$ in place of the derivative costs, under the
assumed derivative bounds for the chart and reference velocity.
\end{corollary}
\begin{proof}
Apply Lemma~\ref{app:finite-spatial-reserve}(b) separately to the
coefficients, initial datum and source, with the corresponding value
of $E$. Their rescaled norms are bounded by $CE$ through order $N'$;
nonoscillatory coefficients are included by taking $H_h=1$.
Changing the normalized basis leaves the longitudinal coefficient of
$V$ bounded by $C\tau_a\ell^{-1}\delta_{q+1}^{1/2}$ and bounds each
transverse coefficient by $C\tau_a\delta_{q+1}^{1/2}/(\vartheta\ell)\le C$.
These bounds also hold for the normalized derivatives. Thus all the
hypotheses of Lemma~\ref{principal:finite-transport} hold, which proves
\eqref{principal:rescaled-profile-output}.

In physical coordinates, every chart spatial derivative divided by $\lambda$
is controlled by a derivative in the rescaled norm. This is immediate
in the longitudinal direction; in a transverse direction the ratio is
$(\lambda\vartheta\ell)^{-1}\le1$. The chart and frame bounds and
Leibniz' rule therefore give the asserted class through $(N'-1,J)$.
For ordinary time derivatives, first apply
Lemma~\ref{app:finite-spatial-reserve}(c) in Lagrangian coordinates,
and then expand $\partial_t=D_t-v\cn$. With $h$ time derivatives,
this uses the reference velocity only through time order $h-1$ and
the solution through time order $h$, as assumed.
\end{proof}

\paragraph{\textbf{The principal phase and frame.}}
In an oriented volume-preserving adapted chart put
\[
 f=\lambda^{-1}\mathfrak a(t,y)\varphi(\lambda y^k),\qquad
 \xi_0^{\rm p}=\curl(f\,dy^\nu).
\]
Assume that $f$ is smooth and supported away from the chart boundary.
Let $X_s^{\rm p}$ be the physical flow of $\xi_0^{\rm p}$, and let
$Y_s$ denote this flow in chart coordinates. The pushforward
$\mathcal U_s^{\rm p}$ acts on scalars by composition with
$Y_s^{-1}$. The profile argument below is $\lambda y^k$.

\begin{lemma}[Action on scalars]\label{phase:prepared-identities}
The principal flow preserves
\begin{equation}\label{phase:prepared-invariance}
 \mathcal U_r^{\rm p}f=f,\qquad
 \mathcal U_r^{\rm p}dy^\nu=dy^\nu,\qquad
 \mathcal U_r^{\rm p}df=df.
\end{equation}
For a smooth scalar $g$ and every $j\ge0$,
\begin{equation}\label{principal:scalar-lie-action}
 \mathcal L_{\xi_0^{\rm p}}^j(g\,dy^\nu)
       =((\xi_0^{\rm p}\cn)^jg)\,dy^\nu,\qquad
 \mathcal L_{\xi_0^{\rm p}}^j\curl(g\,dy^\nu)
       =\curl(((\xi_0^{\rm p}\cn)^jg)\,dy^\nu).
\end{equation}
Moreover
\begin{equation}\label{phase:scalar-cancellation}
 \xi_0^{\rm p}\cn(\lambda_{q+1}y^k)=-\varphi\partial_{y^\zeta}\mathfrak a,
 \qquad
 \xi_0^{\rm p}\cn\mathfrak a=\mathfrak a\varphi'\partial_{y^\zeta}\mathfrak a,
\end{equation}
and hence
\begin{equation}\label{path:intrinsic-phase-drift}
 \partial_r[(\lambda_{q+1}y^k)\circ X_r^{\rm p}]
 =-(\varphi\partial_{y^\zeta}\mathfrak a)\circ X_r^{\rm p}.
\end{equation}
\end{lemma}
\begin{proof}
Since the chart preserves volume,
\begin{equation}\label{phase:scalar-action}
 \iota_{\xi_0^{\rm p}}\mathrm{vol}=df\wedge dy^\nu,
 \qquad
 \xi_0^{\rm p}\cn g
 =\partial_{y^k}f\,\partial_{y^\zeta}g-\partial_{y^\zeta}f\,\partial_{y^k}g .
\end{equation}
Taking $g=f$ or $g=y^\nu$ gives
\[
 \xi_0^{\rm p}\cn f=\xi_0^{\rm p}\cn y^\nu=0.
\]
Cartan's identity consequently gives
\[
 \mathcal L_{\xi_0^{\rm p}}dy^\nu
 =\mathcal L_{\xi_0^{\rm p}}df=0.
\]
The homogeneous Lie transport equation proves
\eqref{phase:prepared-invariance}. Invariance of the one-form factor
then gives the first identity in \eqref{principal:scalar-lie-action}
by the product rule; commutation with exterior differentiation and
volume preservation give the curl identity. Finally, substituting
$f=\lambda_{q+1}^{-1}\mathfrak a\varphi$ into
\eqref{phase:scalar-action}, the two products of derivatives of the amplitude
cancel, giving \eqref{phase:scalar-cancellation}.
The chain rule gives the phase-drift formula.
\end{proof}

The curl product rule gives the directional derivative along $\xi_0^{\rm p}$ in the form
\begin{equation}\label{principal:scalar-action}
 \xi_0^{\rm p}\cn
 =\mathfrak a\varphi'\partial_{y^\zeta}
  +\lambda^{-1}\varphi\,(\nabla\mathfrak a\times dy^\nu)\cn.
\end{equation}
Its coefficients in the original normalized basis are
\[
 -(\partial_{y^\zeta}\mathfrak a)\varphi,\qquad 0,\qquad
 \ell^{-1}\mathfrak a\varphi'
 +(\ell\lambda)^{-1}(\partial_{y^k}\mathfrak a)\varphi.
\]
By \eqref{principal:primitive-sharp}, all three coefficients are bounded
by $C\tau_a\ell^{-1}\delta_{q+1}^{1/2}$. This uses one additional
spatial derivative of $\mathfrak a$, with the same bound on the sum
of the material and magnetic derivative orders. In contrast, the
Euclidean gradient contains $\lambda\mathfrak a\varphi''$, which
may be large. The normalization of the longitudinal-to-transverse
entry removes this factor $\lambda$. The small coefficient bounds
therefore allow us to estimate the phase and displacement by
Lemma~\ref{principal:finite-transport}.

\begin{corollary}[Phase, frame and path-integral bounds]
\label{principal:phase-flow}
Suppose the preceding principal LDF satisfies the coefficient
hypotheses of Lemma~\ref{principal:finite-transport} in the original,
rescaled or ordinary family. Assume these bounds on a neighborhood of
the support invariant under the forward and inverse flows, as is the
case when the LDF vanishes in a surrounding cut off region. Let $Y_s$
be the flow in chart coordinates and set
\[
 \psi_s=\lambda\bigl(y^k\circ Y_s^{-1}-y^k\bigr),\qquad
 d_s^\perp=\ell^{-1}\bigl(y^\perp\circ Y_s^{-1}-y^\perp\bigr).
\]
For integers $0\le J\le N$, assume that the normalized LDF
coefficients are bounded for compositions of spatial, material and magnetic derivatives
whose orders have sum at most $N$, with the sum of the material and magnetic
orders at most $J$. Then
\begin{equation}\label{principal:phase-displacement}
 \sup_{0\le s\le1}
       \bigl(\|\psi_s\|_{N,J}+\|d_s^\perp\|_{N,J}\bigr)
 \le C\tau_a\ell^{-1}\delta_{q+1}^{1/2},
\end{equation}
and the same estimate holds with $Y_s^{-1}$ replaced by $Y_s$, as
well as in the rescaled norm. In particular, the displacement divided
by the rescaled transverse length $\vartheta\ell$ is bounded by
$C\tau_a\ell^{-1}\delta_{q+1}^{1/2}/\vartheta$.
If the same coefficient bounds hold with one additional spatial
derivative, then
\begin{equation}\label{principal:phase-jacobian-bound}
 \sup_{0\le s\le1,\,\sigma\in\{-1,1\}}
 \|D_y(Y_s^\sigma)-\IId\|_{N,J}
 \le C\tau_a\lambda\delta_{q+1}^{1/2},
\end{equation}
with the same bound in the rescaled norm. Under this additional
spatial derivative assumption, the coframe and transverse vector
satisfy the following bounds in the original norm:
\begin{equation}\label{principal:phase-frame-bound}
 \begin{gathered}
 \mathcal U_s^{\rm p}dy^\nu=dy^\nu,\qquad
 \mathcal U_s^{\rm p}dy^k=dy^k+\lambda^{-1}d\psi_s,\\
 \ell\lambda\| (\mathcal U_s^{\rm p}\frac{\partial}{\partial y^\zeta}
                         -\frac{\partial}{\partial y^\zeta})^k\|_{N,J}
 +\|(\mathcal U_s^{\rm p}\frac{\partial}{\partial y^\zeta}
                         -\frac{\partial}{\partial y^\zeta})^\perp\|_{N,J}
       \le C\tau_a\ell^{-1}\delta_{q+1}^{1/2},\\
 \|\partial_{y^\zeta}(\mathcal U_s^{\rm p}y^k-y^k)\|_{N,J}
       \le C(\lambda\ell)^{-1}\tau_a\ell^{-1}\delta_{q+1}^{1/2}.
 \end{gathered}
\end{equation}
For a fixed smooth periodic profile $H$,
\begin{equation}\label{principal:composed-profile-bound}
 \sup_s\|H(\lambda y^k+\psi_s)\|_{N,J}\le C_H,\qquad
 \sup_s\|H(\lambda y^k+\psi_s)-H(\lambda y^k)\|_{N,J}
       \le C_H\tau_a\ell^{-1}\delta_{q+1}^{1/2},
\end{equation}
again with the same estimates in the rescaled norm. The difference
estimate uses derivatives of $H$ through order $N+1$.
For a scalar family $F_r$ and a weight $\omega\in L^1(0,1)$, we also have
\begin{equation}\label{principal:phase-path-average}
 \sup_{0\le s\le1}
 \left\|\int_0^s\omega(r)\mathcal U_r^{\rm p}F_r\,\dd r\right\|_{N,J}
 \le C\int_0^1|\omega(r)|\,\|F_r\|_{N,J}\,\dd r ,
\end{equation}
whenever the right-hand side is finite. In the rescaled norm, one uses
the corresponding bound on the scalar family, either assumed directly
or obtained from Corollary~\ref{principal:rescaled-profile-continuation}.
All additional derivatives above are spatial, so the sum of the material
and magnetic derivative orders remains at most $J$. Interpolation gives
the H\"older estimates using one further integer spatial derivative.
\end{corollary}
\begin{proof}
Scalar composition with $Y_s^{-1}$ solves the homogeneous transport
equation. Subtracting the coordinate functions gives the following
equations for the displacement:
\[
 \begin{aligned}
 (\partial_s+\xi_0^{\rm p}\cn)\psi_s
   &=(\partial_{y^\zeta}\mathfrak a)\varphi,\\
 (\partial_s+\xi_0^{\rm p}\cn)d_s^\nu&=0,\\
 (\partial_s+\xi_0^{\rm p}\cn)d_s^\zeta
   &=-\ell^{-1}\mathfrak a\varphi'
     -(\ell\lambda)^{-1}(\partial_{y^k}\mathfrak a)\varphi,\\
 (\psi_0,d_0^\perp)&=0.
 \end{aligned}
\]
The right-hand sides are, up to sign, the normalized LDF coefficients.
Their assumed bounds and Lemma~\ref{principal:finite-transport} give
\eqref{principal:phase-displacement}. The same argument applies in the
rescaled norm: the source is still bounded by
$\tau_a\ell^{-1}\delta_{q+1}^{1/2}$, and the transverse coefficient
$\tau_a\ell^{-1}\delta_{q+1}^{1/2}/\vartheta\le1$ satisfies the
hypothesis of the transport estimate. Replacing $\xi_0^{\rm p}$ by
$-\xi_0^{\rm p}$ gives the forward-flow estimate. Both signs are
admissible because the cut off region keeps their flows in the domain.

Next apply \eqref{principal:phase-displacement} at $(N+1,J)$.
Multiplying by the coordinate lengths and using $\lambda^{-1}\le\ell$
bounds the displacement by $C\tau_a\delta_{q+1}^{1/2}$.
A spatial derivative contributes at most $\lambda$ in either family,
since $(\vartheta\ell)^{-1}\le\lambda$. This proves
\eqref{principal:phase-jacobian-bound}. It requires the LDF coefficients
through $(N+1,J)$. Thus the primitive $\mathfrak a$ is needed with
spatial, material and magnetic derivative orders summing to at most
$N+2$, while the material and magnetic orders still sum to at most $J$.

The coframe formulas follow from $Y_s^\nu=y^\nu$ and commutation
with exterior differentiation. For the transverse vector, the
pushforward formula gives
\[
 \mathcal U_s^{\rm p}\frac{\partial}{\partial y^\zeta}-\frac{\partial}{\partial y^\zeta}
   =(\partial_{y^\zeta}(Y_s-\mathrm{id}))\circ Y_s^{-1}.
\]
Differentiating the forward displacement estimate in the transverse
direction bounds its longitudinal component by
$C(\ell\lambda)^{-1}\tau_a\ell^{-1}\delta_{q+1}^{1/2}$ and its
transverse components by $C\tau_a\ell^{-1}\delta_{q+1}^{1/2}$.
Composition with the inverse flow preserves these bounds by the scalar
transport estimate. The final estimate in
\eqref{principal:phase-frame-bound} follows from
\[
 \partial_{y^\zeta}(\mathcal U_s^{\rm p}y^k-y^k)
 =\lambda^{-1}\partial_{y^\zeta}\psi_s.
\]

For the profile difference, use the fundamental theorem of calculus:
\[
 H(\lambda y^k+\psi_s)-H(\lambda y^k)
 =\psi_s\int_0^1H'(\lambda y^k+a\psi_s)\,\dd a.
\]
Every positive-order normalized derivative of $\lambda y^k$ vanishes
except its first longitudinal derivative, which equals one. Apply
Leibniz' rule and the chain rule to the last identity, using the bounds
on $\psi_s$ and the derivatives of $H$, to obtain
\eqref{principal:composed-profile-bound}. No bound on the coordinate
itself is needed. For \eqref{principal:phase-path-average}, apply the
scalar transport estimate to each integrand before integrating in $r$;
thus each scalar is estimated at its own transported phase. The same
equations, differentiated in the ordinary family, give the ordinary time
estimates. The extra derivatives of the flow map used here are all spatial.
\end{proof}

The preceding estimate bounds the difference between the flow differential and the identity by
$C\tau_a\lambda\delta_{q+1}^{1/2}$. Where $\mathfrak a$ is constant,
the exact flow is
\[
 Y_s(y)=(y^k,y^\nu,y^\zeta+s\mathfrak a\varphi'(\lambda y^k)).
\]
The phase is fixed, but the longitudinal-to-transverse derivative
can be large. The pushforward of a general vector may consequently
contain a factor $\tau_a\lambda\delta_{q+1}^{1/2}$. Invariant frame
factors permit sharper estimates because only their scalar
coefficients are transported. For a scalar $F$, one derivative
$\xi_0^{\rm p}\cn F$ at order $(N,J)$ uses $F$ through $(N+1,J)$
and the LDF coefficients through $(N,J)$. Two such derivatives use
$F$ through $(N+2,J)$ and the coefficients through $(N+1,J)$.
Estimate~\eqref{principal:scalar-two-actions} supplies the
corresponding small factors, which persist under weighted path
integration. Similarly, the identity
\[
 \mathcal U_r^{\rm p}F-\mathcal U_a^{\rm p}F
 =-\int_a^r\mathcal U_b^{\rm p}(\xi_0^{\rm p}\cn F)\,\dd b
\]
expresses each difference in \eqref{eq:covariance} through one
application of $\xi_0^{\rm p}\cn$, giving two small factors in the
scalar covariance. A field depending on the deformation parameter contributes also
$\partial_rF_r$, as in \eqref{eq:covariance-difference}. For a
one-form $F\,dy^\nu$, invariance of $dy^\nu$ reduces pushforward
to the scalar action on $F$. Formula~\eqref{principal:scalar-lie-action}
then permits curl to be taken after pushforward, at the cost of one
additional spatial derivative in the scalar estimate.

When the first spatial derivative bound is available,
\eqref{principal:phase-displacement} gives
\[
 \begin{aligned}
 \lambda^{-1}\partial_{y^k}(\lambda y^k+\psi_s)
   &=1+O(\tau_a\ell^{-1}\delta_{q+1}^{1/2}),\\
 |\partial_{y^\perp}(\lambda y^k+\psi_s)|
   &\le C\ell^{-1}\tau_a\ell^{-1}\delta_{q+1}^{1/2}.
 \end{aligned}
\]
Since $\lambda\ell\ge1$, the bounded invertible chart differential
gives
\[
 c\lambda\le|\nabla_x(\lambda y^k+\psi_s)|\le C\lambda
\]
when $\tau_a\ell^{-1}\delta_{q+1}^{1/2}$ is sufficiently small,
as guaranteed by Lemma~\ref{iter:parameter-lemma}. After the corrector
pushforward, the exact comparison is
\begin{equation}\label{principal:phase-corrector-composition}
 \mathcal U_s^{\rm c}(\lambda y^k+\psi_r)
       =\lambda\mathcal U_s^{\rm c}y^k+\mathcal U_s^{\rm c}\psi_r ,
\end{equation}
so the principal phase increment is small relative to the corrected
base phase. Its difference from the original phase $\lambda y^k$
need not be small. The additional factor $\varepsilon_\tau$ in the
comparison with the leading velocity term follows from
\eqref{principal:principal-constant-bound}, where the amplitude
and phase derivatives cancel after expansion.

We next apply these estimates to the field increments and to general
path integrals. For invariant frame factors, scalar transport controls
the coefficients, after which we take curls and tensor products.
General vectors require in addition the differential of the flow in
the pushforward formula. Let $(v,B)$ be the local background, with
full and corrector path fields as in
\eqref{path:factorized-field-increments}.

\paragraph{\textbf{Velocity.}}
The affine velocity law and Lemma~\ref{lem:path-transport} give, for
the autonomous principal flow,
\begin{equation}\label{principal:velocity-duhamel}
 (\partial_tX_s^{\rm p})\circ(X_s^{\rm p})^{-1}
   =\int_0^s\mathcal U_r^{\rm p}\partial_t\xi_0^{\rm p}\,\dd r,\qquad
 (\mathcal U_s^{\rm p}-\IId)v
   =\int_0^s\mathcal U_r^{\rm p}[v,\xi_0^{\rm p}]\,\dd r,
\end{equation}
since the two left sides have zero initial value and satisfy,
respectively,
\[
 (\partial_s+\mathcal L_{\xi_0^{\rm p}})u_s=\partial_t\xi_0^{\rm p},
 \qquad
 (\partial_s+\mathcal L_{\xi_0^{\rm p}})u_s=[v,\xi_0^{\rm p}].
\]
The sum of the two equations has source
\[
 \partial_t\xi_0^{\rm p}+[v,\xi_0^{\rm p}]=\mathcal D_t\xi_0^{\rm p}.
\]
An estimate of $\partial_t\xi_0^{\rm p}$ alone would allow size
$\tau_a\lambda\delta_{q+1}^{1/2}$, because an ordinary time derivative
differentiates the oscillatory coordinate. In the sum, this
contribution cancels with the advective derivative in
$[v,\xi_0^{\rm p}]$ by $D_ty^k=0$. The coframe contribution also
vanishes, since $\mathcal D_tdy^\nu=d(D_ty^\nu)=0$ by commutation
of Lie and exterior derivatives. Consequently,
\[
 \mathcal D_t\xi_0^{\rm p}=\curl\bigl((D_tf)\,dy^\nu\bigr),\qquad
 D_tf=\lambda^{-1}(D_t\mathfrak a)\varphi .
\]
Equation~\eqref{principal:primitive-identity} expresses the material
derivative of the primitive as the amplitude multiplied by its temporal
profile, plus the correction.
The factor $\tau_a^{-1}$ from time differentiation cancels the factor
$\tau_a$ in the primitive, so the scalar source for the potential is
bounded by $\lambda^{-1}\delta_{q+1}^{1/2}$ in the normalized norms
under consideration. Since the principal flow preserves $dy^\nu$,
\eqref{path:factorized-field-increments} gives
\[
 v_s-v_s^{\rm c}=\mathcal U_s^{\rm c}\curl(f_{w,s}\,dy^\nu),\qquad
 f_{w,s}=\int_0^s\mathcal U_r^{\rm p}(D_tf)\,\dd r,
\]
and, for $0\le J\le j_1$,
\eqref{principal:phase-path-average} gives
\[
 \sup_s\|f_{w,s}\|_{N,J}\le C\lambda^{-1}\delta_{q+1}^{1/2}.
\]
Taking curl, with one additional spatial derivative, gives the velocity
bound before the corrector pushforward. Applying
\eqref{gg:corrector-class-ranges} then gives
\eqref{principal:fields-bound}. Applying the same argument
to the source in the rescaled and ordinary norms gives the corresponding
estimates. In the ordinary norm, $h$ time derivatives of the field use
$D_t\mathfrak a$ through time order $h$, rather than $\xi_0^{\rm p}$
through time order $h+1$. The material derivatives of the correction
in \eqref{principal:primitive-identity} are bounded by
\eqref{principal:primitive-classes}.

\paragraph{\textbf{Magnetic field.}}
Here $D_By^k=0$ and $D_By^\nu$ is a chart constant, so
$\mathcal L_Bdy^\nu=d(D_By^\nu)=0$ and
\[
 \mathcal L_B\xi_0^{\rm p}=\curl\bigl((D_Bf)\,dy^\nu\bigr),\qquad
 D_Bf=\lambda^{-1}(D_B\mathfrak a)\varphi .
\]
Magnetic differentiation leaves the fast temporal profile unchanged.
Thus the scalar source for the potential has size
$\lambda^{-1}\tau_a\lambda_\parallel\delta_{q+1}^{1/2}$.
In chart coordinates, the local magnetic field has constant components
and no $\frac{\partial}{\partial y^k}$ component. Consequently it does
not differentiate the rapidly oscillating phase, and $D_B\mathfrak a$
is estimated with its magnetic derivative cost. Invariance of the
coframe then gives
\[
 B_s-B_s^{\rm c}=\mathcal U_s^{\rm c}\curl(f_{b,s}\,dy^\nu),\qquad
 f_{b,s}=\int_0^s\mathcal U_r^{\rm p}(D_Bf)\,\dd r,
\]
and, for $0\le J\le j_1$,
\[
 \sup_s\|f_{b,s}\|_{N,J}
 \le C\lambda^{-1}\tau_a\lambda_\parallel\delta_{q+1}^{1/2}
                         \varepsilon_\tau^{-[J+1-j_1]^+},
\]
by \eqref{principal:phase-path-average}. Taking curl requires one
additional spatial derivative. The source also contains one additional
magnetic derivative of $\mathfrak a$, which gives the loss
$\varepsilon_\tau^{-1}$ at $m=j_1$. For ordinary time estimates,
this derivative has the same magnetic cost and requires no further
time differentiation. Applying the rescaled estimates to the same
source retains the factor $\varepsilon_\tau^{-[J+1-j_1]^+}$
for every $0\le J\le j_1$.

Finally, apply the corrector bounds and
\eqref{direct:differentiated-transport} to both estimates. Together
with \eqref{path:factorized-field-increments}, they give the velocity
and magnetic increment bounds relative to the corrected local background,
for every $0\le s\le1$. The argument uses the full $s$-dependent LDF.

\paragraph{\textbf{General path integrals.}}

The same argument applies to any source expressed in invariant
tensor factors. Assume the hypotheses of
Corollary~\ref{principal:phase-flow} in the original, rescaled or
ordinary norm, for the derivative orders in
\eqref{principal:profile-norm}. Let
$E_1,\ldots,E_L$ be smooth tensor fields of one fixed type, acted on by the
corresponding tensor pushforward, with
\[
 \partial_sE_j=0,\qquad \mathcal L_{\xi_0^{\rm p}}E_j=0,
 \qquad j=1,\ldots,L,
\]
and consider a source
\[
 F_r=\sum_{j=1}^L f_{r,j}E_j,\qquad
 \|f_{r,j}\|_{N,J}\le A_j(r),\qquad
 0\le A_j\in L^1(0,1).
\]
Invariance and \eqref{path:factorized-transport} give
\[
 \mathcal U_u^{\rm p}E_j=E_j,\qquad
 \mathcal U_{s\leftarrow r}\mathcal U_r^{\rm c}
 =\mathcal U_s^{\rm c}\mathcal U_{s-r}^{\rm p}.
\]
Thus Duhamel's formula reduces to scalar transport:
\begin{equation}\label{principal:structured-duhamel}
 \int_0^s\mathcal U_{s\leftarrow r}\mathcal U_r^{\rm c}F_r\,\dd r
 =\mathcal U_s^{\rm c}\sum_j c_{s,j}E_j,
 \qquad c_{s,j}=\int_0^s\mathcal U_{s-r}^{\rm p}f_{r,j}\,\dd r.
\end{equation}
Applying Lemma~\ref{principal:finite-transport} to each scalar
coefficient and integrating in $r$ gives
\begin{equation}\label{principal:structured-coefficient-bound}
 \sup_{0\le u\le1}\|\mathcal U_u^{\rm p}f_{r,j}\|_{N,J}\le CA_j(r),
 \qquad
 \sup_{0\le s\le1}\|c_{s,j}\|_{N,J}\le C\|A_j\|_{L^1(0,1)}.
\end{equation}
Corollary~\ref{principal:rescaled-profile-continuation} gives the same
bounds in the rescaled norm; one may also use bounds on the source
assumed directly in that norm. No independence from the path parameter
is required. In particular, under $u=s-r$ the source becomes $F_{s-u}$.

For an oscillatory coefficient $f_{r,j}=a_{r,j}H_j(\lambda y^k)$,
with a fixed smooth periodic profile, the integrand is
\[
 \mathcal U_u^{\rm p}f_{r,j}
 =(a_{r,j}\circ Y_u^{-1})H_j(\lambda y^k+\psi_u).
\]
Since $D_ty^k=D_By^k=0$ before the corrector pushforward,
Corollary~\ref{principal:phase-flow} bounds the phase displacement
$\psi_u$ and the composed profile. Leibniz' rule therefore gives
\eqref{principal:structured-coefficient-bound} with
\[
 A_j(r)=C_{H_j}\|a_{r,j}\|_{N,J},
\]
using the profile through derivative order $N$. The estimate applies
to each integrand with its own phase.

For tensor products, let $E,E'$ be invariant one-forms. Commutation
with curl and preservation of tensor products give
\[
 \mathcal U_u^{\rm p}\bigl(\curl(qE)\otimes\curl(q'E')\bigr)
 =\curl\bigl((\mathcal U_u^{\rm p}q)E\bigr)
       \otimes\curl\bigl((\mathcal U_u^{\rm p}q')E'\bigr);
\]
to estimate each exterior derivative or curl, use the scalar, frame and
LDF bounds through $(N+1,J)$. Only the spatial derivative order increases;
the sum of the material and magnetic orders and the number of ordinary
time derivatives are unchanged. Leibniz' rule bounds products of these
factors, and the corresponding path integral is bounded provided
$A_jA_k\in L^1(0,1)$ for the factors in question.

The tensor estimates also use the assumed derivative bounds for $E_j$.
If $\mathcal D_tE_j=\mathcal L_BE_j=0$, the material and magnetic
Lie derivatives act only on the scalar coefficients. Apply
\eqref{direct:differentiated-transport} and the chart and corrector bounds
to obtain the estimates in the corrected background. The H\"older bound
then follows by interpolation between consecutive integer spatial norms.

In our construction, both $dy^\nu$ and $df$ are invariant under the
principal flow. Only $dy^\nu$ is also invariant under the material and
magnetic Lie derivatives of the local background. For $df$, those
derivatives are $d(D_tf)$ and $d(D_Bf)$, which satisfy the estimates
already proved. For $dy^k$ and $\frac{\partial}{\partial y^\zeta}$,
we use \eqref{principal:phase-frame-bound} instead. Thus the argument
requires invariance only of the specified factors; it does not assume
invariance of the metric or of a general tensor frame.

\begin{corollary}[General vector transport]\label{principal:vector-transport}
Assume the phase and displacement hypotheses of
Corollary~\ref{principal:phase-flow} with one additional spatial derivative
of the LDF coefficients, in the original, rescaled or ordinary
family. For a vector $F$ whose chart components satisfy
$\|F\|_{N,J}\le A_F$,
\begin{equation}
 \sup_s\bigl(\|\mathcal U_s^{\rm p}F\|_{N,J}
           +\|\mathcal U_s^{\rm p}F-F\|_{N,J}\bigr)
 \le C(1+\tau_a\lambda\delta_{q+1}^{1/2})A_F ,
\end{equation}
where every norm is taken in the chosen family. The estimate uses only
the derivatives of $F$ that appear on the left. For
$\int_0^s\mathcal U_{s-r}^{\rm p}F_r\,\dd r$, the analogous bound
has the integral of the source bounds on the right. Under the assumed
frame bounds, these estimates also hold in physical coordinates and
in the corrected background.
\end{corollary}
\begin{proof}
The pushforward formula gives
\begin{equation}\label{adapt:gap-stretching-composition}
 \mathcal U_s^{\rm p}F-F
 =\bigl((D_yY_s-\IId)\circ Y_s^{-1}\bigr)(F\circ Y_s^{-1})
                  +(F\circ Y_s^{-1}-F).
\end{equation}
By \eqref{principal:phase-jacobian-bound}, the matrix $D_yY_s-\IId$
is bounded by $C\tau_a\lambda\delta_{q+1}^{1/2}$. Fix $s$ in this
matrix and apply the scalar transport estimate to each entry and to
each component of $F$. Leibniz' rule bounds the first term by
$C\tau_a\lambda\delta_{q+1}^{1/2}A_F$. For the second term, the
triangle inequality and the scalar composition bound give $CA_F$,
without differentiating $F$ further. Adding the bound for $F$ proves
the assertion for its pushforward. Applying this estimate with parameter
$s-r$ and integrating proves the assertion for the path integral.

The additional derivative assumption is used only in the matrix estimate:
the LDF coefficients are needed through $(N+1,J)$, and hence the
primitive with spatial, material and magnetic derivative orders summing
to at most $N+2$, with material and magnetic orders summing to at most
$J$. The derivatives of $F$ are exactly
those in the asserted norm, with no additional material, magnetic or
ordinary time derivatives.
\end{proof}

This completes the principal path estimates on the local background.
Section~\ref{ssec:background-gap-fields} estimates the additional
errors from its difference with the actual background.
\subsection{Class properties}\label{ssec:class-properties}

In this subsection we establish the estimates for the classes defined in
Section~\ref{ssec:fixed-scales}. These include the extension to higher
spatial derivatives, comparison between backgrounds, and boundedness
of products, smoothing operators and pushforwards.

\paragraph{\textbf{Higher spatial derivatives.}}

The next lemma extends the sharp estimates from order $N$ to order
$N'\ge N$, using the lossy estimates and a larger spatial frequency.
We use this lemma, also recorded as part (vii) of
Lemma~\ref{setup:calculus}, whenever sharp estimates beyond order $N$
are required.

\begin{lemma}[Extension to higher spatial derivatives]
\label{app:finite-spatial-reserve}
Fix $0\le J\le N\le N'\le\infty$ and $E,E'>0$. For the ordinary
estimates in (c), fix an independent integer $H\ge0$. Let $\rho$ be a
pair of derivative costs $\mathrm a_t,\mathrm a_B$, each bounded by
the corresponding component of $\mathrm f$. Assume that
$\ell^{-1}\le\Lambda_1\le\lambda_{q+1}$ satisfies
\begin{equation}\label{app:reserve-inequality}
 \frac{E'}E(\ell\Lambda_1)^{-(N-J)}
 \max\Bigl\{1,
 \frac{\ell^{-1}}{\lambda_{q+1}\delta_{q+1}^{1/2}},
 \frac{\ell^{-1}}{\lambda_{q+1}\delta_{B,q+1}^{1/2}}\Bigr\}^{J}
 \le1.
\end{equation}
\begin{enumerate}
\item[(a)] If $b\in\mathcal C_{N,J}(E;\ell^{-1},\rho)\cap\mathcal K_{N',J}(E')$,
then $b\in\mathcal C_{N',J}(CE;\Lambda_1,\mathrm f)$.
\item[(b)] Let $y^k$ be a coordinate of an adapted chart, annihilated
by the material and magnetic derivatives. Assume that the chart maps
satisfy the estimates of Section~\ref{ssec:material-partition}. Let
$H_h$ be fixed functions with bounded derivatives of every order, and
set $F=\sum_hb_hH_h(\lambda_{q+1}y^k)$, where each coefficient
satisfies (a). With $\vartheta=(\ell\Lambda_1)^{-1}$, the anisotropic
chart norm \eqref{principal:profile-norm} satisfies
$\|F\|_{\bar N,J;\vartheta,\mathrm f}\le CE$ for every $\bar N\le N'$.
In physical coordinates,
$F\in\mathcal C_{N'-1,J}(CE;\lambda_{q+1},\mathrm f)$, where
$N'-1=\infty$ when $N'=\infty$. These conclusions also allow frame
vectors or one-forms in $F$, provided they satisfy the same chart
estimates.
\item[(c)] If $b\in\mathcal O_{N,H}(E;\ell^{-1})\cap\mathcal O_{N',H}(E';\ell^{-1})$,
then $b\in\mathcal O_{N',H}(CE;\Lambda_1)$, and the products $F$ of
(b) lie in $\mathcal O_{N',H}(CE;\lambda_{q+1})$ provided the chart
and any frame factors satisfy the corresponding ordinary estimates
through time order $H$. For this product assertion, the map estimates
are assumed through $r+h\le N'+1$, and the frame estimates through
$r+h\le N'$, with $h\le H$ and with the corresponding H\"older
bounds. No relation between $H$ and $J$ is required.
\end{enumerate}
The constants depend on the number of summands and on the fixed
functions $H_h$.
\end{lemma}
\begin{proof}
(a) For $r+k+m\le N$, the sharp estimate implies the desired bound,
since the quotient of the two bounds is
\[
 (\ell\Lambda_1)^{-r}
 \Bigl(\frac{\mathrm a_t}{\lambda_{q+1}\delta_{q+1}^{1/2}}\Bigr)^k
 \Bigl(\frac{\mathrm a_B}{\lambda_{q+1}\delta_{B,q+1}^{1/2}}\Bigr)^m\le1 .
\]
For derivatives of higher order, we have
\[
 N<r+k+m\le N',\qquad r>N-k-m\ge N-J.
\]
We therefore use the lossy estimate, whose quotient with the desired
bound satisfies
\[
 \frac{E'}E(\ell\Lambda_1)^{-r}
 \Bigl(\frac{\ell^{-1}}{\lambda_{q+1}\delta_{q+1}^{1/2}}\Bigr)^k
 \Bigl(\frac{\ell^{-1}}{\lambda_{q+1}\delta_{B,q+1}^{1/2}}\Bigr)^m\le1
\]
by \eqref{app:reserve-inequality}. Since $\ell^{-\alpha}$ is at most
$(\Lambda_1)_\alpha$, the same argument applies directly to the assumed
H\"older estimates, without using an additional spatial derivative.
This proves (a).

(b) Both transport derivatives annihilate $y^k$ and hence act only on
the coefficient in $b_hH_h(\lambda_{q+1}y^k)$. Applying the spatial
product rule, consider a term in which $i$ of the $p$ longitudinal
derivatives act on $b_h$. All $|\sigma|$ transverse derivatives also
act on this coefficient. Put $d=i+|\sigma|$ and multiply by the
factor from the definition of the chart norm,
\[
 \lambda_{q+1}^{-p}(\vartheta\ell)^{|\sigma|}\mathrm a_t^{-k}\mathrm a_B^{-m},
\]
with the final material and magnetic derivative costs. If
$d+k+m\le N$, the sharp coefficient estimate gives
\[
 CE(\ell\lambda_{q+1})^{-i}\vartheta^{|\sigma|}
 \le CE\vartheta^d.
\]
Otherwise $N-J<d\le N'$, and the lossy estimate gives
\[
 CE(E'/E)\vartheta^d
 \max\{1,\ell^{-1}/\mathrm a_t^{\rm f},\ell^{-1}/\mathrm a_B^{\rm f}\}^{J}
 \le CE.
\]
Here we used \eqref{app:reserve-inequality} and
$\vartheta^{-1}=\ell\Lambda_1$. The remaining longitudinal derivatives
act on the fixed function $H_h$; their contribution is bounded by a
constant depending on the function and derivative order. Summing over
the product-rule terms and the profiles proves the chart estimate.

To pass to physical coordinates, apply the chain rule and the chart
estimates. Each longitudinal derivative contributes a factor
$\lambda_{q+1}$, and each transverse derivative a factor
$\Lambda_1\le\lambda_{q+1}$. Interpolation between consecutive spatial
derivative orders gives the factor $\lambda_{q+1}^{\alpha}$ in the
H\"older estimate. Thus the spatial, material and magnetic derivative
orders in the latter estimate may sum to at most $N'-1$ when
$N'$ is finite.

(c) For the ordinary derivative classes, fix $h\le H$. The same
comparison gives
\[
 \begin{aligned}
 (\ell\Lambda_1)^{-(r+h)}&\le1 &&\text{if }r+h\le N,\\
 (E'/E)(\ell\Lambda_1)^{-(r+h)}&\le1
       &&\text{if }N<r+h\le N',
 \end{aligned}
\]
by \eqref{app:reserve-inequality}. For the products, ordinary time
differentiation acts on the profile argument through
$\partial_ty^k=-v\cn y^k$. The chart bounds and the chain rule
therefore give at most one factor $\lambda_{q+1}$ for each such
derivative. The assumed ordinary chart and frame bounds apply with
the same time index $H$; neither comparison increases it.
\end{proof}

\paragraph{\textbf{Changing transports.}}

The next lemma compares transport and Lie derivatives in two
backgrounds. We retain the number of magnetic derivatives falling on
the scalar factor, which will be needed when its estimates deteriorate
at the highest magnetic derivative orders.

\begin{lemma}[Change of transport derivatives]
\label{aniso:old-frame-operator-change}
Let $D_t=\partial_t+v_*\cn$ and $D_B=B_*\cn$. Assume that
\[
 [D_t,D_B]=0,\qquad J,N\ge1,\quad
 \Lambda\ge1,\quad \mathrm a_t,\mathrm a_B>0.
\]
Suppose the vector fields $u,b$ and the background gradients satisfy
\[
 \begin{aligned}
 \Lambda\|D_t^pD_B^ju\|_r+\|D_t^pD_B^j\nabla v_*\|_r
 &\le C\Lambda^r\mathrm a_t^{p+1}\mathrm a_B^j,\\
 \Lambda\|D_t^pD_B^jb\|_r+\|D_t^pD_B^j\nabla B_*\|_r
 &\le C\Lambda^r\mathrm a_t^p\mathrm a_B^{j+1}.
 \end{aligned}
\]
for $p+j\le J-1$ and $r+p+j\le N-1$. For the assertion concerning
Lie derivatives of tensors, assume in addition, on the same range,
\[
 \|D_t^pD_B^j\nabla u\|_r\le C\Lambda^r\mathrm a_t^{p+1}\mathrm a_B^j,
 \qquad
 \|D_t^pD_B^j\nabla b\|_r\le C\Lambda^r\mathrm a_t^p\mathrm a_B^{j+1}.
\]
Let $H$ be a tensor of fixed type whose Cartesian components satisfy
\[
 \|D_t^pD_B^jH\|_r\le C\Lambda^r\mathrm a_t^p\mathrm a_B^j,
 \qquad p+j\le J,\quad r+p+j\le N.
\]
For every smooth scalar $f$, each component of
\[
 \partial_x^\gamma(\partial_t+\mathcal L_{v_*+u})^k
 \mathcal L_{B_*+b}^{\,m}(fH),
 \qquad |\gamma|=r,\quad k+m\le J,\quad r+k+m\le N,
\]
has the expansion
\begin{equation}\label{aniso:old-frame-normal-form}
 \begin{gathered}
 \sum_{\alpha,p,j}C_{\alpha,p,j}\partial_x^\alpha D_t^pD_B^jf,\\
 p\le k,\quad j\le m,\quad |\alpha|+p+j\le r+k+m,\\
 \|C_{\alpha,p,j}\|_0
 \le C\Lambda^{r-|\alpha|}\mathrm a_t^{k-p}\mathrm a_B^{m-j}.
 \end{gathered}
\end{equation}
The conclusion holds for any ordering of the $k$ operators
$\partial_t+\mathcal L_{v_*+u}$ and the $m$ operators
$\mathcal L_{B_*+b}$. It also holds with the transport derivatives
$\partial_t+(v_*+u)\cn$ and $(B_*+b)\cn$ in place of the Lie
derivatives; in this case the additional assumptions on $\nabla u$
and $\nabla b$ are unnecessary.
\end{lemma}
\begin{proof}
We use the commutator identities
\begin{equation}\label{aniso:normalized-commutators}
 [D_t,D_B]=0,\qquad
 [D_t,\partial_i]=-(\partial_iv_*^a)\partial_a,\qquad
 [D_B,\partial_i]=-(\partial_iB_*^a)\partial_a.
\end{equation}
to put the derivatives in the order
$\partial_x^\alpha D_t^pD_B^j$. Each commutator introduces a
background gradient and reduces by one the number of derivatives
acting on the function. On Cartesian tensor components, the operators
in the conclusion have the form
\[
 D_t+u\cn+\Gamma_t,
 \qquad D_B+b\cn+\Gamma_B,
\]
where $\Gamma_t$ and $\Gamma_B$ act algebraically on the tensor
indices and contain $\nabla(v_*+u)$ and $\nabla(B_*+b)$, respectively.
Their bounds follow from the gradient assumptions with costs
$\mathrm a_t$ and $\mathrm a_B$. Starting from $fH$, apply the
operators successively and use the product rule. The magnetic Lie
derivative of a term $C\partial_x^\alpha D_t^pD_B^jf$ is the sum of
\[
 \begin{aligned}
 &(D_B C)\partial_x^\alpha D_t^pD_B^jf,\qquad
 C\partial_x^\alpha D_t^pD_B^{j+1}f,\qquad
 C[D_B,\partial_x^\alpha]D_t^pD_B^jf,\\
 &b^i(\partial_iC)\partial_x^\alpha D_t^pD_B^jf,\qquad
 b^iC\partial_i\partial_x^\alpha D_t^pD_B^jf,\qquad
 \Gamma_BC\partial_x^\alpha D_t^pD_B^jf.
 \end{aligned}
\]
Only the second term increases $j$. Every factor $b$ is accompanied
by a spatial derivative, so its bound $\mathrm a_B/\Lambda$ gives
the required magnetic derivative cost. The gradient terms have cost
$\mathrm a_B$ as well. The analogous expansion for the material Lie
derivative has cost $\mathrm a_t$, and increases $p$ only when the
material derivative in the original background acts on the scalar.
Induction, followed by the spatial product rule, proves the index
restrictions and powers in \eqref{aniso:old-frame-normal-form}.

We now check that the coefficient bounds use only the stated
hypotheses. One operator is needed to introduce an increment or a
background gradient, leaving at most $k+m-1$ material and magnetic
derivatives combined on that factor. After spatial differentiation,
the sum of its spatial, material and magnetic orders is at most
$r+k+m-1\le N-1$. The tensor terms
also contain first spatial derivatives of the increments, which are
controlled by the additional gradient hypotheses. In contrast, all
derivatives may fall on $H$. The commutations
\eqref{aniso:normalized-commutators} preserve these restrictions, so
the product estimate gives the asserted bound for $C_{\alpha,p,j}$.

For $k+m\le J$ and $r+k+m\le N$, the preceding count requires the
increments and background gradients with spatial, material and
magnetic orders summing to at most $N-1$, of which at most $J-1$
are material or magnetic. For $H$, the corresponding bounds on
these two sums are $N$ and $J$, respectively. The argument
therefore applies also when $N\le J$. By the integer version of
Lemma~\ref{high:component-lie-conversion}, the hypothesis on $H$ may
be given in transport or Lie derivatives in the reference background,
with the same gradient assumptions. No commutation of the new
operators was used, so the proof applies to any ordering, even when
those operators do not commute. The constants depend only on $J,N$,
the tensor type and the constant in the hypotheses.

In particular, the expansion gives the implication
\[
 \begin{gathered}
 \|D_t^pD_B^jf\|_r\le F\Lambda^r\mathrm a_t^p\mathrm a_B^j
 \quad(p+j\le J,\ r+p+j\le N)\\
 \Longrightarrow\quad
 \|(\partial_t+\mathcal L_{v_*+u})^k
       \mathcal L_{B_*+b}^{\,m}(fH)\|_r
 \le CF\Lambda^r\mathrm a_t^k\mathrm a_B^m
 \quad(k+m\le J,\ r+k+m\le N).
 \end{gathered}
\]
Dropping $\Gamma_t,\Gamma_B$ proves the assertion for transport
derivatives of components. If the hypotheses hold in both the integer
and H\"older norms, Leibniz' rule and \eqref{setup:product} give the
same conclusion in $\|\cdot\|_{r+\alpha}$ with the fixed H\"older
factor from the hypotheses. That factor occurs only once in each
product estimate. In particular, for $J\le j_1+1$, setting
$\Lambda=\mathrm a_t=\mathrm a_B=\ell^{-1}$ and taking the
H\"older factor to be $\ell^{-\alpha}$ gives
\[
 f\in\mathcal K_{N,J}(F)
 \quad\Longrightarrow\quad
 fH\in\mathcal K_{N,J}(CF).
\]
The first class is defined in $(v_*,B_*)$ and the second in
$(v_*+u,B_*+b)$, using the transport or Lie derivatives specified
in the conclusion. The coefficient hypotheses are required in both
the integer and H\"older norms.
\end{proof}

\begin{corollary}[Preservation of smallness under a change of transport]
\label{aniso:endpoint-residual-gain}
Assume the coefficient hypotheses of
Lemma~\ref{aniso:old-frame-operator-change}. Fix an integer
$\bar N\ge0$, a constant $F>0$, positive derivative costs
$\mathrm a_{t,0},\mathrm a_{B,0}$, and $d\in\{0,1,2\}$.
Suppose that $J\ge \bar N+d$ and $0<e<1$, with
$\mathrm a_{t,0}/\mathrm a_t\le1$ and
$\mathrm a_{B,0}/\mathrm a_B\le e$. Let $f$ satisfy
\begin{equation}
 \|D_t^pD_B^jf\|_r
 \le CF\Lambda^r\mathrm a_{t,0}^p\mathrm a_{B,0}^j
       e^{\min\{\bar N,J-j-d\}},\qquad
 p+j\le J,\quad r+p+j\le N.
\end{equation}
The exponent in this hypothesis is allowed to be negative.
Then, for $k+m\le J$ and $r+k+m\le N$, we have
\begin{equation}
 \|(\partial_t+\mathcal L_{v_*+u})^k
              \mathcal L_{B_*+b}^{\,m}(fH)\|_r
 \le CF\Lambda^r\mathrm a_t^k\mathrm a_B^m e^{\bar N}.
\end{equation}
The same estimate holds for transport derivatives of components.
If the scalar and coefficient hypotheses hold with one additional
spatial derivative, then the corresponding H\"older estimate holds
with the factor $\Lambda^\alpha$ and the same material and magnetic
derivative orders.
\end{corollary}
\begin{proof}
After division by $F\Lambda^r\mathrm a_t^k\mathrm a_B^m$, a summand of
\eqref{aniso:old-frame-normal-form} is at most
\[
 C(\mathrm a_{t,0}/\mathrm a_t)^p(\mathrm a_{B,0}/\mathrm a_B)^j e^{\min\{\bar N,J-j-d\}}
 \le Ce^{j+\min\{\bar N,J-j-d\}}.
\]
Since $J\ge \bar N+d$, the exponent satisfies, for $0\le j\le J$,
\begin{equation}
 j+\min\{\bar N,J-j-d\}
 =\min\{\bar N+j,J-d\}\ge \bar N.
\end{equation}
Derivatives falling on $H$, an increment or a gradient do not
increase $j$. Thus the preceding inequality bounds every term of
the expansion. Summing over spatial derivative orders $i\le r$ and
using $\Lambda^i\le\Lambda^r$ proves the assertion. Interpolation
between consecutive integer spatial norms proves the H\"older estimate.
\end{proof}

\paragraph{\textbf{Class calculus.}}

We now collect the estimates used for products, derivatives,
multipliers and smoothing, together with the preceding comparison
results.

\begin{lemma}[Class calculus]\label{setup:calculus}
Let the backgrounds and derivative classes be as in
Section~\ref{ssec:fixed-scales}. In the ordinary derivative assertions,
fix an integer $H\ge0$ independently of the material and magnetic
derivative orders. The following
estimates hold for adapted backgrounds, with constants depending only
on the fixed orders.
\begin{enumerate}
\item[(i)] \emph{Products.}
One has
\[
 \begin{aligned}
 \mathcal C_{N,J}(E_1;\Lambda,\rho)
 \cdot\mathcal C_{N,J}(E_2;\Lambda,\rho)
 &\subset\mathcal C_{N,J}(CE_1E_2;\Lambda,\rho),\\
 \mathcal K_{N,J}(E_1')\cdot\mathcal K_{N,J}(E_2')
 &\subset\mathcal K_{N,J}(CE_1'E_2').
 \end{aligned}
\]
The classes $\mathcal O_{N,H}$ satisfy the analogous product
estimates with the same $N,H$. The loss of a product is the product
of the losses of its factors.
\item[(ii)] \emph{Derivatives.}
For $N\ge1$, a spatial derivative gives the mappings
\[
 \begin{aligned}
 \mathcal C_{N,J}(E;\Lambda,\rho)
 &\longrightarrow\mathcal C_{N-1,J}(E\Lambda;\Lambda,\rho),\\
 \mathcal K_{N,J}(E')
 &\longrightarrow\mathcal K_{N-1,J}(E'\ell^{-1}).
 \end{aligned}
\]
For an adapted background, the same bounds hold with the gradient
placed before the material and magnetic derivatives
(Lemma~\ref{high:transport-gradient}). For $N,J\ge1$, the transport
mappings are
\[
 \begin{aligned}
 D_t\colon\mathcal C_{N,J}(E;\Lambda,\rho)
 &\longrightarrow\mathcal C_{N-1,J-1}(E\mathrm a_t;\Lambda,\rho),\\
 D_B\colon\mathcal C_{N,J}(E;\Lambda,\rho)
 &\longrightarrow\mathcal C_{N-1,J-1}(E\mathrm a_B;\Lambda,\rho).
 \end{aligned}
\]
For the lossy class, the corresponding mappings are
\[
 D_t,D_B\colon\mathcal K_{N,J}(E')
 \longrightarrow\mathcal K_{N-1,J-1}(E'\ell^{-1}).
\]
\item[(iii)] \emph{Transport and Lie derivatives.} For an adapted
background, the classes defined by transport and Lie derivatives with
the same parameters coincide up to a fixed constant
(Lemma~\ref{high:component-lie-conversion}).
\item[(iv)] \emph{Change of transport.} If $N=0$ or $J=0$, the sharp
and lossy classes are defined by spatial estimates alone and hence are
independent of the background. For $N,J\ge1$, suppose
two backgrounds are $(\Lambda,\rho)$-adapted through $(N,J)$ and differ by
\[
 \begin{aligned}
 u=v'-v&\in\mathcal C_{N-1,J-1}(C\mathrm a_t\Lambda^{-1};\Lambda,\rho),\\
 b=B'-B&\in\mathcal C_{N-1,J-1}(C\mathrm a_B\Lambda^{-1};\Lambda,\rho),
 \end{aligned}
\]
where the bounds are measured in $(v,B)$. Then the classes
$\mathcal C_{N,J}(\,\cdot\,;\Lambda,\rho)$ defined using the two
backgrounds coincide up to a fixed constant. If in addition
$u,b\in\mathcal K_{N-1,J-1}(C)$, the same holds for the lossy
classes $\mathcal K_{N,J}$
(Lemma~\ref{aniso:old-frame-operator-change}). These conclusions also
hold when only integer derivative estimates are imposed.
\item[(v)] \emph{Order-zero Fourier multipliers.} For an adapted background and
$\Lambda\le\lambda_{q+1}$, each of the order-zero operators
$\mathbb P$, $\mathscr T$, $\mathcal R\curl$ and $\mathcal R\ddiv$
gives the mappings
\[
 \begin{aligned}
 \mathcal C_{N,J}(E;\Lambda,\rho)
 &\longrightarrow
 \mathcal C_{N,J}(C\lambda_{q+1}^{(J+1)\alpha}E;\Lambda,\rho),\\
 \mathcal K_{N,J}(E')
 &\longrightarrow\mathcal K_{N,J}(C\lambda_{q+1}^{(J+1)\alpha}E'),\\
 \mathcal O_{N,H}(E;\Lambda)
 &\longrightarrow
 \mathcal O_{N,H}(C\max\{\ell^{-1},\Lambda\}^{\alpha}E;\Lambda).
 \end{aligned}
\]
The ordinary estimate uses commutation with $\partial_t$ and
boundedness on $C^{r+\alpha}$
(Lemma~\ref{principal:finite-multiplier}).
\item[(vi)] \emph{Spatial smoothing.} For an adapted background,
\[
 \mathcal J_\ell\colon\mathcal C_{N,J}(E;\Lambda,\rho)
 \longrightarrow\mathcal C_{N,J}(CE;\Lambda,\rho)
\]
by Proposition~\ref{moll:differentiated-spatial}. If $N\ge H$
and $\Lambda\le\ell^{-1}$, then also
\[
 \mathcal J_\ell\colon\mathcal O_{N,H}(E;\Lambda)
 \longrightarrow
 \mathcal O_{N,H}(CE;\Lambda)\cap\mathcal O_{\infty,H}(CE;\ell^{-1}).
\]
If $H\ge j_1$ and $v,B\in\mathcal K_{j_1-1}(C)$ as well, the
smoothed function belongs to $\mathcal K(CE)$. The ordinary estimates
retain the same time index $H$, which may be decreased independently
provided $N\ge H$; additional spatial derivatives may be placed on
the convolution kernel. Moreover,
\[
 \|\mathcal J_\ell f-f\|_r\le C\ell^d\|f\|_{r+d}
\]
for $d$ at most one more than the number of vanishing moments
(Proposition~\ref{moll:spatial-basic}).
\item[(vii)] \emph{Higher spatial derivatives.} Let
\[
 b\in\mathcal C_{N,J}(E;\ell^{-1},\rho)\cap\mathcal K_{N',J}(E'),
\]
where $N\le N'\le\infty$ and the costs of $\rho$ do not exceed
the corresponding components of $\mathrm f$. Suppose that
$\ell^{-1}\le\Lambda_1\le\lambda_{q+1}$ satisfies
\begin{equation}\label{setup:reserve-inequality}
 \frac{E'}{E}(\ell\Lambda_1)^{-(N-J)}
 \max\Bigl\{1,\frac{\ell^{-1}}{\lambda_{q+1}\delta_{q+1}^{1/2}},
             \frac{\ell^{-1}}{\lambda_{q+1}\delta_{B,q+1}^{1/2}}\Bigr\}^J\le1 .
\end{equation}
Then
\[
 b\in\mathcal C_{N',J}(CE;\Lambda_1,\mathrm f).
\]
If
\[
 F=\sum_hb_hH_h(\lambda_{q+1}y^k)
\]
is a finite sum with coefficients satisfying the same hypotheses
and fixed bounded profiles depending on an adapted chart coordinate
annihilated by the material and magnetic derivatives, then the chart
estimate of Lemma~\ref{app:finite-spatial-reserve}(b) holds through
order $N'$ with transverse scale $(\ell\Lambda_1)^{-1}$. Constant
profiles $H_h=1$ are allowed. In physical coordinates,
\[
 F\in\mathcal C_{N'-1,J}(CE;\lambda_{q+1},\mathrm f),
\]
where $N'-1=\infty$ when $N'=\infty$. For finite $N'$, the loss
of one derivative is due to H\"older interpolation. Under the same
hypothesis,
\[
 \mathcal O_{N,H}(E;\ell^{-1})\cap\mathcal O_{N',H}(E';\ell^{-1})
 \subset\mathcal O_{N',H}(CE;\Lambda_1)
\]
by Lemma~\ref{app:finite-spatial-reserve}.
\end{enumerate}
\end{lemma}
\begin{proof}
For (i), apply Leibniz' rule and \eqref{setup:product}. In each term,
the spatial, material and magnetic derivative orders add separately
to those of the product. Multiplying the estimates for the factors
therefore gives the asserted amplitude and derivative costs. The
H\"older estimate uses the H\"older norm of one factor and the
integer norm of the other. The same argument proves the assertions
for the lossy and ordinary classes.

For (ii), the definitions give the stated mappings: a spatial
derivative contributes the spatial frequency and reduces the allowed
sum of spatial, material and magnetic derivative orders by one.
A material or magnetic derivative contributes its respective cost and also reduces the allowed sum of
material and magnetic derivative orders by one. If the gradient is
placed before these derivatives, we commute it through using
Lemma~\ref{high:transport-gradient}, with the class indices in the
displayed mappings.

For (iii), apply Lemma~\ref{high:component-lie-conversion}; its
gradient hypotheses are part of the definition of an adapted background.

For (iv), first consider transport derivatives. Apply
Lemma~\ref{aniso:old-frame-operator-change} with $H=1$ to each
Cartesian component. The gradient bounds for $(v,B)$ and the assumed
bounds for $u,b$ give one inclusion for class indices $(N,J)$;
no gradient bounds for the increments are needed for this version of the lemma.
To prove the converse, first apply the same expansion to $u$ and $b$
with class indices $(N-1,J-1)$. The same hypotheses apply to these
indices and give
\[
 \begin{aligned}
 u&\in\mathcal C_{N-1,J-1}(C\mathrm a_t\Lambda^{-1};\Lambda,\rho),\\
 b&\in\mathcal C_{N-1,J-1}(C\mathrm a_B\Lambda^{-1};\Lambda,\rho)
 \end{aligned}
\]
in the background $(v',B')$. If $N=1$ or $J=1$, these are only
spatial estimates and hence are independent of the background. In
all cases, $-u,-b$ satisfy the increment hypotheses for the converse
comparison. We may therefore apply the lemma in $(v',B')$, which is
adapted by assumption, to obtain the reverse inclusion. Part (iii)
in each background then gives the same conclusion for tensor Lie
derivatives.

Each product estimate uses the H\"older factor only once. Thus the
two inclusions preserve both the integer and H\"older estimates;
the integer estimates themselves require only integer hypotheses.
For the lossy classes, use the gradient bounds from
Section~\ref{ssec:fixed-scales}, the assumption
$u,b\in\mathcal K_{N-1,J-1}(C)$, and derivative cost $\ell^{-1}$
in the same argument. The cases $N=0$ or $J=0$ follow from the
definition, since only spatial estimates occur. The case $N=\infty$
follows by applying the preceding argument at each integer value
of the first class index.

For (v), apply Lemma~\ref{principal:finite-multiplier}. The first
assertion and the approximation estimate in (vi) follow from
Propositions~\ref{moll:differentiated-spatial}
and~\ref{moll:spatial-basic}; adaptedness supplies the gradient
hypotheses in each case. To obtain the ordinary estimates at higher
spatial orders, place the spatial derivatives on the convolution
kernel. For $h\le H\le N$ and every $r\ge0$, this gives
\[
 \|\partial_t^h\mathcal J_\ell f\|_r
 \le C\ell^{-r}\|\partial_t^hf\|_0
 \le CE\ell^{-r}\Lambda^h
 \le CE\ell^{-(r+h)}.
\]
One additional spatial derivative of the kernel gives the H\"older
estimate. For the transport conclusion, $H\ge j_1$ supplies all
ordinary time derivatives needed in expanding $D_t^kD_B^m$ with
$k+m\le j_1$. This expansion requires at most
$k+m-1\le j_1-1$ derivatives of the background coefficients. Their
lossy bounds therefore imply membership in $\mathcal K(CE)$. The
ordinary argument preserves the prescribed index $H$. Finally, (vii) follows from
Lemma~\ref{app:finite-spatial-reserve} under
\eqref{setup:reserve-inequality}; the physical H\"older estimate uses
one additional spatial derivative, as in that lemma.
\end{proof}

\paragraph{\textbf{Pushforward of classes.}}\label{ssec:transported-charts}

Let $X_s$ be the flow of an $s$-independent, divergence-free
LDF $\xi$, and write $\mathcal U_s=(X_s)_*$ and
$\mathcal P_X[F]=\int_0^1\mathcal U_sF\,\dd s$.
Given a pair $(v,B)$ satisfying the induction equation, define
\[
 v_s=\partial_tX_s\circ X_s^{-1}+\mathcal U_sv,
 \qquad B_s=\mathcal U_sB.
\]
We first estimate pushforwards in the transported background
$(v_s,B_s)$. The commutation identities reduce this estimate to spatial
derivatives of $\xi$ of order at least one. To compare with the
original background $(v,B)$ and to estimate the pushforward defect,
we also assume bounds on the displacement and its first material and
magnetic Lie variations.

\begin{lemma}[Pushforward of classes]\label{gg:transport-classes}
Fix $N\ge J\ge1$, $\Lambda\ge1$ and
$\rho=(\mathrm a_t,\mathrm a_B)$. Let
$F\in\mathcal C_{N,J}(E;\Lambda,\rho)$, with the class defined in
the background $(v,B)$ using material and magnetic Lie derivatives.
The classes below are also defined using Lie derivatives.
\begin{enumerate}
\item[(i)] Assume that
$\nabla\xi\in\mathcal C_{N,0}(C_0;\Lambda,\rho)$. Then, for
$s\in[0,1]$, the following estimate holds in $(v_s,B_s)$:
\[
 \mathcal U_sF\in\mathcal C_{N,J}(CE;\Lambda,\rho),
 \qquad 0\le s\le1.
\]
\item[(ii)] Suppose $(v,B)$ is $(\Lambda,\rho)$-adapted through $(N+1,J)$ and
\[
 \begin{aligned}
 \xi&\in\mathcal C_{N+1,0}(\epsilon\Lambda^{-1};\Lambda,\rho),
       &&0\le\epsilon\le C_0,\\
 \mathcal D_t\xi
   &\in\mathcal C_{N+1,J-1}
       (\epsilon\mathrm a_t\Lambda^{-1};\Lambda,\rho),\\
 \mathcal L_B\xi
   &\in\mathcal C_{N+1,J-1}
       (\epsilon\mathrm a_B\Lambda^{-1};\Lambda,\rho).
 \end{aligned}
\]
Then the following estimates hold uniformly for $s\in[0,1]$, both
in the original background and in each transported background with
the deformation parameter fixed:
\[
 \begin{aligned}
 \mathcal U_sF,\ \mathcal P_X[F]
   &\in\mathcal C_{N,J}(CE;\Lambda,\rho),\\
 (\mathcal U_s-\IId)F,\ \mathcal P_X[F]-F
   &\in\mathcal C_{N-1,J}(C\epsilon E;\Lambda,\rho).
 \end{aligned}
\]
\end{enumerate}
The constants depend only on the derivative orders, the tensor type,
$\alpha$, $C_0$ and the constants in the adaptedness assumptions.
\end{lemma}
\begin{proof}
For (i), use the commutation identity
\[
 \mathcal D_{t,s}^{k}\mathcal L_{B_s}^{m}\mathcal U_sF
   =\mathcal U_s\mathcal D_t^k\mathcal L_B^mF.
\]
It remains to estimate spatial derivatives of the pushforward, as in
Lemma~\ref{high:natural-transport}. Differentiate its Lie transport
equation in space and integrate along characteristics. Each coefficient
outside the transport term contains at least one spatial derivative of
$\xi$; at spatial order $r$, the largest required order is $r+1$.
These derivatives are bounded by the hypothesis on $\nabla\xi$.

The integer estimates follow from the differentiated equation. For
the H\"older estimate, subtract the equations along two characteristics.
Their separation is controlled by $\|\nabla\xi\|_0$. Apply the
H\"older product estimate, placing either the coefficient or the
transported tensor in the H\"older seminorm. After division by
$\max\{\ell^{-1},\Lambda\}^{\alpha}$, the coefficient terms are
bounded by the integer estimates. The remaining terms contain the
H\"older seminorm of the tensor divided by the same factor. Gronwall's inequality
gives \eqref{high:natural-transport-bound} with the class weights.
Only derivatives of $\xi$ occur in these estimates, so no bound on
$\xi$ itself is needed.

For (ii), first work in the original background. By
Lemma~\ref{lem:path-transport},
\[
 \begin{aligned}
 [\mathcal D_t,\mathcal U_s]
 &=-\int_0^s\mathcal U_{s-u}
       \mathcal L_{\mathcal D_t\xi}\mathcal U_u\,\dd u,\\
 [\mathcal L_B,\mathcal U_s]
 &=-\int_0^s\mathcal U_{s-u}
       \mathcal L_{\mathcal L_B\xi}\mathcal U_u\,\dd u.
 \end{aligned}
\]
The Lie derivative in each integrand has a first-order coefficient of
size $\epsilon\mathrm a_t/\Lambda$ or
$\epsilon\mathrm a_B/\Lambda$. The first-order part differentiates
the tensor once, with spatial cost $\Lambda$. The zeroth-order part
differentiates this coefficient once and satisfies the same bound.
After division by the material or magnetic derivative cost, the
remaining factor is at most $C\epsilon$. Commute successive Lie
derivatives through the pushforwards and apply the spatial assertion
of (i). This proves the estimate by induction on the number of Lie
derivatives. Indeed, each first variation introduced by a commutator
receives at most $J-1$ further Lie derivatives combined, and its
zeroth-order tensor term uses one additional spatial derivative.
Thus it suffices that the spatial, material and magnetic orders on
each first variation sum to at most $N+1$, with at most $J-1$
material and magnetic Lie derivatives combined. No estimate of
$\xi$ with $J$ material and magnetic Lie derivatives combined is
used. Adaptedness allows conversion to transport derivatives
and commutation with spatial derivatives. The H\"older product
estimate again uses its additional factor only once. We have therefore
proved the asserted bounds in the original background.

The affine field identities give
\[
 \begin{aligned}
 c_v=v_s-v&=\int_0^s\mathcal U_u\mathcal D_t\xi\,\dd u,\\
 c_B=B_s-B&=\int_0^s\mathcal U_u\mathcal L_B\xi\,\dd u.
 \end{aligned}
\]
Apply the estimate just proved to the first variations and integrate
in the deformation parameter. With at most $J-1$ material and
magnetic Lie derivatives combined, we obtain in the original background
\[
 \begin{aligned}
 c_v&\in\mathcal C_{N,J-1}
   (C\epsilon\mathrm a_t\Lambda^{-1};\Lambda,\rho),&
 \nabla c_v&\in\mathcal C_{N-1,J-1}
   (C\epsilon\mathrm a_t;\Lambda,\rho),\\
 c_B&\in\mathcal C_{N,J-1}
   (C\epsilon\mathrm a_B\Lambda^{-1};\Lambda,\rho),&
 \nabla c_B&\in\mathcal C_{N-1,J-1}
   (C\epsilon\mathrm a_B;\Lambda,\rho).
 \end{aligned}
\]
Apply Lemma~\ref{aniso:old-frame-operator-change}, in the original
background, to $\nabla v+\nabla c_v$ and $\nabla B+\nabla c_B$.
The preceding estimates show that $(v_s,B_s)$ is adapted through
$(N,J)$. Part (iv) of Lemma~\ref{setup:calculus} now identifies the
classes in the two backgrounds, up to a fixed constant.

To estimate the defect, first use the equation
\[
 (\partial_s+\mathcal L_\xi)(\mathcal U_sF-F)
       =-\mathcal L_\xi F,\qquad
 \mathcal U_0F-F=0.
\]
The right-hand side has amplitude $C\epsilon E$ and spatial
derivative cost $\Lambda$ through order $N-1$, since it uses one
additional derivative of $F$. This proves the spatial estimate.
For material and magnetic Lie derivatives, we use
\[
 \begin{aligned}
 \mathcal D_{t,s}(\mathcal U_sF-F)
   &=(\mathcal U_s-\IId)\mathcal D_tF-\mathcal L_{c_v}F,\\
 \mathcal L_{B_s}(\mathcal U_sF-F)
   &=(\mathcal U_s-\IId)\mathcal L_BF-\mathcal L_{c_B}F.
 \end{aligned}
\]
After repeated differentiation, each term is either a spatial defect
applied to an iterated Lie derivative in the original background, or
contains $\mathcal L_{c_v}$ or $\mathcal L_{c_B}$. The first kind is
bounded by the spatial defect estimate. The increment estimates give
an additional factor $C\epsilon$ for the second kind. Induction and
the product rule therefore prove the defect estimate when the
spatial, material and magnetic orders sum to at most $N-1$, with
at most $J$ material and magnetic Lie derivatives combined.
Transfer this estimate to the original background using
the equivalence of classes, and integrate to obtain the assertions for
$\mathcal P_X$. The same equivalence gives these integral estimates
in every fixed transported background.
\end{proof}

The proof of (ii) also shows that the transported background is adapted
through $(N,J)$, so its transport and Lie derivative classes are
equivalent. The spatial conclusions hold for $J=0$. For
$J\le j_1+1$, setting $\Lambda=\ell^{-1}$ and
$\rho=(\ell^{-1},\ell^{-1})$ in the lemma gives
\[
 F\in\mathcal K_{N,J}(E)
 \quad\Longrightarrow\quad
 \mathcal U_sF\in\mathcal K_{N,J}(CE),
\]
where the first class is defined in $(v,B)$ and the second in
$(v_s,B_s)$, using Lie derivatives. Under the hypotheses of (ii),
we also have
\[
 \begin{aligned}
 \mathcal U_sF,\ \mathcal P_X[F]
   &\in\mathcal K_{N,J}(CE),\\
 (\mathcal U_s-\IId)F,\ \mathcal P_X[F]-F
   &\in\mathcal K_{N-1,J}(C\epsilon E)
 \end{aligned}
\]
in the original background and in every fixed transported background.
The conclusions with first class index $\infty$ follow by taking
$N$ arbitrarily large, provided the hypotheses hold for every such
choice.

There is an analogous statement for ordinary time derivatives, with
an independent integer bound $K$ on the ordinary time order.
Assume that the derivatives of $\xi$ of order at least one satisfy
\[
 \begin{gathered}
 \|\partial_t^h\DD^r\xi\|_0
 +\max\{\ell^{-1},\Lambda\}^{-\alpha}
       \|\partial_t^h\DD^r\xi\|_\alpha
       \le C_0\Lambda^{r+h-1},\\
 1\le r+h\le N+1,\qquad h\le K,
 \end{gathered}
\]
Then the ordinary class is preserved:
\[
 \mathcal U_s\colon\mathcal O_{N,K}(E;\Lambda)
 \longrightarrow\mathcal O_{N,K}(CE;\Lambda).
\]
To see this,
apply $\partial_t^h\DD^r$ to the Lie transport equation and use the
integer and H\"older estimates of
Lemma~\ref{high:natural-transport}. Every commutator differentiates
$\xi$ at least once, and the zeroth-order tensor term uses one
additional spatial derivative. No term contains more than $K$ time
derivatives. Taking $N$ arbitrarily large proves the assertion
when the spatial hypotheses hold at all orders.

Finally, pushforward of a coordinate chart and its dual frames gives
\[
 \begin{gathered}
 y_s^\eta=y^\eta\circ X_s^{-1},\qquad
 dy_s^\eta=\mathcal U_sdy^\eta,\qquad
 \frac{\partial}{\partial y_s^\eta}=\mathcal U_s\frac{\partial}{\partial y^\eta},\\
 \mathcal D_{t,s}y_s^\eta=\mathcal U_sD_ty^\eta,
 \qquad \mathcal L_{B_s}y_s^\eta=\mathcal U_sD_By^\eta.
 \end{gathered}
\]
Part (i) gives the spatial estimates for the frames. Their ordinary
time derivatives are estimated from the flow and inverse-flow
equations; as above, one additional spatial derivative of the LDF
is required. For a finite sum of profiles, we also have
\[
 \mathcal U_s\!\left[\sum_h b_hH_h(\lambda y^k)\right]
       =\sum_h(\mathcal U_sb_h)H_h(\lambda y_s^k).
\]
Thus the pushforward preserves the form of the profiles, with the
phase still annihilated by the transported material and magnetic
derivatives. The chart estimates allow us to apply
Lemma~\ref{app:finite-spatial-reserve} in these coordinates.
\section{An MHD Frobenius Theorem}\label{sec:mhd-frobenius}

We construct volume-preserving coordinates that straighten the
material and magnetic derivatives. The induction equation implies
that these derivatives commute. We may therefore straighten the
magnetic field at the initial time and extend the coordinates by the
velocity flow. The volume normalization will leave the magnetic
coordinate unchanged. We then estimate the charts and use them to
construct a transported partition of unity.

\subsection{Adapted frame conventions}

Let $D_t=\partial_t+v\cn$ and $D_B=B\cn$ denote differentiation
of Cartesian components, and let
$\mathcal D_t=\partial_t+\mathcal L_v$ and $\mathcal L_B$ denote
the corresponding Lie derivatives of tensors. They agree on scalars.
For vectors, the Lie derivative subtracts the field gradient acting
on the vector; for one-forms, it adds the transpose gradient acting
on the one-form. We write $\DD$ for spatial derivatives of Cartesian
components.

All geometric operations below are taken in space at fixed time,
with the Lie derivative and curl conventions of
Section~\ref{ssec:path-conventions}. A one-form is identified with its
Cartesian coefficients when taking a cross product, and a two-form with
its vector representative under the fixed Euclidean volume form.

For a volume-preserving chart $\Psi$ from physical to chart coordinates,
write $M=\DD\Psi$. At fixed time,
$\DD\Psi^{-1}=M^{-1}\circ\Psi^{-1}$. We use $\Psi^*$ for vector
pullback and $\Psi^{1*}$ for one-form pullback:
\begin{equation}
 \Psi^*Z=M^{-1}(Z\circ\Psi),\qquad
 \Psi^{1*}\omega=M^T(\omega\circ\Psi),\qquad \det M=1.
\end{equation}
Define the pushforwards $\Psi_*$ and $\Psi_{1*}$ by the corresponding
pullbacks under $\Psi^{-1}$. The vector formula follows from the
chain rule for $Z\cn f$, and the one-form formula from preservation
of the vector--covector pairing. For a two-form represented by
$\iota_Z\mathrm{vol}$, volume preservation gives
\[
 \Psi^*(\iota_Z\mathrm{vol})=\iota_{\Psi^*Z}\mathrm{vol}.
\]
Consequently, vectors and their associated two-forms have the same
pullback formula, while one-forms transform by the transpose
differential. In each case we use the Lie derivative of the indicated
tensor type.

An oriented constant orthonormal frame is denoted by
$(k,\nu,\zeta)$, with $k\times\nu=\zeta$. On one chart define
\begin{equation}\label{setup:adapted-frame-convention}
 y^{\eta}=\Psi\cdot\eta,\qquad
 d y^{\eta}=\Psi^{1*}\eta,\qquad
 \frac{\partial}{\partial y^{\eta}}=\Psi^*\eta,
 \qquad \eta\in\{k,\nu,\zeta\}.
\end{equation}
The identity $M M^{-1}=\IId$ shows that the coordinate vectors
and differentials are dual. Volume preservation gives
\begin{equation}
 \frac{\partial}{\partial y^{\zeta}}=d y^{k}\times d y^{\nu},
 \qquad \ddiv\frac{\partial}{\partial y^{\zeta}}=0,
 \qquad \partial_{y^{\zeta}}y^{k}=0.
\end{equation}
The divergence of the cross product vanishes by symmetry of the
Hessians. Taking cyclic permutations shows that every coordinate
vector is divergence free. These vectors commute and are the vectors
used in the rank-one tensors of the construction.

\subsection{Frobenius coordinates}

\begin{proposition}[Magnetic coordinates]
\label{frobenius:coordinates}
Fix $0<c_0<C_0$ and a point $(t_0,x_0)$. There are constants
$R>1$ and $\varepsilon_0>0$, depending only on $c_0,C_0$, with the
following property. Let $0<T\le r_{\rm box}$ and let $v,B$ be
smooth vector fields on
\[
 [t_0-T,t_0+T]\times B_{R r_{\rm box}}(x_0),
\]
where the spatial ball is contained in a coordinate neighborhood.
Assume
\begin{equation}\label{frobenius:induction}
 \partial_tB+[v,B]=0,\qquad
 \ddiv v=\ddiv B=0,\qquad
 |v|+|B|\le C_0,\qquad |B(t_0,x_0)|\ge c_0.
\end{equation}
If, in addition, on the same cylinder,
\begin{equation}\label{frobenius:smallness}
 [v]_1\le L_t,\qquad [B]_1\le L_B,\qquad
 r_{\rm box}L_B+TL_t\le\varepsilon_0,
\end{equation}
then there is a smooth chart $\Psi(t,x)$ on
$[t_0-T,t_0+T]\times B_{r_{\rm box}}(x_0)$, volume preserving and
a diffeomorphism onto its image at each fixed time, such that
\begin{equation}\label{frobenius:identities}
 D_t\Psi=0,\qquad D_B\Psi=c^{-1}e_3,\qquad
 \det\DD\Psi=1,\qquad c=|B(t_0,x_0)|^{-1}.
\end{equation}
For any $O\in SO(3)$ satisfying
$OB(t_0,x_0)=|B(t_0,x_0)|e_3$, the chart may be chosen so that
$\Psi(t_0,x_0)=0$ and $\DD\Psi(t_0,x_0)=O$. It then satisfies
\begin{equation}\label{frobenius:closeness}
 \|\DD\Psi-O\|_0+
 \|(\DD\Psi)^{-1}-O^T\|_0
 \le C(r_{\rm box}L_B+TL_t).
\end{equation}
The conclusion also holds on any fixed enlargement of this
cylinder, provided $R$ is increased and $\varepsilon_0$ decreased
by fixed factors.
\end{proposition}

\begin{remark}[Reference coordinates]
\label{part:normalized-chart}
After the orientation-preserving rigid change
$x\mapsto O(x-x_0)$, we have $B(t_0,0)=c^{-1}e_3$. This change
preserves volume, the induction equation, the straightening identities
and all derivative estimates. We therefore construct the chart with
$x_0=0$, $O=\IId$, $\Psi(t_0,0)=0$ and
$\DD\Psi(t_0,0)=\IId$.
We make this choice separately on each chart. Fixed spatial cutoffs,
expressed in the same coordinates, extend along the velocity flow as
$\beta(\Psi)$. Returning to the original coordinates replaces the
reference differential $\IId$ by $O$.
\end{remark}

\begin{proof}[Proof of Proposition~\ref{frobenius:coordinates}]
We use the coordinates of Remark~\ref{part:normalized-chart}.
By the induction equation, for every scalar $f$,
\begin{equation}
 [D_t,D_B]f=(\partial_tB+[v,B])\cn f=0.
\end{equation}

\emph{1. Initial magnetic coordinates.}
Let $Y_r$ be the flow of $cB(t_0,\cdot)$ and define
\begin{equation}
 H(y)=Y_{y^{3}}(y^{1},y^{2},0),\qquad |y^i|<R_0r_{\rm box}.
\end{equation}
Choose $R_0$ sufficiently large. Since $|cB|\le C_0/c_0$, increasing
$R$ by a fixed factor ensures that all trajectories defining $H$
remain in the original ball. To estimate the differential, use the
identities
\[
 \partial_r\DD Y_r=c\DD B(t_0,Y_r)\DD Y_r,
 \qquad \partial_3H=cB(t_0,H),\qquad \DD H(0)=\IId
\]
The first identity and Gronwall's inequality bound the transverse
derivatives of $H$. In the second identity, subtract the value at
the origin to estimate the difference of the derivative in the third
coordinate from its value there. Together with the last identity,
these estimates give
\begin{equation}\label{frobenius:box-close}
 \|\DD H-\IId\|_0\le Cr_{\rm box}L_B,
 \qquad |H(y)-y|\le Cr_{\rm box}L_B|y|.
\end{equation}
Integrating the differential estimate along the segment between $y$ and
$y'$ and taking $\varepsilon_0$ small gives
\[
 |H(y)-H(y')|\ge(1-C\varepsilon_0)|y-y'|
                 \ge\tfrac12|y-y'|.
\]
Thus $H$ is injective. To prove the required inclusion for its image,
fix $|x|\le(C_0+2)r_{\rm box}$. For $R_0$ sufficiently large, the
map $y\mapsto x-(H(y)-y)$ sends the closed cube into its interior
and has Lipschitz constant at most $C\varepsilon_0<1/2$. The
contraction mapping theorem gives a solution of $H(y)=x$, proving
\[
 B_{(C_0+2)r_{\rm box}}(0)\subset H(\{|y^i|<R_0r_{\rm box}\}).
\]

\emph{2. Volume normalization.}
The autonomous flow $Y_r$ preserves volume and transports its
defining vector field. Hence
\[
 \DD H(y)=\DD Y_{y^3}(y^1,y^2,0)
       [e_1,e_2,cB(t_0,y^1,y^2,0)],\qquad \det\DD Y_r=1.
\]
Consequently
\begin{equation}\label{frobenius:initial-density}
 \det\DD H(y)=h(y^1,y^2),\qquad
 h(y^1,y^2)=cB_3(t_0,y^1,y^2,0).
\end{equation}
Here $h(0)=1$, $\|h-1\|_0\le Cr_{\rm box}L_B$
and $[h]_1\le CL_B$. Set
\begin{equation}
 A(y)=\left(y^1,\int_0^{y^2}h(y^1,z)\,\dd z,y^3\right),
 \qquad \Psi_0=A\circ H^{-1}.
\end{equation}
The first and third components of $A$ are unchanged, and its second
component is strictly increasing in $y^2$ since $h\ge1/2$. Thus $A$
is injective. Direct differentiation gives
\[
 \det\DD A=h,\qquad \DD A\,e_3=e_3,\qquad
 A(0)=0,\qquad \DD A(0)=\IId,
\]
and
\begin{equation}\label{frobenius:normalization-close}
 \|\DD A-\IId\|_0\le Cr_{\rm box}L_B,\qquad
 |A(y)-y|\le Cr_{\rm box}L_B|y|.
\end{equation}
As in Step 1, the contraction mapping theorem shows that the image
of $A$ contains a smaller cube of comparable side length. The map
$\Psi_0$ is volume preserving and straightens the magnetic field,
since at $x=H(y)$,
\[
 \det\DD\Psi_0(x)=\frac{\det\DD A(y)}{\det\DD H(y)}=1,
 \qquad
 \DD\Psi_0(x)B(t_0,x)=c^{-1}\DD A(y)e_3=c^{-1}e_3.
\]
Moreover $\Psi_0(0)=0$ and $\DD\Psi_0(0)=\IId$.

\emph{3. Extension along the velocity flow.}
Let $\Phi_{t\leftarrow t_0}$ be the velocity flow. The displacement bound
\[
 |\Phi_{t\leftarrow t_0}(a)-a|\le C_0T
\]
shows that the inverse images of $B_{r_{\rm box}}(0)$ remain in the
domain of $\Psi_0$. We may therefore define
\begin{equation}\label{frobenius:material-chart}
 \Psi(t,x)=\Psi_0\bigl(\Phi_{t_0\leftarrow t}(x)\bigr).
\end{equation}
Then $D_t\Psi=0$. Since $\ddiv v=0$,
\[
 \det\DD\Phi_{t\leftarrow t_0}=1,
 \qquad \det\DD\Psi=1.
\]
To verify the magnetic identity, differentiate the pullback of $B$:
\[
 \begin{aligned}
 \partial_t\left[(\DD\Phi_{t\leftarrow t_0})^{-1}
                  B(t,\Phi_{t\leftarrow t_0})\right]
 &=(\DD\Phi_{t\leftarrow t_0})^{-1}
          (D_tB-D_Bv)(t,\Phi_{t\leftarrow t_0})=0.
 \end{aligned}
\]
The pullback is constant in time and equals its value at $t_0$.
Using the initial straightening identity and the chain rule in
\eqref{frobenius:material-chart}, we obtain
\[
 D_{B(t)}\Psi(t,x)
 =\bigl(D_{B(t_0)}\Psi_0\bigr)(\Phi_{t_0\leftarrow t}(x))
 =c^{-1}e_3.
\]
Thus all identities in \eqref{frobenius:identities} hold throughout
the cylinder.

Finally the differentiated velocity-flow equation gives
\[
 \|\DD\Phi_{t\leftarrow t_0}-\IId\|_0
 +\|\DD\Phi_{t_0\leftarrow t}-\IId\|_0
 \le CTL_t.
\]
Combining this with \eqref{frobenius:box-close} and
\eqref{frobenius:normalization-close} in the chain rule gives
\[
 \|\DD\Psi-\IId\|_0\le C(r_{\rm box}L_B+TL_t).
\]
The right-hand side is small, so the inverse differential is uniformly
bounded. The identity
\[
 M^{-1}-\IId=M^{-1}(\IId-M)
\]
then gives the same estimate for the inverse differential. Returning
to the original coordinates proves \eqref{frobenius:closeness}.
For an enlarged cylinder, increase the initial spatial domain by a
fixed factor to contain the trajectories and coordinate images, and
decrease the smallness constant by a fixed factor. The same construction
and estimates apply.
\end{proof}

\subsection{Chart derivative estimates}

The chart \eqref{frobenius:material-chart} satisfies the following
derivative estimates at any frequency for which the field hypotheses
hold. Spatial and ordinary time derivatives use a common frequency;
material and magnetic derivatives have separate costs.

\begin{lemma}[Flow derivative bounds]
\label{frobenius:finite-flow}
Let $\Phi^{u}_{s\leftarrow r}$ be the flow defined by
\[
 \partial_s\Phi^{u}_{s\leftarrow r}=u_s\circ \Phi^{u}_{s\leftarrow r},
 \qquad \Phi^{u}_{r\leftarrow r}=\IId,
\]
for a vector field $u_s$ on an enlarged spatial domain; the field
may depend on $s$. Assume that all trajectories under consideration
remain in this domain. Fix an integer $N\ge1$ and a frequency
$\Lambda\ge1$, and suppose that
\[
 [u_s]_{j+1}\le C_j\operatorname{Lip}\Lambda^j,
 \qquad 0\le j\le N,
 \qquad |s-r| \operatorname{Lip} \le C_*.
\]
Then
\begin{equation}
 \begin{aligned}
 \|\DD \Phi^{u}_{s\leftarrow r}-\IId\|_0
   &\le C|s-r| \operatorname{Lip} ,\\
 [\Phi^{u}_{s\leftarrow r}]_n
   &\le C_n|s-r| \operatorname{Lip} \Lambda^{n-1},
       &&2\le n\le N+1.
 \end{aligned}
\end{equation}
The same estimates hold for the inverse map
$\Phi^{u}_{r\leftarrow s}$, with constants depending only on $C_*$
and the constants in the hypotheses.

More generally, let $ L_0=1$ and $ L_j\ge1$ satisfy
$ L_i L_j\le L_{i+j}$ for $i+j\le N$. If the bound for the gradient
at order $j$ has the additional factor $ L_j$, then the estimate
for $\DD^n\Phi^{u}_{s\leftarrow r}$ has the additional factor
$ L_{n-1}$.
\end{lemma}

For example,
$ L_j=A^{(j-N_*)_+}$ has this property for any
$A\ge1$ and integer $N_*\ge0$.

Only $|s-r| \operatorname{Lip} $ is required to be bounded. No
smallness of $|s-r|\Lambda$ or of the vector field itself is assumed;
the latter determines the domain enlargement needed to contain the
trajectories.

\begin{proof}
Write $F_s=\Phi^u_{s\leftarrow r}$. The differential satisfies
\[
 \partial_s\DD F_s=\DD u_s(F_s)\DD F_s.
\]
Gronwall's inequality and the bound on the integral of the spatial
gradient over the interval of the flow parameter give
\[
 [F_s]_1\le e^{C|s-r|\operatorname{Lip}},\qquad
 \|\DD F_s-\IId\|_0\le C|s-r|\operatorname{Lip}.
\]
For $n\ge2$, the differentiated equation is
\begin{equation}
 \partial_s\DD^nF_s
 =\DD u_s(F_s)\DD^nF_s+
 \sum_{\substack{2\le p\le n\\n_1+\cdots+n_p=n\\n_i\ge1}}
 C_{n_1,\ldots,n_p}\DD^pu_s(F_s)
       [\DD^{n_1}F_s,\ldots,\DD^{n_p}F_s].
\end{equation}
Assume the estimates through order $n-1$. Since $p\ge2$, every
$n_i\le n-1$, so each flow derivative on the right is bounded by
the induction hypothesis. The bound on the integral of the spatial
gradient over the flow interval then gives
\[
 \left\|\DD^pu_s(F_s)
       [\DD^{n_1}F_s,\ldots,\DD^{n_p}F_s]\right\|_0
 \le C\operatorname{Lip}\Lambda^{p-1}
             \prod_i\Lambda^{n_i-1}
 =C\operatorname{Lip}\Lambda^{n-1}.
\]
Since $\DD^nF_r=0$, variation of constants gives
\[
 [F_s]_n
 \le C\int_{\min(r,s)}^{\max(r,s)}
          e^{C|s-r|\operatorname{Lip}}
          \operatorname{Lip}\Lambda^{n-1}\,\dd t
 \le C_n|s-r|\operatorname{Lip}\Lambda^{n-1}.
\]
This proves the estimate through $n=N+1$, requiring derivatives of
$u_s$ only through order $N+1$. Reversing the flow interval proves
the estimate for the inverse. With the additional loss factors, the
factor in each summand is bounded by
\[
 L_{p-1}\prod_iL_{n_i-1}\le
 L_{p-1+\sum_i(n_i-1)}=L_{n-1}.
\]
Thus the same induction proves the estimate with loss factors. The
example stated after the lemma follows from
\[
 (i-N_*)_++(j-N_*)_+\le(i+j-N_*)_+.
\]
\end{proof}

We apply the lemma to the velocity and magnetic flows separately.
The integrals of their spatial gradients over the respective flow
intervals are small in our applications, although boundedness is
sufficient here. The choice $L_j=A^{(j-N_*)_+}$ with
$A=(\ell\Lambda)^{-1}$ gives the sharp estimates below $N_*$ and
the lossy estimates above it for both flow maps. In particular, no
condition $\tau_a\ell^{-1}\ll1$ is needed.

\begin{proposition}[Chart derivative bounds]
\label{frobenius:estimates}
Let $\Psi$ be the chart constructed in
Proposition~\ref{frobenius:coordinates}, under the hypotheses of that
proposition. Fix an integer $N\ge1$ and a frequency $\Lambda\ge1$
such that $L_t+L_B\le C\Lambda$.
\begin{enumerate}
\item[(i)] \emph{Spatial estimates.}
Assume that the fields satisfy
\begin{equation}
 [v]_{r+1}\le C_rL_t\Lambda^r,\qquad
 [B]_{r+1}\le C_rL_B\Lambda^r,
 \qquad 0\le r\le N.
\end{equation}
on the enlarged cylinder used to construct the chart. Put
$M=\DD\Psi$. Then each matrix field
\begin{equation}\label{frobenius:matrix-family}
 f\in\{M,M^{-1}\},
\end{equation}
satisfies $\|f\|_r\le C_r\Lambda^r$ for $0\le r\le N$.

\item[(ii)] \emph{Ordinary time estimates.}
Fix an integer $k_{\max}\ge1$, independently of the spatial derivative
order. In addition to the hypotheses of (i), assume that
\begin{equation}\label{frobenius:time-input}
 \|\partial_t^a v\|_r\le C_{r,a}\Lambda^{r+a},
 \qquad 0\le a\le  k_{\max}-1,\quad r+a\le N.
\end{equation}
Then every $f$ in \eqref{frobenius:matrix-family} satisfies
\begin{equation}\label{frobenius:ordinary-output}
 \|\partial_t^k f\|_r\le C_{r,k}\Lambda^{r+k},
 \qquad 0\le k\le  k_{\max},\quad r+k\le N.
\end{equation}

These spatial and ordinary time derivative estimates also hold for
the differential of the inverse velocity flow $\Phi_{t_0\leftarrow t}$
and its inverse matrix. The chart and the inverse velocity flow are
volume preserving. The maps themselves satisfy, for
$F=\Psi$, $\Psi(t,\cdot)^{-1}$ or $\Phi_{t_0\leftarrow t}$,
\begin{equation}
 [\partial_t^k F]_r\le C_{r,k}\Lambda^{r+k-1},
 \qquad k\le  k_{\max},\quad 1\le r+k\le N+1.
\end{equation}
The inverse charts are considered on the enlarged coordinate sets
where they are defined. These assertions concern ordinary derivatives
of the three maps. The magnetic identities used in (iii) apply only
to the straightening chart $\Psi$.

If also $0<\alpha<1$ and $\ell\Lambda\le1$, then
\eqref{frobenius:ordinary-output} holds in the norm
$\|\partial_t^kf\|_{r+\alpha}$ with the additional factor
$\ell^{-\alpha}$ whenever $r+k+1\le N$.

\item[(iii)] \emph{Material and magnetic derivatives.}
Assume the hypotheses of (i), and let $0<\alpha<1$ with
$\ell\Lambda\le1$. Fix $j_{\max},N_{\rm char}\ge1$ satisfying
$N_{\rm char}+1\le N$, and choose independent positive derivative
costs $\mathrm a_t,\mathrm a_B$. Suppose that
\begin{equation}\label{frobenius:mixed-input}
 \begin{aligned}
 \|D_t^pD_B^j\DD v\|_r
 &\le C_{r,p,j}\Lambda^r\mathrm a_t^{p+1}\mathrm a_B^j,&
 \|D_t^pD_B^j\DD v\|_{r+\alpha}
 &\le C_{r,p,j}\ell^{-\alpha}\Lambda^r\mathrm a_t^{p+1}\mathrm a_B^j,\\
 \|D_t^pD_B^j\DD B\|_r
 &\le C_{r,p,j}\Lambda^r\mathrm a_t^p\mathrm a_B^{j+1},&
 \|D_t^pD_B^j\DD B\|_{r+\alpha}
 &\le C_{r,p,j}\ell^{-\alpha}\Lambda^r\mathrm a_t^p\mathrm a_B^{j+1},
 \end{aligned}
\end{equation}
for $p+j\le j_{\max}-1$ and $r+p+j+1\le N_{\rm char}$. Then every
$f$ in \eqref{frobenius:matrix-family} satisfies
\begin{equation}\label{frobenius:mixed-output}
 \begin{gathered}
 \|D_t^kD_B^m f\|_r\le C_{r,k,m}\Lambda^r\mathrm a_t^k\mathrm a_B^m,\qquad
 \|D_t^kD_B^m f\|_{r+\alpha}\le C_{r,k,m}\ell^{-\alpha}\Lambda^r\mathrm a_t^k\mathrm a_B^m,\\
 k+m\le j_{\max},\qquad r+k+m\le N_{\rm char}.
 \end{gathered}
\end{equation}
\end{enumerate}
The matrix estimates in (i)--(iii) hold also for fixed polynomial expressions in
the matrices, and for smooth functions of their entries when the
matrices range over a fixed compact subset of the domain of definition.
\end{proposition}

In (ii), the time derivatives of the velocity are required only
through order $k_{\max}-1$ to estimate the chart through order
$k_{\max}$. No ordinary time derivatives of $B$ are needed. The
estimates for the maps are independent of the coordinate origin,
since they contain only derivatives of order at least one.

\begin{proof}
\emph{1. Spatial derivatives.}
Lemma~\ref{frobenius:finite-flow} and \eqref{frobenius:smallness}
give
\[
 [Y_r]_n+[\Phi_{t\leftarrow t_0}]_n
 \le C_n\Lambda^{n-1},\qquad 1\le n\le N+1.
\]
The derivatives of $H$ in the first two coordinates are restrictions
of derivatives of the magnetic flow, so they satisfy the same bound.
For a derivative involving $\partial_3$, use
$\partial_3H=cB(t_0,H)$ and apply the remaining $n-1$ derivatives
to this composition. The chain rule uses derivatives of $B(t_0)$
only through order $n-1$, and gives $[H]_n\le C_n\Lambda^{n-1}$.

By \eqref{frobenius:initial-density},
\begin{equation}
 [h]_m\le C_mL_B\Lambda^{m-1},
 \qquad 1\le m\le N+1.
\end{equation}
For the second component of $A$,
\[
 A_2(y)=\int_0^{y^2}h(y^1,z)\,\dd z.
\]
With $a+b=n$, differentiation gives
\[
 \partial_1^a\partial_2^b A_2=
 \begin{cases}
 \partial_1^a\partial_2^{b-1}h,&b\ge1,\\
 \displaystyle\int_0^{y^2}\partial_1^n h(y^1,z)\,\dd z,&b=0.
 \end{cases}
\]
The other components are affine, and higher derivatives involving
$y^3$ vanish. Estimating the two cases above gives
\begin{equation}
 [A]_n
 \le C_n\bigl(L_B\Lambda^{n-2}
       +r_{\rm box}L_B\Lambda^{n-1}\bigr)
 \le C_n\Lambda^{n-1},\qquad 2\le n\le N+1.
\end{equation}
The last inequality follows from $L_B\le C\Lambda$ for the first
term and $r_{\rm box}L_B\le\varepsilon_0$ for the second. At
$n=N+1$ we require derivatives of $B$ through order $N+1$, as assumed.

For a map $F$ with these bounds and bounded inverse differential,
differentiating $F\circ F^{-1}=\IId$ gives, for $n\ge2$,
\[
 \begin{aligned}
 (\DD F\circ F^{-1})\DD^nF^{-1}
 &=-\sum_{\substack{2\le p\le n\\n_1+\cdots+n_p=n\\n_i\ge1}}
 C_{n_1,\ldots,n_p}(\DD^pF\circ F^{-1})
      [\DD^{n_1}F^{-1},\ldots,\DD^{n_p}F^{-1}],\\
 \Lambda^{p-1}\prod_i\Lambda^{n_i-1}&=\Lambda^{n-1}.
 \end{aligned}
\]
All derivatives of the inverse on the right have order less than
$n$. Multiplying by the bounded inverse differential therefore proves
the bound for $F^{-1}$ by induction. The displayed equality of powers
of $\Lambda$ gives the same bound for compositions. Apply these
observations to $H$, $A$ and the velocity flow, and to their inverses.
We obtain the estimate for
$\Psi= A\circ H^{-1}\circ\Phi_{t_0\leftarrow t}$, and hence the
spatial estimates for $M$, $M^{-1}$ and all the maps in the statement.

\emph{2. Ordinary time derivatives.}
Differentiating $D_t\Psi=0$ in space yields
\begin{equation}
 \partial_tM=-v\cn M-M\DD v.
\end{equation}
For $k\ge1$ and spatial derivative order $r$, Leibniz' rule gives
\[
 \begin{aligned}
 \|\partial_t^kM\|_r
 \le C_{r,k}\sum_{\substack{a+b=k-1\\r_1+r_2\le r}}
 \bigl(&\|\partial_t^av\|_{r_1}
                    \|\partial_t^bM\|_{r_2+1}
       +\|\partial_t^av\|_{r_1+1}
                    \|\partial_t^bM\|_{r_2}\bigr)
 \le C_{r,k}\Lambda^{r+k}.
 \end{aligned}
\]
To justify the last inequality, induct on $k$. Every derivative of
$M$ on the right has time order $b<k$, while each velocity factor
has time order $a\le k-1$. The powers of the frequency add to at
most $r+k$, and \eqref{frobenius:time-input} applies to the velocity
derivatives when $r+k\le N$. Differentiating
$M^{-1}M=\IId$ gives the same estimate for $M^{-1}$. The inverse
velocity flow satisfies the same transport equation with affine
initial data, so its differential satisfies these estimates as well.

For the maps themselves, use
\[
 \partial_t\Psi=-\DD\Psi\,v,\qquad
 \partial_t\Phi_{t_0\leftarrow t}=-\DD(\Phi_{t_0\leftarrow t})v,
 \qquad \partial_t\Psi^{-1}=v\circ\Psi^{-1}.
\]
The last identity follows by differentiating
$\Psi\circ\Psi^{-1}=\IId$ in time. Apply the product rule to the
first two equations and the chain rule to the third. Induction gives
\[
 [\partial_t^kF]_r\le C_{r,k}\Lambda^{r+k-1},
 \qquad k\le k_{\max},\quad 1\le r+k\le N+1,
\]
for each map $F$ in the statement. The induction uses velocity
derivatives with time order at most $k-1$ and sum of spatial and time
derivative orders at most $r+k-1\le N$. The initial normalization
is fixed at $t_0$ and therefore contributes no time derivatives of
$B$.

For the H\"older estimate, interpolate consecutive spatial derivative
bounds using, on the enlarged convex balls,
\[
 [g]_\alpha\le C\|g\|_0^{1-\alpha}[g]_1^\alpha,
 \qquad \ell^\alpha\Lambda^\alpha\le1.
\]
Apply this to every spatial derivative of $\partial_t^kf$ when
$r+k+1\le N$.

\emph{3. Material and magnetic derivatives.}
Spatial differentiation of \eqref{frobenius:identities} gives
\begin{equation}\label{frobenius:matrix-mixeds}
 D_tM=-M\DD v,\qquad D_BM=-M\DD B.
\end{equation}
When $k+m$ is positive, repeated differentiation expresses
$D_t^kD_B^mM$ as a sum of ordered matrix products of the form
\[
 M\prod_{i=1}^p D_t^{k_i}D_B^{m_i}\DD w_i,
 \qquad w_i\in\{v,B\},
\]
where the indices satisfy
\begin{equation}
 \sum_i\bigl(k_i+\mathbf1_{w_i=v}\bigr)=k,
 \qquad
 \sum_i\bigl(m_i+\mathbf1_{w_i=B}\bigr)=m.
\end{equation}
Indeed, applying $D_t$ either introduces a factor $\DD v$ or
increases one $k_i$, while applying $D_B$ either introduces
$\DD B$ or increases one $m_i$. Since $[D_t,D_B]=0$, the derivatives
on every coefficient can be put in the displayed order. The spatial
product rule gives terms with $r_0+\cdots+r_p=r$, each bounded by
\[
 C\Lambda^{r_0}\prod_i
 \left(\Lambda^{r_i}
 \mathrm a_t^{k_i+\mathbf1_{w_i=v}}
 \mathrm a_B^{m_i+\mathbf1_{w_i=B}}\right)
 =C\Lambda^r\mathrm a_t^k\mathrm a_B^m.
\]
At most $k+m-1$ material and magnetic derivatives combined act on
any one field gradient. The remaining restriction in
\eqref{frobenius:mixed-input} holds because
\[
 r_i+k_i+m_i+1\le r+k+m.
\]
For the H\"older estimate, place one factor in $C^{r_i+\alpha}$
and apply the product inequality. This proves
\eqref{frobenius:mixed-output} whenever a transport derivative is
present. With no transport derivatives, use interpolation and the
spatial estimate at one additional order.

Derivatives of the inverse matrix are products of bounded inverse
matrices and derivatives of $M$, with no increase in derivative
order. This proves the assertion for $M^{-1}$. The product and chain
rules then give the matrix estimates in (i)--(iii) for polynomial expressions and
smooth matrix functions on the compact subset specified in the
statement.
\end{proof}

\subsection{An isotropic Lagrangian partition}

\begin{proposition}[Lagrangian partition]
\label{frobenius:partition}
Let $v,B$ be smooth vector fields on
$[t_0-T,t_0+T]\times\mathbb T^3$ satisfying
\eqref{frobenius:induction}, with $|v|+|B|\le C_0$ and
$|B|\ge c_0>0$. Suppose that $0<T\le r$, with $r$ sufficiently
small, and that
\[
 r[B]_1+T[v]_1\le\varepsilon_0.
\]
Then there is a smooth quadratic partition of unity
$\sum_J\chi_J^2=1$ transported by $v$, with volume-preserving magnetic
charts $\Psi_J$ defined on neighborhoods of the respective supports.
Each support has diameter at most $Cr$ in physical and chart
coordinates. The supports have overlap at most eight and admit eight
colors such that the cut off regions of fixed relative width are
disjoint within each color. The partition satisfies
\begin{equation}\label{frobenius:cutoffs}
 D_t\chi_J=0,\qquad
 \|D_B^m\chi_J\|_0\le C_mr^{-m}.
\end{equation}
If the chart and velocity satisfy the estimates in parts (i)--(ii) of
Proposition~\ref{frobenius:estimates} at a frequency
$\Lambda\ge r^{-1}$, then
\begin{equation}\label{frobenius:cutoff-ordinary}
 \|\partial_t^kD_B^m\chi_J\|_j
 \le C_{j,k,m}r^{-m}\Lambda^{j+k},
 \qquad k\le k_{\max},\quad j+k\le N.
\end{equation}
The number $m$ of magnetic derivatives may be fixed arbitrarily,
without increasing the required derivative orders of the fields.
With one additional spatial derivative, the estimate also holds in
the H\"older norm with the factor $\Lambda^\alpha$. The constants
depend only on $c_0,C_0$ and the fixed derivative orders.
\end{proposition}

\begin{proof}
\emph{1. Cover and cut off regions.}
Choose an integer $N_{\rm mesh}$ such that
\[
 d=2\pi/N_{\rm mesh}\asymp r,
\]
where $d/r$ lies between two sufficiently small fixed positive
constants. Take the eight translates of the periodic grid
$d\mathbb Z^3$ by vectors in $\{0,d/2\}^3$, and denote their centers
by $x_J$ and their colors by elements of $\{0,1\}^3$. At each center,
choose a nonnegative bump supported in the cube of side $3d/5$,
uniformly positive on the concentric cube of side $11d/20$, with
derivatives of order $j$ bounded by $C_jd^{-j}$.

Choose $O_JB(t_0,x_J)=|B(t_0,x_J)|e_3$, and express the bump
in these rigid coordinates as $\beta_J$. Construct $\Psi_J$ by
Proposition~\ref{frobenius:coordinates} on a ball of radius $Cr$,
on $[t_0-T,t_0+T]$. Its initial differential is within
$C\varepsilon_0$ of $O_J$, and hence
\begin{equation}
 |\Psi_J(t_0,x)-O_J(x-x_J)|\le C\varepsilon_0r
\end{equation}
on a neighborhood of the bump support. Decrease $\varepsilon_0$ so that
this error is at most $d/100$. Then $\beta_J(\Psi_J(t_0,x))$,
extended by zero, is supported in the concentric physical cube
of side $2d/3$ and is uniformly positive on the cube of side
$21d/40$. The latter cubes cover the torus, since their centers
form the grid $(d/2)\mathbb Z^3$ and $21d/40>d/2$. Thus
\begin{equation}\label{part:original-partition}
 c\le\sum_J\beta_J(\Psi_J(t_0,x))^2\le C.
\end{equation}
Enlarging the support cubes to side $3d/4$ leaves a cut off region of
width $d/24$. At most eight enlarged cubes overlap, because in each
coordinate an interval of length $3d/4<d$ contains at most two centers
from the union of the translated grids. Centers of one color have
spacing $d$, so the enlarged cubes of that color are disjoint.

Transport the supports and cut off regions by the velocity flow.
The bound $T[v]_1\le\varepsilon_0$ controls the differential and its
inverse, so the diameters and relative widths remain comparable to
their initial values. The flow is a common diffeomorphism and hence
preserves overlap and disjointness. Since $T\le r$ and $|v|\le C_0$,
the transported sets remain in the enlarged chart domains. Define
\begin{equation}
 \chi_J(t,x)=\frac{\beta_J(\Psi_J(t,x))}
 {\bigl(\sum_K\beta_K(\Psi_K(t,x))^2\bigr)^{1/2}}.
\end{equation}
Since $D_t\Psi_J=0$, every summand in the denominator is constant
along the velocity flow. The bounds \eqref{part:original-partition}
therefore hold throughout the time interval. The definition gives
$D_t\chi_J=0$ and $\sum_J\chi_J^2=1$. Each support is contained in
a physical ball of radius $Cr$, and the bound on the chart differential
gives the same diameter estimate in chart coordinates.

\emph{2. Derivatives.}
With $c_J=|B(t_0,x_J)|^{-1}$, the exact straightening identity gives
\begin{equation}\label{part:cutoff-exact-magnetic}
 D_B^m[\beta_J(\Psi_J)]
 =c_J^{-m}(\partial_3^m\beta_J)(\Psi_J).
\end{equation}
Apply the magnetic derivatives to the quotient defining the partition.
The product and chain rules give products of the preceding quantities
and negative powers of the sum in the denominator. Its lower bound
and the finite overlap give $\|D_B^m\chi_J\|_0\le C_mr^{-m}$,
proving \eqref{frobenius:cutoffs}. After $j$ spatial derivatives,
the chart matrices are differentiated at most $j-1$ times. Their
bounds and those of the bumps give a factor at most $\Lambda^j$,
since $r^{-1}\le\Lambda$. By
\eqref{part:cutoff-exact-magnetic}, all magnetic derivatives act on
the fixed bumps, so increasing $m$ requires no additional regularity
of either field. Finally, $[D_t,D_B]=0$ implies
\[
 \partial_t(D_B^m\chi_J)=-v\cn D_B^m\chi_J.
\]
Differentiate this equation successively and use the spatial bounds
to obtain \eqref{frobenius:cutoff-ordinary}. The velocity factors
require at most $k-1$ time derivatives, and the sum of the spatial
and time derivative orders on the chart and velocity is at most
$j+k$. One further spatial derivative and interpolation give the
H\"older factor $\Lambda^\alpha$. The cut off regions ensure that
the functions vanish near the chart boundaries, so their zero
extensions are smooth before differentiation.
\end{proof}

\subsection{Frame and phase transport}

\begin{corollary}[Lie invariance of the chart data]
Let $\Psi$ be the chart of Proposition~\ref{frobenius:coordinates},
with normalization $c$. For a constant covector $\nu$ and a constant
vector $\zeta$, set $y^{\nu}=\Psi\cdot\nu$. Then
\begin{equation}
 d y^{\nu}=\Psi^{1*}\nu=M^T\nu,
 \qquad \Psi^*\zeta=M^{-1}\zeta.
\end{equation}
These fields satisfy
\begin{equation}\label{frobenius:natural-frames}
 (\partial_t+\mathcal L_v)d y^{\nu}
 =\mathcal L_B d y^{\nu}=0,
 \qquad
 (\partial_t+\mathcal L_v)(\Psi^*\zeta)
 =\mathcal L_B(\Psi^*\zeta)=0.
\end{equation}
In particular, the dual coordinate frame
\eqref{setup:adapted-frame-convention} satisfies these identities.
If $k\cdot e_3=0$, the coordinate $y^k=\Psi\cdot k$ satisfies
$D_ty^k=D_By^k=0$. Let $\varphi$ be a smooth function of one variable,
let $\lambda>0$, and let $(k,\nu,\zeta)$ be an oriented orthonormal
frame with $k\cdot e_3=0$. Then
\begin{equation}\label{frobenius:curl-profile}
 \curl\!\left[\lambda^{-1}\varphi(\lambda y^k)d y^{\nu}\right]
   =\varphi'(\lambda y^k)\frac{\partial}{\partial y^{\zeta}}
\end{equation}
defines a divergence-free field annihilated by both
$\mathcal D_t$ and $\mathcal L_B$. The estimates of
Proposition~\ref{frobenius:estimates} hold for $d y^{\nu}$,
$\Psi^*\zeta$, their tensor products, and smooth functions of the
frame matrices on fixed compact subsets of their domains in the
space of invertible matrices.
\end{corollary}

\begin{proof}
Lie differentiation commutes with the exterior derivative. Thus
\[
 \begin{aligned}
 (\partial_t+\mathcal L_v)d y^\nu
     &=d(D_ty^\nu)=0,\\
 \mathcal L_Bd y^\nu
     &=d(D_By^\nu)=d(c^{-1}e_3\cdot\nu)=0.
 \end{aligned}
\]
Applying these identities to the constant covectors and using the
product rule for wedge products proves the same invariance for the
pullback of every constant two-form. To obtain the vector identities,
use volume preservation: for $u=v,B$,
\[
 \mathcal L_u(\iota_Z\mathrm{vol})
       =\iota_{[u,Z]}\mathrm{vol},
 \qquad \mathcal L_u\mathrm{vol}=0.
\]
For the material Lie derivative, the same identity contains the
additional term $\partial_tZ$. The identification of two-forms with
vectors therefore proves the vector identities in
\eqref{frobenius:natural-frames}.

For $k\cdot e_3=0$, \eqref{frobenius:identities} gives
$D_ty^k=D_By^k=0$. Since $d^2y^\nu=0$,
\[
 d\bigl[\lambda^{-1}\varphi(\lambda y^k)d y^\nu\bigr]
       =\varphi'(\lambda y^k)d y^k\wedge d y^\nu.
\]
Representing two-forms by vectors, we have
\[
 (M^Tk)\times(M^T\nu)
   =\det(M)M^{-1}(k\times\nu)=\frac{\partial}{\partial y^\zeta},
\]
which proves \eqref{frobenius:curl-profile}. As a curl, this field is
divergence free. Its scalar factor is constant along the velocity and
magnetic flows, while its coordinate vector is invariant under the
corresponding Lie derivatives. The product rule gives
\[
 \mathcal D_t[\varphi'(\lambda y^k)\frac{\partial}{\partial y^\zeta}]=0,
 \qquad \mathcal L_B[\varphi'(\lambda y^k)\frac{\partial}{\partial y^\zeta}]=0.
\]
The derivative estimates follow from the last statement of
Proposition~\ref{frobenius:estimates}.
\end{proof}
\section{Smoothing and Transport Estimates}
\label{sec:analytic-estimates}

This section contains the smoothing and transport estimates used in the
construction. The estimates distinguish between material derivatives
$D_t$ and magnetic derivatives $D_B$, and between the corresponding
Lie derivatives $\mathcal D_t$ and
$\mathcal L_B$ of tensors. We prove the commutator estimates for
arbitrary fixed tensor types. We also consider smoothing along the
two commuting flows, with an independent smoothing length for each
flow.

\subsection{H\"older norms and products}
\label{ssec:holder-calculus}

For an integer $r\ge0$ and $0<\alpha<1$, we define the spatial
H\"older norms and seminorms by
\begin{equation}
 \begin{aligned}
 [f]_r&=\max_{|\boldsymbol\gamma|=r}
                  \|\partial^{\boldsymbol\gamma}f\|_0,
 &\|f\|_r&=\sum_{i=0}^r[f]_i,\\
 [f]_{r+\alpha}&=\max_{|\boldsymbol\gamma|=r}\sup_{x\ne y}
 \frac{|\partial^{\boldsymbol\gamma}f(x)
       -\partial^{\boldsymbol\gamma}f(y)|}
      {\operatorname{dist}_{\mathbb T^3}(x,y)^\alpha},
 &\|f\|_{r+\alpha}&=\|f\|_r+[f]_{r+\alpha}.
 \end{aligned}
\end{equation}
Here $\|f\|_0=\sup_x|f(x)|$; when a time interval is under
consideration, all these quantities include its supremum in time.
The same notation applies componentwise to vectors and tensors.
Unless a spatial derivative is needed to specify the tensor type or
the order of noncommuting operations, we indicate its order by the
norm or seminorm subscript. The spatial multiindex
$\boldsymbol\gamma$ is distinct from the pair of orders $(k,m)$ of
material and magnetic derivatives, respectively. Norms with an explicit space--time subscript are
defined where they occur. Constants in $\lesssim$ may depend on
fixed derivative orders, profiles, geometric separation constants, and the stated
background constants, but are independent of $q$ and the large
parameter $a$. We indicate these dependencies when needed.
The elementary tame product estimate is
\begin{equation}\label{setup:product}
 \|fg\|_{r+\alpha}\le C_r
 \bigl(\|f\|_{r+\alpha}\|g\|_0+
       \|f\|_0\|g\|_{r+\alpha}\bigr).
\end{equation}
Indeed, Leibniz' rule and interpolation bound each differentiated
product by the two terms on the right of \eqref{setup:product}.
For the H\"older seminorm we use the identity
\[
 (fg)(x)-(fg)(y)=f(x)(g(x)-g(y))+g(y)(f(x)-f(y))
\]
on each differentiated product. The chain rule gives the analogous
composition estimate for functions taking values in a fixed compact set.

Only one factor in each product needs to be estimated in a H\"older
norm. In particular, suppose that integer norms are bounded and that
the H\"older bounds contain additional factors $h_1,h_2\ge1$, as in
Section~\ref{ssec:fixed-scales}. The product then has the product of
the amplitudes and the additional factor $\max\{h_1,h_2\}$.
We will also use the boundedness of order-zero singular integrals from
$C^{r+\alpha}$ to $C^{r+\alpha}$; the H\"older factor in the estimate
of the function is unchanged by this bound. For a slow function the
factor is $\ell^{-\alpha}$. In the lower-order operator estimates,
it is absorbed by the smallness of the background gradients. When
such an estimate is used for a source term, its H\"older factor is
included in the subsequent estimates.

\subsection{Spatial Fourier multipliers}
The inverse Laplacian annihilates the zero Fourier mode. Define
\begin{equation}\label{setup:hodge}
 \mathbb P=\IId-\nabla\Delta^{-1}\ddiv,
 \qquad \curl^{-1}=-\Delta^{-1}\curl,
 \qquad \mathscr T=\curl^{-1}\mathbb P\ddiv.
\end{equation}
Here $\curl\curl^{-1}=\mathbb P$ on mean-zero vector fields; constant
vectors are handled separately. The symmetric trace-free inverse
divergence $\mathcal R$ is defined in \eqref{eq:inverse-divergence};
Lemma~\ref{lem:inverse-divergence} proves its symmetry, trace condition,
and right-inverse identity.
We use only symmetry and the right-inverse identity in what follows.
The complete stress also contains explicit quadratic tensors and need
not be trace free. In particular,
\[
 \curl\mathscr T S=\mathbb P\ddiv S,\qquad
 \ddiv\mathcal R\curl\omega=\curl\omega.
\]
The operators $\mathbb P$, $\mathscr T$, and $\mathcal R\curl$ have
order zero, while $\mathcal R$ and $\curl^{-1}$ have order minus one.
The cancellation of their periodic kernels gives boundedness at the
stated $C^{r,\alpha}$ orders. We use these H\"older estimates
throughout; general order-zero singular integrals need not be bounded
on $L^\infty$.

\subsection{Transport estimates}
\label{ssec:component-transport-estimates}

We first prove estimates for equations acting on individual Cartesian
components of vectors and tensors. We state the forward estimates for
$t\ge t_0$; time reversal gives the corresponding estimates before
$t_0$.

\begin{lemma}[H\"older transport estimates]
Let $0<\alpha<1$ and let smooth fields solve
$(\partial_t+a\cn)u=f$, $u(t_0)=u_0$, on a compact interval.
Set
\[
 V(t)=\int_{t_0}^t\|\DD a(s)\|_\alpha\,\dd s.
\]
Then
\begin{equation}\label{high:component-transport-low}
 \|u(t)\|_\alpha
 \le C e^{CV(t)}
 \left(\|u_0\|_\alpha+
              \int_{t_0}^t\|f(s)\|_\alpha\,\dd s\right).
\end{equation}
For every integer $r\ge1$,
\begin{equation}\label{high:component-transport-tame}
 \begin{aligned}
 \|u(t)\|_{r+\alpha}
 \le C_r e^{C_rV(t)}\bigg(&\|u_0\|_{r+\alpha}
       +\int_{t_0}^t\|f(s)\|_{r+\alpha}\,\dd s\\
 &+\int_{t_0}^t
       \|\DD a(s)\|_{r-1+\alpha}\|\DD u(s)\|_\alpha\,\dd s\bigg).
 \end{aligned}
\end{equation}
Using the $C^{r+\alpha}$ bound for the gradient of $a$, one can
replace the last integrand by
$\|\DD a(s)\|_{r+\alpha}\|u(s)\|_\alpha$ in the estimate.
In both versions, the exponential depends only on $\|\DD a\|_\alpha$.
\end{lemma}
\begin{proof}
The flow of $a$ satisfies
\[
 \operatorname{Lip}X_{t\leftarrow s}\le
 \exp(C\int_s^t[a(\sigma)]_1\,\dd\sigma),
 \qquad s\le t,
\]
with the same bound for the inverse flow. The characteristic formula
therefore bounds the supremum norm of the solution. Applying it at two
points and using the Lipschitz bounds for the flow and its inverse gives
\eqref{high:component-transport-low}. For $|\gamma|\le r$,
differentiation of the equation yields
\[
 (\partial_t+a\cn)\partial^\gamma u
 =\partial^\gamma f
  -\sum_{0<\beta\le\gamma}\binom\gamma\beta
       (\partial^\beta a)\cn\partial^{\gamma-\beta}u,
\]
and the product and interpolation estimates give
\[
 \sum_{|\gamma|\le r}
   \|[\partial^\gamma,a\cn]u\|_\alpha
 \le C_r\bigl(
   \|\DD a\|_\alpha\|u\|_{r+\alpha}
   +\|\DD a\|_{r-1+\alpha}\|\DD u\|_\alpha\bigr).
\]
We apply the characteristic estimate to the differentiated equations
and sum over the multiindices. Gronwall's inequality absorbs the term
$\|\DD a\|_\alpha\|u\|_{r+\alpha}$, whereas the second term in
the commutator estimate is integrated as a source. This proves
\eqref{high:component-transport-tame}. To obtain the alternative
estimate, interpolate between $\|\DD a\|_\alpha$ and
$\|\DD a\|_{r+\alpha}$ for the coefficient and between
$\|u\|_\alpha$ and $\|u\|_{r+\alpha}$ for the solution.
When $|\beta|=b$, the higher norms have exponents $(b-1)/r$ and
$(r-b+1)/r$. Their sum is one. Young's inequality therefore gives a
term absorbed by the same Gronwall estimate and the source term
$\|\DD a\|_{r+\alpha}\|u\|_\alpha$.
\end{proof}

\begin{corollary}[Coupled transport]
\label{high:projected-transport}
Let $a_i,c_{ij\nu},f_i$, $1\le i\le m$, be smooth on $I\times\mathbb T^3$, where $I$ is a compact interval containing
$t_0$, and let the finite family of fixed spatial linear operators $T_\nu$ be bounded on $C^{r,\alpha}$ for every $r\ge0$.
For prescribed smooth initial data at $t_0$, the system
\begin{equation}\label{high:coupled-transport-equation}
 (\partial_t+a_i\cn)u_i
       =f_i+\sum_{j,\nu}T_\nu(c_{ij\nu}u_j)
\end{equation}
has a unique smooth solution on $I$. The coefficients $c_{ij\nu}$
may be matrices. If the initial data vanish and the source vanishes
on an interval starting at $t_0$, the solution vanishes throughout
that interval.
Set
\[
 \begin{gathered}
 U_r=\sum_i\|u_i\|_{r+\alpha},\qquad
 F_r=\sum_i\|f_i\|_{r+\alpha},\qquad
 H_r=\sum_i\|\DD a_i\|_{r+\alpha}
               +\sum_{i,j,\nu}\|c_{ij\nu}\|_{r+\alpha},\\
 Q(t)=\int_{t_0}^t H_0(s)\,\dd s.
 \end{gathered}
\]
Then
\begin{align}
 U_0(t)&\le C e^{CQ(t)}
       \left(U_0(t_0)+\int_{t_0}^tF_0(s)\,\dd s\right),
       \label{high:coupled-transport-low}\\
 U_r(t)&\le C_r e^{C_rQ(t)}
       \left(U_r(t_0)+\int_{t_0}^t
             [F_r(s)+H_r(s)U_0(s)]\,\dd s\right),\qquad r\ge1.
       \label{high:coupled-transport-high}
\end{align}
The constants in \eqref{high:coupled-transport-low} and
\eqref{high:coupled-transport-high} depend only on the derivative
order, $\alpha$, the size of the finite system, and the norms of the
fixed operators $T_\nu$ on the H\"older spaces used in the estimate.
The dependence on the coefficients is given by $H_r$ and $Q$.
The estimates use the H\"older norms of $a_i$ of order $r+1+\alpha$
and of $c_{ij\nu},u_i,f_i$ of order $r+\alpha$ only.
\end{corollary}
\begin{proof}
Let $S_{t\leftarrow s}$ be the transport propagator along the respective flows
of $a_i$, and write $\mathcal C(s)u$ for the coupling in
\eqref{high:coupled-transport-equation}. Duhamel's formula is
\begin{equation}
 u(t)=S_{t\leftarrow t_0}u(t_0)+\int_{t_0}^t
                    S_{t\leftarrow s}[f(s)+\mathcal C(s)u(s)]\,\dd s.
\end{equation}
On each $C^{r,\alpha}$, both the propagator and the coupling have
operator norms bounded uniformly on $I$. We iterate Duhamel's formula,
starting from the solution with zero coupling. The correction at the
next iteration is an integral over an ordered time simplex, and its
volume gives the bound
\[
 C(C_r|I|)^n/n!
\]
for the $n$th correction. The series converges in every spatial
H\"older norm, and differentiation of the equation gives smoothness
in time. The same argument applied to the difference of two solutions
proves uniqueness. In particular, the solution is zero on any interval
starting at $t_0$ on which its initial data and source are zero.

For the estimates, the tame product inequality bounds the coupling
in $C^{r,\alpha}$ by $C_r(H_0U_r+H_rU_0)$. We apply
\eqref{high:component-transport-low} when $r=0$ and the version of
\eqref{high:component-transport-tame} with one additional derivative
of the coefficient when $r\ge1$. Summing over $i$ and applying
Gronwall's inequality proves the assertions, with only the lowest
coefficient norms in the exponential. This argument applies in
particular to $T_\nu=\mathbb P$ and $\IId-\mathbb P$.
\end{proof}

\begin{corollary}[Linear projected transport]
\label{high:projected-cauchy}
Let $a,f$ be smooth on $I\times\mathbb T^3$, where $I$ is a compact interval
containing $t_0$, and let $W_0$ be smooth, with
$\ddiv a=\ddiv f=\ddiv W_0=0$. The Cauchy problem
\begin{equation}\label{high:linear-projected-transport}
 \partial_t W+\mathbb P(a\cn W)=f,\qquad W(t_0)=W_0
\end{equation}
has a unique smooth divergence-free solution, whose mean satisfies
\[
 \int_{\mathbb T^3}W(t,x)\,\dd x
 =\int_{\mathbb T^3}W_0(x)\,\dd x
      +\int_{t_0}^t\int_{\mathbb T^3}f(s,x)\,\dd x\,\dd s.
\]
For $r\ge1$, write
$A_r=\|\DD a\|_{r-1+\alpha}$ and
$W_r=\|\DD W\|_{r-1+\alpha}$. Then
\begin{align}
 W_1(t)&\le C e^{C\int_{t_0}^t A_1}
       \left(W_1(t_0)+\int_{t_0}^t\|\DD f(s)\|_\alpha\,\dd s\right),
       \label{high:projected-gradient-low}\\
 W_r(t)&\le C_r e^{C_r\int_{t_0}^t A_1}
       \left(W_r(t_0)+\int_{t_0}^t
       [\|\DD f(s)\|_{r-1+\alpha}+A_r(s)W_1(s)]\,\dd s\right),
       \quad r\ge2.
       \label{high:projected-gradient-high}
\end{align}
Only derivatives of $a,W,f$ through $r+\alpha$ enter these estimates.
In particular, the first gradient bound requires no second derivative
of $a$.
\end{corollary}
\begin{proof}
For arbitrary $W$, divergence freedom of $a$ gives
\[
 \Gamma_aW:=\nabla\Delta^{-1}(\partial_i a_j\,\partial_jW_i)
       =\nabla\Delta^{-1}\partial_j((\partial_i a_j)W_i).
\]
The expression to which the inverse Laplacian is applied has zero
mean. The last identity also writes $\Gamma_aW$ as a sum of the
fixed order-zero operators $\partial_k\Delta^{-1}\partial_j$
applied to $(\partial_i a_j)W_i$. Corollary~\ref{high:projected-transport}
thus gives a smooth solution on $I$ of
\[
 (\partial_t+a\cn)W=f+\Gamma_aW.
\]
Taking divergence, the derivative of the transport coefficient
cancels the divergence of $\Gamma_aW$, so that
\[
 (\partial_t+a\cn)\ddiv W=0.
\]
The divergence vanishes initially, so $\ddiv W=0$ and
\[
 \Gamma_aW=(\IId-\mathbb P)(a\cn W)
          =(\IId-\mathbb P)(W\cn a).
\]
Thus the constructed solution satisfies
\eqref{high:linear-projected-transport}; integration in space gives
the asserted mean. To obtain the gradient estimates, differentiate
the transport equation once:
\[
 (\partial_t+a\cn)\DD W
        =\DD f+\DD\Gamma_aW-(\DD a)\DD W.
\]
Schauder and product estimates give
\[
 \|\DD\Gamma_aW\|_{r-1+\alpha}
 \le C_r(A_1W_r+A_rW_1).
\]
When $r=1$, \eqref{high:component-transport-low} and the estimate
for $\DD\Gamma_aW$ give an integral inequality with coefficient
$A_1$, and Gronwall's inequality proves the first assertion.
For $r\ge2$, we apply \eqref{high:component-transport-tame} to
$\DD W$ at order $r-1$. Interpolation bounds the transport
commutator, as well as the pressure term, by
$C_r(A_1W_r+A_rW_1)$. A second application of Gronwall's inequality
gives the higher gradient estimates.
\end{proof}

\begin{proof}[Proof of Lemma~\ref{prep:ordinary-comparison}]
Recall the local fields of Section~\ref{sec:preparation}. Since
\eqref{prep:mollified-equations} and \eqref{prep:local-equation}
contain the same stress $R_\ell$, their difference satisfies a linear
system with the convolution defects as sources. We construct the
local solution and prove the spatial estimates first. We then estimate
the difference and use the equations to obtain the ordinary time
estimates. We use throughout the inductive hypotheses
\eqref{prep:old-ordinary} and \eqref{prep:old-stress-time}, the
comparisons in \eqref{setup:rate-comparisons}, and
Proposition~\ref{moll:spatial-basic}. In particular, the latter gives
the factor $\varepsilon_\ell^2=(\ell\lambda_q)^{2m_0}$ in the
mollification error.

Fix $n$ and work on $[t_n,t_n+T]$, where $T=4\tau_c$; reversing time gives
the other half of $I_n^{\rm bg}$. Put $f=\mathbb P\ddiv R_\ell$.

\emph{1. Local existence.}
Starting with $Z_{n,0}^\pm=z_\ell^\pm(t_n)$, define the iterates by
\[
 \partial_tZ_{n,j+1}^\pm+
       \mathbb P(Z_{n,j}^\mp\cn Z_{n,j+1}^\pm)=f,
 \qquad Z_{n,j+1}^\pm(t_n)=z_\ell^\pm(t_n).
\]
Corollary~\ref{high:projected-cauchy} gives a smooth solution at each
iteration on $[t_n,t_n+T]$. It preserves divergence freedom
and, since $\int f=0$, the initial mean. We first estimate the gradients by setting
\[
 \operatorname{Lip}_{n,j}=\sum_\pm\|\DD Z_{n,j}^\pm\|_\alpha,
 \qquad F_1=\sup_t\|f(t)\|_{1+\alpha}.
\]
The projected gradient estimate
\eqref{high:projected-gradient-low} gives
\[
 \operatorname{Lip}_{n,j+1}(t)
 \le C\exp\left(C\int_{t_n}^t\operatorname{Lip}_{n,j}(s)\,\dd s\right)
                         (\operatorname{Lip}_{n,0}+TF_1).
\]
Set
\[
 \operatorname{Lip}_{n,\rm cell}=1+\operatorname{Lip}_{n,0}+TF_1.
\]
The data and forcing satisfy
\[
 \begin{aligned}
 \operatorname{Lip}_{n,0}&\le C\ell^{-\alpha}\lambda_q\delta_q^{1/2},\qquad
 F_1\le C\ell^{-\alpha}\lambda_q^{2-\alpha}\delta_{q+1},\\
 T\operatorname{Lip}_{n,\rm cell}
 &\le C\bigl(\tau_c+\tau_c\ell^{-\alpha}\lambda_q\delta_q^{1/2}
       +\tau_c^2\ell^{-\alpha}\lambda_q^{2-\alpha}\delta_{q+1}\bigr)
 =o(1).
 \end{aligned}
\]
The last quantity tends to zero. Thus the preceding gradient
inequality closes by induction, after increasing the large parameter,
and gives
\[
 \sup_j\sup_t\operatorname{Lip}_{n,j}(t)\le 2C\operatorname{Lip}_{n,\rm cell}
\]
throughout $[t_n,t_n+T]$.

For $0<\alpha'<\alpha$, set $d_{n,j+1}^\pm=Z_{n,j+1}^\pm-Z_{n,j}^\pm$.
Using the identity
\[
 (\IId-\mathbb P)(a\cn d)=(\IId-\mathbb P)(d\cn a)
\]
in the difference of the iteration equations gives
\[
 (\partial_t+Z_{n,j}^\mp\cn)d_{n,j+1}^\pm
  =(\IId-\mathbb P)(d_{n,j+1}^\pm\cn Z_{n,j}^\mp)
      -\mathbb P(d_{n,j}^\mp\cn Z_{n,j}^\pm),
 \qquad d_{n,j+1}^\pm(t_n)=0.
\]
The right-hand side contains no derivative of the differences:
each derivative falls on an iterated field. We may therefore apply
\eqref{high:coupled-transport-low} with exponent $\alpha'$ and
treat the last term as the source. We obtain
\[
 \sup_t\sum_\pm\|d_{n,j+1}^\pm\|_{\alpha'}
 \le CT\operatorname{Lip}_{n,\rm cell}e^{CT\operatorname{Lip}_{n,\rm cell}}
                   \sup_t\sum_\pm\|d_{n,j}^\pm\|_{\alpha'}
 \le\tfrac12\sup_t\sum_\pm\|d_{n,j}^\pm\|_{\alpha'},
\]
where the last inequality follows by increasing $a$.

\emph{2. Spatial estimates and smoothness.}
For $r\ge1$, put
\[
 \operatorname{Lip}_{n,j,r}=\sum_\pm\|\DD Z_{n,j}^\pm\|_{r-1+\alpha}.
\]
For $r\ge2$, \eqref{high:projected-gradient-high} and the uniform
bound on the lowest gradient norm give a linear integral inequality,
with coefficient $C_r\operatorname{Lip}_{n,\rm cell}$, for
$\operatorname{Lip}_{n,j,r},\operatorname{Lip}_{n,j+1,r}$.
Iterating this inequality gives
\[
 \sup_j\operatorname{Lip}_{n,j,r}(t)
 \le C_re^{C_rT\operatorname{Lip}_{n,\rm cell}}
       \left(\operatorname{Lip}_{n,0,r}
               +\int_{t_n}^{t_n+T}\|f(s)\|_{r+\alpha}\,\dd s\right);
\]
the $j$th iterated integral is bounded by
$(C_rT\operatorname{Lip}_{n,\rm cell})^j/j!$. No smallness of the
higher norms is needed. By interpolation, these uniform estimates
and the convergence in Step 1 imply convergence in every spatial
norm. The equations then give convergence of the time derivatives.
The pressures in the two limiting equations agree, since their
Poisson equations have the same right-hand side:
\[
 \partial_i(z_{\ell,n}^-)_j\partial_j(z_{\ell,n}^+)_i
 =\partial_i(z_{\ell,n}^+)_j\partial_j(z_{\ell,n}^-)_i.
\]
The limit therefore solves \eqref{prep:local-equation}; the difference
estimate from Step 1 gives uniqueness. To make the spatial bounds
explicit, convolution gives, for every $r$,
\[
 \begin{gathered}
 \operatorname{Lip}_{n,0,r}
 \le C_r\ell^{-\alpha}\lambda_q^{\min\{r,\rgood\}}
      \ell^{-[r-\rgood]^+}\delta_q^{1/2}
 \le C_r\ell^{-\alpha}\lambda_q\ell^{-(r-1)}\delta_q^{1/2},
 \\
 \|f\|_{r+\alpha}
 \le C_r\ell^{-\alpha}\lambda_q^{\min\{r+1,\rgood\}-\alpha}
      \ell^{-[r+1-\rgood]^+}\delta_{q+1},
 \end{gathered}
\]
and the forcing integral is bounded by the same quantity as the initial
data: the ratio of their bounds is at most
\[
 \tau_c\lambda_q^{1-\alpha}\frac{\delta_{q+1}}{\delta_q^{1/2}}
 =\lambda_q^{-\alpha}\varepsilon_{q+1}^{\gamma_\ell}
                       \frac{\delta_{q+1}}{\delta_q}\le1.
\]
Substitution in the preceding estimate proves
\eqref{prep:whole-patch-fields} and its lossy version for $k=0$.
The restriction $r+3\le\rgood$ leaves three additional spatial
derivatives, accounting for the divergence in the forcing and the
spatial derivative used for H\"older interpolation.

\emph{3. Spatial comparison.}
We now subtract the equations for $z_\ell^\pm$ and
$z_{\ell,n}^\pm$. Cancellation of the common stress leaves the two
convolution defects as sources:
\begin{equation}
 \begin{aligned}
 (\partial_t+z_{\ell,n}^\mp\cn)\Delta_n^\pm
 ={}&\mathbb P\ddiv R_\ell^c\pm\curl M_\ell
 -\mathbb P(\Delta_n^\mp\cn z_\ell^\pm)\\
 &+\nabla\Delta^{-1}\ddiv
                  (\Delta_n^\pm\cn z_{\ell,n}^\mp),
 \qquad \Delta_n^\pm(t_n)=0.
 \end{aligned}
\end{equation}
Corollary~\ref{high:projected-transport} applies to this equation.
The fixed multipliers are $\mathbb P$ and $\IId-\mathbb P$; the
coefficients are the gradients
$\DD z_{\ell,n}^\pm,\DD z_\ell^\pm$, and the sources are
$\mathbb P\ddiv R_\ell^c\pm\curl M_\ell$. The spatial estimates give
$H_r\le C\ell^{-\alpha}\lambda_q^{r+1}\delta_q^{1/2}$ for
$r+1\le\rgood-3$, while $TH_0=o(1)$ controls the exponential.
The quadratic convolution estimate for the sources, followed by
\eqref{high:coupled-transport-low}--\eqref{high:coupled-transport-high},
then gives
\[
 \|\Delta_n^\pm\|_{r+\alpha}
 \le C\tau_c\ell^{-\alpha}\varepsilon_\ell^2\lambda_q^{r+1}\delta_q
 =C\ell^{-\alpha}\varepsilon_\ell^2\lambda_q^r
                         \varepsilon_{q+1}^{\gamma_\ell}\delta_q^{1/2}.
\]
Here the coefficients require spatial derivatives of the local field
through order $r+1$. These are bounded by
\eqref{prep:whole-patch-fields}, since $r+1\le\rgood-3$.
For the source, taking the divergence and using one additional spatial
derivative for H\"older interpolation requires the old fields through
spatial order $r+2m_0+3$, which is at most $\rgood$ by the
restriction $r\le\rgood-2m_0-3$ in the comparison estimate.
Adding the mollification error proves
\eqref{prep:complete-comparison} for $k=0$. Both errors tend to
zero; hence $\|z_{\ell,n}^\pm\|_0\le C$ and $|B_{\ell,n}|$ is
bounded away from zero. This argument uses the approximation estimate
and does not require positivity of the convolution kernel.

For $G_n^\pm$, the lossy estimate follows by the triangle inequality
and interpolation between consecutive integer spatial estimates.
The old fields satisfy \eqref{prep:old-ordinary} through
$\rgood+4$, with $\lambda_q\le\ell^{-1}$; the local fields satisfy
the lossy version of \eqref{prep:whole-patch-fields}, restated in
\eqref{prep:local-field-classes}.

\emph{4. Independent time derivatives.}
For $k\ge1$ the projected equation gives the exact recursion
\[
 \partial_t^kz_{\ell,n}^\pm
 =\partial_t^{k-1}f
 -\sum_{h=0}^{k-1}\binom{k-1}{h}\mathbb P
       (\partial_t^hz_{\ell,n}^\mp\cn
                  \partial_t^{k-1-h}z_{\ell,n}^\pm).
\]
After applying spatial derivatives, the numbers of spatial and
ordinary time derivatives in each product sum to $r+k$. We apply the
induction hypothesis to the factor with the most derivatives. Each
other differentiated field contributes a factor at most
$\ell^{-\alpha}\delta_q^{1/2}\le1$, and undifferentiated fields are
bounded by $C$. This proves \eqref{prep:whole-patch-fields} by
induction through $k\le j_1+1$: the source requires time derivatives
of the stress only through $k-1\le j_1$. The lossy assertion follows
by the same induction from the lossy spatial bounds. For the difference
of the fields, we use the identity
\begin{align*}
 \partial_t^k\Delta_n^\pm
 ={}&\mathbb P\ddiv\partial_t^{k-1}R_\ell^c
                     \pm\curl\partial_t^{k-1}M_\ell\\
 &-\sum_{h=0}^{k-1}\binom{k-1}{h}\mathbb P
 \bigl(\partial_t^hz_{\ell,n}^\mp\cn
                         \partial_t^{k-1-h}\Delta_n^\pm
       +\partial_t^h\Delta_n^\mp\cn
                         \partial_t^{k-1-h}z_\ell^\pm\bigr).
\end{align*}
Each factor now has fewer than $k$ ordinary time derivatives, and the
numbers of spatial and time derivatives sum to $r+k$. The source
terms obey the same estimate because $\tau_c^{-1}\le\lambda_q$:
the inverse time factor in the spatial comparison is bounded by the
spatial frequency. Induction proves
\eqref{prep:correction-comparison}, and adding the estimate in
\eqref{moll:spatial-accuracy} gives \eqref{prep:complete-comparison}.
\end{proof}

\subsection{Comparison of ordinary, transport and Lie derivative estimates}
\label{ssec:finite-derivative-conversion}

The following three lemmas compare ordinary, transport, and Lie
derivatives in the integer norms $\|\cdot\|_r$. They also hold in
$\|\cdot\|_{r+\alpha}$ with a common fixed factor $\varkappa\ge1$
on the right of every H\"older hypothesis and conclusion, provided
the corresponding integer bounds hold. Indeed, the proofs use only
Leibniz' rule and \eqref{setup:product}, whose terms contain just one
H\"older norm. For the classes in Section~\ref{ssec:fixed-scales},
this factor is $\varkappa=\max\{\ell^{-1},\Lambda\}^{\alpha}$.

\begin{lemma}[Comparison of transport derivatives]
\label{high:ordinary-transport-comparison}
Let $D_t=\partial_t+v\cn$, $D_B=B\cn$ and
$D_t'=\partial_t+v'\cn$, $D_B'=B'\cn$. Fix integers $N\ge1$,
$K\ge0$, and $\Lambda\ge1$. For $h\le K$ and $r+h\le N-1$,
assume $\|\partial_t^h w\|_{r}\le C\Lambda^{r+h}$ for
$w=v,B,v',B'$, and the same bound with an additional factor
$e\ge0$ for $w=v'-v,B'-B$. Assume also
\[
 \|\partial_t^h f\|_{r}\le F\Lambda^{r+h},\qquad
 \|\partial_t^h(f'-f)\|_{r}\le E\Lambda^{r+h}
 \qquad h\le K,\quad r+h\le N.
\]
Then, for $k\le K$ and $r+k+m\le N$,
\begin{equation}
 \|(D_t')^k(D_B')^m f'-D_t^kD_B^mf\|_{r}
       \le C_{N,K}(E+eF)\Lambda^{r+k+m}.
\end{equation}
The same estimate holds for any ordering of the derivatives when
primed and unprimed operators occupy corresponding positions. No
commutation relation is needed, and ordinary time derivatives are
required only through order $k$.
\end{lemma}
\begin{proof}
For $d=k+m$ operators in corresponding positions,
\[
 L_1'\cdots L_d'-L_1\cdots L_d
 =\sum_{a=1}^d L_1'\cdots L_{a-1}'(L_a'-L_a)
                                      L_{a+1}\cdots L_d.
\]
Apply this identity to $f$ and add
$L_1'\cdots L_d'(f'-f)$. In the first sum, Leibniz' rule leaves a
derivative of a coefficient difference in every term; in the added
term it leaves a derivative of the difference of the functions.
After spatial differentiation, the numbers of spatial and ordinary
time derivatives in each product sum to at most $r+d$. Those on any
coefficient sum to at most $r+d-1$. Moreover, only the $k$ material
derivatives can produce ordinary time derivatives. The hypotheses
therefore give a factor $eF$ for the first sum and $E$ for the added
term. The product inequality gives the bound
$C(eF+E)\Lambda^{r+d}$.
\end{proof}

\begin{lemma}[Commuting a spatial gradient]
\label{high:transport-gradient}
Let $D_t=\partial_t+v\cn$ and $D_B=B\cn$ commute. Fix integers
$J,N\ge1$, $\Lambda\ge1$, amplitudes $a,b\ge0$ and positive
derivative costs $\mathrm a_t,\mathrm a_B$, with $\Lambda a\le \mathrm a_t$, $\Lambda b\le \mathrm a_B$.
Suppose, for $k+m\le J$ and $r+k+m\le N-1$,
\[
 \|D_t^kD_B^mv\|_{r+1}\le C a\Lambda^{r+1}\mathrm a_t^k\mathrm a_B^m,
 \qquad
 \|D_t^kD_B^mB\|_{r+1}\le C b\Lambda^{r+1}\mathrm a_t^k\mathrm a_B^m.
\]
Then, for $k+m\le J$ and $r+k+m\le N-1$,
\begin{equation}\label{high:transport-gradient-coefficients}
 \begin{aligned}
 \|D_t^kD_B^m\nabla v\|_{r}
       &\le C'a\Lambda^{r+1}\mathrm a_t^k\mathrm a_B^m,\\
 \|D_t^kD_B^m\nabla B\|_{r}
       &\le C'b\Lambda^{r+1}\mathrm a_t^k\mathrm a_B^m.
 \end{aligned}
\end{equation}
For $k+m\le J$ and $r+k+m\le N-1$,
\[
 \|D_t^kD_B^mF\|_{r+1}
 \le F_0\Lambda^{r+1}\mathrm a_t^k\mathrm a_B^m
\]
implies
\[
 \|D_t^kD_B^m\nabla F\|_{r}
 \le C'F_0\Lambda^{r+1}\mathrm a_t^k\mathrm a_B^m;
\]
the latter bound conversely
implies the same bound for $\nabla(D_t^kD_B^mF)$.
The comparison for a general field also holds if
\eqref{high:transport-gradient-coefficients} is assumed directly.
When the orders of the material and magnetic derivatives sum to
$k+m$, the commutator terms require these derivatives of the
background gradients only when their orders sum to at most $k+m-1$.
\end{lemma}
\begin{proof}
Use
\[
 [D_t,\partial_i]=-(\partial_iv^j)\partial_j,
 \qquad [D_B,\partial_i]=-(\partial_iB^j)\partial_j
\]
and argue by induction on $k+m$. Moving the gradient to the left
gives the leading term $\nabla D_t^kD_B^mF$. Every commutator term
is a product of a differentiated background gradient and derivatives
of $\nabla F$ with strictly fewer material and magnetic derivatives
combined. A commutation with the material derivative gives a gradient
bounded by $\Lambda a\le \mathrm a_t$; a commutation with the magnetic
derivative gives one bounded by $\Lambda b\le \mathrm a_B$. Leibniz'
rule distributes the remaining derivatives between the two factors.
Their material orders add to the original material order after counting
the first commutation, and the same is true of their magnetic orders.
The product estimate therefore gives the asserted bound.

For $v,B$, this is a simultaneous induction. Every background
gradient in a commutator has fewer material and magnetic derivatives
combined than the gradient being estimated, so the induction
hypothesis applies. Commutativity, $[D_t,D_B]=0$, allows us to put
them in the stated order. Including the spatial derivatives, each
term satisfies $r+k+m+1\le N$. This proves the estimates for the
background gradients and hence the comparison for a general field.
To prove the converse, solve the same identity for
$\nabla D_t^kD_B^mF$. The commutator terms require at most
$k+m-1$ material and magnetic derivatives combined on a background
gradient, and the same induction bounds them for $k+m\le J$ and
$r+k+m\le N-1$.
\end{proof}

\begin{lemma}[Transport and Lie derivatives]
\label{high:component-lie-conversion}
Let $D_t=\partial_t+v\cn$ and $D_B=B\cn$ commute. Fix a tensor
type, $J,N\ge1$, $\Lambda\ge1$ and positive derivative costs $\mathrm a_t,\mathrm a_B$.
Assume
\begin{equation}\label{high:conversion-strain-input}
 \begin{aligned}
 \|D_t^kD_B^m\nabla v\|_{r}&\le C\Lambda^r\mathrm a_t^{k+1}\mathrm a_B^m,\\
 \|D_t^kD_B^m\nabla B\|_{r}&\le C\Lambda^r\mathrm a_t^k\mathrm a_B^{m+1},
 \end{aligned}
 \qquad k+m\le J-1,\quad r+k+m\le N-1.
\end{equation}
The families
\begin{equation}
 \begin{aligned}
 \|D_t^kD_B^mF\|_{r}&\le F_0\Lambda^r\mathrm a_t^k\mathrm a_B^m,\\
 \|\mathcal D_t^k\mathcal L_B^mF\|_{r}
                      &\le F_0\Lambda^r\mathrm a_t^k\mathrm a_B^m
 \end{aligned}
 \qquad k+m\le J,\quad r+k+m\le N
\end{equation}
are equivalent for every ordering, up to a constant depending only
on the fixed orders, tensor type, and the constant in the gradient
hypotheses. Consequently, for a $(\Lambda,\rho)$-adapted background
through $(N,J)$ in the sense of Section~\ref{ssec:fixed-scales}, the
transport and Lie definitions of $\mathcal C_{N,J}(F_0;\Lambda,\rho)$
give the same class up to a fixed constant. With the corresponding
lossy gradient hypotheses, the same equivalence holds for
$\mathcal K_{N,J}(F_0')$. If those hypotheses hold at every spatial
order, the equivalence also holds for $\mathcal K_J(F_0')$.
\end{lemma}
\begin{proof}
In Cartesian components,
\[
 \mathcal D_t=D_t+Z_t,\qquad \mathcal L_B=D_B+Z_B,
\]
where $Z_t,Z_B$ are linear in $\nabla v,\nabla B$ and act on the
tensor indices. Expand a composition of Lie derivatives using these
identities and Leibniz' rule. Its leading term is the corresponding
composition of material and magnetic derivatives of $F$. Each other
term contains differentiated background gradients and strictly fewer
material and magnetic derivatives combined on the tensor. Assign
$(1,0)$ to $Z_t$ and $(0,1)$ to $Z_B$, counting each as the derivative
it replaces. The degrees of the factors in each term then sum to
$(k,m)$. A background gradient has at most $k+m-1$ material and
magnetic derivatives combined; after spatial differentiation, the
sum of its spatial, material and magnetic orders is at most
$r+k+m-1$. Thus \eqref{high:conversion-strain-input} and the
product estimate give the bound for the Lie derivatives.

The expansion is triangular in the sum of the material and magnetic
derivative orders on the tensor. Solving it for $D_t^kD_B^mF$
and inducting on $k+m$ proves the converse. Since $[D_t,D_B]=0$ also implies
$[\mathcal D_t,\mathcal L_B]=0$, the estimates hold for every
ordering. For $\mathcal K_{N,J}(F_0')$, and for
$\mathcal K_J(F_0')$ under the hypotheses at every spatial order,
observe that replacing a material
or magnetic derivative by its zeroth-order tensor action introduces
one spatial derivative of the corresponding field. Thus the numbers
of spatial, material and magnetic derivatives in each product still
sum to at most $r+k+m$. The same product estimate proves the assertion.
\end{proof}
\subsection{Tensor pushforward estimates}

The following estimate for tensor pushforwards will be applied to the
corrector path average. We include an inhomogeneous term in the Lie
equation. In Cartesian components, this equation is a transport system
whose zeroth-order coefficients are the entries of the gradient of
the LDF.

\begin{lemma}[Tensor pushforward]
\label{high:natural-transport}
Fix an integer $N\ge0$, $0<\alpha<1$, $\Lambda\ge1$ and a
tensor type. Let $\xi_s$ be smooth and divergence free for
$0\le s\le1$, and suppose
\[
 \sup_s\|\xi_s\|_{r+\alpha}
       \le \epsilon\Lambda^{r-1},\qquad 0\le r\le N+1,
 \qquad 0\le\epsilon\le C_0.
\]
For this lemma write
\[
 \|F\|_{N,\Lambda}=\max_{0\le r\le N}
 \Lambda^{-r}\|F\|_{r+\alpha}.
\]
If
$(\partial_s+\mathcal L_{\xi_s})F_s=H_s$, then
\begin{equation}\label{high:natural-transport-bound}
 \sup_s\|F_s\|_{N,\Lambda}
 \le C\left(\|F_0\|_{N,\Lambda}
              +\int_0^1\|H_s\|_{N,\Lambda}\,\dd s\right).
\end{equation}
For a homogeneous solution,
\begin{equation}\label{high:natural-transport-defect}
 \sup_s\|F_s-F_0\|_{N,\Lambda}
              \le C\epsilon\|F_0\|_{N+1,\Lambda}.
\end{equation}
The constants depend only on the fixed orders, tensor type,
$\alpha$ and $C_0$. The same result holds with ordinary time
derivatives: suppose the coefficient estimate for
$\partial_t^k\xi_s$ is $\epsilon\Lambda^{r+k-1}$ on
$r+k\le N+1$, $k\le K$. Then bounds of size
$A\Lambda^{r+k}$ for the initial tensor and source give the same
bounds for the solution on $r+k\le N$, $k\le K$. The estimates
require at most $K$ time derivatives of the LDF, initial tensor,
and source.

For $N\ge2$, there is also a tame estimate. Assume the displayed
LDF bound without loss through $r=2$, and with an additional
factor $M\ge1$ at the remaining orders through $r=N+1$. Then
\begin{equation}\label{high:natural-transport-tame}
 \begin{aligned}
 \sup_s\|F_s\|_{N,\Lambda}
 \le C\bigg(&\|F_0\|_{N,\Lambda}
       +\int_0^1\|H_s\|_{N,\Lambda}\,\dd s\\
 &+\epsilon M\Big[\|F_0\|_{1,\Lambda}
       +\int_0^1\|H_s\|_{1,\Lambda}\,\dd s\Big]\bigg),
 \end{aligned}
\end{equation}
where $C$ is independent of $M$. In particular, the higher spatial norms of the coefficients enter
linearly; only their lower spatial H\"older norms occur in the
exponential. The same assertions hold with $\|\cdot\|_r$ in place
of $\|\cdot\|_{r+\alpha}$. Thus, for the classes of
Section~\ref{ssec:fixed-scales}, pushforward by an LDF satisfying
these hypotheses preserves the integer and H\"older bounds, with the
same amplitude. The difference from the initial tensor satisfies the
corresponding bound with the additional factor $\epsilon$.
\end{lemma}
\begin{proof}
In Cartesian components write
\[
 \mathcal L_{\xi_s}F=\xi_s\cn F+Z_sF,
\]
where $Z_s$ is linear in $\DD\xi_s$, with one term for each tensor index.
Spatial differentiation therefore gives, for
$|\boldsymbol\gamma|=r$,
\[
 \begin{aligned}
 (\partial_s+\xi_s\cn)\partial^{\boldsymbol\gamma}F_s
 ={}&\partial^{\boldsymbol\gamma}H_s
 -\partial^{\boldsymbol\gamma}(Z_sF_s)\\
 &-\sum_{0<\boldsymbol\beta\le\boldsymbol\gamma}
 \binom{\boldsymbol\gamma}{\boldsymbol\beta}
 (\partial^{\boldsymbol\beta}\xi_s)\cn
       \partial^{\boldsymbol\gamma-\boldsymbol\beta}F_s.
 \end{aligned}
\]
The coefficient hypotheses and the product estimate bound the
zeroth-order tensor action and the transport commutators as follows:
\[
 \Lambda^{-r}\left(
 \|\partial^{\boldsymbol\gamma}(Z_sF_s)\|_\alpha
 +\sum_{0<\boldsymbol\beta\le\boldsymbol\gamma}
 \|(\partial^{\boldsymbol\beta}\xi_s)\cn
       \partial^{\boldsymbol\gamma-\boldsymbol\beta}F_s\|_\alpha
 \right)
 \le C_N\epsilon\|F_s\|_{N,\Lambda}.
\]
Indeed, a transport commutator contains spatial derivatives of orders
$b$ and $r-b+1$ of the LDF and tensor, respectively. A term from
the zeroth-order action has the corresponding orders $b+1$ and
$r-b$. The assumed powers of the frequency therefore give exponent
$r$ in either case, and require the LDF only through spatial order
$r+1$.

We apply \eqref{high:component-transport-low} to these differentiated
components, divided by their frequency weights. The bound
$\|\DD\xi_s\|_\alpha\le C\epsilon$ controls the transport flow, and
the displayed inequality bounds the couplings by $C_N\epsilon$.
Consequently, \eqref{high:coupled-transport-low} gives
\eqref{high:natural-transport-bound}.

For $G_s=F_s-F_0$ in the homogeneous case,
\[
 \begin{gathered}
 (\partial_s+\mathcal L_{\xi_s})G_s=-\mathcal L_{\xi_s}F_0,
 \qquad G_0=0,\\
 \|\mathcal L_{\xi_s}F_0\|_{N,\Lambda}
       \le C\epsilon\|F_0\|_{N+1,\Lambda}.
 \end{gathered}
\]
Applying \eqref{high:natural-transport-bound} with this source and
zero initial data proves \eqref{high:natural-transport-defect}.

For the ordinary time estimates, differentiation by
$\partial_t^k\partial^{\boldsymbol\gamma}$ produces transport
commutators of the form
\[
 (\partial_t^a\partial^{\boldsymbol\beta}\xi_s)\cn
       \partial_t^{k-a}\partial^{\boldsymbol\gamma-\boldsymbol\beta}F_s,
 \quad a+|\boldsymbol\beta|\ge1,
\]
together with the terms from the tensor action, which contain one
additional spatial derivative of $\xi_s$. In each product the
frequency exponents sum to $r+k$. The spatial and ordinary time
orders sum to at most $r+k$ on the tensor and $r+k+1$ on the LDF;
the ordinary time order on either factor is at most $k\le K$.
The hypotheses thus bound every coefficient of the differentiated
system after division by the frequency weights, and
\eqref{high:coupled-transport-low} applies as before.

For the tame estimate, apply the result with $N=1$ first. This uses
the LDF only through $r=2$ and bounds the lowest tensor norms
independently of the higher coefficient norms. At higher spatial
orders, the tame product and commutator estimates give
\[
 C\epsilon\|F_s\|_{N,\Lambda}
       +C\epsilon M\|F_s\|_{1,\Lambda}.
\]
To see this, the highest derivative of the LDF in the tensor action
has spatial order $N+1$ and multiplies the undifferentiated tensor.
In a transport commutator it has order $N$ and multiplies one spatial
derivative of the tensor. Division by the powers of $\Lambda$
leaves the stated bound for these two products; interpolation gives
the same bound for the intermediate products.

We now use \eqref{high:component-transport-low} for the differentiated
system. The lowest tensor norm has already been estimated and may be
included in the source. Gronwall's inequality absorbs only
$C\epsilon\|F_s\|_{N,\Lambda}$, so $M$ does not occur in the
exponential. This proves \eqref{high:natural-transport-tame} and the
claimed linear dependence on the highest spatial norms of the LDF.
\end{proof}

\subsection{Singular integral commutators}

We next prove commutator estimates for spatial Fourier multipliers
acting on components and on tensors. A commutation with a transport
derivative differentiates the kernel once and introduces a difference
of the transport vector field. This difference vanishes to first
order on the diagonal and compensates for the derivative of the
kernel. For iterated commutators, we use the full distributional
kernel, including its diagonal part.

Fix $0<\alpha<1$, with spatial norms $\|\cdot\|_{r+\alpha}$ on
$\mathbb T^3$. All identities are first proved for smooth functions
and then extended under the stated regularity assumptions. We
consider spatial operators $T$ which, at each fixed time, have the
componentwise representation
\begin{equation}\label{high:cz-kernel-representation}
 Tf(x)=c_Tf(x)+\operatorname{p.v.}
          \int_{\mathbb T^3}K(x-y)f(y)\,\dd y.
\end{equation}
Here $K$ is smooth away from the origin and agrees near the origin,
up to a smooth function, with a homogeneous kernel $K_0$ of degree $-3$
satisfying $\int_{\mathbb S^2}K_0=0$.
We use only this hypothesis on the kernel. It is satisfied by the
order-zero multipliers occurring here, including $\mathbb P$,
$\mathcal R\ddiv$, $\mathcal R\curl$ and
$\nabla\Delta^{-1}\ddiv$. Changing the zero Fourier mode changes
only the smooth part of the kernel. The constant part commutes with
component transport derivatives, although it also contributes to tensor
Lie commutators.

The representation follows from the Euclidean Green functions
\[
 G_1^{\mathbb R^3}(x)=-\frac1{4\pi|x|},\qquad
 G_2^{\mathbb R^3}(x)=-\frac{|x|}{8\pi},\qquad
 \Delta G_1^{\mathbb R^3}=\delta_0,\quad
 \Delta G_2^{\mathbb R^3}=G_1^{\mathbb R^3}.
\]
The unit flux of $\nabla G_1^{\mathbb R^3}$ through a sphere proves
the first distributional identity. Ordinary differentiation proves
the second away from zero; the vanishing boundary flux excludes a
point mass at the origin. A cutoff $\chi$ equal to one near zero in
a periodic coordinate ball therefore gives, for $a=1,2$,
\[
 \Delta^a(\chi G_a^{\mathbb R^3})=\delta_0+r_a,
 \qquad r_a\in C^\infty(\mathbb T^3),\quad \int r_a\,\dd x=-1.
\]
The equation
\[
 \Delta^aw_a=-|\mathbb T^3|^{-1}-r_a
\]
has a right-hand side of zero mean. Inverting $\Delta^a$ on the
nonzero Fourier modes gives a periodic smooth solution: after any
number of spatial differentiations, its Fourier coefficients still
decay rapidly. We choose the additive constant so that the periodic
Green function has zero mean and obtain
\[
 G_a^{\mathbb T^3}=\chi G_a^{\mathbb R^3}+w_a,
 \qquad \Delta^aG_a^{\mathbb T^3}=\delta_0-|\mathbb T^3|^{-1}.
\]
The operators under consideration are constant combinations of the
identity, second derivatives of $G_1^{\mathbb T^3}$ and fourth
derivatives of $G_2^{\mathbb T^3}$; this includes
$\mathscr T=\curl^{-1}\mathbb P\ddiv$.
Their singular kernels have degree $-3$. Writing the last derivative
as the divergence of a homogeneous vector field of degree $-2$ gives
\[
 \log(R/r)\int_{\mathbb S^2}K_0\,\dd S
 =\int_{r<|x|<R}\ddiv V\,\dd x=0,
\]
because the two boundary fluxes agree. Thus the spherical mean vanishes. Integration by parts outside
$B(0,r)$ identifies the full distribution as
$\operatorname{p.v.}K+c_T\delta_0$. Only the constant Taylor
coefficient of the test function contributes in the boundary limit;
every other term contains a positive power of $r$ and tends to zero.
Adding the smooth periodic correction proves
\eqref{high:cz-kernel-representation}.

The same Green functions yield integer-norm bounds for the
operators of order $-1$.
Expanding \eqref{eq:inverse-divergence} gives
\[
 (\mathcal Rf)_{ij}
 =\partial_i\Delta^{-1}f_j+\partial_j\Delta^{-1}f_i
 -\tfrac12\partial_i\partial_j\Delta^{-2}\ddiv f
 -\tfrac12\delta_{ij}\Delta^{-1}\ddiv f.
\]
Its kernel consists of first derivatives of ${G_1^{\mathbb T^3}}$ and third
derivatives of ${G_2^{\mathbb T^3}}$. These have size $O(|x|^{-2})$ at the
origin and are smooth elsewhere, hence are integrable on the
torus. The same holds for $\curl^{-1}$. Convolution and
commutation with ordinary derivatives therefore give
\begin{equation}\label{high:integrable-hodge}
 \|\mathcal Ru\|_r+\|\curl^{-1}u\|_r\le C_r\|u\|_r,
 \qquad r\in\mathbb N_0.
\end{equation}
The operators annihilate the constant mode by convention. The
integer-norm boundedness asserted here will suffice for the small
remainder in the finite expansion below; no derivative gain in these
norms is asserted.

For the commutator calculation, the order-zero kernel is the full distribution
\[
 \mathcal K=c_T\delta_0+\operatorname{p.v.}K,\qquad Tf=\mathcal K*f.
\]
For $T=\mathcal R$ or $\curl^{-1}$, it is the integrable
distributional kernel just described. The same difference-kernel
notation and transport identities apply in both cases.

For smooth functions, a spatial multiindex $\boldsymbol\gamma$ of length
$m$, and $\boldsymbol f=(f_1,\ldots,f_m)$, define
\begin{equation}\label{high:cz-difference-kernel}
 T_{\boldsymbol\gamma,\boldsymbol f}u(x)
 =\left\langle\partial^{\boldsymbol\gamma}\mathcal K(w),
   \prod_{i=1}^m\bigl(f_i(x)-f_i(x-w)\bigr)u(x-w)\right\rangle.
\end{equation}
At $m=0$ this is $T$. Away from the diagonal its kernel is
\[
 k(x,y)=\prod_i(f_i(x)-f_i(y))\,
                         \partial^{\boldsymbol\gamma}K(x-y).
\]
In the order-zero case, for $w\to0$ the test factor in this
distributional pairing satisfies
\[
 \prod_{i=1}^m(f_i(x)-f_i(x-w))u(x-w)
 =u(x)\prod_{i=1}^m(\nabla f_i(x)\cdot w)+O(|w|^{m+1}).
\]
After integration by parts outside a small ball, only the degree-$m$
polynomial contributes to the boundary limit. Hence the difference
between the distributional kernel and its principal-value part is
multiplication by a fixed contraction of $\prod_i\nabla f_i(x)$;
it contains no derivative of $u$. The derivatives of $c_T\delta_0$
give terms of the same form. We include these diagonal terms in the
following estimate.

\begin{lemma}[Difference kernels]
For an order-zero $T$, $|\boldsymbol\gamma|=m$, every integer
$r\ge0$, and $0<\alpha<1$,
\begin{equation}\label{high:cz-kernel-estimate}
 \|T_{\boldsymbol\gamma,\boldsymbol f}u\|_{r+\alpha}
 \le C\sum_{r_0+\cdots+r_m=r}
       \|u\|_{r_0+\alpha}\prod_{i=1}^m\|\nabla f_i\|_{r_i+\alpha}.
\end{equation}
For $T=\mathcal R$ or $\curl^{-1}$, the corresponding integer estimate is
\begin{equation}\label{high:integrable-difference-estimate}
 \|T_{\boldsymbol\gamma,\boldsymbol f}u\|_r
 \le C\sum_{r_0+\cdots+r_m=r}
       \|u\|_{r_0}\prod_{i=1}^m\|\nabla f_i\|_{r_i}.
\end{equation}
Both types of kernel also satisfy spatial and transport identities.
For time-dependent coefficients and a divergence-free field $z$, set
$\mathcal A=\partial_t+z\cn$. Then
\begin{align}
 \partial_kT_{\boldsymbol\gamma,\boldsymbol f}
 &=T_{\boldsymbol\gamma,\boldsymbol f}\partial_k
   +\sum_iT_{\boldsymbol\gamma,(f_1,\ldots,\partial_kf_i,\ldots,f_m)},
                                      \label{high:cz-spatial-rule}\\
 [\mathcal A,T_{\boldsymbol\gamma,\boldsymbol f}]
 &=\sum_{k=1}^3T_{\boldsymbol\gamma+e_k,(\boldsymbol f,z_k)}
   +\sum_iT_{\boldsymbol\gamma,(f_1,\ldots,\mathcal Af_i,\ldots,f_m)}.
                                      \label{high:cz-closure}
\end{align}
The constant depends only on $T,m,r,\alpha$. The identities hold for
smooth functions and extend in distributions under the regularity
assumptions of the displayed estimate.
\end{lemma}
\begin{proof}
We first treat the order-minus-one kernels by absolute integration.
Their singular part has degree
$-2$, so, for $|\boldsymbol\gamma|=m$,
\[
 \left|\partial^{\boldsymbol\gamma}\mathcal K(x-y)
       \prod_{i=1}^m(f_i(x)-f_i(y))\right|
 \le C|x-y|^{-2}\prod_{i=1}^m\|\nabla f_i\|_0.
\]
The distributional product agrees with this integrable kernel.
Indeed, integration by parts outside a ball of radius $\epsilon$
leaves boundary terms of size $O(\epsilon)$, since the product of
differences vanishes to order $m$. The smooth periodic part satisfies
the same bound. Differentiating at fixed $x-y$ differentiates only
the function and the coefficient differences. Leibniz' rule therefore
proves \eqref{high:integrable-difference-estimate}. We prove the
identities for both kernel types in Step~2, after establishing the
H\"older estimate for order-zero kernels.

\emph{1. Cancellation and the zeroth-order estimate.}
For a homogeneous polynomial $P$ of degree $m=|\boldsymbol\gamma|$,
integration by parts on an annulus gives
\begin{equation}\label{high:cz-spherical-moments}
 \int_{\mathbb S^2}P(\omega)\partial^{\boldsymbol\gamma}K_0(\omega)
                       \,\dd S
 =(-1)^m(\partial^{\boldsymbol\gamma}P)
                   \int_{\mathbb S^2}K_0(\omega)\,\dd S=0.
\end{equation}
At each integration by parts the boundary integrand has degree
$-2$, so the inner and outer fluxes agree. The volume integrands
have degree $-3$ and hence the same logarithmic radial integral.
After $m$ integrations, the zero spherical mean of the original
kernel proves the identity.

We prove the estimate for $r=0$ by induction on $m$. When $m=0$,
$T1$ is constant, and the argument below applies without coefficient
differences. For the induction step, choose $j$ with $\gamma_j>0$
and integrate the kernel distribution by parts:
\[
 T_{\boldsymbol\gamma,\boldsymbol f}1
 =-\sum_{i=1}^m
 T_{\boldsymbol\gamma-e_j,\boldsymbol f\setminus f_i}(\partial_jf_i),
 \qquad
 \|T_{\boldsymbol\gamma,\boldsymbol f}1\|_\alpha
       \le C\prod_i\|\nabla f_i\|_\alpha.
\]
The minus sign follows from
\[
 \partial_{w_j}(f_i(x)-f_i(x-w))=\partial_jf_i(x-w).
\]
The bound for $T_{\boldsymbol\gamma,\boldsymbol f}1$ follows by applying
the induction hypothesis to each summand. For a general function, we
subtract its value at $x$ and write
\begin{equation}\label{high:cz-pv-split}
 T_{\boldsymbol\gamma,\boldsymbol f}u(x)
 =\int k(x,y)(u(y)-u(x))\,\dd y
                       +u(x)T_{\boldsymbol\gamma,\boldsymbol f}1(x).
\end{equation}
For smooth $u$, the factor paired with the kernel in the first term
vanishes through degree $m$; hence every diagonal term is zero.
The integral is also absolutely convergent for $u\in C^\alpha$, as
\[
 |k(x,y)|\le C|x-y|^{-3}\prod_i[f_i]_1,
 \qquad
 \int_0^1 t^{\alpha-1}\,\dd t=\alpha^{-1}.
\]
This proves its supremum bound.

For the H\"older seminorm, let $d=|x-x'|$. The contributions from
the union of $B(x,2d)$ and $B(x',2d)$ are bounded by absolute
integration:
\[
 Cd^\alpha\|u\|_\alpha\prod_i\|\nabla f_i\|_\alpha.
\]
On the complement, we separate the difference of the kernels from
the difference of the functions and estimate the first term as follows:
\[
 \begin{aligned}
 &k(x,y)(u(y)-u(x))-k(x',y)(u(y)-u(x'))\\
 &\quad=(k(x,y)-k(x',y))(u(y)-u(x))
                     +k(x',y)(u(x')-u(x)),\\
 |k(x,y)-k(x',y)|
 &\le Cd|x-y|^{-4}\prod_i[f_i]_1,
 \qquad d\int_{2d}^1t^{\alpha-2}\,\dd t\le C_\alpha d^\alpha.
 \end{aligned}
\]
The second term contains a truncated kernel integral. To bound it
uniformly, Taylor expansion at $x'$ gives
\[
 \prod_i(f_i(x')-f_i(y))
 =\prod_i\bigl(\nabla f_i(x')\cdot(x'-y)\bigr)
 +O\!\left(|x'-y|^{m+\alpha}\prod_i\|\nabla f_i\|_\alpha\right).
\]
By \eqref{high:cz-spherical-moments}, the product of the polynomial
and the homogeneous kernel derivative integrates to zero on each
centered annulus. The remainder is integrable.
Replacing the two excluded balls by a centered ball changes the
integral only on an annulus whose two radii are comparable to $d$.
The absolute integral of the kernel on that annulus is uniformly
bounded. Using
\[
 |u(x')-u(x)|\le d^\alpha\|u\|_\alpha,
\]
we obtain the required H\"older estimate for the second term. The
smooth periodic correction satisfies the same estimate, as does the
multiplication term in \eqref{high:cz-pv-split}. This completes the
proof of \eqref{high:cz-kernel-estimate} when $r=0$.

\emph{2. Spatial and transport derivatives.}
At fixed $w$, a spatial derivative in
\eqref{high:cz-difference-kernel} acts on the function or one
coefficient difference. This proves \eqref{high:cz-spatial-rule}; its $r$-fold
iteration and the estimate just proved give
\[
 \|T_{\boldsymbol\gamma,\boldsymbol f}u\|_{r+\alpha}
 \le C\sum_{r_0+\cdots+r_m=r}
      \|u\|_{r_0+\alpha}\prod_i\|\nabla f_i\|_{r_i+\alpha}.
\]
For the commutator, integrate by parts in the term in which the
transport derivative acts on the function. Since $\ddiv z=0$, the
kernel and the coefficient differences are acted on by
$\partial_t+z(x)\cn_x+z(y)\cn_y$. We have
\[
 \begin{aligned}
 (\partial_t+z(x)\cn_x+z(y)\cn_y)(f_i(x)-f_i(y))
   &=\mathcal Af_i(x)-\mathcal Af_i(y),\\
 (z(x)\cn_x+z(y)\cn_y)\partial^{\boldsymbol\gamma}\mathcal K(x-y)
   &=\sum_k(z_k(x)-z_k(y))
                         \partial^{\boldsymbol\gamma+e_k}\mathcal K(x-y).
 \end{aligned}
\]
Leibniz' rule gives \eqref{high:cz-closure}, including the diagonal
part of the kernel. To extend the identities to the stated regularity,
take smooth approximations and use the uniform estimates. Formula
\eqref{high:cz-pv-split} identifies their limits in distributions,
so this passage does not require convergence in the full $C^\alpha$
norm.
\end{proof}

\begin{proposition}[Tensor commutators]
\label{high:cz-iterated}
Let $z_1,\ldots,z_J$ be smooth divergence-free fields, and set
$\mathcal A_i=c_i\partial_t+z_i\cn$, where $c_i$ is constant. For two fixed Cartesian tensor types
$E,F$, write $\mathscr A_i^E=c_i\partial_t+\mathcal L_{z_i}|_E$ and
$\mathscr A_i^F=c_i\partial_t+\mathcal L_{z_i}|_F$.
Let $T$ act componentwise between these types. For a finite sequence
$\iota=(i_1,\ldots,i_m)$ define $S_\varnothing=T$ and successively
\[
 S_{(i_1,\ldots,i_a)}
 =\mathscr A_{i_a}^F S_{(i_1,\ldots,i_{a-1})}
             -S_{(i_1,\ldots,i_{a-1})}\mathscr A_{i_a}^E.
\]
The products of Lie derivatives are ordered as
$\mathscr A_\iota=\mathscr A_{i_m}\cdots\mathscr A_{i_1}$,
and analogously for $\mathcal A_\iota$;
$\iota|P$ keeps the entries with positions in $P$ in this order.
Then
\begin{equation}\label{high:cz-ordered-expansion}
 \mathscr A_\iota^F(Tu)
 =\sum_{P\subset\{1,\ldots,m\}}
                  S_{\iota|P}\mathscr A_{\iota|P^c}^E u.
\end{equation}
For an order-zero $T$, $m\ge1$ and every integer $r\ge0$, one has
\begin{equation}\label{high:cz-iterated-estimate}
 \|S_\iota u\|_{r+\alpha}
 \le C\sum_{p=1}^m
 \sum_{\substack{|\nu_1|+\cdots+|\nu_p|=m-p\\r_0+\cdots+r_p=r}}
 \sum_{b_1,\ldots,b_p=1}^J
 \|u\|_{r_0+\alpha}\prod_{a=1}^p
                  \|\DD\mathcal A_{\nu_a}z_{b_a}\|_{r_a+\alpha}.
\end{equation}
For $T=\mathcal R$ or $\curl^{-1}$, the same expansion and estimate
hold with every $r_a+\alpha$ replaced by $r_a$, including $a=0$.
The entries of each sequence $\nu_a$ belong to $\{1,\ldots,J\}$;
empty sequences are allowed. Constants depend only on the fixed
derivative orders, $T$, $\alpha$, and the tensor types. No commutation of the fields is assumed.
\end{proposition}

For scalars this is the commutator estimate for transport derivatives.
The choice $(c_i,z_i)=(1,0)$ includes an ordinary time derivative,
while $(1,v)$ and $(0,B)$ give the material and magnetic Lie
derivatives, respectively.
\begin{proof}
For every operator $S$,
\[
 \mathscr A_i^F(Su)
 =\bigl(\mathscr A_i^FS-S\mathscr A_i^E\bigr)u
                                      +S\mathscr A_i^Eu.
\]
Repeated application gives \eqref{high:cz-ordered-expansion}: the
positions in $P$ correspond to commutators, and those in $P^c$ to
derivatives of the function. No operators are interchanged, so the
order within each subsequence is preserved.

Write $\mathscr A_i^E=\mathcal A_i+Z_i^E$ and
$\mathscr A_i^F=\mathcal A_i+Z_i^F$. The matrices $Z_i$ are linear
in $\DD z_i$; on vectors and one-forms they are respectively
$-\DD z_i$ and $(\DD z_i)^T$. Hence the first commutator is
\begin{equation}\label{high:cz-natural-first}
 [\mathcal A_i,T]+Z_i^FT-TZ_i^E.
\end{equation}
More generally, each term applied to $u$ has the form
\[
 a\,T_{\boldsymbol\gamma,\boldsymbol f}(bu),
\]
where $\boldsymbol f$ consists of components of
$\mathcal A_\nu z_j$ and $a,b$ are products of
$\DD\mathcal A_\nu z_j$. We count each such coefficient with degree
$|\nu|+1$. Induction on the number $m$ of commutators gives
\begin{equation}\label{high:cz-count}
 \sum_{\text{coefficient factors}}(|\nu|+1)=m,
 \qquad |\boldsymbol\gamma|=\#\boldsymbol f.
\end{equation}
For the induction step, the next commutator applied to $u$ is
\[
 \begin{aligned}
 &(\mathcal A_i a)\,T_{\boldsymbol\gamma,\boldsymbol f}(bu)
 +a\,[\mathcal A_i,T_{\boldsymbol\gamma,\boldsymbol f}](bu)
 +a\,T_{\boldsymbol\gamma,\boldsymbol f}((\mathcal A_i b)u)\\
 &\qquad+Z_i^Fa\,T_{\boldsymbol\gamma,\boldsymbol f}(bu)
              -a\,T_{\boldsymbol\gamma,\boldsymbol f}(bZ_i^Eu).
 \end{aligned}
\]
For the middle term we use \eqref{high:cz-closure}. For derivatives
of the matrix factors we use
\begin{equation}\label{high:cz-gradient-transport}
 \mathcal A_i\partial_k h
 =\partial_k(\mathcal A_i h)-(\partial_kz_i)\cn h
\end{equation}
on each factor. Differentiating a coefficient increases its degree
by one, and multiplication by a new gradient introduces a factor of
degree one. In \eqref{high:cz-gradient-transport}, the second term
has two gradient factors whose degrees add to that of the first term.
Likewise, in the kernel commutator, each new coefficient difference
is accompanied by one kernel derivative. Thus both identities in
\eqref{high:cz-count} hold with $m$ increased by one. The calculation
uses the full kernel throughout and therefore also includes its
diagonal part.

We apply \eqref{high:cz-kernel-estimate} to each term of the
expansion, using the spatial product rule for the matrix factors.
For $p$ coefficient factors, \eqref{high:cz-count} yields
\[
 \sum_{a=1}^p|\nu_a|=m-p,
 \qquad
 \|a\,T_{\boldsymbol\gamma,\boldsymbol f}(bu)\|_{r+\alpha}
 \le C\sum_{r_0+\cdots+r_p=r}
 \|u\|_{r_0+\alpha}\prod_{a=1}^p
             \|\DD\mathcal A_{\nu_a}z_{b_a}\|_{r_a+\alpha}.
\]
Summation proves \eqref{high:cz-iterated-estimate}. Since $p\ge1$,
each coefficient has at most $m-1$ transport derivatives. The order
of all matrix products, including their position relative to $T$,
has been preserved. For operators of order minus one, use
\eqref{high:integrable-difference-estimate} instead of
\eqref{high:cz-kernel-estimate}. The identical expansion and derivative
count give the integer estimate.
\end{proof}

For the material and magnetic derivatives the first commutators are
\[
 [D_t,T]=\sum_iT_{e_i,(v_i)},\qquad
 [D_B,T]=\sum_iT_{e_i,(B_i)},
\]
and, for example,
\[
 [D_B,[D_t,T]]
 =\sum_iT_{e_i,(D_Bv_i)}
  +\sum_{i,j}T_{e_i+e_j,(v_i,B_j)}.
\]
A magnetic commutation therefore either applies a magnetic derivative
to an existing coefficient or introduces a gradient of $B$. Velocity
gradients arise only from commutations with the material derivative.
Counting these two kinds of commutation separately in
\eqref{high:cz-count} gives the following estimates.

\begin{lemma}[Material and magnetic derivatives of Fourier multipliers]
\label{principal:finite-multiplier}
Fix integers $\bar N,J$, a spatial frequency $\Lambda\ge1$, a
H\"older factor $\varkappa\ge\Lambda^\alpha$, and
positive material and magnetic derivative costs $\mathrm a_t,\mathrm a_B$. Assume, for
$r+k+m\le \bar N+J-1$ and $k+m\le J-1$,
\begin{equation}
 \begin{aligned}
 \|D_t^kD_B^m\nabla v\|_{r+\alpha}
   &\le C\varkappa\Lambda^{r}\mathrm a_t^{k+1}\mathrm a_B^m,\\
 \|D_t^kD_B^m\nabla B\|_{r+\alpha}
   &\le C\varkappa\Lambda^{r}\mathrm a_t^k\mathrm a_B^{m+1}.
 \end{aligned}
\end{equation}
Suppose $[D_t,D_B]=0$ and assume for the tensor under consideration
\begin{equation}
 \|\mathcal D_t^k\mathcal L_B^m u\|_{r+\alpha}
 \le F\varkappa\Lambda^{r}\mathrm a_t^k\mathrm a_B^m,
 \qquad k+m\le J,\quad r+k+m\le \bar N+J.
\end{equation}
For any of the order-zero Fourier multipliers above, between its stated
tensor types,
\begin{equation}\label{high:cz-scaled-output}
 \|\mathcal D_t^k\mathcal L_B^m(Tu)\|_{r+\alpha}
 \le CF\varkappa^{k+m+1}\Lambda^{r}\mathrm a_t^k\mathrm a_B^m,
 \qquad r+k+m\le \bar N+J,\quad k+m\le J.
\end{equation}
If instead $T=\mathcal R$ or $\curl^{-1}$ and the hypotheses hold
in integer norms with every factor $\varkappa$ omitted, then
\begin{equation}\label{high:integrable-hodge-mixed}
 \|\mathcal D_t^k\mathcal L_B^m(Tu)\|_r
 \le CF\Lambda^r\mathrm a_t^k\mathrm a_B^m,
 \qquad r+k+m\le\bar N+J,\quad k+m\le J.
\end{equation}
For the order-zero operators, the choices $\Lambda=\lambda_{q+1}$,
$\varkappa=\lambda_{q+1}^{\alpha}$, $\mathrm a_t=\tau_a^{-1}$ and
$\mathrm a_B=\lambda_\parallel$ yield
\begin{equation}
 \|\mathcal D_t^k\mathcal L_B^m(Tu)\|_{r+\alpha}
 \le CF\lambda_{q+1}^{r+(k+m+1)\alpha}\tau_a^{-k}
       \lambda_\parallel^m.
\end{equation}
The same bounds hold for component material and magnetic derivatives.
For a composition of $k+m$ such derivatives, the tensor is
differentiated at most $k+m$ times in these directions combined,
and a coefficient at most $k+m-1$ times. If an ordinary time
derivative is also applied, it acts on these factors by Leibniz' rule;
the corresponding ordinary time estimates must be assumed separately.
For the classes of Section~\ref{ssec:fixed-scales}, a
$(\Lambda,\rho)$-adapted background through $(N,J)$, with
$\Lambda\le\lambda_{q+1}$, satisfies the hypotheses with
\[
 \varkappa=\max\{\ell^{-1},\Lambda\}^{\alpha}
 \le\lambda_{q+1}^\alpha.
\]
Consequently,
\[
 \begin{aligned}
 T\colon\mathcal C_{N,J}(F;\Lambda,\rho)
 &\longrightarrow
 \mathcal C_{N,J}(CF\lambda_{q+1}^{(J+1)\alpha};\Lambda,\rho),\\
 T\colon\mathcal K_{N,J}(F')
 &\longrightarrow\mathcal K_{N,J}(CF'\lambda_{q+1}^{(J+1)\alpha}).
 \end{aligned}
\]
These mappings give rule (v) of Lemma~\ref{setup:calculus}, since
the H\"older bound controls the corresponding integer norm. The
mapping between the $\mathcal K_{N,J}$ classes follows with
$\varkappa=\ell^{-\alpha}$.
\end{lemma}
\begin{proof}
Assign degree $(1,0)$ to each material derivative and degree
$(0,1)$ to each magnetic derivative. In the proof of
Proposition~\ref{high:cz-iterated}, we can then apply
\eqref{high:cz-count} separately to the two derivative orders.
A difference of $D_t^aD_B^bv$ has degree $(a+1,b)$; a difference
of $D_t^aD_B^bB$ has degree $(a,b+1)$. The same assignment is
made to their gradients. The identities
\[
 [D_t,\partial_i]=-(\partial_iv)\cn,\qquad
 [D_B,\partial_i]=-(\partial_iB)\cn
\]
preserve these degrees. It follows that in every term of
\eqref{high:cz-ordered-expansion}, the degrees of the tensor and all
coefficients add to $(k,m)$. Distribute $r$ spatial derivatives by
Leibniz' rule. For a term with $p$ coefficient factors,
\eqref{high:cz-kernel-estimate}, the bound for the tensor, and the
bounds for the $p$ gradients give
\[
 CF\varkappa^{p+1}\Lambda^{r}\mathrm a_t^k\mathrm a_B^m,\qquad p\le k+m.
\]
Summing proves \eqref{high:cz-scaled-output}. Each coefficient arises
from at least one commutator and hence has at most $k+m-1$ further
material and magnetic derivatives combined. The sum of its spatial,
material and magnetic orders is at most $r+k+m-1\le\bar N+J-1$,
and its material and magnetic orders sum to at most $J-1$.
These are precisely the coefficient bounds assumed in the statement.
Ordinary time differentiation commutes with $T$; applying it to the expansion differentiates the
same factors by Leibniz' rule. The integer assertion follows from the
integer estimate in Proposition~\ref{high:cz-iterated}, with the
same count of material and magnetic derivatives and no H\"older
factors.
\end{proof}

The hypotheses on $\nabla B$ are needed already for the first
magnetic commutator; a bound only on $D_BB$ would not suffice.
For order-minus-one operators, \eqref{high:integrable-hodge-mixed}
provides the corresponding bound directly in integer norms.

\subsection{Spatial moment kernels}
\label{ssec:spatial-mollification}

We will use spatial mollifiers whose moments vanish up to a fixed
order. The vanishing moments improve both the approximation estimate
and the estimate for the quadratic commutator. We give the construction
and the estimates here, and then estimate material and magnetic
derivatives of the mollified functions.

All spatial norms in this subsection are taken on $\mathbb T^3$, using
periodic extensions in convolution integrals. Fix an integer $M\ge0$. There is a real smooth kernel
$\rho$ supported in the unit ball with
\begin{equation}\label{moll:spatial-moments}
 \int\rho=1,\qquad \int y^\gamma\rho(y)\,\dd y=0
 \quad(1\le |\gamma|\le M).
\end{equation}
For the construction of $\rho$, take a radial smooth unit-mass
bump $\rho_0$ in the unit ball and distinct dilation factors
$0<c_i<1$, for $0\le i\le\lfloor M/2\rfloor$. The Vandermonde
system determines the coefficients in
\[
 \sum_i b_i c_i^{2j}=\begin{cases}1,&j=0,\\0,&1\le j\le\lfloor M/2\rfloor,
                         \end{cases}
 \qquad \rho(y)=\sum_i b_ic_i^{-3}\rho_0(y/c_i).
\]
Each moment of degree $2j$ is multiplied by $\sum_i b_ic_i^{2j}$,
and radial symmetry makes the odd moments vanish. Thus
\eqref{moll:spatial-moments} holds, although the kernel need not be
nonnegative. Its rescalings and the associated convolutions are
\[
 \rho_{\ell}(y)=\ell^{-3}\rho(y/\ell),\qquad
 \mathcal J_{\ell}f(x)=\int\rho_{\ell}(y)f(x-y)\,\dd y,
 \qquad 0<\ell\le1.
\]
Constants below depend on the fixed kernel and derivative orders, not
on $\ell$.

\begin{proposition}[Spatial mollification]
\label{moll:spatial-basic}
For integers $r,a\ge0$ and $1\le d\le M+1$,
\begin{equation}\label{moll:spatial-accuracy}
 \|\mathcal J_\ell f\|_{r+a}\le C\ell^{-a}\|f\|_r,
 \qquad
 \|\mathcal J_\ell f\|_{r+\alpha}\le C\ell^{-\alpha}\|f\|_r,
 \qquad
 \|(\mathcal J_\ell-\IId)f\|_r\le C\ell^d\|f\|_{r+d}.
\end{equation}
For either $d=2$ or $3\le d\le M+1$,
\begin{equation}
 \|\mathcal J_\ell f\,\mathcal J_\ell g-\mathcal J_\ell(fg)\|_r
 \le C\ell^d\sum_{a+b=r+d-2}\|\nabla f\|_a\|\nabla g\|_b.
\end{equation}
The case $d=2$ uses only unit mass and therefore also holds when
$M=0$. The product estimate applies componentwise to tensors;
ordinary time derivatives act on the two factors by Leibniz's rule.
The smoothing and approximation estimates hold with $r$ replaced
by $r+\alpha$. The H\"older form of the quadratic estimate has
$\alpha$ on exactly one norm in each product, summed over the two
choices.
\end{proposition}
\begin{proof}
Commutation with ordinary derivatives and the kernel scaling give
\[
 [\mathcal J_\ell f]_{r+a}
 \le C\|\DD^a\rho_\ell\|_{L^1}\|f\|_r
 \le C\ell^{-a}\|f\|_r.
\]
Summing the lower derivatives gives the full integer-norm bound.
Interpolation between $a=0$ and $a=1$ yields the H\"older factor
$\ell^{-\alpha}$. For the approximation error, the constant Taylor term cancels with
the unsmoothed function. All nonconstant terms below degree $d$
integrate to zero by the moment conditions. The integral remainder
therefore gives
\[
 \begin{aligned}
 (\mathcal J_\ell-\IId)f(x)
 &=\int\rho_\ell(y)\frac{(-1)^d}{(d-1)!}
    \int_0^1(1-s)^{d-1}(y\cn)^df(x-sy)\,\dd s\,\dd y,\\
 \|(\mathcal J_\ell-\IId)f\|_r
 &\le C\int|\rho_\ell(y)|\,|y|^d\,\dd y\,\|f\|_{r+d}
 \le C\ell^d\|f\|_{r+d}.
 \end{aligned}
\]
These prove \eqref{moll:spatial-accuracy}.

Unit mass yields the exact quadratic identity
\begin{equation}\label{moll:quadratic-kernel}
 \begin{aligned}
 &\mathcal J_\ell(fg)(x)-(\mathcal J_\ell f)(x)(\mathcal J_\ell g)(x)\\
 &\quad=\frac12\iint\rho_\ell(y)\rho_\ell(y')
 [f(x-y)-f(x-y')][g(x-y)-g(x-y')]\,\dd y\,\dd y'.
 \end{aligned}
\end{equation}
To estimate \eqref{moll:quadratic-kernel}, we first record a bound for
a product of differences. Write each of $n$ differences as the
integral of a first derivative over the segment joining its arguments.
This gives $n$ factors of the translation variables. Expand the
remaining derivative factors to degree $s-1$ in those variables.
The remainder is bounded by
\begin{equation}
 C\ell^{n+s}\sum_{a_1+\cdots+a_n=s}
                         \prod_{i=1}^n\|\nabla f_i\|_{a_i}.
\end{equation}
This follows from the integral Taylor formula on the segment from
zero to the translation vectors and $s$ applications of Leibniz'
rule. For $s=0$, it is the bound on the product before expansion.

Apply this bound to \eqref{moll:quadratic-kernel} after taking $r$
spatial derivatives, with $n=2$ and $s=d-2$. The subtracted
monomials have degrees between $2$ and $d-1$. Each integral is a
product of kernel moments, at least one of positive degree and all of
degree at most $d-1\le M$. It therefore vanishes. If $d=2$, no
monomial has been subtracted, and unit mass alone suffices. Hence
\[
 \|\mathcal J_\ell f\,\mathcal J_\ell g-\mathcal J_\ell(fg)\|_r
 \le C\ell^d\sum_{a+b=r+d-2}\|\nabla f\|_a\|\nabla g\|_b.
\]
Ordinary time derivatives commute with convolution and act on the
products by Leibniz' rule. For the H\"older bounds, take a spatial
difference of the integral remainder and telescope each product. Each
term contains a H\"older difference of one factor and integer norms
of the others. This proves the remaining assertions with the same
numbers of spatial and ordinary time derivatives.
\end{proof}

\subsection{Differentiated smoothing}

We next estimate material and magnetic derivatives of spatial
convolutions. As in the singular integral commutators, each derivative
of the kernel is multiplied by a coefficient difference. The
resulting kernel has an absolute integral bounded independently of
the smoothing length, which gives the following estimate.

\begin{proposition}[Material and magnetic derivatives of spatial convolutions]
\label{moll:differentiated-spatial}
Let $D_t=\partial_t+v\cn$ and $D_B=B\cn$ commute, with
$\ddiv v=\ddiv B=0$. Fix positive derivative costs $\mathrm a_t,\mathrm a_B$, a frequency
$\Lambda\ge1$ and derivative orders $N,J$. Suppose
\[
 \begin{aligned}
 \|D_t^kD_B^mT\|_r&\le F\Lambda^r\mathrm a_t^k\mathrm a_B^m,\\
 \|D_t^kD_B^m\nabla v\|_r&\le C\Lambda^r\mathrm a_t^{k+1}\mathrm a_B^m,\\
 \|D_t^kD_B^m\nabla B\|_r&\le C\Lambda^r\mathrm a_t^k\mathrm a_B^{m+1}.
 \end{aligned}
\]
The first estimate is assumed for $k+m\le J$, $r+k+m\le N$,
and the last two for $k+m\le J-1$, $r+k+m+1\le N$.
Then the spatial convolution $\mathcal J_\ell$ with the kernel of
\eqref{moll:spatial-moments} satisfies, using only smoothness and unit
mass,
\begin{equation}\label{moll:spatial-retained}
 \|D_t^kD_B^m\mathcal J_\ell T\|_r
 \le CF\Lambda^r\mathrm a_t^k\mathrm a_B^m,\qquad k+m\le J,\quad r+k+m\le N.
\end{equation}
The corresponding Lie derivative estimate replaces $D_t,D_B$ in
the estimates for the tensor and its convolution by
$\mathcal D_t,\mathcal L_B$, with the same component estimates for
the field gradients. The corresponding spatial H\"older estimates
hold with the same H\"older factor. After $h$ ordinary time
derivatives, the same proof gives products in which the ordinary time
orders on the tensor and background gradients sum to $h$; bounds for
these derivatives are assumed separately.
To estimate $m$ magnetic derivatives of the convolution, one needs
at most that many magnetic derivatives of the tensor. In the classes of
Section~\ref{ssec:fixed-scales}, a $(\Lambda,\rho)$-adapted
background therefore gives a map $\mathcal J_\ell$ from
$\mathcal C_{N,J}(F;\Lambda,\rho)$
into $\mathcal C_{N,J}(CF;\Lambda,\rho)$, which is part of rule (vi)
of Lemma~\ref{setup:calculus}.
\end{proposition}
\begin{proof}
For a smooth kernel put
\[
 T^{\rm sm}_{\gamma,\boldsymbol f}u(x)
 =\int\partial^\gamma\rho_\ell(y)
       \prod_i(f_i(x)-f_i(x-y))u(x-y)\,\dd y.
\]
If $D=c\partial_t+z\cn$, where $c$ is constant and $\ddiv z=0$,
integration by parts gives
\begin{equation}
 [D,T^{\rm sm}_{\gamma,\boldsymbol f}]
 =\sum_jT^{\rm sm}_{\gamma+e_j,(\boldsymbol f,z_j)}
  +\sum_iT^{\rm sm}_{\gamma,(f_1,\ldots,Df_i,\ldots,f_p)}.
\end{equation}
The term $[z(x)-z(x-y)]\cn f_i(x-y)$ obtained by transporting a
coefficient difference cancels the corresponding term from integration
by parts. Thus every new kernel derivative is accompanied by a new
coefficient difference. If their number is $p$, applying the mean
value formula to the differences gives $p$ factors of $|y|$.
These compensate for the kernel derivatives, since
\[
 \int|\nabla^p\rho_\ell(y)|\,|y|^p\,\dd y\le C_p.
\]
Count material and magnetic derivatives as in
Lemma~\ref{principal:finite-multiplier}. If the function has degree
$(k_0,m_0)$ and the gradient factors have degrees $(k_i,m_i)$,
including the commutations which introduced them, then
\begin{equation}\label{moll:transport-finite-bound}
 (k_0,m_0)+\sum_i(k_i,m_i)=(k,m),\qquad
 r_0+\sum_ir_i=r.
\end{equation}
The function and coefficient factors are bounded, respectively, by
\[
 F\Lambda^{r_0}\mathrm a_t^{k_0}\mathrm a_B^{m_0},
 \qquad
 \prod_i C\Lambda^{r_i}\mathrm a_t^{k_i}\mathrm a_B^{m_i}.
\]
By \eqref{moll:transport-finite-bound}, their product is bounded by
the right-hand side of \eqref{moll:spatial-retained}. The number of
terms depends only on the derivative orders. Summation proves the
component estimate.

For tensor Lie derivatives write $\mathcal L_z=D_z+Z_z$, where $Z_z$ is
linear in $\nabla z$. The exact identity
\[
 [\mathcal L_z,\rho_\ell*]
 =[D_z,\rho_\ell*]+Z_z(\rho_\ell*)-(\rho_\ell*)Z_z
\]
gives an expansion in the same difference kernels, with products
of differentiated gradients on either side. The count in
Proposition~\ref{high:cz-iterated} is unchanged: the spatial,
material and magnetic orders on a gradient factor sum to at most
$r+k+m-1$, so the hypotheses bound every such factor.
Finally, $\partial_t$ commutes with convolution and differentiates
the factors in this expansion by Leibniz' rule. Taking a spatial
difference of a product gives a H\"older seminorm on one factor at
a time. The assumed ordinary time and H\"older bounds therefore
prove the corresponding assertions.
\end{proof}

\subsection{Smoothing along the velocity and magnetic flows}

Fix a background pair $(v,B)$ satisfying
\begin{equation}
 \partial_tB+[v,B]=0,\qquad [D_t,D_B]=0,\qquad
 [\mathcal D_t,\mathcal L_B]=0.
\end{equation}
All transport and smoothing operators in this subsection refer to
this pair. Let $\Phi_a^t$ and $\Phi_b^B$ be the spacetime flows of
$(1,v)$ and $(0,B)$, respectively. The second flow preserves
physical time. By the induction equation the flows commute wherever
both compositions are defined. Let $\rho$ be a smooth kernel on $\mathbb R$
supported in $[-1,1]$, with unit mass and vanishing moments through
order $P-1$, and put $\rho_s(u)=s^{-1}\rho(u/s)$. For positive
lengths $s_t,s_B$, define
\begin{align}
 \mathcal J_{s_t,s_B}^{D}f
 &=\iint\rho_{s_t}(a)\rho_{s_B}(b)
             f\circ\Phi_a^t\Phi_b^B\,\dd a\,\dd b,
 \\
 \mathcal J_{s_t,s_B}^{\mathcal L}T
 &=\iint\rho_{s_t}(a)\rho_{s_B}(b)
             (\Phi_a^t\Phi_b^B)^*T\,\dd a\,\dd b.
\end{align}
In the second formula, pullback acts on spatial tensor indices;
a spatial covector acquires no $dt$ component. The definitions agree
on overlapping charts, provided that the full flow segments in the
integrals are contained in the domains of the charts.

In a Lagrangian chart adapted to the magnetic field, that field has
components $c^{-1}e_3$. The tensor components $T^\Psi$ therefore
satisfy
\begin{equation}\label{moll:mixed-translation}
 (\mathcal J_{s_t,s_B}^{\mathcal L}T)^\Psi(t,y)
 =\iint\rho_{s_t}(a)\rho_{s_B}(b)
            T^\Psi(t+a,y+bc^{-1}e_3)\,\dd a\,\dd b.
\end{equation}
Indeed, $\Psi$ satisfies $D_t\Psi=0$ and
$D_B\Psi=c^{-1}e_3$. The two flows are thus translations in these
coordinates. Their differentials are the identity, which proves the
formula for each tensor type.

We write $\mathcal J_{s_t,s_B}$ for the component or tensor
smoother associated with the derivative pair under consideration.
The exponentials below denote composition or spatial tensor pullback
along these flows.

\begin{proposition}[Smoothing along the velocity and magnetic flows]
\label{moll:mixed-calculus}
Write $(\mathscr D_t,\mathscr D_B)=(D_t,D_B)$ for components, or
$(\mathcal D_t,\mathcal L_B)$ for tensor Lie derivatives. If
$0\le k_0\le k$ and $0\le m_0\le m$, then
\begin{equation}\label{moll:mixed-transfer}
 \begin{aligned}
 \mathscr D_t^k\mathscr D_B^m\mathcal J_{s_t,s_B}T
 ={}&(-1)^{k-k_0+m-m_0}\iint
 \rho_{s_t}^{(k-k_0)}(a)\rho_{s_B}^{(m-m_0)}(b)\\
 &\qquad\cdot e^{a\mathscr D_t}e^{b\mathscr D_B}
             \mathscr D_t^{k_0}\mathscr D_B^{m_0}T\,\dd a\,\dd b.
 \end{aligned}
\end{equation}
For $\nu=r$ or $r+\alpha$, consequently,
\begin{equation}\label{moll:mixed-overflow}
 \|\mathscr D_t^k\mathscr D_B^m\mathcal J_{s_t,s_B}T\|_\nu
 \le Cs_t^{-(k-k_0)}s_B^{-(m-m_0)}
 \sup_{|a|\le s_t,|b|\le s_B}
 \|e^{a\mathscr D_t}e^{b\mathscr D_B}
        \mathscr D_t^{k_0}\mathscr D_B^{m_0}T\|_\nu.
\end{equation}
For $1\le d\le P$, one also has
\begin{equation}\label{moll:mixed-accuracy}
 \begin{aligned}
 \|(\mathcal J_{s_t,s_B}-\IId)T\|_\nu
 \le C_d\sum_{i+j=d}s_t^is_B^j
 \sup_{|a|\le s_t,|b|\le s_B}
       \|e^{a\mathscr D_t}e^{b\mathscr D_B}
                      \mathscr D_t^i\mathscr D_B^jT\|_\nu.
 \end{aligned}
\end{equation}
The approximation estimate after applying
$\mathscr D_t^k\mathscr D_B^m$ uses material and magnetic derivatives
of the tensor of orders $(k+i,m+j)$, $i+j=d$, respectively,
throughout the indicated flow segments. When these derivatives are
not among the hypotheses, \eqref{moll:mixed-transfer} instead
places derivatives on the kernels, at a cost of $s_t^{-1}$ or
$s_B^{-1}$ per derivative.
\end{proposition}
\begin{proof}
Commutativity gives
\[
 \partial_ae^{a\mathscr D_t}e^{b\mathscr D_B}
 =\mathscr D_te^{a\mathscr D_t}e^{b\mathscr D_B},
\]
and similarly in $b$. Integrate by parts $k-k_0$ times in $a$ and
$m-m_0$ times in $b$, leaving material and magnetic orders
$(k_0,m_0)$, respectively, on the tensor. There are no boundary
terms because the kernels are compactly supported. This proves
\eqref{moll:mixed-transfer}, including its sign. Since
\[
 \|\rho_s^{(h)}\|_1=s^{-h}\|\rho^{(h)}\|_1,
\]
the norm estimate \eqref{moll:mixed-overflow} follows.

For the approximation error, Taylor-expand the map
\[
 \theta\longmapsto
 e^{\theta a\mathscr D_t}e^{\theta b\mathscr D_B}T
\]
at zero to degree $d-1$. Every nonconstant term has positive
degree less than $d$ in at least one of $a,b$ and therefore
integrates to zero by the moment conditions. The difference from
the constant term is consequently the integral of the Taylor
remainder
\[
 \frac1{(d-1)!}\int_0^1(1-\theta)^{d-1}
 e^{\theta a\mathscr D_t}e^{\theta b\mathscr D_B}
 (a\mathscr D_t+b\mathscr D_B)^dT\,\dd\theta.
\]
Expand the $d$th power and apply Minkowski's inequality to obtain
\eqref{moll:mixed-accuracy}. The material and magnetic derivatives
commute with both flows, so the argument also applies to
$\mathscr D_t^k\mathscr D_B^mT$ when the tensor derivatives of
material order $k+i$ and magnetic order $m+j$, with $i+j=d$, are
bounded on the flow segments in the statement. The alternative
estimate follows by applying
\eqref{moll:mixed-transfer} to place the remaining derivatives
on the kernels.
\end{proof}

\paragraph{Ordinary derivatives and domains.}
The ordinary derivatives are controlled by the chart representation.
To estimate spatial order $r$ and ordinary time order $a$, assume
the bounds for the chart and its inverse when the spatial and
ordinary time orders sum to at most $r+a$ and the ordinary time
order is at most $a$, with a factor $\Lambda$ for each spatial
derivative and $\Lambda$ for each ordinary time derivative.
For tensors, assume the same bounds for the frame matrices.
Then the chain rule in \eqref{moll:mixed-translation} gives
\begin{equation}
 \|\partial_t^a\mathcal J_{s_t,s_B}T\|_r
 \le CF\Lambda^{r+a}
\end{equation}
whenever the corresponding bound for the tensor holds throughout
the domain traversed by the flows. For a scalar, taking spatial and ordinary time derivatives whose
orders sum to $r+a$ uses chart derivatives with at most the same
sum of orders. Tensor components also contain frame factors. These
require one additional spatial derivative of the chart, but no
additional ordinary time derivative. In fact, a chain-rule term with
$p$ spatial and ordinary time derivatives combined on the tensor,
including the accompanying chart factors, is bounded by
\[
 F\Lambda^p\prod_{i=1}^p\Lambda^{d_i-1}
 =F\Lambda^{r+a},
 \qquad \sum_i d_i=r+a.
\]
The frame factors satisfy the same identity once their additional
spatial derivative is included. Interpolation requires just one
further spatial derivative.

For the convolution to be defined, the time interval must extend
by $s_t$ in both directions. In chart coordinates the velocity
flow fixes the spatial point, while the magnetic flow changes only
the third coordinate, by at most $Cs_B$. The domain of the background
must contain these full segments. In physical coordinates, the
displacement bound shows that spatial cutoff regions of width
$C(s_t\|v\|_0+s_B\|B\|_0)$ also suffice.
For the Euclidean contravariant two-tensor $\IId$ appearing in the
stress, we have
\[
 \mathcal L_B\IId=-2\sym\nabla B,\qquad
 \mathcal D_t\IId=-2\sym\nabla v.
\]
Consequently, even a volume-preserving chart need not make magnetic
pullback smoothing preserve the identity tensor.

\begin{lemma}[Estimates for smoothed amplitudes]
\label{aniso:amplitude-overflow}
Fix $\Lambda\ge1$, $F>0$, $\mathrm a_t,\mathrm a_B>0$, $0<e\le1$ and an
integer $J\ge0$. Suppose the material and magnetic Lie derivatives of the tensor satisfy the bounds
$F\Lambda^r\mathrm a_t^p\mathrm a_B^m$ through $p+m\le J$ on
every flow segment in the convolution, and the flow and frame maps
have the assumed spatial bounds through order $r$.
Set $s_t=e/\mathrm a_t$ and $s_B=e/\mathrm a_B$. For each estimate
below, assume the bounds for spatial, material and magnetic
derivatives of the tensor of respective orders $(r,p_0,m_0)$, where $m_0=\min(m,J)$ and
$p_0=\min(p,J-m_0)$, together with the spatial composition bounds
in \eqref{moll:mixed-overflow}. Then
\begin{equation}\label{aniso:flow-overflow-bound}
 \|\mathcal D_t^p\mathcal L_B^m\mathcal J_{s_t,s_B}T\|_r
 \le CF\Lambda^r\mathrm a_t^p\mathrm a_B^m e^{-[p+m-J]^+}.
\end{equation}
The same statement holds for components, for assumed H\"older
bounds with their H\"older factor, and for the lossy bounds for
compositions of spatial, material and magnetic derivatives.

The smoothed tensor also yields scalar coefficient estimates in the
same chart $\Psi$. Write $\mathcal J=\mathcal J_{s_t,s_B}^{\mathcal L}$,
with chart components evaluated at $y=\Psi(t,x)$.
Let $g$ be smooth on a fixed compact set containing the normalized
chart components $F^{-1}(\mathcal JT)^\Psi$. If a slow cutoff $\chi$
has bounds $C\Lambda^r\mathrm a_t^p\mathrm a_B^m$ through spatial
order $r$, material order $p$ and magnetic order $m$, set
$a=\sqrt F\,\chi\,g(F^{-1}(\mathcal JT)^\Psi)$. Then
\begin{equation}\label{aniso:amplitude-overflow-bound}
 \|D_t^pD_B^ma\|_r
 \le C\sqrt F\,\Lambda^r\mathrm a_t^p\mathrm a_B^m e^{-[p+m-J]^+}.
\end{equation}
For the square-root maps in the geometric decomposition, the compact
set is assumed to lie where the coefficient functionals are positive.
The chart and cutoff bounds are assumed at all derivative orders in
the conclusion. For the tensor, bounds are required only when the sum
of the numbers of material and magnetic Lie derivatives is at most $J$.
\end{lemma}
\begin{proof}
With $m_0=\min(m,J)$ and $p_0=\min(p,J-m_0)$, we have
$p_0+m_0=\min(p+m,J)$. Formula~\eqref{moll:mixed-overflow}
therefore uses only the assumed derivatives of the tensor. All
remaining material and magnetic derivatives act on the corresponding
kernels and give the factor
\[
 \mathrm a_t^{p_0}\mathrm a_B^{m_0}(\mathrm a_t/e)^{p-p_0}(\mathrm a_B/e)^{m-m_0}
 =\mathrm a_t^p\mathrm a_B^m e^{-[p+m-J]^+},
\]
which, together with the spatial bounds for the flow and frame, proves
\eqref{aniso:flow-overflow-bound}. No higher material or magnetic
Lie derivatives of the tensor occur.

In the chart, the Lie derivatives act on the tensor components as
constant-coefficient derivatives. The chain rule for $g$ thus gives
products of derivatives of the normalized components. Denote by
$d_1,\ldots,d_h$ the respective sums of the material and magnetic
orders on these factors. Then
\[
 \sum_i[d_i-J]^+\le[\sum_i d_i-J]^+.
\]
The material orders on the factors add to the material order being
estimated, and the magnetic orders add to its magnetic order. Since
$0<e\le1$, the displayed inequality bounds the product of the
negative powers of $e$ by the claimed power. Applying Leibniz' rule
to $\chi$ and including $\sqrt F$ proves
\eqref{aniso:amplitude-overflow-bound}. Spatial differentiation uses
the product rule and the assumed chart estimates. The H\"older and
lossy estimates follow from the corresponding product bounds.
\end{proof}

\subsection{Localized antiderivatives}

The next lemma constructs a localized approximate antiderivative
with respect to the material derivative. Successive terms in a finite
sum cancel the derivatives of the slow amplitude, leaving an explicit
error. We estimate the antiderivative both under the original
derivative bounds and under the bounds obtained by placing additional
derivatives on the smoothing kernels.

\begin{lemma}[Localized approximate antiderivative]
\label{principal:material-primitive}
Let $D_t=\partial_t+v\cn$ and $D_B=B\cn$ commute on scalars.
Fix $N\ge1$, $\tau_a>0$ and smooth unit-periodic profiles
$\alpha^{[0]}=\alpha,\ldots,\alpha^{[N]}$ satisfying
$(\alpha^{[j]})'=\alpha^{[j-1]}$ and vanishing in a common open
cutoff region in each period. For a smooth scalar $a$, the profiles
are evaluated at $t/\tau_a$, and we set
\[
 \mathfrak a=\sum_{j=0}^{N-1}(-1)^j\tau_a^{j+1}\alpha^{[j+1]}D_t^ja,
 \qquad a^{\rm c}=(-1)^{N-1}\tau_a^N\alpha^{[N]}D_t^Na.
\]
Then $D_t\mathfrak a=\alpha a+a^{\rm c}$. The supports of
$\mathfrak a,a^{\rm c},D_t\mathfrak a$ stay in
$\supp a$ and outside the common fast-time cutoff region.

Suppose $F>0$, $\Lambda\ge1$, $\mathrm a_t,\mathrm a_B>0$, $\tau_a \mathrm a_t\le1$,
and $\|D_t^pD_B^ma\|_r\le F\Lambda^r\mathrm a_t^p\mathrm a_B^m$
for derivatives of the amplitude through material order $k+N$,
magnetic order $m$ and spatial order $r$. Then
\begin{align}
 \|D_t^kD_B^m\mathfrak a\|_r
 &\le CF\Lambda^r\tau_a^{1-k}\mathrm a_B^m,\notag\\
 \|D_t^kD_B^ma^{\rm c}\|_r
 &\le CF\Lambda^r\tau_a^{-k}\mathrm a_B^m(\tau_a \mathrm a_t)^N,
 \\
 \|D_t^kD_B^m(\mathfrak a-\tau_a\alpha^{[1]}a)\|_r
 &\le CF\Lambda^r\tau_a^{1-k}\mathrm a_B^m(\tau_a \mathrm a_t).
\end{align}
Alternatively, for an integer $J\ge0$, suppose the amplitude satisfies
\begin{equation}\label{principal:amplitude-overflow-input}
 \|D_t^pD_B^ma\|_r
 \le CF\Lambda^r\mathrm a_t^p\mathrm a_B^m
                   (\tau_a \mathrm a_t)^{-[p+m-J]^+}
\end{equation}
through material order $k+N$, magnetic order $m$ and spatial order
$r$. The resulting estimates are
\begin{equation}\label{principal:primitive-overflow}
 \begin{aligned}
 \|D_t^kD_B^m\mathfrak a\|_r
 &\le CF\Lambda^r\tau_a^{1-k}\mathrm a_B^m
                          (\tau_a \mathrm a_t)^{-[m-J]^+},\\
 \|D_t^kD_B^ma^{\rm c}\|_r
 &\le CF\Lambda^r\tau_a^{-k}\mathrm a_B^m
                          (\tau_a \mathrm a_t)^{\min\{N,J-m\}},\\
 \|D_t^kD_B^m(\mathfrak a-\tau_a\alpha^{[1]}a)\|_r
 &\le CF\Lambda^r\tau_a^{1-k}\mathrm a_B^m
                          (\tau_a \mathrm a_t)^{\min\{1,J-m\}}.
 \end{aligned}
\end{equation}
The exponents may be negative when $m>J$. The assertions use only
the stated derivatives of the amplitude and the corresponding
derivatives of the profiles. Directly assumed H\"older bounds retain their factor;
a H\"older bound obtained by interpolation requires one additional
spatial derivative.
\end{lemma}
\begin{proof}
For $1\le j<N$, the coefficient of
$\alpha^{[j]}D_t^ja$ in $D_t\mathfrak a$ is
$(-1)^{j-1}\tau_a^j+(-1)^j\tau_a^j=0$. Hence all interior terms cancel and the remaining terms are
$\alpha a$ and $a^{\rm c}$. Differentiation does not enlarge the support of a smooth function,
so every term is supported in the support of the amplitude. Each term also vanishes in the common fast-time
cutoff region of the profiles. This proves both support assertions.

Since $D_B\alpha^{[j]}=0$, a term with $b$ of the $k$ material
derivatives applied to the amplitude in the $j$th summand of
$\mathfrak a$ is bounded by
\[
 \tau_a^{j+1}\tau_a^{-(k-b)}\mathrm a_t^{j+b}
 =\tau_a^{1-k}(\tau_a \mathrm a_t)^{j+b}.
\]
Spatial and magnetic derivatives contribute $\Lambda^r\mathrm a_B^m$.
For $a^{\rm c}$, the same calculation has $j=N$ and no leading factor
$\tau_a$. Subtracting $\tau_a\alpha^{[1]}a$ removes the summand
with $j=0$. Since $\tau_a \mathrm a_t\le1$, summing these finite
expansions gives the first three estimates.

For the alternative bounds, \eqref{principal:amplitude-overflow-input}
follows from Lemma~\ref{aniso:amplitude-overflow} after composition
with a smooth function and multiplication by a cutoff. Substitution of
\eqref{principal:amplitude-overflow-input} into the same product
calculation gives the exponent
\[
 j+b-[j+b+m-J]^+=\min\{j+b,J-m\}.
\]
Its smallest value is $\min\{0,J-m\}$ for $\mathfrak a$,
$\min\{N,J-m\}$ for $a^{\rm c}$, and $\min\{1,J-m\}$ for the
difference $\mathfrak a-\tau_a\alpha^{[1]}a$. When $N=1$, this
difference vanishes. Summing the finite Leibniz expansion proves
\eqref{principal:primitive-overflow}. The $k$ material derivatives have already contributed the factor
$\tau_a^{-1}$ per derivative in the calculation above; they
therefore cause no further loss in the displayed powers.
The H\"older assertions follow from the H\"older product estimate
under the same material and magnetic derivative hypotheses.
\end{proof}

\subsection{Oscillatory inverse divergence}
\label{sec:inverse-divergence}

\paragraph{The symmetric inverse divergence.}

We begin with the symmetric, trace-free inverse of the tensor
divergence that enters the stress identities of
Section~\ref{sec:adapted-perturbation} and the estimates of
Section~\ref{sec:endpoint-stress}. The inverse Laplacian has Fourier
multiplier $-|k|^{-2}$ for $k\in\mathbb Z^3\setminus\{0\}$ and
annihilates the zero mode; the Leray projection $\mathbb P$ is as in
\eqref{setup:hodge}.
For a smooth mean-zero vector field $f$, put $u=\Delta^{-1}f$ and define
\begin{equation}\label{eq:inverse-divergence}
 \mathcal R f=\frac14\bigl[\DD\mathbb P u+(\DD\mathbb P u)^T\bigr]
       +\frac34\bigl[\DD u+(\DD u)^T\bigr]
       -\frac12(\ddiv u)\IId.
\end{equation}
\begin{lemma}\label{lem:inverse-divergence}
The tensor $\mathcal Rf$ is symmetric and trace free, and
$\ddiv\mathcal Rf=f$.
\end{lemma}
\begin{proof}
The definition gives $(\mathcal R f)^T=\mathcal R f$. Since
$\ddiv\mathbb P u=0$ and $\operatorname{tr}\IId=3$,
\[
 \operatorname{tr}(\mathcal R f)
 =\tfrac12\ddiv\mathbb P u+\tfrac32\ddiv u
        -\tfrac32\ddiv u=0.
\]
For the tensor divergence of Section~\ref{ssec:path-conventions},
$\ddiv(\DD u)=\Delta u$ and
$\ddiv((\DD u)^T)=\nabla\ddiv u$. Using $\Delta u=f$ and the definition
of $\mathbb P$, we obtain
\begin{align*}
 \ddiv\mathcal R f
 &=\tfrac14\Delta\mathbb P u
       +\tfrac34\Delta u+\tfrac14\nabla\ddiv u\\
 &=\tfrac14\mathbb P f+\tfrac34f+\tfrac14(f-\mathbb P f)=f.
\end{align*}
\end{proof}

\paragraph{Oscillatory tensors.}

We seek a small symmetric tensor with the same divergence as a given
oscillatory tensor tangent to the phase level sets. A finite expansion
in the chart gives such a tensor up to a remainder, to which we apply
the inverse divergence above. Fix a volume-preserving chart and a
reference pair $(v,B)$ defining all transport derivatives below.
In the application, the relevant coordinates are $y_{I,1}^{k_I}$
and $y_{I,1}^{\zeta_I}$ for one index $I$ at $s=1$, and the pair is
$(v_{\ell,I,1}^{\rm c},B_{\ell,I,1}^{\rm c})$. The required chart
estimates are those following \eqref{gg:corrected-chart-bounds},
obtained from Lemma~\ref{gg:transport-classes} for the charts before
the corrector pushforward. In the notation of the lemma, the
invariance identities are
\[
 D_ty^k=D_By^k=0,\qquad
 \mathcal D_t\frac{\partial}{\partial y^\zeta}
 =\mathcal L_B\frac{\partial}{\partial y^\zeta}=0.
\]
Thus material and magnetic derivatives do not differentiate the
oscillatory phase. We assume ordinary time estimates separately.
We also assume global estimates for the reference gradients, since
the inverse divergence of the remainder is nonlocal.

\begin{lemma}[Symmetric inverse divergence for oscillatory tensors]
\label{osc:parametrix}
Let $H$ be smooth, $2\pi$-periodic and mean zero, and let $f$ have
compact support in a chart domain. Set
\[
 Q=f(\partial_{y^\zeta}\otimes\partial_{y^\zeta})
                              H(\lambda_{q+1}y^k).
\]
Suppose that, on the range $r+k+m\le N$, $k+m\le j_1$,
\begin{equation}
 \|D_t^kD_B^mf\|_r
 \le CF\ell^{-r}\tau_a^{-k}
          \lambda_\parallel^m.
\end{equation}
Assume also the independent ordinary time bounds
\[
 \|\partial_t^af\|_r\le CF\ell^{-r-a},
 \qquad a\le j_1,\quad r+a\le N.
\] The frame and coframe are annihilated by both Lie derivatives
and satisfy the ordinary time estimates
\[
 \|\partial_t^a\frac{\partial}{\partial y^\zeta}\|_r
 +\|\partial_t^a dy^k\|_r\le C\ell^{-r-a},
 \qquad a\le j_1,\quad r+a\le N,
 \qquad |dy^k|\ge c>0.
\]
The reference fields are smooth and divergence free on the whole
torus, satisfy $\partial_tB+[v,B]=0$, and are bounded in $C^0$.
On a chart neighborhood of the support their gradients obey
\[
 \begin{aligned}
 \|D_t^kD_B^m\nabla v\|_r
 &\le C\ell^{-r}\tau_a^{-k-1}
       \lambda_\parallel^m,\\
 \|D_t^kD_B^m\nabla B\|_r
 &\le C\ell^{-r}\tau_a^{-k}
       \lambda_\parallel^{m+1},
 \end{aligned}
 \qquad k+m\le j_1-1,\quad r+k+m\le N-1.
\]
The reference fields have independent ordinary component bounds
$C\ell^{-r-a}$ for $a\le j_1-1$, $r+a\le N$. Suppose also that
the same bounds for material and magnetic derivatives of the gradients hold on the whole torus with
$\ell^{-r}$ replaced by $\lambda_{q+1}^r$, for $k+m\le j_1-1$
and $r+k+m\le N-1$.
The bounds for $f$, the frame and coframe, and the ordinary
component bounds for the reference fields are required only on the
specified chart neighborhood of the support. Every product containing
the amplitude has compact support there and extends smoothly by zero.
Choose an integer $\rfg$ for which
\begin{equation}\label{osc:reserve}
 \lambda_{q+1}(\ell\lambda_{q+1})^{-\rfg}\le1.
\end{equation}
There is a symmetric $S$ such that
\begin{equation}\label{osc:exact-divergence}
 \ddiv S=\ddiv Q,
\end{equation}
and, for $k+m\le j_1$ and $r+k+m+\rfg+3\le N$,
\begin{equation}\label{osc:integer-estimate}
 \|\mathcal D_t^k\mathcal L_B^mS\|_r
 \le CF(\ell\lambda_{q+1})^{-1}\lambda_{q+1}^r
       \tau_a^{-k}\lambda_\parallel^m.
\end{equation}
The independent ordinary time estimate is
\[
 \|\partial_t^aS\|_r
 \le CF(\ell\lambda_{q+1})^{-1}\lambda_{q+1}^{r+a},
 \qquad a\le j_1,\quad r+a+\rfg+3\le N.
\]
If the ordinary bounds for the amplitude, frame and coframe also hold
for $a\le j_1+1$, $r+a\le N$, and the reference fields have the
stated ordinary component bounds for $a\le j_1$, $r+a\le N$, then
the same ordinary conclusion holds for $a\le j_1+1$ and
$r+a+\rfg+3\le N$. This does not extend the material or magnetic
derivative range. For the continuation to higher spatial orders
below, the corresponding higher-order ordinary hypotheses must hold
through time order $j_1$, or through $j_1+1$ for this conditional
extension.
If the preceding hypotheses also hold with $r$ replaced by $r+1$,
one also has
\begin{equation}\label{osc:holder-estimate}
 \|\mathcal D_t^k\mathcal L_B^mS\|_{r+\alpha}
 \le CF(\ell\lambda_{q+1})^{-1}\lambda_{q+1}^{r+\alpha}
       \tau_a^{-k}\lambda_\parallel^m.
\end{equation}
If the coefficients satisfy the $\mathcal K_{N',j_1}$ hypotheses
and the inequality in Lemma~\ref{app:finite-spatial-reserve}, with
the amplitudes and spatial range specified there, then the same
$S$ satisfies the following estimates at each spatial order for
which the global reference gradients obey the corresponding material
and magnetic derivative bounds with spatial factor
$\lambda_{q+1}^r$.
Alternatively, assume all the coefficient, frame and gradient bounds
above with the material and magnetic derivative factors replaced by
$\lambda_{q+1}\delta_{q+1}^{1/2}$ and
$\lambda_{q+1}\delta_{B,q+1}^{1/2}$, respectively. These bounds refer to the same material and magnetic
derivatives, and the spatial factor for the coefficients remains
$\ell^{-r}$. Under these alternative hypotheses the conclusion
holds for $r+k+m+\rfg+3\le N$. Extending it beyond this range
requires the separate higher-order estimates and the same inequality
as above. In both cases,
\begin{equation}\label{osc:final-continuation}
 \begin{aligned}
 \|\mathcal D_t^k\mathcal L_B^mS\|_r
 &\le CF(\ell\lambda_{q+1})^{-1}\lambda_{q+1}^r
   (\lambda_{q+1}\delta_{q+1}^{1/2})^k
                         (\lambda_{q+1}\delta_{B,q+1}^{1/2})^m,\\
 \|\partial_t^aS\|_r
 &\le CF(\ell\lambda_{q+1})^{-1}\lambda_{q+1}^{r+a}.
 \end{aligned}
\end{equation}
Without the hypotheses for spatial order $r+1$, only the integer
estimate at spatial order $r$ is asserted. The tensor vanishes
at every time when $f$ vanishes identically in space. The conclusion also holds for a
sum with bounded overlap whose coefficients satisfy the estimates with
respect to the same reference transports.
\end{lemma}
\begin{proof}
\emph{1. Exact construction.}
The chart identities give
\[
 |\nabla y^k|\ge c>0,\qquad
 \ddiv\frac{\partial}{\partial y^\zeta}=0,\qquad
 \partial_{y^\zeta}y^k=0.
\] For a vector $u$, put
\begin{equation}\label{osc:algebraic-inverse}
 \mathcal S_{\nabla y^k}u
 =\frac{u\otimes\nabla y^k+\nabla y^k\otimes u}{|\nabla y^k|^2}
 -\frac{u\cdot\nabla y^k}{|\nabla y^k|^4}
                         \nabla y^k\otimes\nabla y^k.
\end{equation}
This symmetric tensor is a right inverse to contraction with
$\nabla y^k$, since
\[
 (\mathcal S_{\nabla y^k}u)\nabla y^k=u.
\]
The tangency of $\frac{\partial}{\partial y^\zeta}$ to the level
sets of $y^k$ removes the derivative of the oscillatory profile
from the divergence, giving
\begin{equation}\label{osc:tangent-divergence}
 \ddiv Q=F_0H(\lambda_{q+1}y^k),\qquad
 F_0=(\partial_{y^\zeta}f)\frac{\partial}{\partial y^\zeta}
       +f(\frac{\partial}{\partial y^\zeta}\cn)\frac{\partial}{\partial y^\zeta}.
\end{equation}
Let $H_0=H$ and let $H_{j+1}$ be the mean-zero periodic primitive
of $H_j$. Define $S_j=\mathcal S_{\nabla y^k}F_j$ and
$F_{j+1}=\ddiv S_j$. Set
\begin{equation}\label{osc:finite-sum}
 T_{\rfg}=\sum_{j=0}^{\rfg-1}
 (-1)^j\lambda_{q+1}^{-j-1}S_jH_{j+1}(\lambda_{q+1}y^k).
\end{equation}
The product rule gives
\[
 \ddiv(\lambda_{q+1}^{-j-1}S_jH_{j+1})
 =\lambda_{q+1}^{-j}F_jH_j
   +\lambda_{q+1}^{-j-1}F_{j+1}H_{j+1},
\]
where all profiles are evaluated at $\lambda_{q+1}y^k$.
The alternating signs cancel consecutive terms, leaving
\begin{equation}
 \ddiv Q-\ddiv T_{\rfg}
 =(-1)^{\rfg}\lambda_{q+1}^{-\rfg}
          F_{\rfg}H_{\rfg}(\lambda_{q+1}y^k)
 =e_{\rfg}.
\end{equation}
All local products extend smoothly by zero, so $e_{\rfg}$ is a
difference of divergences on the torus and has zero spatial mean.
We may therefore apply the inverse divergence to it. The tensor
\begin{equation}
 S=T_{\rfg}+\mathcal R e_{\rfg}
\end{equation}
is symmetric with divergence given by \eqref{osc:exact-divergence}.

\emph{2. Estimates for material and magnetic Lie derivatives.}
The induction equation implies $[D_t,D_B]=0$ and commutation of the
corresponding Lie derivatives. Their commutators with spatial
derivatives and their action on tensor indices are controlled by the
gradient hypotheses. In particular, the magnetic derivative requires
the full gradient of $B$. We use the following identity for a
divergence-free vector $z$ and a symmetric tensor $U$:
\begin{equation}\label{osc:divergence-commutator}
 \mathcal L_z(\ddiv U)-\ddiv(\mathcal L_zU)
                              =U_{ij}\partial_i\partial_jz.
\end{equation}
To verify it, expand in Cartesian components. The terms with one
derivative on each factor cancel, and those containing
$\partial_j\ddiv z$ vanish. The remaining term is the contraction
of $U$ with the Hessian of $z$. Taking $z=B$ gives the magnetic
identity; the material identity follows because $\partial_t$
commutes with divergence. Although $dy^k$ is invariant under both
Lie derivatives, raising its index in \eqref{osc:algebraic-inverse}
also differentiates the Euclidean metric. This produces factors of
$\nabla v$ or $\nabla B$, bounded at the respective material and
magnetic derivative costs. Apply the product rule to the algebraic
inverse and use \eqref{osc:divergence-commutator} at each recursive
divergence. The coefficient $F_0$ requires one spatial derivative
of the amplitude or frame, and each subsequent divergence requires
one more. Induction on $j$ gives
\begin{equation}\label{osc:recursive-coefficients}
 \begin{aligned}
 &\|\mathcal D_t^k\mathcal L_B^mS_j\|_r
    +\ell\|\mathcal D_t^k\mathcal L_B^mF_{j+1}\|_r\\
 &\qquad\le CF\ell^{-r-j-1}\tau_a^{-k}
       \lambda_\parallel^m.
 \end{aligned}
\end{equation}
The ordinary time estimate, obtained separately by Leibniz' rule,
has factor $\ell^{-r-a-j-1}$. At step $j$ we have used $j+1$
additional spatial derivatives of the coefficients, without increasing
either the material or the magnetic derivative order. The three
additional derivatives in the hypothesis allow for the initial
divergence, the spatial derivatives of the background gradients in
the commutators, and one spatial derivative for interpolation.

Since $D_ty^k=D_By^k=0$, no derivative of the phase occurs when
we apply material or magnetic Lie derivatives. A spatial derivative
costs at most $\lambda_{q+1}$ by the chain rule. An ordinary time
derivative has the same cost, by the independent ordinary estimates
and the identity $\partial_ty^k=-v\cn y^k$. Thus the coefficient
bound gives, for the $j$th term of \eqref{osc:finite-sum},
\[
 CF(\ell\lambda_{q+1})^{-j-1}\lambda_{q+1}^r
     \tau_a^{-k}\lambda_\parallel^m.
\]
By the assumed inequality, $(\ell\lambda_{q+1})^{-1}$ is at most
one. The finite sum is therefore bounded by the first term times a
constant depending on $\rfg$, which proves the estimates for
$T_{\rfg}$.

From \eqref{osc:recursive-coefficients} and phase invariance, the
remainder $e_{\rfg}$ of the finite expansion satisfies
\[
 \|\mathcal D_t^k\mathcal L_B^m e_{\rfg}\|_r
 \le CF\ell^{-\rfg-1}\lambda_{q+1}^{r-\rfg}
          \tau_a^{-k}\lambda_\parallel^m.
\]
Apply \eqref{high:integrable-hodge-mixed} with
$\Lambda=\lambda_{q+1}$, $\mathrm a_t=\tau_a^{-1}$ and
$\mathrm a_B=\lambda_\parallel$, taking amplitude
$CF\ell^{-\rfg-1}\lambda_{q+1}^{-\rfg}$ for the remainder.
The global gradient hypotheses give the coefficient bounds required
by this nonlocal commutator estimate. Hence
\[
 \|\mathcal D_t^k\mathcal L_B^m(\mathcal R e_{\rfg})\|_r
 \le CF\ell^{-\rfg-1}\lambda_{q+1}^{r-\rfg}
          \tau_a^{-k}\lambda_\parallel^m.
\]
We verify the derivative ranges in this application. The spatial
order of the remainder and the sum of its material and magnetic
orders are bounded by $(r,k+m)$, respectively. The recursion supplies
these bounds whenever $r+k+m+\rfg+3\le N$. Each coefficient in
the commutator has at most $k+m-1\le j_1-1$ material and magnetic
derivatives combined. Including its spatial derivatives, a background
gradient is differentiated at most $r+k+m-1\le N-1$ times. Thus
the global hypotheses suffice. The ratio of this remainder bound to
\eqref{osc:integer-estimate} is
$\lambda_{q+1}(\ell\lambda_{q+1})^{-\rfg}$, which is at most one by
\eqref{osc:reserve}. For ordinary time derivatives,
\eqref{high:integrable-hodge} applies directly because the kernel
commutes with time differentiation. The independent bound on
$e_{\rfg}$ is $CF\ell^{-\rfg-1}\lambda_{q+1}^{r+a-\rfg}$, whose
ratio to the asserted ordinary time estimate is again bounded by
\eqref{osc:reserve}. For the conditional ordinary endpoint assertion,
the same spatial recursion applies. Ordinary time derivatives of the
amplitude, frame and coframe are used through order $a$, while
$\partial_ty^k=-v\cn y^k$ requires ordinary derivatives of the
reference velocity only through order $a-1$. Thus the extended
ordinary hypotheses give the estimate through $a\le j_1+1$,
without changing the material or magnetic derivative restrictions.
The inverse divergence commutes with ordinary time differentiation,
so the remainder has the same extended range. Interpolation between
consecutive integer spatial norms proves \eqref{osc:holder-estimate}.

\emph{3. The remaining derivative estimates.}
Divide $S_j$ by $F\ell^{-j-1}$ and apply
Lemma~\ref{app:finite-spatial-reserve} to the resulting coefficient
bounds. The lemma also applies to
$F_{\rfg}/(F\ell^{-\rfg-1})$. The recursion uses $j+1$ additional
spatial derivatives for $S_j$ and $\rfg+1$ for $F_{\rfg}$.
The assumed higher-order bounds and the inequality in that lemma
therefore give the asserted estimates. For ordinary time derivatives, use part (c) with the same time index as
in the higher-order ordinary hypotheses, including the extended
index in the conditional endpoint assertion.
For the remainder, apply \eqref{high:integrable-hodge-mixed} with
the higher-order coefficient bounds and the global gradient bounds,
using the same material and magnetic derivative factors. The global
gradients have spatial factor $\lambda_{q+1}^r$. The sum of their
material and magnetic orders is at most one less than the sum for
the tensor being estimated. The inequality also absorbs the
additional coefficient factor beyond the sharp spatial range.
The remaining factor is
\[
 \lambda_{q+1}(\ell\lambda_{q+1})^{-\rfg},
\]
which is at most one by \eqref{osc:reserve}, independently of
the material and magnetic orders. This proves
\eqref{osc:final-continuation} under the stated higher-order
hypotheses. The global gradients need no spatial bound with factor
$\ell^{-r}$.

Under the alternative derivative hypotheses, repeat Step 2 with
those factors in each coefficient estimate. The construction itself
does not change: it gives the same tensors $S_j$, the same finite
sum, and the same remainder $e_{\rfg}$. The reference transports
still annihilate the phase. The global gradient hypotheses have the
same derivative factors, so \eqref{osc:reserve} again bounds the
remainder. This proves the alternative assertion for $k+m\le j_1$,
$r+k+m+\rfg+3\le N$, and for the ordinary time orders specified in
the statement.

For a family with bounded overlap, sum the local tensors $T_{\rfg}$
and remainders $e_{\rfg}$ first, and apply $\mathcal R$ to the sum
of the remainders. Bounded overlap controls both sums, and the same
nonlocal estimate then applies. Finally, the construction uses only
spatial operations at each fixed time. If $f$ is identically zero
in space at that time, every recursive coefficient is zero, and so
are $e_{\rfg}$ and $S$. This proves the assertion on time support.
\end{proof}
\section{Forced Galbrun Equation}
\label{sec:galbrun}

Fix a smooth background $(v,B,p,R)$ satisfying \eqref{eq:mhd}.
The Galbrun equation is the momentum equation linearized in the Lie
algebra of volume-preserving diffeomorphisms. For a deformation
with initial LDF $\xi=\xi_0$, the tangent fields are
\[
 w=\left.\partial_sv_s\right|_{s=0}=\mathcal D_t\xi,\qquad
 b=\left.\partial_sB_s\right|_{s=0}=\mathcal L_B\xi,
\]
by \eqref{eq:path-linear-fields}, and satisfy the linearized induction
equation. By \eqref{eq:displacement-force}, the linearized momentum
equation with prescribed force is
\[
 \mathcal F_{v,B}(w,b)+\nabla\pi
 =(D_t^2-D_B^2)\xi-\DD(D_tv-D_BB)\xi+\nabla\pi
 =\ddiv\mathsf F.
\]
The same displacement determines both field increments. The principal
part of its equation is $D_t^2-D_B^2$, with $D_t$ the material
derivative and $D_B$ the magnetic derivative.

To impose $\ddiv\xi=0$, we write $\xi=\curl\Theta$. We solve an equation
for the one-form potential $\Theta$ and then determine the pressure by
solving a Poisson equation. The estimates use an explicit approximate
solution obtained by integrating the oscillatory profiles in time.
Subtracting it leaves a small remainder. Since the explicit potential
is supported in the time support of the source coefficients, this
decomposition also improves the estimates in the cut off regions.

All operators below refer to the fixed background. In the main
construction we apply the result on each background interval and sum
the resulting localized potentials $\Theta^{\rm c}$ and cutoff stresses
$R^{\mathrm{cut}}$.

We use the transport derivatives $D_t,D_B$ on components and the
Lie derivatives $\mathcal D_t,\mathcal L_B$ on tensors. The induction
equation gives $[D_t,D_B]=0$ on components and
$[\mathcal D_t,\mathcal L_B]=0$ on tensors. Although the
characteristic system is written with
$\mathcal A^\pm=D_t\pm D_B$, we estimate material and magnetic derivatives separately.

\subsection{The linearized equation and solvability}

The operator on Cartesian one-form components is
\begin{equation}\label{galbrun:operator}
 \mathscr G\Theta=(D_t^2-D_B^2)\Theta
       +\mathscr H_tD_t\Theta+\mathscr H_BD_B\Theta+\mathscr H_2\Theta,
\end{equation}
with the lower-order operators defined in
\eqref{galbrun:separate-lower-operators} and \eqref{mom:H2} below.
Put $z^\pm=v\pm B$ and use the operators of \eqref{setup:hodge}.
For a one-form $Y$ and a vector field $z$ we write $Y\times\nabla z$
for the matrix
\begin{equation}\label{mom:Kdefinition}
 (Y\times\nabla z)_{ij}=(Y\times\nabla z_i)_j
 =\sum_{k,l}\varepsilon_{jkl}Y_k\partial_lz_i,
\end{equation}
whose $i$th row is the cross product of $Y$ with the gradient of the
component $z_i$. The vector field $z$ is differentiated once, and
no derivative falls on $Y$. Antisymmetry of $\varepsilon_{jkl}$ gives
\begin{equation}\label{mom:Kdiv}
 \ddiv(Y\times\nabla z)=(\curl Y)\cn z.
\end{equation}
Indeed, differentiating $Y$ gives the right-hand side, while the term
containing the Hessian of $z_i$ vanishes by symmetry against
$\varepsilon_{jkl}$.

\begin{proposition}[The potential equation]\label{mom:lower-structure}
Define the lower-order operators in \eqref{galbrun:operator} by
\begin{equation}\label{galbrun:separate-lower-operators}
 \begin{aligned}
 \mathscr H_tY&=2(\DD v)^TY+2\mathscr T(Y\times\nabla v),\\
 \mathscr H_BY&=-2(\DD B)^TY-2\mathscr T(Y\times\nabla B),
 \end{aligned}
\end{equation}
\begin{equation}\label{mom:H2}
 \mathscr H_2\Theta=\DD(D_tv-D_BB)^T\Theta
 +2\mathscr T\bigl(((\DD v)^T\Theta)\times\nabla v
                         -((\DD B)^T\Theta)\times\nabla B\bigr).
\end{equation}
Let $\xi=\curl\Theta$ and define $w,b$ by \eqref{eq:linear-fields}.
Then
\begin{equation}\label{mom:projected-G}
 \curl\mathscr G\Theta=\mathbb P\mathcal F_{v,B}(w,b).
\end{equation}
Consequently, if $\mathsf S$ is symmetric and
$\curl\mathscr G\Theta=\mathbb P\ddiv \mathsf S$, there is a unique
mean-zero pressure $\pi$ such that
\begin{equation}\label{mom:linear-system}
 \mathcal F_{v,B}(w,b)+\nabla\pi=\ddiv \mathsf S,
\end{equation}
determined by
\begin{equation}\label{mom:pressure}
 \Delta\pi=\ddiv\ddiv \mathsf S-2\ddiv\bigl(w\cn v-b\cn B\bigr).
\end{equation}
\end{proposition}
\begin{proof}
Commutation of curl with Lie derivatives gives $\curl(\partial_t+\mathcal L_v)\Theta=w$ and
$\curl\mathcal L_B\Theta=b$. Using
$[\partial_t+\mathcal L_v,\mathcal L_B]=0$ and \eqref{mom:Kdiv},
\begin{align*}
 \curl\mathscr A^-\mathscr A^+\Theta
 &=(\partial_t+\mathcal L_v)w-\mathcal L_Bb\\
 &=\mathcal F_{v,B}(w,b)-2(w\cn v-b\cn B),
\end{align*}
\[
 \begin{aligned}
 &\curl\,2\mathscr T\left(
  (D_t\Theta+(\DD v)^T\Theta)\times\nabla v
 -(D_B\Theta+(\DD B)^T\Theta)\times\nabla B\right)\\
 &\qquad=2\mathbb P(w\cn v-b\cn B).
 \end{aligned}
\]
Apply the Leray projection to the first identity and add the second.
The left-hand side of the first identity is already divergence free,
so the sum equals $\mathbb P\mathcal F_{v,B}(w,b)$. To identify its
potential, use
\[
 \mathcal A^-((\DD z^+)^T)
 =\DD(\mathcal A^-z^+)^T-(\DD z^-)^T(\DD z^+)^T.
\]
Substituting this identity into $\mathscr A^-\mathscr A^+\Theta$
cancels the quadratic matrix product and leaves
\[
 \mathcal A^-\mathcal A^+\Theta
 +(\DD z^-)^T\mathcal A^+\Theta
 +(\DD z^+)^T\mathcal A^-\Theta
 +\DD(D_tv-D_BB)^T\Theta,
\]
and the first-order terms combine as
\[
 (\DD z^-)^T\mathcal A^+\Theta+(\DD z^+)^T\mathcal A^-\Theta
 =2(\DD v)^TD_t\Theta-2(\DD B)^TD_B\Theta.
\]
Expanding the two cross-product terms gives
\eqref{galbrun:separate-lower-operators} and \eqref{mom:H2}, proving
\eqref{mom:projected-G}.

Since $w$ has zero mean and the other force terms are divergences,
$\int\mathcal F_{v,B}(w,b)=0$. Thus
\[
 \mathbb P(\ddiv\mathsf S-\mathcal F_{v,B}(w,b))=0
\]
implies
\[
 \ddiv\mathsf S-\mathcal F_{v,B}(w,b)=\nabla\pi,
 \qquad \int\pi=0.
\]
Incompressibility gives
\[
 \ddiv\mathcal F_{v,B}(w,b)=2\ddiv(w\cn v-b\cn B),
\]
which gives \eqref{mom:pressure}. The right-hand side has zero mean,
and hence determines $\pi$ uniquely under the stated normalization.
\end{proof}

Thus a remainder in the potential equation,
$\mathscr G\Theta=\mathscr T\mathsf S+Y$, is represented in the
momentum equation by the additional stress $\mathcal R\curl Y$.
This additional stress has zero double divergence and therefore does
not affect the pressure equation. For the estimates below, observe
that the coefficient of $D_B\Theta$ contains $\nabla B$. The
background momentum equation gives
\begin{equation}\label{mom:null-acceleration}
 D_tv-D_BB=\ddiv R-\nabla p .
\end{equation}

\paragraph{\textbf{Hypotheses and estimates.}}

Throughout this appendix the scales satisfy
\begin{equation}
 \tau_a\le\tau_c,\qquad \lambda_\parallel\le\tau_c^{-1},
 \qquad \tau_a^{-1}\le\lambda_q\le\ell^{-1},
\end{equation}
together with
\[
 \ell^{-\alpha}\varepsilon_\tau\le1,\qquad
 \ell^{1-\alpha}\lambda_q\le1,\qquad
 \ell^{-\alpha}\delta_q^{1/2}\le1.
\]
These comparisons follow from the choices in
Section~\ref{sec:parameter-choice}. The sharp estimates below use
bounds for derivatives of the background and source up to order $\rprep$.

\begin{assumption}[Background and forcing assumptions]
\label{galbrun:linear-assumptions}
On a time interval of length at most $C_{\rm win}\tau_c$, assume
the following bounds in the classes of Section~\ref{ssec:fixed-scales},
defined using the background transports.
\begin{enumerate}
\item[(i)] \emph{Background estimates.} With $\rprep\ge j_1+12$,
\begin{equation}\label{galbrun:background-low}
 \begin{aligned}
 \nabla v&\in\mathcal C_{\rprep-1}(C\lambda_{q+1}^{-3\alpha}\tau_c^{-1};\lambda_q,\mathrm c),\\
 \nabla B&\in\mathcal C_{\rprep-1}(C\lambda_{q+1}^{-3\alpha}\lambda_\parallel;\lambda_q,\mathrm c),\\
 \nabla(\ddiv R-\nabla p)&\in\mathcal C_{\rprep-1}(C\lambda_{q+1}^{-6\alpha}\tau_c^{-2};\lambda_q,\mathrm c),
 \end{aligned}
\end{equation}
so that the background is $(\lambda_q,\mathrm c)$-adapted through
$(\rprep,j_1+1)$. The factors $\lambda_{q+1}^{-3\alpha}$ and
$\lambda_{q+1}^{-6\alpha}$ absorb the loss $\ell^{-\alpha}$
in the H\"older estimates and the factor $\lambda_{q+1}^{\alpha}$
for each coefficient in a commutator with a spatial multiplier
(Lemma~\ref{principal:finite-multiplier}), since
\[
 \ell^{-\alpha}\lambda_{q+1}^{-3\alpha}\lambda_{q+1}^{\alpha}
 \le\lambda_{q+1}^{-\alpha}.
\]
\item[(ii)] \emph{Ordinary derivatives of the background.} The fields are bounded,
$\|v\|_0+\|B\|_0\le C$.
We also assume the ordinary derivative estimates
\[
 \|\partial_t^h(v,B)\|_{r+\alpha}
 \lesssim\ell^{-\alpha}\lambda_q^{r+h}\delta_q^{1/2},
 \qquad h\le j_1+1,\qquad 1\le r+h\le\rprep.
\]
The estimates below hold for $k\le j_1+1$ in the first line,
$k\le j_1$ in the second, and every $r\ge0$:
\begin{align}
 \|\partial_t^k\nabla v\|_{r+\alpha}
 +\|\partial_t^k\nabla B\|_{r+\alpha}
 &\lesssim\lambda_q\ell^{-r-k-2\alpha}\delta_q^{1/2},
 \label{galbrun:crude-gradient-input}\\
 \|\partial_t^k\nabla(\ddiv R-\nabla p)\|_{r+\alpha}
 &\lesssim\lambda_q^2\ell^{-r-k-4\alpha}\delta_q.
 \label{galbrun:crude-acceleration-input}
\end{align}
\item[(iii)] \emph{Source.} The tensor source has the form
\begin{equation}
 \mathsf F(t,x)=\sum_{b\in\mathscr B}g_b(t/\tau_a)\mathsf S_b(t,x),
\end{equation}
where $\mathscr B$ is a fixed finite set and the $g_b$ are fixed smooth
mean-zero unit-periodic scalar functions. The slow tensors
$\mathsf S_b$ have a common compact time support strictly after the
initial time $t_-$ and satisfy
\begin{equation}\label{galbrun:source-bounds}
 \mathsf S_b\in\mathcal C_{\rprep}(M_{\rm src};\lambda_q,\mathrm c)
 \cap\mathcal K(M_{\rm src})\cap\mathcal O_{\infty,j_1}(M_{\rm src};\ell^{-1}).
\end{equation}
\end{enumerate}
The constants may depend on $C_{\rm win}$, the low coefficient norms
needed to estimate composition with the characteristic flows, and
the derivatives and mean-zero periodic primitives of the fixed $g_b$
used in the proofs. They are independent of the iteration stage.
The H\"older estimate for ordinary source derivatives in
\eqref{galbrun:source-bounds} reads
\begin{equation}\label{galbrun:crude-source-bounds}
 \sum_{b\in\mathscr B}\|\partial_t^k\mathsf S_b\|_{r+\alpha}
 \lesssim\ell^{-r-k-\alpha}M_{\rm src},\qquad k\le j_1,
\end{equation}
in the corresponding H\"older norms.
\end{assumption}

Proposition~\ref{prep:background-verification} verifies (i) and (ii)
for the local backgrounds through derivative order $\rprep$, and
Section~\ref{sec:adapted-perturbation} verifies (iii) for the forcing
\eqref{amplitude:centered-source} with $M_{\rm src}=C\delta_{q+1}$.

\begin{theorem}[Linear forced Galbrun equation]
\label{galbrun:linear-theorem}
Under Assumption~\ref{galbrun:linear-assumptions}, the Cauchy problem
\begin{equation}\label{galbrun:fixed-data-solve}
 \mathscr G\Theta=\mathscr T\mathsf F,
 \qquad \Theta(t_-)=D_t\Theta(t_-)=0,
\end{equation}
has a unique smooth solution on the given interval, with the following
properties. All material and magnetic derivatives are taken with
respect to the fixed background.
\begin{enumerate}
\item[(a)] \emph{Sharp bounds.} We have
\begin{equation}\label{galbrun:potential-bounds}
 \Theta\in\mathcal C_{\rprep-10,j_1-5}(C\ell^{-\alpha}\tau_a^2M_{\rm src};\lambda_q,\mathrm a),
\end{equation}
\begin{equation}\label{galbrun:first-state-bounds}
 \begin{aligned}
 D_t\Theta&\in\mathcal C_{\rprep-10,j_1-5}(C\ell^{-\alpha}\tau_aM_{\rm src};\lambda_q,\mathrm a),\\
 D_B\Theta&\in\mathcal C_{\rprep-10,j_1-5}(C\ell^{-\alpha}\tau_a^2\lambda_\parallel M_{\rm src};\lambda_q,\mathrm a).
 \end{aligned}
\end{equation}
In particular
\begin{equation}\label{galbrun:low-deformation}
 \|\nabla\curl\Theta\|_\alpha
 \lesssim\tau_a^2\lambda_q^2\ell^{-\alpha}M_{\rm src}.
\end{equation}
\item[(b)] \emph{Coarse bounds.}
\begin{equation}\label{galbrun:collarstate}
 \begin{aligned}
 \Theta&\in\mathcal C_{\rprep-2}(C\ell^{-\alpha}\tau_c^2M_{\rm src};\lambda_q,\mathrm a),\\
 D_t\Theta,\,D_B\Theta&\in\mathcal C_{\rprep-2}(C\ell^{-\alpha}\tau_cM_{\rm src};\lambda_q,\mathrm a).
 \end{aligned}
\end{equation}
\item[(c)] \emph{$\mathcal K$ bounds and ordinary derivative estimates.} For every $r\ge0$,
\begin{align}
 \|\partial_t^k\Theta\|_{r+\alpha}
 &\lesssim\tau_c^2M_{\rm src}\ell^{-r-k-\alpha},
 &&k\le j_1+2,
 \label{galbrun:crude-potential-bounds}\\
 \|\partial_t^kD_t\Theta\|_{r+\alpha}
 +\|\partial_t^kD_B\Theta\|_{r+\alpha}
 &\lesssim\tau_cM_{\rm src}\ell^{-r-k-\alpha},
 &&k\le j_1+1;
 \label{galbrun:crude-first-state-bounds}
\end{align}
The corresponding estimates in integer norms give
\[
 \begin{aligned}
 \Theta&\in\mathcal K(C\tau_c^2M_{\rm src})
       \cap\mathcal O_{\infty,j_1+2}(C\tau_c^2M_{\rm src};\ell^{-1}),\\
 D_t\Theta,D_B\Theta&\in\mathcal K(C\tau_cM_{\rm src})
       \cap\mathcal O_{\infty,j_1+1}(C\tau_cM_{\rm src};\ell^{-1}).
 \end{aligned}
\]
The bounds for compositions of
material and magnetic derivatives follow by expanding these operators
in ordinary derivatives and using assumption (ii).
\item[(d)] \emph{Localization.} Let $\tilde\eta=1$ near the common
source support, with $|\tilde\eta^{(h)}|\lesssim_h\tau_c^{-h}$ and
derivatives supported in the cut off regions before and after the
source support. Then $\Theta^{\rm c}=\tilde\eta\Theta$ is compactly supported in time, and its
additional source is
\begin{equation}\label{galbrun:cutoff-identity}
 \begin{aligned}
 \mathscr C^{\rm cut}[\tilde\eta,\Theta]
 &:=\mathscr G(\tilde\eta\Theta)-\tilde\eta\mathscr G\Theta\\
 &=\tilde\eta''\Theta+2\tilde\eta'D_t\Theta+\tilde\eta'\mathscr H_t\Theta,
 \end{aligned}
\end{equation}
which vanishes in the cut off region before the source support. The cutoff term and the
symmetric cutoff stress
$R^{\mathrm{cut}}=\mathcal R\curl\mathscr C^{\rm cut}[\tilde\eta,\Theta]$
lie in $\mathcal C_{\rprep-2}(C\ell^{-\alpha}M_{\rm src};\lambda_q,\mathrm a)$.
For every $1\le N\le j_1-2$,
\begin{equation}\label{galbrun:genericoutgoingtailA}
 R^{\mathrm{cut}}\in
 \mathcal C_{\rprep-N-7,\,j_1-N-2}(C_N\ell^{-\alpha}\varepsilon_\tau^NM_{\rm src};\lambda_q,\mathrm a),
\end{equation}
together with
\begin{equation}\label{galbrun:genericoutgoingtailT}
 \mathscr C^{\rm cut}[\tilde\eta,\Theta],\,
 R^{\mathrm{cut}}\in\mathcal K(CM_{\rm src})\cap\mathcal O_{\infty,j_1+1}(CM_{\rm src};\ell^{-1}).
\end{equation}
Both $\Theta$ and $\Theta^{\rm c}$ satisfy (a)--(c).
\item[(e)] \emph{Pressure.} The mean-zero pressures of $\Theta$ and
$\Theta^{\rm c}$ satisfy Corollary~\ref{galbrun:split-pressure}.
\item[(f)] \emph{Estimates for the iteration.} With the parameter choices of Lemmas~\ref{iter:parameter-lemma}
and \ref{aniso:parameter-budgets}, the bounds in (a) extend to all
$k+m\le j_1$ when the background derivatives are measured in
the classes used in the iteration:
\begin{equation}\label{galbrun:final-rate-classes}
 \begin{aligned}
 \Theta&\in\mathcal C_{\rprep-10}(C\ell^{-\alpha}\tau_a^2M_{\rm src};\lambda_q,\mathrm f),\\
 D_t\Theta&\in\mathcal C_{\rprep-10}(C\ell^{-\alpha}\tau_aM_{\rm src};\lambda_q,\mathrm f),\\
 D_B\Theta&\in\mathcal C_{\rprep-10}(C\ell^{-\alpha}\tau_a^2\lambda_\parallel M_{\rm src};\lambda_q,\mathrm f),\\
 R^{\mathrm{cut}}&\in\mathcal C_{\rprep-j_0-7}(C\ell^{-\alpha}\varepsilon_\tau^{j_0}M_{\rm src};\lambda_q,\mathrm f),
 \end{aligned}
\end{equation}
and likewise for $\Theta^{\rm c}$. The losses of the $\mathcal K$ bounds in
(c) and (d), compared with the estimates above, are
$\varepsilon_\tau^{-2}$, $\varepsilon_\tau^{-1}$,
$\varepsilon_\tau^{-2}\varepsilon_{q+1}^{-\gamma_\parallel}$ and $\varepsilon_\tau^{-j_0}$.
\end{enumerate}
\end{theorem}

We first solve the equation for a general smooth source. We then
use time oscillations to prove the interior estimates and estimate
the error introduced by the time cutoff.

\paragraph{\textbf{Existence and uniqueness.}}

\begin{lemma}[Vanishing initial data]\label{galbrun:zero-incoming}
For every smooth one-form $\mathbb F$, the problem
\begin{equation}\label{galbrun:galbrunzeroincoming}
 \mathscr G\Theta=\mathbb F,\qquad
 \Theta(t_0)=D_t\Theta(t_0)=0
\end{equation}
has a unique smooth solution on every compact subinterval $I\ni t_0$
of the background's time domain, and it vanishes on every interval
containing the initial time on which the source vanishes. With
$U^\pm=D_t\Theta\pm D_B\Theta$, the triple $(\Theta,U^+,U^-)$ is
the unique solution of the coupled transport system
\begin{equation}\label{galbrun:firstordergalbrun}
 \begin{gathered}
 (D_t\mp D_B)U^\pm
 =\mathbb F-\tfrac12(\mathscr H_t+\mathscr H_B)U^+
              -\tfrac12(\mathscr H_t-\mathscr H_B)U^-
              -\mathscr H_2\Theta,\\
 D_t\Theta=\tfrac12(U^++U^-),\qquad
 (\Theta,U^+,U^-)(t_0)=(0,0,0),
 \end{gathered}
\end{equation}
in which $U^+$, $U^-$ and $\Theta$ are transported by the vector
fields $v-B$, $v+B$ and $v$.
\end{lemma}
\begin{proof}
\emph{1. The half-wave system.}
We apply Corollary~\ref{high:projected-transport} to
\eqref{galbrun:firstordergalbrun}. By
\eqref{galbrun:separate-lower-operators} and \eqref{mom:H2}, the two
half-wave equations have right-hand side $\mathbb F$ plus terms
$T(cU)$. Here $U$ is one of the unknowns, $c$ is a smooth matrix
formed from $\DD v$, $\DD B$ and $\DD(D_tv-D_BB)$, and $T$ is either
the identity or $\mathscr T$ from \eqref{setup:hodge}. The latter is
bounded on $C^{r,\alpha}(\mathbb T^3)$ for every $r\ge0$, so these
terms do not differentiate the unknowns. The transport velocities
are $v$, $v-B$ and $v+B$, and the source is $(0,\mathbb F,\mathbb F)$.
The corollary gives a unique solution $(\Theta,U^+,U^-)$ on $I$.
Estimates \eqref{high:coupled-transport-low} and
\eqref{high:coupled-transport-high} apply at every spatial order:
the higher coefficient norms occur in the inhomogeneous terms and
only low norms enter the exponential. Thus the solution is smooth,
with no loss of derivatives in these estimates.
It has zero initial data and vanishes on every interval containing
the initial time on which the source vanishes.

\emph{2. The second-order problem.}
Subtract the two transport equations in
\eqref{galbrun:firstordergalbrun} and use its equation for
$D_t\Theta$. Since $[D_t,D_B]=0$,
\[
 \begin{aligned}
 D_t\left(\tfrac12(U^+-U^-)-D_B\Theta\right)
 &=\tfrac12D_B(U^++U^-)-D_BD_t\Theta=0.
 \end{aligned}
\]
Since the expression has zero initial data,
$\tfrac12(U^+-U^-)=D_B\Theta$. The equation for the potential gives
$\tfrac12(U^++U^-)=D_t\Theta$, and hence
$U^\pm=D_t\Theta\pm D_B\Theta$. Adding the two half-wave equations yields
\[
 (D_t^2-D_B^2)\Theta+\mathscr H_tD_t\Theta+\mathscr H_BD_B\Theta
 +\mathscr H_2\Theta=\mathbb F,
\]
with $\Theta(t_0)=D_t\Theta(t_0)=0$, which is
\eqref{galbrun:galbrunzeroincoming}.
Conversely, if $\Theta$ solves \eqref{galbrun:galbrunzeroincoming}
on $I$, then $U^\pm=D_t\Theta\pm D_B\Theta$ satisfy
\[
 \begin{aligned}
 (D_t\mp D_B)U^\pm
 &=(D_t^2-D_B^2)\Theta\pm[D_t,D_B]\Theta\\
 &=\mathbb F-\mathscr H_tD_t\Theta-\mathscr H_BD_B\Theta-\mathscr H_2\Theta,
 \end{aligned}
\]
which is \eqref{galbrun:firstordergalbrun} because
$D_t\Theta=\tfrac12(U^++U^-)$ and $D_B\Theta=\tfrac12(U^+-U^-)$, and
$U^\pm(t_0)=0$ since $D_B$ is spatial and $\Theta(t_0)=0$. Thus
$(\Theta,U^+,U^-)$ coincides with the solution of Step~1. This proves
uniqueness for \eqref{galbrun:galbrunzeroincoming} and the asserted
vanishing wherever the source vanishes on an interval containing the
initial time.
\end{proof}

\subsection{Interior estimates}

\paragraph{\textbf{Decomposition of the solution.}}
The source profiles oscillate on the scale $\tau_a$, while their
coefficients and the background vary on the scale $\tau_c$.
Integrating each profile twice suggests a potential of size
$\tau_a^2M_{\rm src}$. Derivatives of the coefficients produce errors,
which we cancel by adding further periodic primitives. This gives
the decomposition
\begin{equation}\label{galbrun:solution-splitting}
 \Theta=\Theta_{\mathrm{exp},N}+Z_N,\qquad
 \mathscr G\Theta_{\mathrm{exp},N}=\mathbb F+E_N,\qquad
 \mathscr GZ_N=-E_N,
\end{equation}
where $\mathbb F=\mathscr T\mathsf F$ and
$\Theta_{\mathrm{exp},N}$ is the finite sum in
Lemma~\ref{galbrun:finite-normal-form}. The remainder $Z_N$ has zero
initial data. Each cancellation contributes
$\varepsilon_\tau=\tau_a/\tau_c$ to the estimate for $E_N$.
Only the decomposition depends on the truncation order; the solution
$\Theta$ is independent of $N$.

Taking $N=2$ suffices for the interior estimates. The explicit potential
has size $\tau_a^2M_{\rm src}$, and the Cauchy estimate gives $Z_2$
the same bound because
$\tau_c^2\varepsilon_\tau^2M_{\rm src}=\tau_a^2M_{\rm src}$.
Applying a material derivative to a profile gives size
$\tau_aM_{\rm src}$. Magnetic derivatives act only on the coefficients
and cost $\lambda_\parallel$.

In either cut off region, the explicit potential vanishes and the
solution equals $Z_N$. We may therefore continue the expansion to
obtain a smaller cutoff error, provided the required coefficient and
source bounds are assumed. The derivatives of the cutoff contribute
$\tau_c^{-1}$, cancelling the factors from the Cauchy estimate and
leaving $\varepsilon_\tau^N$. Before the source support the solution
is zero by uniqueness. No vanishing between the fast principal
profiles is required.

We prove these statements in the H\"older norms $\|\cdot\|_{r+\alpha}$.
At derivative orders for which the expansion cannot be estimated from
the assumed bounds, we use the coarse Cauchy estimate. The final proof
shows that the derivative costs in the iteration absorb its losses.

\paragraph{\textbf{Coefficient and source bounds.}}

\begin{lemma}[Lower-order operator bounds]\label{mom:cell-smallness}
The separate operators satisfy
\begin{align}
 \|\mathscr H_tY\|_{r+\alpha}
 &\le C\sum_{a+b=r}\|\nabla v\|_{a+\alpha}\|Y\|_{b+\alpha},\notag\\
 \|\mathscr H_BY\|_{r+\alpha}
 &\le C\sum_{a+b=r}\|\nabla B\|_{a+\alpha}\|Y\|_{b+\alpha},
 \\
 \|\mathscr H_2Y\|_{r+\alpha}
 &\le C\sum_{a+b=r}
  \|\nabla(\ddiv R-\nabla p)\|_{a+\alpha}\|Y\|_{b+\alpha}\notag\\
 &\quad+C\sum_{a+b+c=r}
 (\|\nabla v\|_{a+\alpha}\|\nabla v\|_{b+\alpha}
  +\|\nabla B\|_{a+\alpha}\|\nabla B\|_{b+\alpha})\|Y\|_{c+\alpha}.
\end{align}
The corresponding estimates for material and magnetic derivatives
follow by Leibniz' rule and the commutator bounds for $\mathscr T$ in
Proposition~\ref{high:cz-iterated}.
Under \eqref{galbrun:background-low}, the differentiated coefficients
satisfy \eqref{galbrun:galbrunmixedcoefficientinput} below.
\end{lemma}
\begin{proof}
The product estimate and the order-zero bound for $\mathscr T$ give
\[
 \|(\DD u)^TY\|_{r+\alpha}+\|\mathscr T(Y\times\nabla u)\|_{r+\alpha}
 \le C\sum_{a+b=r}\|\nabla u\|_{a+\alpha}\|Y\|_{b+\alpha}.
\]
This gives the first two bounds by
\eqref{galbrun:separate-lower-operators}. For \eqref{mom:H2}, apply the
product estimate twice to $((\DD u)^TY)\times\nabla u$ and use
\eqref{mom:null-acceleration} for the remaining coefficient. To
estimate material and magnetic derivatives, iterate
\[
 D_t(TY)=TD_tY+[D_t,T]Y,\qquad
 D_B(TY)=TD_BY+[D_B,T]Y.
\]
For $k+m\ge1$ and $r\ge0$, Proposition~\ref{high:cz-iterated}
bounds the commutator terms arising from $D_t^kD_B^m(TY)$ in
$C^{r+\alpha}$. Each background gradient is differentiated at
most $k+m-1$ times by material and magnetic derivatives in all;
the sum of its spatial, material, and magnetic derivative orders
is at most $r+k+m-1$. These are the coefficient derivatives needed
when applying the proposition.
\end{proof}

Lemma~\ref{mom:cell-smallness}, the commutator rule of
Proposition~\ref{high:cz-iterated} and the background bounds
\eqref{galbrun:background-low} give the following factors in the coefficient
estimates after $a$ material and $b$ magnetic derivatives:
\begin{equation}\label{galbrun:galbrunmixedcoefficientinput}
 \begin{array}{c|c}
 \text{operator}&\text{coefficient factor}\\ \hline
 \mathscr H_t&\tau_c^{-1-a}\lambda_\parallel^b\\
 \mathscr H_B&\tau_c^{-a}\lambda_\parallel^{b+1}\\
 \mathscr H_2&\tau_c^{-2-a}\lambda_\parallel^b
 \end{array}
\end{equation}
In each product, the remaining derivatives act on the one-form.
The tame product estimate places the high spatial norm on one factor,
and Lemma~\ref{high:component-lie-conversion} gives the analogous Lie
derivative estimates. The small factors in
\eqref{galbrun:background-low} absorb the H\"older losses for the
coefficients, so only the factor $\ell^{-\alpha}$ from the one-form
remains. We state the resulting bounds on ranges
$r+k+m+s\le \rprep$, allowing $s$ additional spatial derivatives
of the source and coefficients. The $\mathcal K$ estimates hold for
every spatial order $r$, with $k+m\le j_1$.

\begin{lemma}[The forcing in the potential equation]
Set
\begin{equation}\label{galbrun:derivedgalbrunforcing}
 \mathbb F_{b,0}=\mathscr T\mathsf S_b,
 \qquad \mathbb F=\mathscr T\mathsf F
       =\sum_{b\in\mathscr B}g_b(t/\tau_a)\mathbb F_{b,0}.
\end{equation}
Under Assumption~\ref{galbrun:linear-assumptions}, for $k+m\le j_1$
and $r+k+m\le \rprep$,
\begin{equation}\label{galbrun:Fsharp}
 \sum_{b\in\mathscr B}\|D_t^kD_B^m\mathbb F_{b,0}\|_{r+\alpha}
 \lesssim\lambda_q^r\ell^{-\alpha}M_{\rm src}\tau_c^{-k}\lambda_\parallel^m .
\end{equation}
Every $\mathbb F_{b,0}$ is supported in the common time support of the
source coefficients, and the forcing satisfies
\begin{equation}\label{galbrun:physical-source-bounds}
 \|D_t^kD_B^m\mathbb F\|_{r+\alpha}
 \lesssim\lambda_q^r\ell^{-\alpha}M_{\rm src}\tau_a^{-k}\lambda_\parallel^m .
\end{equation}
\end{lemma}
\begin{proof}
The multiplier acts only in $x$, so it commutes with
$g_b(t/\tau_a)$ and preserves time support. Commuting $D_B$ with
$\mathscr T$ gives $[D_B,\mathscr T]=[B\cn,\mathscr T]$, while
$D_t$ gives $[D_t,\mathscr T]=[v\cn,\mathscr T]$.
Apply the commutator estimate to a term with $a,c$ material and
magnetic derivatives on the operator and $k-a,m-c$ on the tensor.
The coefficient and source bounds give
\[
 \tau_c^{-a}\lambda_\parallel^c
 \tau_c^{-(k-a)}\lambda_\parallel^{m-c}
 =\tau_c^{-k}\lambda_\parallel^m.
\]
The tame kernel estimate bounds the spatial norm with the highest
derivatives on a single factor. Summing over $\mathscr B$ proves
\eqref{galbrun:Fsharp}. Finally, Leibniz' rule gives
\[
 D_t^kD_B^m\mathbb F
 =\sum_{b\in\mathscr B}\sum_{h=0}^k\binom kh
    \tau_a^{-h}g_b^{(h)}(t/\tau_a)D_t^{k-h}D_B^m\mathbb F_{b,0},
\]
so $\tau_c^{-1}\le\tau_a^{-1}$ gives
\eqref{galbrun:physical-source-bounds}. There are no magnetic derivatives of the profiles, since they depend
only on time.
\end{proof}

\paragraph{\textbf{The explicit potential.}}

For a fixed mean-zero unit-periodic profile $g$, let $g^{[j]}$ be its
$j$th mean-zero periodic primitive, with $g^{[0]}=g$. Since
$D_tg(t/\tau_a)=\tau_a^{-1}g'(t/\tau_a)$ and $D_Bg(t/\tau_a)=0$, the
spatial lower-order operators give, for a slow one-form $F$,
\begin{equation}\label{galbrun:galbrunnormalformpreview}
 \mathscr G\{g(t/\tau_a)F\}
 =\tau_a^{-2}g''(t/\tau_a)F
 +\tau_a^{-1}g'(t/\tau_a)(2D_t+\mathscr H_t)F
   +g(t/\tau_a)\mathscr GF .
\end{equation}
Spatial Fourier multipliers commute with the profiles. Their periodic
primitives need not be compactly supported, since each product is supported
in the time support of its slow coefficient and therefore vanishes near
the initial time.

\begin{lemma}[Finite antiderivative expansion]\label{galbrun:finite-normal-form}
Under Assumption~\ref{galbrun:linear-assumptions}, fix an integer
$N\ge1$ and let $\mathbb F_{b,0}$ be defined by
\eqref{galbrun:derivedgalbrunforcing}. Set $\mathbb F_{b,-1}=0$ and define
\begin{equation}\label{galbrun:recursivegalbrunterms}
 \mathbb F_{b,n}=-(2D_t+\mathscr H_t)\mathbb F_{b,n-1}-\mathscr G\mathbb F_{b,n-2},\qquad 1\le n\le N,
\end{equation}
and
\begin{equation}\label{galbrun:physical-normal-form}
 \Theta_{\mathrm{exp},N}
 =\sum_{b\in\mathscr B}\sum_{n=0}^{N-1}
       \tau_a^{n+2}g_b^{[n+2]}(t/\tau_a)\mathbb F_{b,n}.
\end{equation}
Then
\begin{equation}
 \mathscr G\Theta_{\mathrm{exp},N}=\mathbb F+E_N,
\end{equation}
where
\begin{equation}\label{galbrun:recursivegalbrunresidual}
 E_N=\sum_{b\in\mathscr B}
 \left\{-\tau_a^Ng_b^{[N]}(t/\tau_a)\mathbb F_{b,N}
 +\tau_a^{N+1}g_b^{[N+1]}(t/\tau_a)\mathscr G\mathbb F_{b,N-1}\right\}.
\end{equation}
All these expressions have the common time support of the source
coefficients. For each $0\le n\le N$ and all nonnegative integers
$r,k,m$ satisfying $k+m+n+1\le j_1$ and
$r+k+m+n+4\le\rprep$,
\begin{equation}\label{galbrun:recursivegalbrunprofilebound}
 \sum_{b\in\mathscr B}\|D_t^kD_B^m\mathbb F_{b,n}\|_{r+\alpha}
 \lesssim_N\lambda_q^r\ell^{-\alpha}M_{\rm src}\tau_c^{-n-k}\lambda_\parallel^m,
\end{equation}
and for
\begin{equation}\label{galbrun:normal-form-range}
 k+m+N+2\le j_1,\qquad r+k+m+N+5\le \rprep
\end{equation}
the residual satisfies
\begin{equation}\label{galbrun:recursivegalbrunEprofilebound}
 \|D_t^kD_B^mE_N\|_{r+\alpha}
 \lesssim_N\lambda_q^r\ell^{-\alpha}M_{\rm src}\varepsilon_\tau^N\tau_a^{-k}\lambda_\parallel^m .
\end{equation}
Each summand of \eqref{galbrun:physical-normal-form} has size
$\tau_a^2M_{\rm src}\varepsilon_\tau^n$, with a factor $\tau_a^{-1}$ for each material derivative
and $\lambda_\parallel$ for each magnetic derivative.
\end{lemma}
\begin{proof}
Apply \eqref{galbrun:galbrunnormalformpreview} to the summand with
index $n$. The resulting terms are
\[
 \begin{gathered}
 \tau_a^ng_b^{[n]}\mathbb F_{b,n},\qquad
 \tau_a^{n+1}g_b^{[n+1]}(2D_t+\mathscr H_t)\mathbb F_{b,n},\\
 \tau_a^{n+2}g_b^{[n+2]}\mathscr G\mathbb F_{b,n},
 \end{gathered}
\]
where the profiles are evaluated at $t/\tau_a$. Sum over the indices.
For $1\le n\le N-1$, the recurrence
\eqref{galbrun:recursivegalbrunterms} cancels the terms of the same
power. The zeroth power is $\mathbb F$; at power $N$ the coefficient
is $-\mathbb F_{b,N}$ multiplied by $\tau_a^Ng_b^{[N]}$, and the term
of power $N+1$ gives the second term of $E_N$. This also holds for
$N=1$, since $\mathbb F_{b,-1}=0$. Differentiation and spatial
Fourier multipliers preserve time support, proving the support
assertion.

We prove the coefficient bound by induction in the recurrence.
In addition to the $r+k+m$ derivatives being estimated,
$2D_t+\mathscr H_t$ requires one derivative and contributes
$C\tau_c^{-1}$. For the five terms
$D_t^2$, $D_B^2$, $\mathscr H_tD_t$, $\mathscr H_BD_B$ and
$\mathscr H_2$ in $\mathscr G$, the corresponding factors are
\[
 \tau_c^{-2},\qquad \lambda_\parallel^2,\qquad \tau_c^{-2},
 \qquad \lambda_\parallel^2,\qquad \tau_c^{-2}.
\]
Each is bounded by $C\tau_c^{-2}$ and requires at most two additional
derivatives of the source and coefficients. Leibniz' rule preserves
the separate material and magnetic estimates. In each product, use
\eqref{galbrun:galbrunmixedcoefficientinput} for $D_t^aD_B^c$
of the coefficient and the induction hypothesis for
$D_t^{k-a}D_B^{m-c}$ of the one-form. The magnetic factors from
these derivatives multiply to $\lambda_\parallel^m$, in addition
to the factors displayed above for the operator itself.

It follows by induction that $2D_t+\mathscr H_t$ and $\mathscr G$
in the recurrence for $\mathbb F_{b,n}$ produce at most $n$
additional differentiations of the source and background. The
estimate thus requires at most $n$ additional derivatives, including
those of intermediate coefficients, and gives the factor
$\tau_c^{-n}$. Applying the tame product estimate at each step proves
\eqref{galbrun:recursivegalbrunprofilebound}. Only one factor requires
the high spatial norm; one further spatial derivative suffices for
the H\"older estimate.

Substitution in \eqref{galbrun:recursivegalbrunresidual} bounds the two
coefficients by $M_{\rm src}\varepsilon_\tau^N$ and
$M_{\rm src}\varepsilon_\tau^{N+1}$. The second requires at most $N+1$
additional derivatives, together with one spatial derivative for the
H\"older bound, all allowed by \eqref{galbrun:normal-form-range}.
A material derivative of a profile costs $\tau_a^{-1}$, while its
magnetic derivative vanishes. Leibniz' rule now proves
\eqref{galbrun:recursivegalbrunEprofilebound} and the bound for each
summand of the potential. The constants depend only on the fixed
primitives through order $N+1$ and finitely many of their ordinary
derivatives.
\end{proof}

\paragraph{\textbf{The remainder equation.}}
The following Cauchy estimate applies both to the remainder in
\eqref{galbrun:solution-splitting} and to the original forcing. In the
latter case it gives the coarse bounds without requiring additional
material or magnetic derivatives of the source.

\begin{lemma}[Estimates for the Cauchy problem]\label{galbrun:slow-window}
Assume the scale and background hypotheses of this appendix, including
the bound $C_{\rm win}\tau_c$ on the length of the time interval and
Assumption~\ref{galbrun:linear-assumptions}(i).
Let $s$ be a nonnegative integer denoting the additional spatial
derivatives of the source and coefficients described above, and let $M\ge0$.
Suppose $G$ is a smooth one-form that vanishes in a neighborhood of the
initial time and, for $k+m\le j_1$ and $r+k+m+s\le\rprep$, satisfies
\[
 \|D_t^kD_B^mG\|_{r+\alpha}
 \le C\lambda_q^r\ell^{-\alpha}M\tau_a^{-k}\lambda_\parallel^m .
\]
Let $Z$ be the solution of
\[
 \mathscr GZ=G,\qquad Z(t_-)=D_tZ(t_-)=0,
\]
given by Lemma~\ref{galbrun:zero-incoming}. Then, for
$k+m\le j_1$ and $r+k+m+s+2\le\rprep$,
\begin{equation}\label{galbrun:collarstate-general}
 \begin{gathered}
 \|D_t^kD_B^mZ\|_{r+\alpha}
 \lesssim\lambda_q^r\ell^{-\alpha}\tau_c^2M\tau_a^{-k}\lambda_\parallel^m,\\
 \begin{aligned}
 &\|D_t^kD_B^mD_tZ\|_{r+\alpha}
 +\|D_t^kD_B^mD_BZ\|_{r+\alpha}\\
 &\qquad\lesssim\lambda_q^r\ell^{-\alpha}\tau_cM\tau_a^{-k}\lambda_\parallel^m .
 \end{aligned}
 \end{gathered}
\end{equation}
\end{lemma}

The endpoint $k+m=j_1$ requires at most $j_1$ material and magnetic
derivatives in all on the source and coefficients. The bound for
the first magnetic derivative has size $\tau_cM$. The improved size
$\tau_c^2M\lambda_\parallel$ follows by applying the potential bound
at $(k,m+1)$ whenever $k+m+1\le j_1$ and
$r+k+m+s+3\le\rprep$.
\begin{proof}
In the half-wave formulation \eqref{galbrun:firstordergalbrun}
for $Z$, set $U^\pm=D_tZ\pm D_BZ$. Multiplying the potential
by $\tau_c^{-1}$ gives
\[
 \begin{aligned}
 D_t(\tau_c^{-1}Z)&=(2\tau_c)^{-1}(U^++U^-),\\
 (D_t\mp D_B)U^\pm
 &=G-\tfrac12(\mathscr H_t+\mathscr H_B)U^+\\
 &\quad-\tfrac12(\mathscr H_t-\mathscr H_B)U^-
       -\tau_c\mathscr H_2(\tau_c^{-1}Z).
 \end{aligned}
\]
The operators in \eqref{galbrun:separate-lower-operators} and
\eqref{mom:H2} are products with matrix coefficients, possibly followed
by a fixed spatial multiplier. We can therefore apply
Corollary~\ref{high:projected-transport} with velocities $v,v-B,v+B$.
Its coefficient norms satisfy
\[
 H_r\lesssim\tau_c^{-1}\lambda_q^r\ell^{-\alpha},\qquad
 \int_{t_-}^{t_+}H_0\,\dd t\lesssim1.
\]
Here $\lambda_\parallel\le\tau_c^{-1}$, and the normalization
replaces $\mathscr H_2$ by $\tau_c\mathscr H_2$ in the system.
Equations~\eqref{high:coupled-transport-low}--\eqref{high:coupled-transport-high}
therefore give \eqref{galbrun:collarstate-general} at $k=m=0$.
The estimates for the differentiated coefficients require the two
additional spatial derivatives stated in the lemma.

For the higher estimates, apply $D_t^kD_B^m$ to the system and argue
by induction on the sum of the two orders. The principal part is
unchanged, and every coefficient commutator contains an unknown
with fewer than $k+m$ material and magnetic derivatives in all.
By \eqref{galbrun:galbrunmixedcoefficientinput},
its coefficient contributes $\tau_c^{-1}$, cancelling the factor
$\tau_cM$ from the induction hypothesis. If $a,b$ material and
magnetic derivatives fall on that coefficient, the remaining ratio is
\[
 \frac{\tau_c^{-a}\lambda_\parallel^b
       \tau_a^{-(k-a)}\lambda_\parallel^{m-b}}
      {\tau_a^{-k}\lambda_\parallel^m}
 =\varepsilon_\tau^a\le1.
\]
The induction hypothesis for orders below $k+m$ therefore bounds
the right-hand side by
\[
 C\lambda_q^r\ell^{-\alpha}M\tau_a^{-k}\lambda_\parallel^m.
\]
The derivatives on the coefficient and the unknown sum to $r+k+m$,
so the assumptions permit the same transport estimate at each step.
All differentiated unknowns have zero initial data. Integration gives
a factor $\tau_c$ and proves the bounds for $D_tZ,D_BZ$ and, after
rescaling, for $Z$. Since $D_t\mp D_B$ remains on the left, the
estimate for $U^\pm$ requires no additional magnetic derivative of
the source or coefficients. This proves the endpoint $k+m=j_1$
using at most $j_1$ material and magnetic derivatives in all.
\end{proof}

Applied to $G=\mathbb F$ with $s=0$ and $M=M_{\rm src}$, the lemma
gives \eqref{galbrun:collarstate}.

\paragraph{\textbf{Estimates using time oscillations.}}
We first estimate the solution on the whole time interval. The
estimates in the cut off regions will then follow from the vanishing
of the explicit potential there.

\begin{proposition}[Improved Galbrun bounds]\label{galbrun:fast-interior}
The solution of \eqref{galbrun:fixed-data-solve} obeys
\eqref{galbrun:potential-bounds} and \eqref{galbrun:first-state-bounds}
on
\begin{equation}\label{galbrun:interior-range}
 k+m+5\le j_1,\qquad r+k+m+10\le \rprep;
\end{equation}
the potential bound alone holds for $k+m+4\le j_1$ and
$r+k+m+9\le \rprep$. The additional magnetic derivative needed to estimate $D_B\Theta$ accounts for the narrower range in \eqref{galbrun:interior-range}.
\end{proposition}
\begin{proof}
Take $N=2$ in \eqref{galbrun:solution-splitting}. For the remainder
$Z_2$, the source $-E_2$ has size $M_{\rm src}\varepsilon_\tau^2$
with magnetic derivative cost $\lambda_\parallel$. Applying
Lemma~\ref{galbrun:slow-window} with zero initial data gives
\[
 \|D_t^kD_B^mZ_2\|_{r+\alpha}
 \lesssim\lambda_q^r\ell^{-\alpha}
 \tau_c^2\varepsilon_\tau^2M_{\rm src}\tau_a^{-k}\lambda_\parallel^m.
\]
Since $\tau_c^2\varepsilon_\tau^2=\tau_a^2$, this is the desired
potential bound. For the first material derivative of $Z_2$, the
same lemma gives the factor $\tau_c\varepsilon_\tau^2\le\tau_a$.
For the magnetic derivative, apply the potential bound for $Z_2$
with magnetic order $m+1$, giving $\tau_a^2\lambda_\parallel$.

The two summands $n=0,1$ in $\Theta_{\mathrm{exp},2}$ satisfy the
same estimates: their size is at most $\tau_a^2M_{\rm src}$ and the
derivative costs are $\tau_a^{-1}$ and $\lambda_\parallel$ for
material and magnetic derivatives. Their first derivatives therefore
have sizes $\tau_aM_{\rm src}$ and
$\tau_a^2M_{\rm src}\lambda_\parallel$. Adding the bounds for $Z_2$
proves the assertion. With $N=2$,
\eqref{galbrun:normal-form-range} gives the residual estimate for
$k+m+4\le j_1$ and $r+k+m+7\le\rprep$. Applying
Lemma~\ref{galbrun:slow-window} requires two further spatial
derivatives of the source and coefficients.
Allowing one further derivative for the first derivatives of the
potential gives \eqref{galbrun:interior-range}.
\end{proof}

The Lie derivatives $\mathcal D_t\Theta$ and $\mathcal L_B\Theta$
satisfy the bounds for $D_t\Theta$ and $D_B\Theta$, respectively.
On one-forms,
\[
 \mathcal D_t\Theta=D_t\Theta+(\nabla v)^T\Theta,\qquad
 \mathcal L_B\Theta=D_B\Theta+(\nabla B)^T\Theta.
\]
In the first formula the added product is smaller than
$\tau_aM_{\rm src}$ by the factor $\tau_a/\tau_c$; in the second it
satisfies the magnetic bound. Lemma~\ref{high:component-lie-conversion},
with costs $\tau_a^{-1}$ and $\lambda_\parallel$, then bounds the
higher Lie derivatives. Its coefficient estimates require fewer
material and magnetic derivatives in all than the derivative being
estimated. Commuting curl with these Lie derivatives gives the
corrector field estimates. Taking two spatial derivatives of the
potential gives \eqref{galbrun:low-deformation}.

For ordinary time derivatives, the identity
$\partial_t=D_t-v\cn$, the ordinary background bounds and the bounds
for the additional spatial derivatives give, with $m=0$, for
$k+5\le j_1$ and $r+k+10\le\rprep$,
\begin{align}
 \|\partial_t^k\Theta\|_{r+\alpha}
 &\lesssim\lambda_q^{r+k}\ell^{-\alpha}\tau_a^2M_{\rm src},
 \\
 \|\partial_t^kD_t\Theta\|_{r+\alpha}
 &\lesssim\lambda_q^{r+k}\ell^{-\alpha}\tau_aM_{\rm src},\notag\\
 \|\partial_t^kD_B\Theta\|_{r+\alpha}
 &\lesssim\lambda_q^{r+k}\ell^{-\alpha}\tau_a^2M_{\rm src}\lambda_\parallel .
\end{align}
To prove these estimates, expand one $\partial_t$ at each induction
step. The resulting terms require either one spatial derivative or
one material derivative, and the ordinary background estimates bound
the coefficients. Since $\tau_a^{-1}\le\lambda_q$, both terms cost
at most $\lambda_q$. The identity
$\partial_t\Theta=D_t\Theta-v\cn\Theta$ also gives the next time
derivative of the potential with bound
$\lambda_q^{r+k+1}\ell^{-\alpha}\tau_a^2M_{\rm src}$ whenever the
estimates for the first material derivative, potential and coefficients
hold at the required orders. This last step requires no additional
magnetic derivative.

\paragraph{\textbf{Higher spatial and ordinary time derivatives.}}

\begin{lemma}[Higher ordinary derivative estimates]\label{galbrun:whole-window-lemma}
Under the background bounds (i)--(ii) of
Assumption~\ref{galbrun:linear-assumptions} and the ordinary source bound
\eqref{galbrun:crude-source-bounds}, the solution satisfies
\eqref{galbrun:crude-potential-bounds} and
\eqref{galbrun:crude-first-state-bounds}, also after multiplication by
$\tilde\eta$. Here the bound for $D_B\Theta$ is weaker than
\eqref{galbrun:first-state-bounds}.
\end{lemma}
\begin{proof}
Leibniz' rule gives
\[
 \partial_t^k\mathbb F
 =\sum_{b}\sum_{h=0}^k\binom kh\tau_a^{-h}
   g_b^{(h)}(t/\tau_a)\mathscr T(\partial_t^{k-h}\mathsf S_b).
\]
The order-zero H\"older bound for $\mathscr T$ applies directly to
\eqref{galbrun:crude-source-bounds}. Since
\[
 \tau_a^{-h}\ell^{-r-(k-h)-\alpha}\le\ell^{-r-k-\alpha},
\]
we obtain
\[
 \|\partial_t^k\mathbb F\|_{r+\alpha}
 \lesssim\ell^{-r-k-\alpha}M_{\rm src}
\]
for $0\le k\le j_1$ and every $r\ge0$. Apply the argument of
Lemma~\ref{galbrun:slow-window} with $Z=\Theta$ and $G=\mathbb F$.
For the unknowns $\tau_c^{-1}\Theta$ and
$U^\pm=D_t\Theta\pm D_B\Theta$, the coefficient norms in
\eqref{galbrun:crude-gradient-input}--\eqref{galbrun:crude-acceleration-input}
satisfy

\[
 \tau_cH_r
 \lesssim\ell^{-r}\{1+\tau_c\lambda_q\ell^{-2\alpha}\delta_q^{1/2}
                   +\tau_c^2\lambda_q^2\ell^{-4\alpha}\delta_q\}
 \lesssim\ell^{-r}.
\]
The linear term in braces is bounded by
$\varepsilon_{q+1}^{\gamma_\ell-2\gamma_S}\le1$ and the quadratic
term by its square; here $2\gamma_S<\gamma_\ell$ follows from
\eqref{iter:small-gammaS}. Estimates
\eqref{high:coupled-transport-low}--\eqref{high:coupled-transport-high}
therefore prove the claim for $k=0$. They place the higher coefficient
norms, with factor $\ell^{-r-\alpha}$, in the inhomogeneous terms,
and the exponential depends only on $\tau_cH_0$. Assumption (ii)
provides these bounds for every spatial derivative, including the
extra derivative used in the H\"older estimate.

For the time derivatives, the same assumption gives
\[
 \|\partial_t^h(v,B)\|_{r+1+\alpha}
 \lesssim\ell^{-r-h-1-\alpha}
\]
for $r\ge0$ and $0\le h\le j_1+1$. For positive spatial orders
this follows from
$\lambda_q\ell^{-\alpha}\delta_q^{1/2}\le\ell^{-1}$; for time
derivatives without spatial differentiation, use the ordinary field
bounds in (ii). The fields themselves are bounded in $C^0$. Hence
the tame product estimate
\[
 \|a\cn U\|_{r+\alpha}
 \lesssim\|a\|_0\|U\|_{r+1+\alpha}+\|a\|_{r+1+\alpha}\|U\|_0
\]
requires one additional spatial derivative and no further H\"older
loss. By the equation and $[D_t,D_B]=0$, we have
\[
 \begin{aligned}
 \partial_t\Theta&=D_t\Theta-v\cn\Theta,\\
 \partial_t(D_t\Theta)
 &=D_B^2\Theta+\mathbb F-\mathscr H_tD_t\Theta
                         -\mathscr H_BD_B\Theta\\
 &\quad-\mathscr H_2\Theta-v\cn(D_t\Theta),\\
 \partial_t(D_B\Theta)&=D_BD_t\Theta-v\cn(D_B\Theta).
 \end{aligned}
\]
Apply $k-1$ ordinary time derivatives and use the product estimate
above. The source $\mathbb F$ and the coefficient
$\DD(\ddiv R-\nabla p)^T$ in $\mathscr H_2$ are differentiated
in time at most $k-1$ times. Each transport term requires one
additional spatial derivative, while $\tau_c^{-1}\le\ell^{-1}$ and
$\tau_c^{-2}\le\ell^{-2}$ bound the lower-order terms. Induction
with $r+k$ fixed proves the assertions for $D_t\Theta$ and
$D_B\Theta$ up to $k=j_1+1$. The first identity gives the potential
bound up to $j_1+2$ without a further source derivative. Finally,
Leibniz' rule proves the same estimate for $\tilde\eta\Theta$, since
$|\tilde\eta^{(h)}|\lesssim\ell^{-h}$.
\end{proof}

\paragraph{Bounds in integer norms.}
We first prove uniform bounds for ordinary time derivatives and then
interpolate with the spatial H\"older estimates. The ordinary
background estimates give
\[
 \|\partial_t^a(v,B)\|_0\lesssim\lambda_q^a,
 \qquad
 \|\partial_t^a(v,B)\|_{r+1+\alpha}
       \lesssim\lambda_q^{r+a+1}
\]
for $0\le a\le j_1+1$, with $r+a+1\le\rprep$ in the second
estimate. At $a=0$, the first estimate is the assumed uniform bound
for the fields. Expand $\partial_t=D_t-v\cn$ in the coarse
estimates and apply the preceding tame product inequality at each
step. This gives
\[
 \begin{aligned}
 \|\partial_t^h\Theta\|_{r+\alpha}
   &\lesssim\ell^{-\alpha}\tau_c^2M_{\rm src}\lambda_q^{r+h},\\
 \|\partial_t^h(D_t\Theta,D_B\Theta)\|_{r+\alpha}
   &\lesssim\ell^{-\alpha}\tau_cM_{\rm src}\lambda_q^{r+h},
 \end{aligned}
 \qquad h\le j_1,\qquad r+h\le\rprep-4.
\]
There is no additional H\"older loss in the coefficient products,
since every positive-order ordinary derivative of a background field
contributes $\ell^{-\alpha}\delta_q^{1/2}\le1$.
To estimate one more time derivative from the equation, we use
\[
 \begin{aligned}
 \|(\partial_t^a\mathscr H_t)Y\|_{r+\alpha}
 +\|(\partial_t^a\mathscr H_B)Y\|_{r+\alpha}
 &\lesssim\tau_c^{-1}\lambda_q^a
       \bigl(\|Y\|_{r+\alpha}+\lambda_q^r\|Y\|_\alpha\bigr),\\
 \|(\partial_t^a\mathscr H_2)Y\|_{r+\alpha}
 &\lesssim\tau_c^{-2}\lambda_q^a
       \bigl(\|Y\|_{r+\alpha}+\lambda_q^r\|Y\|_\alpha\bigr).
 \end{aligned}
\]
Here $\partial_t^a\mathscr H$ means differentiation of its
coefficients, with the fixed spatial multipliers left unchanged.
The first line holds for $a\le j_1+1$, $r+a+1\le\rprep$; the
second is used for $a\le j_1$, $r+a+3\le\rprep$.
Indeed, the ordinary derivatives of the background gradients in the
first line satisfy
\[
 \begin{aligned}
 \ell^{-\alpha}\lambda_q^{r+a+1}\delta_q^{1/2}
 &=\tau_c^{-1}\lambda_q^{r+a}
       \ell^{-\alpha}\varepsilon_{q+1}^{\gamma_\ell}\\
 &\le\tau_c^{-1}\lambda_q^{r+a}
       \ell^{-\alpha}\lambda_{q+1}^{-3\alpha}.
 \end{aligned}
\]
For $\nabla(\ddiv R-\nabla p)$, expand ordinary derivatives using
the sharp material and magnetic estimates to obtain the factor
$\tau_c^{-2}\lambda_q^{r+a}
\ell^{-\alpha}\lambda_{q+1}^{-6\alpha}$.
The small factors in both cases absorb the coefficient H\"older
losses. Leibniz' rule gives the same bound for the products of
gradients in $\mathscr H_2$, and boundedness of the fixed multipliers
on H\"older spaces proves the two operator inequalities.

Expanding \eqref{galbrun:physical-source-bounds} also gives
$\|\partial_t^h\mathbb F\|_{r+\alpha}
\lesssim\ell^{-\alpha}M_{\rm src}\lambda_q^{r+h}$ for $h\le j_1$
and the spatial orders needed above. Apply $j_1$ ordinary time
derivatives to the equations for $D_t\Theta,D_B\Theta$. The
transport terms contain one more spatial derivative, already
bounded, and the coefficient inequalities control the remaining
terms: $\tau_c^{-1}$ is bounded by $\lambda_q$, and the potential
in the $\mathscr H_2$ term has an additional factor $\tau_c$.
No time derivative of the source beyond $j_1$ is needed. We obtain
the bounds for $D_t\Theta$ and $D_B\Theta$ at $h=j_1+1$, $r=0$.
Next use $\partial_t\Theta=D_t\Theta-v\cn\Theta$ to bound the
potential at $h=j_1+1$, $r\le1$, and then at $h=j_1+2$, $r=0$.
The condition $\rprep\ge j_1+12$ ensures that all coefficient and
source estimates used in these steps lie in the assumed ranges.

For $h\ge1$, these ordinary derivative estimates imply the asserted
uniform bounds, since
\[
 \ell^{-\alpha}(\ell\lambda_q)^h\le\varepsilon_{q+1}^{h\gamma_\ell-\gamma_S}\le1.
\]
At $h=0$, the sharp bounds instead give those amplitudes using
$\tau_a/\tau_c=\varepsilon_\tau$ and
$\varepsilon_\tau\ell^{-\alpha}\le1$; the bound for $D_B\Theta$
also uses $\tau_c\lambda_\parallel\le1$.
For $F=\Theta$ and $h\le j_1+2$, or for
$F=D_t\Theta,D_B\Theta$ and $h\le j_1+1$, interpolate these
uniform bounds with the ordinary H\"older estimates
\eqref{galbrun:crude-potential-bounds}--\eqref{galbrun:crude-first-state-bounds}.
For every integer $r\ge1$, this gives
\[
 \|\partial_t^hF\|_r
 \lesssim
 \|\partial_t^hF\|_0^{\alpha/(r+\alpha)}
 \|\partial_t^hF\|_{r+\alpha}^{r/(r+\alpha)}.
\]
The resulting spatial factor is $\ell^{-r}$, with the same amplitude.
This proves the ordinary derivative bounds in (c). To pass to
material and magnetic derivatives, expand the two operators in
ordinary derivatives. The background gradient estimates cost
$\ell^{-1}$ per derivative in integer norms, because
\[
 \lambda_q\ell^{1-2\alpha}\delta_q^{1/2}
 =(\ell^{1-\alpha}\lambda_q)(\ell^{-\alpha}\delta_q^{1/2})\le1.
\]
The ordinary field estimates bound the coefficient time derivatives
in the uniform norm. Leibniz' rule therefore gives
\[
 \Theta\in\mathcal K(C\tau_c^2M_{\rm src}),\qquad
 D_t\Theta,D_B\Theta\in\mathcal K(C\tau_cM_{\rm src}),
\]
with $k+m\le j_1$ and every $r\ge0$, completing (c).

\subsection{Estimates in the cut off regions}

To obtain a potential compactly supported in time, multiply the
Cauchy solution by a cutoff equal to one near the source support.
The resulting error is supported in the cut off regions before and
after that support. We estimate it using the remainder in the
antiderivative expansion.

\paragraph{\textbf{The cutoff error.}}

Let $\tilde\eta$ be as in Theorem~\ref{galbrun:linear-theorem}(d), with
\begin{equation}
 |\tilde\eta^{(h)}|\lesssim_h\tau_c^{-h}.
\end{equation}
The identity \eqref{galbrun:cutoff-identity} follows from
\eqref{galbrun:galbrunnormalformpreview} with the slow profile
$\tilde\eta$; here $D_B\tilde\eta=0$, so there is no magnetic derivative
of the time cutoff.

\begin{proposition}[Cutoff stress]\label{galbrun:fast-collar}
The cutoff term and the cutoff stress
$R^{\mathrm{cut}}=\mathcal R\curl\mathscr C^{\rm cut}[\tilde\eta,\Theta]$
vanish in the cut off region before the source support and satisfy the assertions of
Theorem~\ref{galbrun:linear-theorem}(d) on the ranges stated there.
\end{proposition}
\begin{proof}
For the first bound in (d), apply Lemma~\ref{galbrun:slow-window}
with $M=M_{\rm src}$ in the cutoff identity. Each of the three terms
has size $M_{\rm src}$ with the stated derivative costs, since cutoff
derivatives cost $\tau_c^{-1}\le\tau_a^{-1}$ and $\mathscr H_t$ has
size $\tau_c^{-1}$. The small factors in
\eqref{galbrun:background-low} absorb the coefficient H\"older loss.
Applying the commutator kernel estimate for $\mathcal R\curl$ gives
the same bound for $R^{\mathrm{cut}}$.

For \eqref{galbrun:genericoutgoingtailA}, fix $1\le N\le j_1-2$
and nonnegative integers $r,k,m$ satisfying
\[
 k+m+N+2\le j_1,\qquad r+k+m+N+7\le\rprep.
\]
All the slow coefficients in
Lemma~\ref{galbrun:finite-normal-form} vanish in both cut off
regions. The decomposition \eqref{galbrun:solution-splitting}
therefore gives
\begin{equation}
 \mathscr C^{\rm cut}[\tilde\eta,\Theta]
 =\mathscr C^{\rm cut}[\tilde\eta,Z_N]
               \quad\hbox{on }\operatorname{supp}\tilde\eta'.
\end{equation}
Before the source support the solution is zero by uniqueness. After
it, use \eqref{galbrun:normal-form-range} and
\eqref{galbrun:recursivegalbrunEprofilebound} to estimate the source
$-E_N$ with amplitude $\varepsilon_\tau^NM_{\rm src}$. These bounds
require $N+5$ additional spatial derivatives. Apply
Lemma~\ref{galbrun:slow-window} up to $j_1-N-2$ material and magnetic
derivatives in all, with $s=N+5$ and $M=\varepsilon_\tau^NM_{\rm src}$.
Its two further spatial derivatives are allowed by
$r+k+m+N+7\le\rprep$. The potential and its first material derivative
then have sizes $\tau_c^2\varepsilon_\tau^NM_{\rm src}$ and
$\tau_c\varepsilon_\tau^NM_{\rm src}$, respectively. Substitution
in the cutoff identity gives the factors
\[
 \tau_c^{-2}\tau_c^2\varepsilon_\tau^N,\qquad
 \tau_c^{-1}\tau_c\varepsilon_\tau^N,\qquad
 \tau_c^{-1}\tau_c^{-1}\tau_c^2\varepsilon_\tau^N,
\]
all equal to $\varepsilon_\tau^N$. Leibniz' rule gives the same gain
for the derivatives: cutoff material derivatives cost
$\tau_c^{-1}\le\tau_a^{-1}$, cutoff magnetic derivatives vanish,
and magnetic derivatives of $\mathscr H_t$ cost $\lambda_\parallel$
by \eqref{galbrun:galbrunmixedcoefficientinput}. The commutator
kernel estimate for $\mathcal R\curl$ now proves
\eqref{galbrun:genericoutgoingtailA}.

For \eqref{galbrun:genericoutgoingtailT}, substitute the ordinary
estimates in (c) in the cutoff identity. After $k\le j_1+1$
ordinary time derivatives, Leibniz' rule bounds every product in
$C^{r+\alpha}$ by $M_{\rm src}\ell^{-r-k-\alpha}$. Here we use
$\tau_c^{-1}\le\ell^{-1}$, $\|\mathscr H_t\|\lesssim\tau_c^{-1}$
and the coefficient estimates from
\eqref{galbrun:crude-gradient-input}. The same bound holds for the
stress by boundedness of
$\mathcal R\curl$ on H\"older spaces.

For the integer norms, first use the ordinary derivative estimates
proved above in the cutoff identity to obtain
\[
 \|\partial_t^h\mathscr C^{\rm cut}[\tilde\eta,\Theta]\|_\alpha
 +\|\partial_t^hR^{\rm cut}\|_\alpha
 \lesssim\ell^{-\alpha}M_{\rm src}\lambda_q^h,
 \qquad h\le j_1+1.
\]
At the highest time order, the ordinary gradient estimate bounds
$\partial_t^{j_1+1}\mathscr H_t$. Cutoff derivatives cost
$\tau_c^{-1}\le\lambda_q$, and $\mathcal R\curl$ commutes with
ordinary time differentiation. Since the source vanishes on the
support of the cutoff derivatives, the equations there require no
further time derivative of the forcing.

For $h\ge1$, use
$\ell^{-\alpha}(\ell\lambda_q)^h\le\varepsilon_{q+1}^{h\gamma_\ell-\gamma_S}\le1$
to obtain the uniform bound with amplitude $M_{\rm src}$. For $h=0$,
the preceding estimate with $N=1$ and
$\varepsilon_\tau\ell^{-\alpha}\le1$ gives the same bound after the
source support; before it the error vanishes. Interpolation with
the spatial H\"older estimates gives the integer spatial bounds for
every $r\ge0$ and $h\le j_1+1$. Expanding $D_t^kD_B^m$ in
ordinary derivatives for $k+m\le j_1$ then proves
\eqref{galbrun:genericoutgoingtailT}.
\end{proof}

\paragraph{\textbf{Pressure before and after localization.}}

\begin{corollary}[Pressure bounds]\label{galbrun:split-pressure}
Suppose $\mathscr G\Theta=\mathscr T\mathsf F$ with zero initial
data and the source and coefficient bounds of
Lemma~\ref{galbrun:slow-window}, with source size $M$ and $s$
additional spatial derivatives as in that lemma. Then the pressure
in \eqref{mom:linear-system} satisfies, for $k+m\le j_1$ and
$r+k+m+s+4\le \rprep$,
\begin{equation}\label{galbrun:split-pressure-bound}
 \|D_t^kD_B^m\pi\|_{r+\alpha}
 \lesssim\lambda_q^r\ell^{-\alpha}M\tau_a^{-k}\lambda_\parallel^m .
\end{equation}
The same bound holds for the
pressure associated with $\Theta^{\rm c}=\tilde\eta\Theta$. Under the
ordinary derivative hypotheses of
Lemma~\ref{galbrun:whole-window-lemma}, the same pressures satisfy
\begin{equation}\label{galbrun:ordinary-pressure-bound}
 \|\partial_t^k\pi\|_{r+\alpha}
 \lesssim M\ell^{-r-k-\alpha},\qquad
 k\le j_1,\quad r\ge0 .
\end{equation}
\end{corollary}

Estimate \eqref{galbrun:split-pressure-bound} holds up to the endpoint
of the assumed material and magnetic derivative range, without
requiring one more magnetic derivative of the source. The additional
stress in the equation for $\Theta^{\rm c}$ is
$\mathcal R\curl\mathscr C^{\rm cut}[\tilde\eta,\Theta]$; its double
divergence is zero.
\begin{proof}
On one-forms,
\[
 \mathcal D_t\Theta=D_t\Theta+(\nabla v)^T\Theta,\qquad
 \mathcal L_B\Theta=D_B\Theta+(\nabla B)^T\Theta.
\]
Commutation of curl with Lie derivatives gives $w=\curl\mathcal D_t\Theta$ and
$b=\curl\mathcal L_B\Theta$.
Equations~\eqref{mom:Kdiv} and \eqref{mom:pressure} therefore yield the exact
order-zero formula
\[
 \pi=\Delta^{-1}\ddiv\ddiv
       \{\mathsf F-2(\mathcal D_t\Theta)\times\nabla v
                    +2(\mathcal L_B\Theta)\times\nabla B\},
\]
where the inverse Laplacian acts on mean-zero functions. We estimate
the three terms inside braces before applying the order-zero
operator. Lemma~\ref{galbrun:slow-window} bounds $D_t\Theta$ and
$D_B\Theta$ with amplitude $\tau_cM$. The same bound holds for
$\mathcal D_t\Theta$ and $\mathcal L_B\Theta$, since the extra
products have sizes
\[
 \tau_c^{-1}\tau_c^2M=\tau_cM,\qquad
 \lambda_\parallel\tau_c^2M\le\tau_cM.
\]
Consequently, the products in the pressure formula have sizes bounded by
\[
 \|\nabla v\|_\alpha\tau_cM\lesssim M,\qquad
 \|\nabla B\|_\alpha\tau_cM\lesssim M.
\]
To bound material and magnetic derivatives, apply Leibniz' rule and
\eqref{galbrun:galbrunmixedcoefficientinput} to the products. The tame
estimate requires a high spatial norm of only one factor.
Proposition~\ref{high:cz-iterated} estimates the commutators with
$\Delta^{-1}\ddiv\ddiv$, proving
\eqref{galbrun:split-pressure-bound}. At $k+m=j_1$, the coarse
bound for $\mathcal L_B\Theta$ still applies without an extra
magnetic derivative of the coefficients.

For the localized potential,
\[
 \begin{aligned}
 \mathcal D_t\Theta^{\rm c}
 &=\tilde\eta\mathcal D_t\Theta+\tilde\eta'\Theta,\\
 \mathcal L_B\Theta^{\rm c}&=\tilde\eta\mathcal L_B\Theta.
 \end{aligned}
\]
The extra term has size $\tau_c^{-1}\tau_c^2M=\tau_cM$.
Leibniz' rule preserves this bound on the stated range: material
derivatives of the cutoff cost $\tau_c^{-1}\le\tau_a^{-1}$ and its
magnetic derivatives vanish. Moreover, the cutoff stress contributes
nothing to the pressure equation, because
\[
 \ddiv\ddiv(\mathcal R\curl\mathscr C^{\rm cut})
 =\ddiv\curl\mathscr C^{\rm cut}=0.
\]
Hence the pressure formula uses $\mathsf F$,
$\mathcal D_t\Theta^{\rm c}$ and $\mathcal L_B\Theta^{\rm c}$,
with no additional cutoff term.

For \eqref{galbrun:ordinary-pressure-bound}, combine the ordinary
derivative bounds for $D_t\Theta,D_B\Theta$ in
\eqref{galbrun:crude-first-state-bounds} with
\eqref{galbrun:crude-gradient-input}. This gives
\[
 \|\partial_t^k\mathcal D_t\Theta^{\rm c}\|_{r+\alpha}
 +\|\partial_t^k\mathcal L_B\Theta^{\rm c}\|_{r+\alpha}
 \lesssim\tau_cM\ell^{-r-k-\alpha},
\]
since the added products in the Lie derivatives contain the small
relative factor
\[
 \tau_c\lambda_q\ell^{-2\alpha}\delta_q^{1/2}=o(1),
\]
while $\tilde\eta'\Theta$ satisfies the same bound directly.
In the pressure formula, ordinary derivatives commute with the
spatial multipliers. The source is bounded by $M\ell^{-r-k-\alpha}$,
and each product
\[
 (\mathcal D_t\Theta^{\rm c})\times\nabla v,\qquad
 (\mathcal L_B\Theta^{\rm c})\times\nabla B
\]
has the same bound multiplied by the small factor above. All terms
in the Leibniz expansion with $r$ spatial and $k$ ordinary time
derivatives are covered by the assumed estimates for $k\le j_1$
and every $r\ge0$. The order-zero H\"older bound therefore proves \eqref{galbrun:ordinary-pressure-bound}
for both potentials.
\end{proof}

\paragraph{\textbf{Proof of the theorem.}}

\begin{proof}[Proof of Theorem~\ref{galbrun:linear-theorem}]
Existence and uniqueness follow from Lemma~\ref{galbrun:zero-incoming}.
The forcing estimates and finite antiderivative expansion used in
Proposition~\ref{galbrun:fast-interior} give
\[
 \|D_t^kD_B^m\Theta\|_{r+\alpha}
 \lesssim\lambda_q^r\ell^{-\alpha}\tau_a^2M_{\rm src}
             \tau_a^{-k}\lambda_\parallel^m
\]
and the corresponding bounds for its first material and magnetic
derivatives on \eqref{galbrun:interior-range}.
These prove (a), including \eqref{galbrun:low-deformation}, since the
H\"older estimates also bound the integer norms.
Lemma~\ref{galbrun:slow-window} applied to $\mathbb F$ gives (b), and
Lemma~\ref{galbrun:whole-window-lemma} and the integer derivative
argument following its proof give (c). The cutoff identity
and Proposition~\ref{galbrun:fast-collar} give (d). In particular,
Leibniz' rule and $|\tilde\eta^{(h)}|\lesssim\tau_c^{-h}\le\tau_a^{-h}$
show that $\tilde\eta\Theta$ inherits (a)--(c). Finally,
Proposition~\ref{mom:lower-structure} and
Corollary~\ref{galbrun:split-pressure} give (e).

For (f), fix $k+m\le j_1$ and $r+k+m\le \rprep-10$.
If $k+m\le j_1-5$, apply (a), since its derivative costs are no
larger than those in (f). If $k+m\ge j_1-4$, apply (b). Its largest
loss relative to (a) occurs for the first magnetic derivative:
\begin{equation}
 \frac{\tau_c}{\tau_a^2\lambda_\parallel}
 =\varepsilon_\tau^{-2}\varepsilon_{q+1}^{-\gamma_\parallel},
\end{equation}
where we used
$\tau_c\lambda_\parallel=\varepsilon_{q+1}^{\gamma_\parallel}$.
The potential and its material derivative lose only
$\varepsilon_\tau^{-2}$ and $\varepsilon_\tau^{-1}$,
respectively. By \eqref{aniso:corrector-final-comparison}, the larger
derivative costs in (f) absorb all three losses when $k+m\ge j_1-4$.
For the source and coefficients, this uses at most $j_1$ material
and magnetic derivatives combined.

For the cutoff error with $k+m$ material and magnetic derivatives
in all, choose
\[
 N=\min\{j_0,\max\{0,j_1-k-m-2\}\}.
\]
If $N\ge1$, then $k+m+N+2\le j_1$, so the finite antiderivative
estimate applies. If $N=0$, use the first bound in (d). Both cases
include the asserted spatial derivative range, because
\[
 \rprep-N-7\ge \rprep-j_0-7.
\]
Using \eqref{aniso:collar-depth-comparison}, the bound with factor
$\varepsilon_\tau^N$ implies the required factor
$\varepsilon_\tau^{j_0}$ with the derivative costs in (f). Finally,
comparison of (c) and (d) with \eqref{galbrun:final-rate-classes} gives
the stated losses for the $\mathcal K$ bounds.
\end{proof}
\bibliographystyle{amsplain}
\bibliography{refs_onsager}

\end{document}